\documentclass[12pt]{article}
\pdfoutput=1
\usepackage{amssymb,amsmath,indentfirst}
\usepackage{srcltx}
\usepackage{graphicx}
\usepackage[all]{xypic}
\newcounter{lemma}[section]
\newcommand{\lemma}{\par\refstepcounter{lemma}%
{\bf Lemma \arabic{section}.\arabic{lemma}.} }

\newcommand{\predlo}{\par\refstepcounter{lemma}%
{\bf Proposition \arabic{section}.\arabic{lemma}.} }
\newcommand{\theorem}{\par\refstepcounter{lemma}%
{\bf Theorem \arabic{section}.\arabic{lemma}.} }
\newcommand{\remark}{\par\refstepcounter{lemma}%
{\bf Remark \arabic{section}.\arabic{lemma}.} }

\input amssym.def
\input amssym.tex

\title{
{\bf On the complement to a real\\ caustic germ of type $E_6$}
\author{Vyacheslav D.~Sedykh}}

\begin{document}

\date{}
\maketitle

\begin{abstract}
We prove that the complement to the caustic of a stable Lagrangian map germ of type $E_6^\pm$ has seven connected components, six of which are contractible and one is homotopy equivalent to a circle. The inverse image of the noncontractible component under this map has three connected components. The restriction of the map to one of them is a two-sheeted covering; the restriction to each of the other two is a diffeomorphism. 

\medskip
{\bf Key words:} Lagrangian map, caustic, simple stable singularities, multisingularities.

\medskip
{\bf Mathematics Subject Classification:} 14P25, 14Q30, 53D12, 57R45, 58K15
\end{abstract}

\section{Introduction}

\begin{figure}[hb]
\begin{center}
\includegraphics[width=10cm]{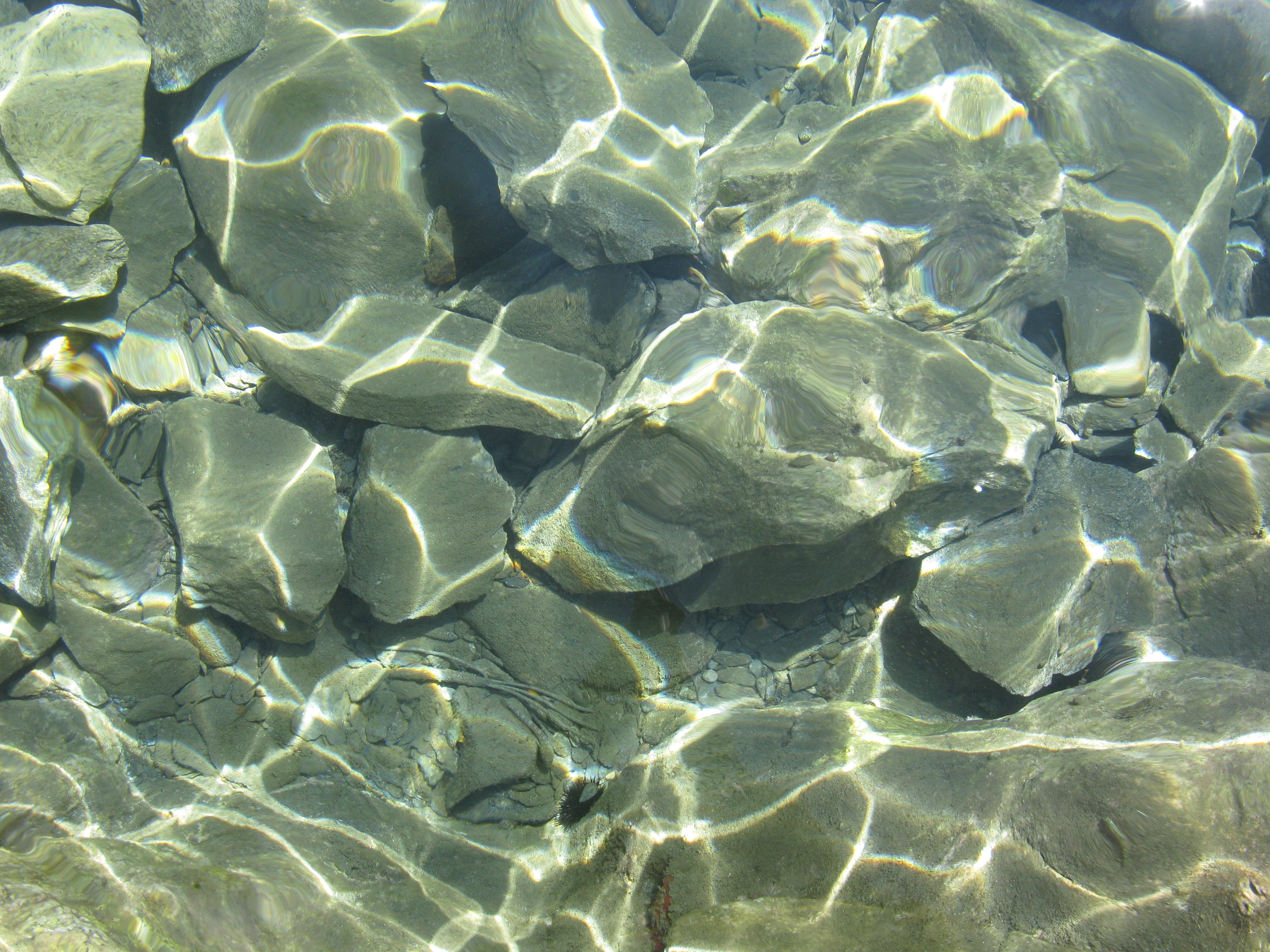}
\caption{A caustic at the bottom of the sea (Crete, 2018).}
\label{dno}
\end{center}
\end{figure}

A caustic is the set of critical values of a Lagrangian map (see \cite{Arn96}, \cite{Sed-kniga-2021}). A light caustic can be seen, for example, on a sunny day at the bottom of the sea with slight waves (Fig. \ref{dno}). Points of this caustic glow brighter than the others. Evolutes of plane curves, focal sets of hypersurfaces and other envelopes of systems of rays of different nature are also caustics.

The caustic of a generic Lagrangian map is a singular hypersurface. The singular points of this hypersurface are described by Arnold's theorem on Lagrangian singularities.
 Namely simple (that is, having zero modality) stable germs of a Lagrangian map $f:L\rightarrow V$ of a smooth manifold $L$ into a smooth manifold $V$ of the same dimension $n$ are Lagrangian equivalent to the germs at the origin of the map
\begin{equation}
\mathbb{R}^n\to \mathbb{R}^n,\quad
(t,q)\mapsto \left(
-\frac{\partial S(t,q)}{\partial t},q\right),\quad t=(t_1,\dots,t_k),\quad q=(q_{k+1},\dots,q_n)
\label{grad-map}
\end{equation}
defined by the function $S=S(t,q)$ of one of the following types corresponding to positive integers $\mu\leq n+1$:
\begin{center}
\begin{tabular}{rll}
$A_{\mu}^\pm:$ & $S=\pm t_1^{\mu+1}+q_{\mu-1}t_1^{\mu-1}+...+q_2t_1^2$, & $\mu\geq1$;\medskip \\
$D_{\mu}^{\pm}:$ & $S=t_1^2t_2\pm t_2^{\mu-1}+q_{\mu-1}t_2^{\mu-2}+...+q_3t_2^2$, & $\mu\geq4$;\medskip \\
$E_6^\pm:$ & $S=t_1^3\pm t_2^4+q_5t_1t_2^2+q_4t_1t_2+q_3t_2^2$, & $\mu=6$;\medskip \\
$E_7:$ & $S=t_1^3+t_1t_2^3+q_6t_1^2t_2+q_5t_1^2+q_4t_1t_2+q_3t_2^2$, & $\mu=7$;\medskip \\
$E_8:$ & $S=t_1^3+t_2^5+q_7t_1t_2^3+q_6t_1t_2^2+q_5t_2^3+q_4t_1t_2+q_3t_2^2$, & $\mu=8$.\\
\end{tabular}
\end{center}

The equivalence class of a Lagrangian map germ at a critical point with respect to Lagrangian equivalence is called a (Lagrangian) singularity. A generic Lagrangian map to a manifold of dimension $n\leq 5$ can have only simple stable singularities. For $n>5$ there are singularities that have functional moduli depending on $n$ variables. These singularities cannot be removed by a small (Lagrangian) deformation of the Lagrangian map.

The type of the function $S$ determines the type of a Lagrangian map germ and its singularity. The number $\mu$ is called the degree of the singularity; the number $\mu-1$ is called its codimension. If $\mu$ is even or $\mu=1$, then germs of types $A_{\mu}^{+}$ and $A_{\mu}^{-}$ are Lagrangian equivalent (their types are denoted by $A_{\mu}$). In other cases, the germs of the types listed above are pairwise Lagrangian non-equivalent.

Let $y$ be an arbitrary point of the target space $V$ of a generic proper Lagrangian map $f$ with simple stable singularities. We consider the unordered set of the
symbols from Arnold's theorem that are the types of germs of $f$ at the preimages
of $y$. The formal commutative product $\mathcal{A}$ of these symbols is called the type of the multisingularity of $f$ at the point $y$ (or the type of a monosingularity if $y$ has only one preimage). If $f^{-1}(y)=\emptyset$, then $\mathcal{A}=\mathbf{1}$. The types of multisingularities of a generic map $f$ belong to the free Abelian multiplicative semigroup $\mathbb{S}^+$ with the unity $\mathbf{1}$ and generators
$$
A_1,\quad A_{2k},A_{2k+1}^+,A_{2k+1}^- (k=1,2,\dots),\quad D_{\mu}^+,D_{\mu}^- (\mu=4,5,\dots), \quad E_6^+, E_6^-, E_7, E_8.
$$

The type of a multisingularity of a Lagrangian map $f$ at a critical value $y$ determines the type of its caustic germ at this point. Germs of a caustic are diffeomorphic if and only if their types either coincide or differ in the number of factors $A_1$ and the signs in the superscript of factors of the forms $A_{2k+1}^\pm,D_{2k+1}^\pm,E_6^\pm$. Therefore the factors $A_1$ and the superscript of the symbols of the mentioned forms are often not written if we are talking about the type of a caustic germ (but not about the type of a multisingularity of a Lagrangian map).

The set $\mathcal{A}_f$ of points $y\in V$ at which $f$ has a multisingularity of type $\mathcal{A}$ is a smooth submanifold of the space $V$. It is called the manifold of multisingularities of type $\mathcal{A}$. The codimension of $\mathcal{A}_f$ in $V$ is equal to the sum of the codimensions of singularities of $f$ at all preimages of an arbitrary point $y\in \mathcal{A}_f$. This sum is called the codimension of a multisingularity of type $\mathcal{A}$ and is denoted by $\mathrm{codim}\,\mathcal{A}$. The sum of the degrees of singularities at the preimages is called the degree of the multisingularity and is denoted by $\deg\mathcal{A}$.

Let $\delta=\pm1$. In what follows, $X_{\mu}^\delta$ means $X_{\mu}^+$ if $\delta=+1$, and $X_{\mu}^-$ if $\delta=-1$. In the paper \cite{Sed2015} we studied the topology of manifolds of multisingularities for Lagrangian germs of types $A_\mu^\delta$ and $D_\mu^\delta$. In particular from \cite[Theorem 7.1]{Sed2015}, it follows that the total number of connected components of the complement to the caustic of the map $(\ref{grad-map})$ with a singularity of type $A_{\mu}^\delta$ at the origin is equal to
$$
\left[\frac{\mu}2\right]+1,
$$
where $[\cdot]$ is the integer part of a number. All these components are contractible.

From \cite[Theorem 7.8]{Sed2015}, it follows that the total number of connected components of the complement to the caustic of the map $(\ref{grad-map})$ with a singularity of type $D_{\mu}^\delta$ at the origin is equal to:
$$
\begin{aligned}
\frac{k^2+3k-2}2& \quad\mbox{ if }\quad \mu=2k,\delta=+1,\\
\frac{k^2+k}2& \quad\mbox{ if }\quad \mu=2k,\delta=-1,\\
\frac{k^2+3k}2& \quad\mbox{ if }\quad \mu=2k+1,\delta=\pm1.\\
\end{aligned}
$$

\noindent
Each of these components is either homotopy equivalent to a circle $S^1$ or contractible. The number of non-contractible components is equal to
$$
\frac{k^2-3k+2}2,\quad \frac{k^2-k}2,\quad \frac{k^2-k}2,
$$
respectively.

We proved in \cite{Sed2023} that the complement to the image of the map $(\ref{grad-map})$ with a singularity of type $E_6^\delta$ at the origin is contractible. Now we describe the topology of the connected components of the complement to the caustic of the map $(\ref{grad-map})$ in its image.

\medskip
\theorem\label{th} {\it Let $f$ be a Lagrangian map given by $(\ref{grad-map})$. Assume that it has a singularity of type $E_6^\delta$ at the origin. Then the complement to the caustic of the map $f$ has seven connected components: two connected components of the manifold of multisingularities of each type $A_1^2,A_1^4,A_1^6$ and the complement to the image. All of them are contractible except for one connected component of type $A_1^4$. This non-contractible component is homotopy equivalent to a circle. Its inverse image under $f$ has three connected components. The restriction of $f$ to one of them is a two-sheeted covering; the restriction to each of the other two is a diffeomorphism.}

\medskip
This statement was announced in \cite{Sed-umn}. The proof is given below. We consider the case $n=5$ and coordinates $(t_1,t_2,S_3,S_4,q_5)$ in the source space where $S_3=\frac{\partial^2 S}{\partial t_2^2}, S_4=\frac{\partial^2 S}{\partial t_1\partial t_2}$. To prove Theorem \ref{th}, we study sections of the inverse image of the caustic of the map $f$ by two-dimensional surfaces $(t_1,S_4,q_5)=\mathrm{const}$. Attempts to study sections of caustic of type $E_6$ have been made previously; for example, see \cite{Cal}.

\medskip
\remark Recently V. A. Vassiliev \cite[Theorem 2]{Vas-23} calculated the number of connected components of the complements to caustic germs of types $E_7$ and $E_8$: $10$ and $15$ respectively. The question about the number of non-contractible components among them remains open. Vassiliev gives only a lower bound for this number: $\geq5$ and $\geq7$ respectively.

\medskip
\remark Theorem \ref{th} and results of the papers \cite{Sed2015}--\cite{Sed2021}, \cite{Sed2023} allowed us to discover new properties of caustics in five-dimensional spaces (see \cite{Sed2025}). Other results of various authors who have studied the topology of Lagrangian maps can be found, for example, in \cite{Vas-88}.

\section{The complement to the caustic inverse image}

Let us fix the Lagrangian map
$$
f:\mathbb{R}^{5}\to \mathbb{R}^5,\quad
(t,q)\mapsto
\left(q_1,q_2,q\right),\quad q_1=-\frac{\partial S}{\partial t_1}(t,q),\quad q_2=-\frac{\partial S}{\partial t_2}(t,q),
$$
where
$$t=(t_1,t_2),\quad q=(q_3,q_4,q_5),\quad
S=S(t,q)=t_1^3+\delta t_2^4+q_5t_1t_2^2+q_4t_1t_2+q_3t_2^2.
$$
It is quasihomogeneous and has a singularity of type $E_6^\delta$ at the origin.

Let
$$
S_3=\frac{\partial^2 S}{\partial t_2^2}(t,q)=2(6\delta t_2^2+q_5t_1+q_3),\quad
S_4=\frac{\partial^2 S}{\partial t_1\partial t_2}(t,q)=2q_5t_2+q_4\,,
$$
$$
H_1=H_1(t,q)=6t_1S_3-S_4^2,\quad
H_2=H_2(t,q)=144\delta t_1^2t_2-S_3S_4-6t_1S_4q_5,
$$
$$
H_3=H_3(t,q)=48\delta t_1^3-(S_3+2t_1q_5)^2.
$$
Then $t_1,t_2,S_3,S_4,q_5$ are new quasihomogeneous coordinates in $\mathbb{R}^{5}$ with weights
\begin{equation}
w(t_1)=4,\quad w(t_2)=3,\quad w(S_3)=w(q_3)=6,\quad w(S_4)=w(q_4)=5,\quad w(q_5)=2.
\label{vesa}
\end{equation}

\medskip
\lemma\label{raspad-E6} \rm{(\cite{Sed2018})} {\it The singularity of the map $f$ at a point $(t_1,t_2,S_3,S_4,q_5)\neq 0$ has the following type:

\medskip
{\rm1)}
$A_1$ if $H_1\neq0$;

\medskip
{\rm2)}
$A_2$ if either $t_1=0$, $S_4=0$, $S_3\neq0$ or $H_1=0$, $H_2\neq0$;

\medskip
{\rm3)}
$A_3^\pm$ if $H_1=H_2=0$, $\pm t_1H_3>0$;

\medskip
{\rm4)}
$A_4$ if $H_1=H_2=H_3=0$, $S_4\neq0$;

\medskip
{\rm5)}
$A_5^\pm$ if $t_1=\delta q_5^2/12$, $t_2=0$, $S_3=0$, $S_4=0$, $\mp \delta q_5>0$;

\medskip
{\rm6)}
$D_4^\pm$ if $t_1=0$, $S_3=0$, $S_4=0$, $\pm(q_5^3+108t_2^2)>0$;

\medskip
{\rm7)}
$D_5^\delta$ if $t_1=0$, $S_3=0$, $S_4=0$, $q_5^3+108t_2^2=0$,
$t_2\neq0$.}

\medskip
The partition $\mathfrak{S}$ of the source space $\mathbb{R}^{5}$ of the map $f$ into connected components (strata) of the manifolds of the form $\mathcal{A}^f=f^{-1}(\mathcal{A}_f),\mathcal{A}\in\mathbb{S}^+$ is a finite semialgebraic Whitney
stratification (see \cite{Mather}). Connected components of the manifold $\mathcal{A}^f$ are referred to as strata of type $\mathcal{A}$. The union $\Sigma$ of strata of codimension $1$ and higher is a hypersurface in $\mathbb{R}^{5}$. It is called the codimension $1$ skeleton of the stratification $\mathfrak{S}$.

The skeleton $\Sigma$ is the inverse image of the caustic under the map $f$. It consists of critical points of this map and non-critical preimages of its critical values. The set of critical points of $f$ is given by the equation $H_1(t,q)=0$.

Let us fix an arbitrary $(t,q)\in \mathbb{R}^5$ and assume that there is a vector $v=(v_1,v_2)\in\mathbb{R }^2\setminus\{0\}$ such that
\begin{equation}
f(t+v,q)=f(t,q).
\label{peresech-1}
\end{equation}

\medskip
\lemma\label{peresecheniya-E6} \rm{(\cite{Sed2018})} {\it The equation $(\ref{peresech-1})$ with respect to $v$ is equivalent to the system of equations
\begin{equation}
3v_1^2+q_5v_2^2+6t_1v_1+S_4v_2=0, \quad
4\delta v_2^3+2q_5v_1v_2+12\delta t_2v_2^2+S_4v_1+S_3v_2=0
\label{peresech-2}
\end{equation}
with respect to $v_1,v_2$.}

\medskip
The system (\ref{peresech-2}) defines a smooth five-dimensional submanifold $\Gamma$ of the direct product $(\mathbb{R}^2\setminus\{0\})\times\mathbb{R}^5$. The restriction of the projection
$$
\mathbb{R}^2\times\mathbb{R}^5\rightarrow\mathbb{R}^5,\quad (v;t,q)\mapsto(t,q)
$$
to $\Gamma$ is a smooth map
$$
\Psi: \Gamma\rightarrow\mathbb{R}^5.
$$
The set of critical values of this map consists of non-critical preimages of critical values of the map $f$ and critical preimages of self-intersection points of its caustic.

\medskip
\lemma\label{l-krit-psi} {\it The set of critical points of the map $\Psi$ is the intersection of $\Gamma$ with the hypersurface in $\mathbb{R}^2\times\mathbb{R}^5$ given by the equation}
\begin{equation}
6(t_1+v_1)(S_3+12\delta v_2^2+24\delta t_2v_2+2q_5v_1)=(S_4+2q_5v_2)^2.
\label{psi-1}
\end{equation}

{\sc Proof.} From Lemma \ref{raspad-E6} it follows that critical points of the map $\Psi$ satisfy the equation $H_1(t+v,q)=0$. $\Box$\footnote{Here and below the symbol $\Box$ denotes the end of the proof.}

\medskip
A critical point of the map $\Psi$ is called a fold point if the local algebra of the germ $\Psi$ at this point is two-dimensional.

\medskip
\lemma\label{vyr-krit-psi} {\it The set of critical points of the map $\Psi$ that are not fold points is the intersection of $\Gamma$ with the surface in $\mathbb{R}^2\times\mathbb{R}^5$ given by the equations $(\ref{psi-1})$ and}
\begin{equation}
144\delta (t_1+v_1)^2(t_2+v_2)=(S_3+12\delta v_2^2+24\delta t_2v_2+2q_5v_1)(S_4+2q_5v_2)+6(t_1+v_1)(S_4+2q_5v_2)q_5.
\label{psi-sborki}
\end{equation}

{\sc Proof.} Critical points of $\Psi$ that are not fold points satisfy the equation $H_2(t+v,q)=0$ by Lemma \ref{raspad-E6}. $\Box$

\medskip
Let us return to the system of equations (\ref{peresech-2}). If $v_2=0,v_1\neq0$, then this system implies
\begin{equation}
S_4=0,\quad t_1=-\frac{v_1}2.
\label{v1ne0}
\end{equation}
If $v_2\neq0$, then there exists $u\in\mathbb{R}$ such that $v_1=uv_2$. Using (\ref{peresech-2}), we get:
\begin{equation}
S_3=-4\delta v_2^2-(u(q_5-3u^2)+12\delta t_2)v_2+6t_1u^2, \quad
S_4=-(q_5+3u^2)v_2-6t_1u.
\label{peresech-4}
\end{equation}

\medskip
\lemma\label{structure-sigma} {\it The skeleton $\Sigma$ is the union of the following five subsets of $\mathbb{R}^5$:
\begin{equation}
6t_1S_3=S_4^2;
\label{crit-tochki}
\end{equation}
\begin{equation}
S_3=4q_5t_1, S_4=0;
\label{peresech-3}
\end{equation}
the image of the map
$
(u,t_1,S_4,q_5)\mapsto(t_1,t_2,S_3,S_4,q_5)
$
given by the formulas
\begin{equation}
\begin{aligned}
\delta t_2 &=\frac{(S_4+6t_1u)(q_5+3u^2)^2}{72(uS_4-t_1(q_5-3u^2))}+\frac{2\delta (S_4+6t_1u)}{3(q_5+3u^2)}-\frac{u(q_5+u^2)}{4},\\
S_3 &=\frac{(S_4+6t_1u)^2(q_5+3u^2)}{6(uS_4-t_1(q_5-3u^2))}+4\delta \frac{(S_4+6t_1u)^2}{(q_5+3u^2)^2}-2u(S_4+3t_1u);\\
\end{aligned}
\label{psi-2}
\end{equation}
the intersections of the cylinder over the Whitney umbrella $S_4^2+12q_5t_1^2=0$ with the hypersurfaces
\begin{equation}
S_3t_1^6=9\delta\left(t_1^3t_2+\frac{S_4^3}{432\delta}\right)^2+\frac{S_4^2t_1^5}{6}
\label{psi-3}
\end{equation}
and} 
\begin{equation}
S_3S_4^2+72\delta t_1^2\left(S_4t_2+2t_1^2\right)=0.
\label{psi-4}
\end{equation}

{\sc Proof.} The set of critical points of the map $f$ is given by (\ref{crit-tochki}). The system (\ref{psi-1}),(\ref{v1ne0}) in the case $v_2=0,v_1\neq0$ determines (\ref{peresech-3}).

Let $v_2\neq0$ and $q_5+3u^2\neq0$. Then the second formula of (\ref{peresech-4}) implies
$$
v_2=-\frac{S_4+6t_1u}{q_5+3u^2}.
$$
Substituting this expression for $v_2$ in (\ref{psi-1}) and the first formula of (\ref{peresech-4}), we get (\ref{psi-2}) if $t_1+uv_2\neq0$. The last is true in the considered case. Indeed, if $t_1+uv_2=0$, then $S_4+2q_5v_2=0$ by (\ref{psi-1}). These two equalities and the second formula of (\ref{peresech-4}) imply $v_2(q_5+3u^2)=0$, which contradicts the conditions $v_2\neq0$ and $q_5+3u^2\neq0$.

Now let $v_2\neq0$ and $q_5+3u^2=0$. Then $S_4=-6t_1u$ by the second formula of (\ref{peresech-4}). Hence
\begin{equation}
S_4+2q_5v_2=-6u(t_1+uv_2).
\label{ravenstvo}
\end{equation}
If $t_1+uv_2\neq0$, then (\ref{ravenstvo}),(\ref{psi-1}) and the first formula of (\ref{peresech-4}) imply:
$$
t_2=\frac{u^3}{2\delta}-\frac23v_2,\quad
S_3=4\delta v_2^2+6t_1u^2.
$$
In the case $t_1\neq0$, these formulas determine the intersection of (\ref{psi-3}) and $S_4^2+12q_5t_1^2=0$. If $t_1+uv_2=0$, then the first formula of (\ref{peresech-4}) implies
$
S_3=-4\delta v_2(v_2+3t_2).
$
In the case $t_1\neq0$, this gives a parametrization for the intersection of (\ref{psi-4}) and $S_4^2+12q_5t_1^2=0$. Finally, if $t_1=0$, then $S_4=0$ and we get a part of (\ref{crit-tochki}). $\Box$

\medskip
There is also the following description of the hypersurface $\Sigma\subset\mathbb{R}^{5}$.

\medskip
\lemma\label{proizv-psi} {\it The skeleton $\Sigma$ is the intersection of the caustic of the Lagrangian map
$$
F:\mathbb{R}^{7}\to \mathbb{R}^7,\quad
(v,t,q)\mapsto
\left(\lambda_1,\lambda_2,t,q\right),\quad
\lambda_1=-\frac{\partial \Phi}{\partial v_1}(v,t,q),\quad
\lambda_2=-\frac{\partial \Phi}{\partial v_2}(v,t,q),
$$
$$
\Phi=\Phi(v,t,q)=v_1^3+\delta v_2^4+q_5v_1v_2^2+S_4v_1v_2+4\delta t_2v_2^3+3t_1v_1^2+\frac12S_3v_2^2
$$
with the subspace $\lambda_1=\lambda_2=0$ of $\mathbb{R}^{7}$. The map $F$ has a singularity of type $E_6^\delta$ at the origin.}

{\sc Proof.} The first statement follows from (\ref{peresech-2}). The second one follows from Arnold's theory on the classification of Lagrangian singularities (see \cite{Arn96},\cite{Sed-kniga-2021}). Namely the family
$
\Phi(v,t,q)+\lambda_1v_1+\lambda_2v_2
$
of functions of $v$ with the parameter $(\lambda_1,\lambda_2,t,q)$ is a generating family of the Lagrangian map $F$. The germ of this family at zero is an $R^+$-versal deformation of the germ $v_1^3+\delta v_2^4$. Such a deformation is $R^+$-equivalent to the deformation
$
v_1^3+\delta v_2^4+q_5v_1v_2^2+S_4v_1v_2+S_3v_2^2+\lambda_1v_1+\lambda_2v_2
$
with the same base. $\Box$

\medskip
The complement $\mathbb{R}^{5}\setminus\Sigma$ to the skeleton $\Sigma$ consists of points $(t,q)$ such that $f(t,q)$ does not belong to the caustic of the map $f$. This implies:

\medskip
\lemma\label{kratnost} {\it Let $\mathcal{A}=A_1^k$, where $k\neq0$ and $\mathcal{A}_f\neq\emptyset$. Then the manifold $\mathcal{A}^f$ is the total space of a locally trivial bundle over $\mathcal{A}_f$ with fibre consisting of $k$ points.}

\medskip
The topology of connected components of the manifolds $(A_1^k)^f$ is described below. To obtain such a description, we study families of sections of the hypersurface $\Sigma$ by planes of various dimensions.

Let us consider the hyperplane $S_4=\mathrm{const}$ in $\mathbb{R}^{5}$. Its partition $\mathfrak{S}_{S_4}$ into connected components of the intersections with strata of $\mathfrak{S}$ is a finite semialgebraic Whitney strati\-fication. Connected components of the intersection with the manifold $\mathcal{A}^f$ are referred to as strata of type $\mathcal{A}$. Since the map $f$ is quasihomogeneous, we see that any hyperplane $S_4=\mathrm{const}\neq0$ transversally intersects every stratum of $\mathfrak{S}$ that does not lie entirely in the hyperplane $S_4=0$. This implies:

\medskip
\lemma\label{sech} {\it Stratifications $\mathfrak{S}_{S_4}$ for non-zero values of $S_4$ having the same sign are diffeo\-morphic.}

\medskip
Let $\Sigma_{S_4}$ be the codimension $1$ skeleton of the stratification $\mathfrak{S}_{S_4}$. The reflection
$$
\mathbb{R}^5\rightarrow\mathbb{R}^5,\quad (t_1,t_2,S_3,S_4,q_5)\mapsto (-t_1,-t_2,-S_3,-S_4,q_5)
$$
takes $\Sigma_{S_4}$ for $\delta=-1$ to $\Sigma_{-S_4}$ for $\delta=1$. Therefore, we can assume that $\delta=1$.

\medskip
{\remark Since the function $S$ does not change when the signs of the variables $t_2$ and $q_4$ are simultaneously reversed, we see that $\Sigma$ is invariant under reflection
\begin{equation}
\mathbb{R}^5\rightarrow\mathbb{R}^5,\quad (t_1,t_2,S_3,S_4,q_5)\mapsto (t_1,-t_2,S_3,-S_4,q_5).
\label{reflect}
\end{equation}
In particular, $\Sigma_{-S_4}$ can be obtained by reflection of $\Sigma_{S_4}$ with respect to the hyperplane $t_2=0$.

\medskip
Consider the family of two-dimensional planes in $\mathbb{R}^5$ given by the equation $(S_4,t_1,q_5)=\mathrm{const}$. The partition $\mathfrak{S}_{S_4,t_1,q_5}$ of such a plane into connected components of the intersections with strata of $\mathfrak{S}$ is a finite semialgebraic Whitney strati\-fication. Connected components of the intersection with the manifold $\mathcal{A}^f$ are referred to as strata of type $\mathcal{A}$.

Let $\Sigma_{S_4,t_1,q_5}$ be the codimension $1$ skeleton of the stratification $\mathfrak{S}_{S_4,t_1,q_5}$. It is a curve in $\mathbb{R}^2=\{(t_2,S_3)\}$. From Lemma \ref{raspad-E6} it follows that for every fixed $S_4$ there is an open everywhere dense set of points in $\mathbb{R}^2=\{(t_1,q_5)\}$ (generic points), for which (\ref{crit-tochki}) defines a straight line, the straight line given by (\ref{peresech-3}) does not coincide with {\rm(\ref{crit-tochki})}, and the curve $\Sigma_{S_4,t_1,q_5}$ has finitely many singular points. Moreover, the singular points of $\Sigma_{S_4,t_1,q_5}$ are from the following list:

1) semicubical cusps and transversal intersection points of two smooth branches of the curve given by formulas (\ref{psi-2});

2) points of a transversal intersection and a simple tangency (with order $1$) of a smooth branch of curve (\ref{psi-2}) with lines (\ref{crit-tochki}) and (\ref{peresech-3});

3) points of a simple tangency of line (\ref{peresech-3}) with a smooth branch of the closure of curve $(\ref{psi-2})$ (as a subset of $\mathbb{R}^2$).

These singular points are called {\it ordinary}. A point $(t_1,q_5)$ is said to be a {\it bifurcation point} for $\Sigma_{S_4,t_1,q_5}$ if in any its neighborhood there is a generic point $({\tilde t}_1,{\tilde q}_5)$, for which $\Sigma_{S_4,t_1,q_5}$ and $\Sigma_{S_4,{\tilde t}_1,{\tilde q}_5}$ have different sets of singular points. Namely either $\Sigma_{S_4,t_1,q_5}$ has a non-ordinary singular point, or a singular point of the curve $\Sigma_{S_4,{\tilde t}_1,{\tilde q}_5}$ tends to infinity as $({\tilde t}_1,{\tilde q}_5)\rightarrow(t_1,q_5)$. The set of bifurcation points for $\Sigma_{S_4,t_1,q_5}$ is called the {\it bifurcation diagram} and is denoted by $B_{S_4}$.

The bifurcation diagram $B_{S_4}$ consists of (semi)algebraic curves. Diagrams $B_{S_4}$ for nonzero values of $S_4$ having the same sign are diffeomorphic by Lemma \ref{sech}.
Since $\Sigma$ is invariant under reflection (\ref{reflect}), we see that $B_{S_4}$ does not change when the sign of $S_4$ is reversed.

\section{The bifurcation diagram $B_{S_4},S_4\neq0$}

The bifurcation diagram $B_{S_4}$ for $S_4\neq0$ is shown in Fig. \ref{razbienie} separately in each quarter of the coordinate $(\delta t_1,q_5)$-plane. The coordinate axes belong to the diagram. The other curves defining the diagram are described in Lemmas \ref{samoperesechS4-2} -- \ref{samoperesechS4-5} below. They are shown in Fig. \ref{razbienie} up to homeomorphism. This homeomorphism is a diffeomorphism everywhere except in the neighborhood of the bold point $(\ref{spoint4})$; for clarity, we made a non-zero angle between the branches $\xi_5^+$ and $\xi_5^-$ of the semicubical parabola at this point.

\begin{figure}
\begin{center}
\begin{tabular}{cc}
\includegraphics[width=8cm]{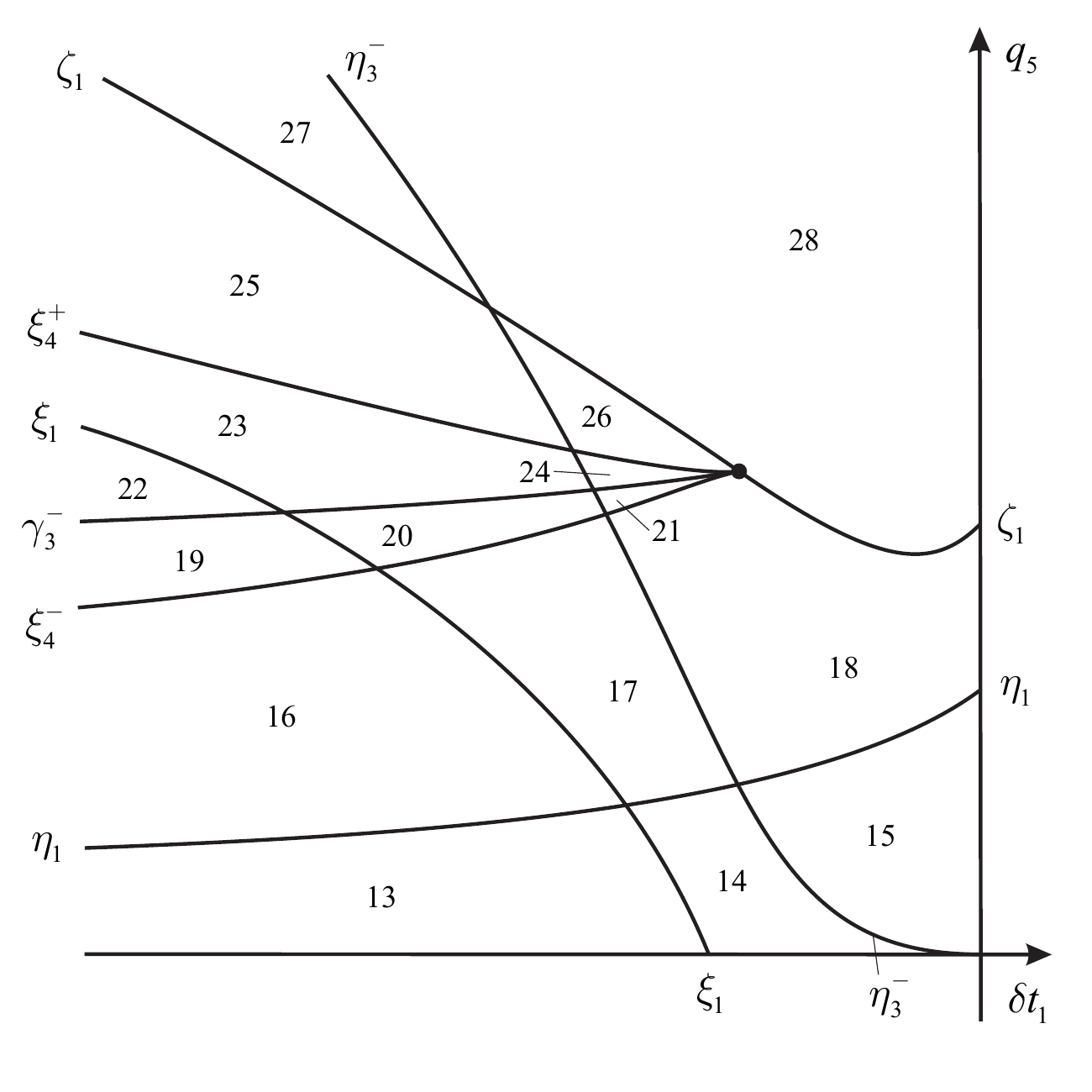}&\includegraphics[width=8cm]{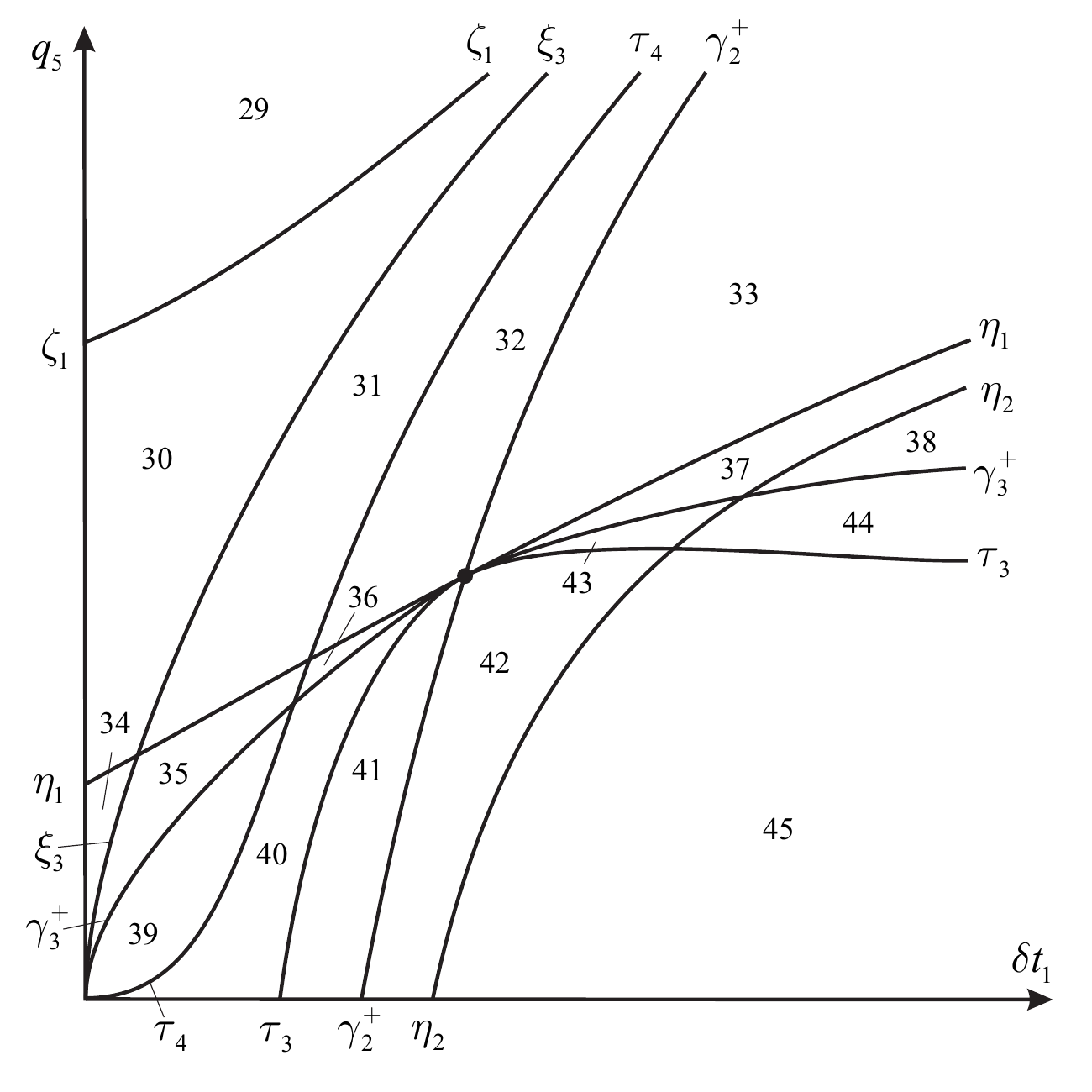}\\
$\bullet$ is point $(\ref{spoint1})$&$\bullet$ is point $(\ref{spoint2})$\\
\\
\\
\includegraphics[width=8cm]{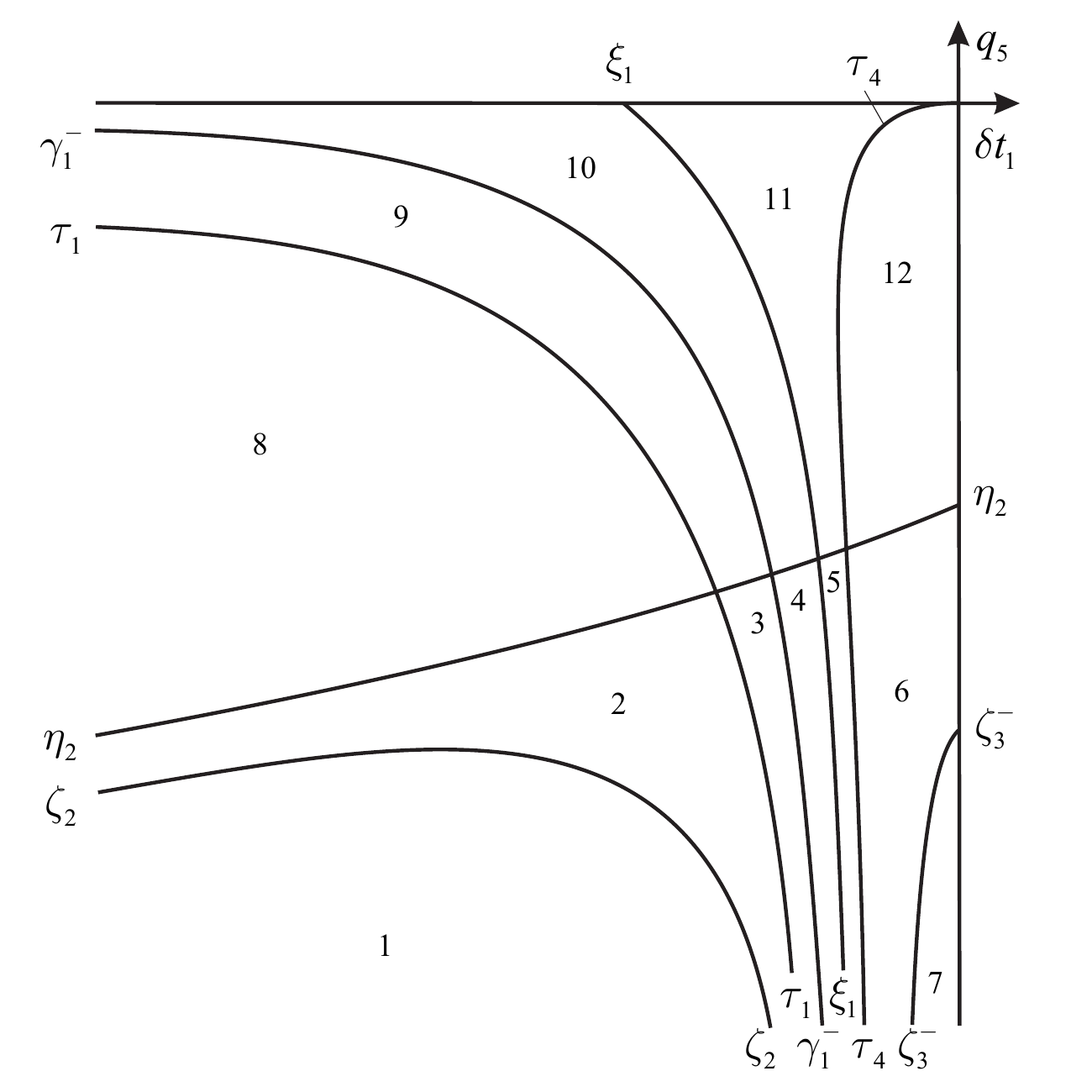}&\includegraphics[width=8cm]{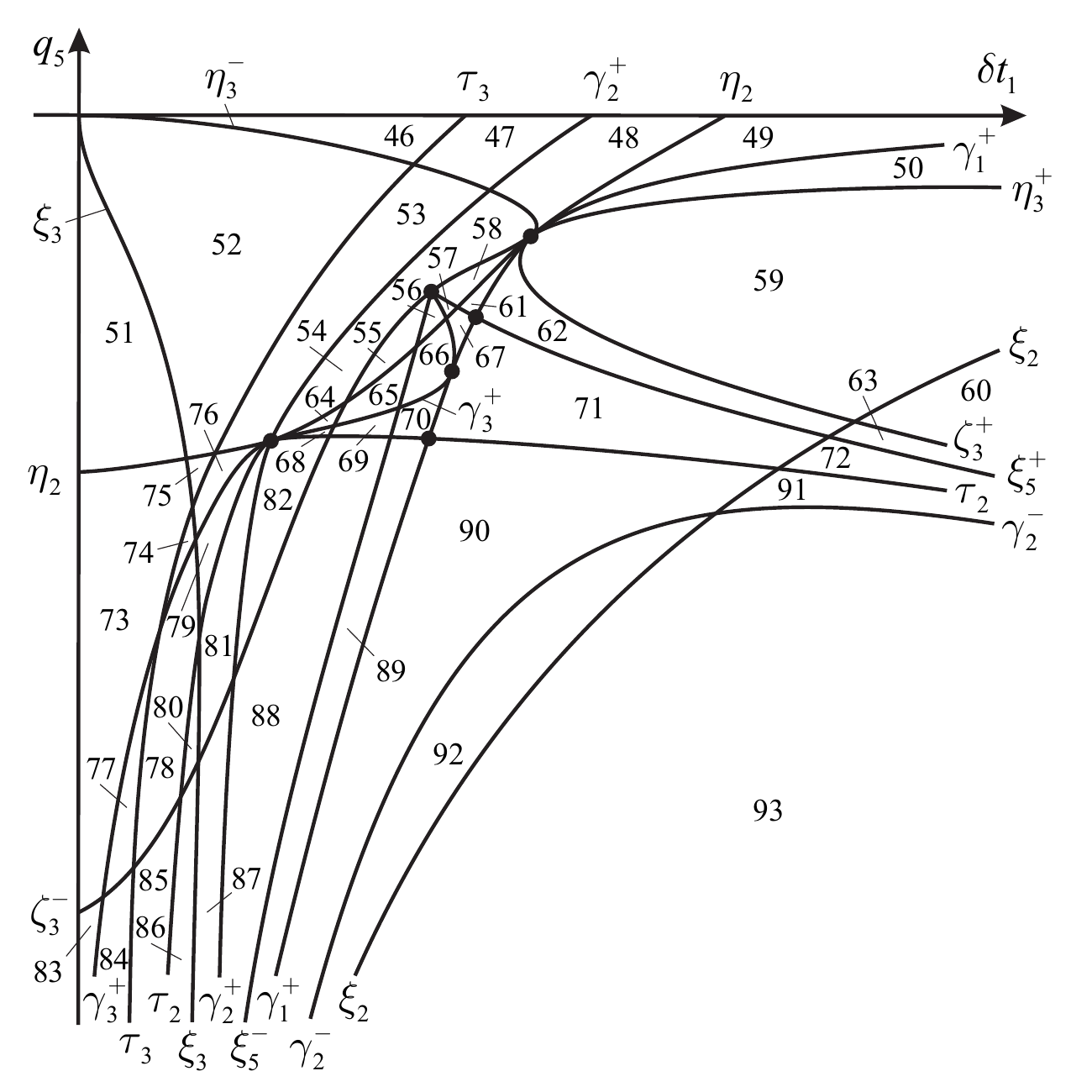}\\
&$\bullet$ are points $(\ref{spoint6})$,$(\ref{spoint8})$,$(\ref{spoint5})$,$(\ref{trper-xi5-g1})$,$(\ref{spoint4})$,$(\ref{vozvrat})$\\
&(in increasing order of their $q_5$ coordinates)\\
\\
\end{tabular}
\caption{The bifurcation diagram $B_{S_4},S_4\neq0$.}
\label{razbienie}
\end{center}
\end{figure}

\medskip
\lemma\label{samoperesechS4-2} {\it Let us fix $S_4\neq0$ and arbitrary $t_1,q_5$.

{\rm 1)} If $S_4^2+12q_5t_1^2\neq0$, then the curve $(\ref{psi-2})$ is closed {\rm(}as a subset of $\mathbb{R}^2${\rm)}. Its closure in the case
\begin{equation}
S_4^2+12q_5t_1^2=0
\label{giperbola2}
\end{equation}
is the union of the curve with the limit point
\begin{equation}
t_2=-\frac{S_4^4+1728\delta t_1^5}{432\delta t_1^3S_4},\quad S_3=\frac{S_4^4+864\delta t_1^5}{6t_1S_4^2}
\label{spec2}
\end{equation}
as $u\rightarrow-\frac{S_4}{6t_1}$.

{\rm 2)} The points at which the curve $(\ref{psi-2})$ is not immersed are determined by the roots $u\in\mathbb{R}$ of the polynomial
$$
N(u)=\left((q_5+9u^2)S_4^2-18t_1u(q_5-u^2)S_4+24q_5^2t_1^2\right)(q_5+3u^2)^3+288\delta\left(uS_4-t_1(q_5-3u^2)\right)^3.
$$

{\rm 3)} If $(t_1,q_5)$ belongs to hyperbola $(\ref{giperbola2})$ but does not coincide with the point
\begin{equation}
\delta t_1=\sqrt[5]{\frac{S_4^4}{432}},\quad q_5=-\sqrt[5]{\frac{3S_4^2}{4}},
\label{vozvrat}
\end{equation}
then the closure of the curve $(\ref{psi-2})$ is smooth at $(\ref{spec2})$. If $(t_1,q_5)\equiv(\ref{vozvrat})$, then {\rm(\ref{spec2})} is a semicubical cusp of this closure\footnote{The symbol $\equiv$ means the coincidence.}. In both cases, the line $(\ref{psi-4})$ is tangent to the closure of $(\ref{psi-2})$ at $(\ref{spec2})$.

{\rm 4)} Let $(t_1,q_5)$ be a point of the complement to the union of the following sets in $\mathbb{R}^2$: the lines $t_1=0$ and $q_5=0$; hyperbola $(\ref{giperbola2})$; the set $\zeta$ of points each of which satisfies the equation
\begin{equation}
\begin{aligned}
&144q_5^{12}t_1^4+24q_5^{11}S_4^2t_1^2+q_5^{10}(S_4^4+6912\delta t_1^5)-63360\delta q_5^9S_4^2t_1^3\\
&+240\delta q_5^8t_1(29S_4^4+432\delta t_1^5)
+3188160q_5^7S_4^2t_1^4+1070640q_5^6S_4^4t_1^2\\
&+60q_5^5S_4^2(1375S_4^4-407808\delta t_1^5)
-34560\delta q_5^4t_1^3(277S_4^4+432\delta t_1^5)\\
&-5760\delta q_5^3S_4^2t_1(125S_4^4+12528\delta t_1^5)-20736q_5^2t_1^4(2945S_4^4+6912\delta t_1^5)\\
&-15552q_5S_4^2t_1^2(625S_4^4-9984\delta t_1^5)-144(3125S_4^8-131328\delta S_4^4t_1^5+2985984t_1^{10})=0
\label{sing-points-1}
\end{aligned}
\end{equation}
and determines the polynomial $N(u)$ having a multiple real root. Then the points at which the curve {\rm(\ref{psi-2})} is not immersed are semicubical cusps. The number of cusps depends only on a connected component of the mentioned complement.}

{\sc Proof.} This is the result of direct calculations. In particular, the equations $q_5=0$ and (\ref{sing-points-1}) are determined by the discriminant of the polynomial $N(u)$ with respect to $u$. The inclusion of the line $t_1=0$ into the bifurcation diagram is due to the fact that a cusp of the curve $(\ref{psi-2})$ disappears at infinity when the point $(t_1,q_5)$ hits the $q_5$ axis. $\Box$

\medskip
\remark For any fixed $S_4\neq0$, there exists an isolated real solution $t_1=-\frac{\delta}2\sqrt[5]{9S_4^4}$, $q_5=-\sqrt[5]{3S_4^2}$ of (\ref{sing-points-1}). The polynomial $N(u)$ for these $t_1,q_5$ does not have multiple real roots.

\medskip
We denote the branches of hyperbola (\ref{giperbola2}) by
$$
\gamma_1^-: \delta t_1<0\quad \mbox{ and }\quad\gamma_1^+: \delta t_1>0.
$$

\medskip
\lemma\label{samoperesechS4-2-1} {\it Let $S_4\neq0$. Then the set $\zeta\subset\mathbb{R}^2=\{(t_1,q_5)\}$ from Lemma {\rm\ref{samoperesechS4-2}} consists of the following three curves.

{\rm 1)} The graph of a smooth positive function
$
q_5=\zeta_1(\delta t_1), t_1\in\mathbb{R}.
$
The curve $\zeta_1$ transversally intersects the $q_5$ axis at $q_5=\sqrt[5]{150(123\sqrt{5}-275)S_4^2}$.

{\rm 2)} The graph of a smooth negative function
$
q_5=\zeta_2(\delta t_1), \delta t_1<0.
$
The curve $\zeta_2$ lies in the domain\footnote{A domain in $\mathbb{R}^{n}$ is a connected open subset.} that is bounded by the curve $\gamma_1^-$ and does not contain the origin.

{\rm 3)} A curve $\zeta_3$ with one degenerate cusp {\rm (\ref{vozvrat})}. The cusp is of type $5/2$ and corresponds to the root $u=-\frac{S_4}{6t_1}$ of the polynomial $N(u)$ {\rm(}with multiplicity $7${\rm)}. The curve $\zeta_3$ is tangent to $\gamma_1^+$ at {\rm (\ref{vozvrat})}. The cusp divides $\zeta_3$ into the graphs of two smooth functions
$
\delta t_1=\zeta_3^\pm(q_5), q_5<-\sqrt[5]{\frac{3S_4^2}{4}}.
$
The curve $\zeta_3^+$ lies in the domain that is bounded by the curve $\gamma_1^+$ and does not contain the origin. The curve $\zeta_3^-$ transversally intersects the $q_5$ axis at $q_5=-\sqrt[5]{150(123\sqrt{5}+275)S_4^2}$ and lies in the domain that is bounded by {\rm(\ref{giperbola2})}.}

{\sc Proof.} All formulas are obtained by direct calculations. To determine the type of the singular point (\ref{vozvrat}) of the curve $\zeta_3$ it is necessary to make the following change of coordinates:
$$
t_1\mapsto t_1+\delta \sqrt[5]{\frac{S_4^4}{432}}
, \quad q_5\mapsto q_5-\sqrt[5]{\frac{3S_4^2}{4}}+2\delta\sqrt[5]{\frac{324}{S_4^2}}\,t_1
-\frac{18}{S_4}\sqrt[5]{\frac{18}{S_4}}\,t_1^2.
$$
The equation (\ref{sing-points-1}) in new coordinates is equivalent to
\begin{equation}
q_5^2+\frac{32768\delta}{3S_4^3}\sqrt[5]{\frac{2}{27S_4}}\,t_1^5+\dots=0,
\label{puiso1}
\end{equation}
where dots denote a polynomial in $t_1,q_5$ consisting of monomials with quasi-degree greater than $10$ in the filtration $w(t_1)=2,w(q_5)=5$. Puiseux series in $t_1$ for a germ of the complex extension of the curve (\ref{puiso1}) at zero starts at $t_1^{5/2}$ (see \cite{Shaf}). $\Box$

\medskip
\lemma\label{samoperesechS4-3} {\it Let us fix $S_4\neq0$ and arbitrary $t_1,q_5$.

{\rm 1)} The increments of $t_2$ and $S_3$ in $(\ref{psi-2})$ corresponding to the increment $\Delta u\neq0$ of $u$ are rational functions of $u$,$t_1$,$S_4$,$q_5$,$\Delta u$; the numerators $\Delta_{t_2}$ and $\Delta_{S_3}$ of these functions are divisible by $\Delta u$. The resultant of the polynomials $\frac{1}{\Delta u}\Delta_{t_2}$ and $\frac{1}{\Delta u}\Delta_{S_3}$ with respect to $\Delta u$ is a polynomial $P$ in $u$; its coefficients are polynomials in $t_1,S_4,q_5$. The degree of the polynomial $P(u)$ is equal to $18$ if $q_5\neq0$; $16$ if $q_5=0,t_1\neq0$; and $15$ if $t_1=q_5=0$. Each self-intersection point of the curve $(\ref{psi-2})$ defines the pair of different real roots of the polynomial $P(u)$.

{\rm 2)} Let $(t_1,q_5)$ be a point of the complement to the union of the following sets in $\mathbb{R}^2$: the line $q_5=0$ and hyperbola $(\ref{giperbola2})$; the set $\eta$ of points each of which satisfies the equation
\begin{equation}
\begin{aligned}
&4096q_5^{21}+221184\delta q_5^{19}t_1+4478976q_5^{17}t_1^2+1009152q_5^{16}S_4^2+38651904\delta q_5^{15}t_1^3
\\
&-4962816\delta q_5^{14}S_4^2t_1+76329216q_5^{13}t_1^4+148210560q_5^{12}S_4^2t_1^2
\\
&+3888q_5^{11}(1771S_4^4-167616\delta t_1^5)-
 621831168\delta q_5^{10}S_4^2t_1^3
\\
&-15552\delta q_5^9t_1(1585S_4^4+84096\delta t_1^5)+1285372800q_5^8S_4^2t_1^4
\\
&+15552q_5^7t_1^2(11485S_4^4+539136\delta t_1^5)+17496q_5^6S_4^2(625S_4^4-170016\delta t_1^5)
\\
&-559872\delta q_5^5t_1^3(1705S_4^4+20736\delta t_1^5)-104976\delta q_5^4S_4^2t_1(625S_4^4+72192\delta t_1^5)
\\
&-104976q_5^3t_1^4(11335S_4^4-49152\delta t_1^5)-962280q_5^2S_4^2t_1^2(125S_4^4+1152\delta t_1^5)
\\
&-2187q_5S_4^4(3125S_4^4+396672\delta t_1^5)-23328\delta S_4^2t_1^3(3125S_4^4-23328\delta t_1^5)=0
\end{aligned}
\label{per-points-1}
\end{equation}
and determines the polynomial $P(u)$ having two different real roots corresponding to a self-intersection point of curve $(\ref{psi-2})$, moreover at least one of these roots is multiple. Then the number of self-intersection points of this curve depends only on a connected component of the mentioned complement.}

{\sc Proof.} The polynomial $P(u)$ is calculated explicitly but it is too large to present here. The equations $q_5=0$, (\ref{giperbola2}) and (\ref{per-points-1}) are determined by the discriminant of this polynomial with respect to $u$. In particular, some self-intersection points of curve $(\ref{psi-2})$ disappears at infinity when the point $(t_1,q_5)$ hits the $t_1$ axis. $\Box$

\medskip
\lemma\label{samoperesechS4-3-1} {\it Let $S_4\neq0$. Then the set $\eta\subset\mathbb{R}^2=\{(t_1,q_5)\}$ from Lemma {\rm\ref{samoperesechS4-3}} consists of the following three curves.

{\rm 1)} The graph of a smooth positive function
$
q_5=\eta_1(\delta t_1), t_1\in\mathbb{R}.
$
The curve $\eta_1$ transversally intersects the $q_5$ axis at
$$
q_5=\frac12\sqrt[5]{3\left(\sqrt{3\left(272763 + 111368\sqrt{6}\right)}-657-252\sqrt{6}\right)S_4^2}.
$$

{\rm 2)} The graph of a smooth function
$
q_5=\eta_2(\delta t_1), t_1\in\mathbb{R}.
$
The curve $\eta_2$ transversally intersects the $t_1$ axis at $\delta t_1=\frac{5}{6}\sqrt[5]{\frac{S_4^4}{3}}$ and the $q_5$ axis at
$$
q_5=-\frac12\sqrt[5]{3\left(\sqrt{3\left(272763 + 111368\sqrt{6}\right)}+657+252\sqrt{6}\right)S_4^2}.
$$

{\rm 3)} The curve $\eta_3$ with one degenerate cusp {\rm (\ref{vozvrat})} of type $5/2$. The curve $\eta_3$ is tangent to the curves $\gamma_1^+$ and $\eta_2$ at {\rm(\ref{vozvrat})}. The cusp divides $\eta_3$ into the graph of a smooth positive function
$
\delta t_1=\eta_3^+(q_5), -\sqrt[5]{\frac{3S_4^2}{4}}<q_5<0
$
and the graph of a function
$
\delta t_1=\eta_3^-(q_5), q_5>-\sqrt[5]{\frac{3S_4^2}{4}}.
$
The curve $\eta_3^+$ lies in the domain that is bounded by the curve $\gamma_1^+$ and does not contain the origin. The curve $\eta_3^-$ is smooth and lies in the domain that is bounded by the curves $\gamma_1^-$ and $\eta_2$. It is tangent to the $t_1$ axis at zero {\rm(}with order $2${\rm)} and $q_5\eta_3^-(q_5)<0$ for all $q_5>-\sqrt[5]{\frac{3S_4^2}{4}}, q_5\neq0$.}

{\sc Proof.} This is a result of direct calculations as well. To determine the type of the singular point (\ref{vozvrat}) of the curve (\ref{per-points-1}) it is necessary to make the following change of coordinates:
$$
t_1\mapsto t_1+\delta \sqrt[5]{\frac{S_4^4}{432}}
, \quad q_5\mapsto q_5-\sqrt[5]{\frac{3S_4^2}{4}}+2\delta\sqrt[5]{\frac{324}{S_4^2}}\,t_1
-\frac{18}{S_4}\sqrt[5]{\frac{18}{S_4}}\,t_1^2
-\frac{1296\delta}{S_4^2}t_1^3.
$$
The equation (\ref{per-points-1}) in new coordinates is equivalent to
\begin{equation}
q_5^3-\delta q_5t_1^5+t_1^8+\dots=0,
\label{puiso2}
\end{equation}
where dots denote a polynomial in $t_1,q_5$ consisting of monomials with quasi-degree greater than $16$ in the filtration $w(t_1)=2,w(q_5)=5$. A germ (\ref{puiso2}) at zero is diffeomorphic to a germ of the curve
\begin{equation}
q_5^3-\delta q_5t_1^5+t_1^8=0
\label{puiso2-1}
\end{equation}
by Arnold's theorem on the normal form of a semi-quasihomogeneous function (see.~\cite[Theorem 7.2]{Arn74}). The curve (\ref{puiso2-1}) has two real branches at zero. One of them is smooth and the other one is singular. The smooth branch is $\eta_2$ and the singular one is $\eta_3$. Puiseux series in $t_1$ for a germ of the complex extension of the curve $\eta_3$ at zero starts at $t_1^{5/2}$. $\Box$

\medskip
\lemma\label{samoperesechS4-4} {\it Let us fix $S_4\neq0$ and arbitrary $t_1,q_5$.

{\rm 1)} The common points of the line {\rm(\ref{crit-tochki})} and the curve $(\ref{psi-2})$ are determined by the roots $u\in\mathbb{R}$ of the polynomial
$$
Q(u)=(uS_4-2q_5t_1)(q_5+3u^2)^2-24\delta t_1(uS_4-t_1(q_5-3u^2))
$$
and the value $u=-\frac{S_4}{6t_1}$.

{\rm 2)} If $(t_1,q_5)\notin(\ref{giperbola2})$, then {\rm(\ref{crit-tochki})} is tangent to {\rm(\ref{psi-2})} at the point
\begin{equation}
t_2=\frac{S_4(S_4^2+36q_5t_1^2)}{864\delta t_1^3},\quad S_3=\frac{S_4^2}{6t_1}
\label{spec1}
\end{equation}
as $u=-\frac{S_4}{6t_1}$. The order of this tangency is equal to $1$ if $(S_4^2+12q_5t_1^2)^2\neq1728\delta t_1^5$; in the case
\begin{equation}
(S_4^2+12q_5t_1^2)^2=1728\delta t_1^5,
\label{vozvrat-points}
\end{equation}
the point {\rm(\ref{spec1})} is a semicubical cusp of the curve {\rm(\ref{psi-2})}.

{\rm 3)} If $(t_1,q_5)\in(\ref{giperbola2})$, then parabola given by {\rm(\ref{psi-3})} is tangent to {\rm(\ref{crit-tochki})} at {\rm(\ref{spec1})}. The straight line given by {\rm(\ref{psi-4})} and parabola {\rm(\ref{psi-3})} are tangent to the closure of the curve {\rm(\ref{psi-2})} at {\rm(\ref{spec2})}. This tangency is simple if $(t_1,q_5)\not\equiv(\ref{vozvrat})$.

{\rm 4)} Let $(t_1,q_5)$ be a point of the complement to the union of the following sets in $\mathbb{R}^2$: the lines $t_1=0$ and $q_5=0$; hyperbola {\rm(\ref{giperbola2})}; the curves {\rm(\ref{vozvrat-points})} and
\begin{equation}
\begin{aligned}
&331776q_5^5t_1^8+69120q_5^4S_4^2t_1^6+432q_5^3t_1^4(11S_4^4-9216\delta t_1^5)\\
&+120q_5^2S_4^2t_1^2(S_4^4-10368\delta t_1^5)+q_5S_4^4(S_4^4-51840\delta t_1^5)
-288\delta S_4^2t_1^3(S_4^4-23328\delta t_1^5)=0;
\label{dobavka2}
\end{aligned}
\end{equation}
the set $\xi$ of points each of which satisfies the equation
\begin{equation}
\begin{aligned}
&1728q_5^9t_1^6+432q_5^8S_4^2t_1^4+36q_5^7t_1^2(S_4^4-432\delta t_1^5)+q_5^6S_4^2(S_4^4+10368\delta t_1^5)\\
&+108\delta q_5^5t_1^3(S_4^4-540\delta t_1^5)-36\delta q_5^4S_4^2t_1(2S_4^4+4239\delta t_1^5)+648q_5^3t_1^4(13S_4^4-72\delta t_1^5)\\
&+432q_5^2S_4^2t_1^2(4S_4^4-945\delta t_1^5)-202176\delta q_5S_4^4t_1^5-432\delta S_4^2t_1^3(32S_4^4+729\delta t_1^5)=0
\label{per-points-3}
\end{aligned}
\end{equation}
and determines the polynomial $Q(u)$ having a multiple real root; the set $\gamma_3$ of points each of which satisfies the equation
\begin{equation}
\delta t_1(q_5^5-12S_4^2)+2q_5^2S_4^2=0
\label{dobavka1}
\end{equation}
and determines the polynomials $P(u)$ {\rm(}from Lemma $\ref{samoperesechS4-3}${\rm)} and $Q(u)$ having two common real roots corresponding to a self-intersection point of the curve $(\ref{psi-2})$. Then the number of common points of {\rm(\ref{crit-tochki})} and $(\ref{psi-2})$ depends only on a connected component of the mentioned complement.}

{\sc Proof.} The equations $q_5=0$, (\ref{dobavka2}) and (\ref{dobavka1}) are determined by the resultant of the polynomials $P(u)$ and $Q(u)$ with respect to $u$. The equation (\ref{dobavka2}) follows also from the equality $P(u)=0$ for $u=-\frac{S_4}{6t_1}$. The equation (\ref{per-points-3}) is determined by the discriminant of $Q(u)$ with respect to $u$. $\Box$

\medskip
The next two statements are the result of direct calculations.

\medskip
\lemma\label{samoperesechS4-4-0} {\it The equation {\rm(\ref{vozvrat-points})} for $S_4\neq0$ determines the graphs of two smooth functions
$
\gamma_2^\pm: q_5=\pm\sqrt{12\delta t_1}-\frac{S_4^2}{12t_1^2}, \delta t_1>0.
$
The curve $\gamma_2^+$ transversally intersects the $t_1$ axis at $\delta t_1=\sqrt[5]{\frac{S_4^4}{1728}}$; the curve $\eta_1$ at the point
\begin{equation}
\delta t_1=\frac{\sqrt[5]{144\left(\sqrt{5}-2\right)^2S_4^4}}{6\left(3-\sqrt{5}\right)},\quad q_5=\sqrt[5]{12 \left(\sqrt{5}-2\right)S_4^2};
\label{spoint2}
\end{equation}
and the curve $\eta_2$ at the point}
\begin{equation}
\delta t_1=\frac{\sqrt[5]{144\left(\sqrt{5}+2\right)^2S_4^4}}{6\left(\sqrt{5}+3\right)},\quad q_5=-\sqrt[5]{12\left(\sqrt{5}+2\right)S_4^2}.
\label{spoint6}
\end{equation}

\medskip
\lemma\label{samoperesechS4-4-2} {\it The equation {\rm(\ref{dobavka2})} for $S_4\neq0$ determines four curves.

{\rm 1)} The graph of a smooth negative function
$
q_5=\tau_1(\delta t_1), \delta t_1<0.
$
The curve $\tau_1$ lies in the domain that is bounded by the curves $\gamma_1^-$ and $\zeta_2$.

{\rm 2)} The graph of a smooth negative function
$
q_5=\tau_2(t_1), \delta t_1>0.
$
The curve $\tau_2$ has a simple tangency with the curve $\eta_2$ at {\rm(\ref{spoint6})} and lies in the domain that is bounded by the curves $\eta_2$ and $\gamma_2^-$. It transversally intersects the curve $\gamma_1^+$ at the point
\begin{equation}
\delta t_1=\frac13 \sqrt[5]{\frac{S_4^4}{12}},\quad q_5=-\frac32 \sqrt[5]{\frac{9S_4^2}{2}}.
\label{spoint8}
\end{equation}

{\rm 3)} The graph of a smooth function
$
q_5=\tau_3(\delta t_1), \delta t_1>0.
$
The curve $\tau_3$ transversally intersects the $t_1$ axis at $\delta t_1=\frac{1}{6}\sqrt[5]{\frac{S_4^4}{3}}$, has a simple tangency with the curve $\eta_1$ at {\rm(\ref{spoint2})} and lies in the domain that is bounded by the curve $\gamma_1^+$ and the $q_5$ axis.

{\rm 4)} The graph of a function
$
\delta t_1=\tau_4(q_5), q_5\in\mathbb{R}.
$
The curve $\tau_4$ is smooth and lies in the domain that is bounded by the curves $\gamma_1^-$ and $\tau_3$. It is tangent to the $t_1$ axis at the origin {\rm(}with order $2${\rm)}, and $q_5\tau_4(q_5)>0$ for all $q_5\neq0$.}

\medskip
\lemma\label{samoperesechS4-4-1} {\it
Let $S_4\neq0$. Then the set $\xi\subset\mathbb{R}^2=\{(t_1,q_5)\}$ from Lemma {\rm\ref{samoperesechS4-4}} consists of five curves.

{\rm 1)} The graph of a smooth function
$
q_5=\xi_1(\delta t_1), \delta t_1<0.
$
The curve $\xi_1$ transversally intersects the $t_1$ axis at $\delta t_1=-\frac{2}{3}\sqrt[5]{\frac{S_4^4}{3}}$ and lies in the domain that is bounded by the curves $\gamma_1^-$ and $\tau_4$.

{\rm 2)} The graph of a smooth function
$
q_5=\xi_2(\delta t_1), \delta t_1>0.
$
The curve $\xi_2$ lies in the domain that is bounded by the curve $\gamma_1^+$ and does not contain the origin.

{\rm 3)} The graph of a function
$
\delta t_1=\xi_3(q_5), q_5\in\mathbb{R}
$, 
where $\xi_3(q_5)>0$ for all $q_5\neq0$. The curve $\xi_3$ has a singularity of type $11/3$ at the origin and is smooth at all other points. It is tangent to the $q_5$ axis at the origin and lies in the domain that is bounded by the curve $\gamma_2^+$.

{\rm 4)} The curve $\xi_4$ with one semicubical cusp
\begin{equation}
\delta t_1=-\sqrt[5]{\frac{(13\sqrt{5}+29)S_4^4}{108}},\quad q_5=\sqrt[5]{6(\sqrt{5}+1)S_4^2}.
\label{spoint1}
\end{equation}
This point lies on the curve $\zeta_1$ and divides $\xi_4$ into the graphs of two smooth functions
$$
q_5=\xi_4^\pm(t_1),\quad \delta t_1<-\sqrt[5]{\frac{(13\sqrt{5}+29)S_4^4}{108}}.
$$
The curves $\xi_4^+$ and $\xi_4^-$ lies in the domain that is bounded by the curves $\zeta_1$ and $\eta_1$.

{\rm 5)} The curve $\xi_5$ with one semicubical cusp
\begin{equation}
\delta t_1=\sqrt[5]{\frac{(13\sqrt{5}-29)S_4^4}{108}},\quad q_5=-\sqrt[5]{6(\sqrt{5}-1)S_4^2}.
\label{spoint4}
\end{equation}
This point lies on the curve $\zeta_3^-$ and divides $\xi_5$ into the graphs of two smooth functions
$$
t_1=\xi_5^\pm(q_5),\quad q_5<-\sqrt[5]{6(\sqrt{5}-1)S_4^2}.
$$
The curves $\xi_5^+$ and $\xi_5^-$ lies in the domain that is bounded by the curves $\gamma_2^+,\gamma_2^-$. The curve $\xi_5^+$ transversally intersects $\gamma_1^+$ at the point}
\begin{equation}
\delta t_1=\frac{1}{3}\sqrt[5]{\frac{S_4^4}{6}},\quad q_5=-\frac{3}{2}\sqrt[5]{\frac{9S_4^2}{8}}.
\label{trper-xi5-g1}
\end{equation}

{\sc Proof.} To determine the type of the singular point $(0,0)$ of the curve (\ref{per-points-3}) it is necessary to make the following change of coordinates:
$
t_1\mapsto t_1+\delta \frac{q_5^2}{24}, q_5\mapsto q_5.
$
The equation (\ref{per-points-3}) in new coordinates is equivalent to
\begin{equation}
\delta t_1^3-\delta t_1q_5^9-q_5^{11}+\dots=0,
\label{puiso3}
\end{equation}
where dots denote a polynomial in $t_1,q_5$ consisting of monomials with quasi-degree greater than $38$ in the filtration $w(t_1)=11,w(q_5)=3$. The curve (\ref{puiso3}) has one real branch $\xi_3$ at the origin. Puiseux series in $q_5$ for a germ of the complex extension of the curve $\xi_3$ at zero starts at $q_5^{11/3}$. The rest is the result of direct calculations. $\Box$

\medskip
\lemma\label{samoperesechS4-5} {\it Let $S_4\neq0$. Then the set $\gamma_3\subset\mathbb{R}^2=\{(t_1,q_5)\}$ from Lemma {\rm\ref{samoperesechS4-4}} consists of two smooth curves:
$$
\gamma_3^-: q_5\in\left(\sqrt[5]{12S_4^2},\sqrt[5]{6(\sqrt{5}+1)S_4^2}\right),
\quad
\gamma_3^+: q_5\in\left(-\infty,-\sqrt[5]{6(\sqrt{5}-1)S_4^2}\right)\cup\left(0,\sqrt[5]{12S_4^2}\right).
$$
The curve $\gamma_3^-$ is bounded by the point {\rm(\ref{spoint1})}. The curve $\gamma_3^+$ consists of two connected components. The first component is bounded by the origin and passes through the point {\rm(\ref{spoint2})} touching the curve $\eta_1$. The second one is bounded by the point {\rm(\ref{spoint4})} and passes through the point {\rm(\ref{spoint6})} touching the curve $\eta_2$. It is also tangent to the curve $\gamma_1^+$ at the point
\begin{equation}
\delta t_1=\sqrt[5]{\frac{S_4^4}{1728}}, \quad q_5=-\sqrt[5]{12S_4^2}.
\label{spoint5}
\end{equation}
All three tangencies mentioned above are simple. The curve $(\ref{psi-2})$ corresponding to the point $(\ref{spoint5})$ passes at $u=0$ through the intersection point
\begin{equation}
t_2=-\frac{S_4^4+864\delta t_1^5}{432\delta t_1^3S_4},\quad S_3=\frac{S_4^2}{6t_1}
\label{spec-kas}
\end{equation}
of lines {\rm(\ref{crit-tochki})} and {\rm(\ref{psi-4})}.}

{\sc Proof.} If $(t_1,q_5)$ satisfies the equation (\ref{dobavka1}), then the polynomials $P(u)$ from Lemma \ref{samoperesechS4-3} and $Q(u)$ from Lemma \ref{samoperesechS4-4} are divisible by the polynomial
$$
3(q_5^5-12S_4^2)u^2+12q_5^3S_4u+q_5(q_5^5+12S_4^2).
$$
The discriminant of this polynomial with respect to $u$ must be positive. $\Box$

\medskip
\lemma\label{peresech-izolir} {\it Let $S_4\neq0$ and $(t_1,q_5)\in\gamma_1^+$. Then the line $(\ref{psi-4})$ intersects the curve $(\ref{psi-2})$ at two different points:
\begin{equation}
t_2=\frac{S_4^4+48S_4^2\sqrt{3\delta t_1^5}-864\delta t_1^5}{432\delta t_1^3S_4},
\quad
S_3=-\frac{S_4^2+48\sqrt{3\delta t_1^5}}{6t_1},
\label{peresech1}
\end{equation}
\begin{equation}
t_2=\frac{S_4^4-48S_4^2\sqrt{3\delta t_1^5}-864\delta t_1^5}{432\delta t_1^3S_4},
\quad
S_3=-\frac{S_4^2-48\sqrt{3\delta t_1^5}}{6t_1}.
\label{peresech2}
\end{equation}
The point $(\ref{spec-kas})$ lies between $(\ref{peresech1})$ and $(\ref{spec2})$. The point $(\ref{peresech2})$ lies between $(\ref{peresech1})$ and $(\ref{spec-kas})$ if $(t_1,q_5)$ belongs to the curve
\begin{equation}
\gamma_1^+: 0<\delta t_1<\sqrt[5]{\frac{S_4^4}{1728}}\quad (\mbox{bounded by $(\ref{spoint5})$});
\label{gam1+b}
\end{equation}
the point $(\ref{peresech2})$ lies between $(\ref{spec-kas})$ and $(\ref{spec2})$ if $(t_1,q_5)$ belongs to any from the curves
\begin{equation}
\gamma_1^+: \sqrt[5]{\frac{S_4^4}{1728}}
<\delta t_1<\sqrt[5]{\frac{S_4^4}{432}}\quad (\mbox{between $(\ref{spoint5})$ and $(\ref{vozvrat})$}),
\label{gamma1+part}
\end{equation}
\begin{equation}
\gamma_1^+: \delta t_1>\sqrt[5]{\frac{S_4^4}{432}}\quad (\mbox{bounded by $(\ref{vozvrat})$}).
\label{gam1+t}
\end{equation}
Finally, $(\ref{peresech2})\equiv(\ref{spec-kas})$ if $(t_1,q_5)\equiv(\ref{spoint5})$, and $(\ref{peresech2})\equiv(\ref{spec2})$ if $(t_1,q_5)\equiv(\ref{vozvrat})$.}

\medskip
This statement is the result of simple calculations.

\section{Perestroikas of skeletons $\Sigma_{S_4,t_1,q_5},S_4\neq0$}

The diagram $B_{S_4},S_4\neq0$ divides $\mathbb{R}^2=\{(t_1,q_5)\}$ into $93$ domains (connected components of the complement to $B_{S_4}$). These domains are enumerated as shown in Fig. \ref{razbienie}. Every domain of the partition is contractible. The skeletons $\Sigma_{S_4,t_1,q_5}$ corresponding to points $(t_1,q_5)$ from the same domain are diffeomorphic. The topology of $\Sigma_{S_4,t_1,q_5}$ is changed at the moment when $(t_1,q_5)$ passes into a neighboring domain, moving along a smooth curve that transversally intersects the boundary of this domain at a non-singular point. This change is said to be perestroika.

There are three types of perestroikas. The following six perestroikas have the {\it first type}.

\medskip
${\cal P}_1$) The perestroika of the section of swallowtail by a plane when the plane passes through the swallowtail vertex transversally to the tangent line at the vertex to the closure of the cuspidal edge. When choosing a suitable perestroika direction, two semicubical cusps and a transversal intersection point of two smooth branches of the curve $(\ref{psi-2})$ appear; one two-dimensional stratum of the stratification $\mathfrak{S}_{S_4,t_1,q_5}$ appears as well. The perestroika ${\cal P}_1$ happens when we cross the curves $\zeta_1,\zeta_2,\zeta_3^\pm$.

\medskip
${\cal P}_2$) The passage of a smooth branch of the curve $(\ref{psi-2})$ through a semicubical cusp transversally to the tangent line to the curve at this point. When choosing a suitable perestroika direction, two transversal intersection points of two smooth branches of the curve $(\ref{psi-2})$ appear; one two-dimensional stratum of $\mathfrak{S}_{S_4,t_1,q_5}$ appears as well; one two-dimensional stratum is divided into two ones. The perestroika ${\cal P}_2$ happens when we cross the curves $\eta_1,\eta_2,\eta_3^\pm$.

\medskip
${\cal P}_3$) The passage of a smooth branch of the curve $(\ref{psi-2})$ through a point of the simple tangency of its other smooth branch with the line $(\ref{crit-tochki})$ transversally to this line. Here, one two-dimensional stratum of $\mathfrak{S}_{S_4,t_1,q_5}$ disappears, and another one appears. The perestroika ${\cal P}_3$ happens when we cross the curves $\tau_1, \tau_2, \tau_3, \tau_4$.

\medskip
${\cal P}_4$) The passage of the line $(\ref{crit-tochki})$ through a semicubical cusp of the curve $(\ref{psi-2})$ transversally to the tangent line to the curve a this point. When choosing a suitable perestroika direction, two transversal intersection points of the line $(\ref{crit-tochki})$ with smooth branches of the curve appear; one two-dimensional stratum of $\mathfrak{S}_{S_4,t_1,q_5}$ appears, and one of two-dimensional strata is divided into two ones. The perestroika ${\cal P}_4$ happens when we cross the curves $\xi_1,\xi_2,\xi_3,\xi_4^\pm,\xi_5^\pm$.

\medskip
${\cal P}_5$) The passage of the line $(\ref{crit-tochki})$ through a transversal intersection point of two smooth branches of the curve $(\ref{psi-2})$ transversally to both intersecting branches. Here, two transversal intersection points of the line with these branches are transposed; one two-dimensional stratum of $\mathfrak{S}_{S_4,t_1,q_5}$ disappears, and another one appears. The perestroika ${\cal P}_5$ happens when we cross the curves $\gamma_3^\pm$.

\medskip
${\cal P}_6$) The passage of a simple tangency point of the line $(\ref{crit-tochki})$ with the curve $(\ref{psi-2})$ through a semicubical cusp that is not a self-intersection point of this curve. Here, one two-dimensional stratum of $\mathfrak{S}_{S_4,t_1,q_5}$ disappears, and another one appears. The perestroika ${\cal P}_6$ happens when we cross the curves $\gamma_2^\pm$.

\medskip
The {\it second type} perestroikas happens when we cross hyperbola (\ref{giperbola2}). There are three such perestroikas.

\begin{figure}[h]
\begin{center}
\begin{tabular}{ccc}
\includegraphics[width=4cm]{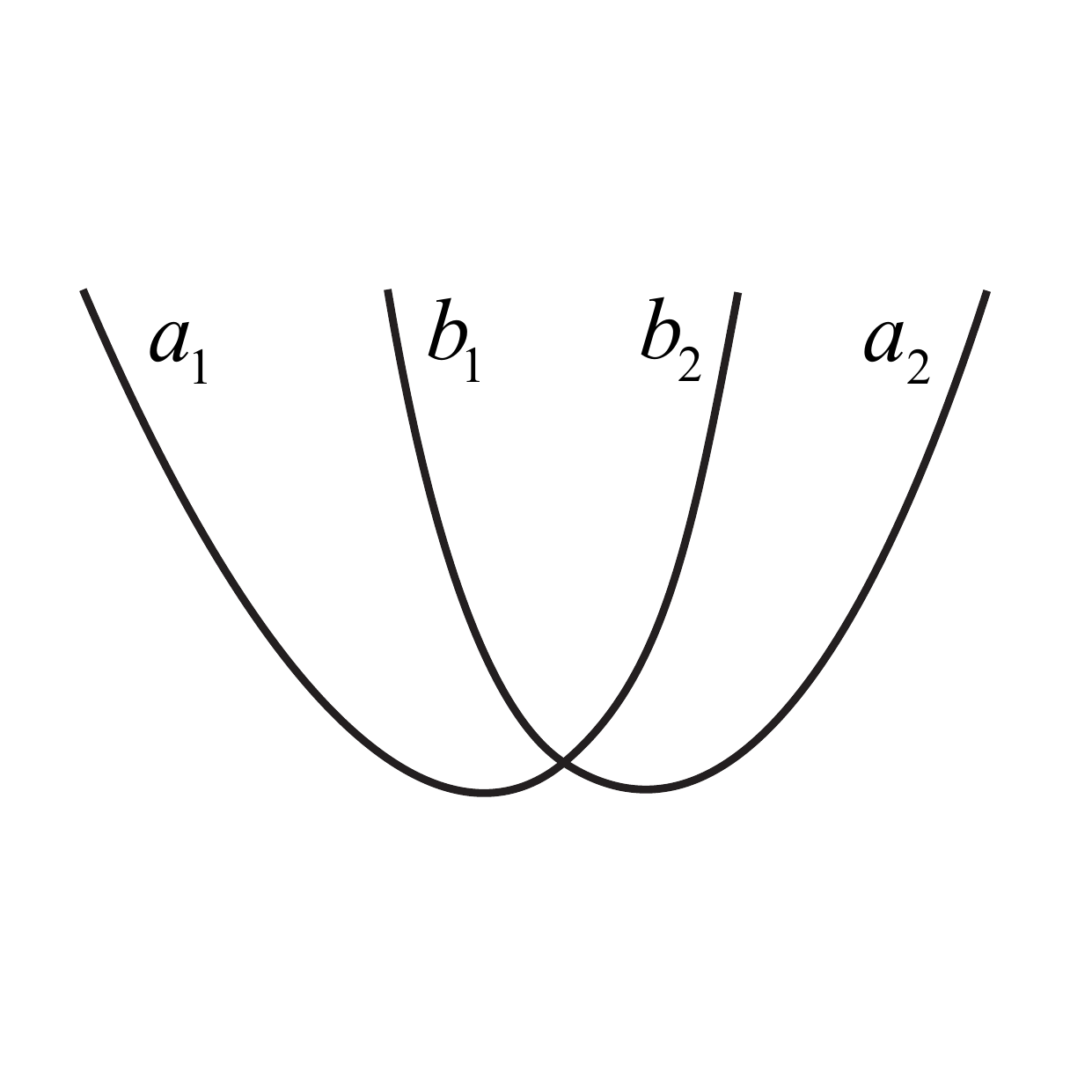}&
\includegraphics[width=4cm]{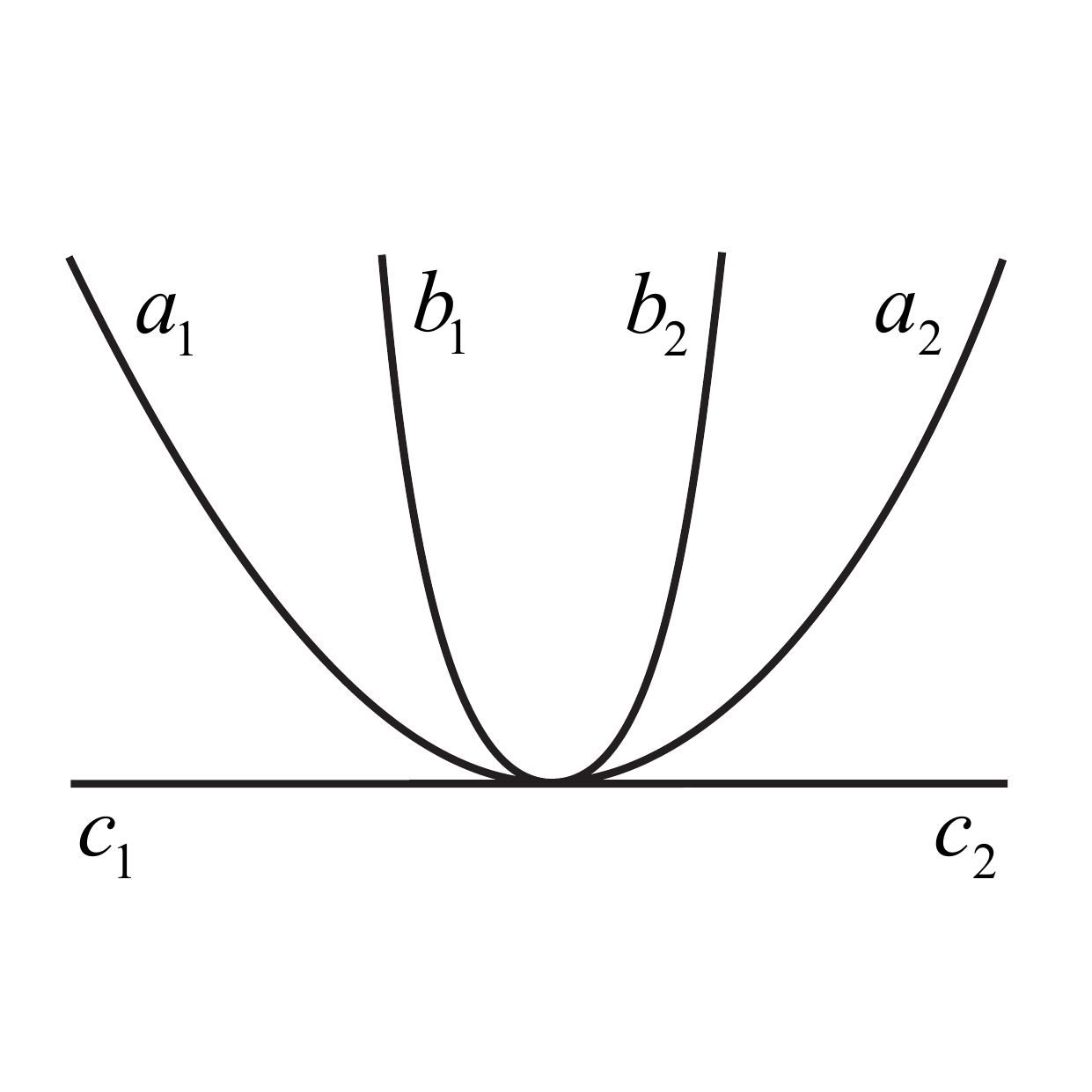}&
\includegraphics[width=4cm]{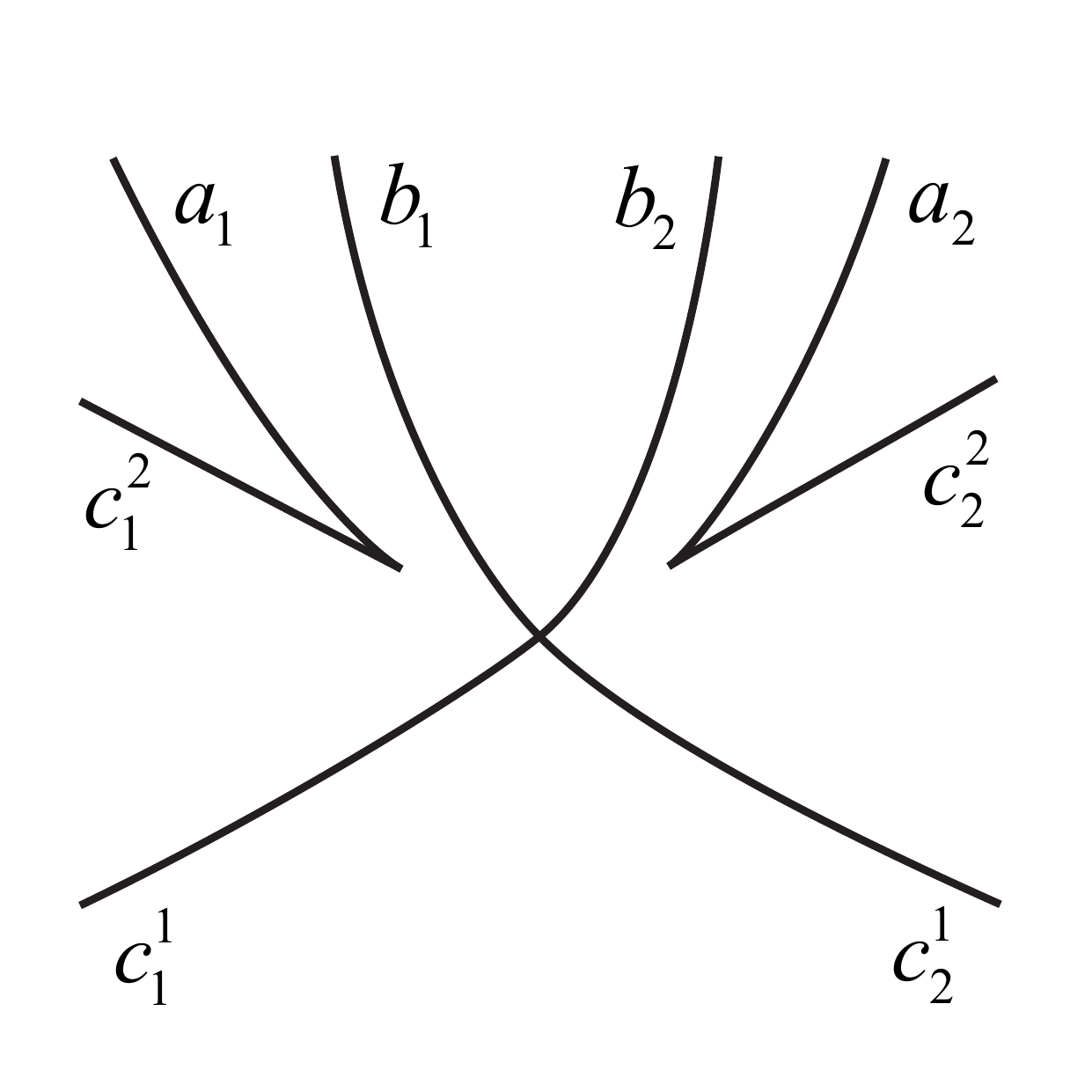}\\
\end{tabular}
\caption{The perestroika ${\cal P}_7$.}
\label{perestr3-4}
\end{center}
\end{figure}

\medskip
${\cal P}_7$) The perestroika in a neighbourhood of the point (\ref{spec2}) when we cross the curve $\gamma_1^-$. It is shown as a triptych in Fig. \ref{perestr3-4} up to diffeomorphism.

The middle picture corresponds to the intersection point $p\in\gamma_1^-$. The curves (\ref{psi-2}), (\ref{psi-3}) and (\ref{psi-4}) are represented by a smooth curve $a_1a_2$, parabola $b_1b_2$ and a straight line $c_1c_2$, respectively. The point of their tangency is (\ref{spec2}). It divides every curve into two parts: $a_1$ and $a_2$, $b_1$ and $b_2$, $c_1$ and $c_2$. The picture on the left corresponds to a point $p'$, which is close to $p$ and lies in the domain that is bounded by the curve $\gamma_1^-$ and does not contain the origin. Two smooth transversally intersecting curves $a_1b_2$ and $b_1a_2$ represent the curve (\ref{psi-2}). The intersection point divides every of them into two curves continuously changing as $p'\rightarrow p$. These curves and their limits are denoted by the same symbols.

The picture on the right corresponds to a point $p''$, which is close to $p$ and lies in the domain that is bounded by the curve $\gamma_1^-$ and contains the origin. Two smooth transversally intersecting curves $b_1c_2^1,c_1^1b_2$ and two singular curves $a_1c_1^2,c_2^2a_2$ with semicubical cusps represent the curve (\ref{psi-2}). The intersection point and the cusps divides every of them into two curves continuously changing as $p''\rightarrow p$. The limit of curves $c_i^1$ and $c_i^2$ is the ray $c_i, i=1,2$.

When choosing a suitable perestroika direction, two semicubical cusps of the curve (\ref{psi-2}) appear; one two-dimensional stratum of $\mathfrak{S}_{S_4,t_1,q_5}$ is divided into three ones. Global changes are also taking place. Namely one two-dimensional stratum lying on the other side of the line $(\ref{crit-tochki})$ with respect to the point (\ref{spec2}) is also divided into three ones.

\begin{figure}[h]
\begin{center}
\begin{tabular}{ccc}
\includegraphics[width=4cm]{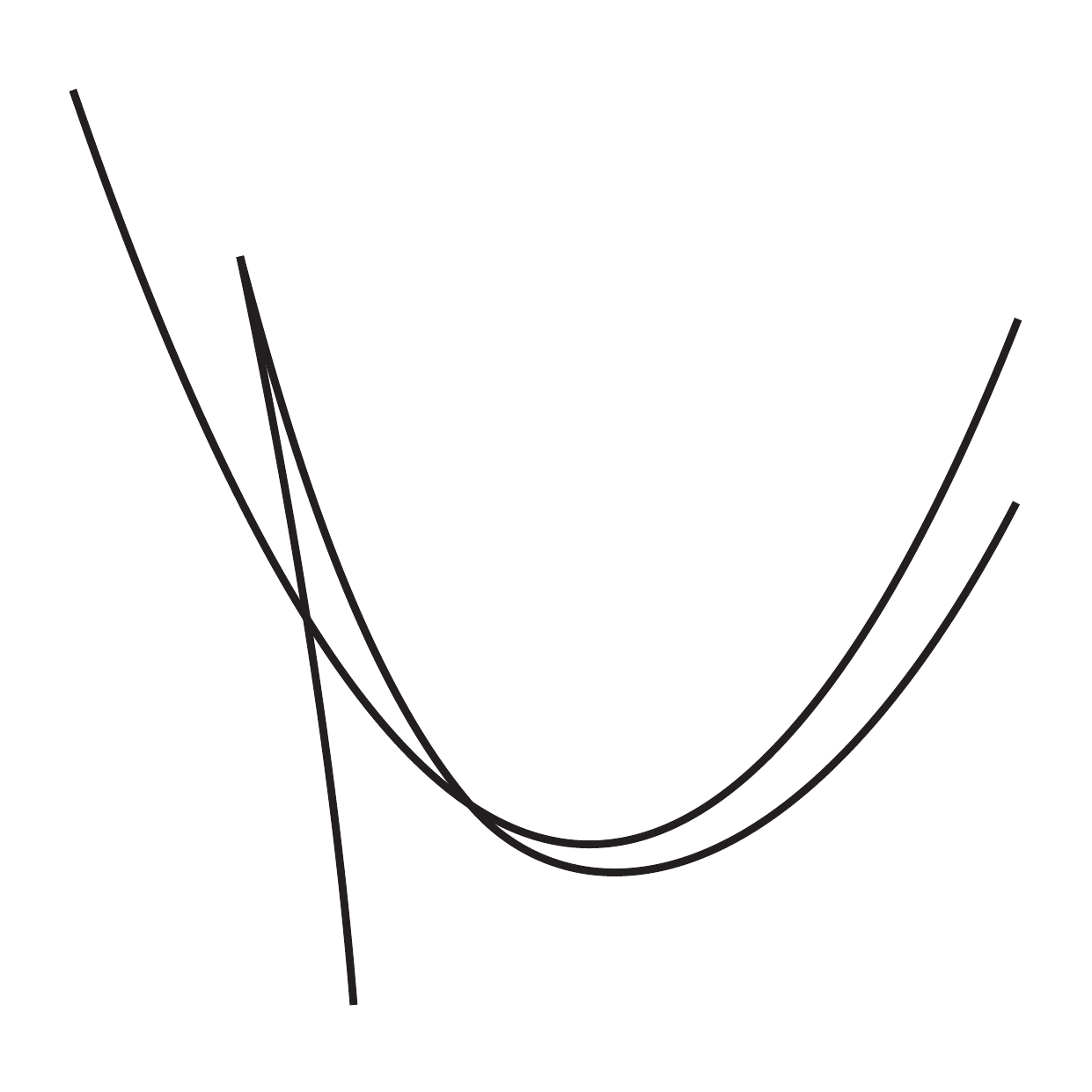}&
\includegraphics[width=4cm]{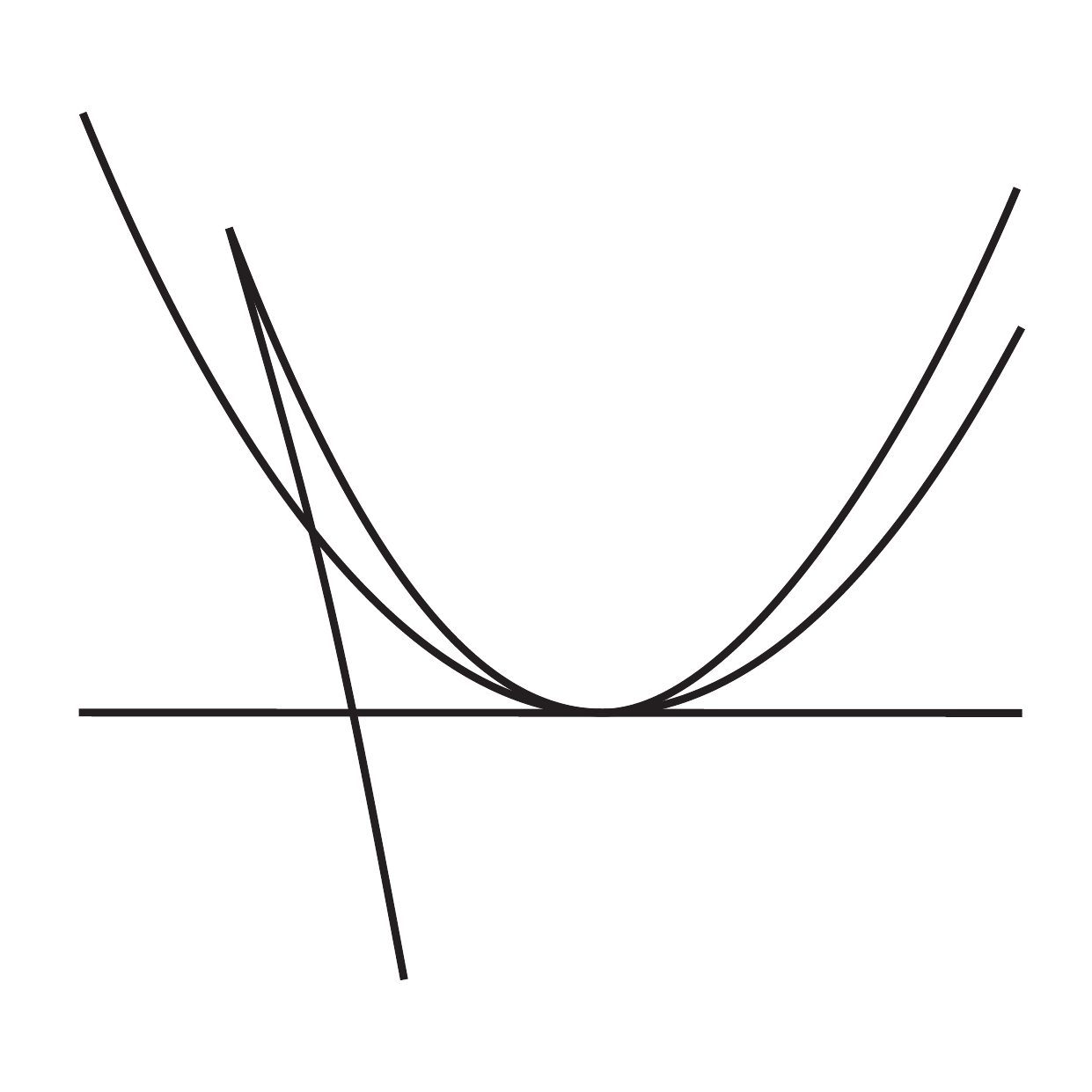}&
\includegraphics[width=4cm]{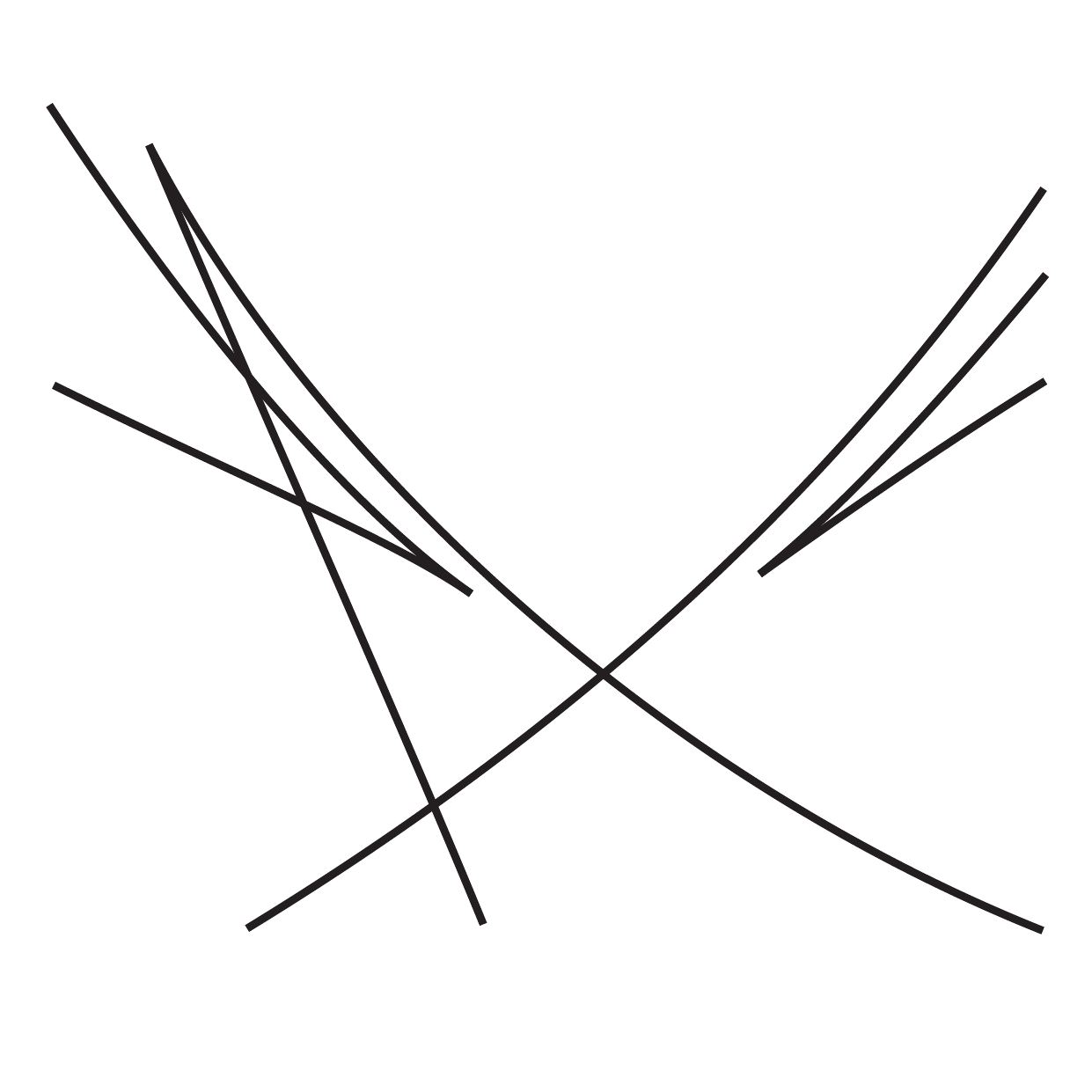}\\
\includegraphics[width=4cm]{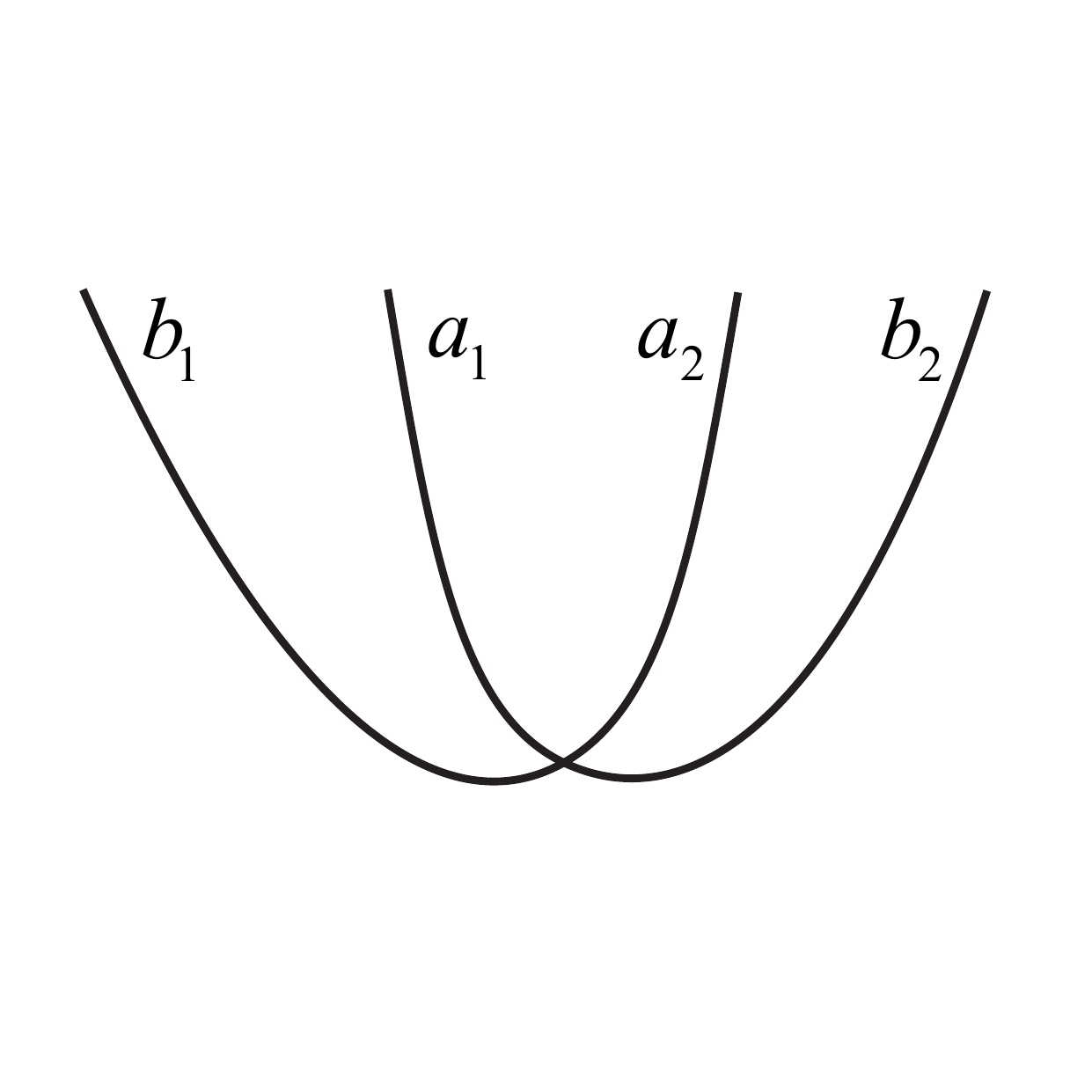}&
\includegraphics[width=4cm]{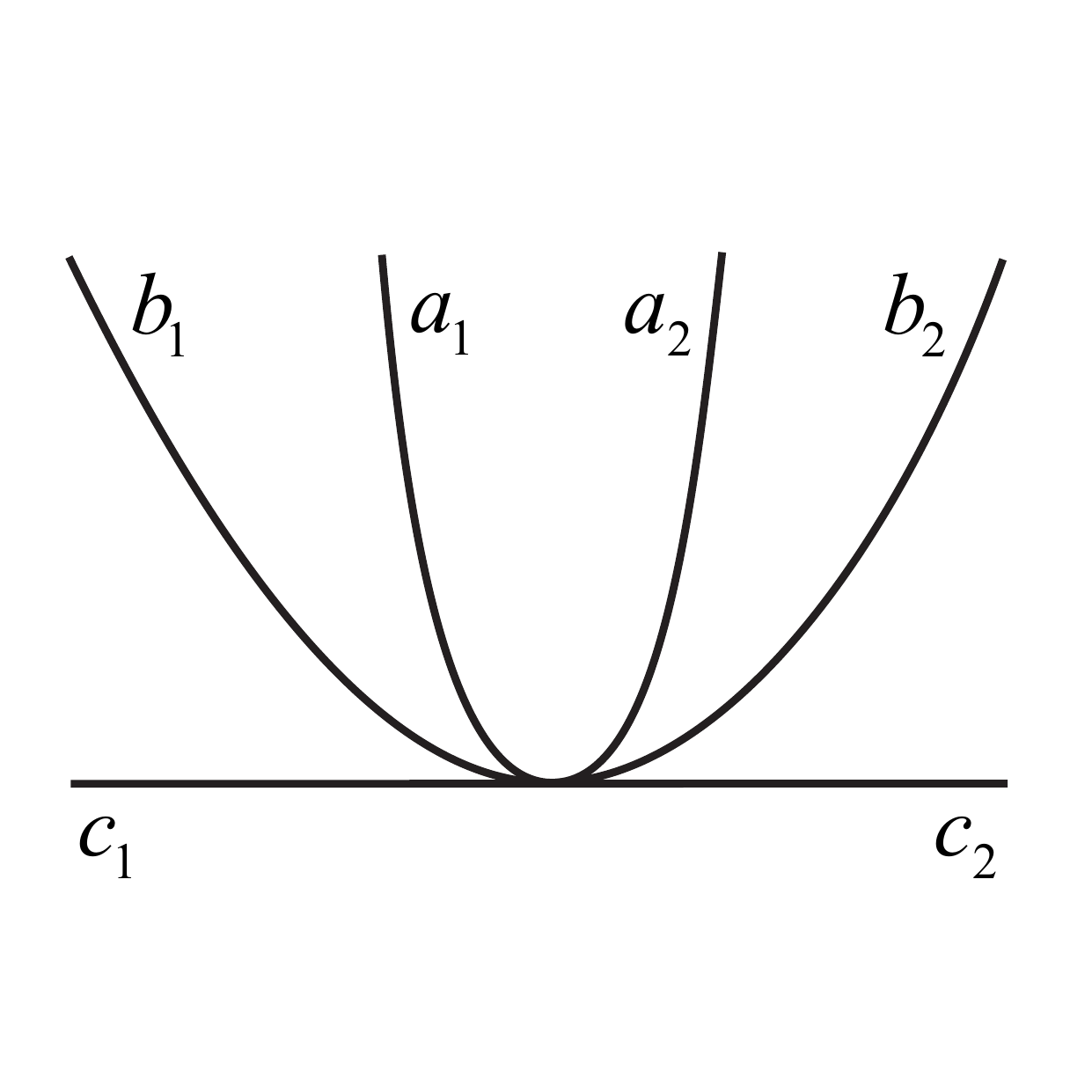}&
\includegraphics[width=4cm]{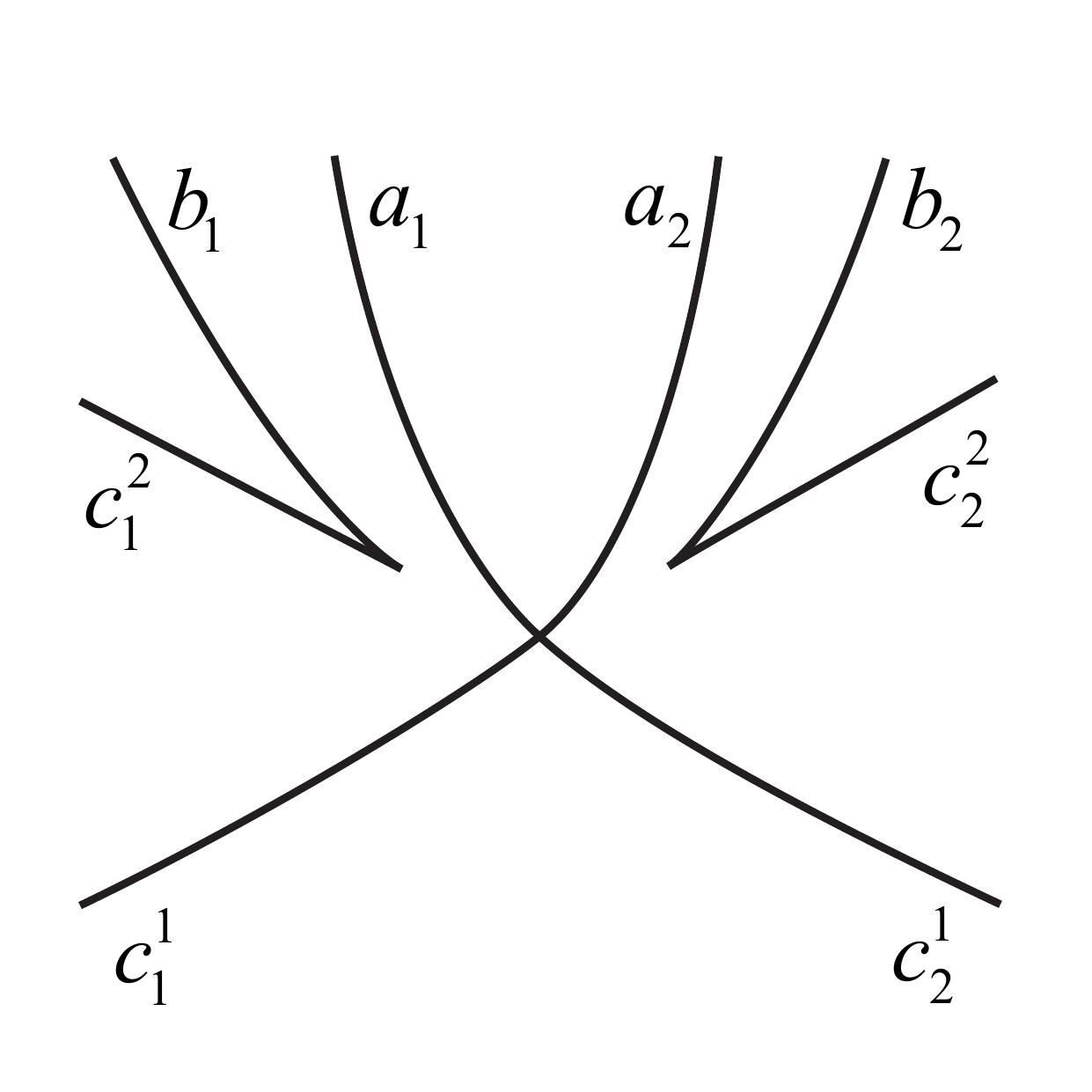}\\
\end{tabular}
\caption{The perestroika ${\cal P}_8$.}
\label{perestr50-49}
\end{center}
\end{figure}

\medskip
${\cal P}_8$) The perestroika in a neighbourhood of the point (\ref{spec2}) when we cross the curve (\ref{gam1+t}). It is shown by two triptychs in Fig. \ref{perestr50-49} up to diffeomorphism. The first triptych shows the perestroika in a neighbourhood of (\ref{spec2}) that contains a nonvanishing cusp of the curve (\ref{psi-2}). The second one demonstrate the perestroika in a sufficiently small neighbourhood of (\ref{spec2}). The description of this triptych is similar to the case ${\cal P}_7$. When choosing a suitable perestroika direction, two semicubical cusps of the curve (\ref{psi-2}) appear; one two-dimensional stratum of $\mathfrak{S}_{S_4,t_1,q_5}$ is divided into three ones. However, more complicated global changes are taking place.

The line (\ref{psi-4}) transversally intersects two smooth branches of the curve (\ref{psi-2}) at points lying on opposite sides of the line $(\ref{crit-tochki})$. One of these branches is shown in the first triptych as a branch of a semicubical parabola with nonvanishing cusp. It lies on the same side of the line $(\ref{crit-tochki})$ as the point (\ref{spec2}) (the second branch is not shown in Fig. \ref{perestr50-49}). Hence, there are three more two-dimensional strata, each of which is divided into three ones.

\begin{figure}
\begin{center}
\begin{tabular}{ccc}
\includegraphics[width=4cm]{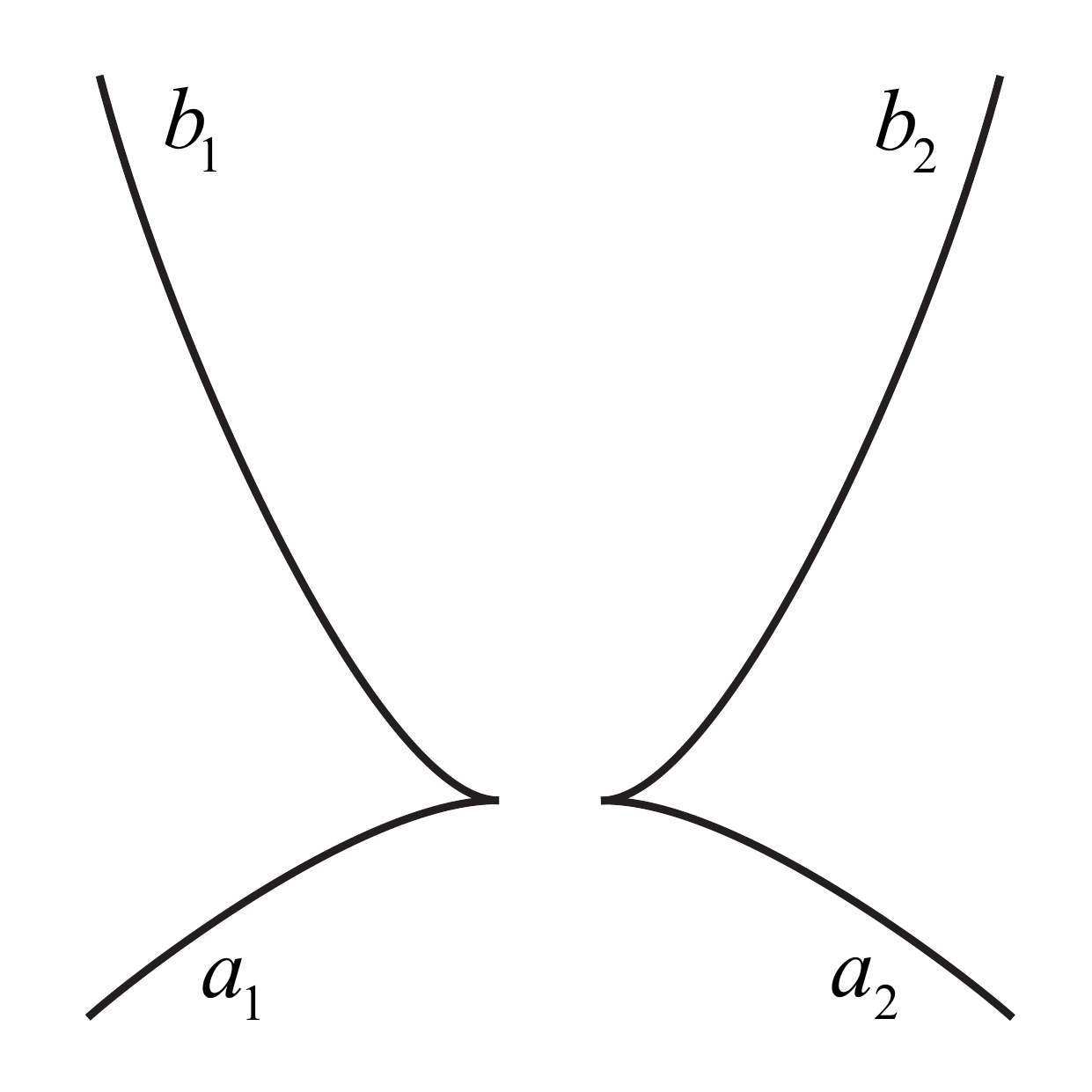}&
\includegraphics[width=4cm]{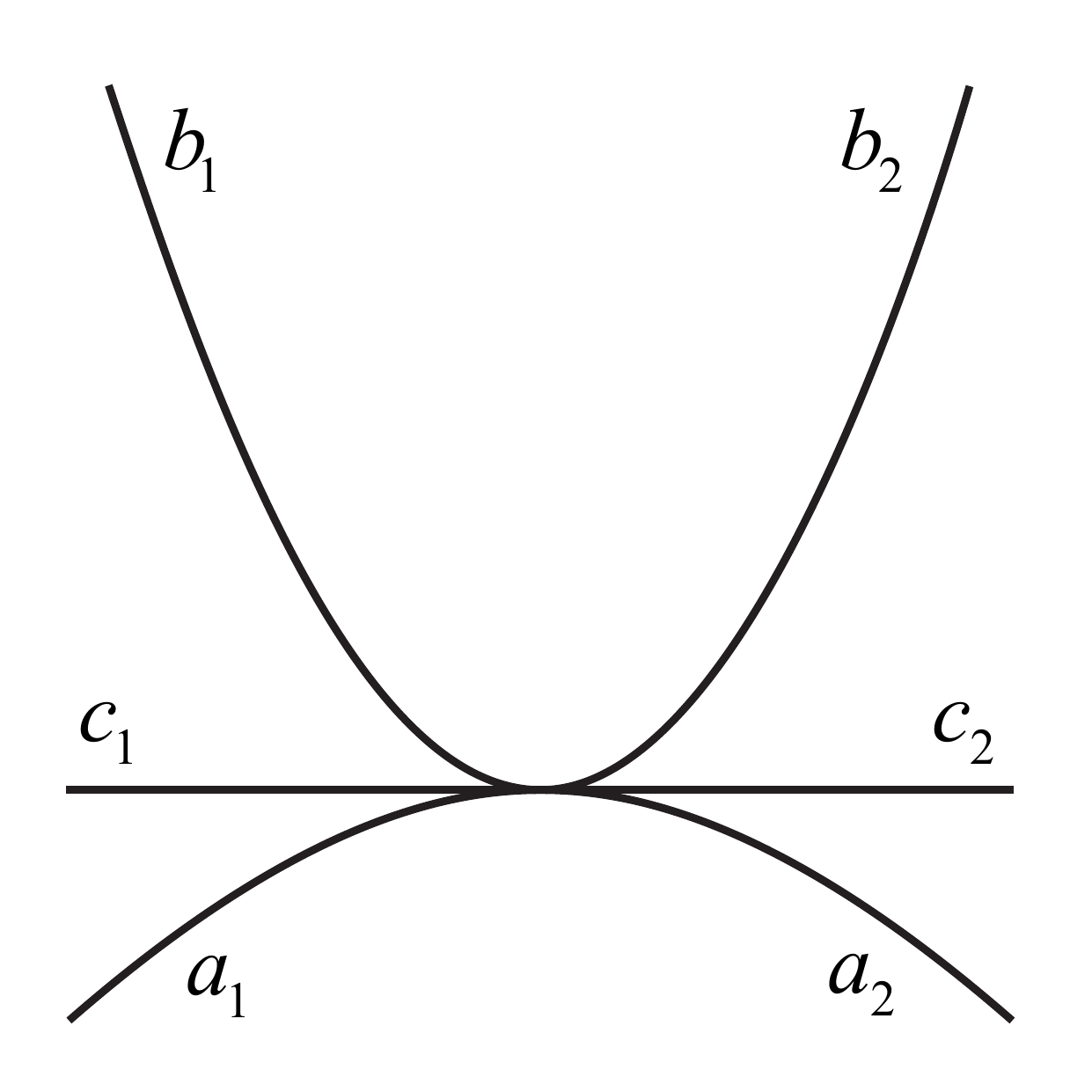}&
\includegraphics[width=4cm]{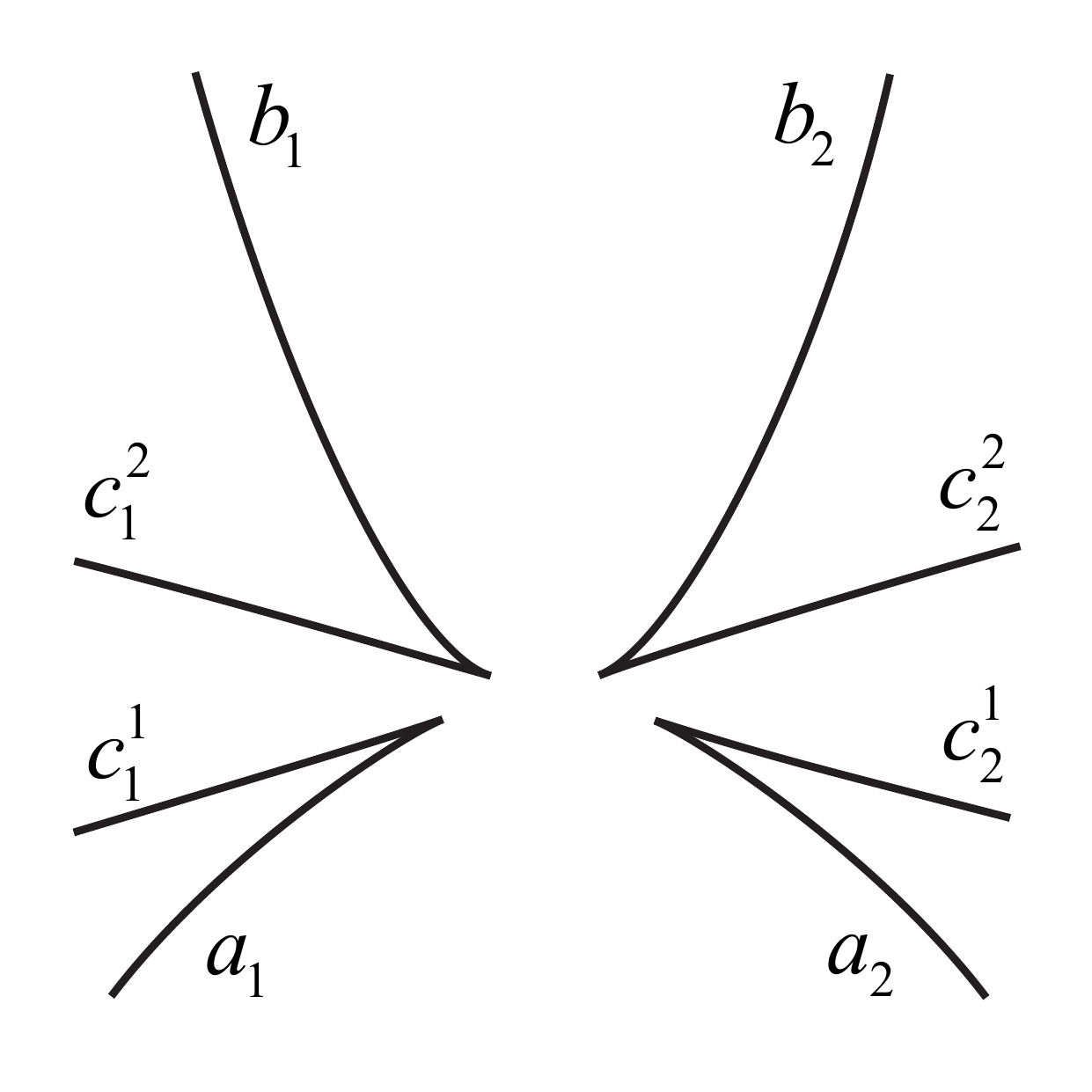}\\
\includegraphics[width=4cm]{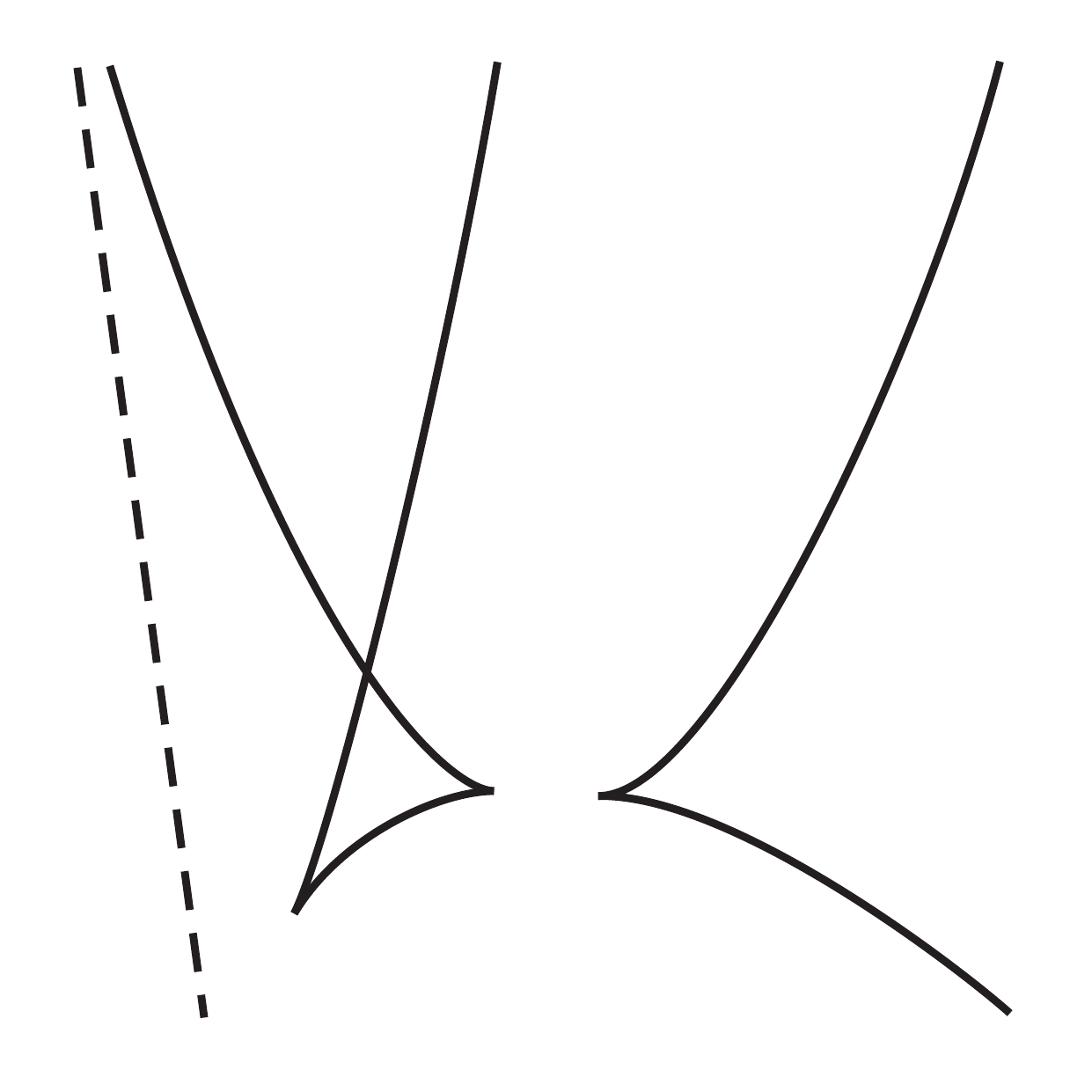}&
\includegraphics[width=4cm]{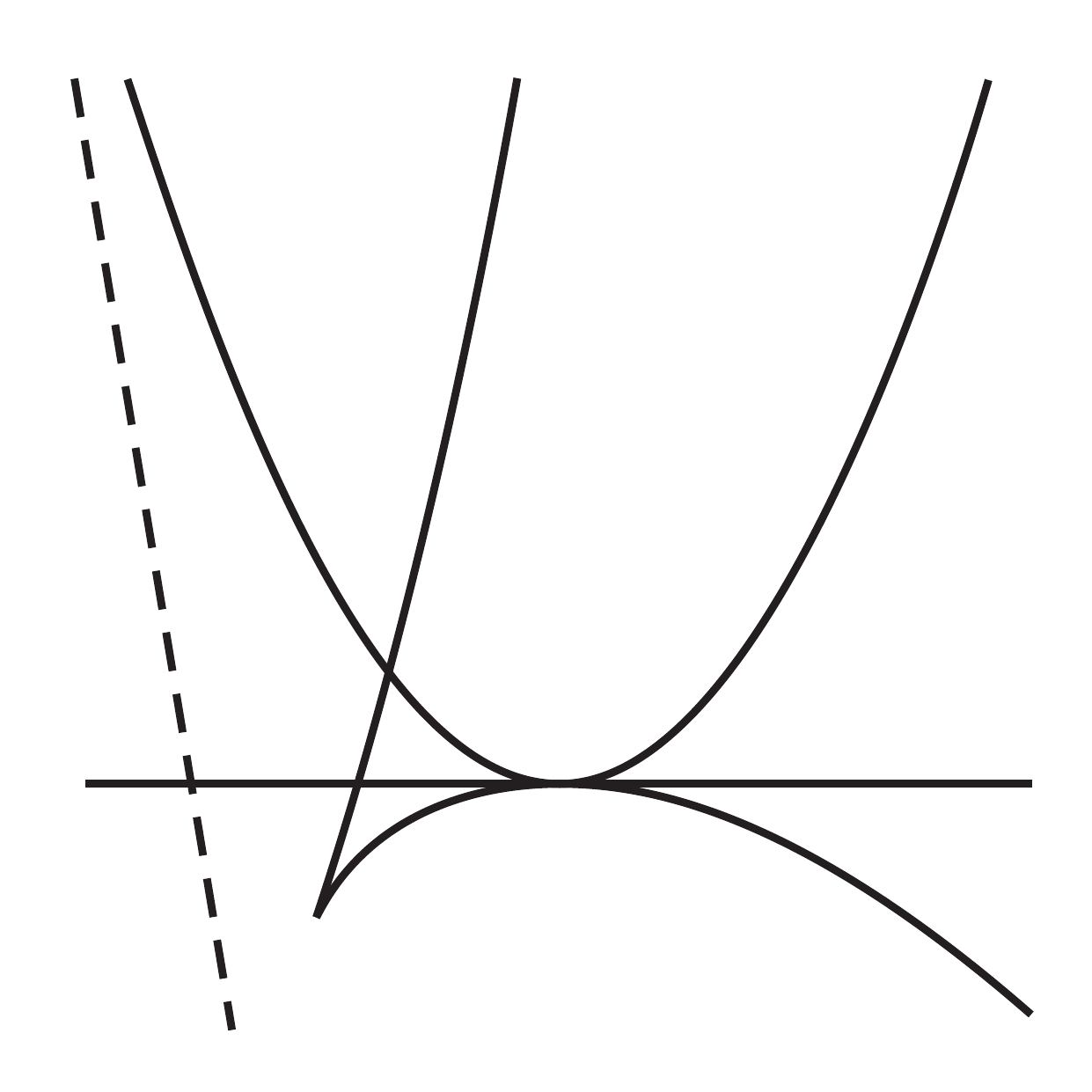}&
\includegraphics[width=4cm]{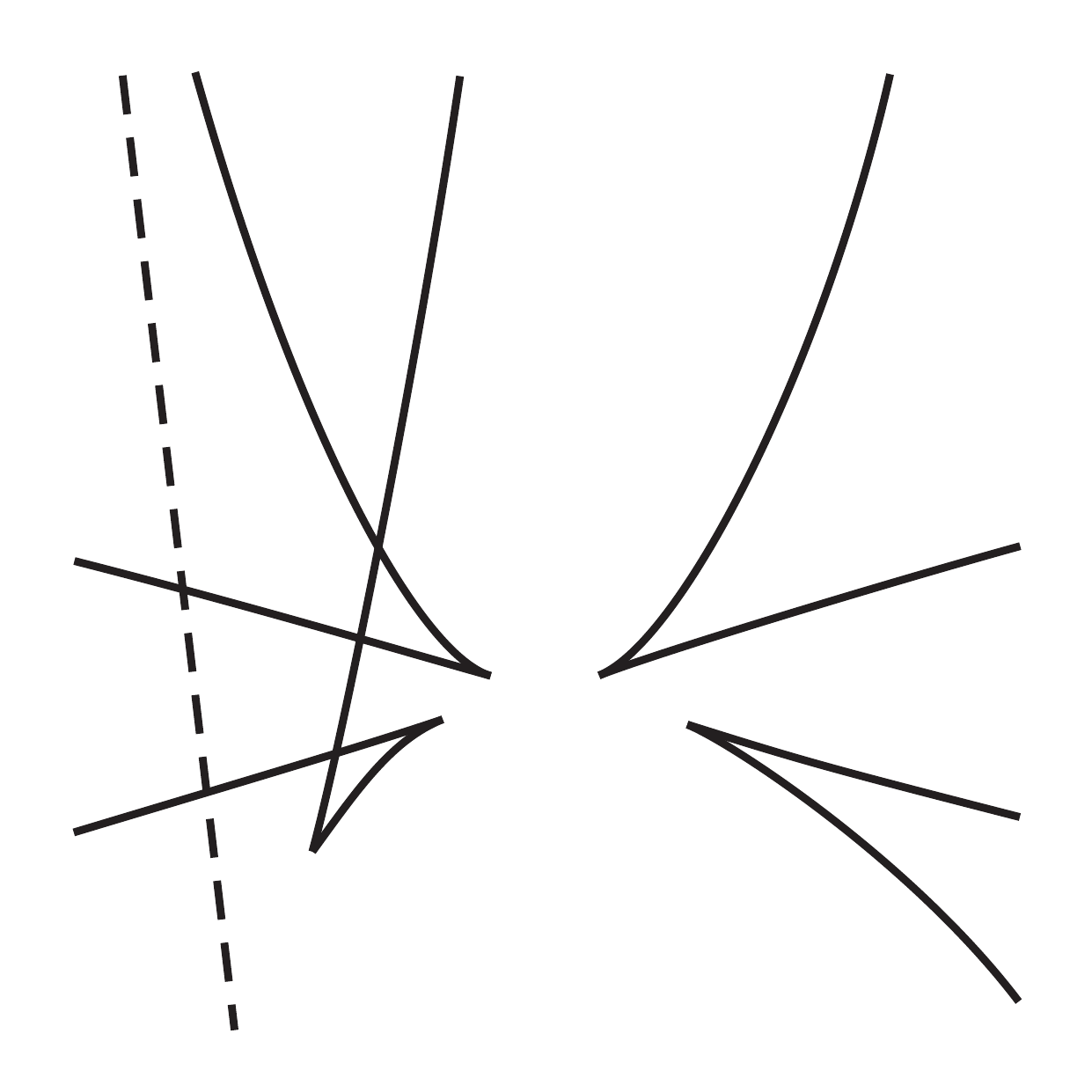}\\
\includegraphics[width=4cm]{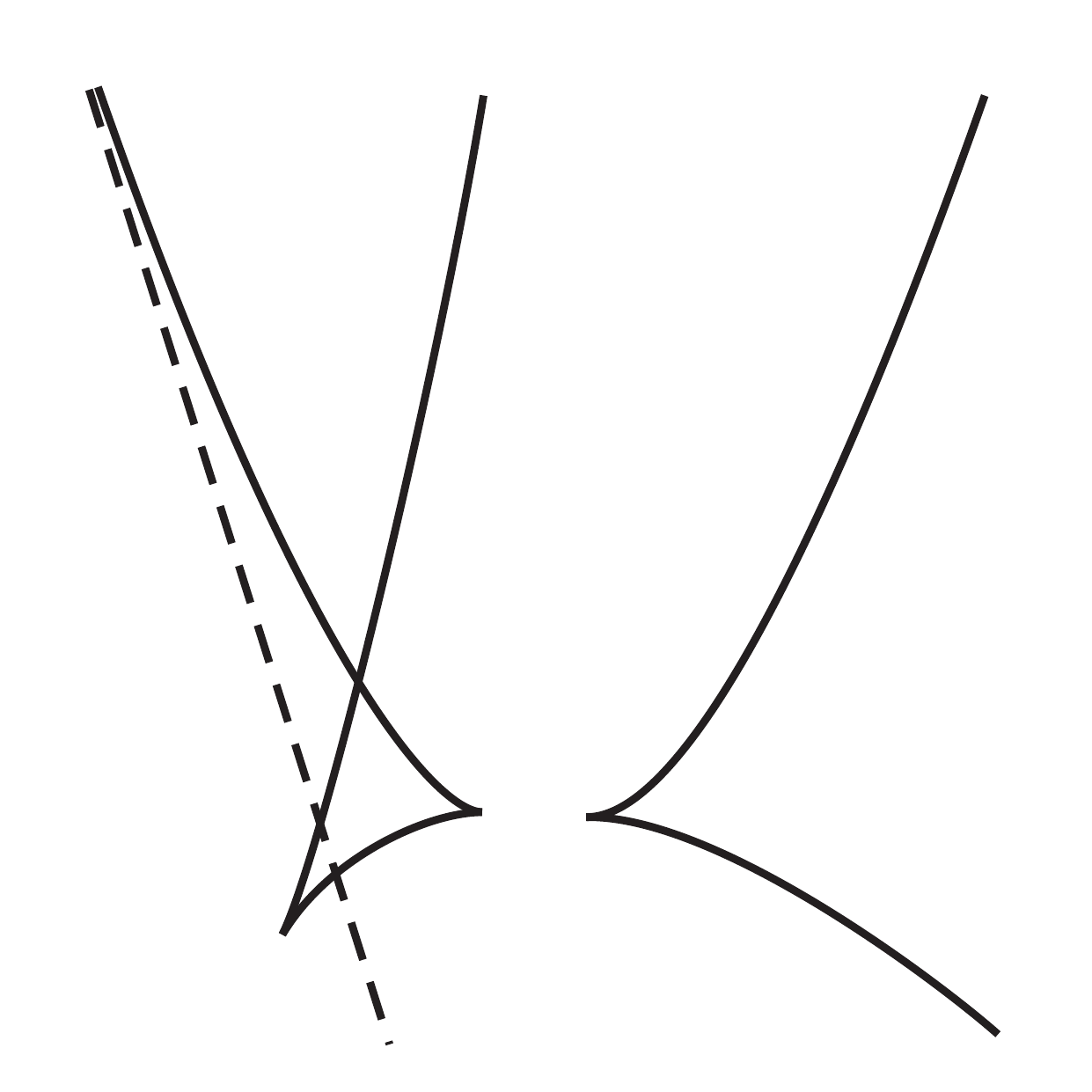}&
\includegraphics[width=4cm]{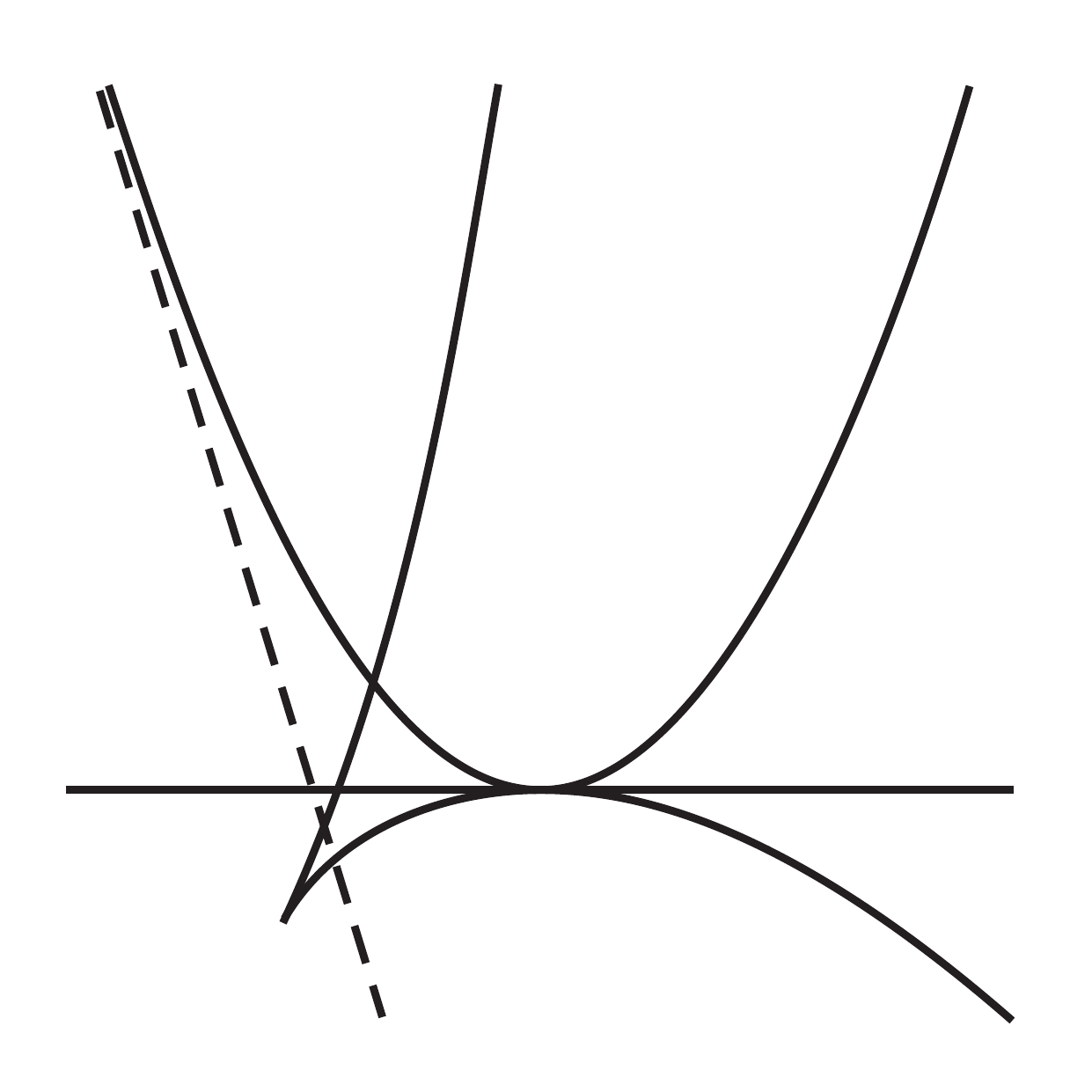}&
\includegraphics[width=4cm]{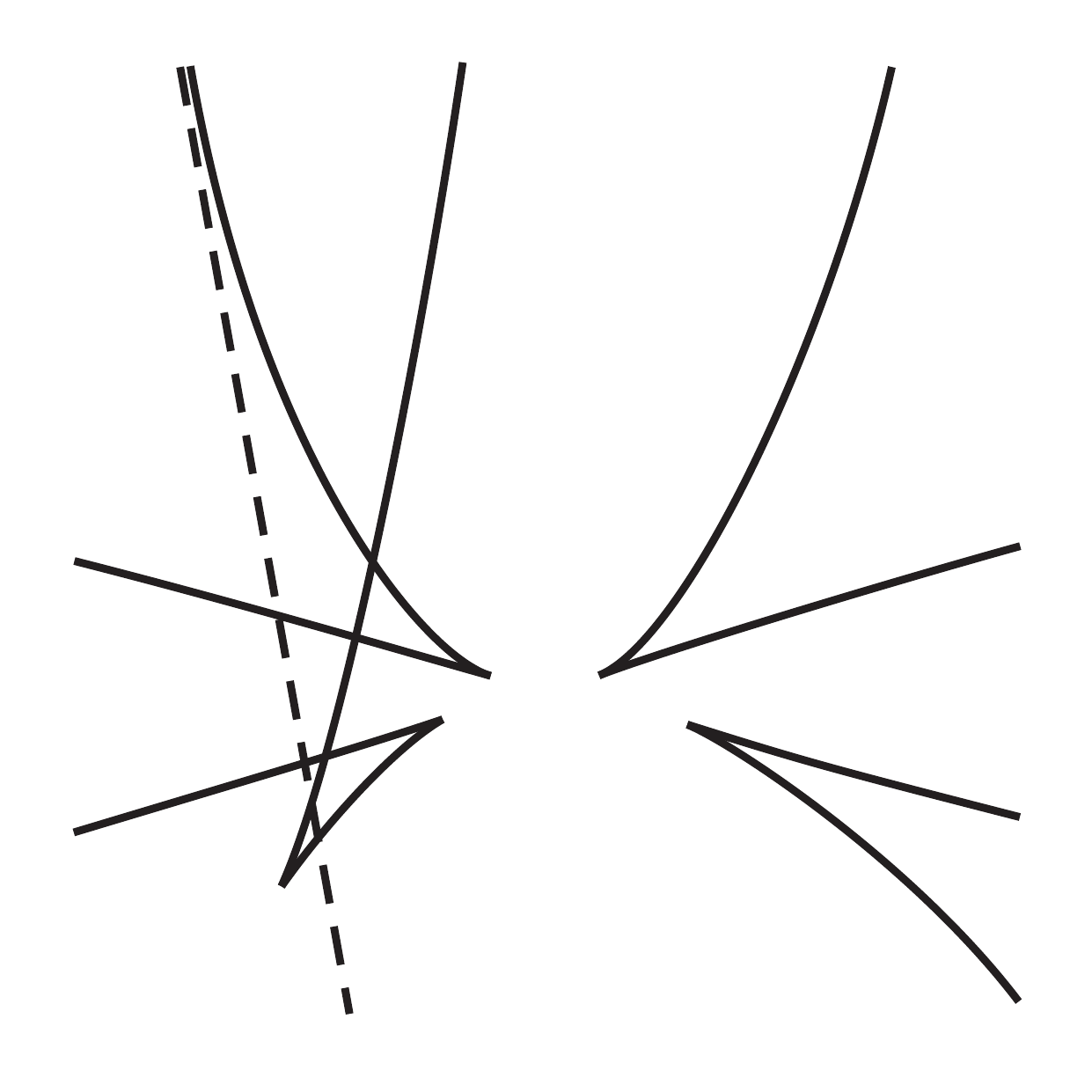}\\
\includegraphics[width=4cm]{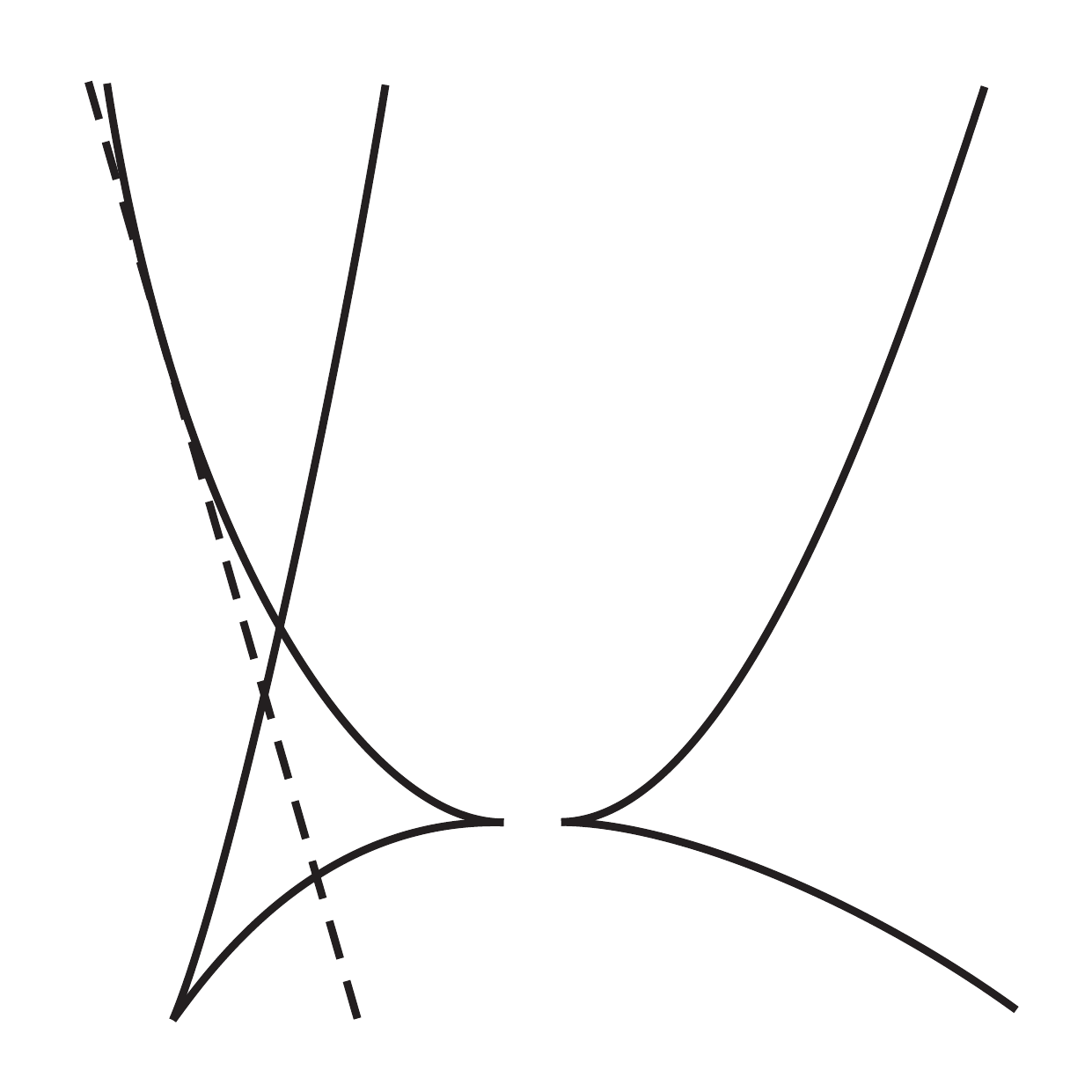}&
\includegraphics[width=4cm]{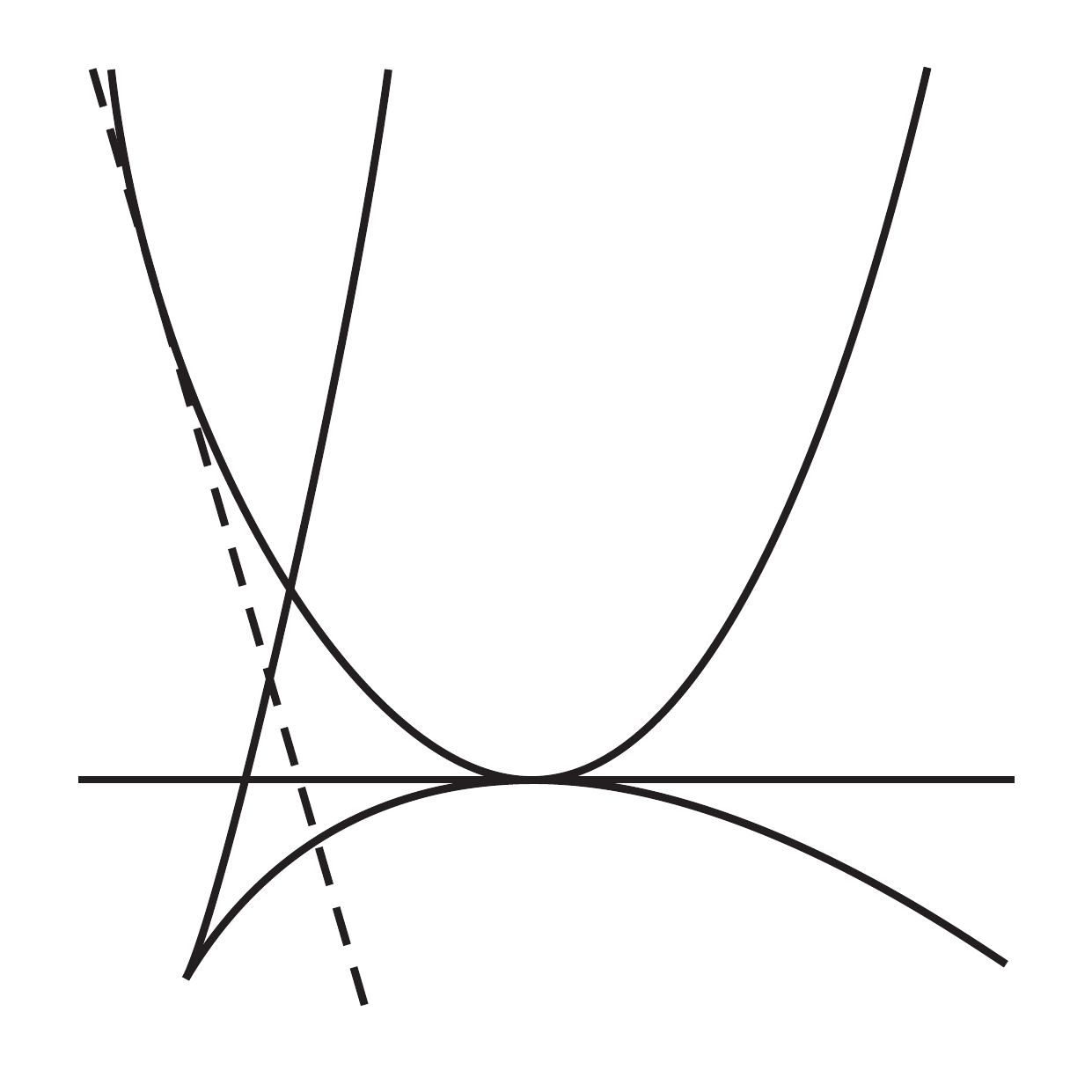}&
\includegraphics[width=4cm]{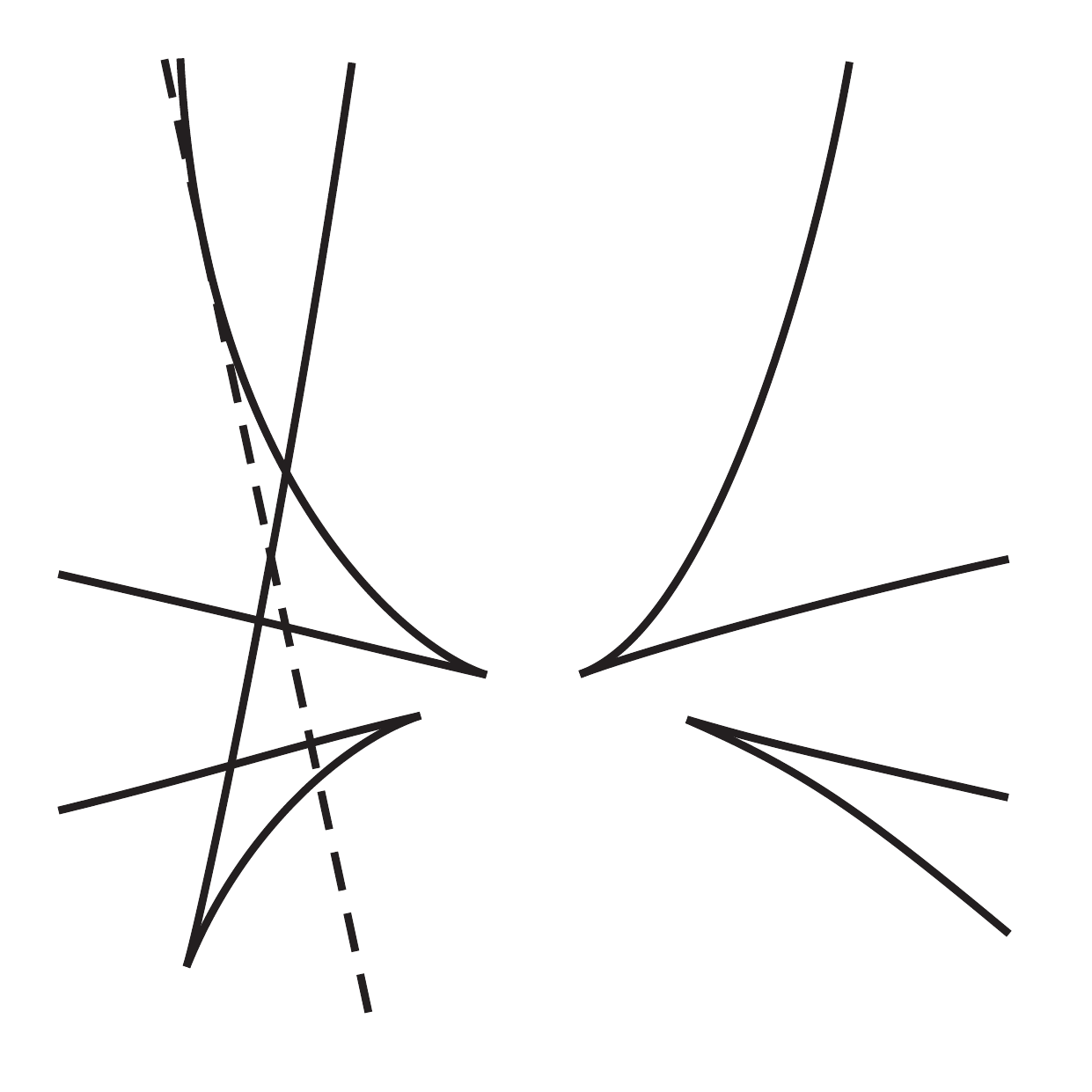}\\
\includegraphics[width=4cm]{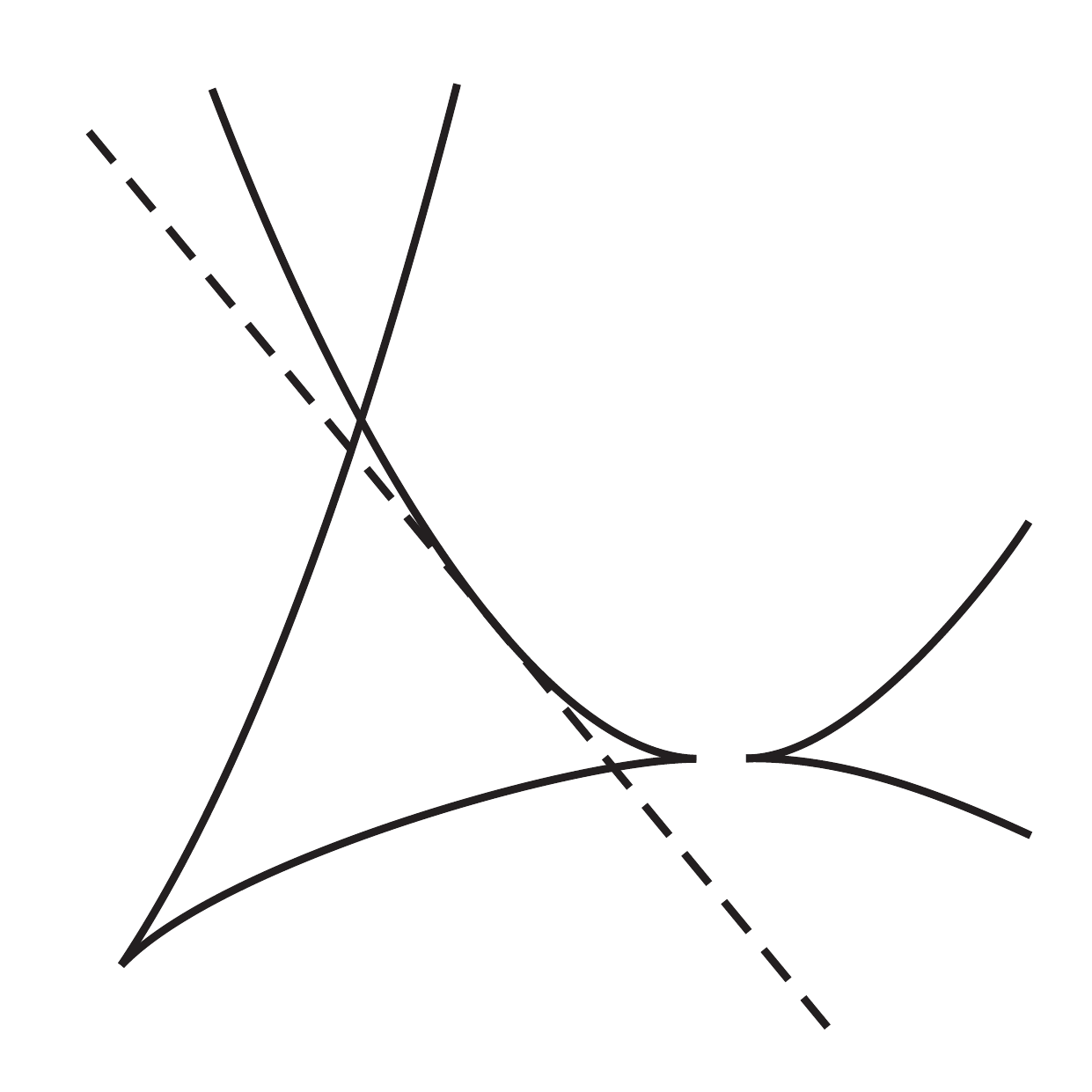}&
\includegraphics[width=4cm]{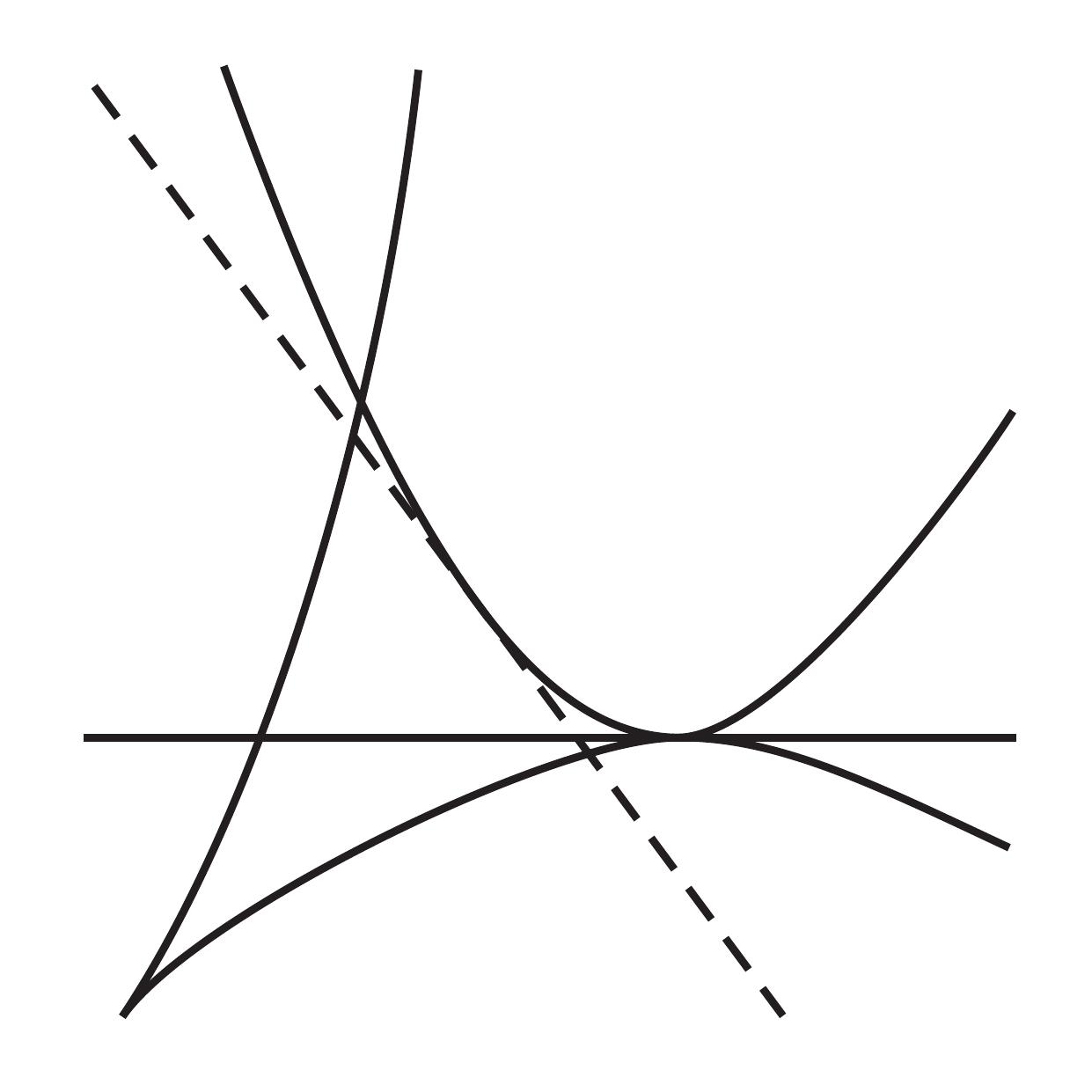}&
\includegraphics[width=4cm]{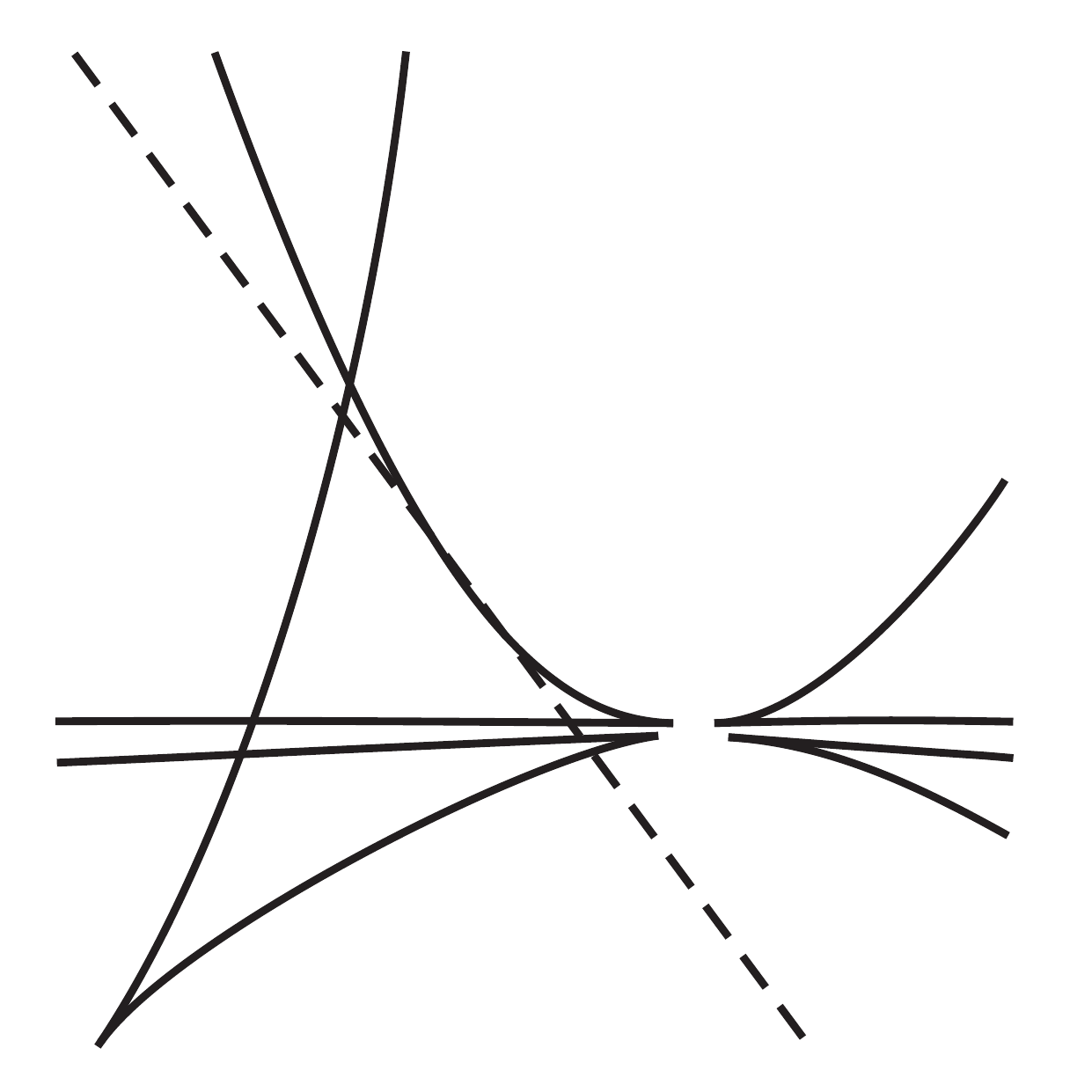}\\
\end{tabular}
\caption{The perestroika ${\cal P}_9$.}
\label{perestr62-61}
\end{center}
\end{figure}

\medskip
${\cal P}_9$) The perestroika in a neighbourhood of the point (\ref{spec2}) when we cross the curve
\begin{equation}
\gamma_1^+: \delta t_1<\sqrt[5]{\frac{S_4^4}{432}}\quad (\mbox{bounded by (\ref{vozvrat})}).
\label{gamma1+bv}
\end{equation}
It is shown by five triptychs in Fig. \ref{perestr62-61} up to diffeomorphism. The first triptych shows the perestroika in a sufficiently small neighbourhood of (\ref{spec2}). The others demonstrate it in a neighbourhood of (\ref{spec2}) that contains a nonvanishing cusp of the curve (\ref{psi-2}). Namely the second, third and fourth triptychs show the perestroika when we cross the curve $\gamma_1^+$ between the points (\ref{vozvrat}) and (\ref{trper-xi5-g1}), (\ref{trper-xi5-g1}) and (\ref{spoint5}), (\ref{spoint5}) and (\ref{spoint8}), respectively. The fifth triptych demonstrate the perestroika when we cross the curve $\gamma_1^+: \delta t_1<\frac13 \sqrt[5]{\frac{S_4^4}{12}}$ (bounded by (\ref{spoint8})). The line $(\ref{crit-tochki})$ is shown dashed.

The description of the first triptych is similar to the case ${\cal P}_7$. When choosing a suitable perestroika direction, two new semicubical cusps of the curve (\ref{psi-2}) appear; each of two-dimensional strata of $\mathfrak{S}_{S_4,t_1,q_5}$ having only one vanishing cusp in its boundary (there are two such strata) is divided into two ones. In addition, the line (\ref{psi-4}) transversally intersects two smooth branches of the curve (\ref{psi-2}). One of them is shown as a branch of a semicubical parabola with non-vanishing cusp. The second branch is not shown in Fig. \ref{perestr62-61}. It intersects (\ref{psi-4}) at a point lying on the other side of the line $(\ref{crit-tochki})$ with respect to the point (\ref{spec2}). Hence, there are three more two-dimensional strata, each of which is divided into three ones.

\medskip
The {\it third type} perestroikas happens in a neighbourhood of infinity\footnote{A neighbourhood of infinity in $\mathbb{R}^n$ is the complement to a compact subset.} in $\mathbb{R}^2=\{(t_2,S_3)\}$ when we cross the coordinate axes of $\mathbb{R}^2=\{(t_1,q_5)\}$. We describe these perestroikas on the plane with coordinates $\widetilde{t}_2,\widetilde{S}_3$ and the origin $\infty$ using the inversion
$$
\widetilde{t}_2=\frac{t_2}{t_2^2+S_3^2},\quad \widetilde{S}_3=\frac{S_3}{t_2^2+S_3^2}
$$
centered at zero of the original coordinates $t_2,S_3$.

The equation (\ref{crit-tochki}) in the coordinates $\widetilde{t}_2,\widetilde{S}_3$ defines a circle with a punctured point $\infty$. The tangent line to this circle at $\infty$ is given by the equation $\widetilde{S}_3=0$. The closure $\gamma$ of the curve $(\ref{psi-2})$ as a subset of $\mathbb{R}^2=\{(\widetilde{t}_2,\widetilde{S}_3)\}$ is an algebraic curve passing through $\infty$.

\medskip
\lemma\label{per-infty} {\it Let $t_1q_5\neq0,S_4^2+12q_5t_1^2>0$. Then a germ $\gamma$ at $\infty$ has: 

{\rm 1)} two smooth branches that are transversal to each other and to the lines $\widetilde{t}_2=0,\widetilde{S}_3=0$; 

{\rm 2)} one singular branch that has semicubical cusp and the tangent line $\widetilde{S}_3=0$ at $\infty$; 

{\rm 3)} two singular branches with semicubical cusps and the tangent line $\widetilde{t}_2=0$ at $\infty$ if $q_5<0$; the branches of these two semicubical parabolas alternate {\rm(}each domain bounded by one of parabolas contains exactly one branch of another one{\rm)}.}

{\sc Proof.} The slope $\kappa$ of the line $\widetilde{S}_3=\kappa\widetilde{t}_2$ intersecting the curve $\gamma$ at a point $(\widetilde{t}_2,\widetilde{S}_3)\neq0$ satisfies the equation
$$
q_5\kappa^2(q_5\kappa^2-12\delta S_4\kappa-432t_1^2)=41472\delta\widetilde{t}_2(q_5S_4t_1^3+\dots),
$$
where dots denote a polynomial in $\kappa,\widetilde{t}_2$ that vanishes at $\kappa=\widetilde{t}_2=0$ (its coefficients are polynomials in $t_1,S_4,q_5$). If the line $\widetilde{S}_3=\kappa\widetilde{t}_2$ is tangent to $\gamma$ at $\infty$, then $\kappa$ is a real root of the equation
$$
\kappa^2(q_5\kappa^2-12\delta S_4\kappa-432t_1^2)=0.
$$
In the case under consideration, this equation has two simple non-zero roots and one zero root with multiplicity $2$. The first two statements of the lemma follow from this.

Now we take a line $\widetilde{t}_2=r\widetilde{S}_3$ with the slope $r$ close to zero. Let $\widetilde{S}_3=9\delta r^2+s$, and the line $\widetilde{t}_2=r\widetilde{S}_3$ intersect $\gamma$ at $(\widetilde{t}_2,\widetilde{S}_3)\neq0$. Then
$$
s^2+87480q_5^3r^6+\dots=0,
$$
where dots denote a polynomial in $r,s$ consisting of monomials with quasi-degree greater than $6$ in the filtration $w(r)=1,w(s)=3$ (its coefficients are polynomials in $t_1,S_4,q_5$ as well). For $q_5\neq0$ this equation has solutions $(r,s)\in\mathbb{R}^2$ in a neighborhood of the origin if and only if $q_5<0$. This implies the third statement of the lemma. $\Box$

\medskip
There are two third type perestroikas.

\begin{figure}[h]
\begin{center}
\begin{tabular}{ccc}
\includegraphics[width=4cm]{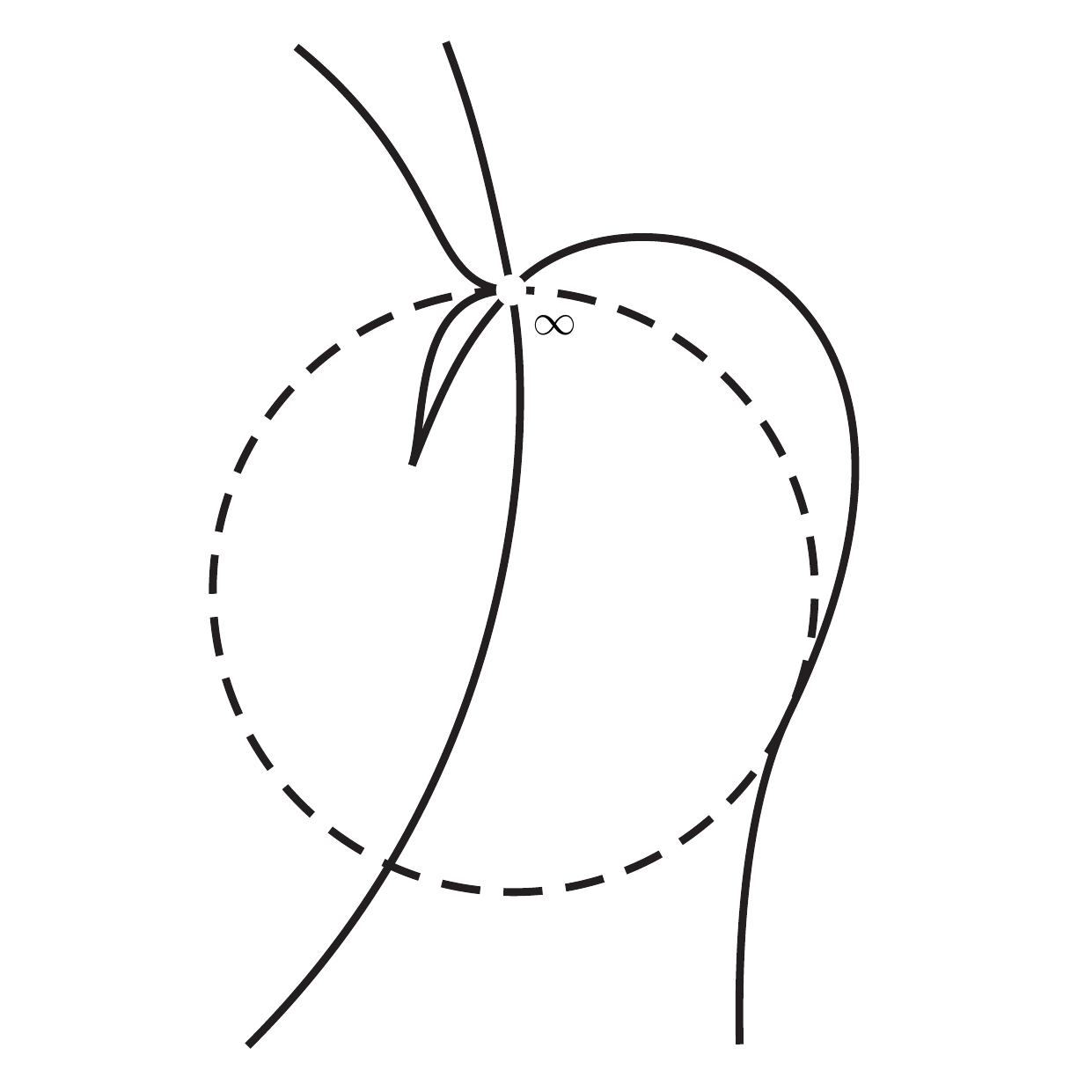}&
\includegraphics[width=4cm]{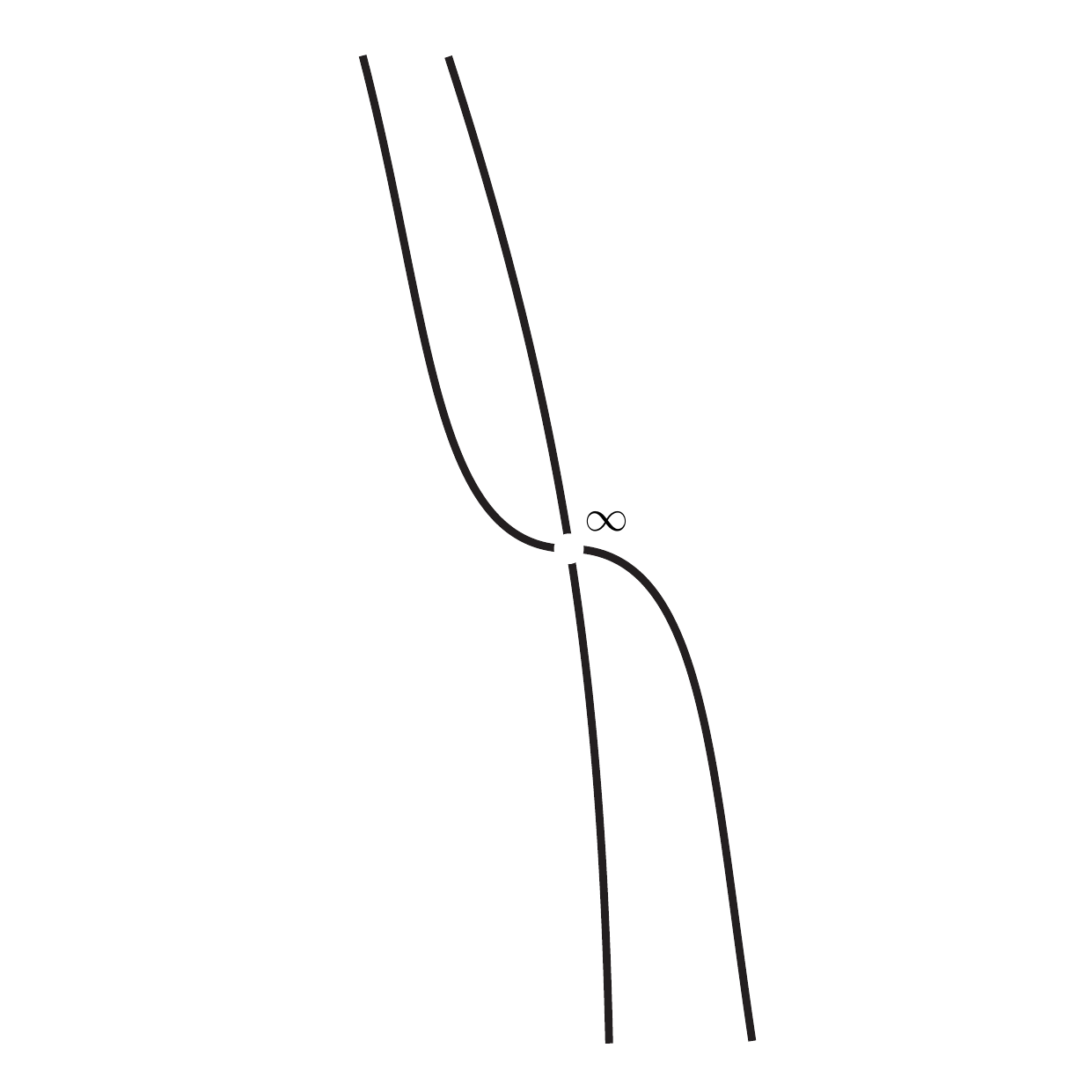}&
\includegraphics[width=4cm]{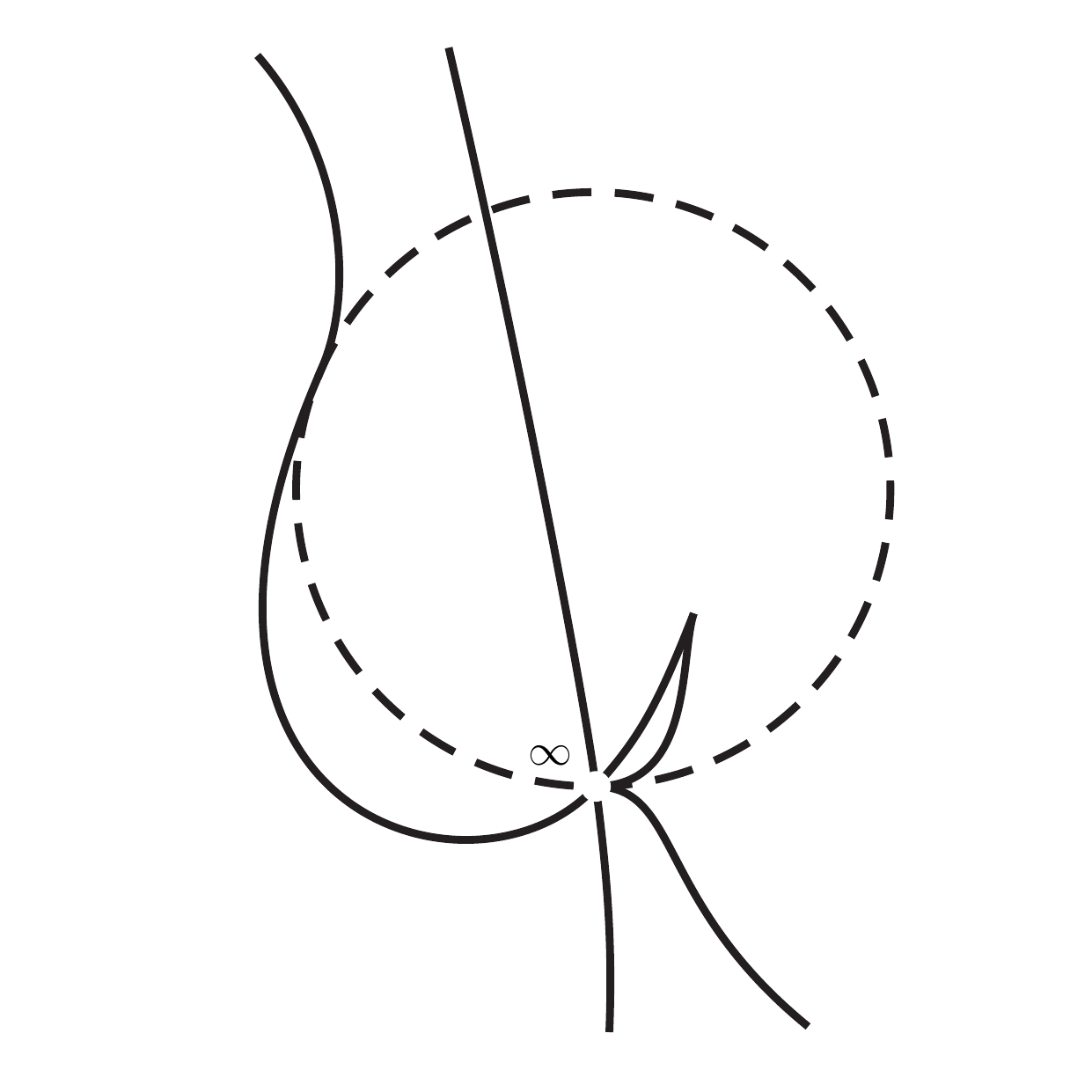}\\
\includegraphics[width=4cm]{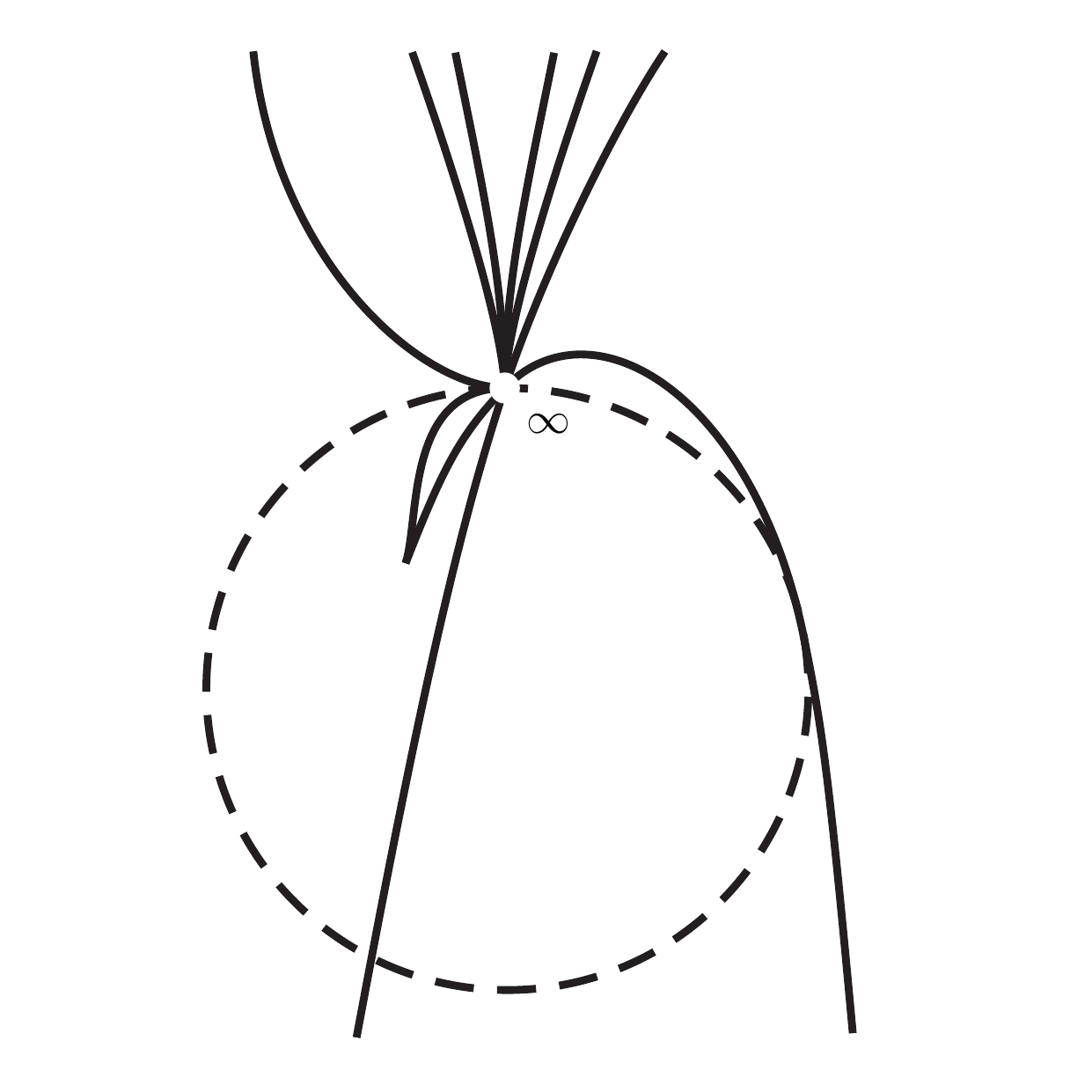}&
\includegraphics[width=4cm]{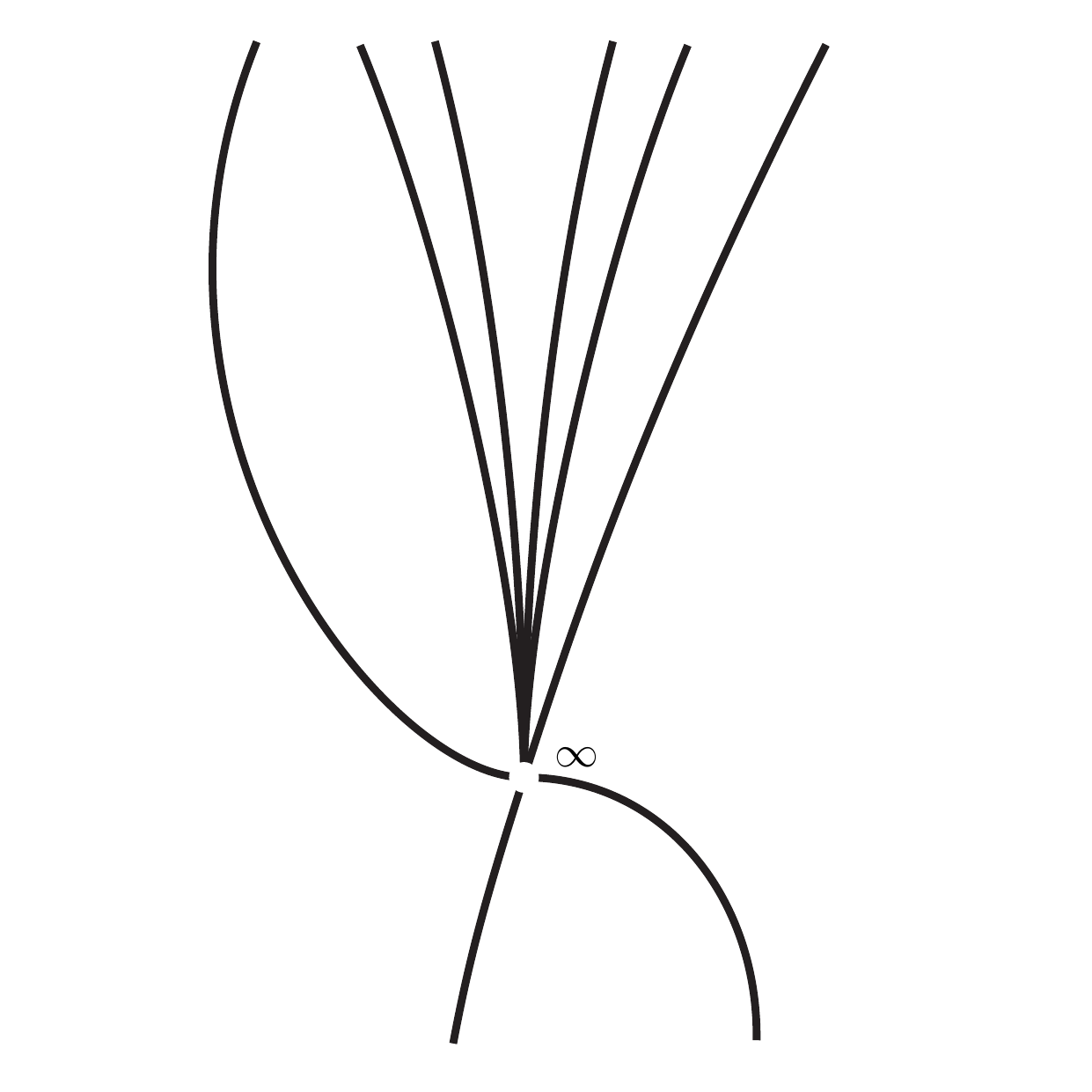}&
\includegraphics[width=4cm]{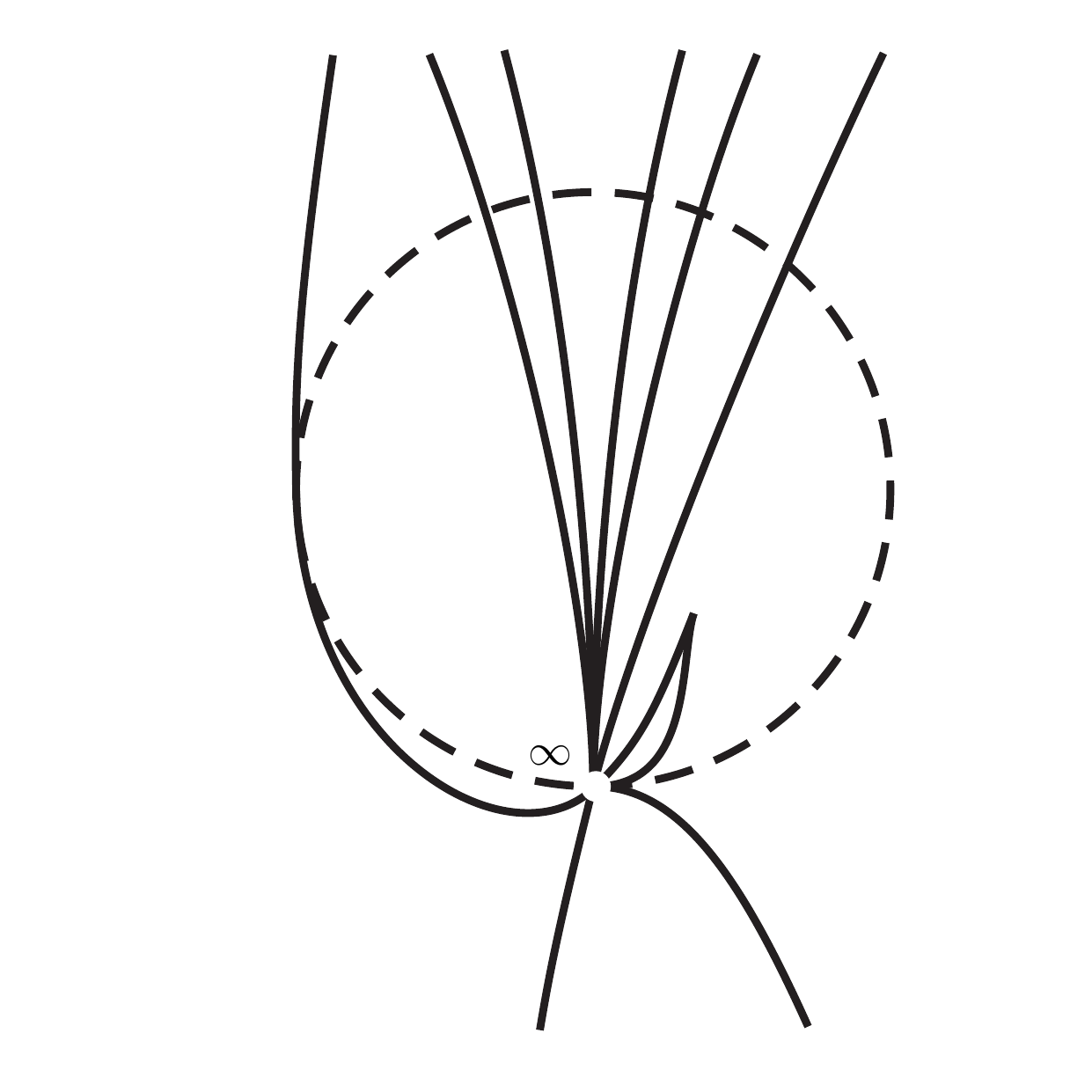}\\
\end{tabular}
\caption{The perestroika ${\cal P}_{10}$.}
\label{perestr28-29}
\end{center}
\end{figure}

\begin{figure}[h]
\begin{center}
\begin{tabular}{ccc}
\includegraphics[width=4cm]{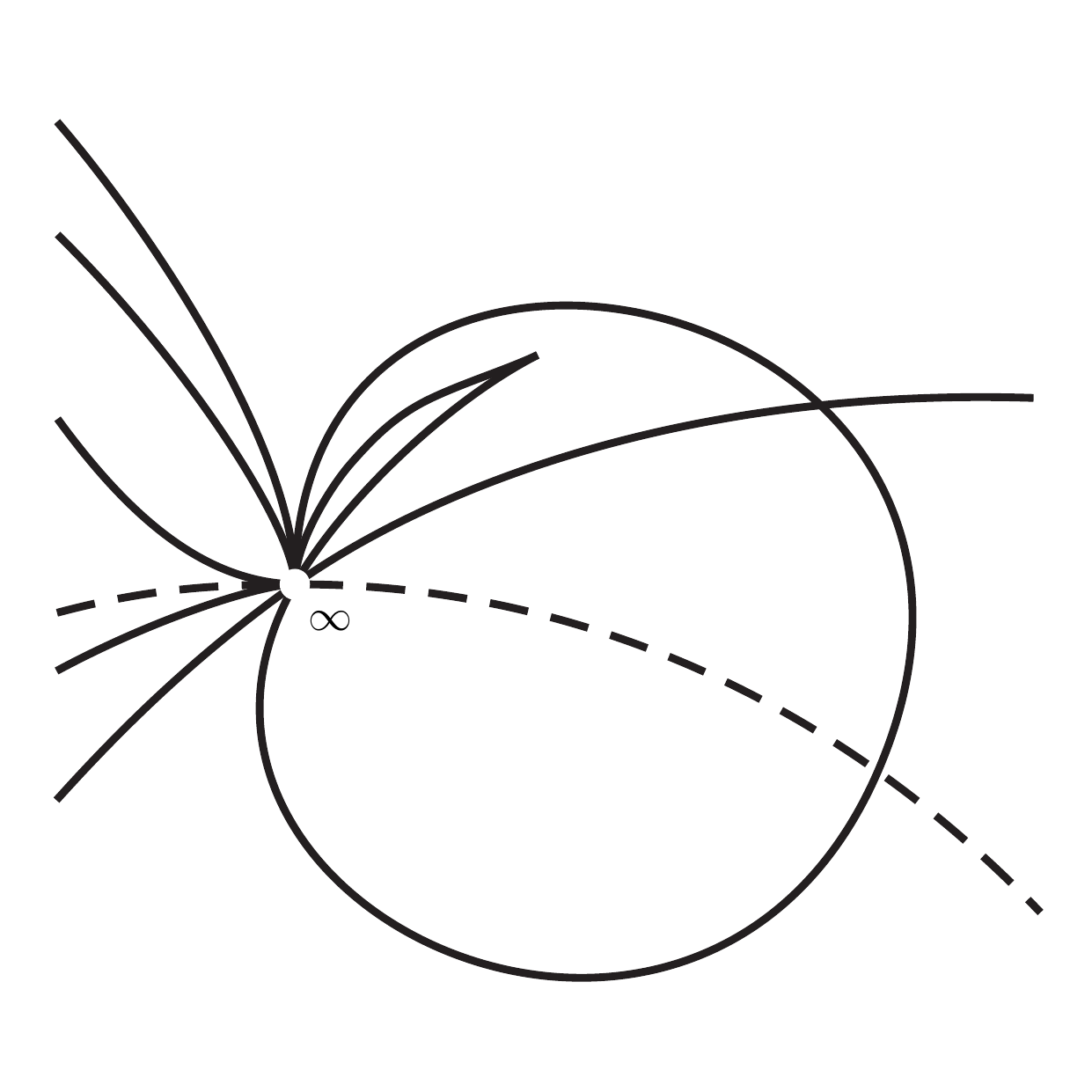}&
\includegraphics[width=4cm]{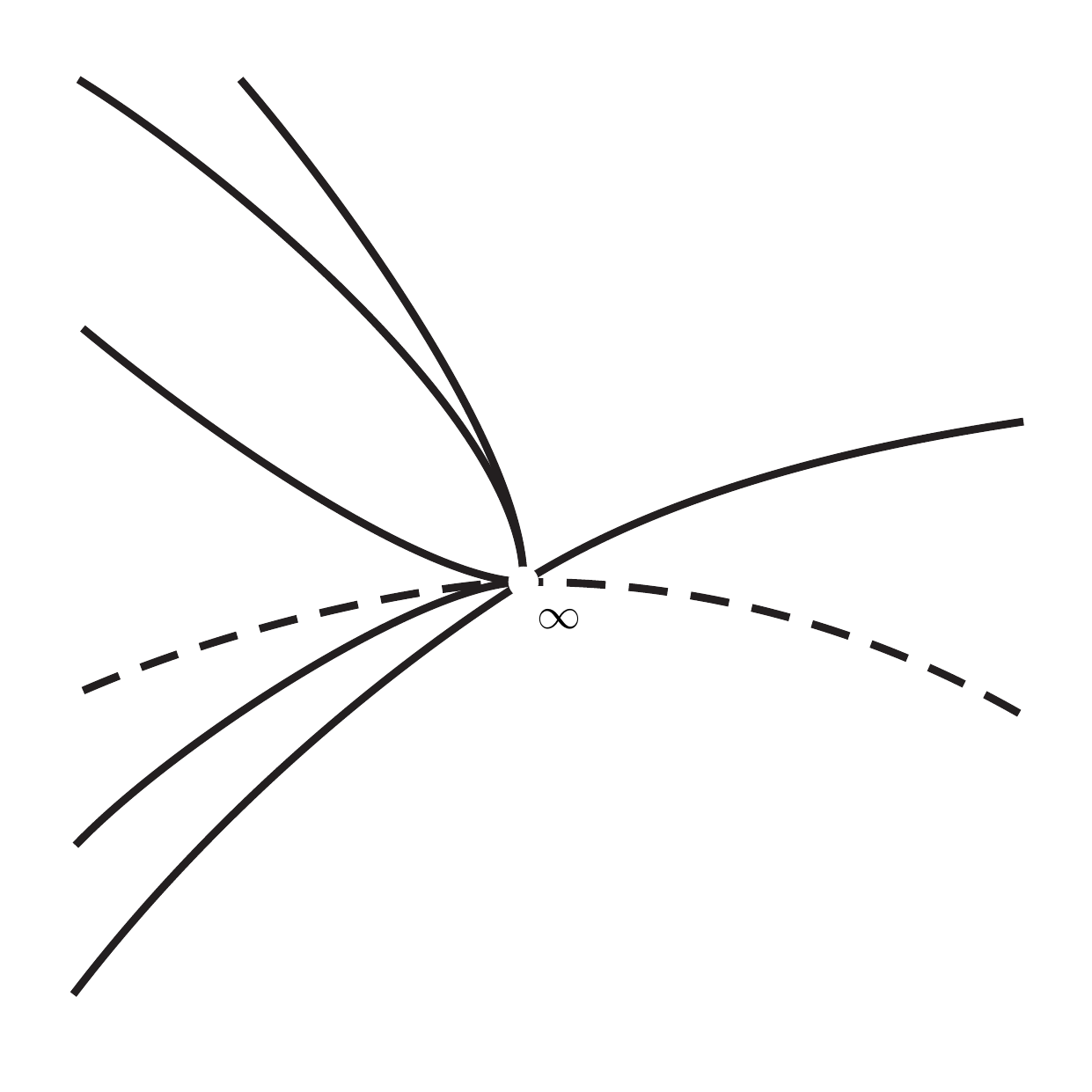}&
\includegraphics[width=4cm]{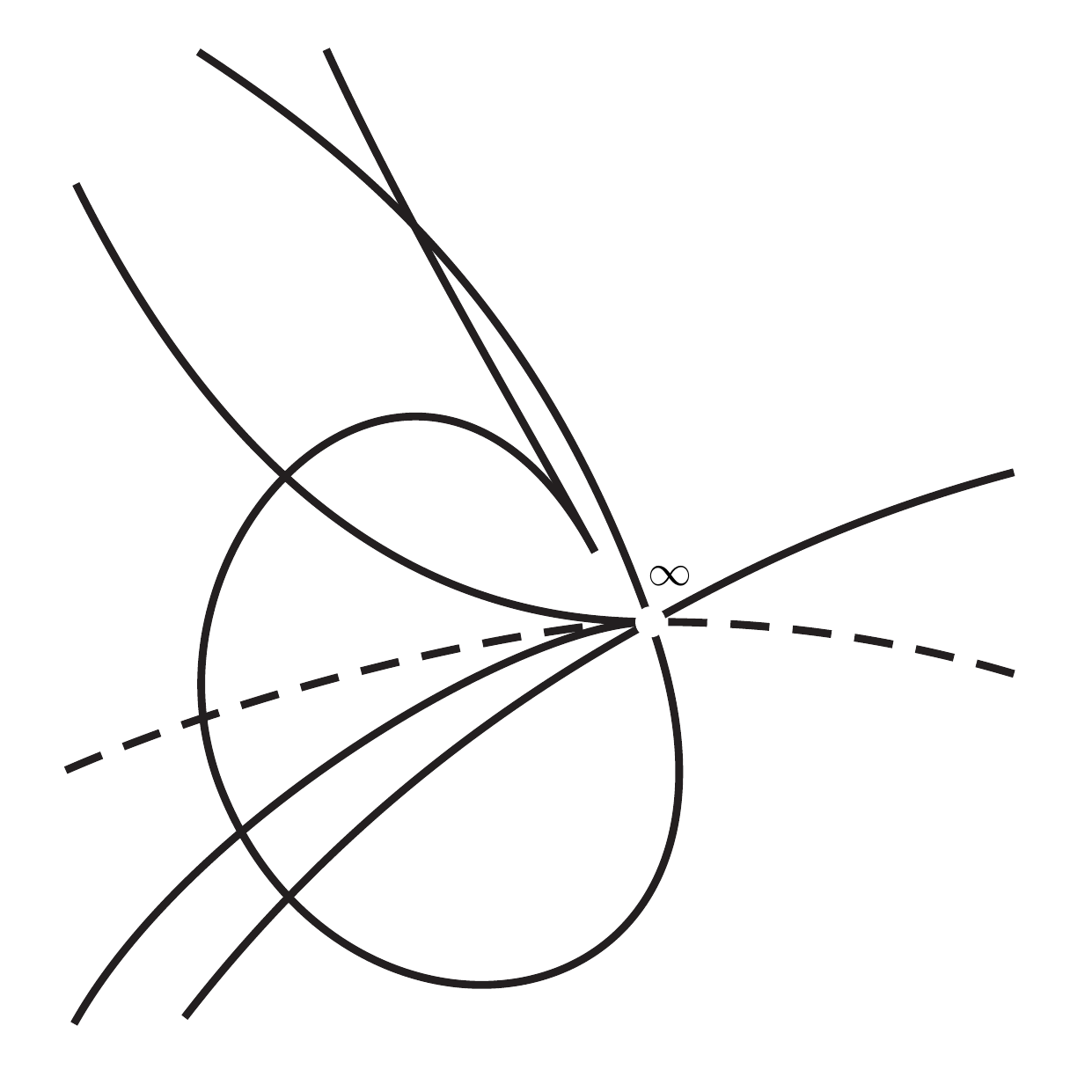}\\
\includegraphics[width=4cm]{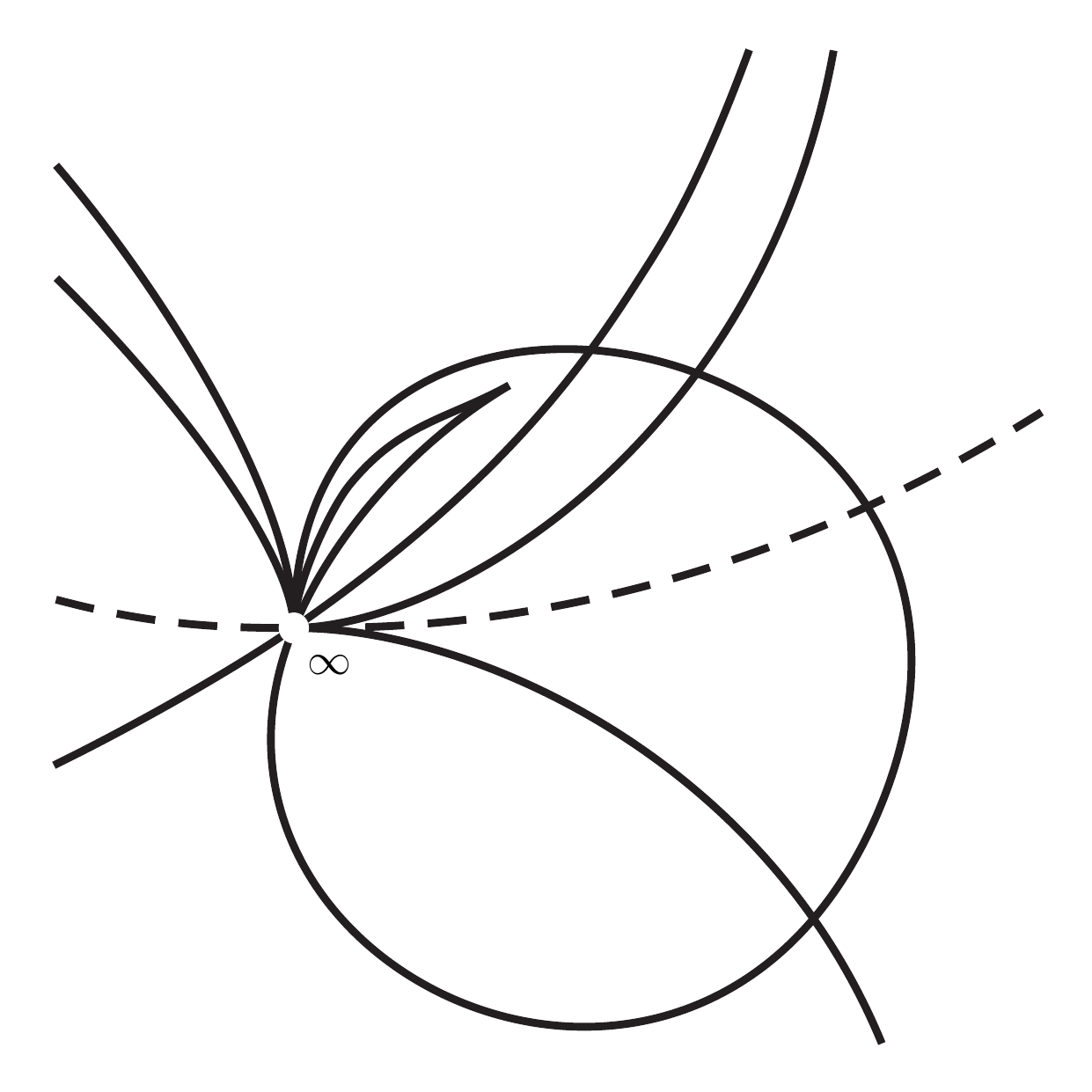}&
\includegraphics[width=4cm]{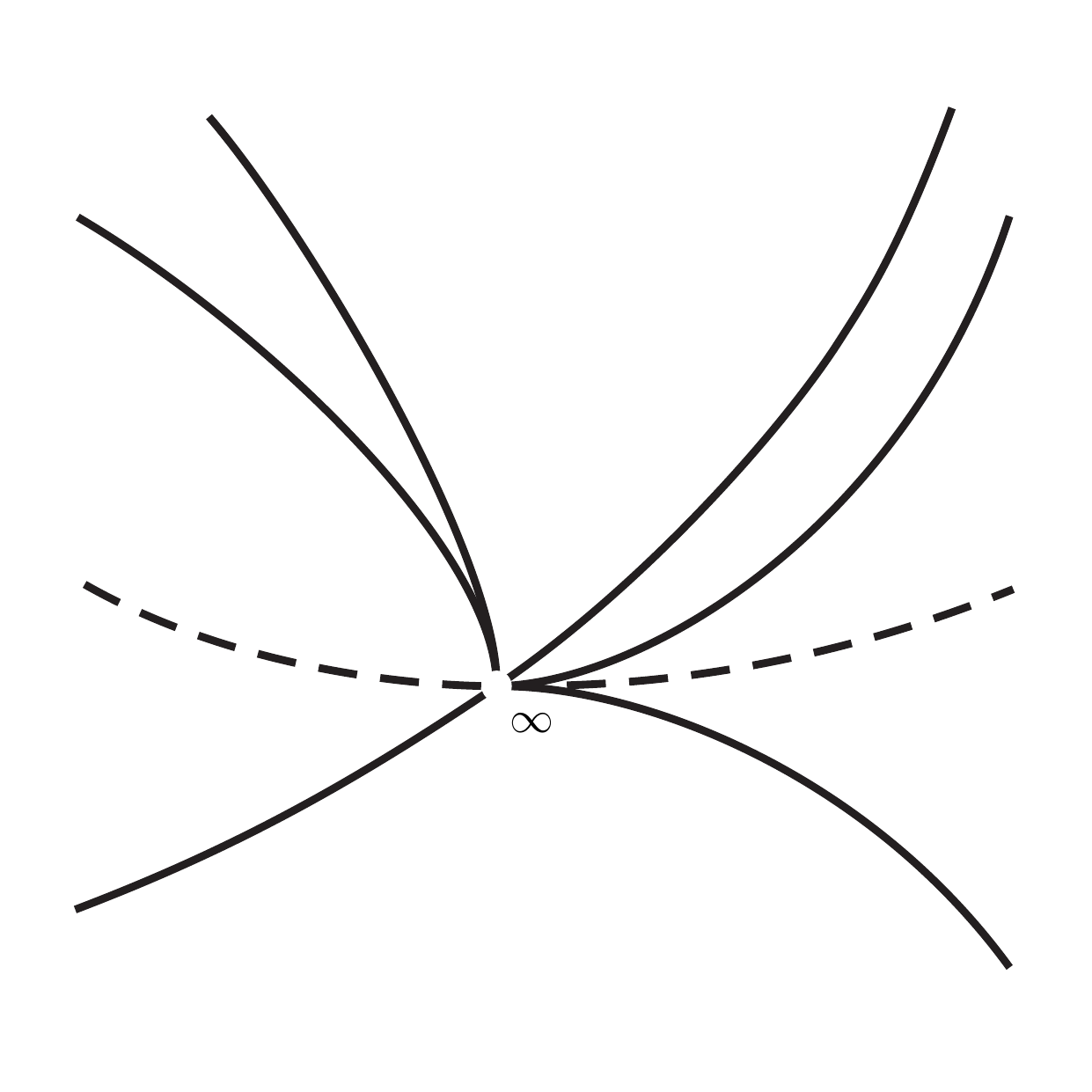}&
\includegraphics[width=4cm]{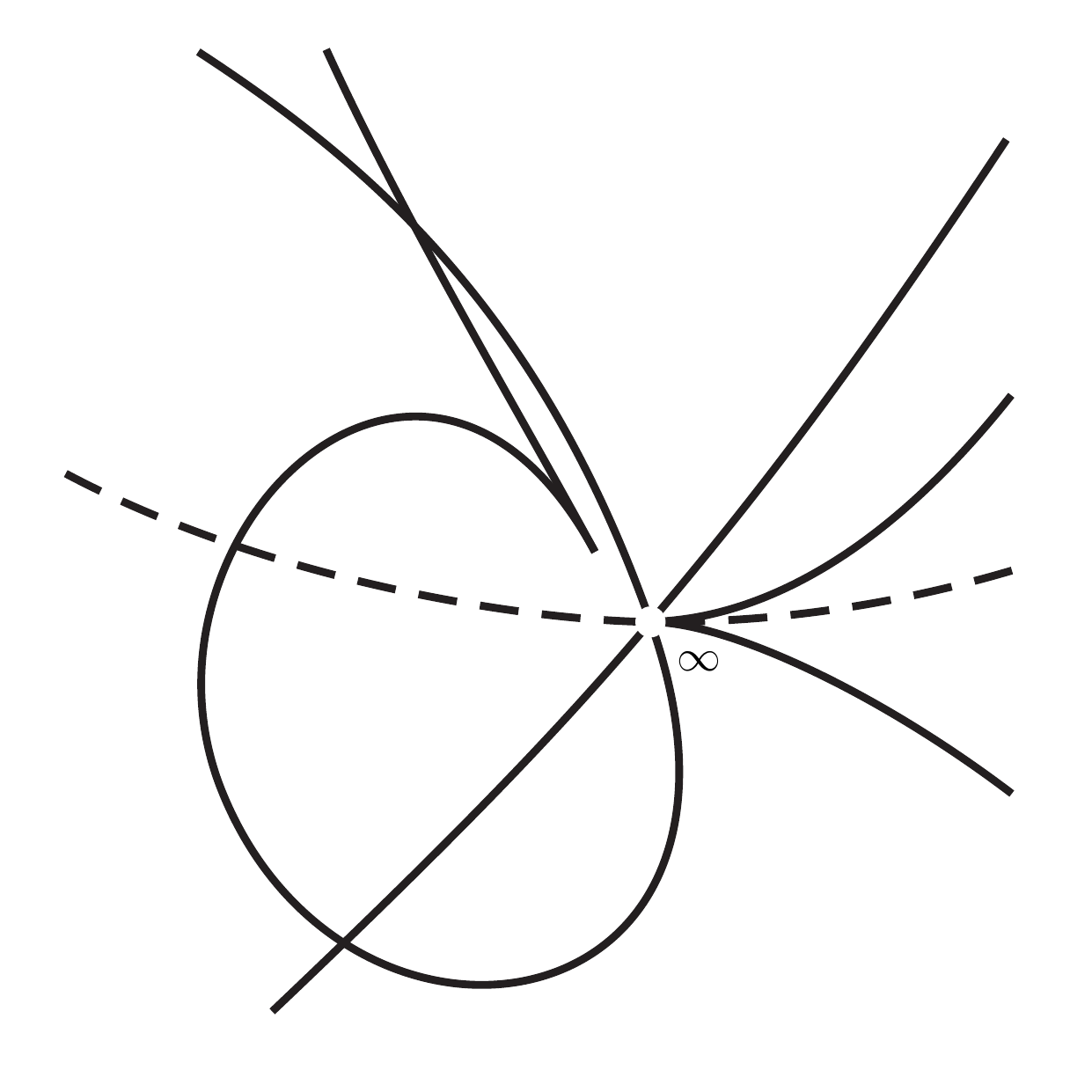}\\
\end{tabular}
\caption{The perestroika ${\cal P}_{11}$.}
\label{perestr11-14}
\end{center}
\end{figure}

\medskip
${\cal P}_{10}$) The perestroika in a neighbourhood of infinity when we cross the $q_5$ axis. It is shown in coordinates $\widetilde{t}_2,\widetilde{S}_3$ by two triptychs in Fig. \ref{perestr28-29} up to diffeomorphism. The left-to-right direction corresponds to the increasing of the parameter $\delta t_1$. The line (\ref{crit-tochki}) is represented by a dashed circle. The point $\infty$ is the punctured dot.

The first triptych shows the perestroika when we cross the ray
\begin{equation}
t_1=0 : q_5>0.
\label{luchq5+}
\end{equation}
In the left and right figures we see two smooth branches of the curve $\gamma$, which transversally intersect the circle at $\infty$. One of them transversally intersects the circle at one more point. The other has a simple tangency with this circle. It is also a branch of a semicubical parabola with cusp lying inside the disk bounded by the circle. This cusp together with the circle and the point of tangency of the circle with $\gamma$ disappear at $\infty$ as $\delta t_1\rightarrow-0$ (the radius of the circle tends to zero). They appear again for $\delta t_1>0$ but on the other side of the tangent line $\widetilde{S}_3=0$.

The second triptych demonstrates the perestroika when we cross the ray
\begin{equation}
t_1=0 : q_5<0.
\label{luchq5-}
\end{equation}
This case differs from the previous one by two additional non-vanishing singular branches of $\gamma$, each of which has a semicubical cusp at $\infty$.

\medskip
${\cal P}_{11}$) The perestroika in a neighbourhood of infinity when we cross the $t_1$ axis. It is shown in coordinates $\widetilde{t}_2,\widetilde{S}_3$ by two triptychs in Fig. \ref{perestr11-14} up to diffeomorphism. The left-to-right direction corresponds to the increasing of the parameter $q_5$. The line (\ref{crit-tochki}) is also represented by a dashed circle. The point $\infty$ is the punctured dot.

The first triptych shows the perestroika when we cross the ray
\begin{equation}
q_5=0: \delta t_1<0.
\label{lucht1-}
\end{equation}
In the left figure we see a loop with one singular point. It is obtained from a smooth arc of curve (\ref{psi-2}) by adding $\infty$. The loop enters $\infty$ on the one hand as a branch of a semicubical parabola that is transversal to the circle and on the other as a smooth branch of the curve $\gamma$. This smooth branch  is extended inside the disk bounded by the loop to a semicubical cusp. The loop together with this cusp disappear at $\infty$ as $q_5\rightarrow-0$. A loop appears again for $q_5>0$ but now it has two singular points: a semicubical cusp and a transversal intersection point of two smooth branches. The point $\infty$ is non-singular for this loop.

The second triptych demonstrates the perestroika when we cross the ray
\begin{equation}
q_5=0 : \delta t_1>0.
\label{lucht1+}
\end{equation}
This case differs from the previous one in the location of a singular branch of the curve $\gamma$ that touches the dashed circle at the point $\infty$. Namely in the left figure it twice transversally intersects the smooth part of the vanishing loop having one singular point.

\section{Type $A_1^k$ strata of the stratification $\mathfrak{S}_{S_4},S_4\neq0$}

The skeletons $\Sigma_{S_4,t_1,q_5},S_4\neq0$ for points $(t_1,q_5)$ of domains $1-93$ in Fig. \ref{razbienie} and for some singular points of the diagram $B_{S_4}$ are shown in Fig. \ref{zona1-15} -- \ref{sing-sigma} up to diffeomorphism. The line (\ref{crit-tochki}) is dashed. The $\delta t_2$ axis is horizontal. The $\delta S_3$ axis is vertical and is directed upward. If $\delta S_4<0$, then the $\delta t_2$ axis is directed to the right. Changing the sign of $S_4$ leads to a mirror reflection of $\Sigma_{S_4,t_1,q_5}$ with respect to the vertical axis.

\medskip
\remark The skeletons $\Sigma_{S_4,t_1,q_5}$ in Fig. \ref{sing-sigma} for points (\ref{spoint1}) and (\ref{spoint4}) have one non-ordinary singular point $(t_2,S_3)$, which is determined by formulas (\ref{psi-2}) at
\begin{equation}
u=-\delta\sqrt[5]{\frac{(47+21\sqrt{5})S_4}{9}}
\label{u-36}
\end{equation}
and
\begin{equation}
u=-\delta\sqrt[5]{\frac{(47-21\sqrt{5})S_4}{9}},
\label{u-37}
\end{equation}
respectively. It is a transversal intersection point of curve (\ref{psi-2}) and line (\ref{crit-tochki}). The curve (\ref{psi-2}) has a singularity of type $4/3$ at this point.

\begin{figure}
\begin{center}
\begin{tabular}{ccc}
1)\includegraphics[width=4cm]{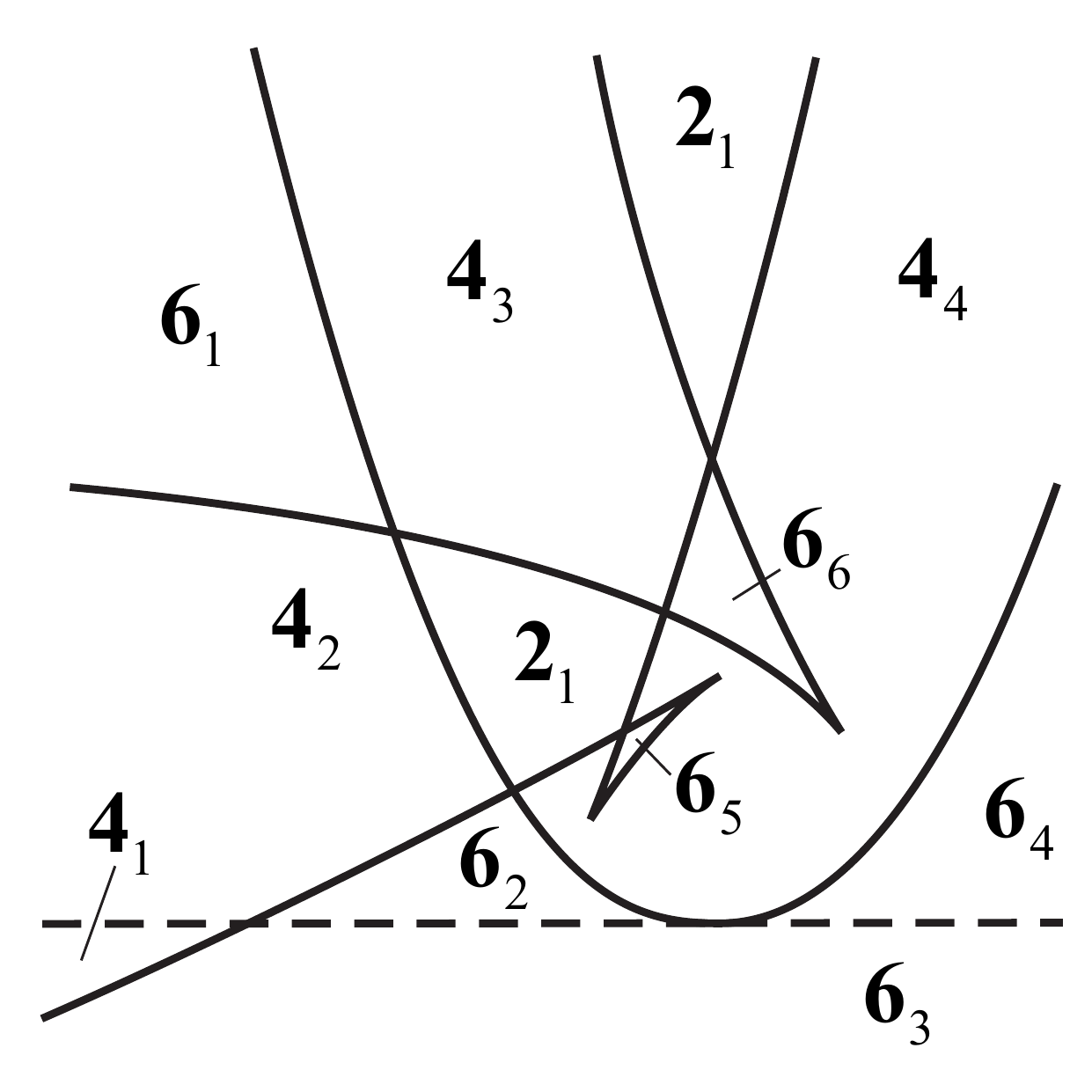}&
2)\includegraphics[width=4cm]{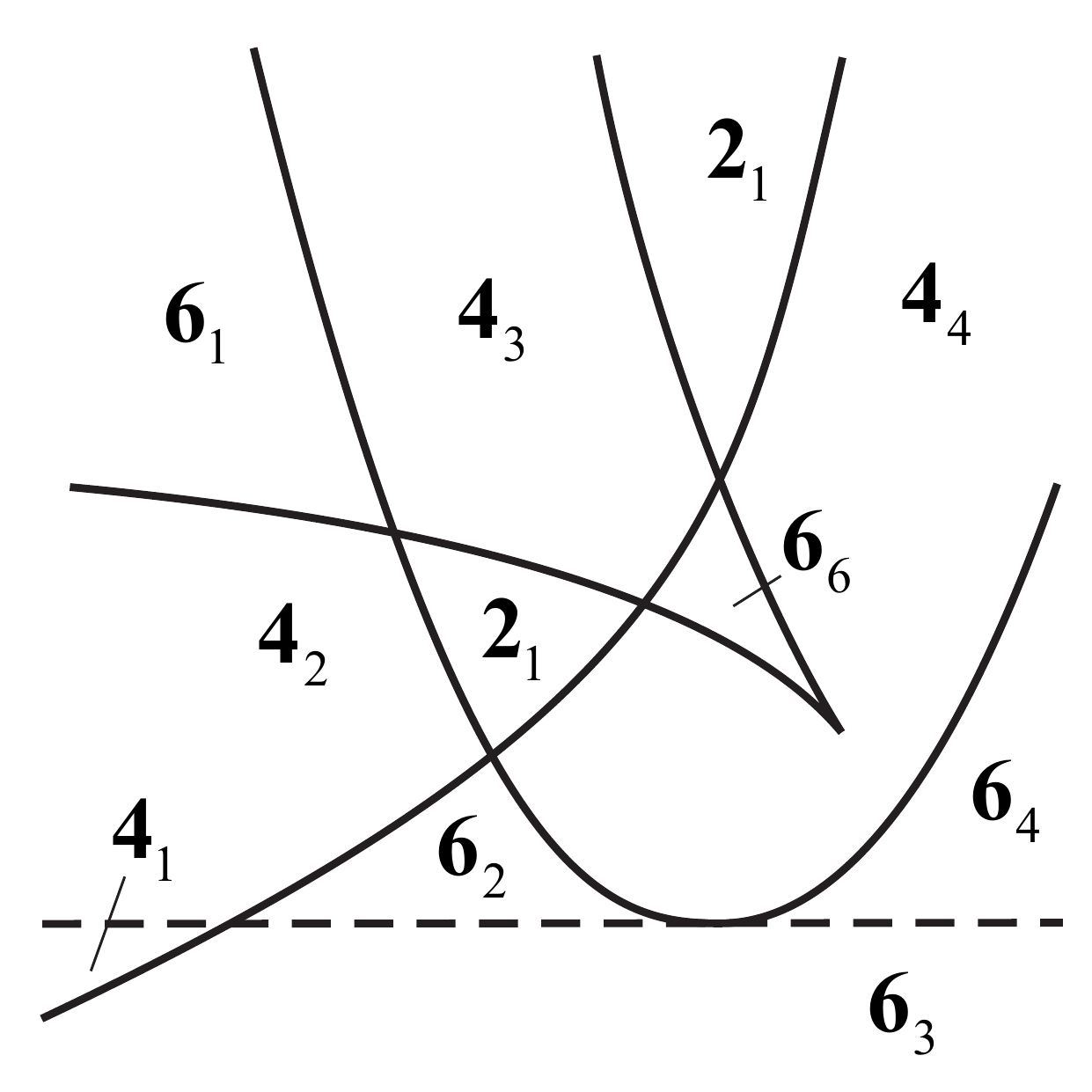}&
3)\includegraphics[width=4cm]{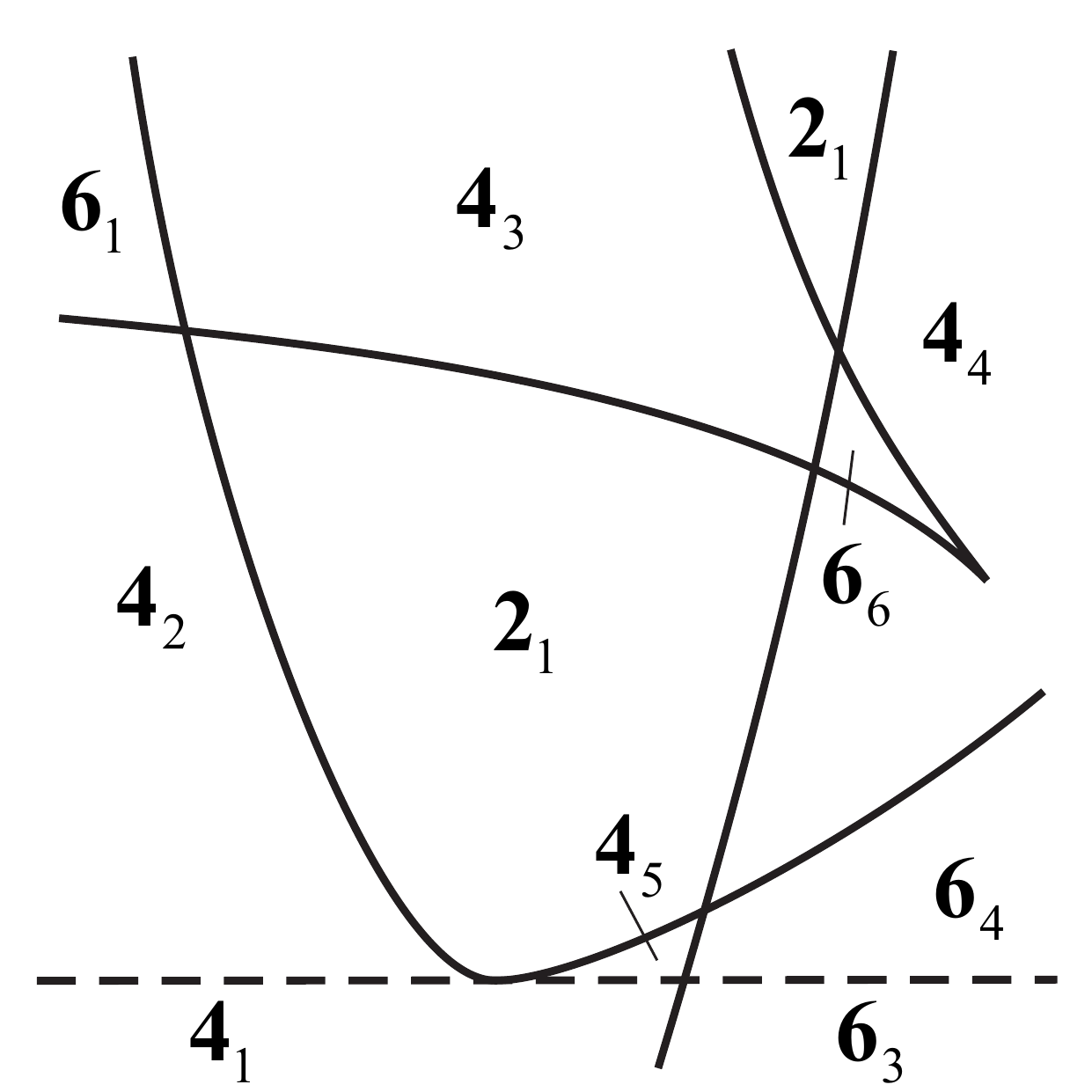}\\
4)\includegraphics[width=4cm]{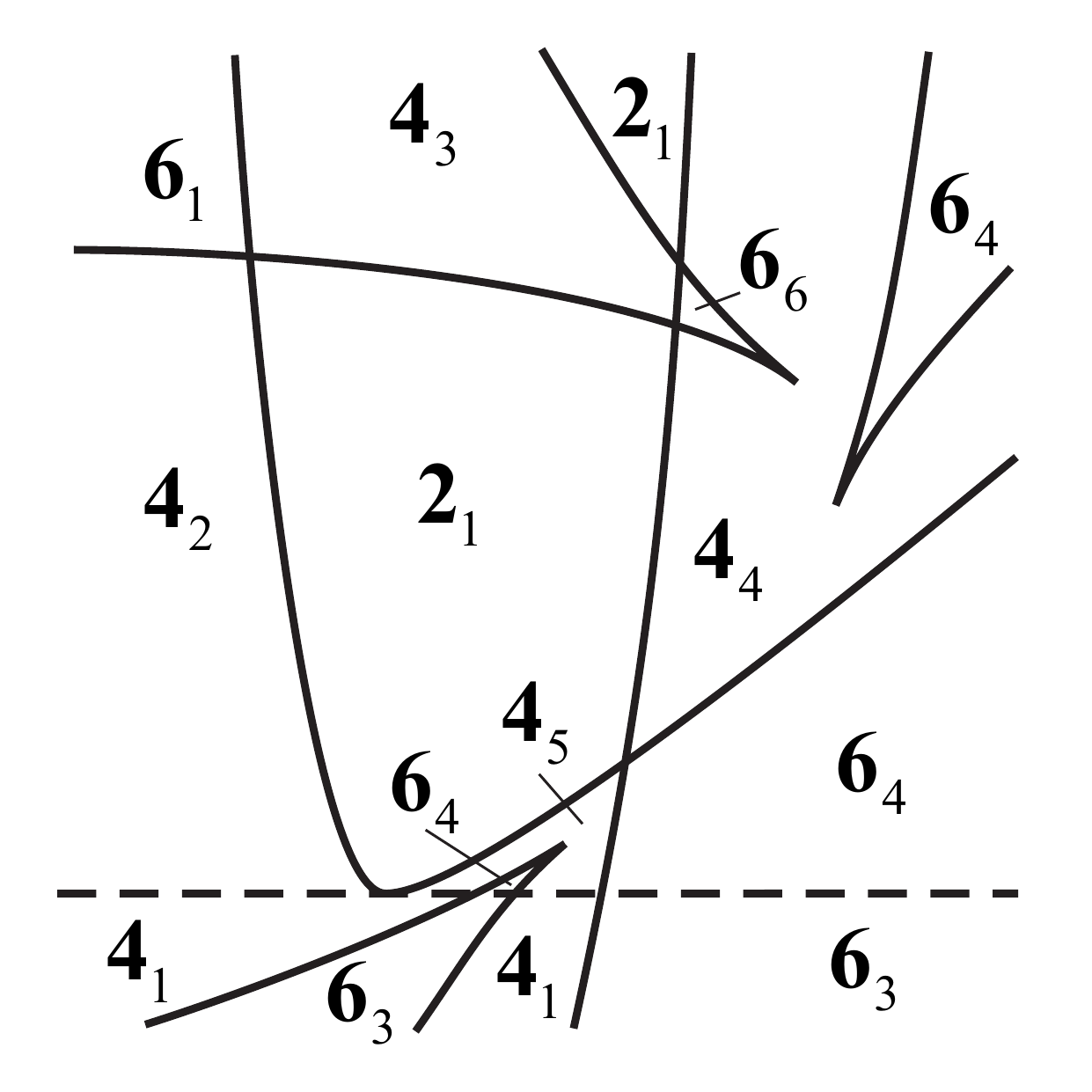}&
5)\includegraphics[width=4cm]{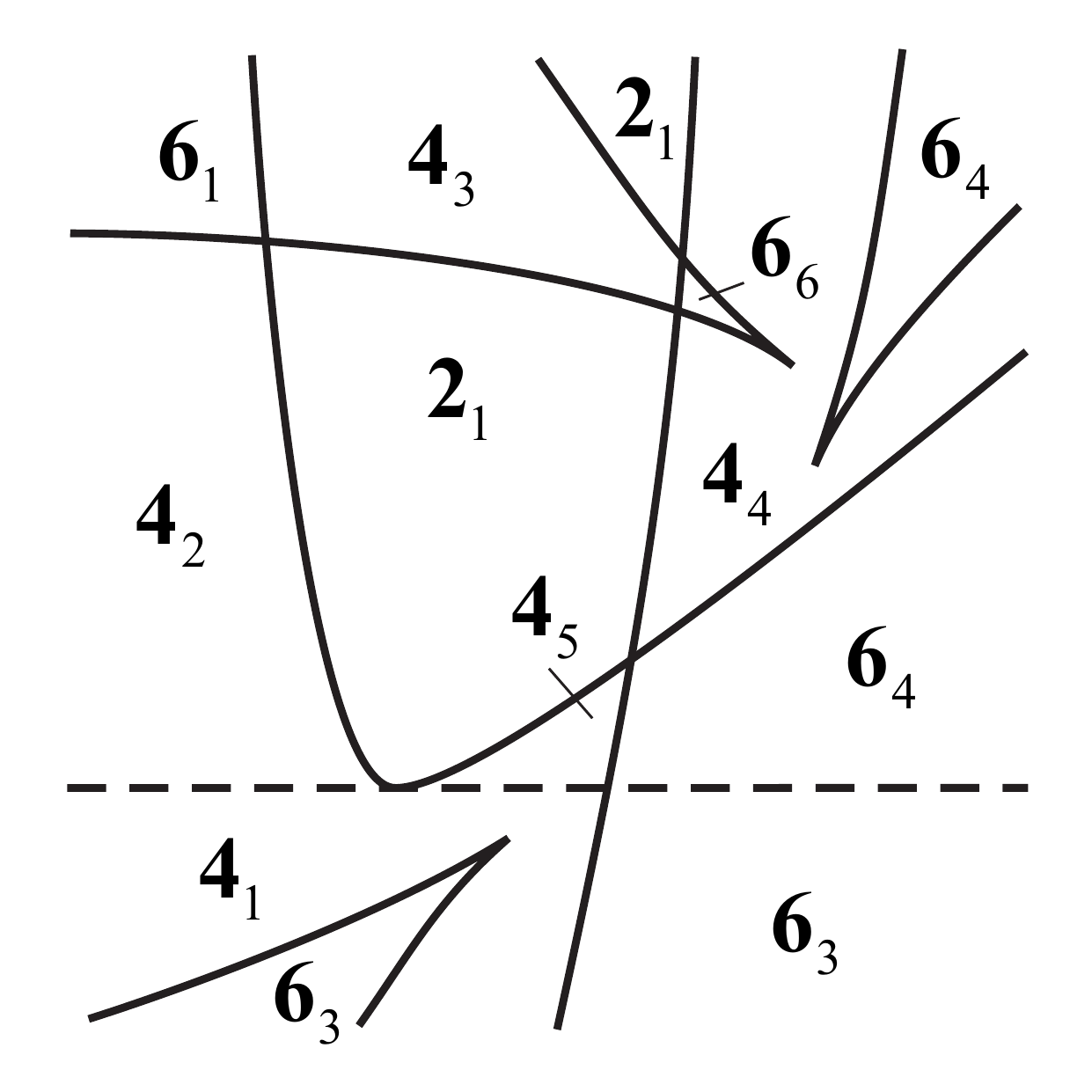}&
6)\includegraphics[width=4cm]{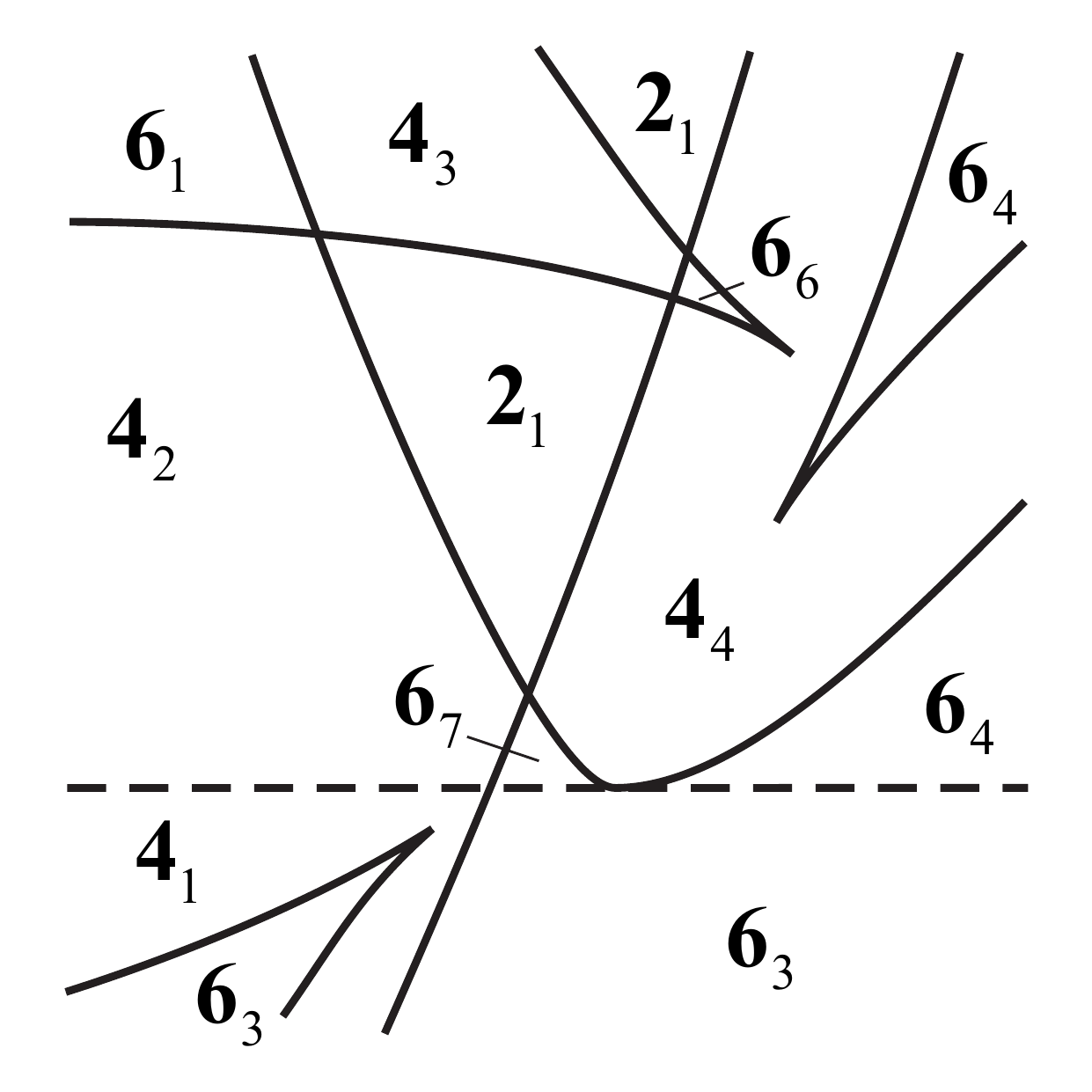}\\
7)\includegraphics[width=4cm]{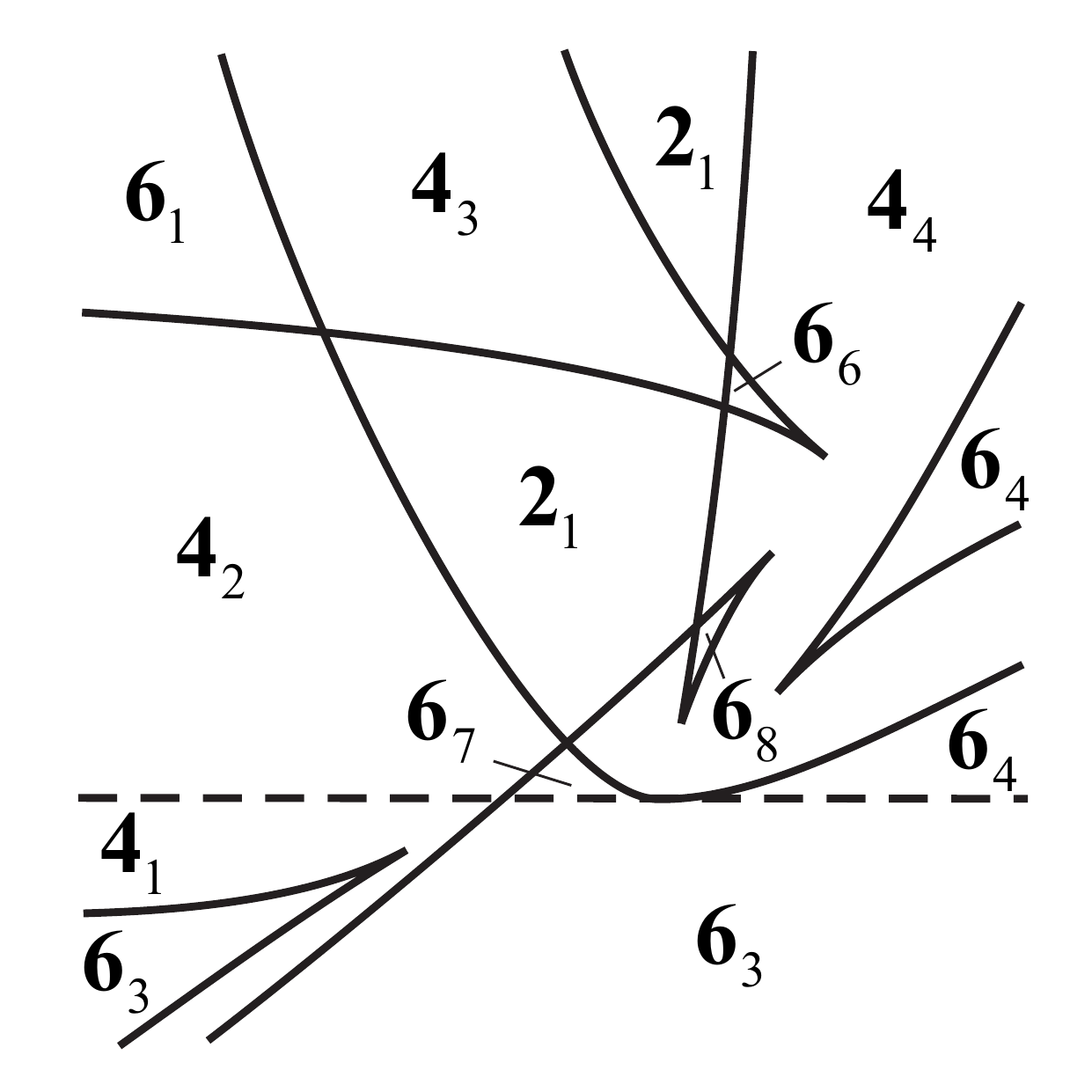}&
8)\includegraphics[width=4cm]{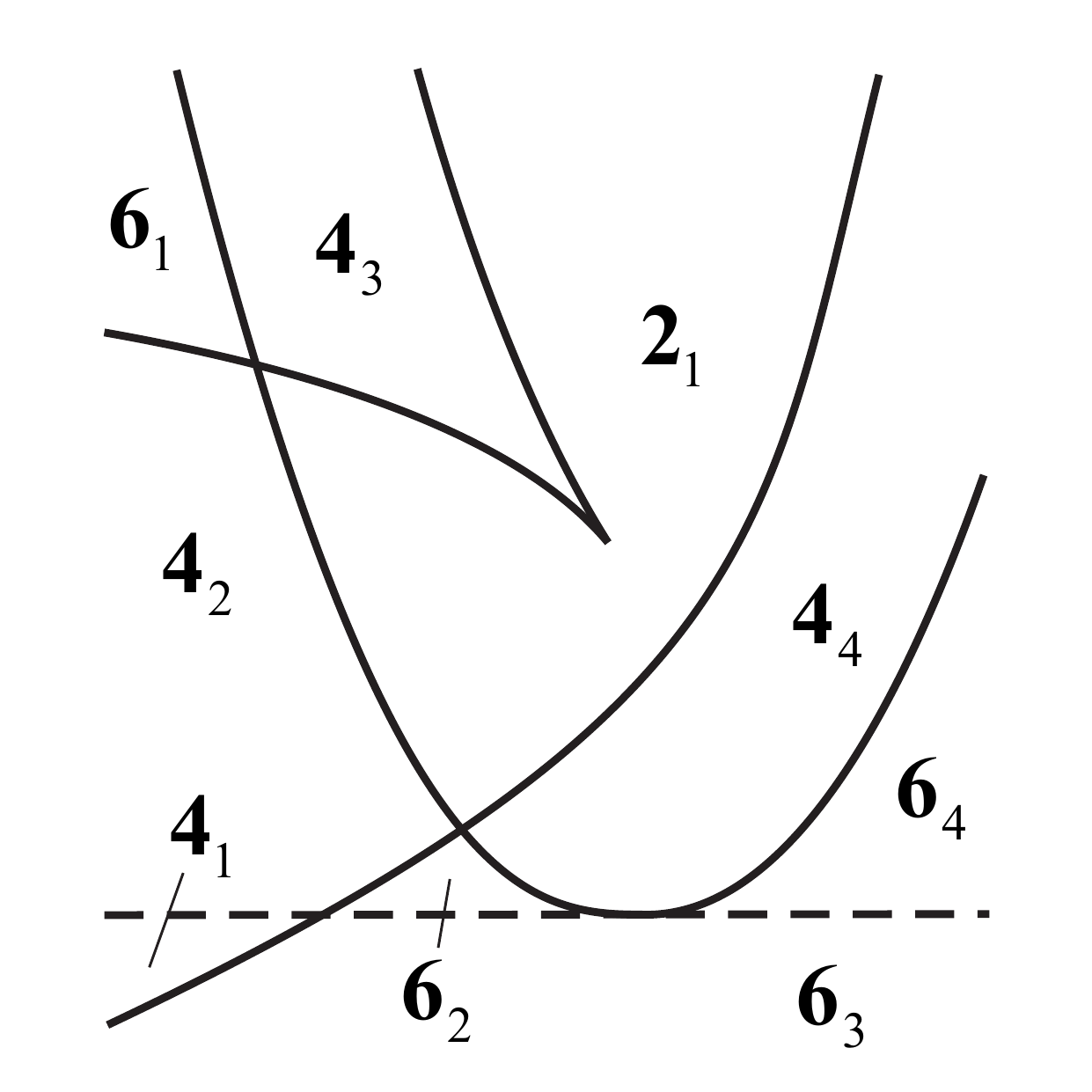}&
9)\includegraphics[width=4cm]{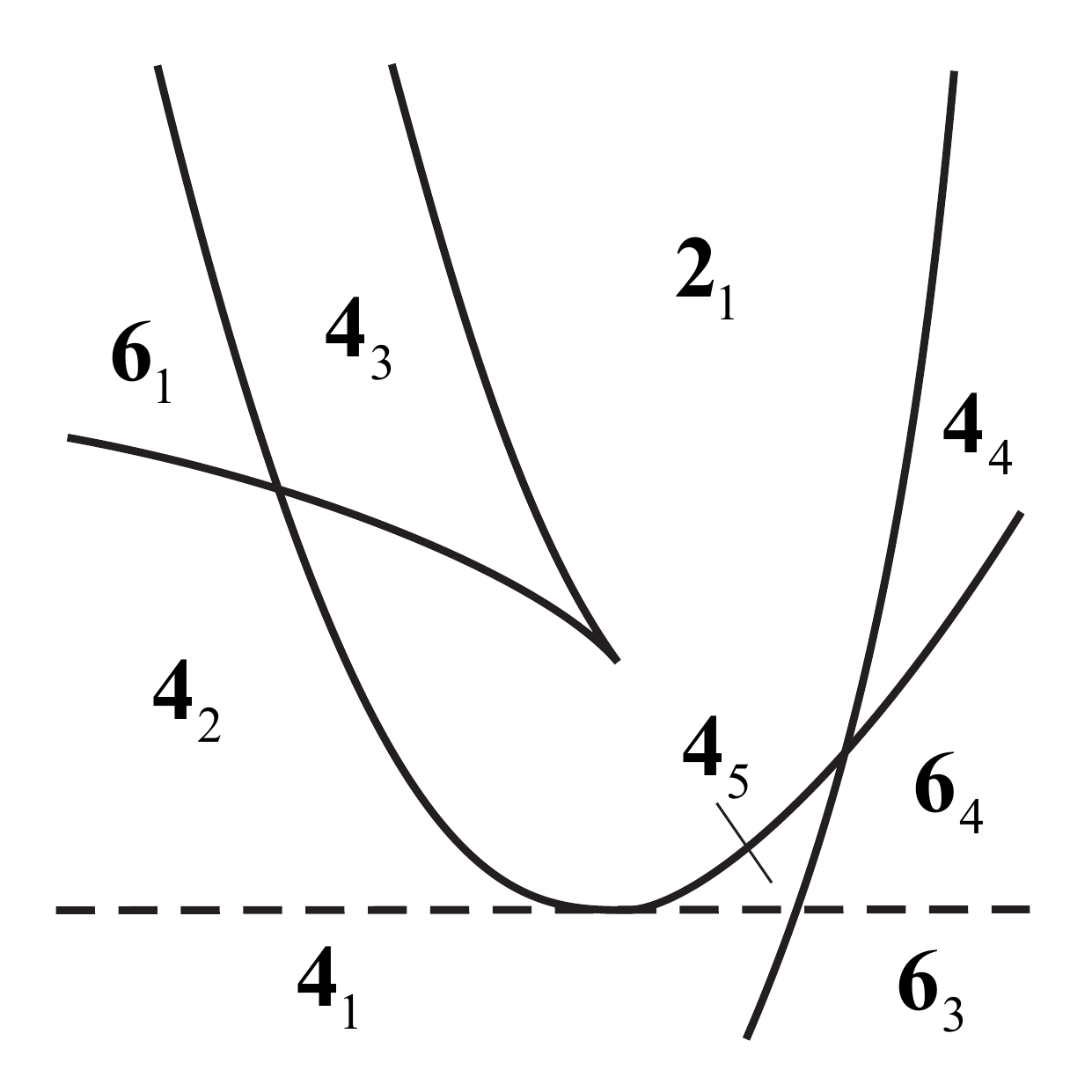}\\
10)\includegraphics[width=4cm]{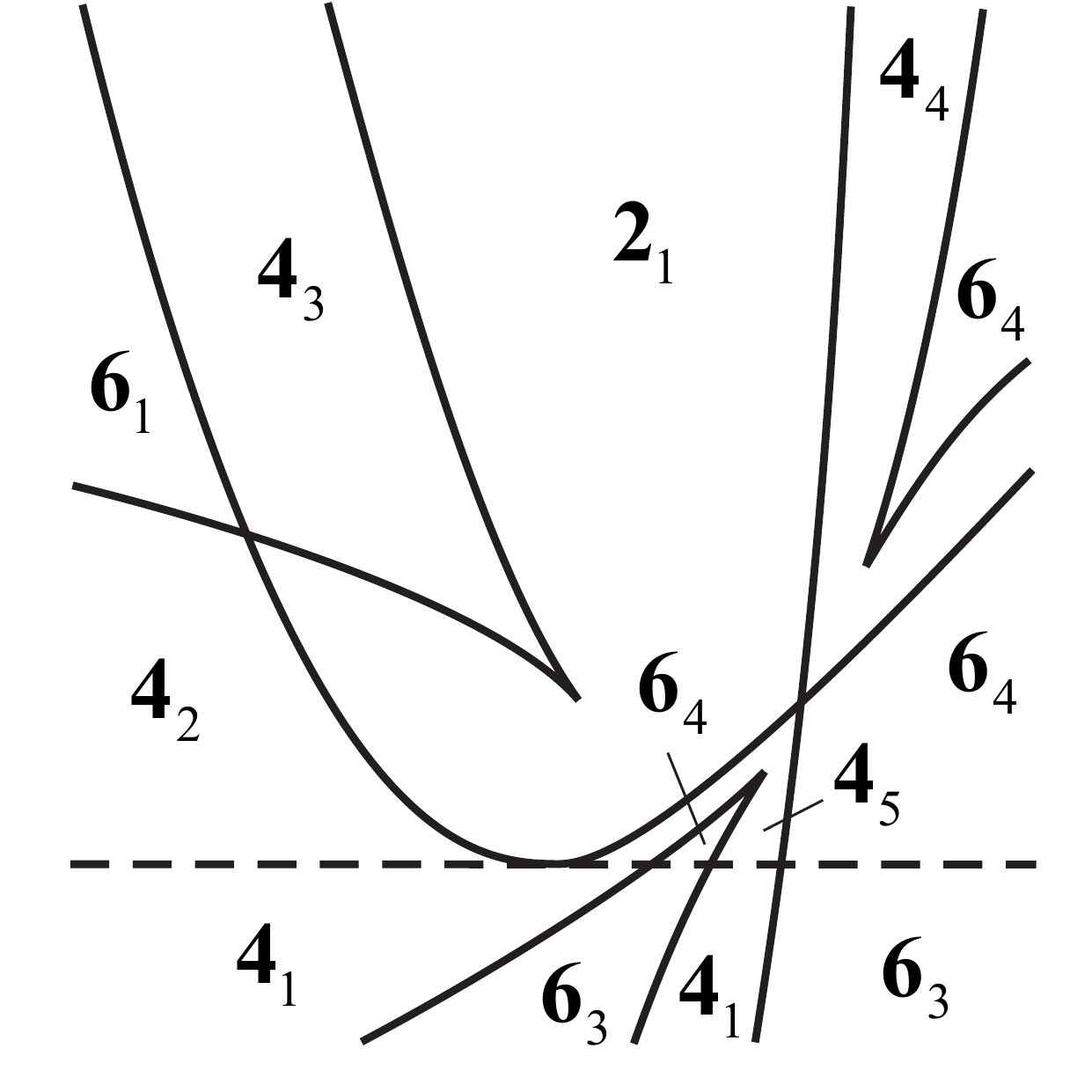}&
11)\includegraphics[width=4cm]{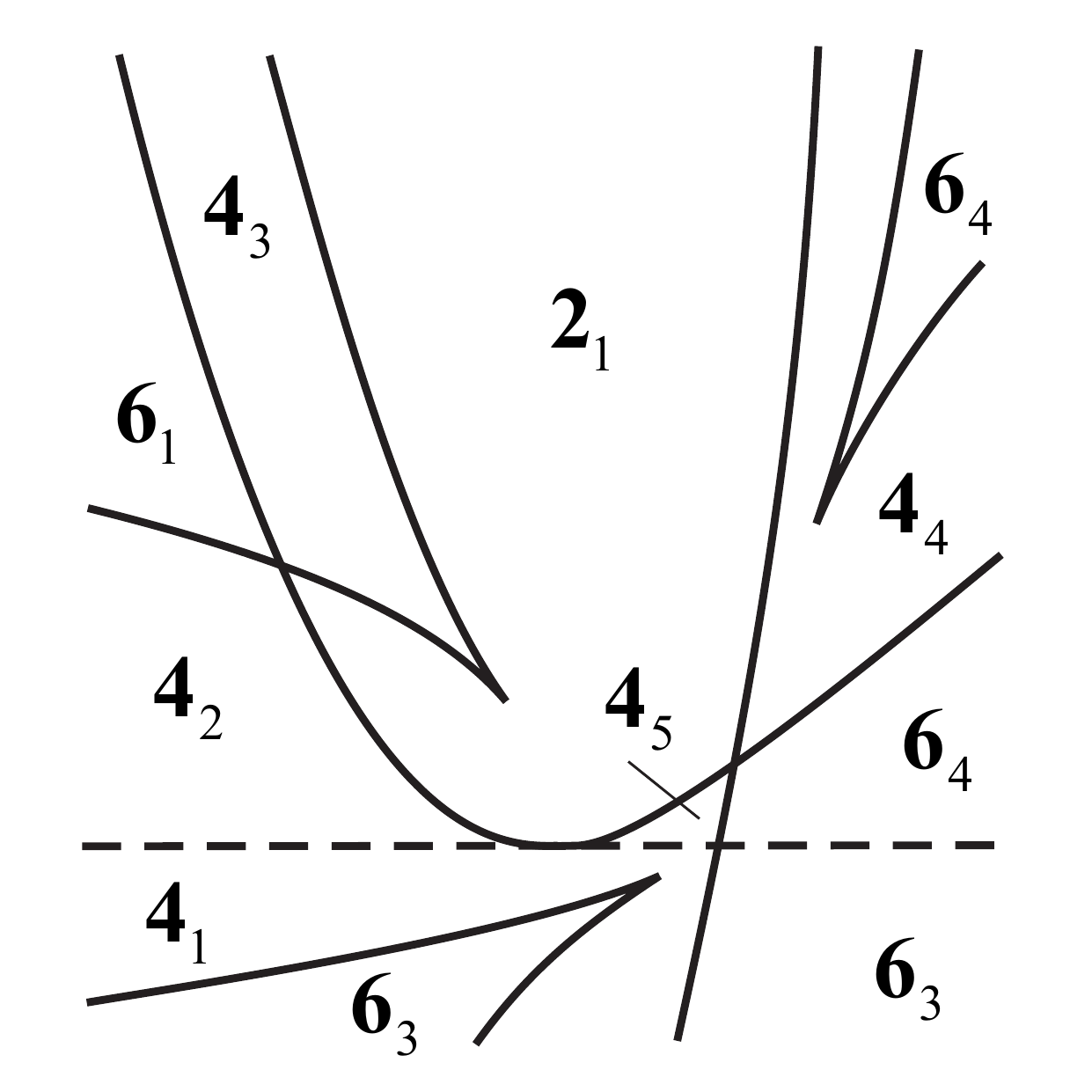}&
12)\includegraphics[width=4cm]{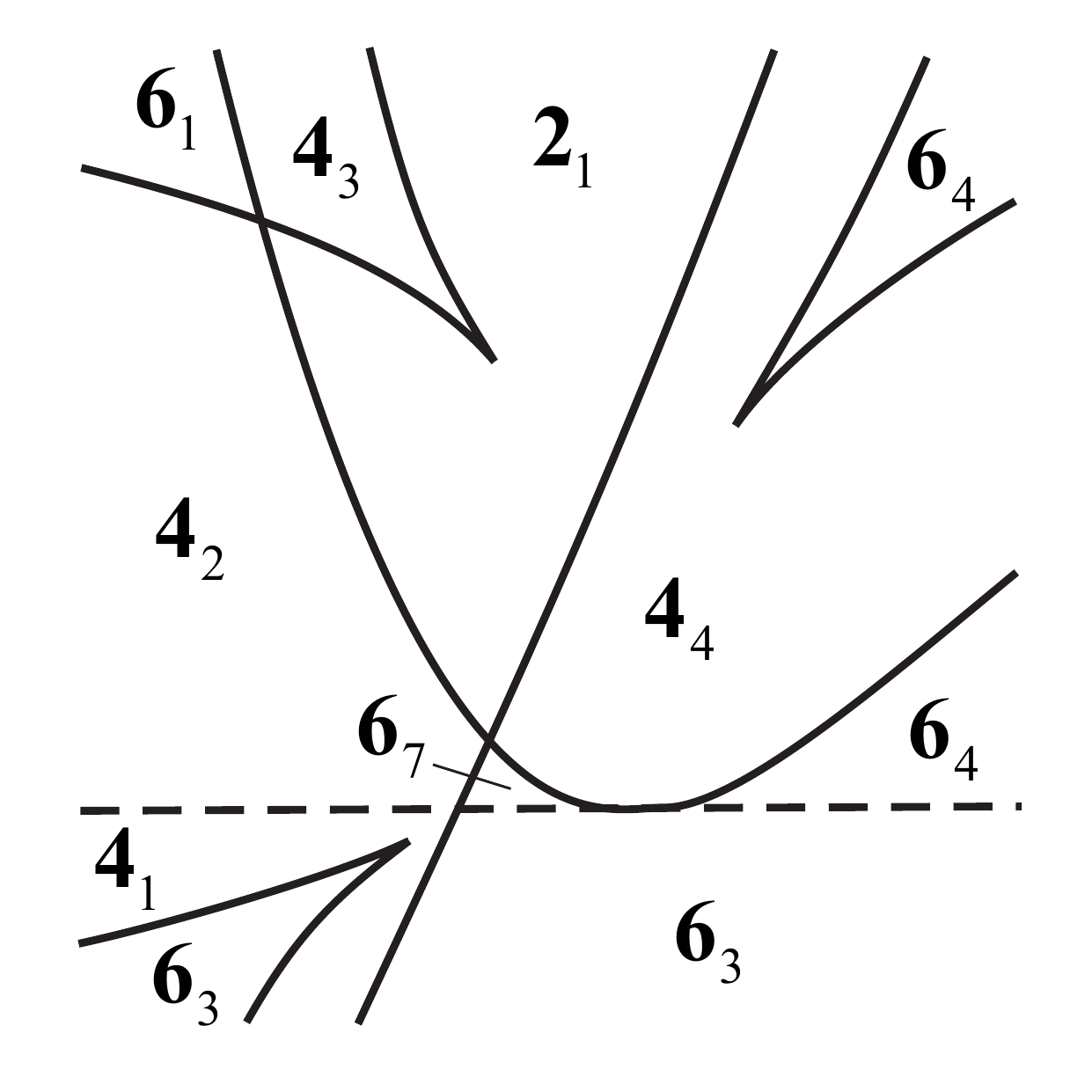}\\
13)\includegraphics[width=4cm]{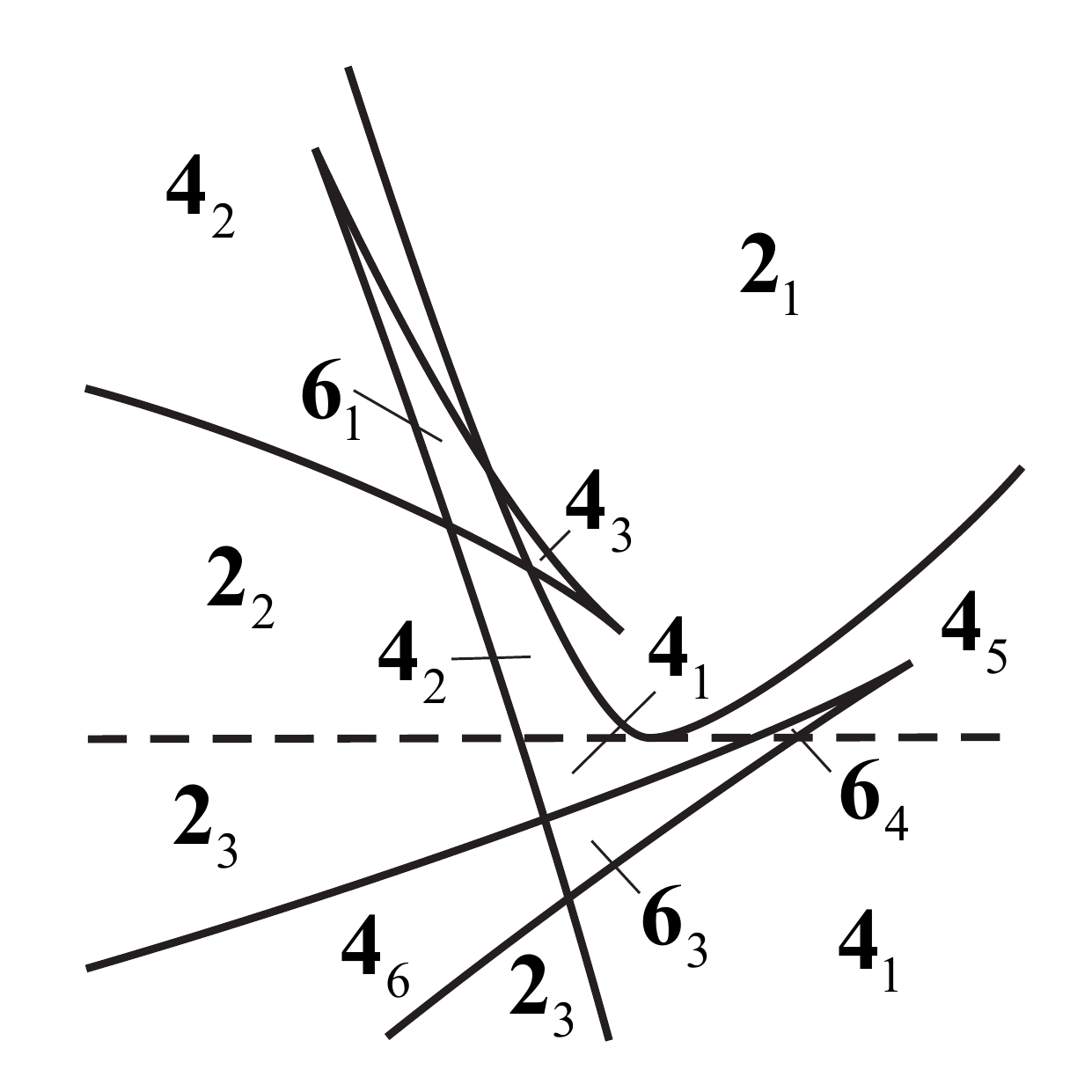}&
14)\includegraphics[width=4cm]{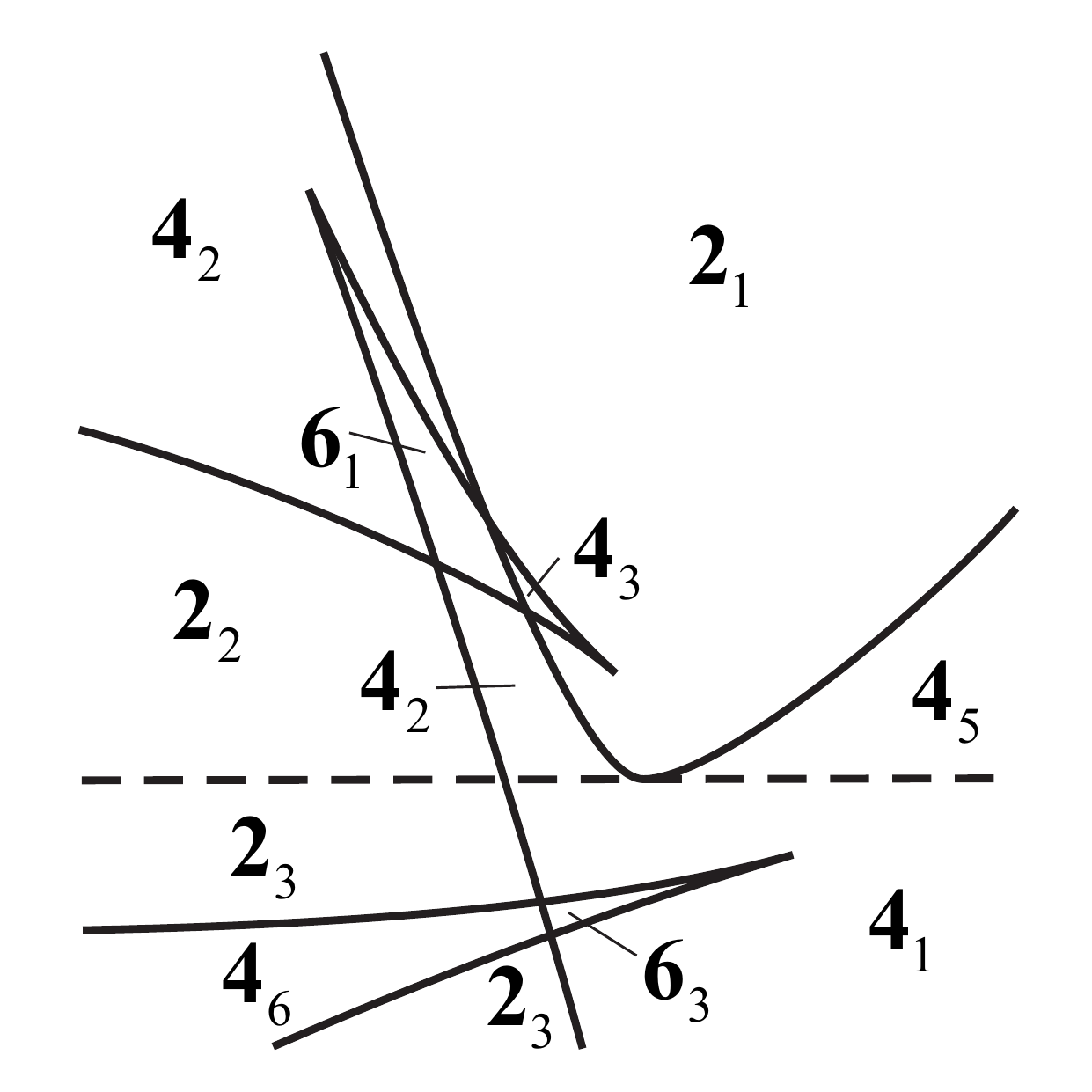}&
15)\includegraphics[width=4cm]{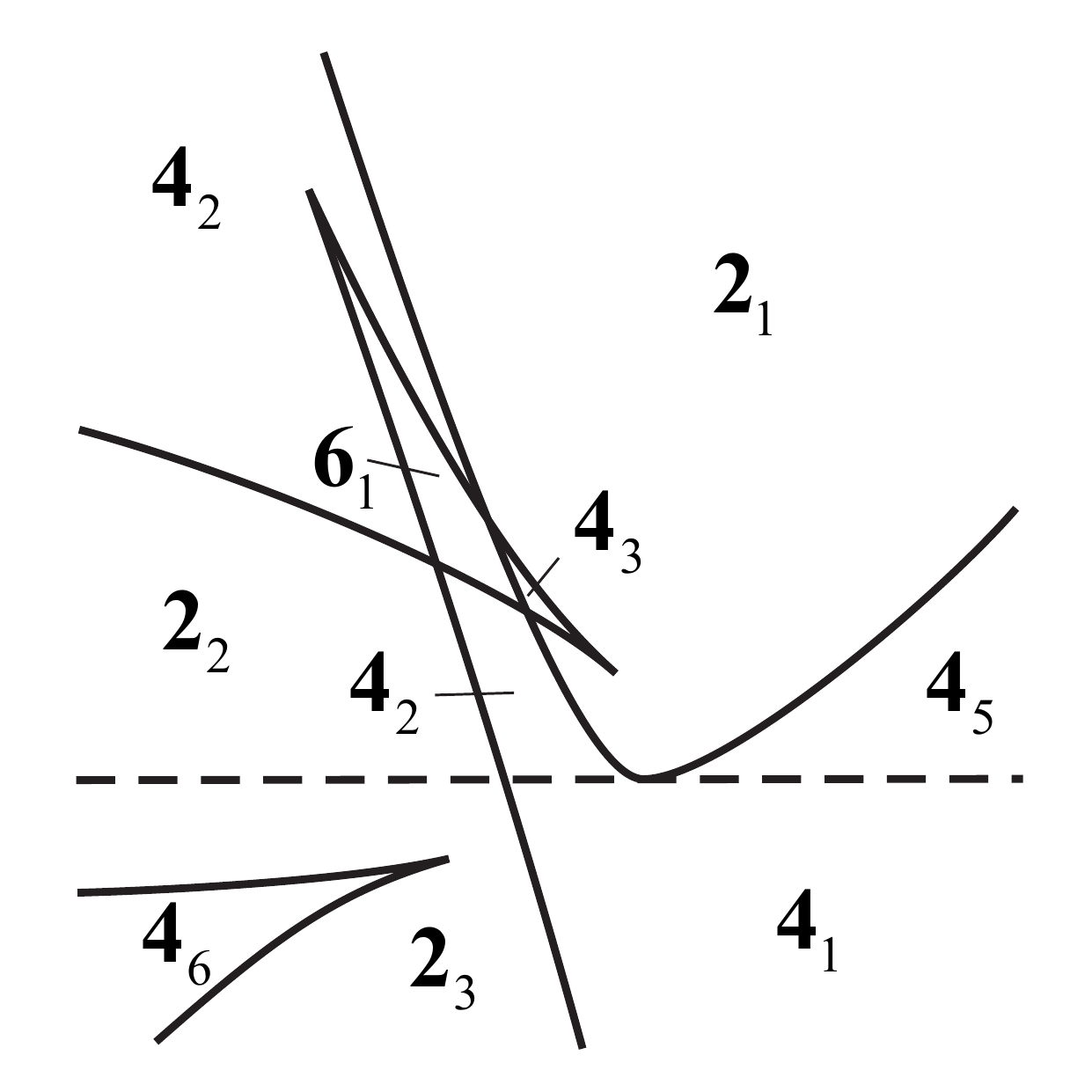}\\
\end{tabular}
\caption{The skeletons $\Sigma_{S_4,t_1,q_5},S_4\neq0$ for domains $1 - 15$.}
\label{zona1-15}
\end{center}
\end{figure}

\begin{figure}
\begin{center}
\begin{tabular}{ccc}
16)\includegraphics[width=4cm]{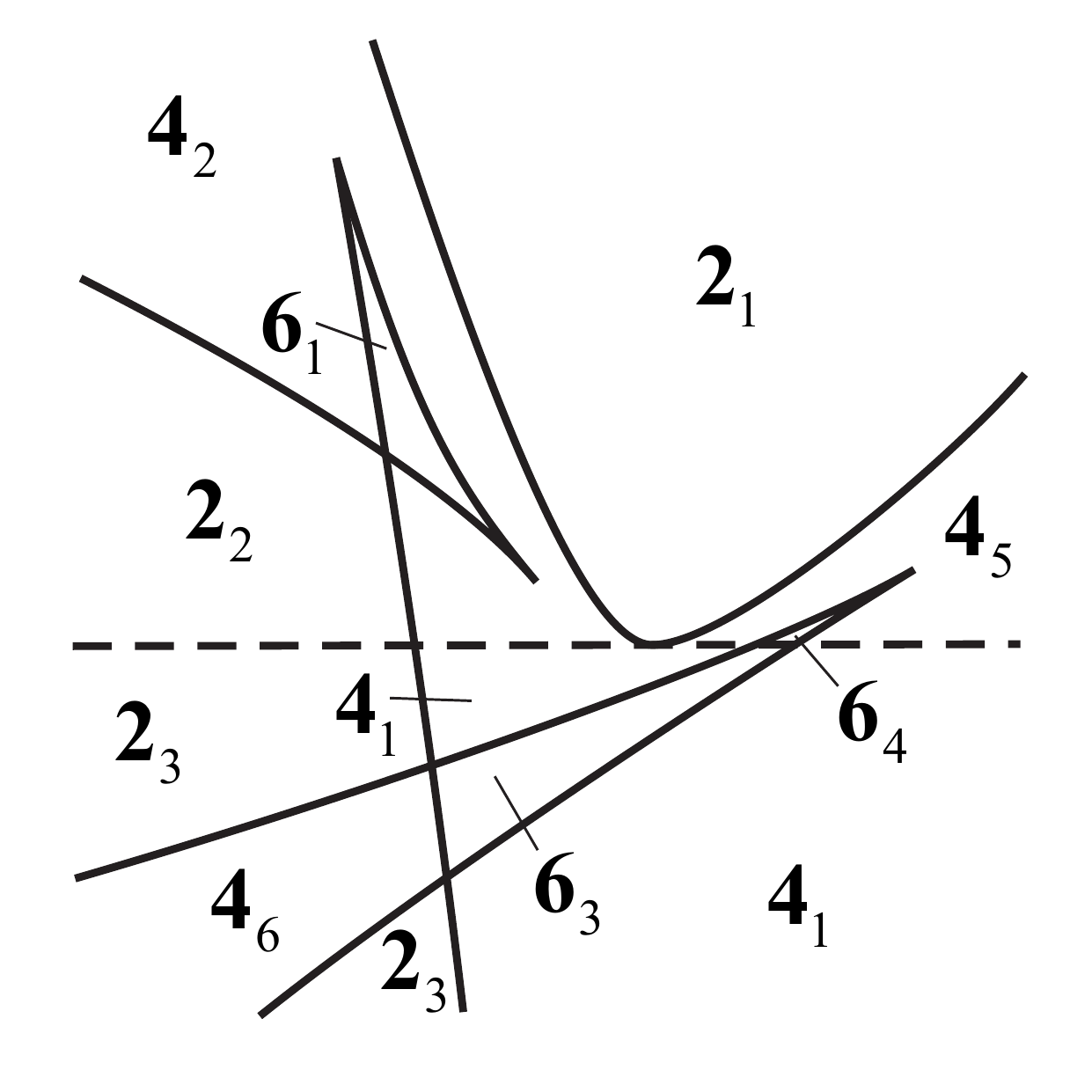}&
17)\includegraphics[width=4cm]{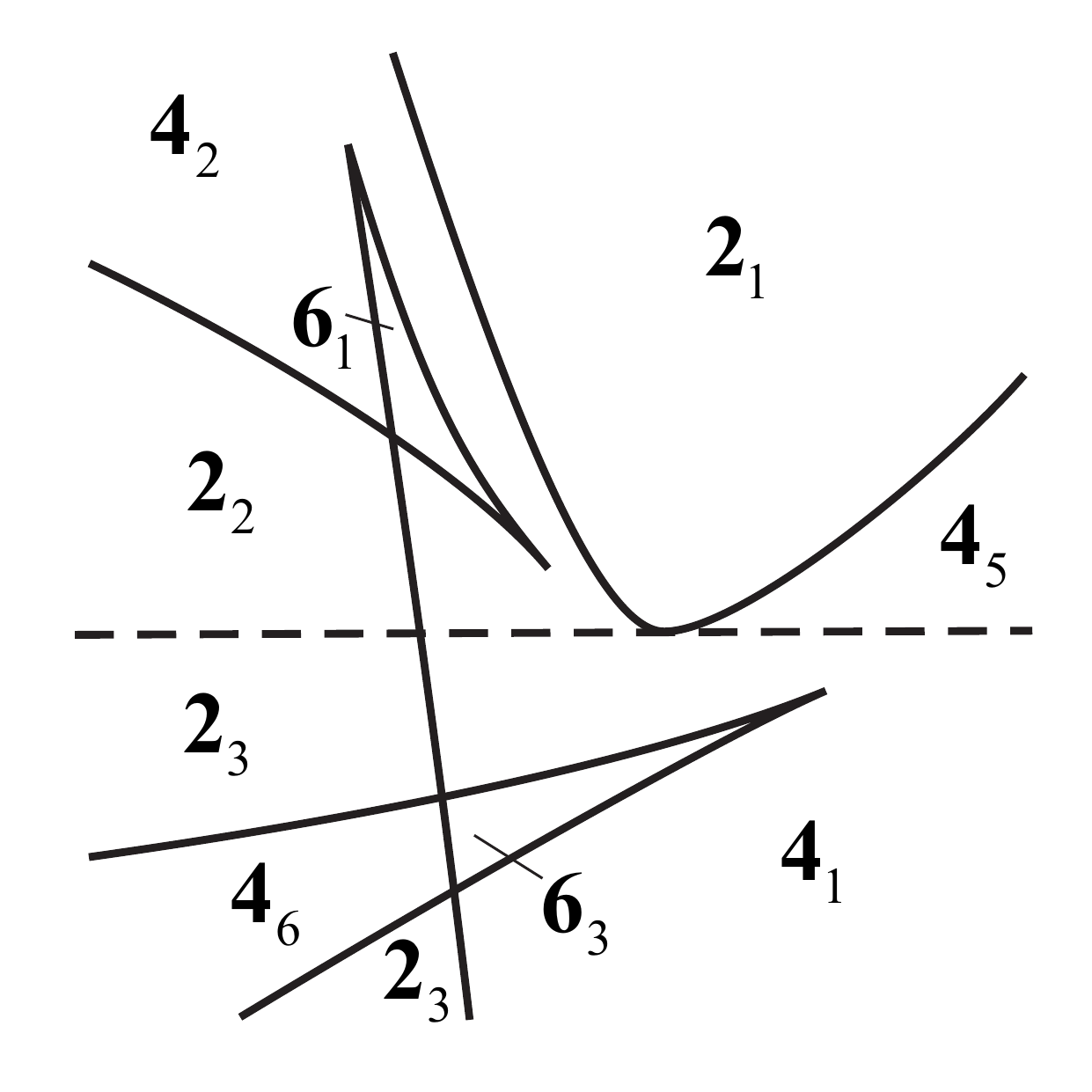}&
18)\includegraphics[width=4cm]{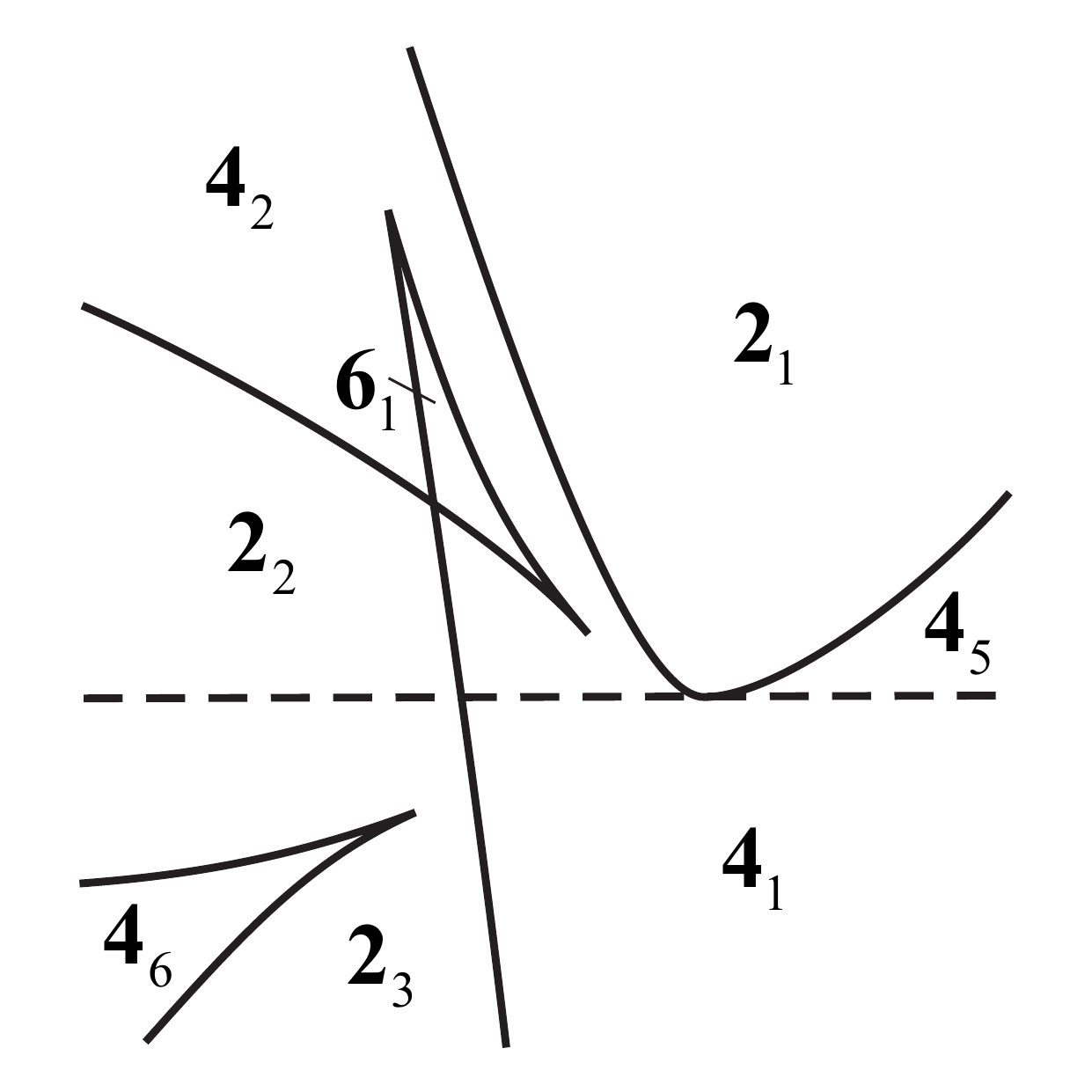}\\
19)\includegraphics[width=4cm]{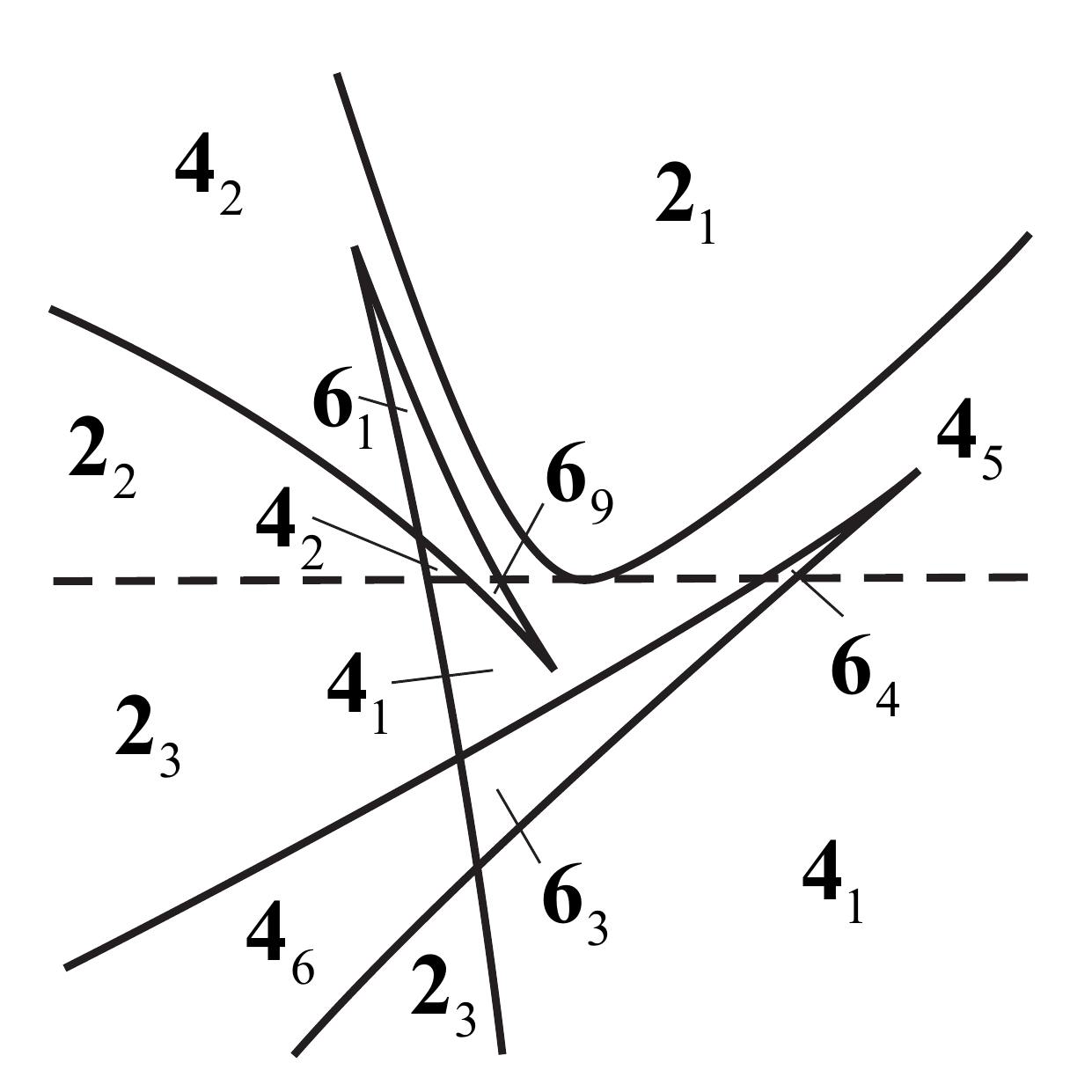}&
20)\includegraphics[width=4cm]{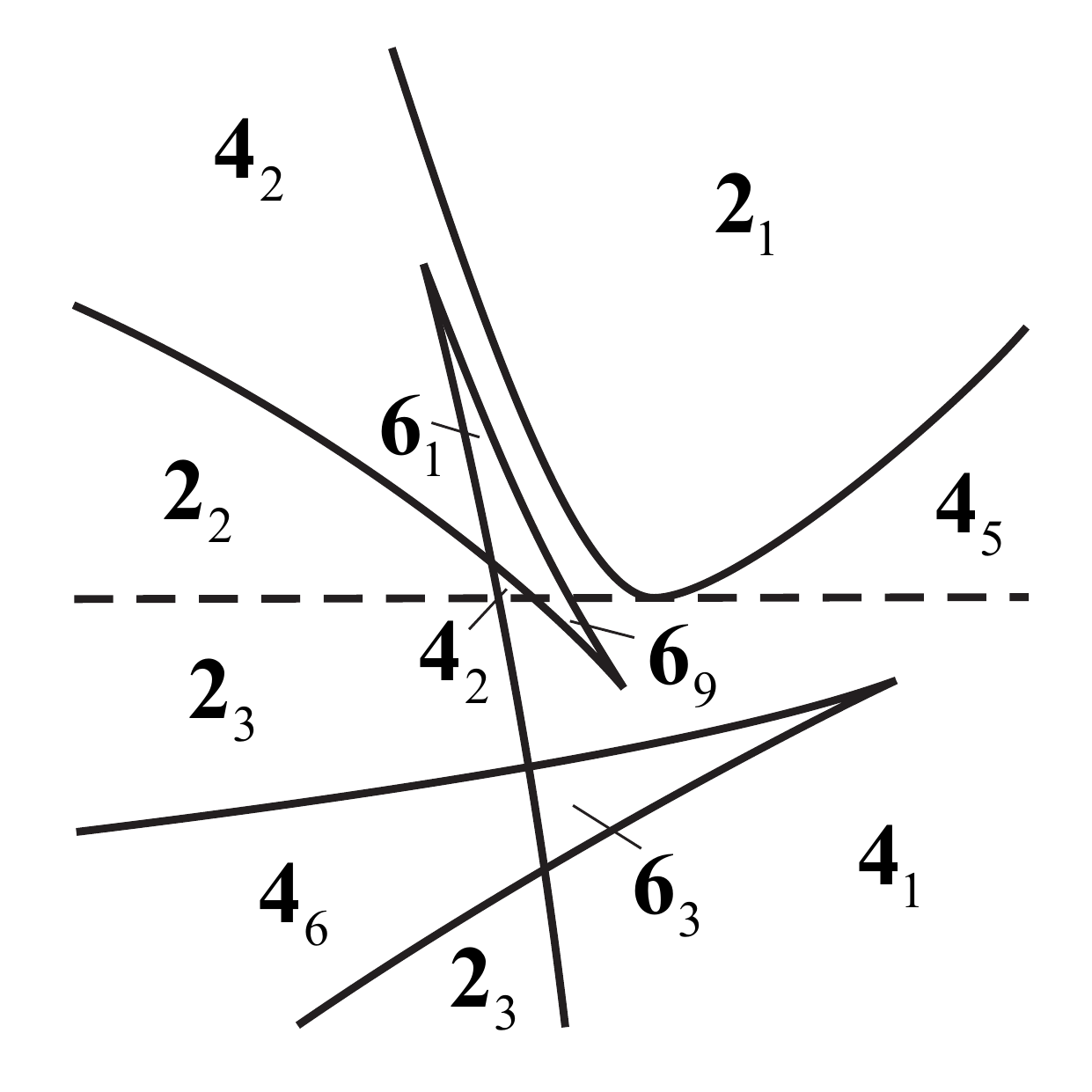}&
21)\includegraphics[width=4cm]{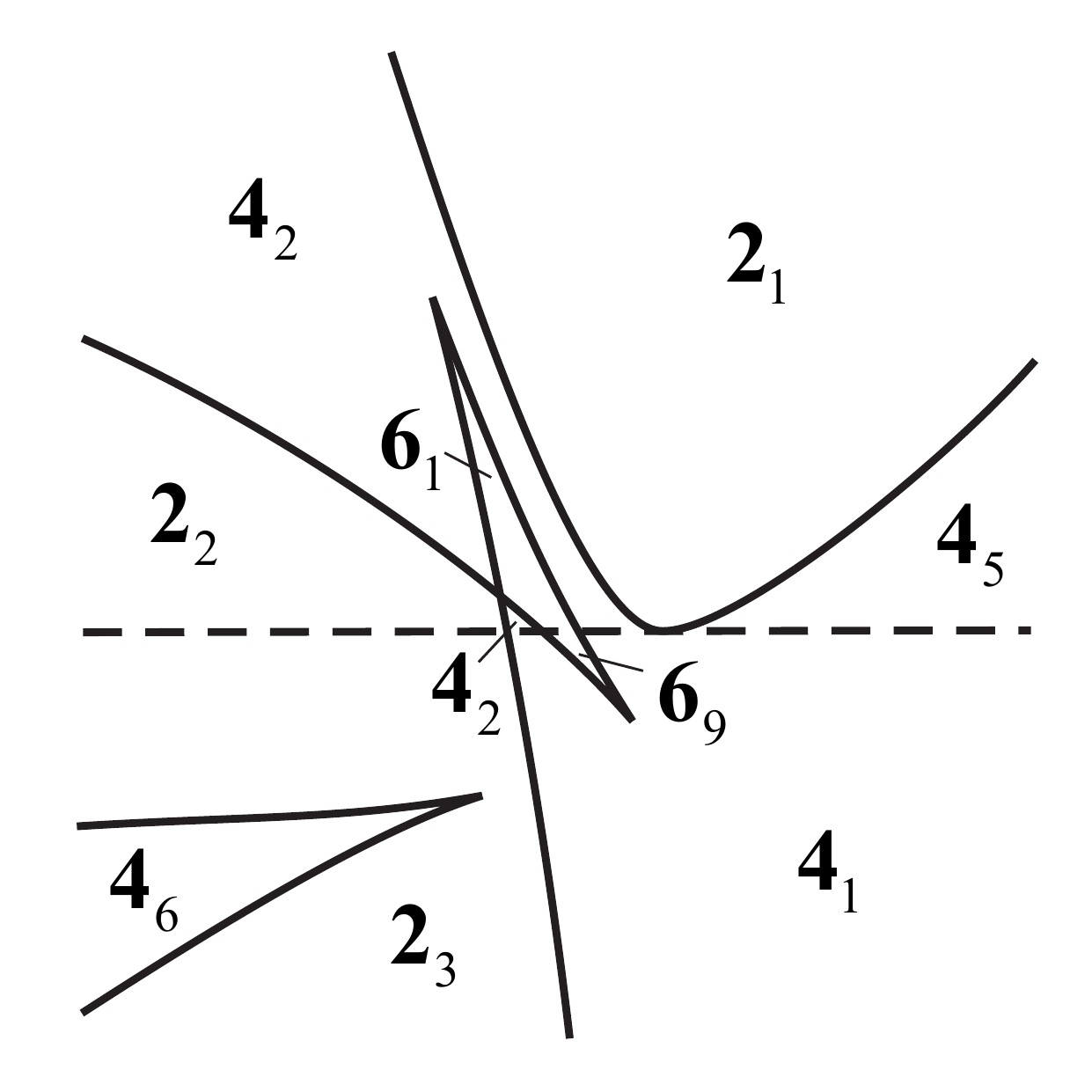}\\
22)\includegraphics[width=4cm]{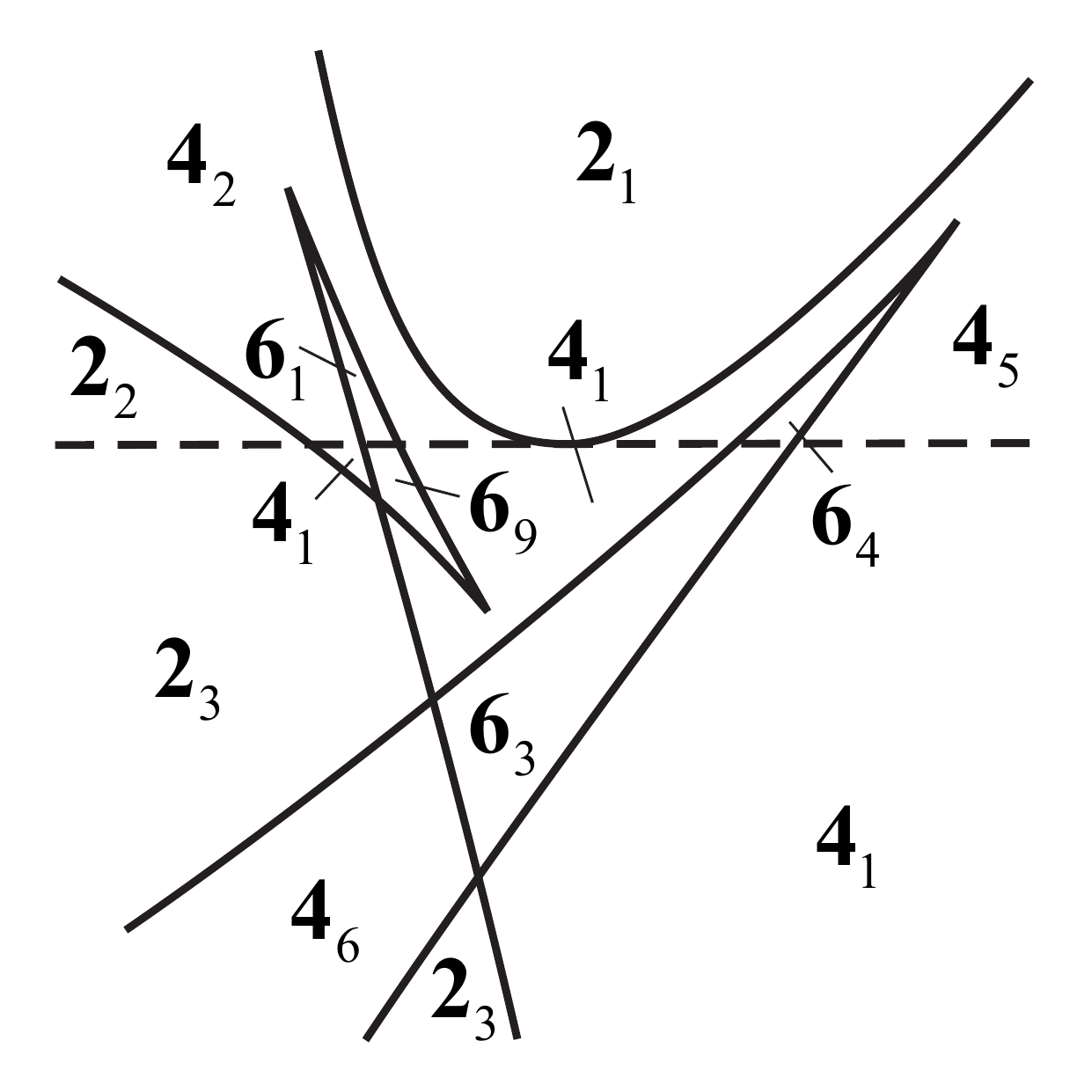}&
23)\includegraphics[width=4cm]{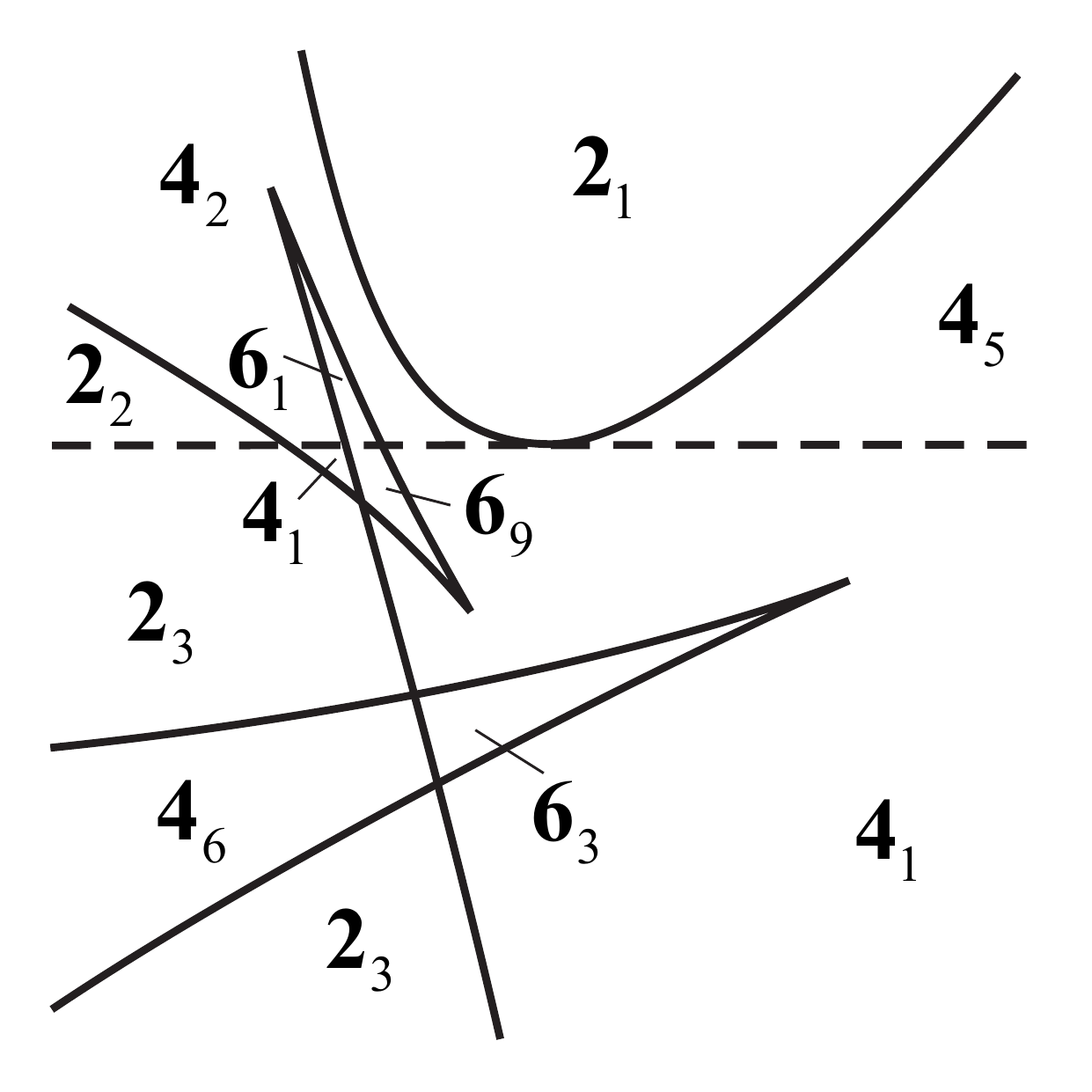}&
24)\includegraphics[width=4cm]{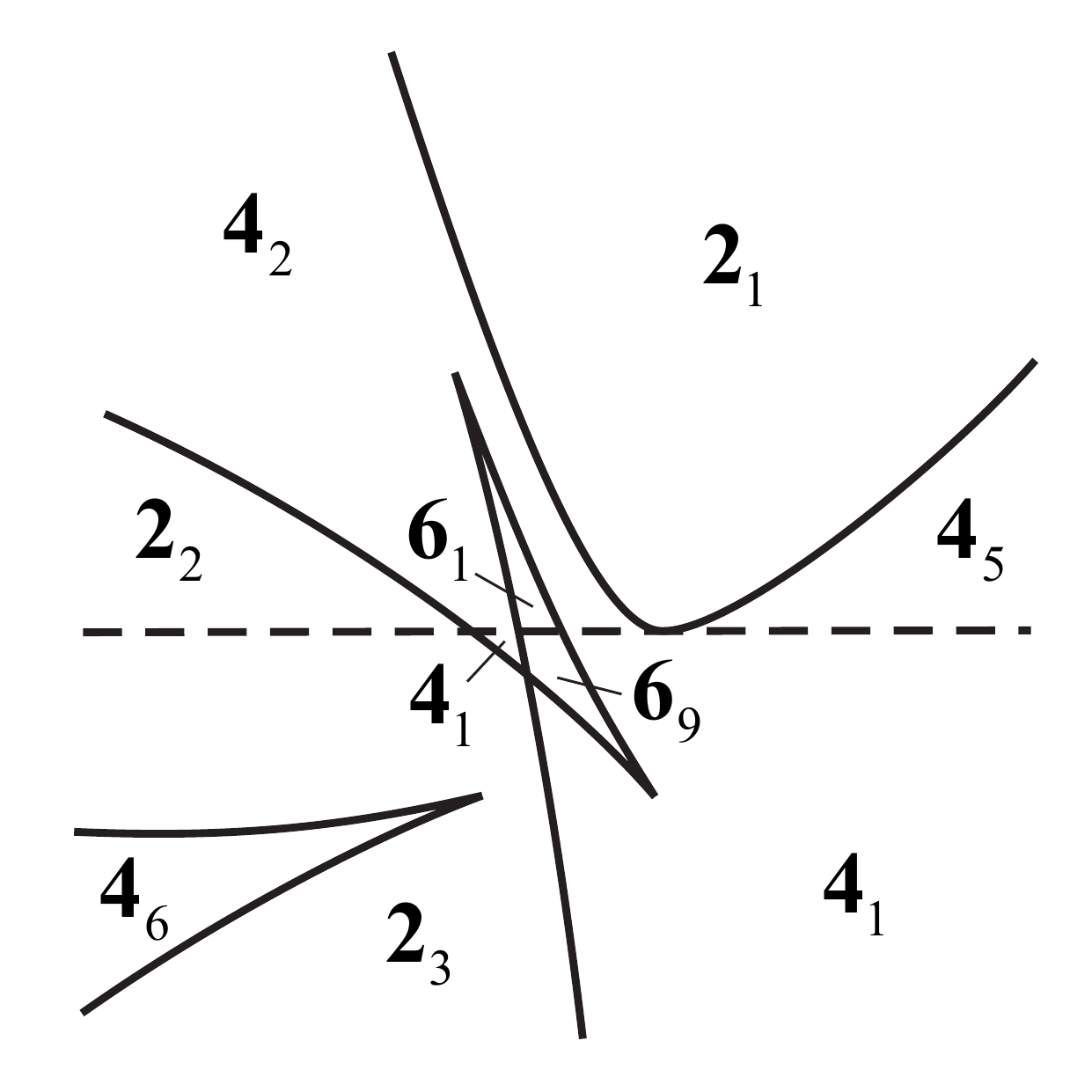}\\
25)\includegraphics[width=4cm]{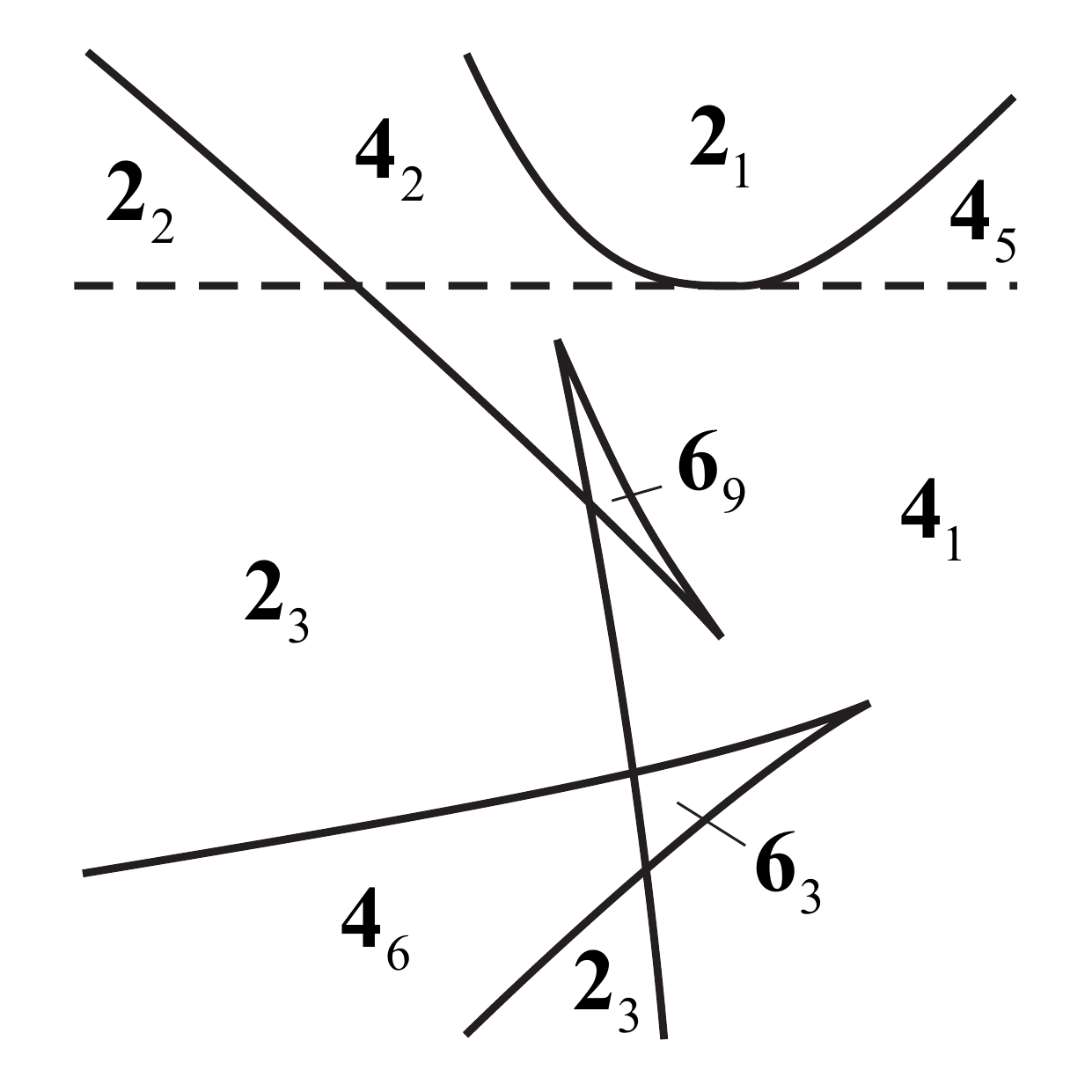}&
26)\includegraphics[width=4cm]{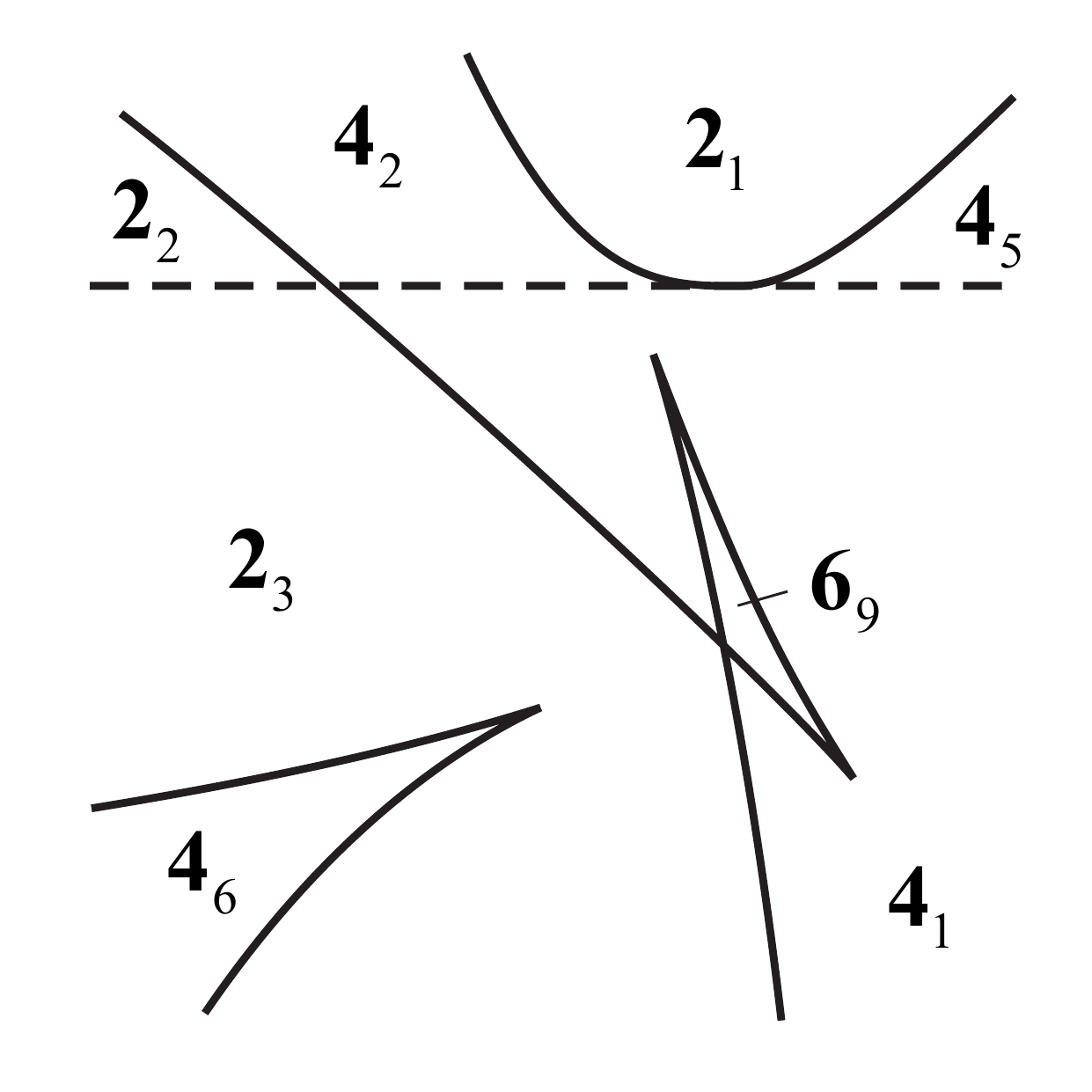}&
27)\includegraphics[width=4cm]{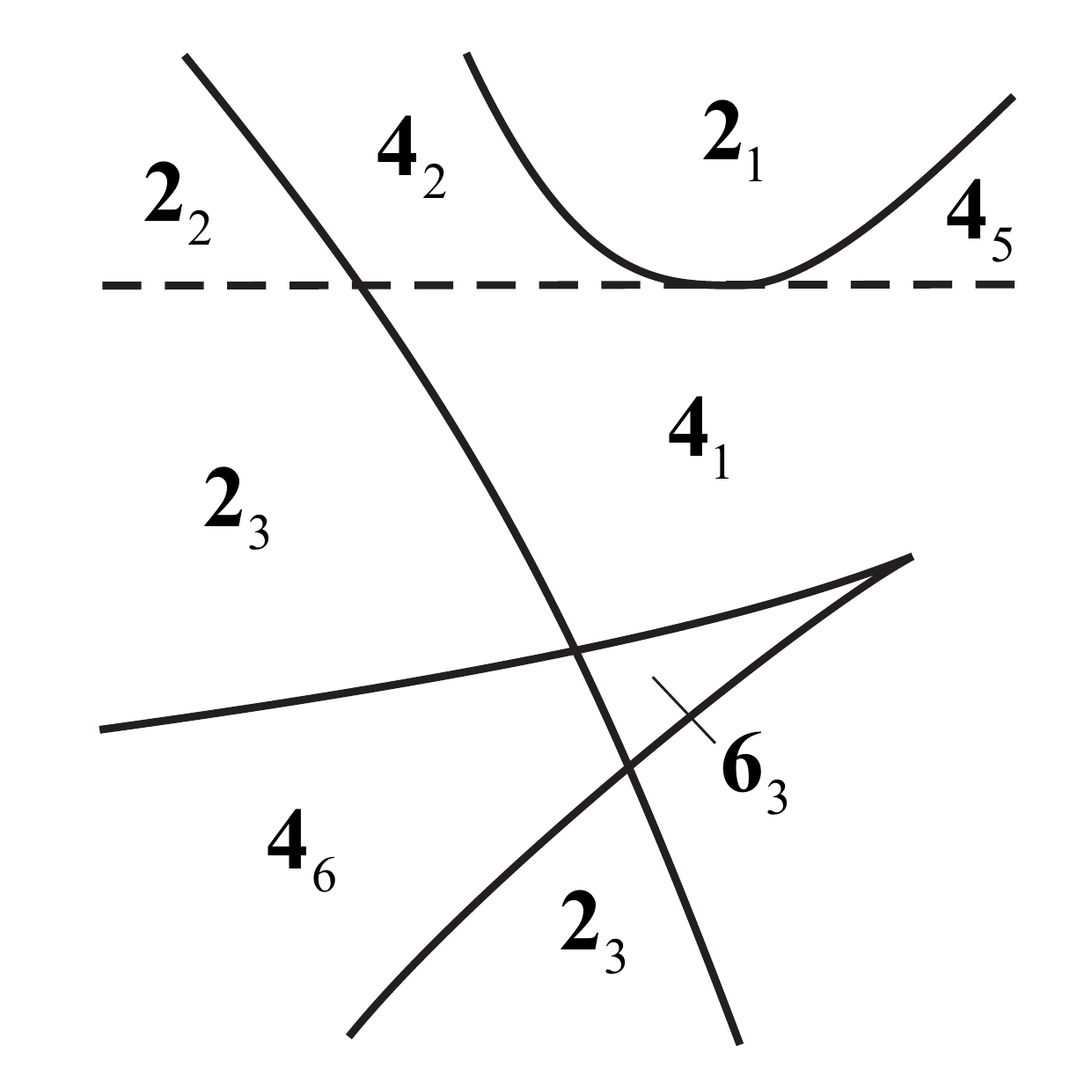}\\
28)\includegraphics[width=4cm]{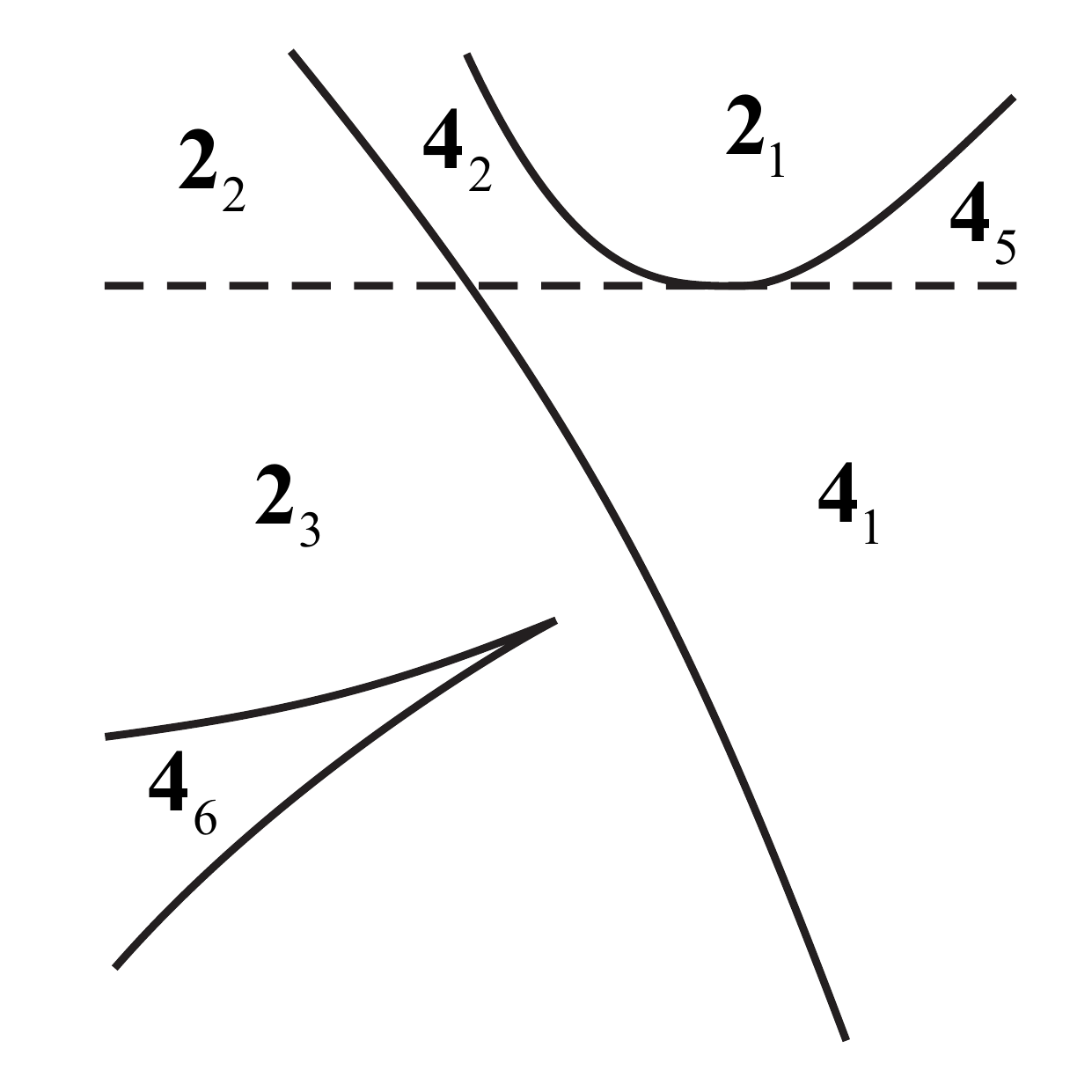}&
29)\includegraphics[width=4cm]{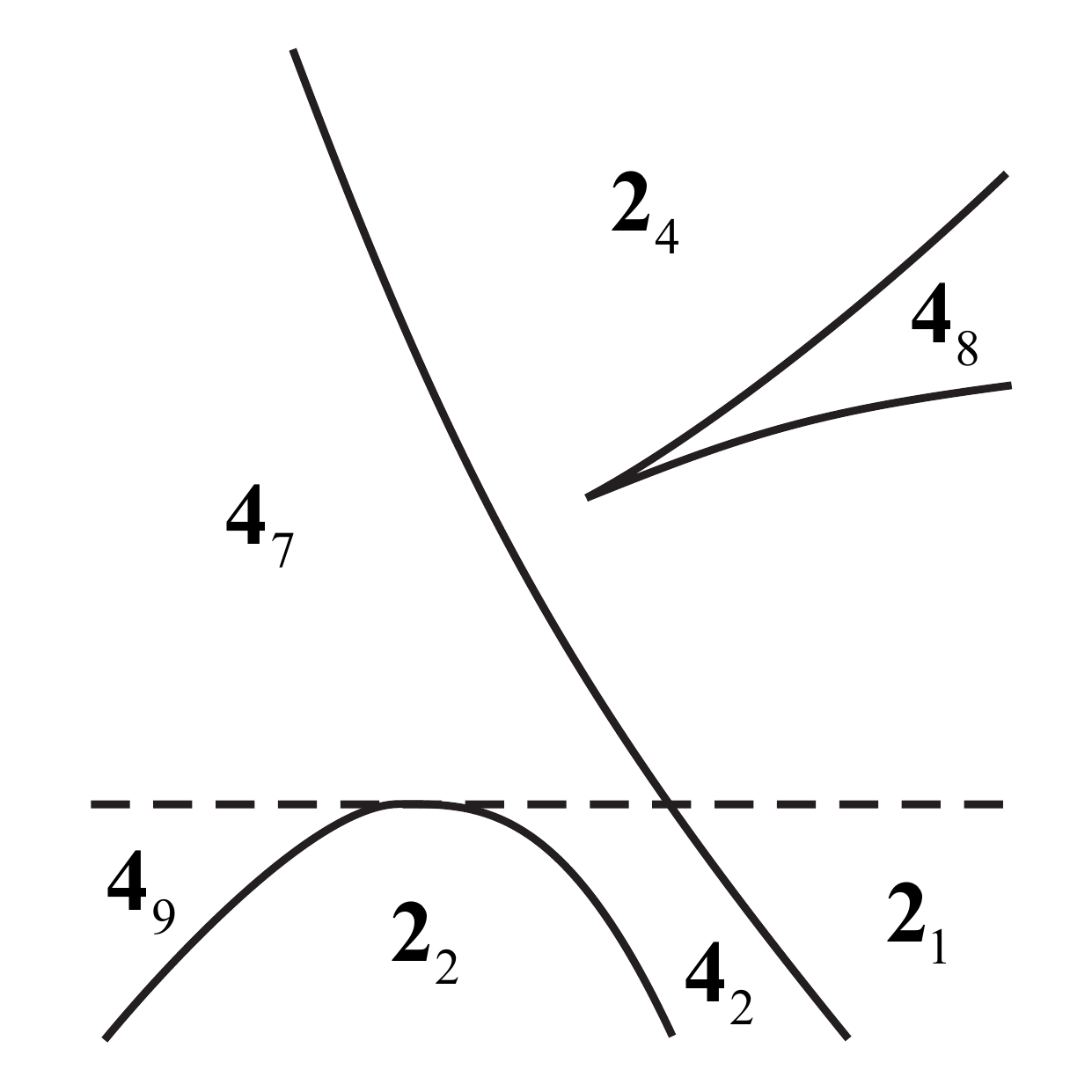}&
30)\includegraphics[width=4cm]{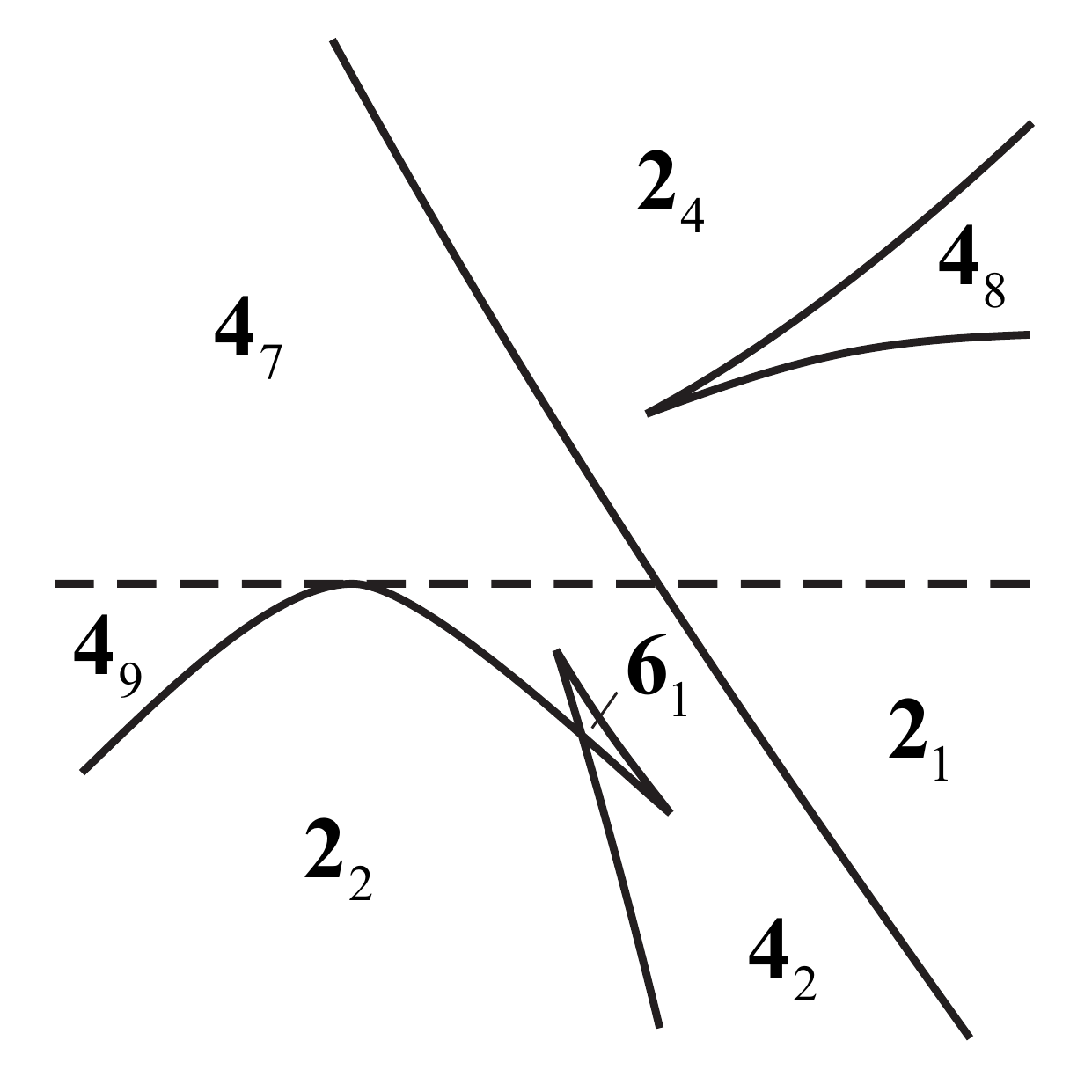}\\
\end{tabular}
\caption{The skeletons $\Sigma_{S_4,t_1,q_5},S_4\neq0$ for domains $16 - 30$.}
\label{zona16-30}
\end{center}
\end{figure}

\begin{figure}
\begin{center}
\begin{tabular}{ccc}
31)\includegraphics[width=4cm]{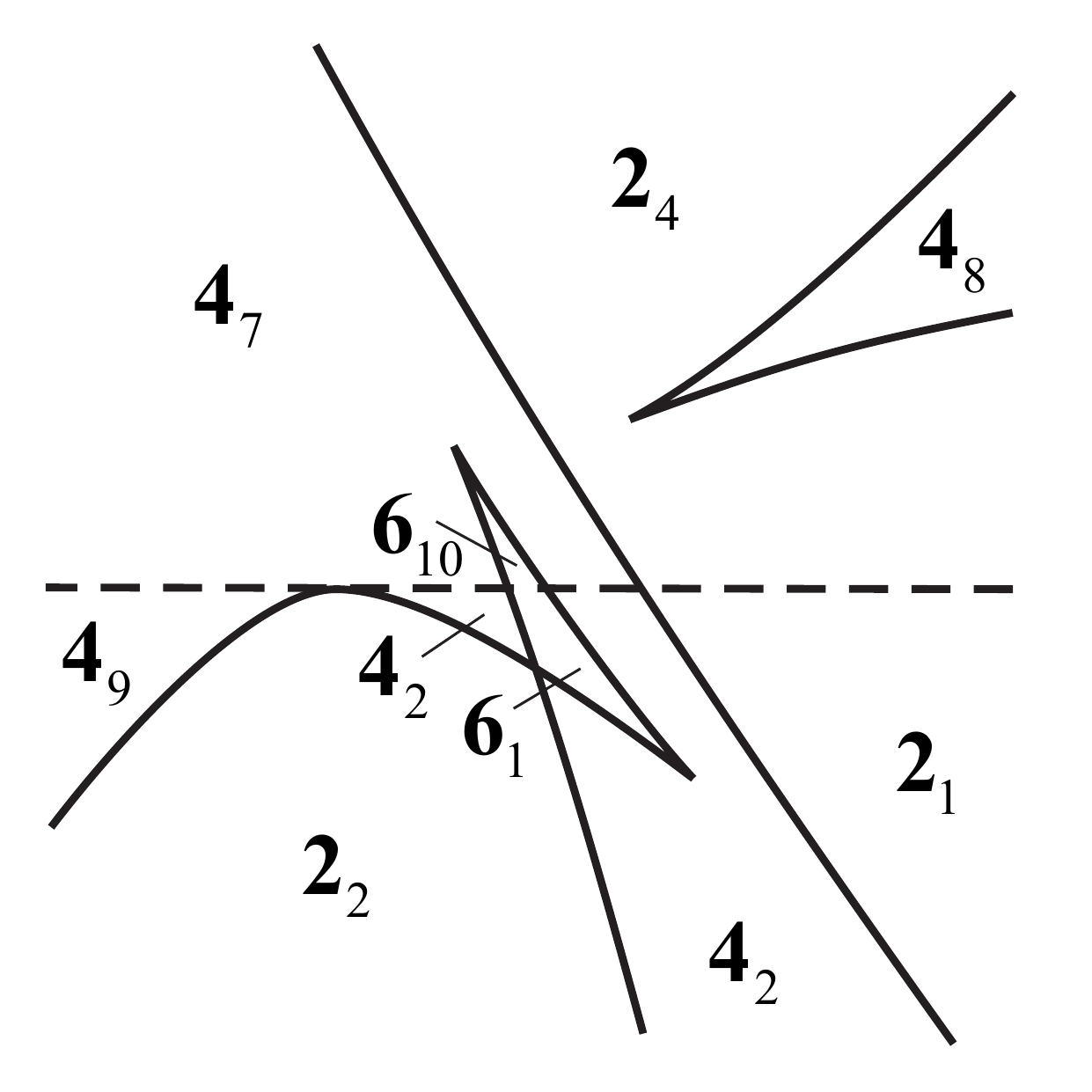}&
32)\includegraphics[width=4cm]{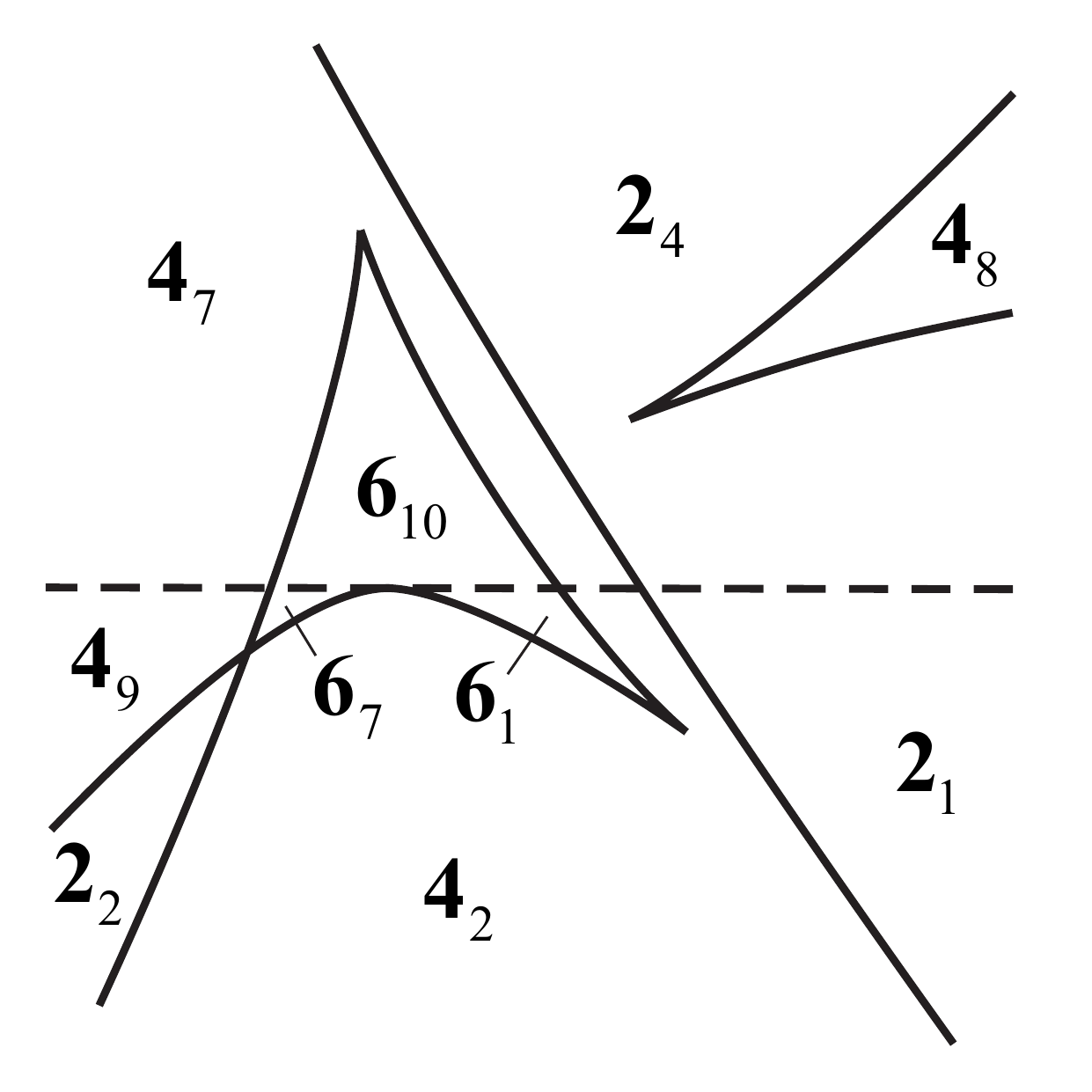}&
33)\includegraphics[width=4cm]{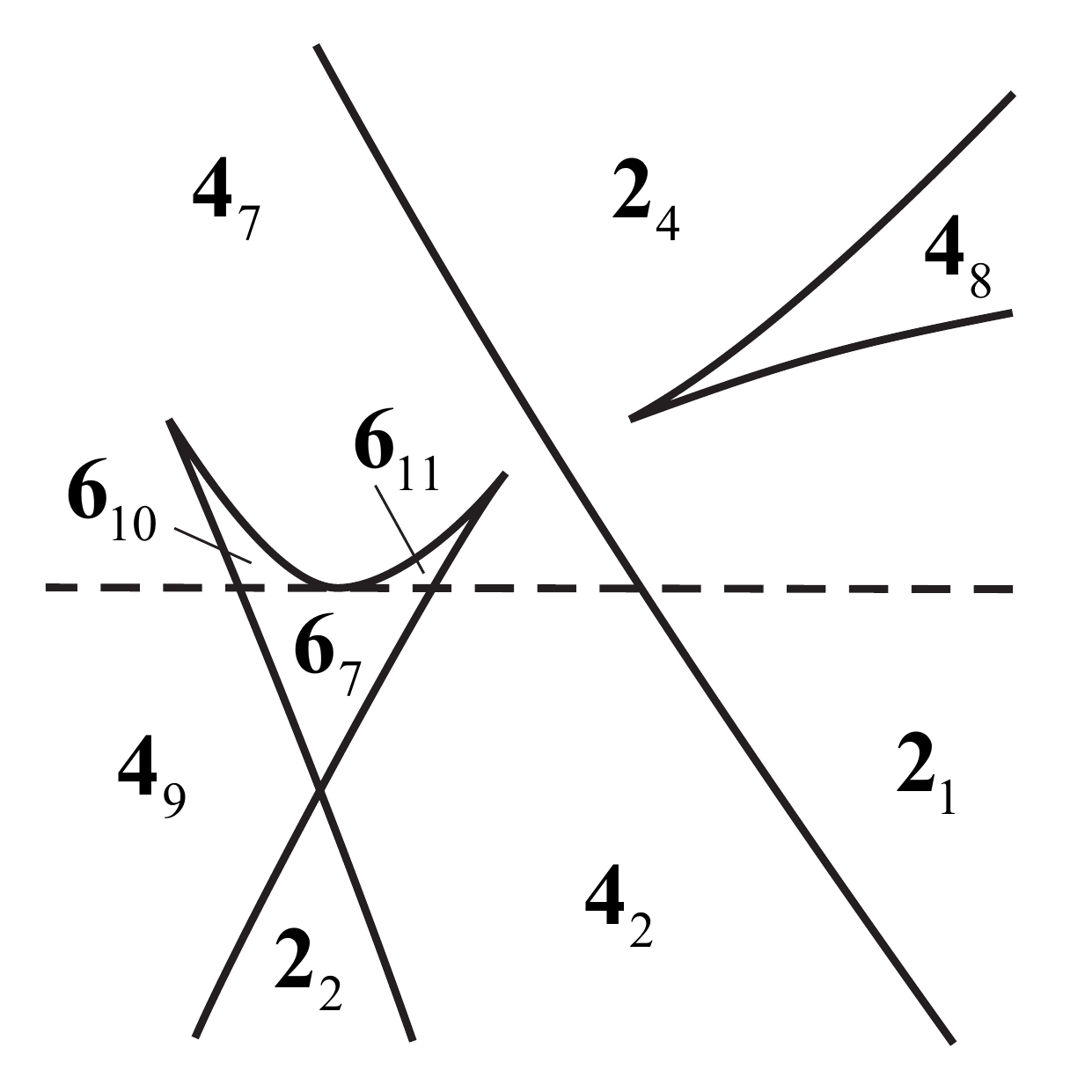}\\
34)\includegraphics[width=4cm]{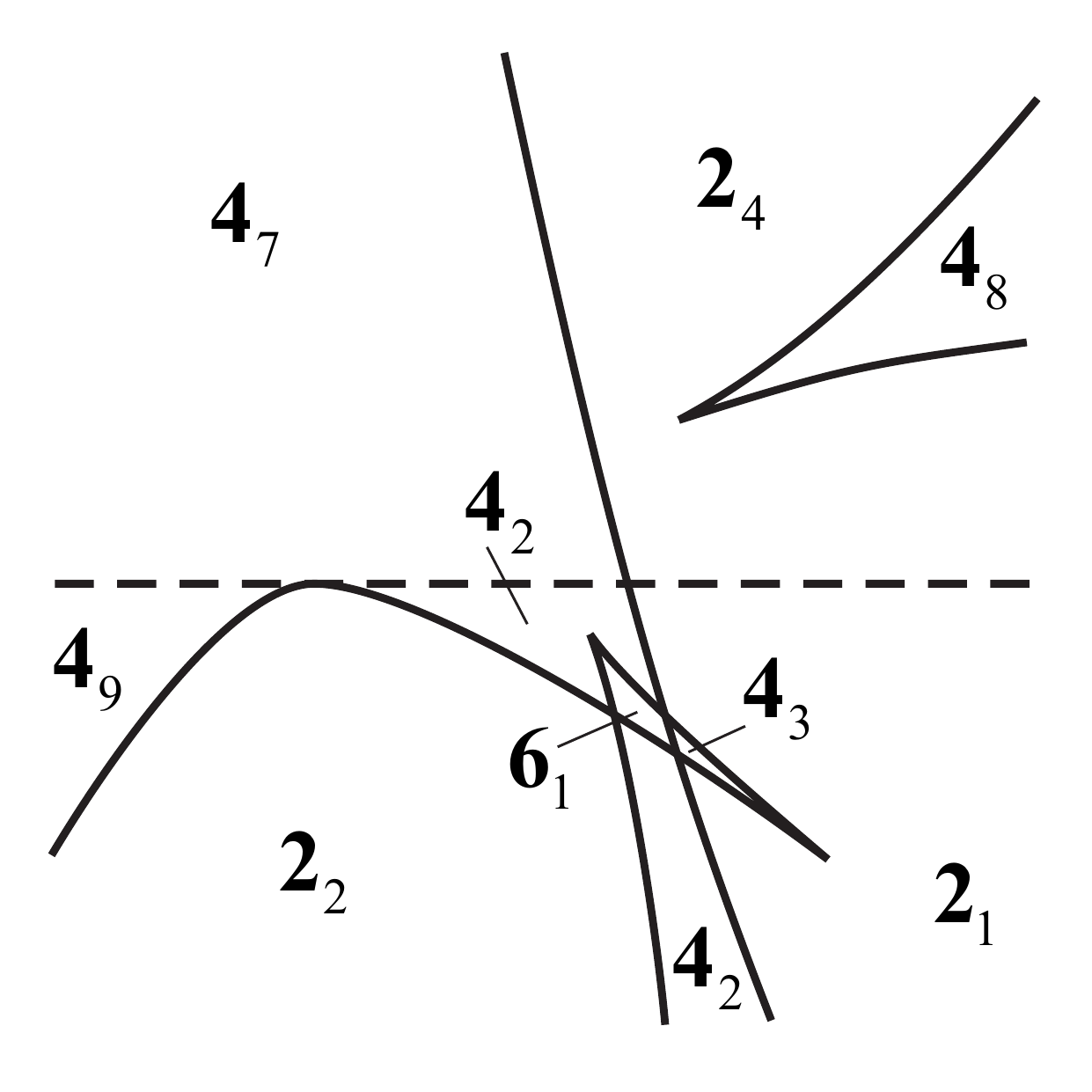}&
35)\includegraphics[width=4cm]{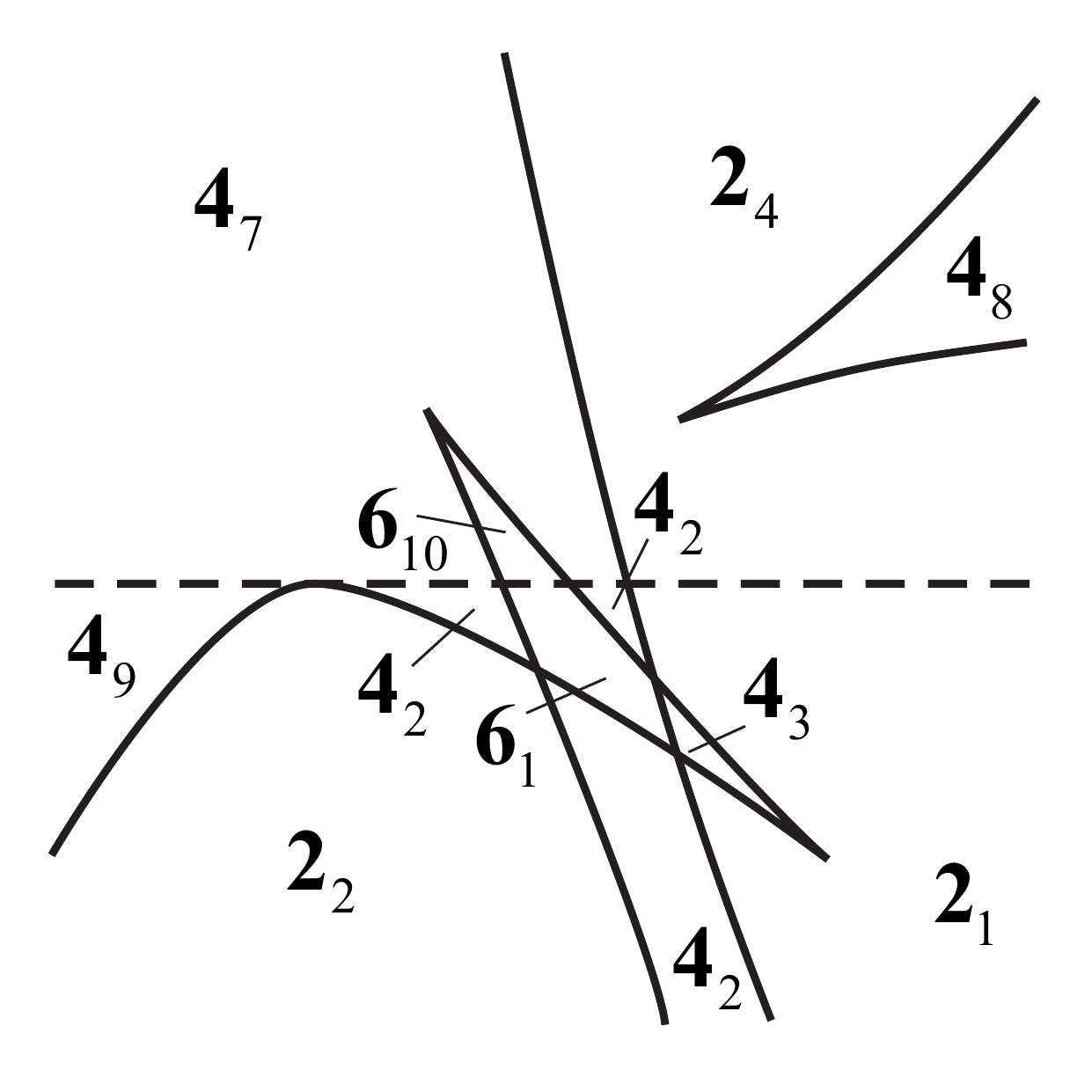}&
36)\includegraphics[width=4cm]{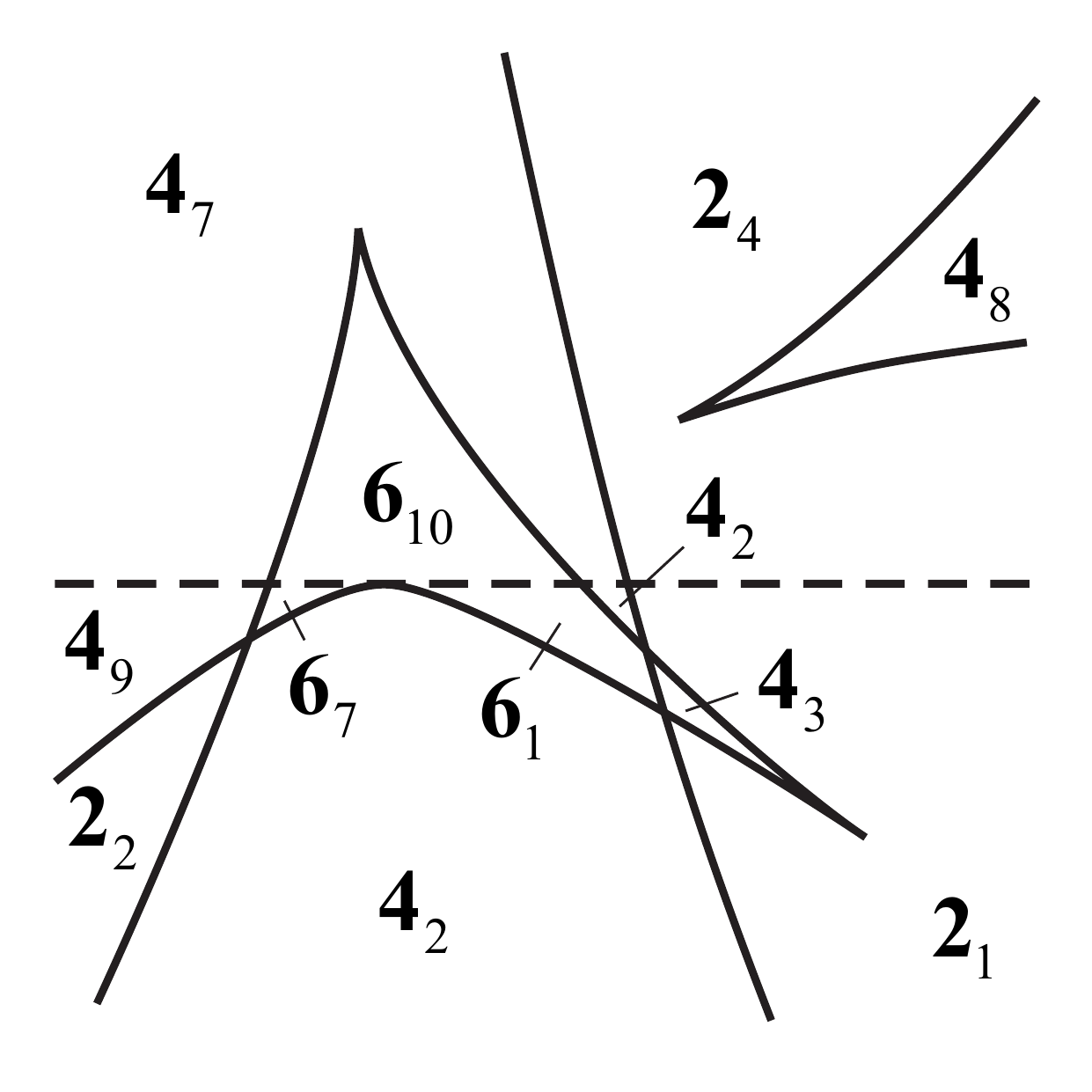}\\
37)\includegraphics[width=4cm]{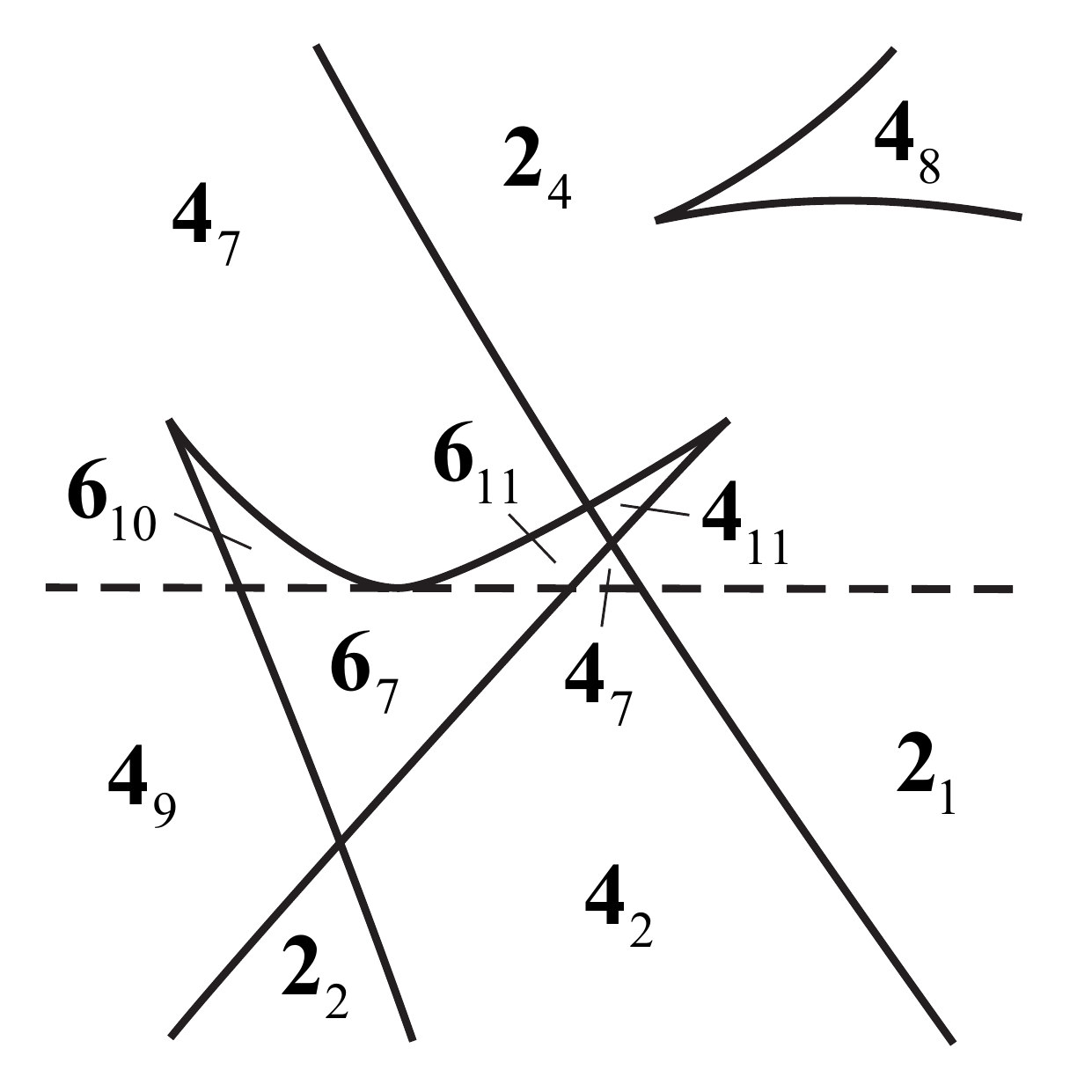}&
38)\includegraphics[width=4cm]{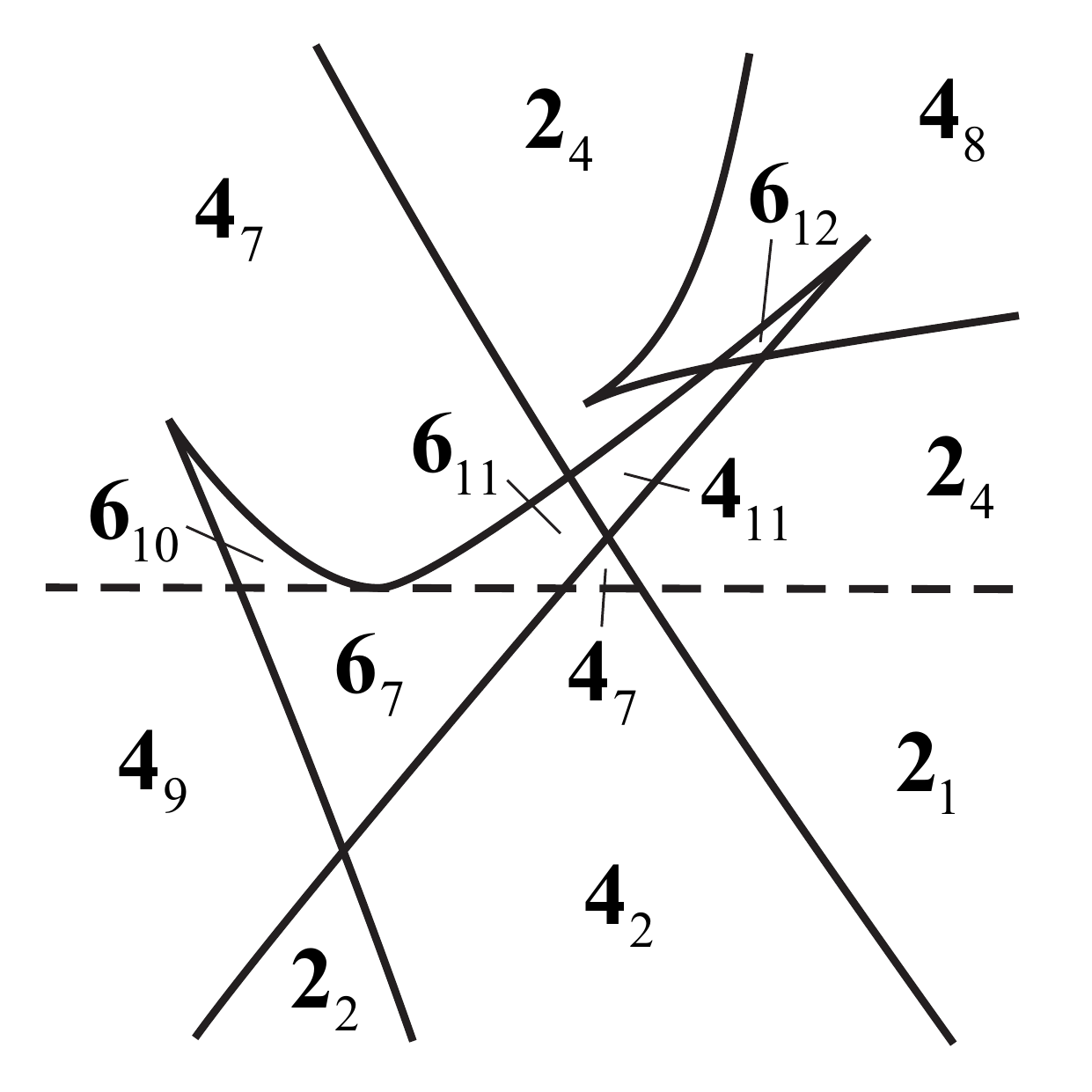}&
39)\includegraphics[width=4cm]{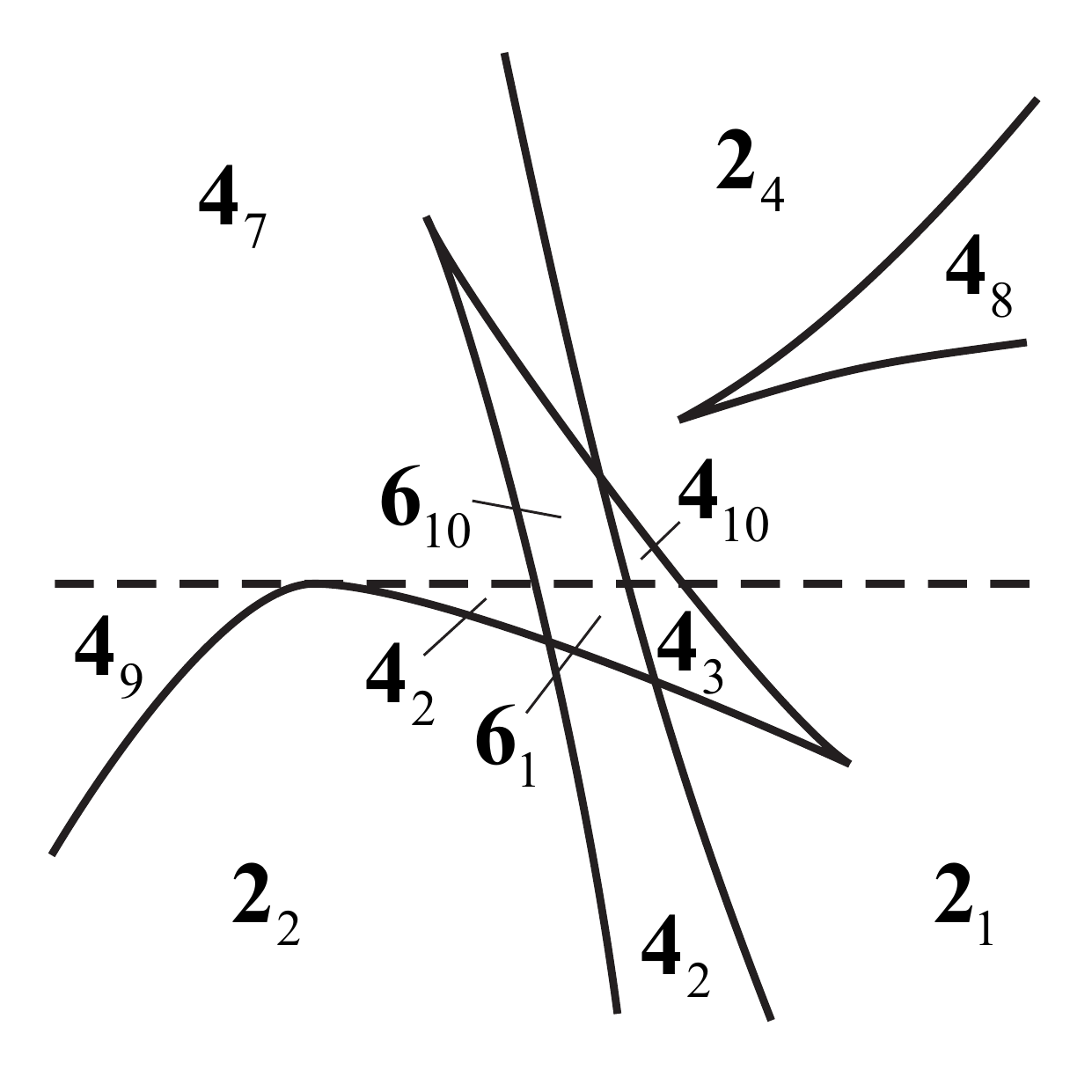}\\
40)\includegraphics[width=4cm]{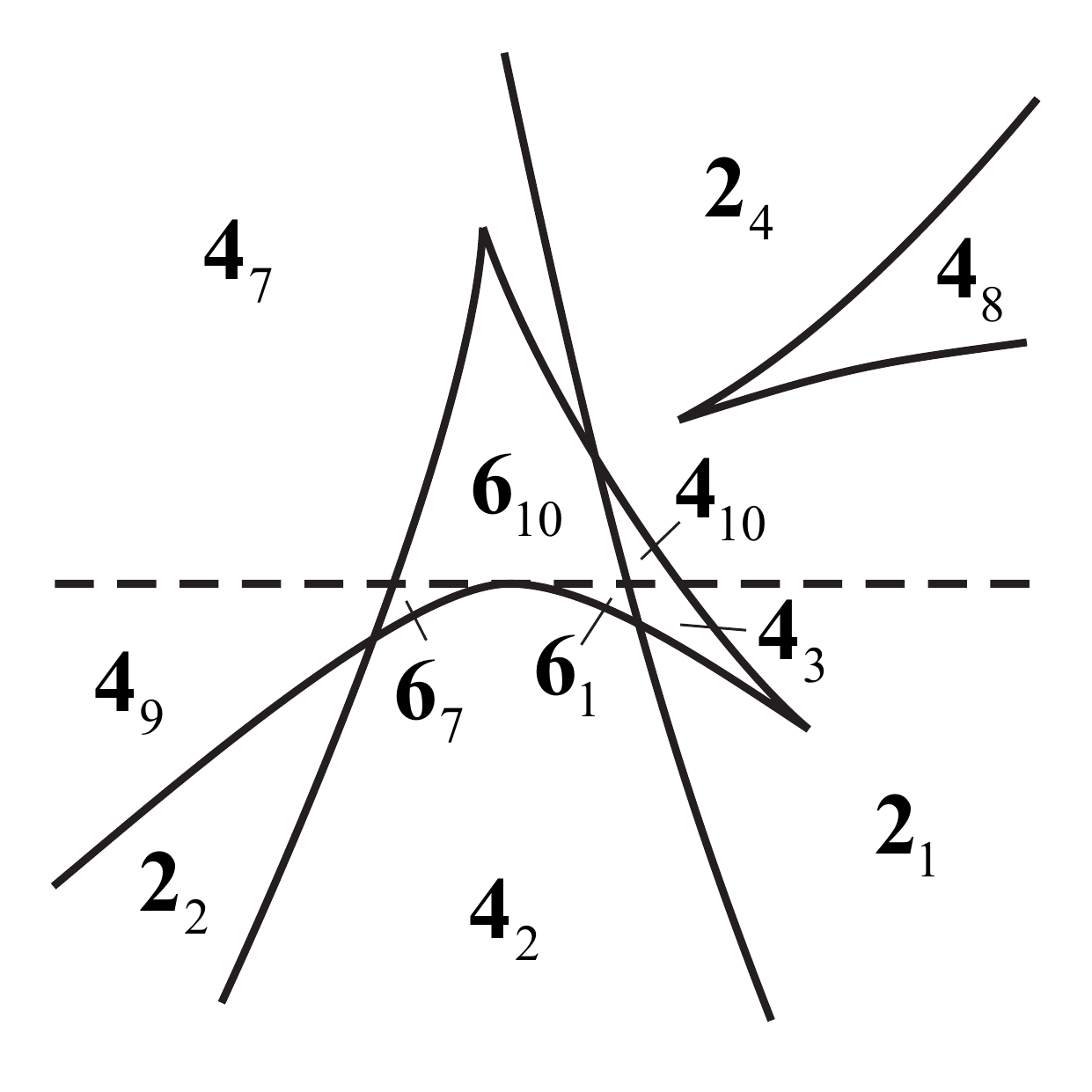}&
41)\includegraphics[width=4cm]{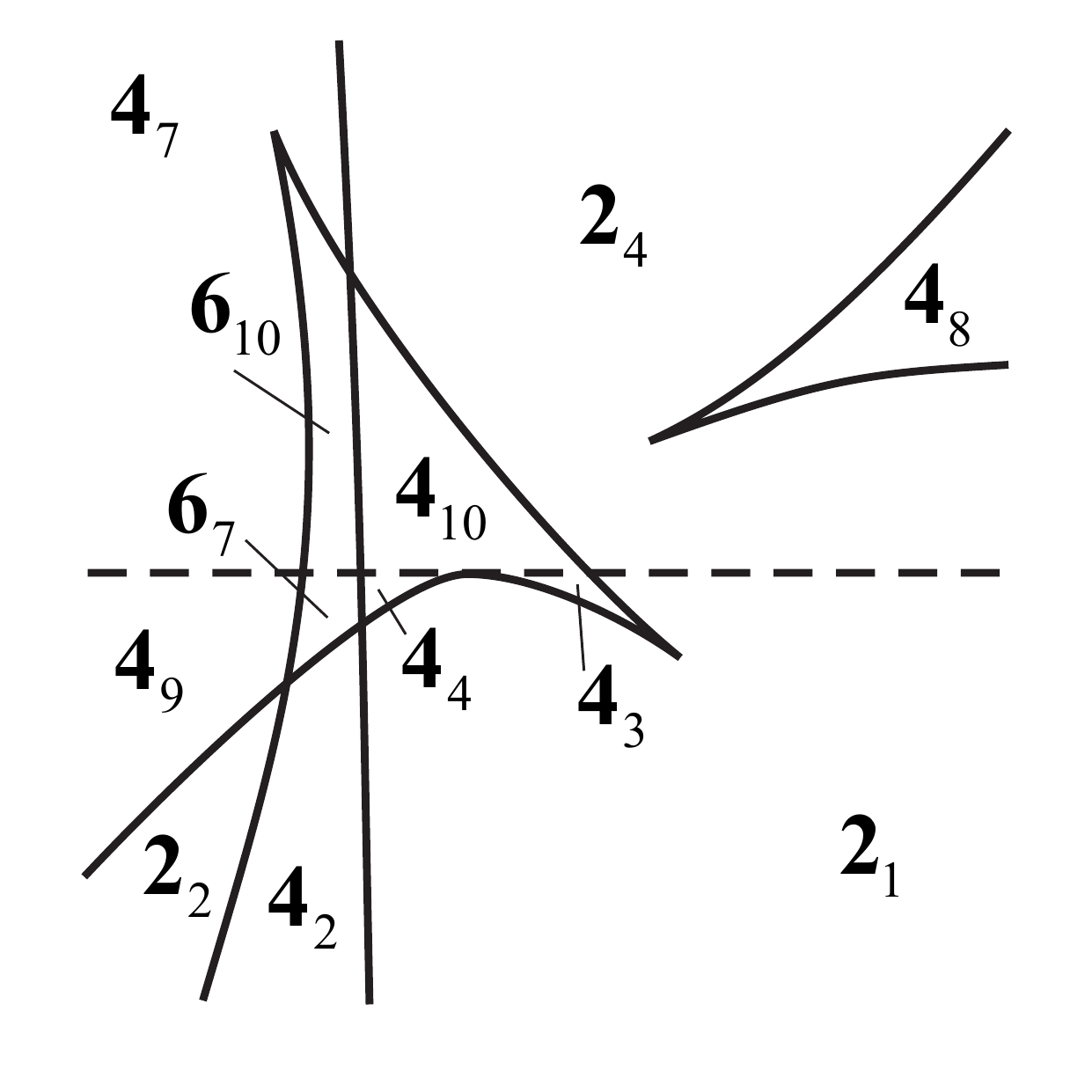}&
42)\includegraphics[width=4cm]{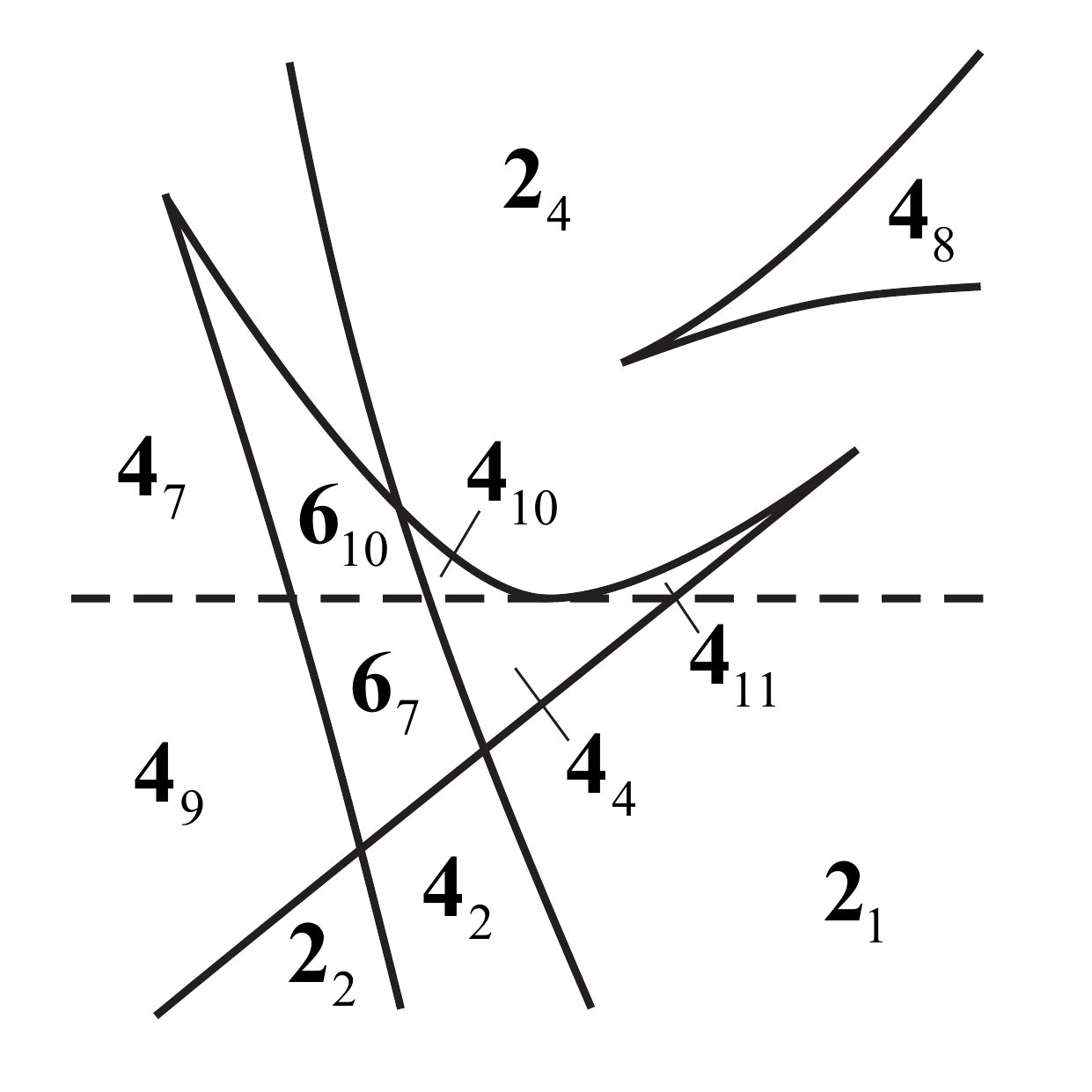}\\
43)\includegraphics[width=4cm]{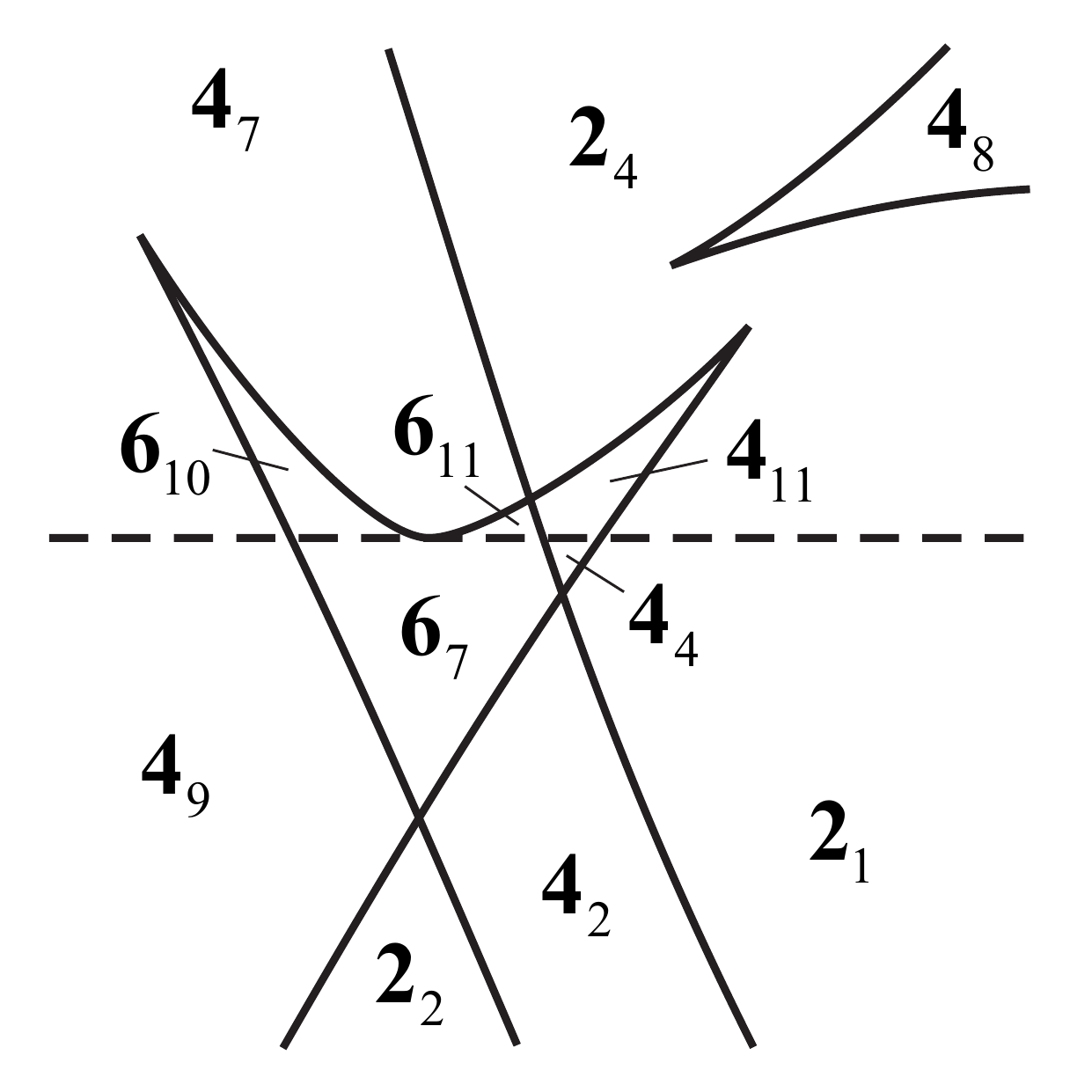}&
44)\includegraphics[width=4cm]{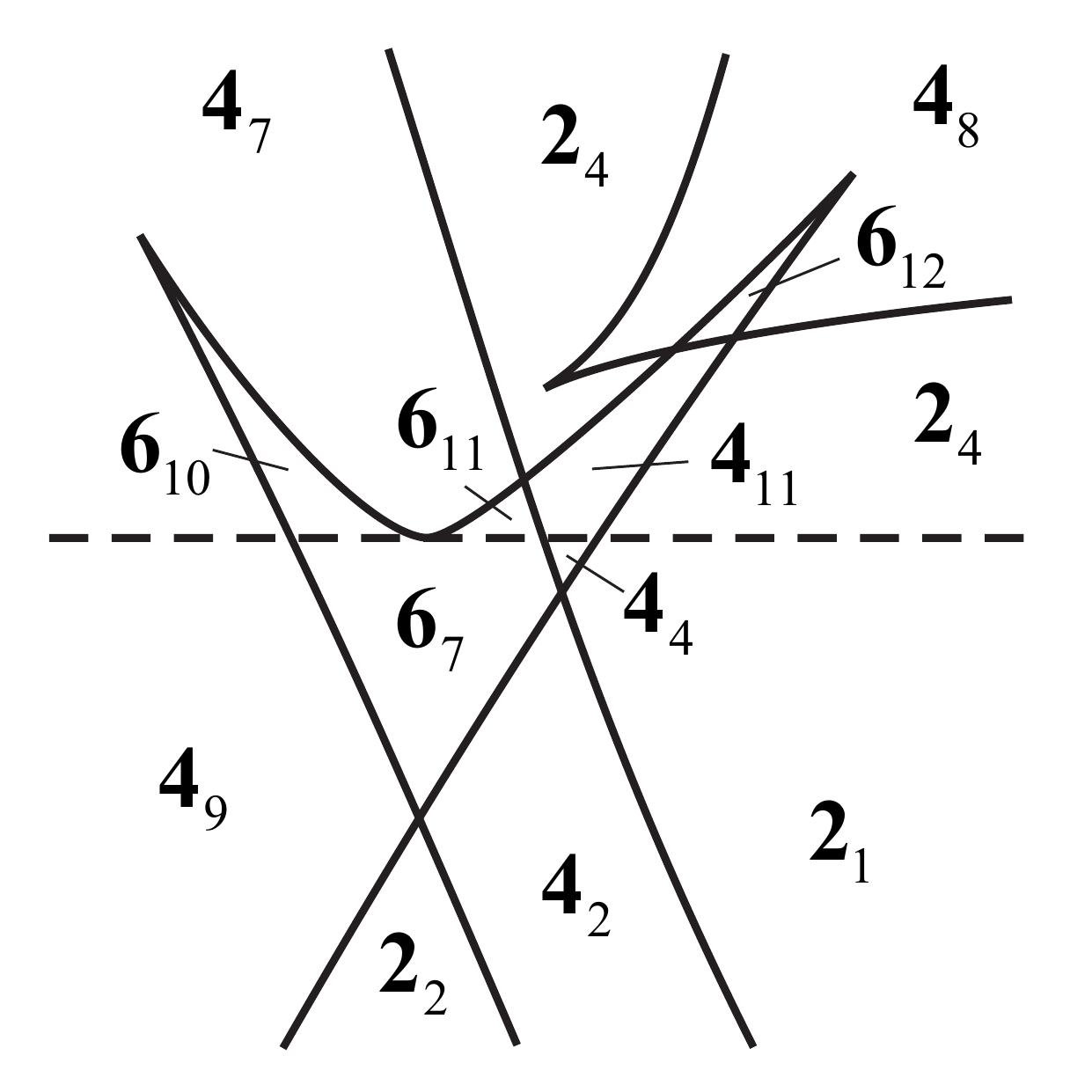}&
45)\includegraphics[width=4cm]{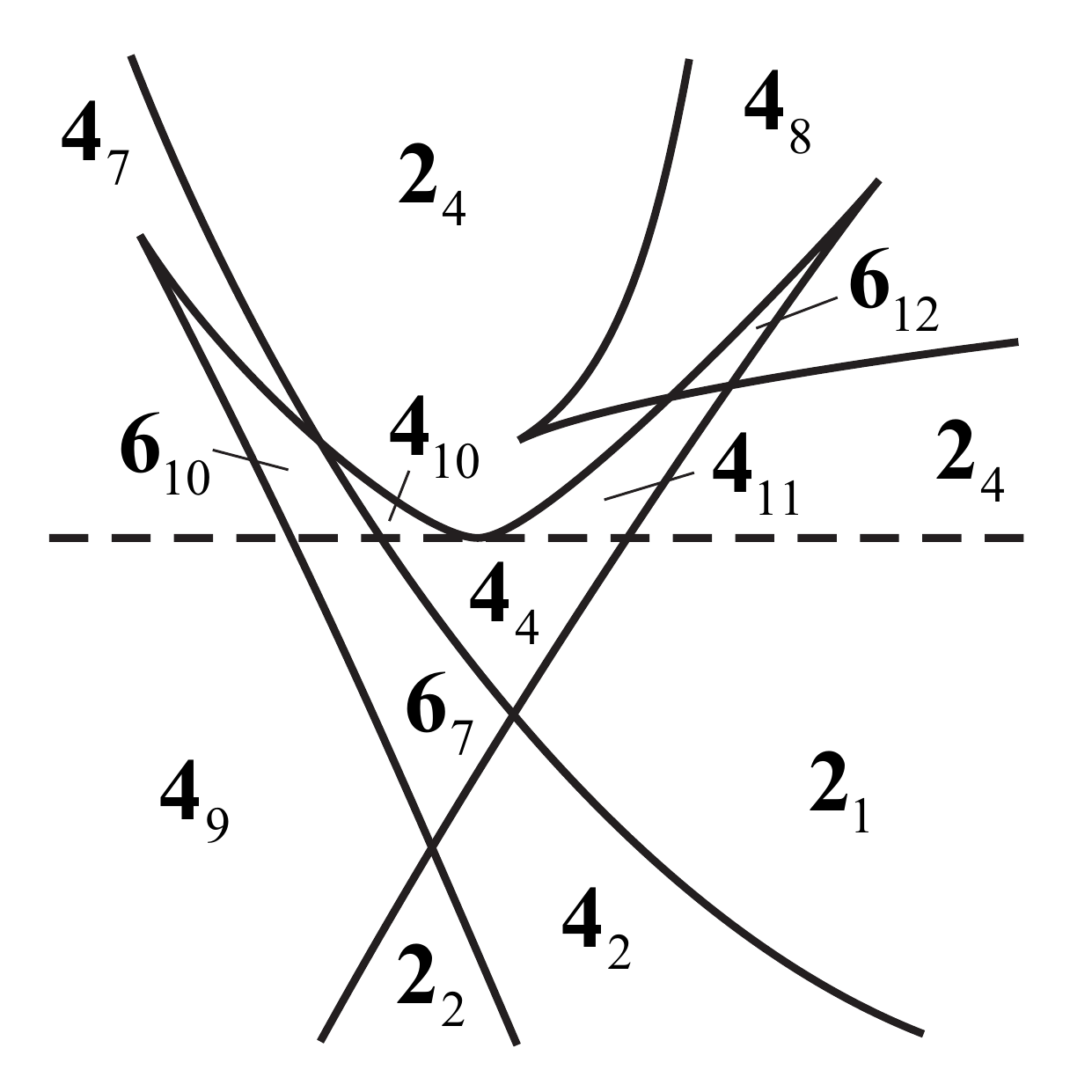}
\end{tabular}
\caption{The skeletons $\Sigma_{S_4,t_1,q_5},S_4\neq0$ for domains $31 - 45$.}
\label{zona31-45}
\end{center}
\end{figure}

\begin{figure}
\begin{center}
\begin{tabular}{ccc}
46)\includegraphics[width=4cm]{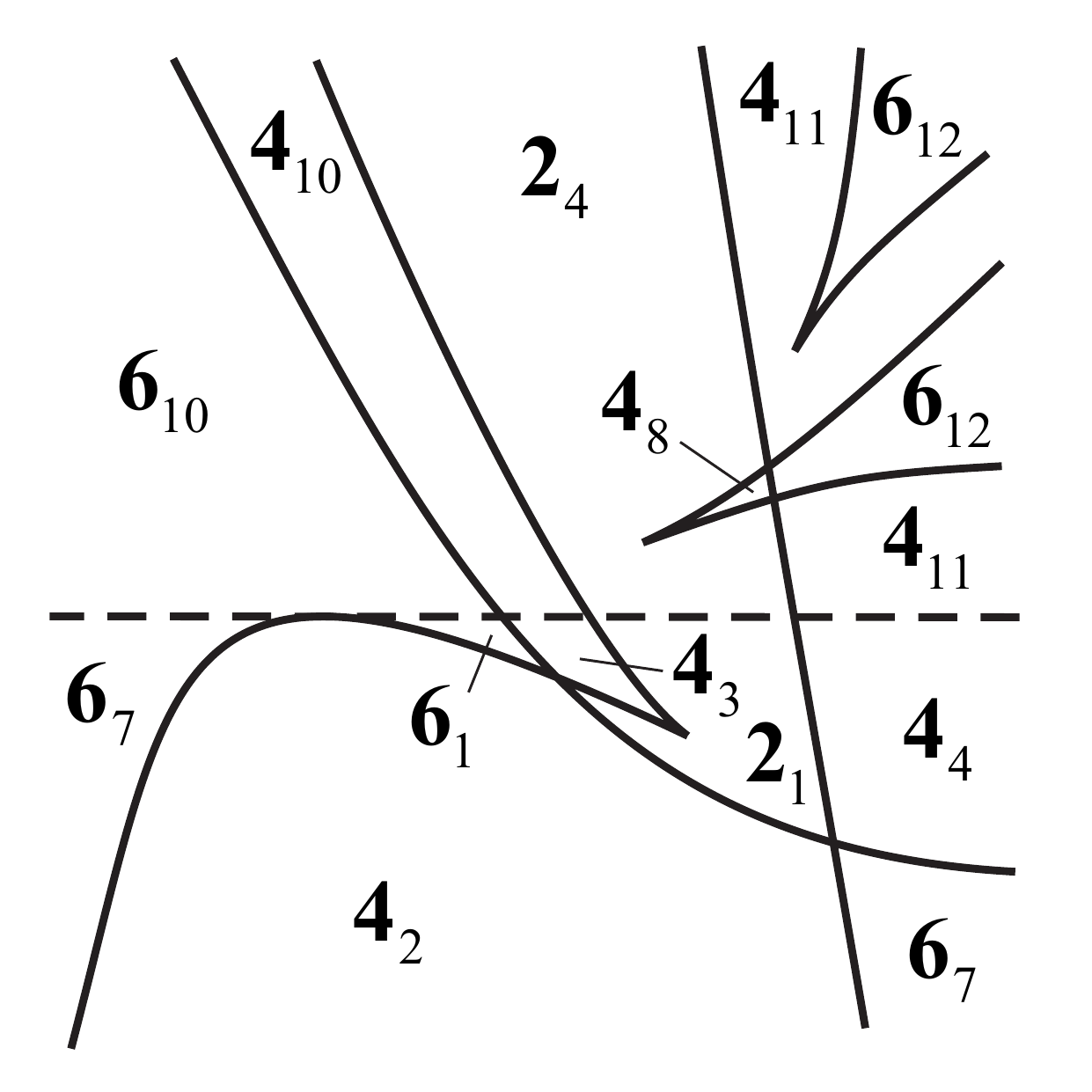}&
47)\includegraphics[width=4cm]{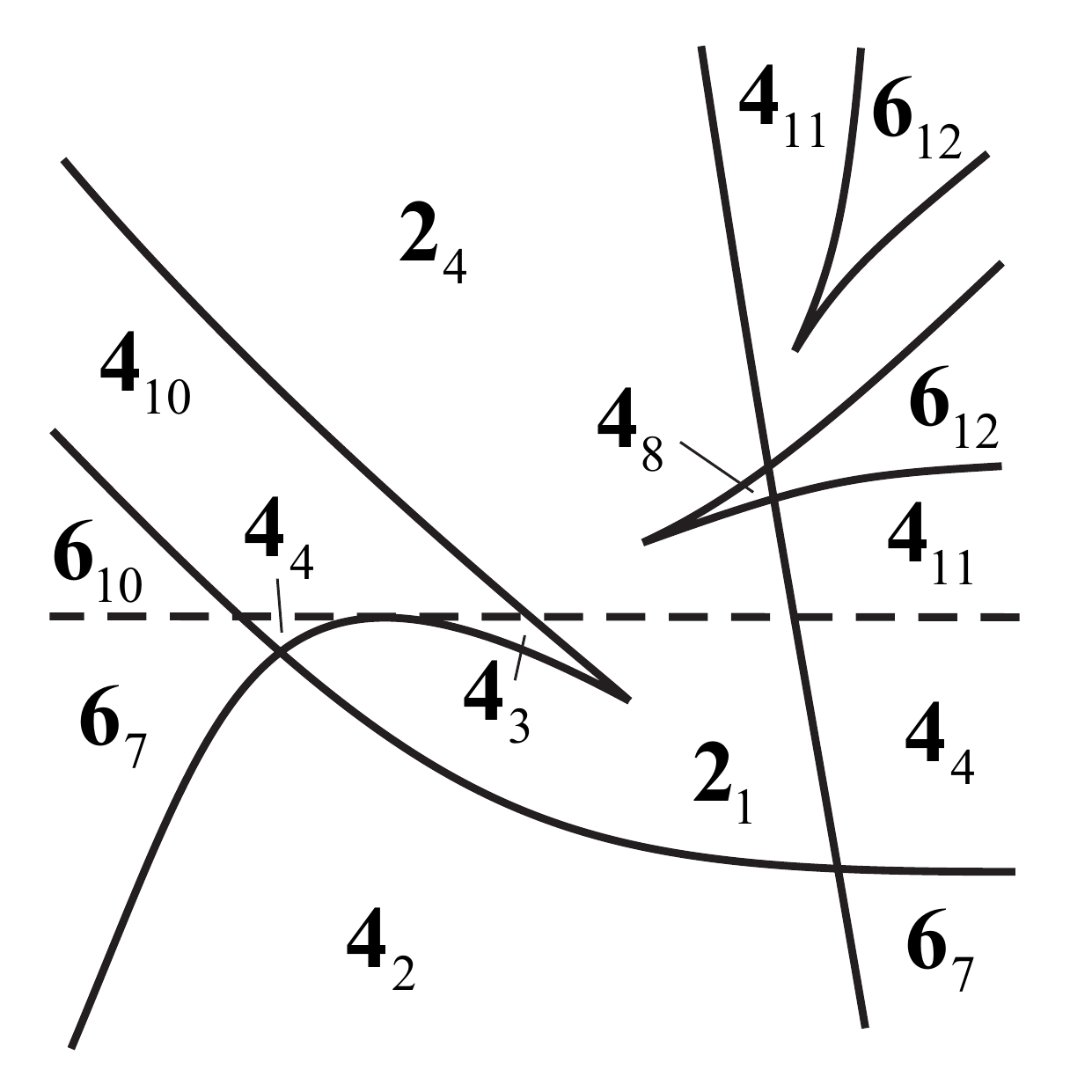}&
48)\includegraphics[width=4cm]{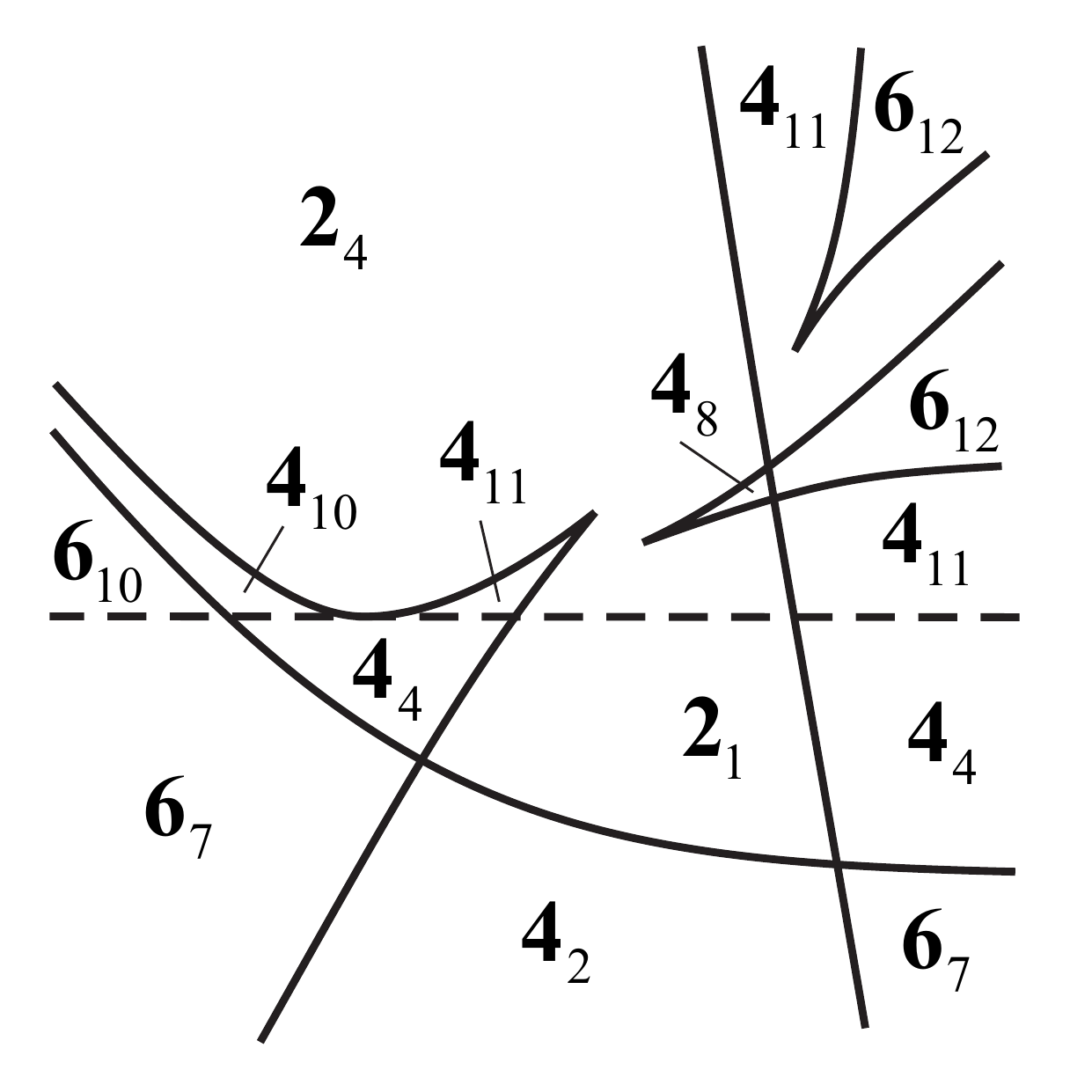}\\
49)\includegraphics[width=4cm]{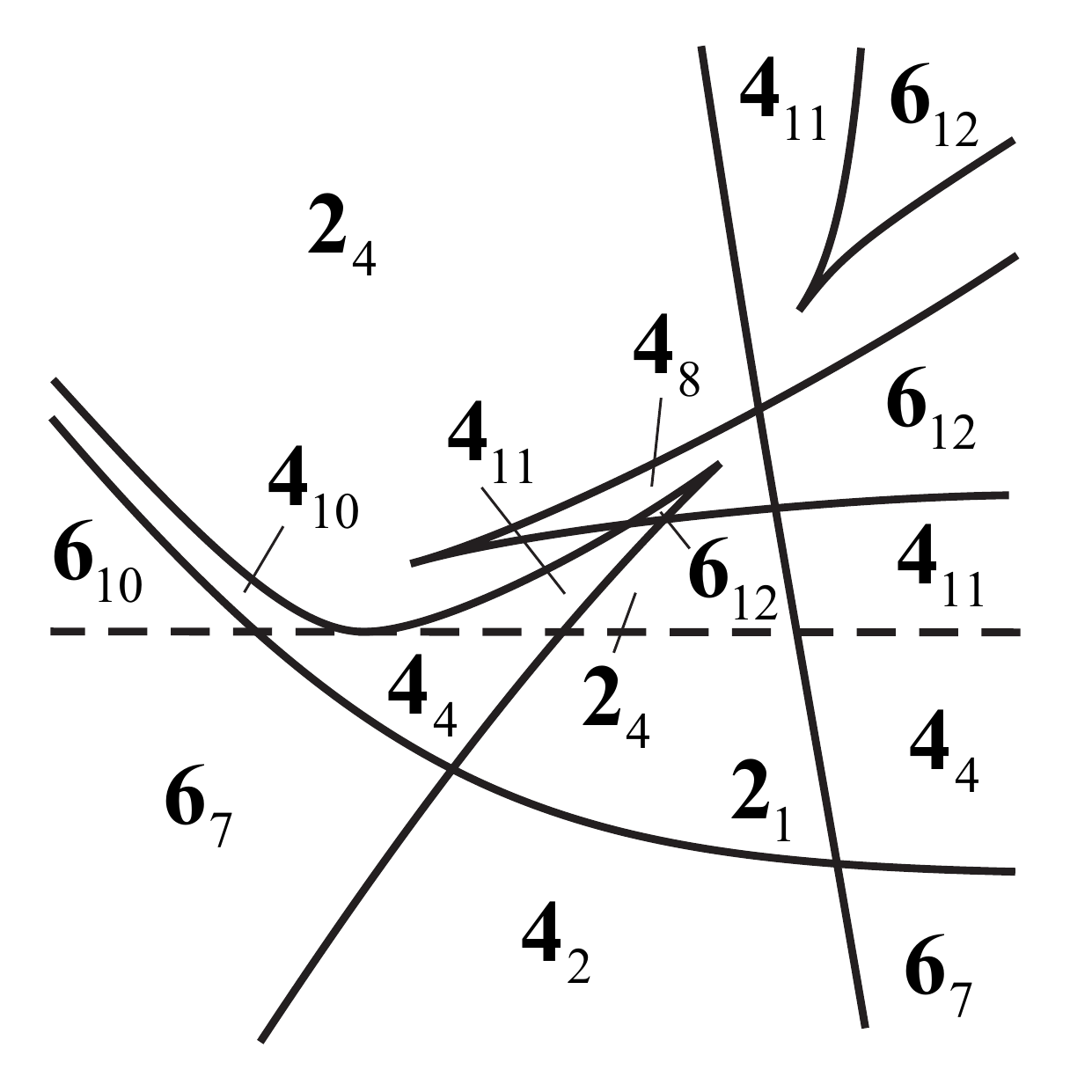}&
50)\includegraphics[width=4cm]{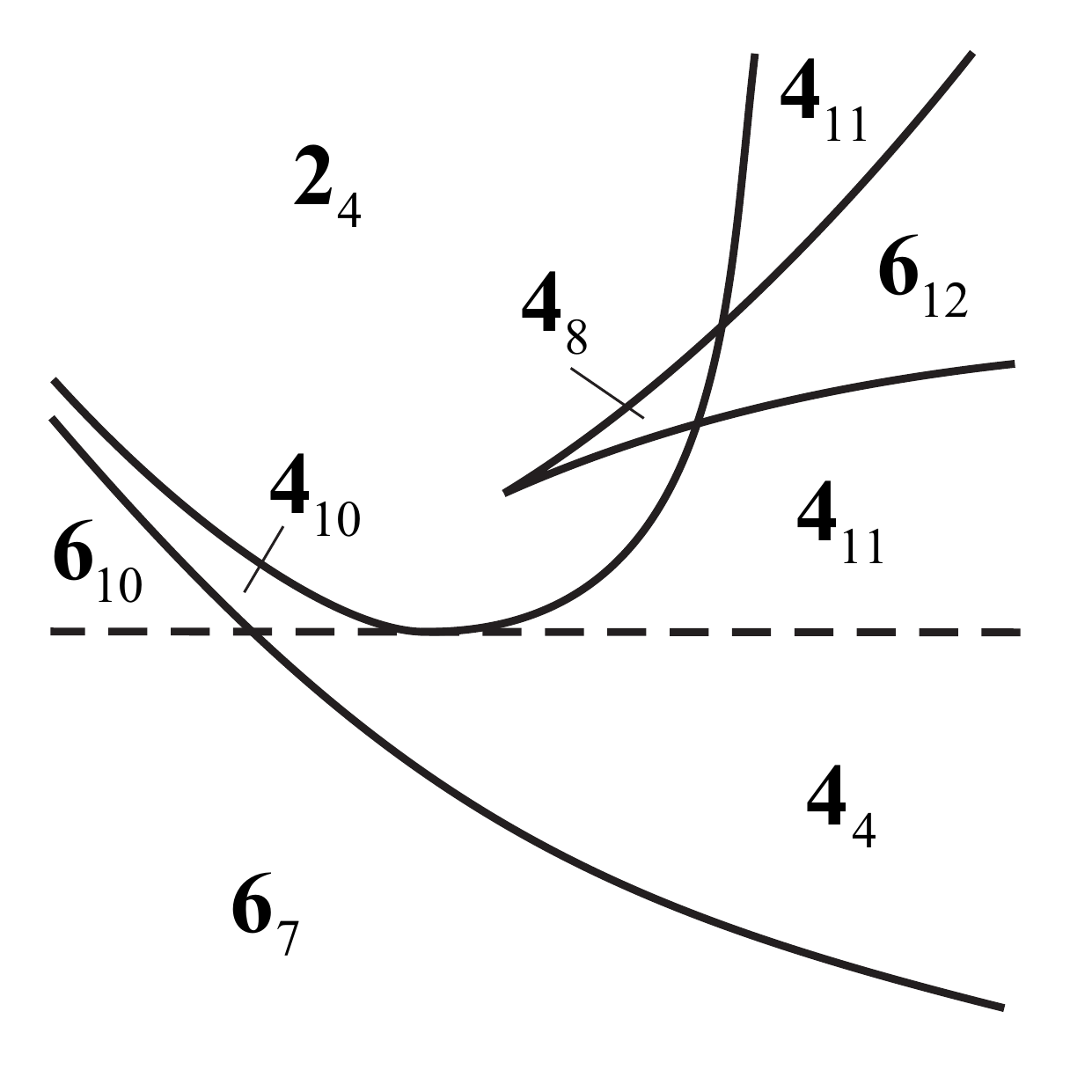}&
51)\includegraphics[width=4cm]{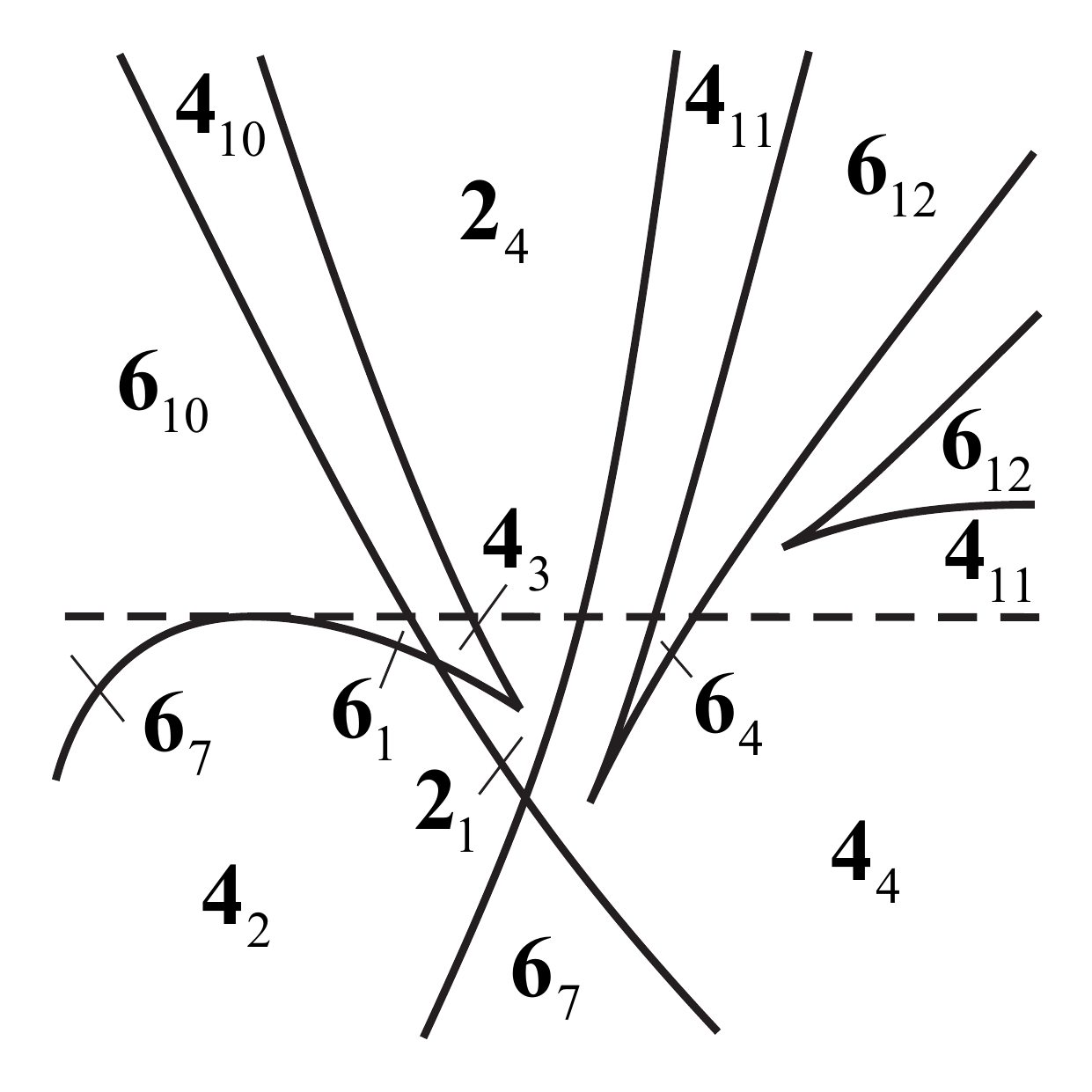}\\
52)\includegraphics[width=4cm]{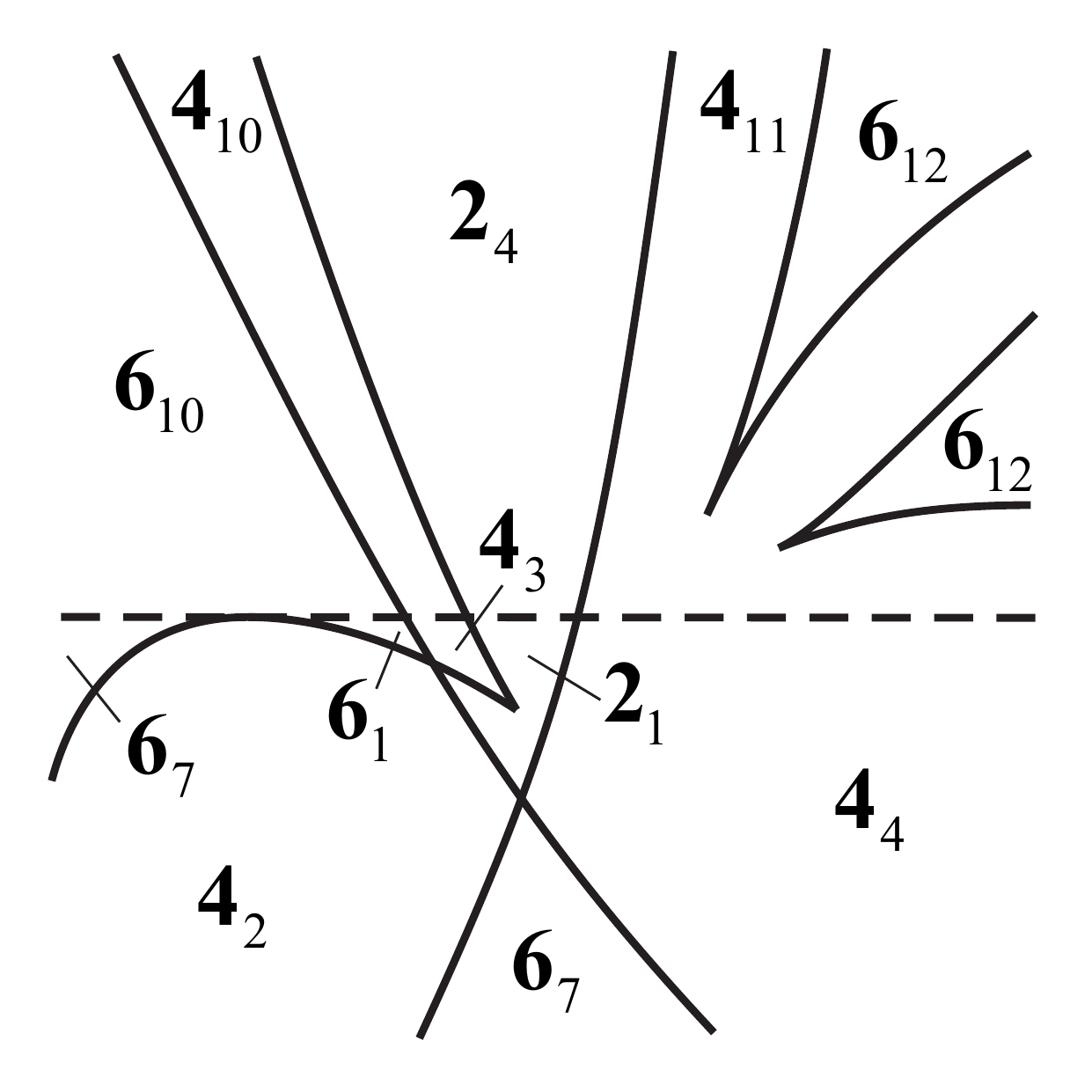}&
53)\includegraphics[width=4cm]{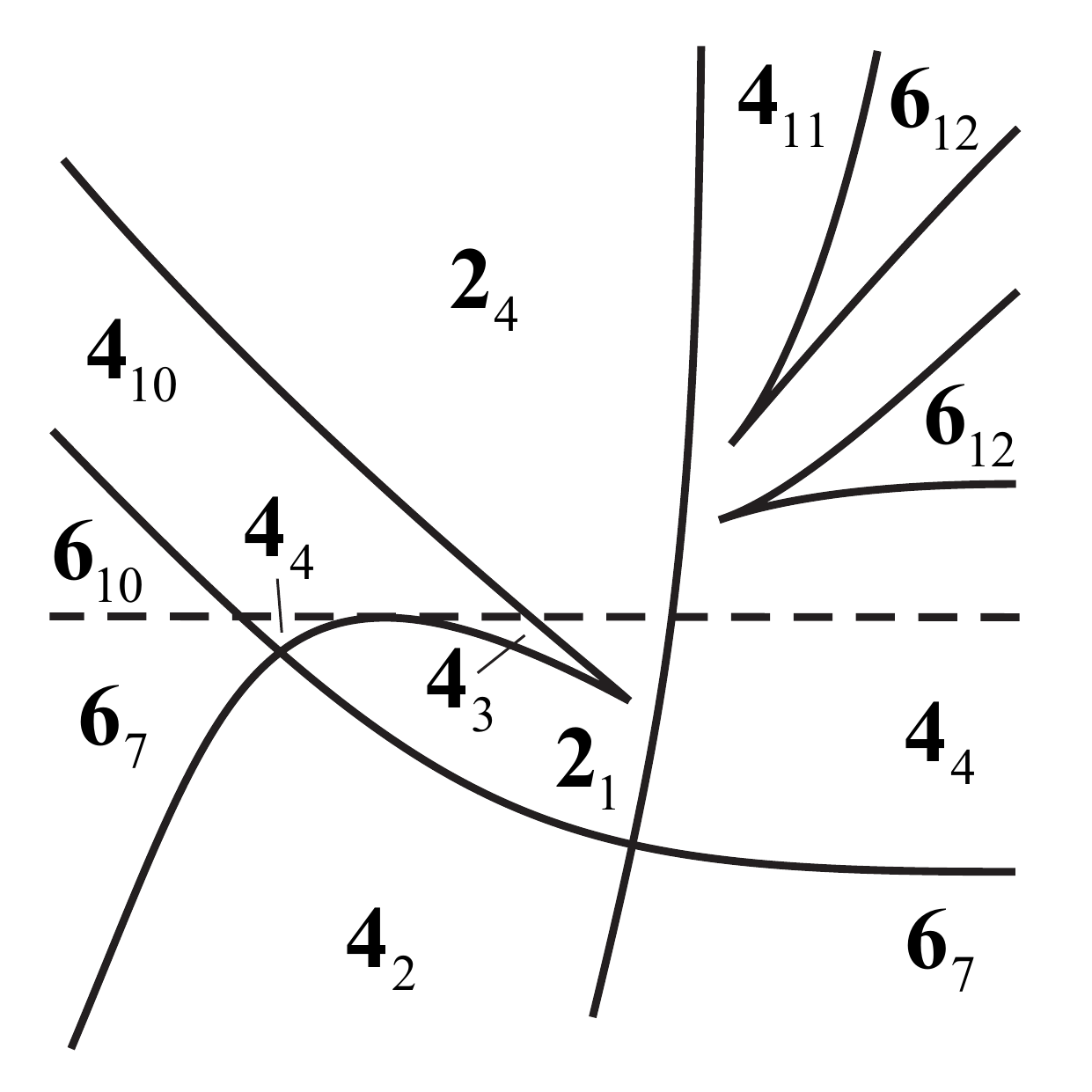}&
54)\includegraphics[width=4cm]{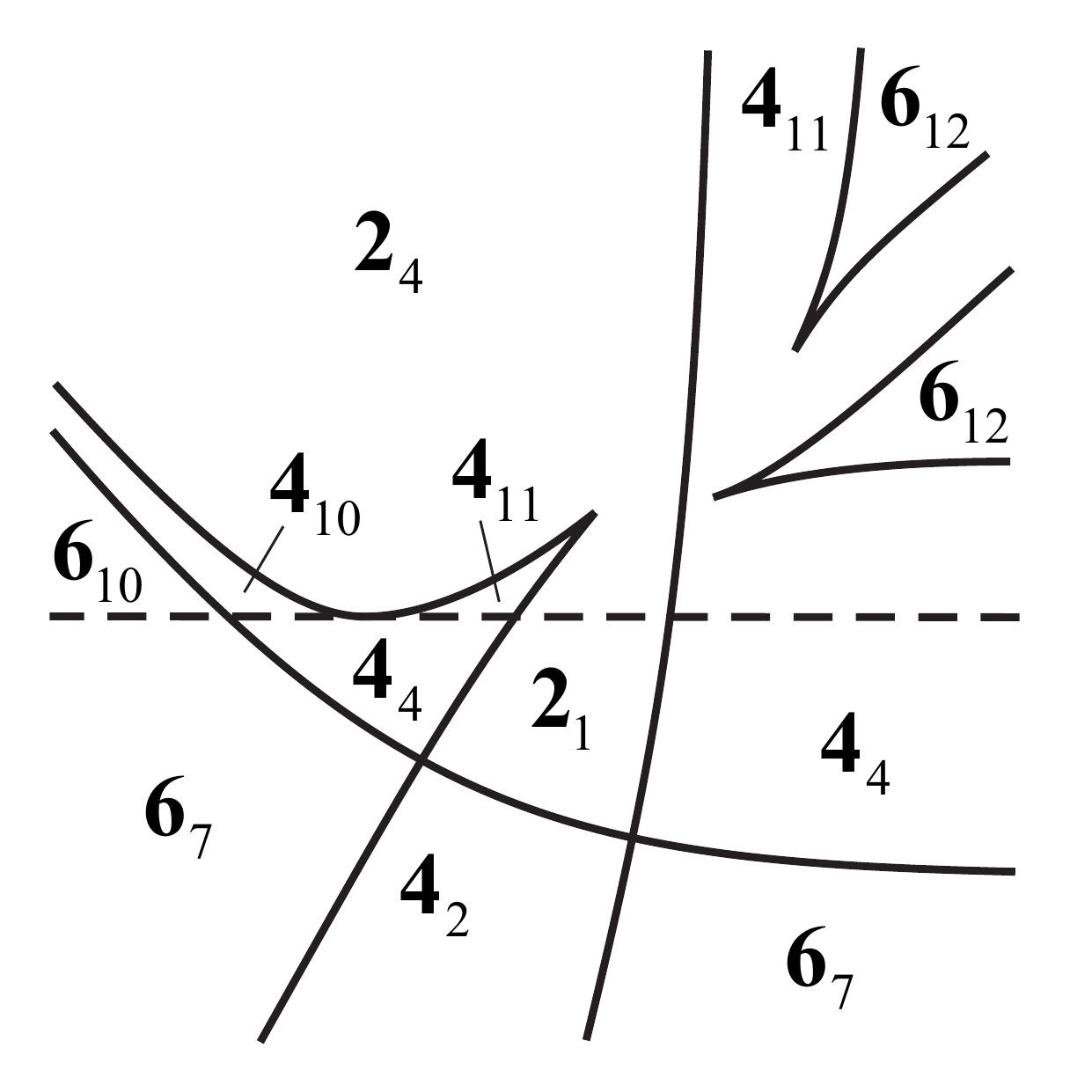}\\
55)\includegraphics[width=4cm]{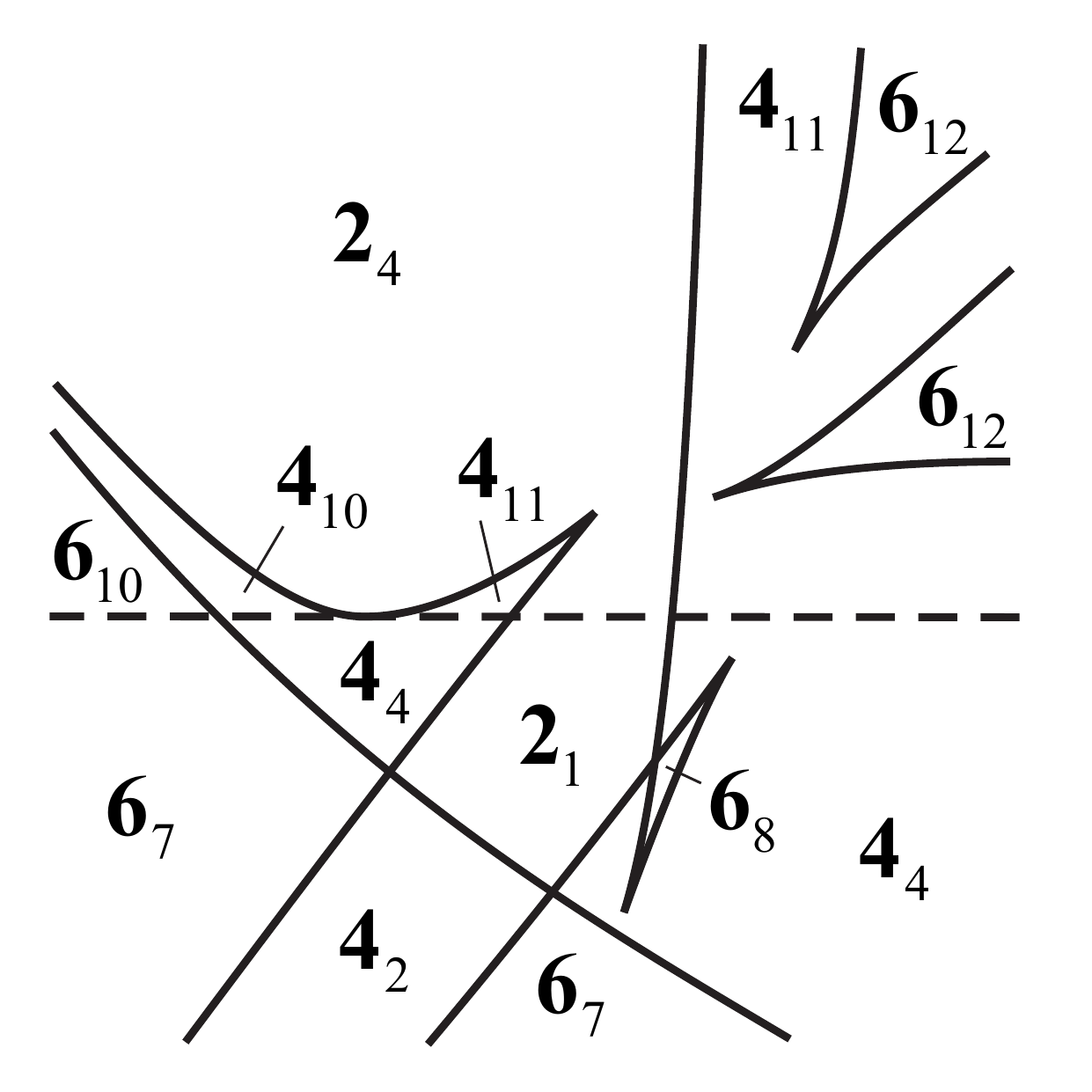}&
56)\includegraphics[width=4cm]{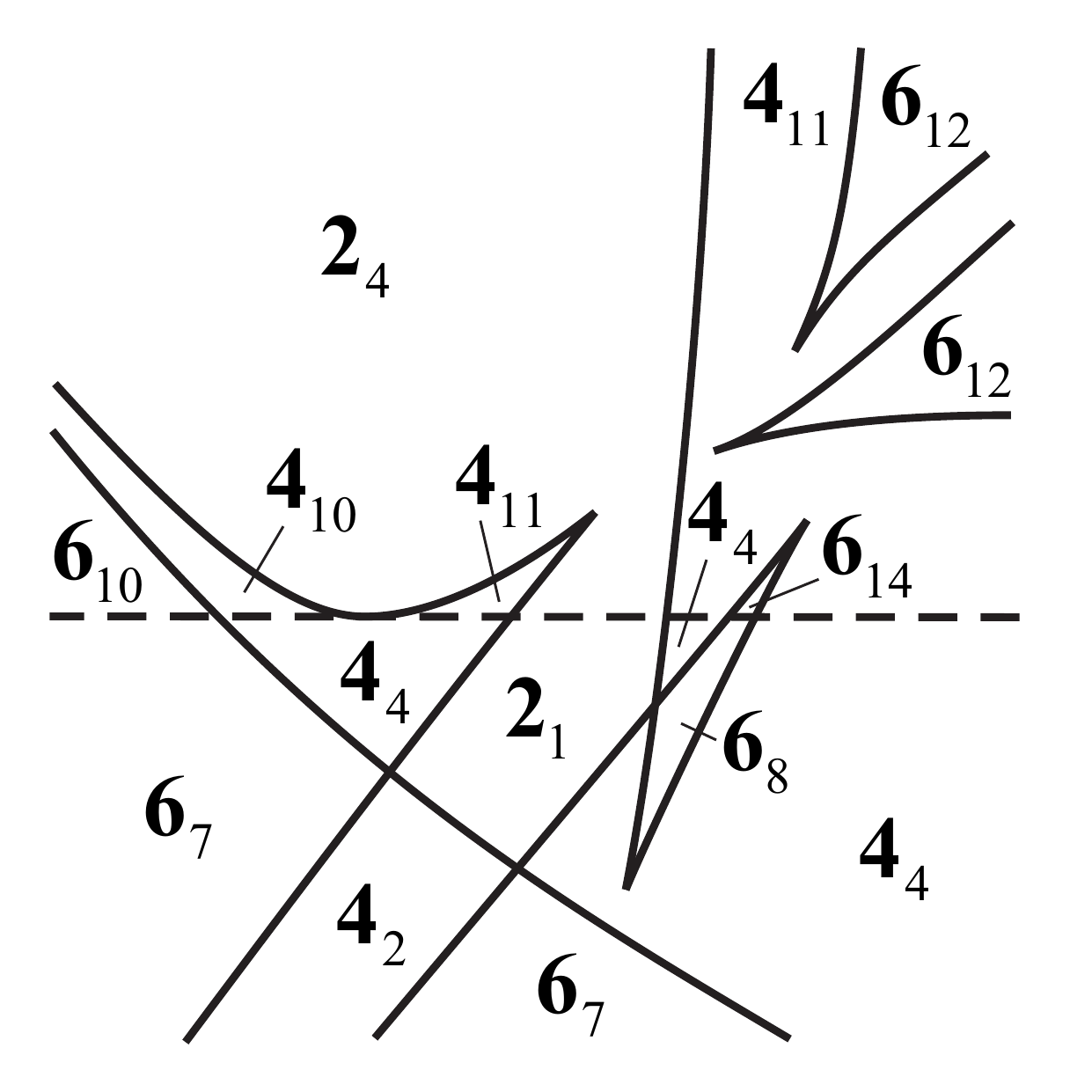}&
57)\includegraphics[width=4cm]{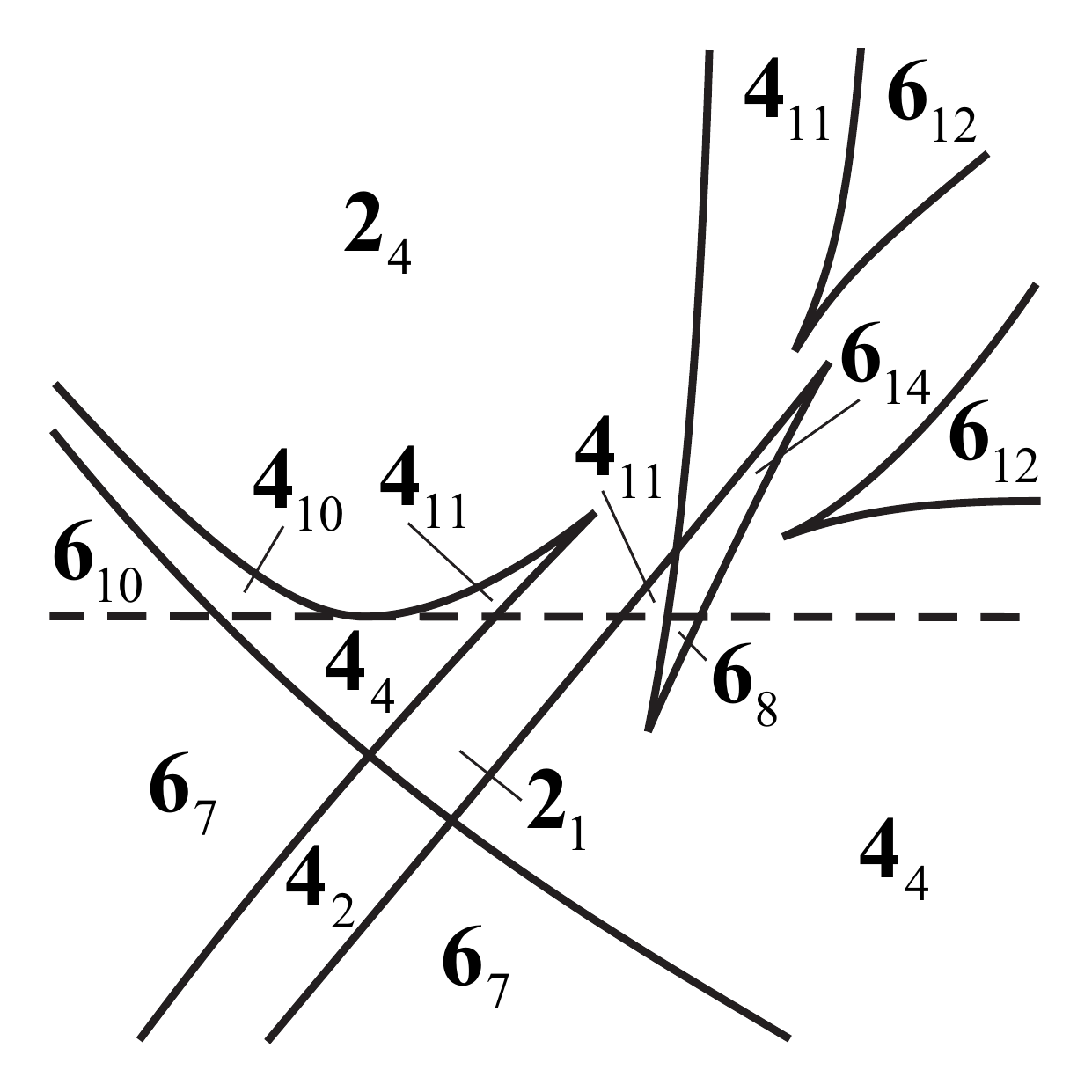}\\
58)\includegraphics[width=4cm]{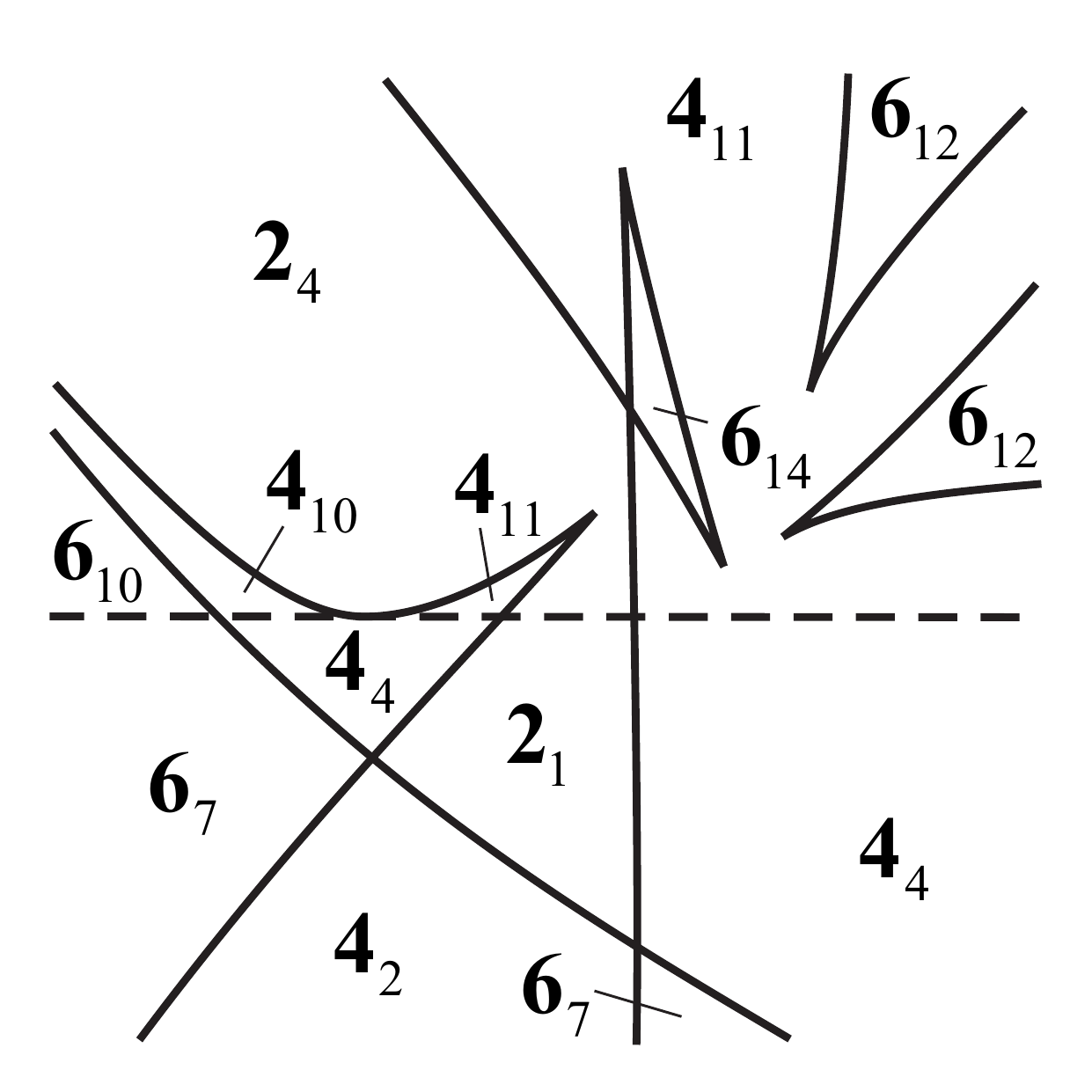}&
59)\includegraphics[width=4cm]{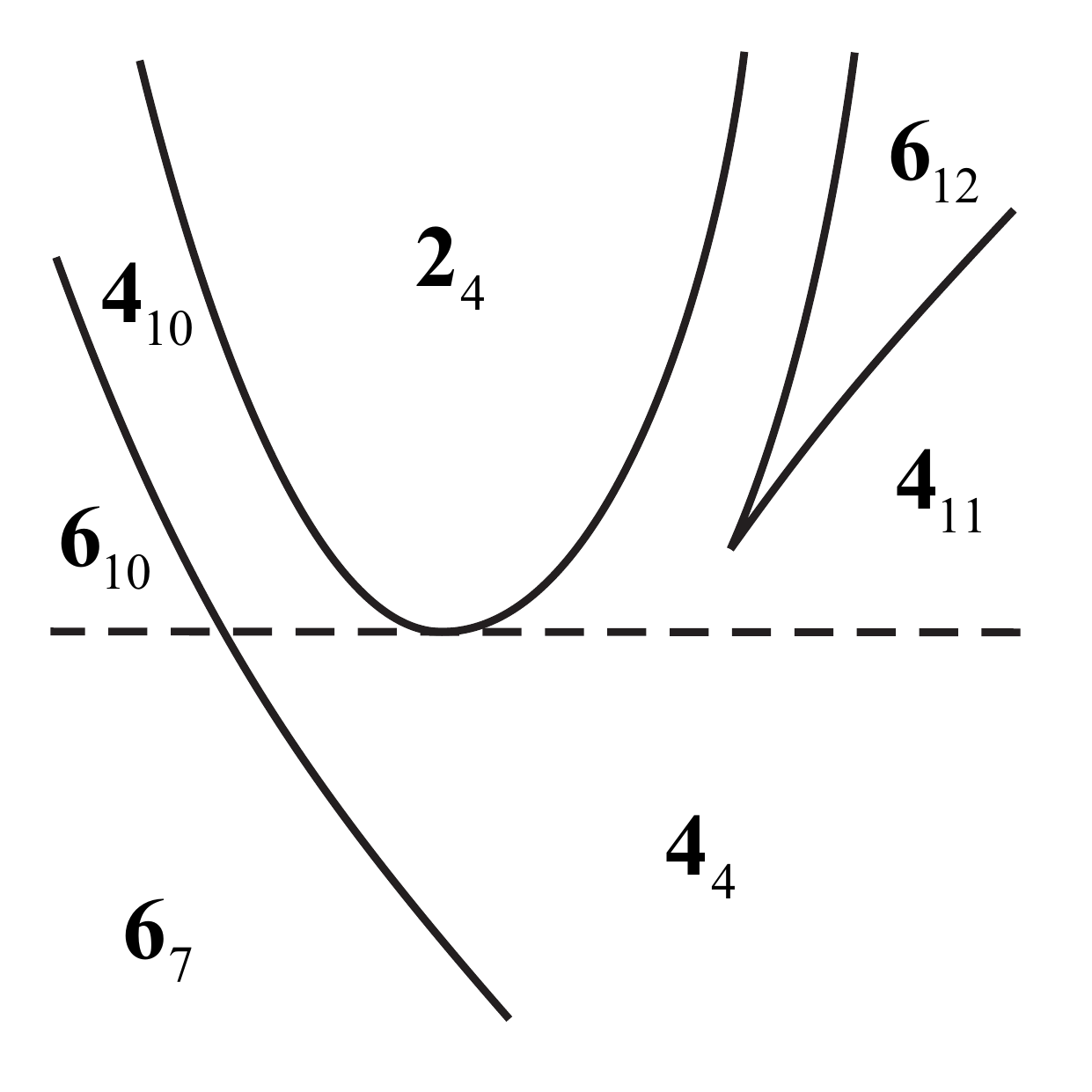}&
60)\includegraphics[width=4cm]{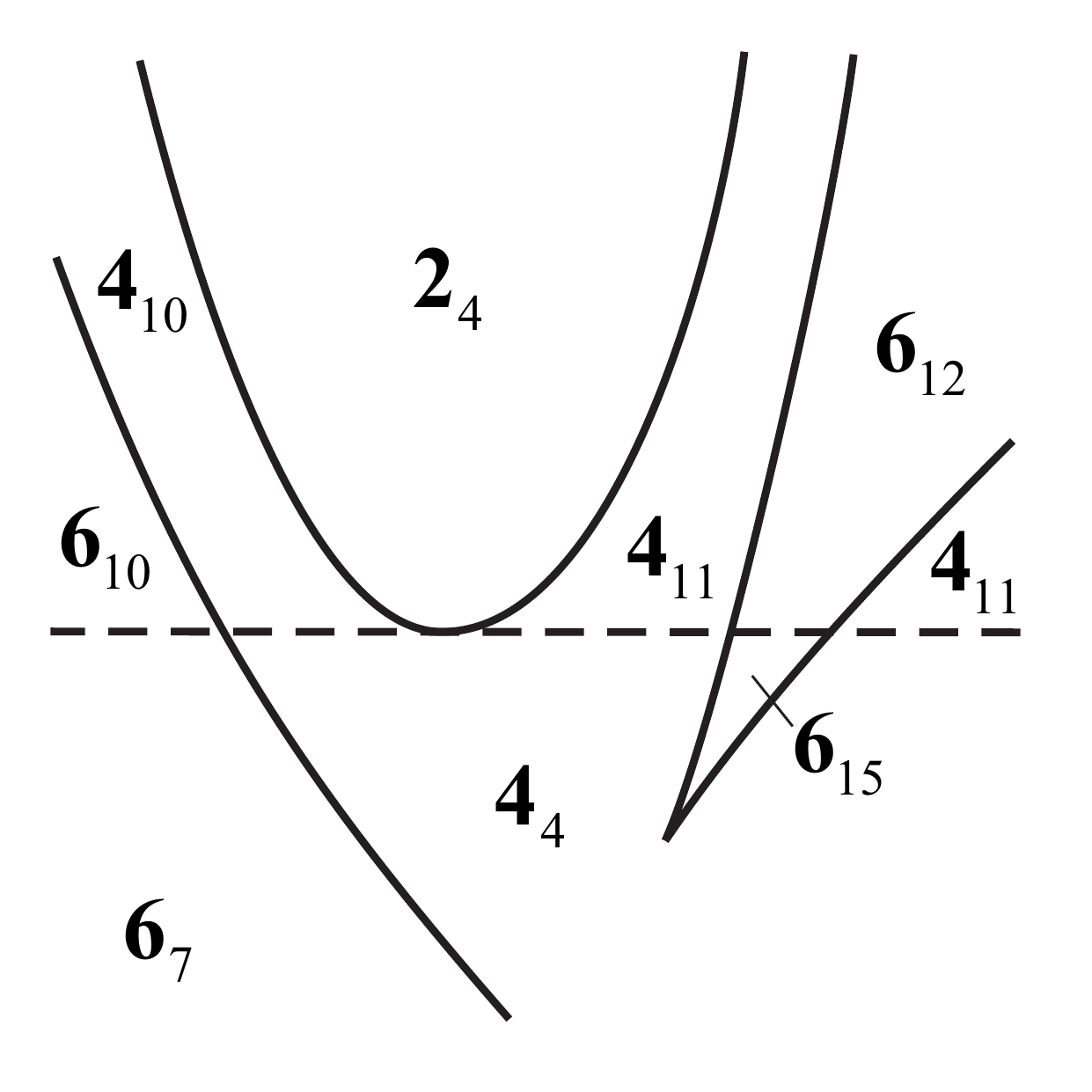}
\end{tabular}
\caption{The skeletons $\Sigma_{S_4,t_1,q_5},S_4\neq0$ for domains $46 - 60$.}
\label{zona46-60}
\end{center}
\end{figure}

\begin{figure}
\begin{center}
\begin{tabular}{ccc}
61)\includegraphics[width=4cm]{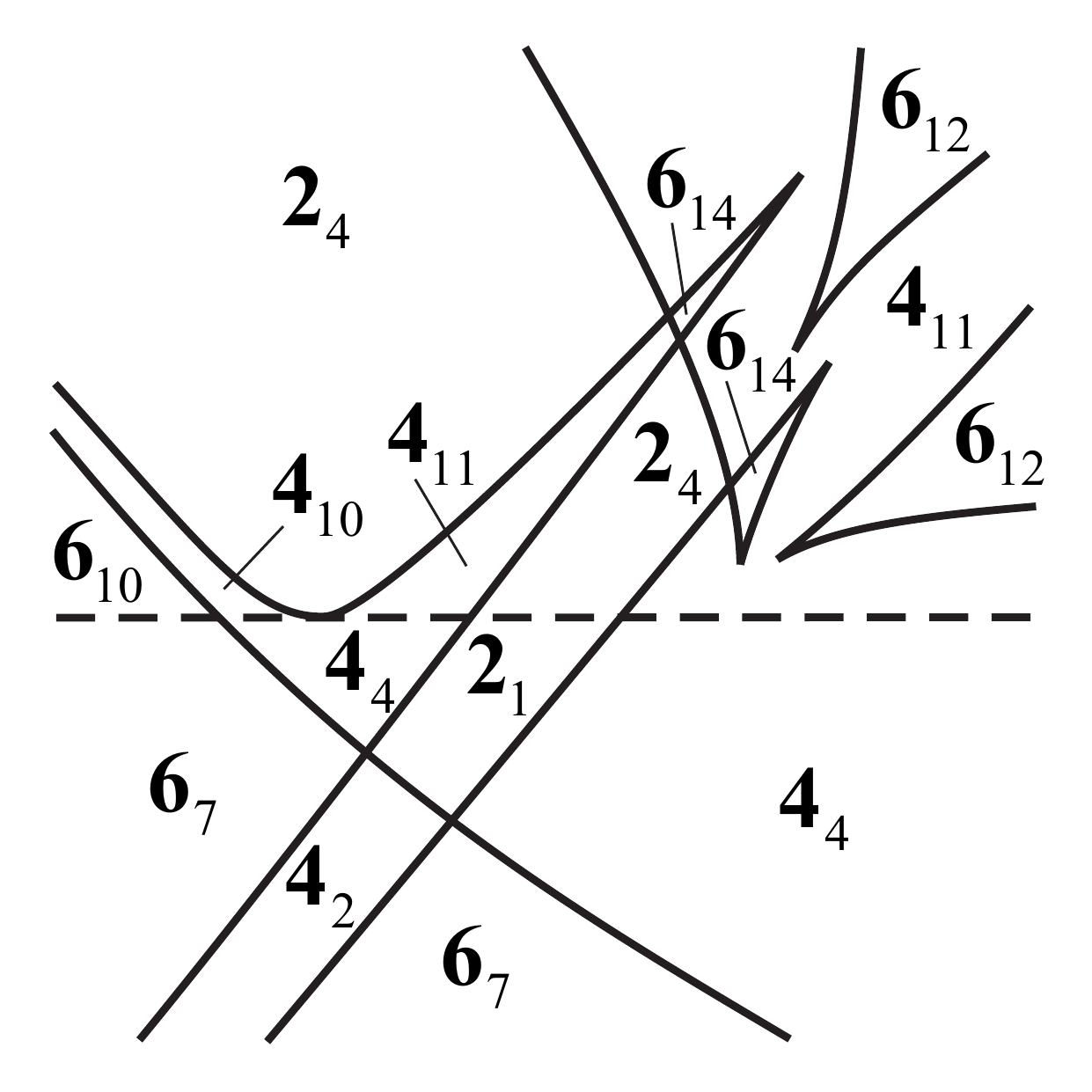}&
62)\includegraphics[width=4cm]{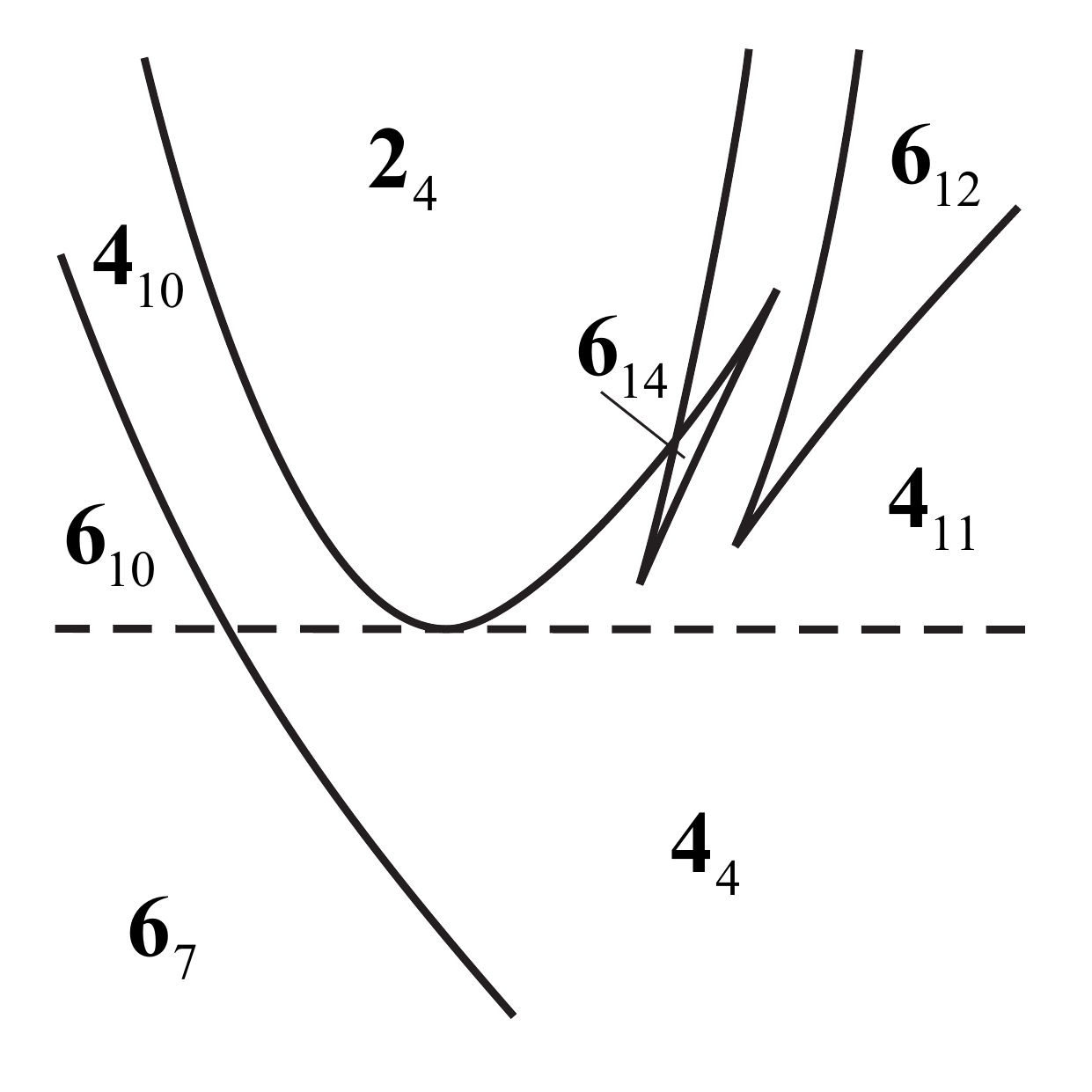}&
63)\includegraphics[width=4cm]{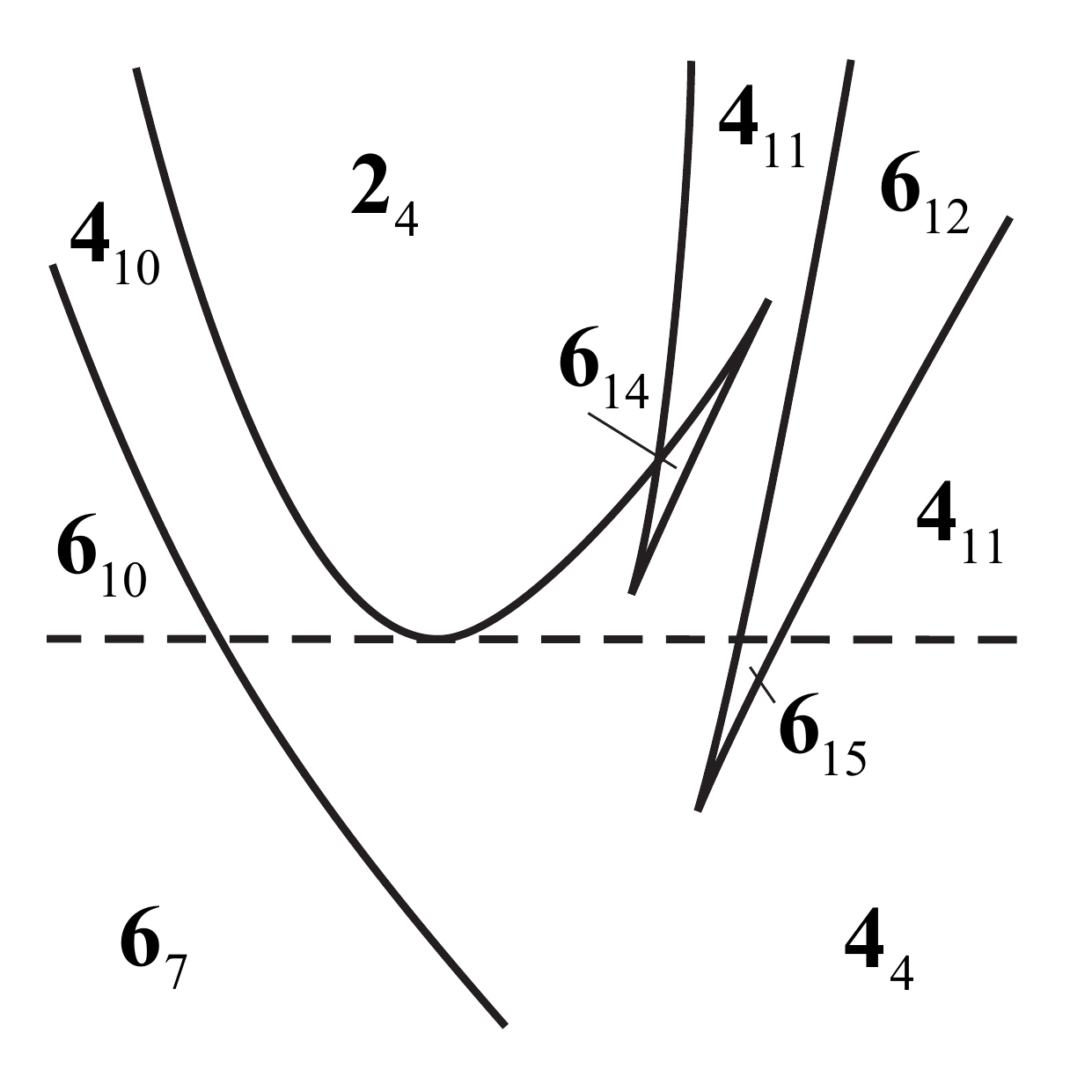}\\
64)\includegraphics[width=4cm]{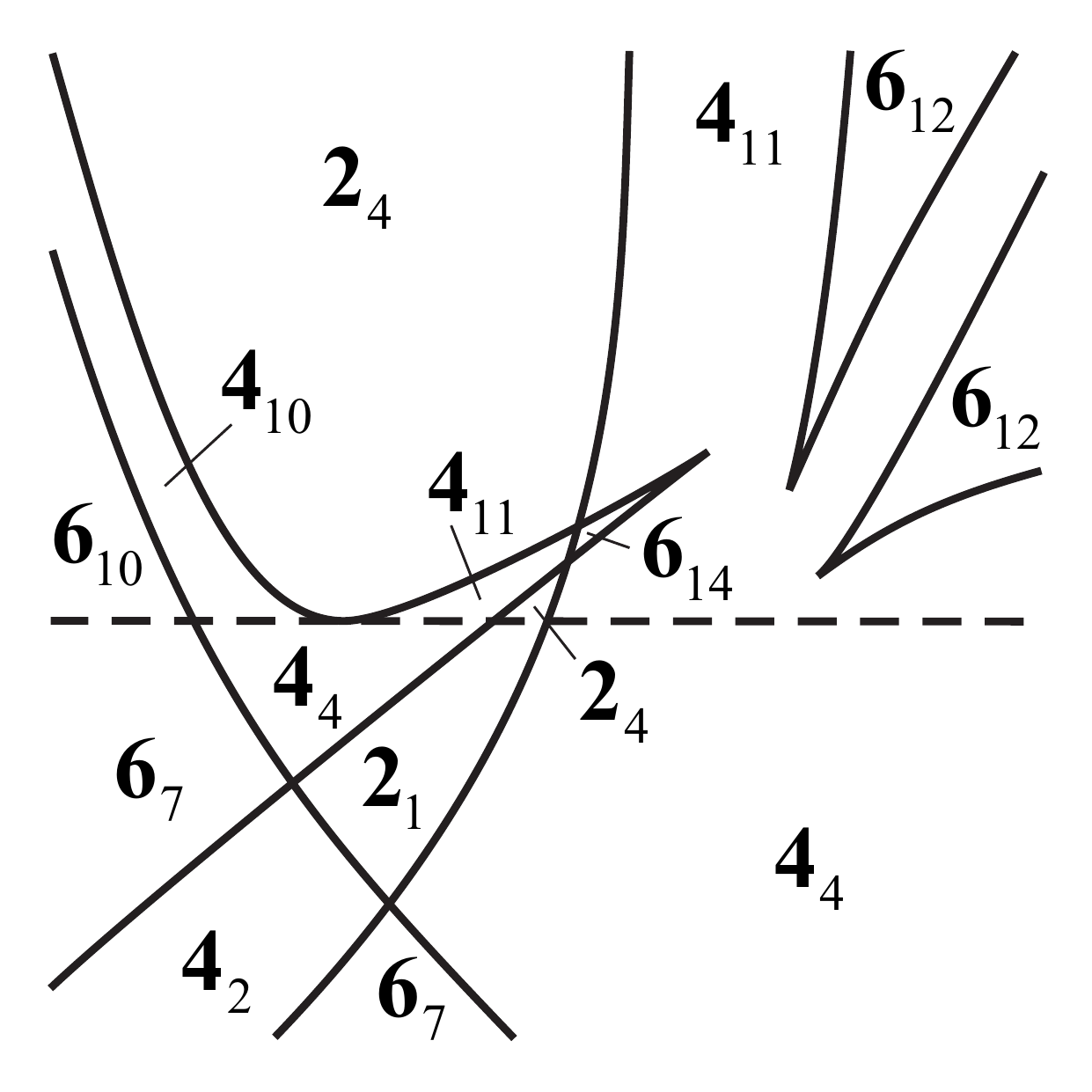}&
65)\includegraphics[width=4cm]{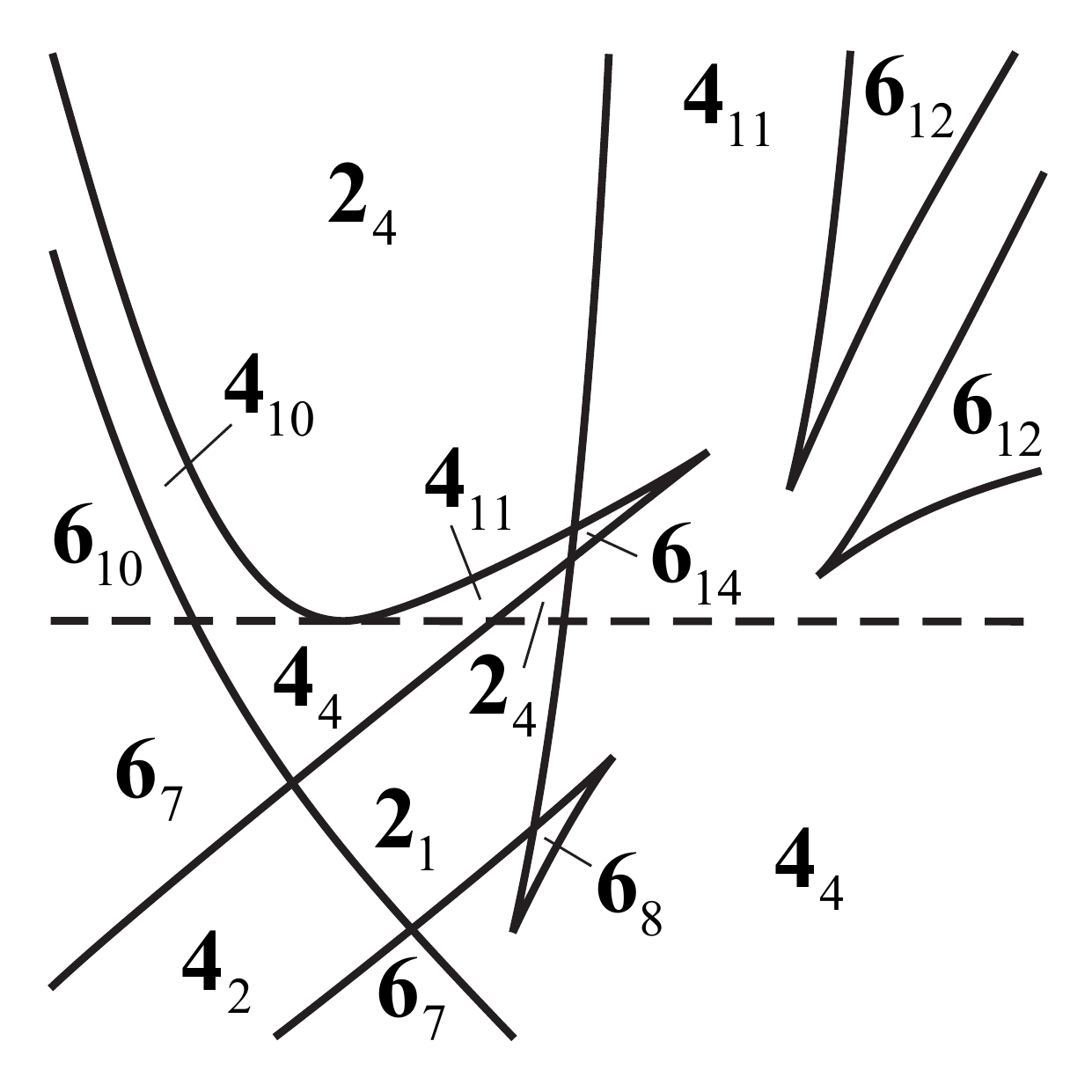}&
66)\includegraphics[width=4cm]{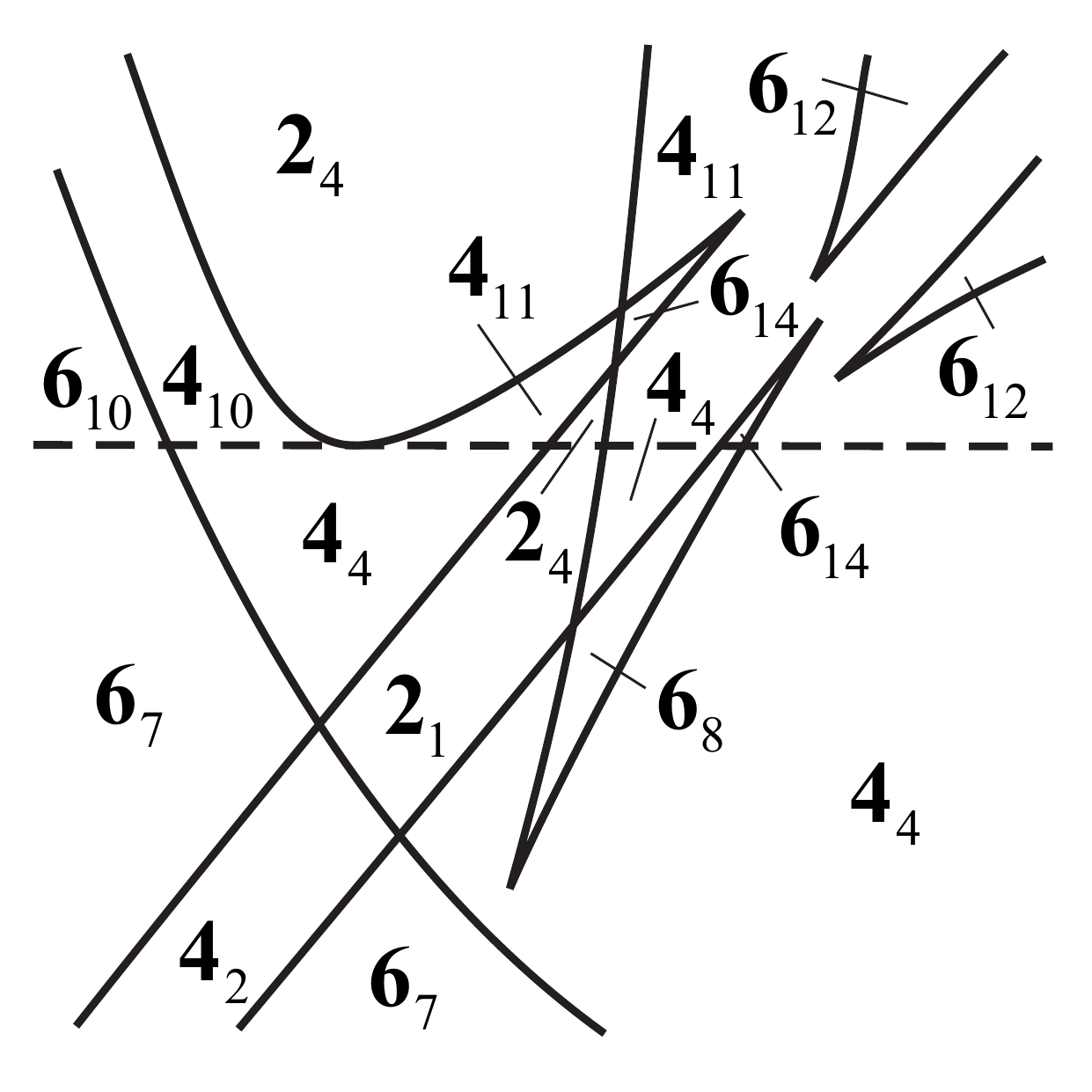}\\
67)\includegraphics[width=4cm]{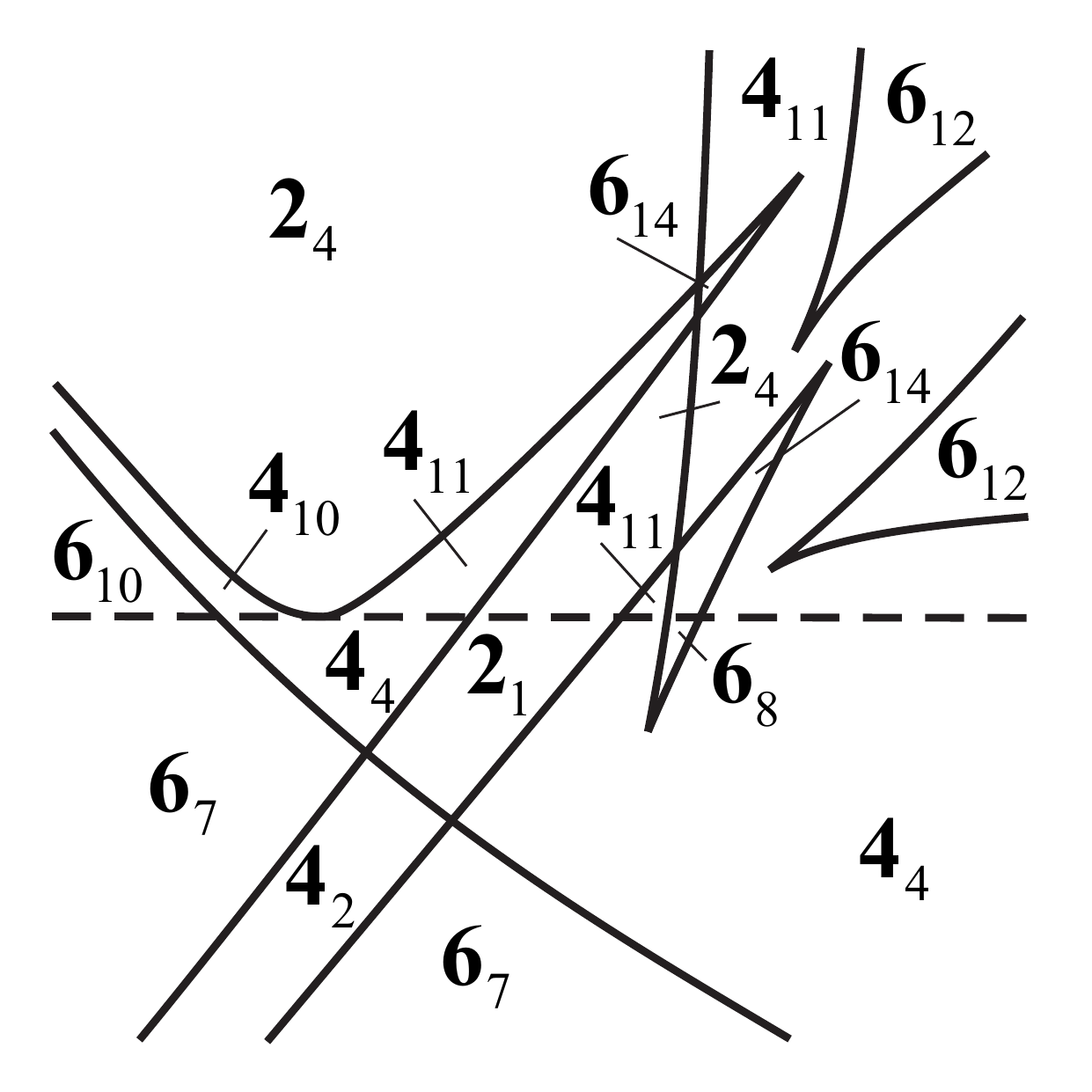}&
68)\includegraphics[width=4cm]{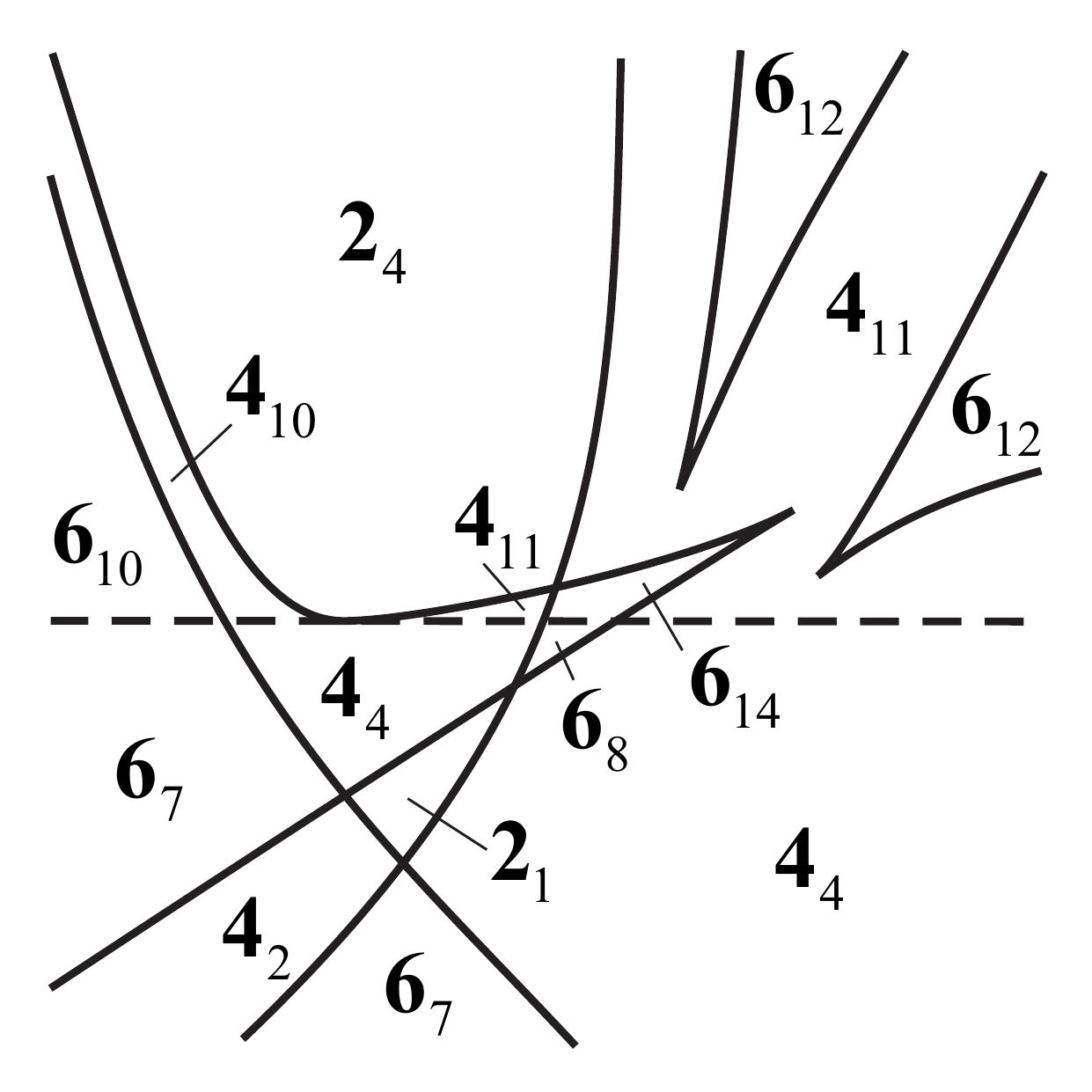}&
69)\includegraphics[width=4cm]{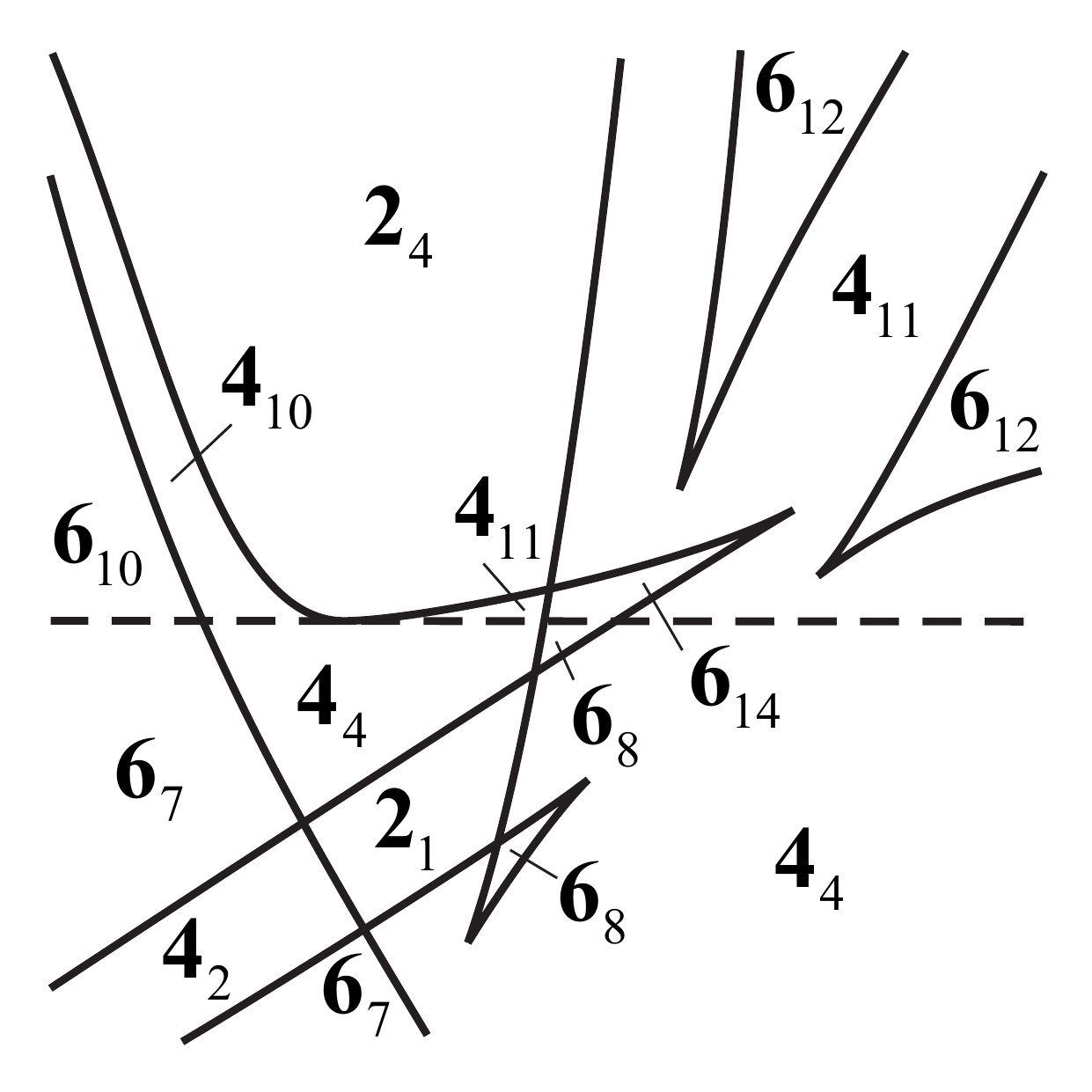}\\
70)\includegraphics[width=4cm]{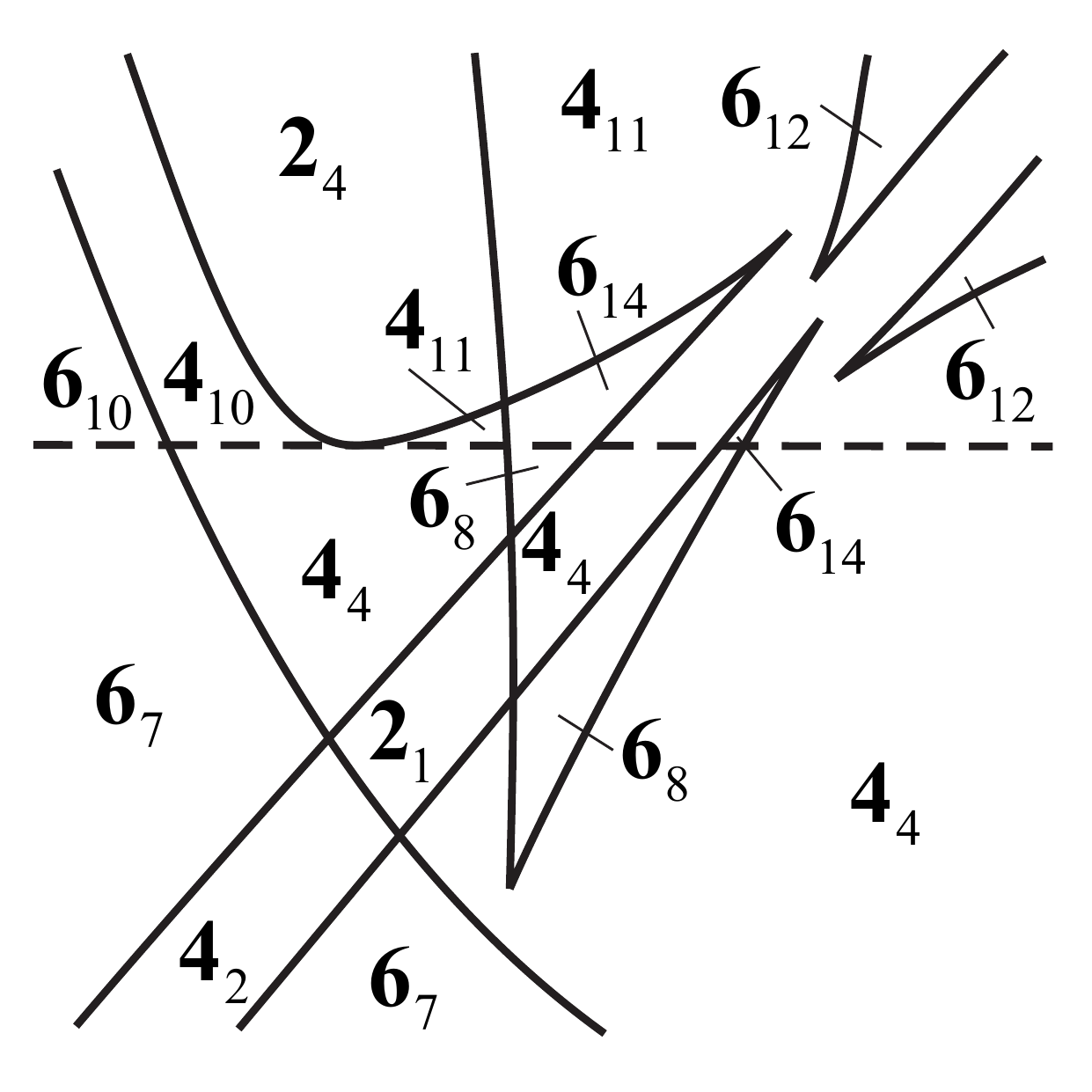}&
71)\includegraphics[width=4cm]{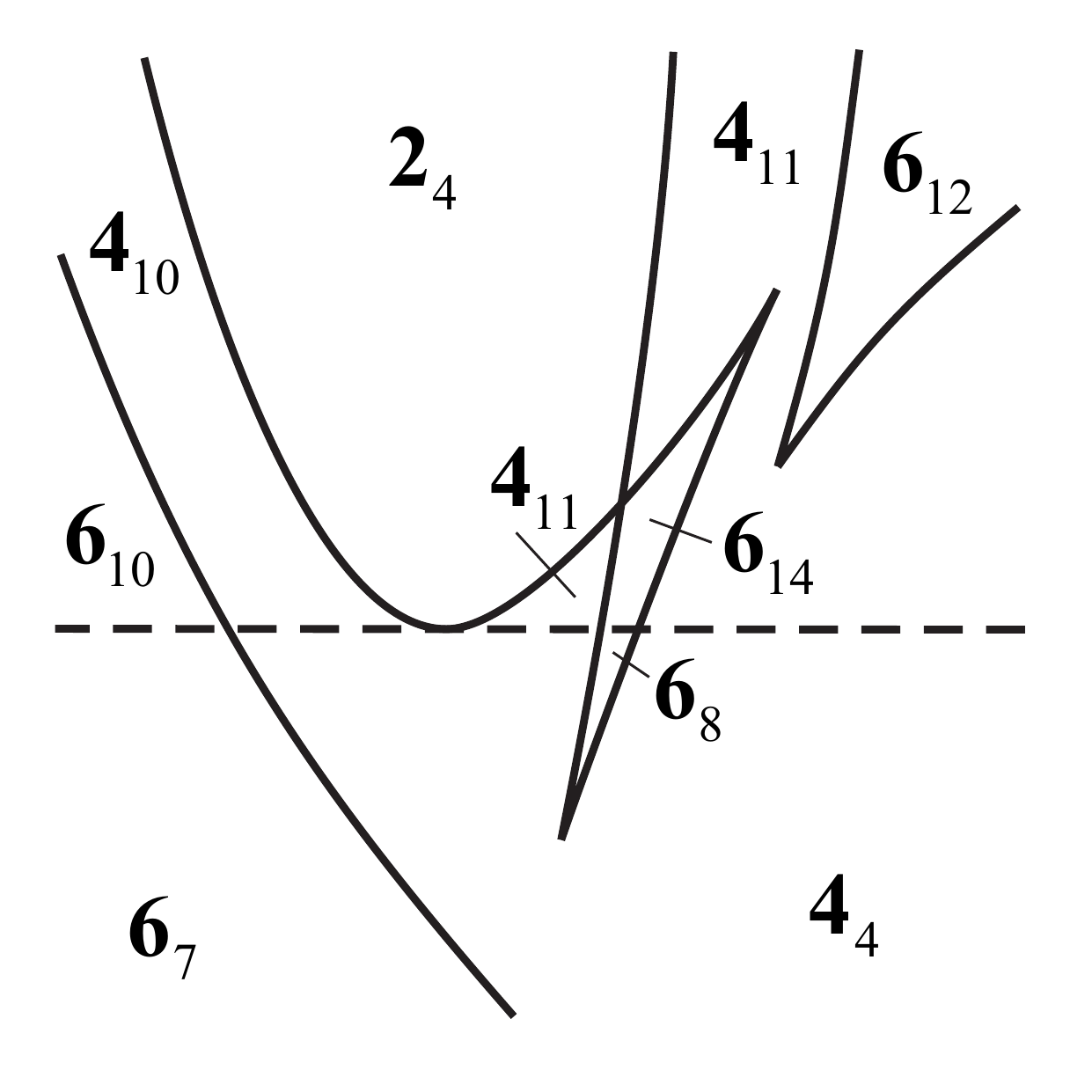}&
72)\includegraphics[width=4cm]{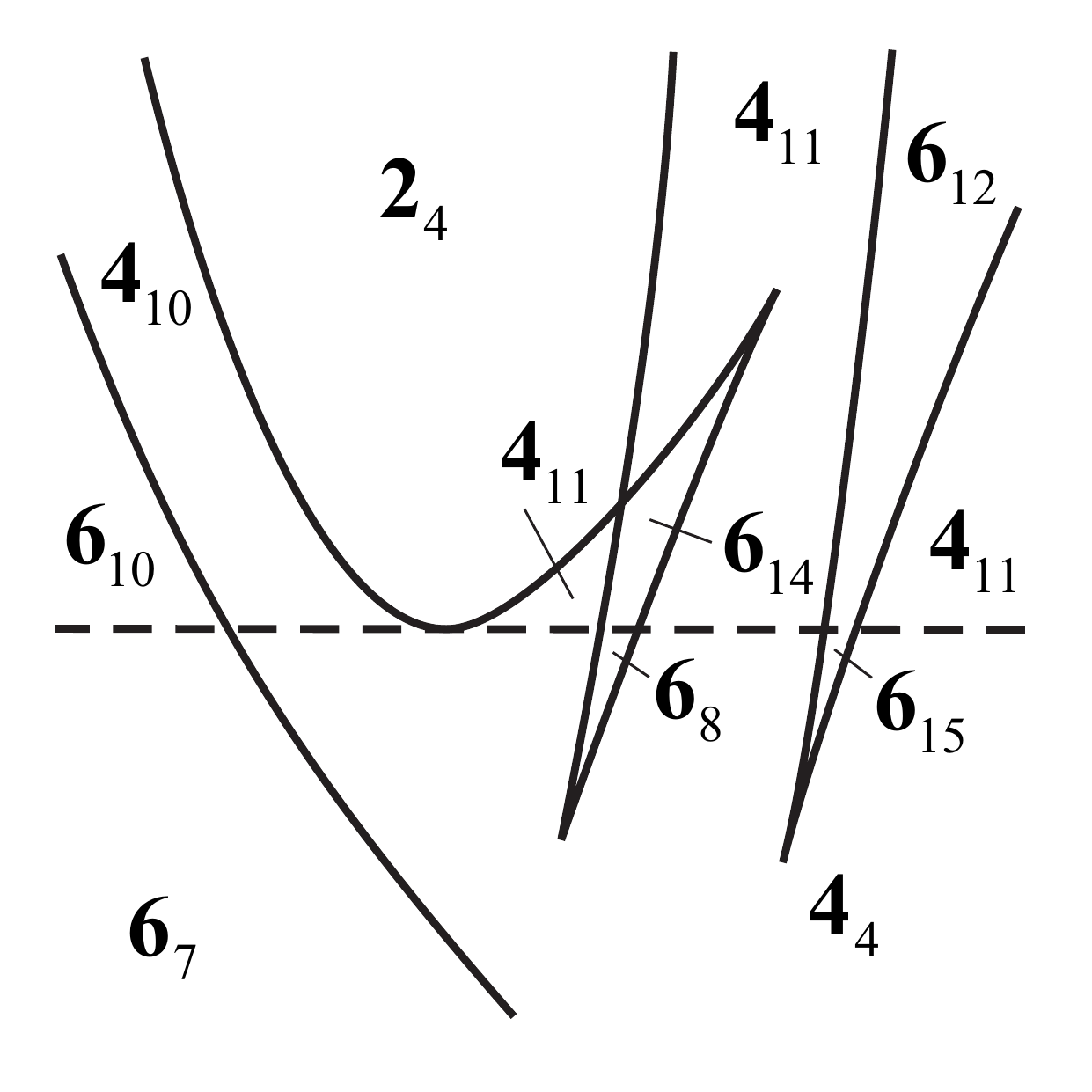}\\
73)\includegraphics[width=4cm]{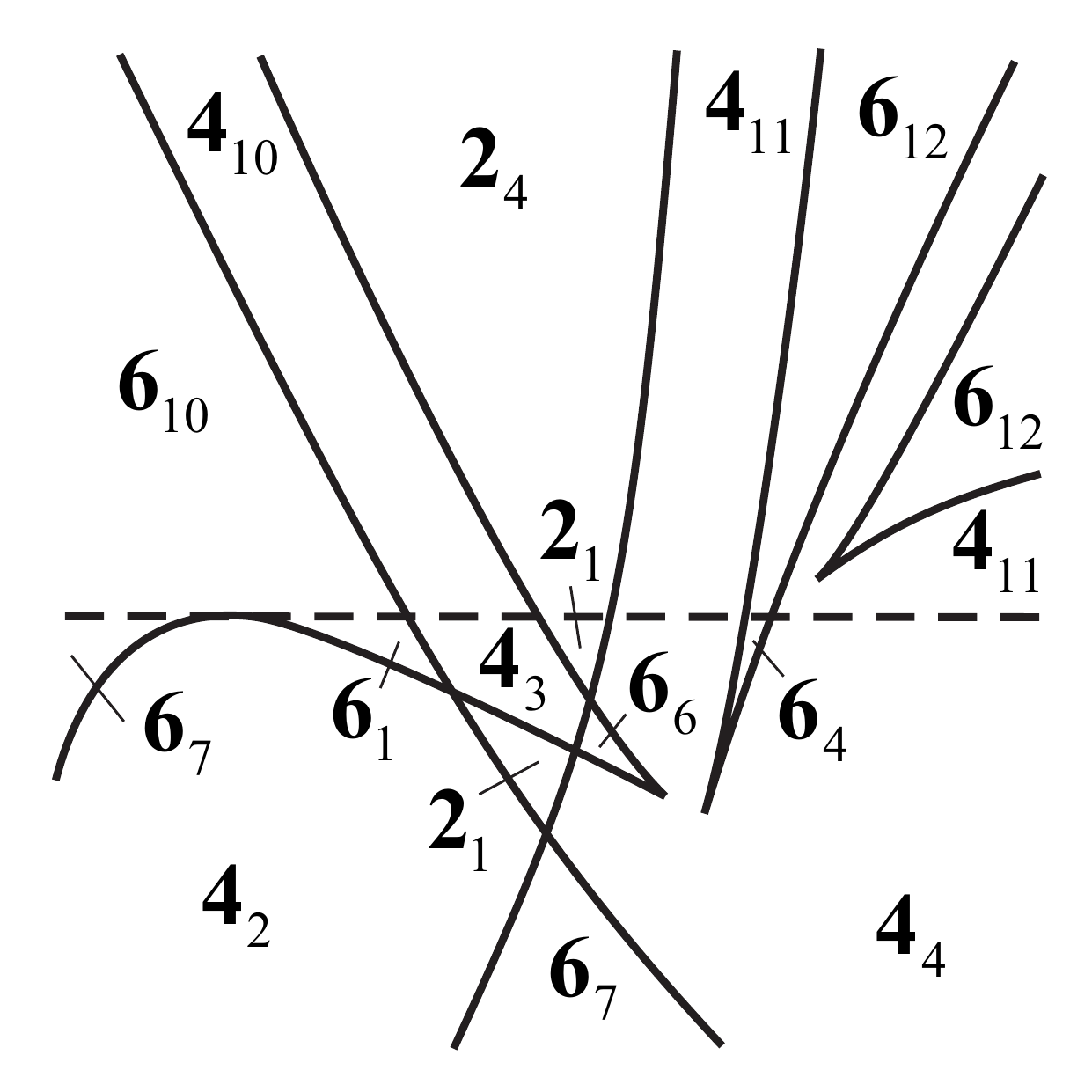}&
74)\includegraphics[width=4cm]{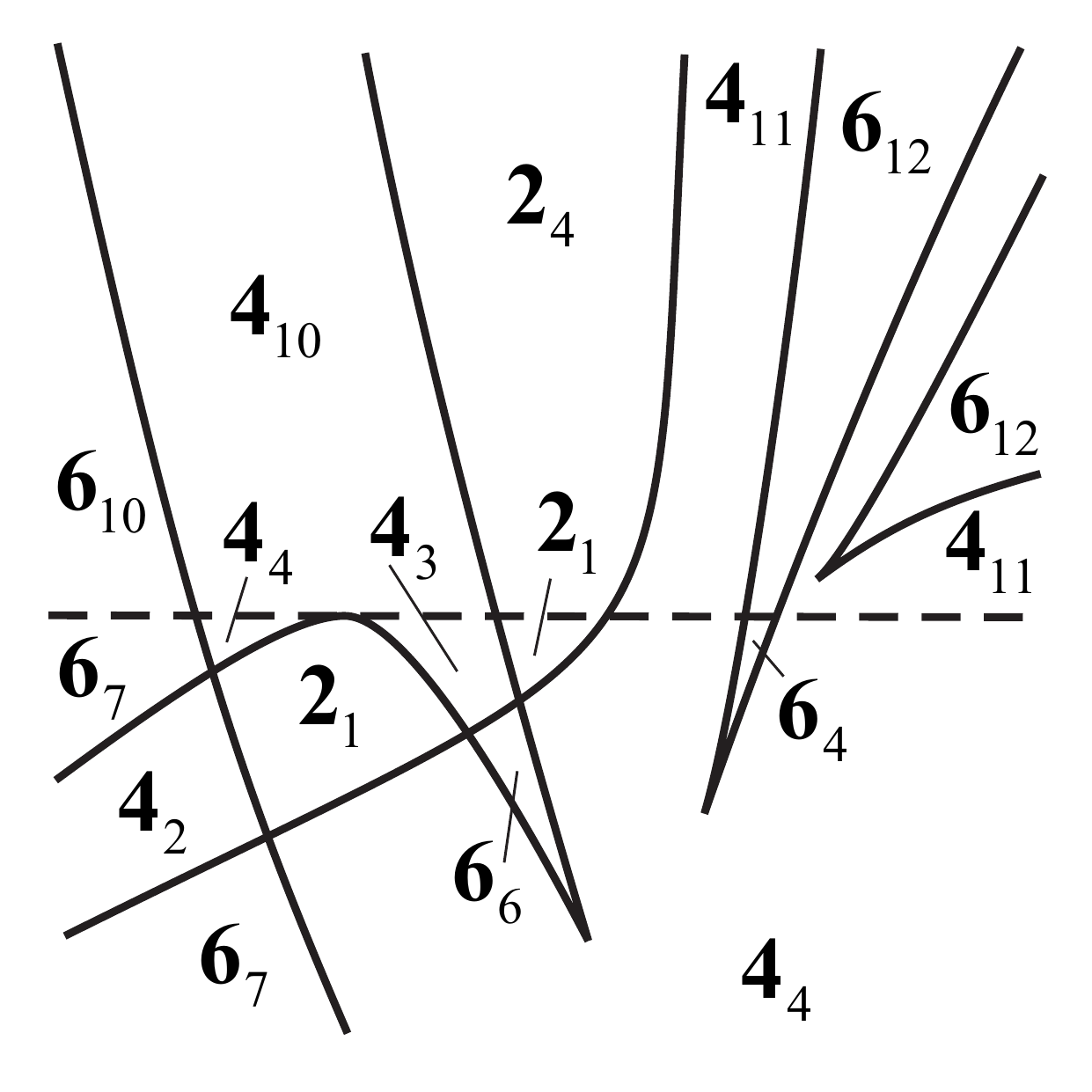}&
75)\includegraphics[width=4cm]{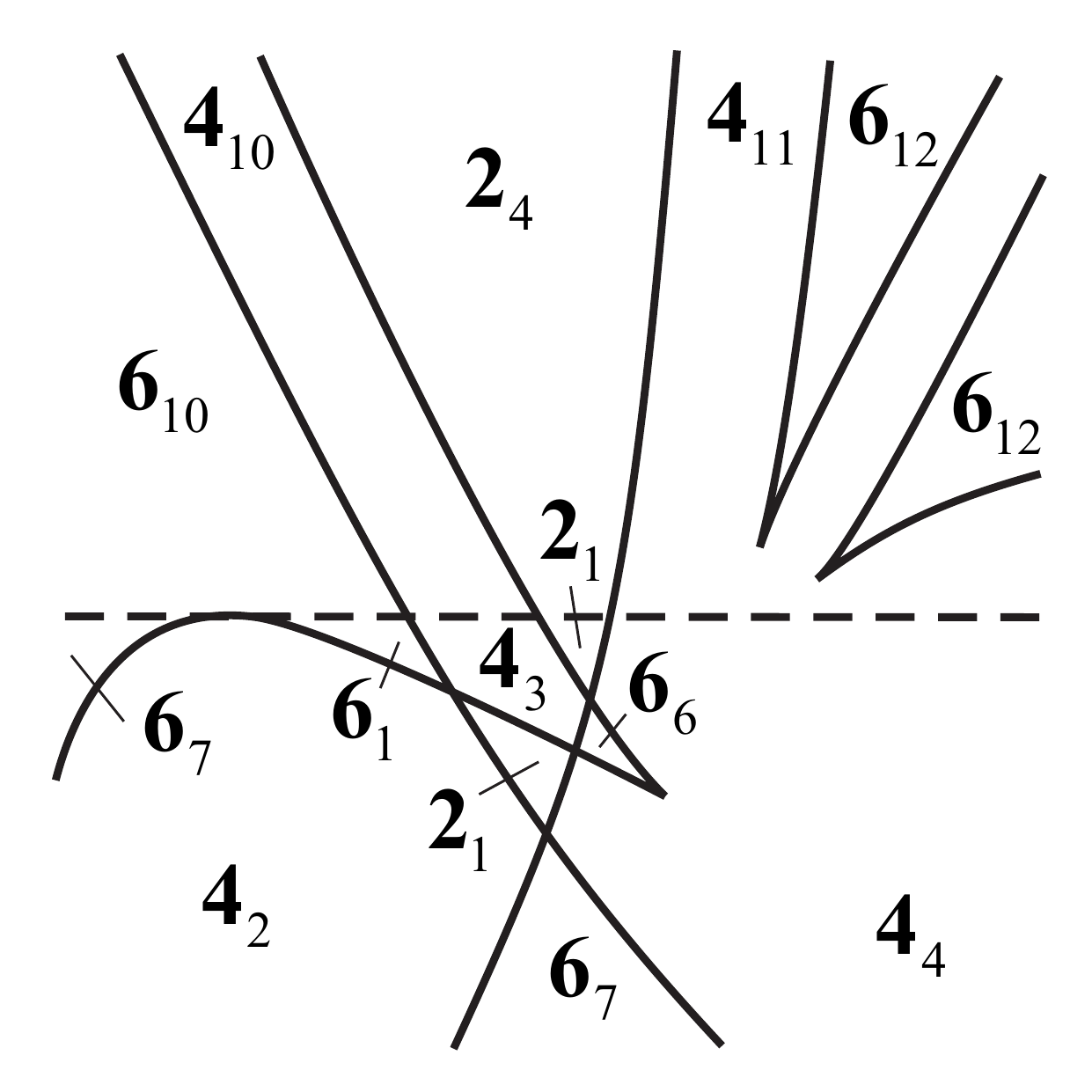}
\end{tabular}
\caption{The skeletons $\Sigma_{S_4,t_1,q_5},S_4\neq0$ for domains $61 - 75$.}
\label{zona61-75}
\end{center}
\end{figure}

\begin{figure}
\begin{center}
\begin{tabular}{ccc}
76)\includegraphics[width=4cm]{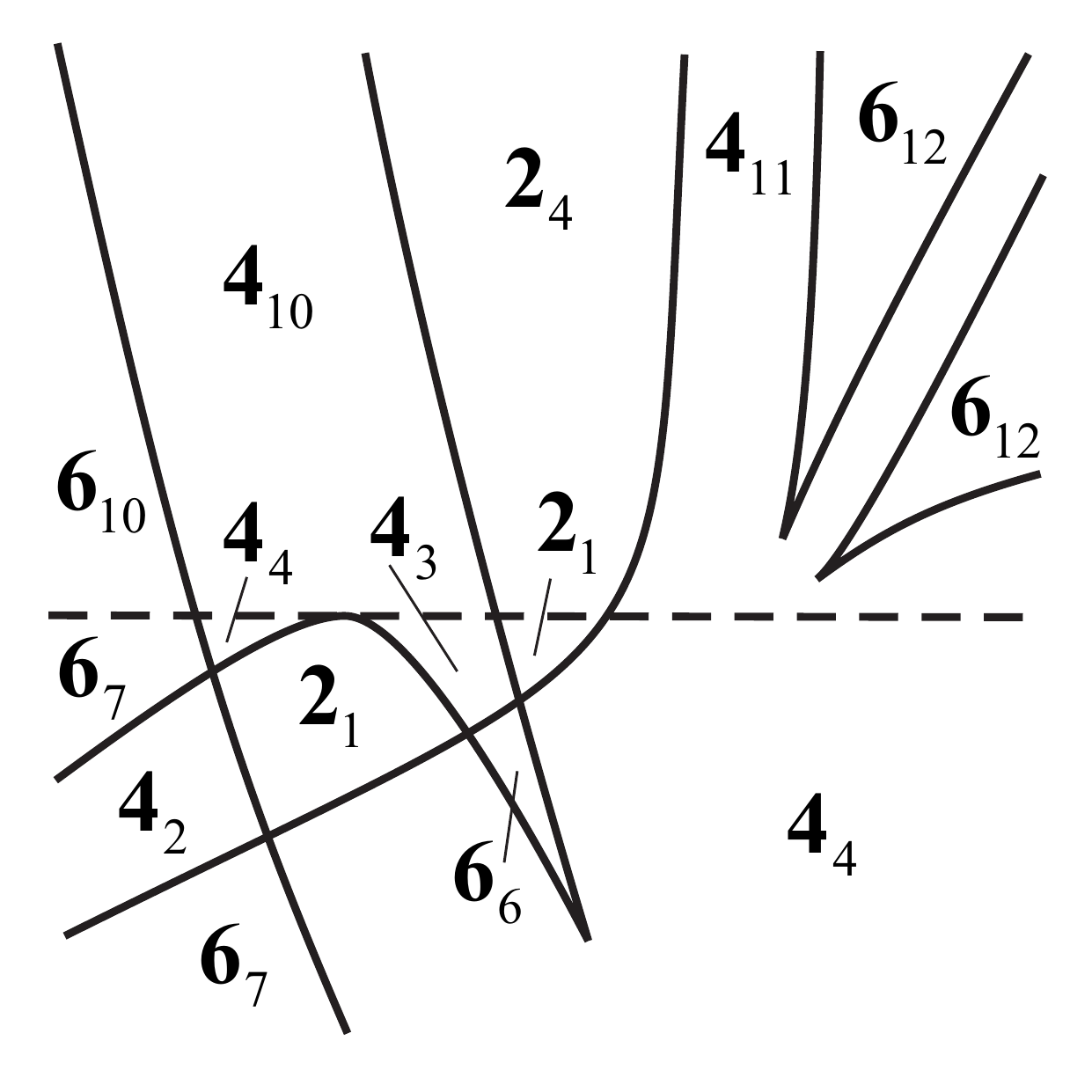}&
77)\includegraphics[width=4cm]{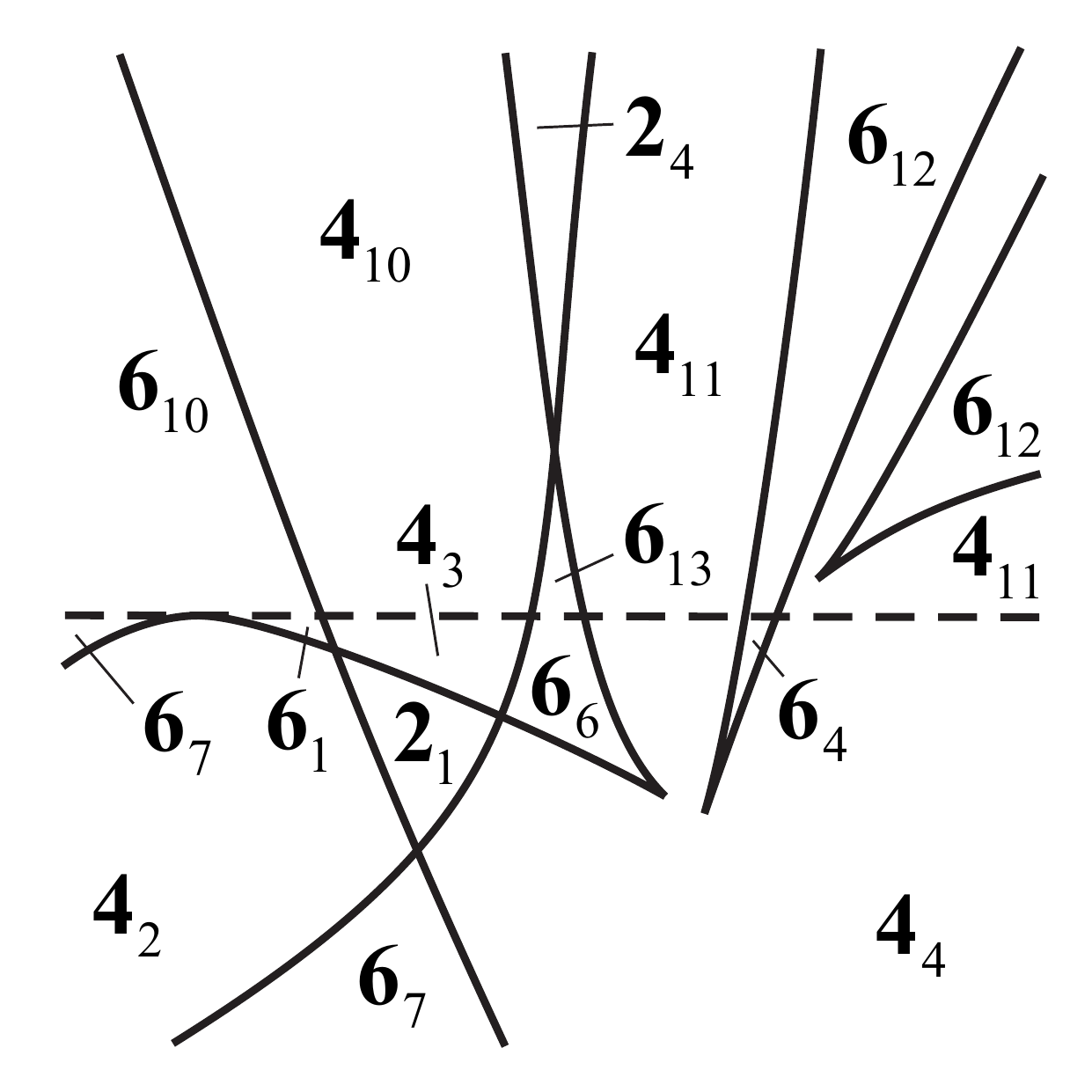}&
78)\includegraphics[width=4cm]{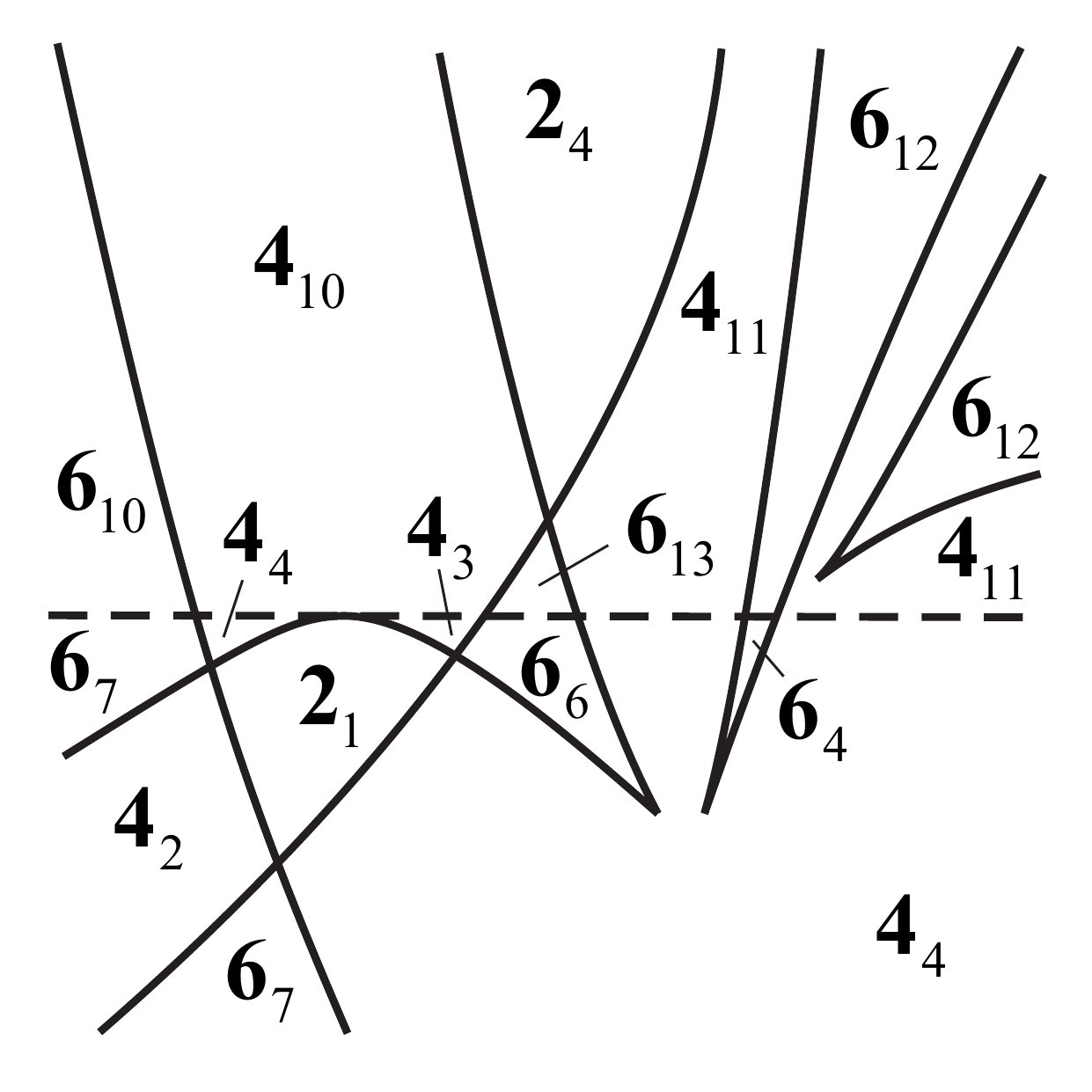}\\
79)\includegraphics[width=4cm]{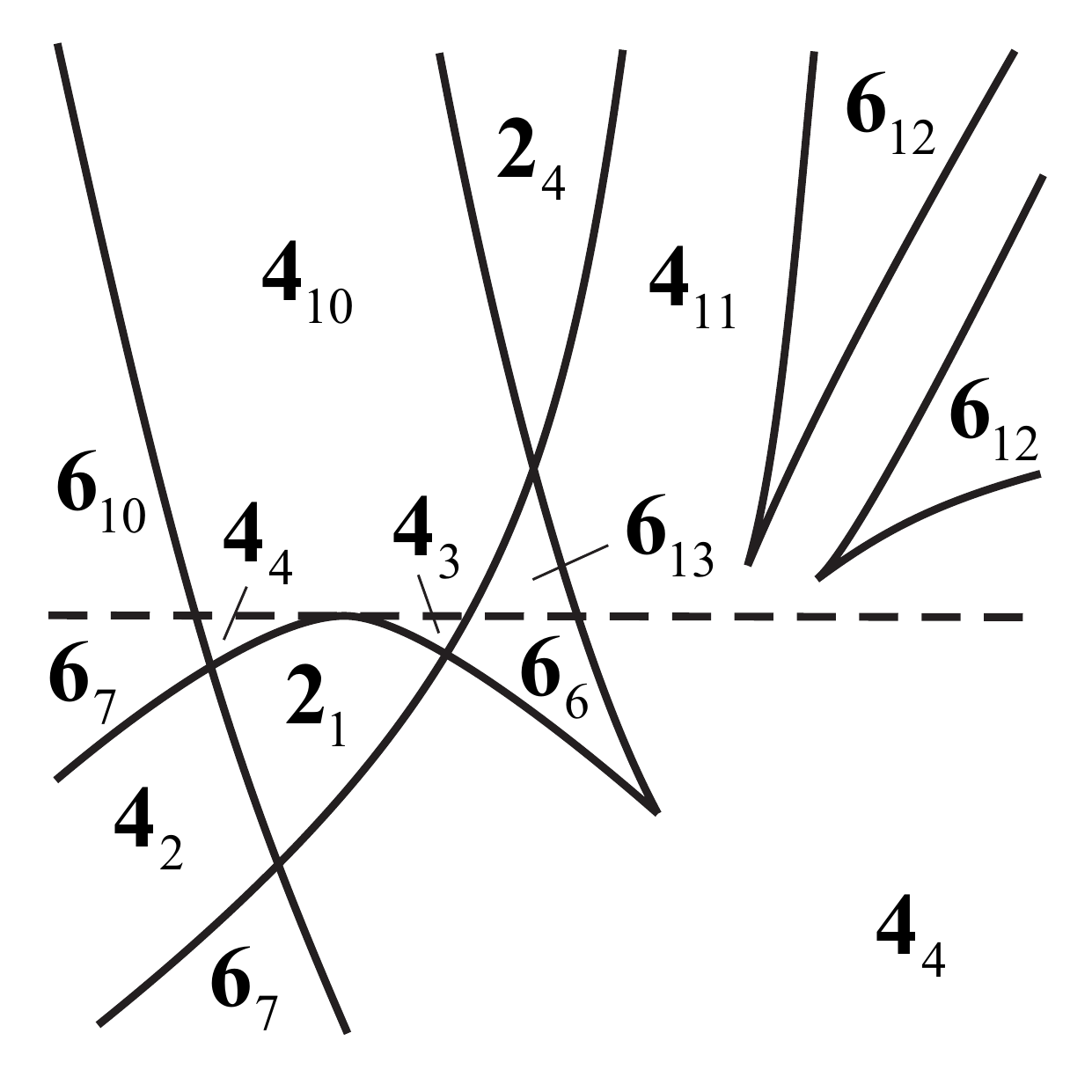}&
80)\includegraphics[width=4cm]{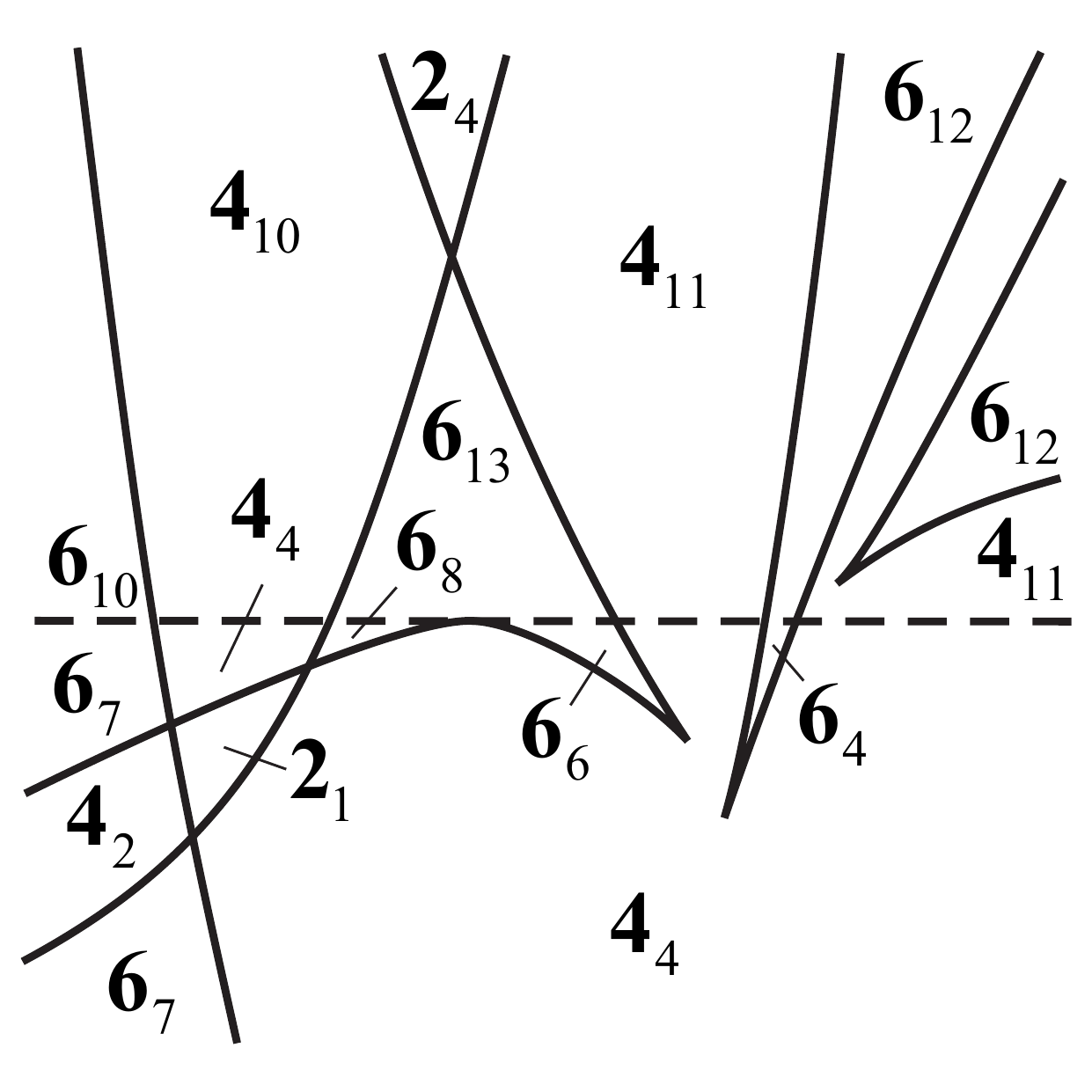}&
81)\includegraphics[width=4cm]{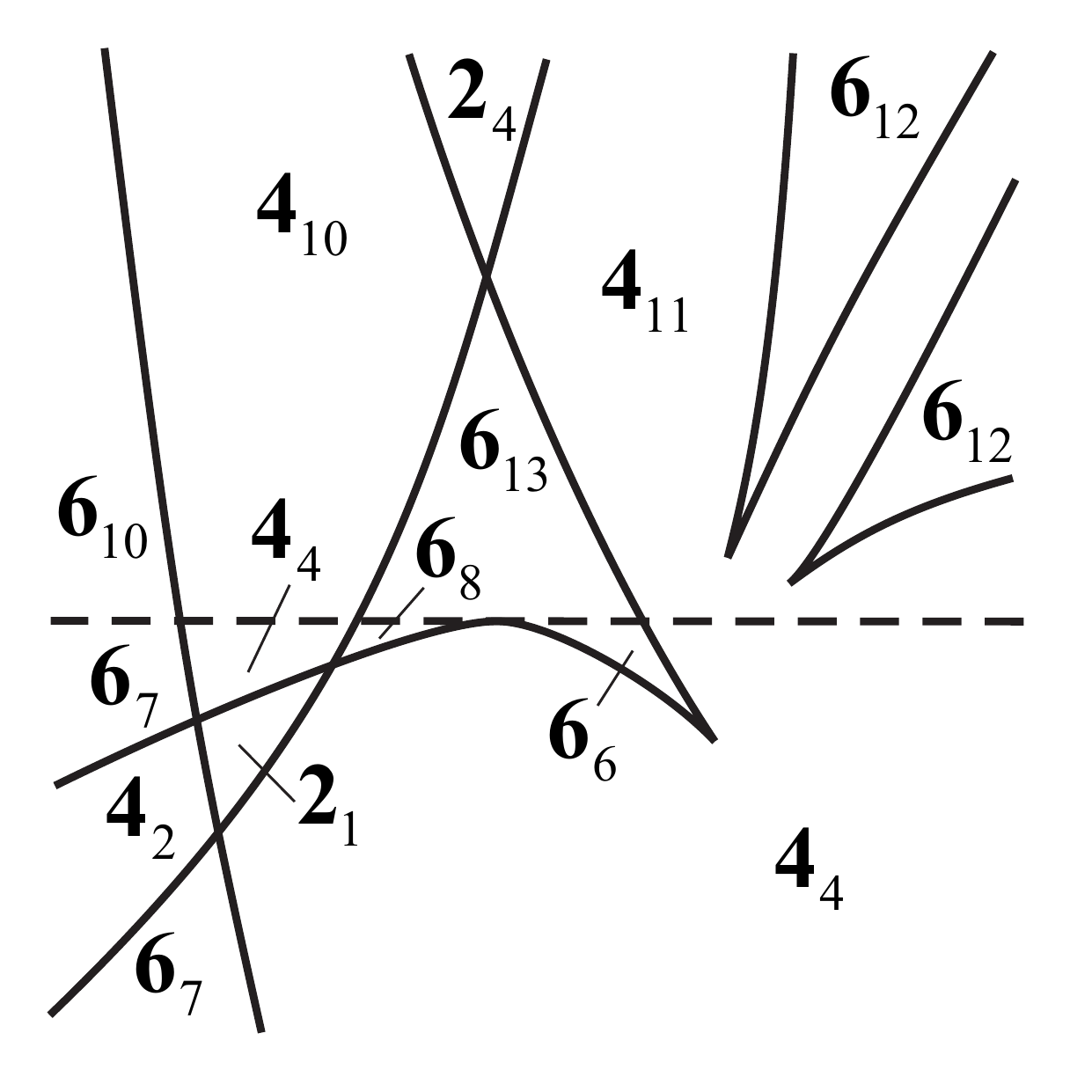}\\
82)\includegraphics[width=4cm]{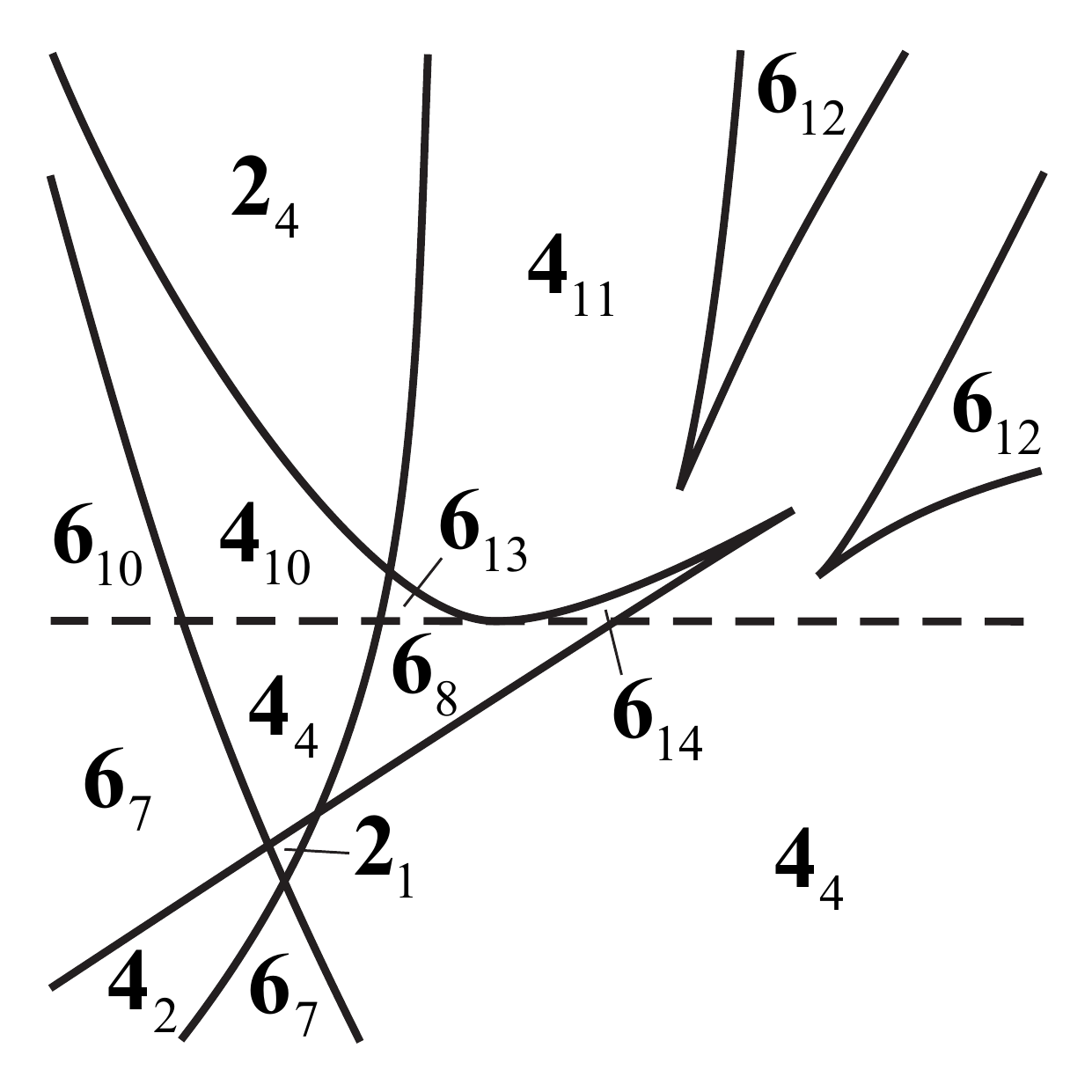}&
83)\includegraphics[width=4cm]{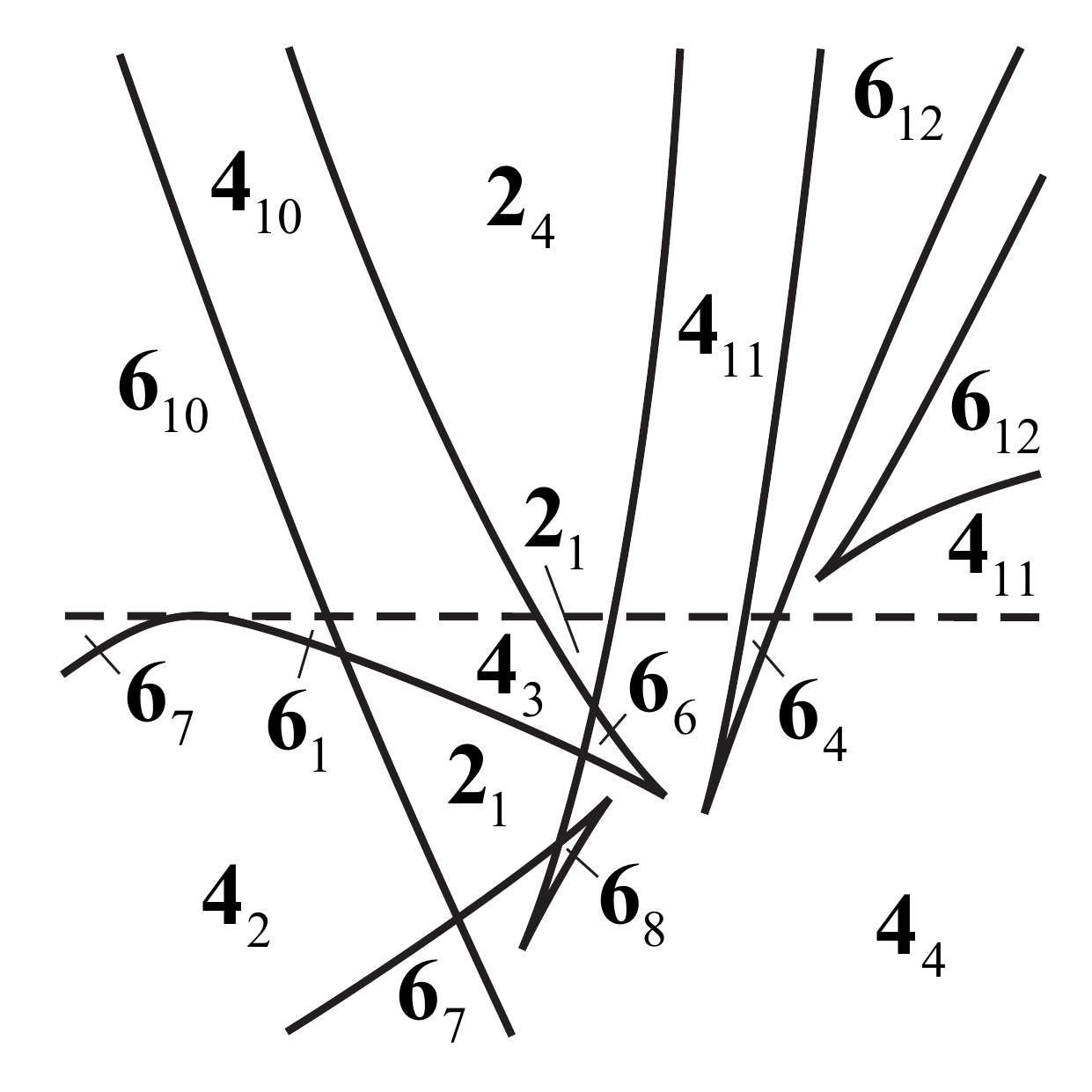}&
84)\includegraphics[width=4cm]{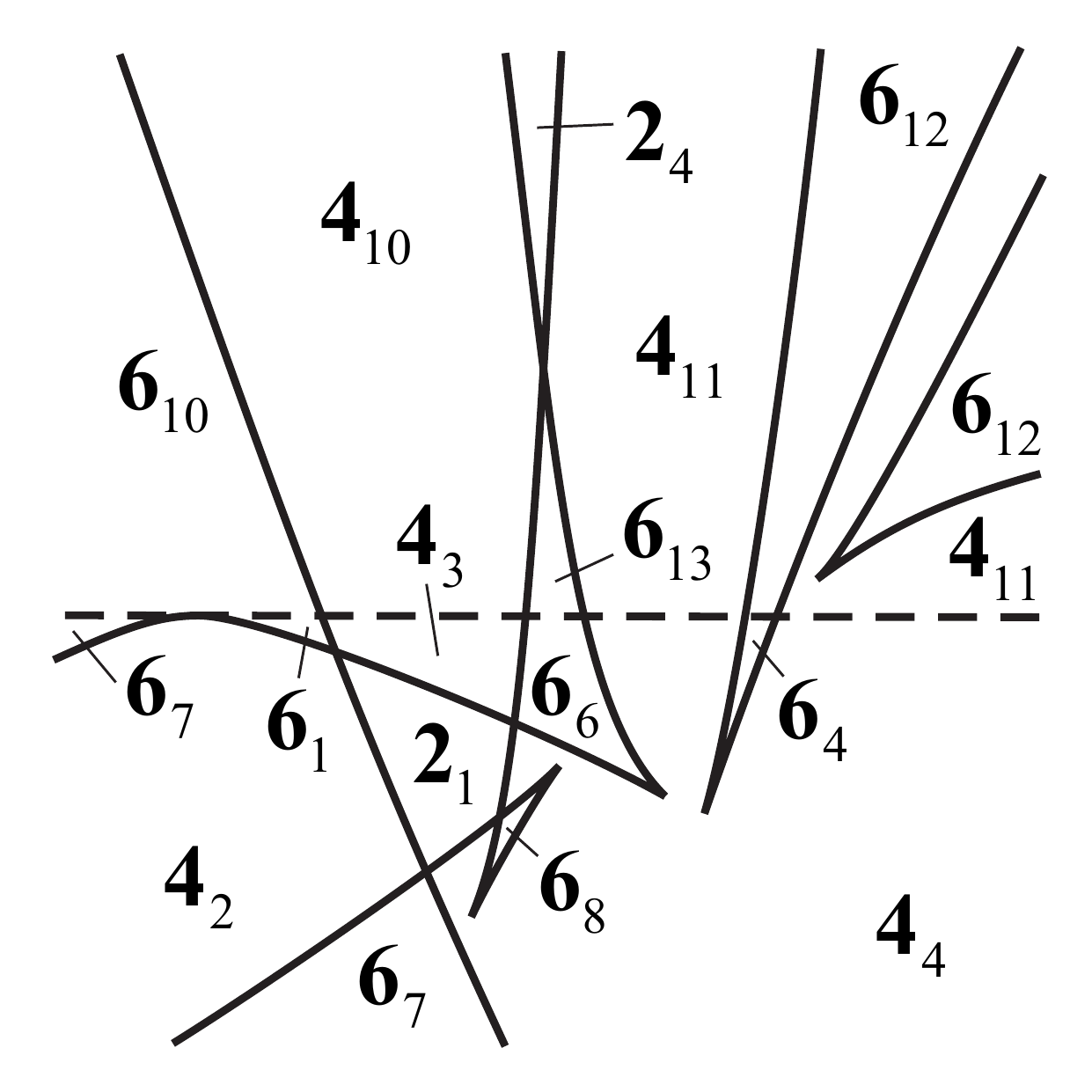}\\
85)\includegraphics[width=4cm]{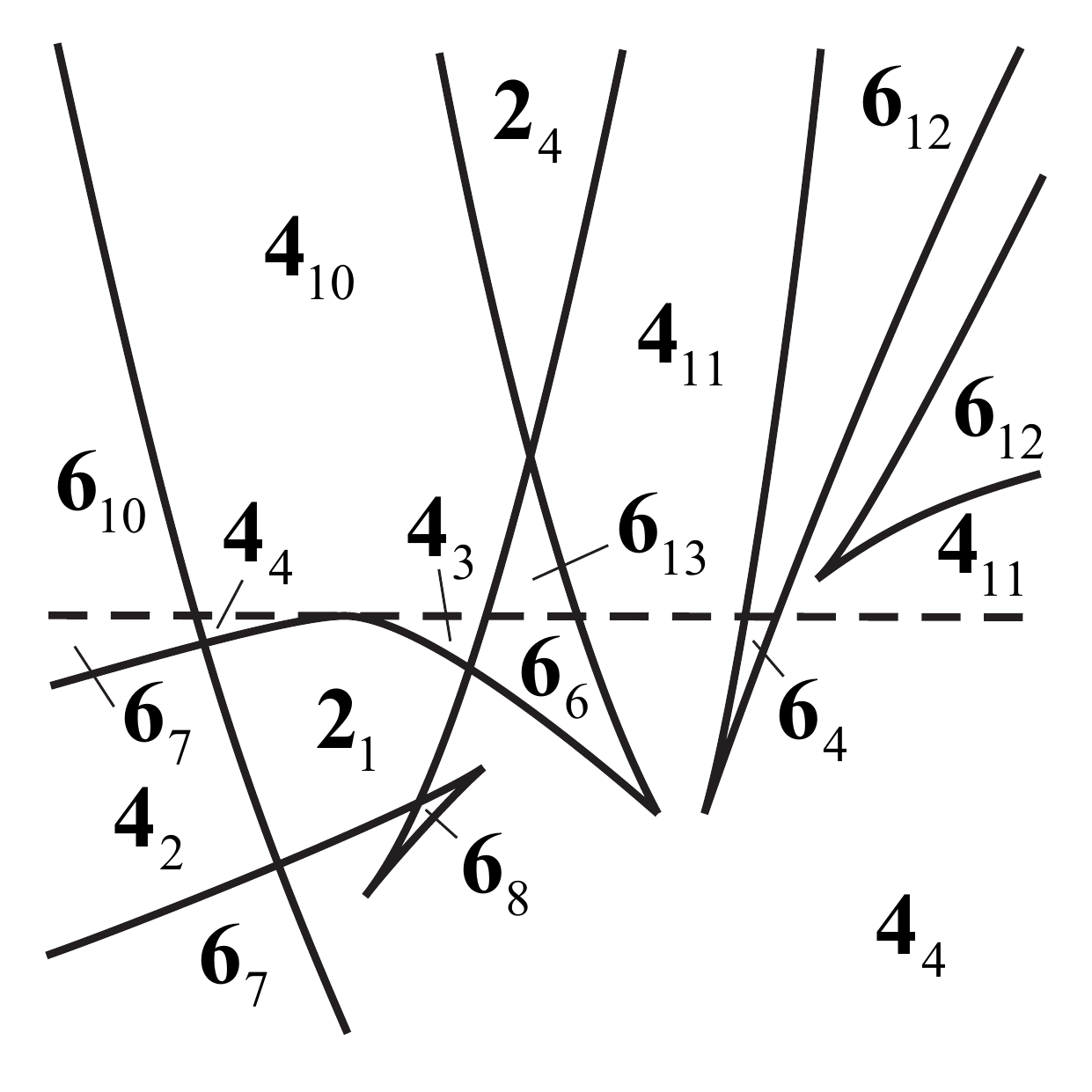}&
86)\includegraphics[width=4cm]{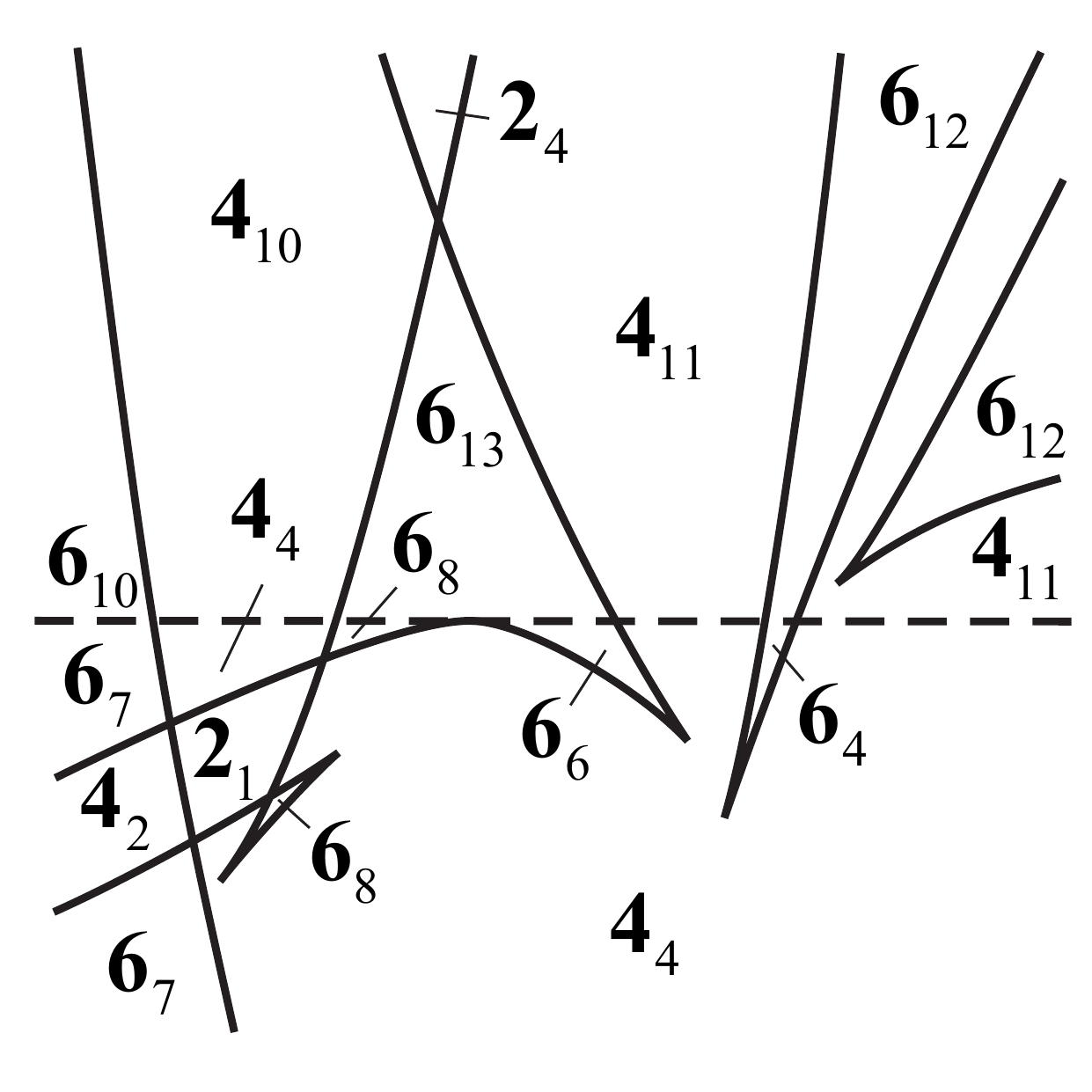}&
87)\includegraphics[width=4cm]{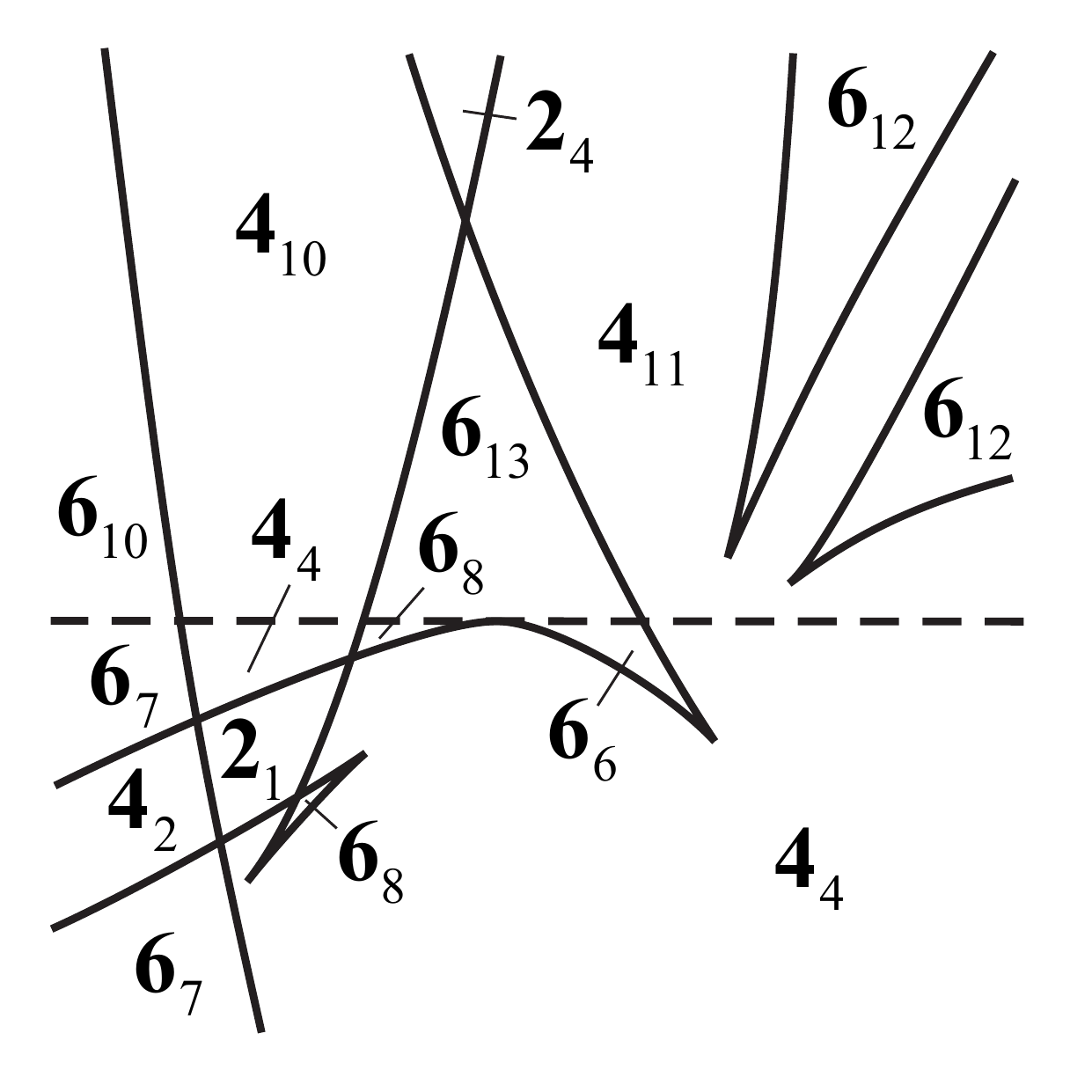}\\
88)\includegraphics[width=4cm]{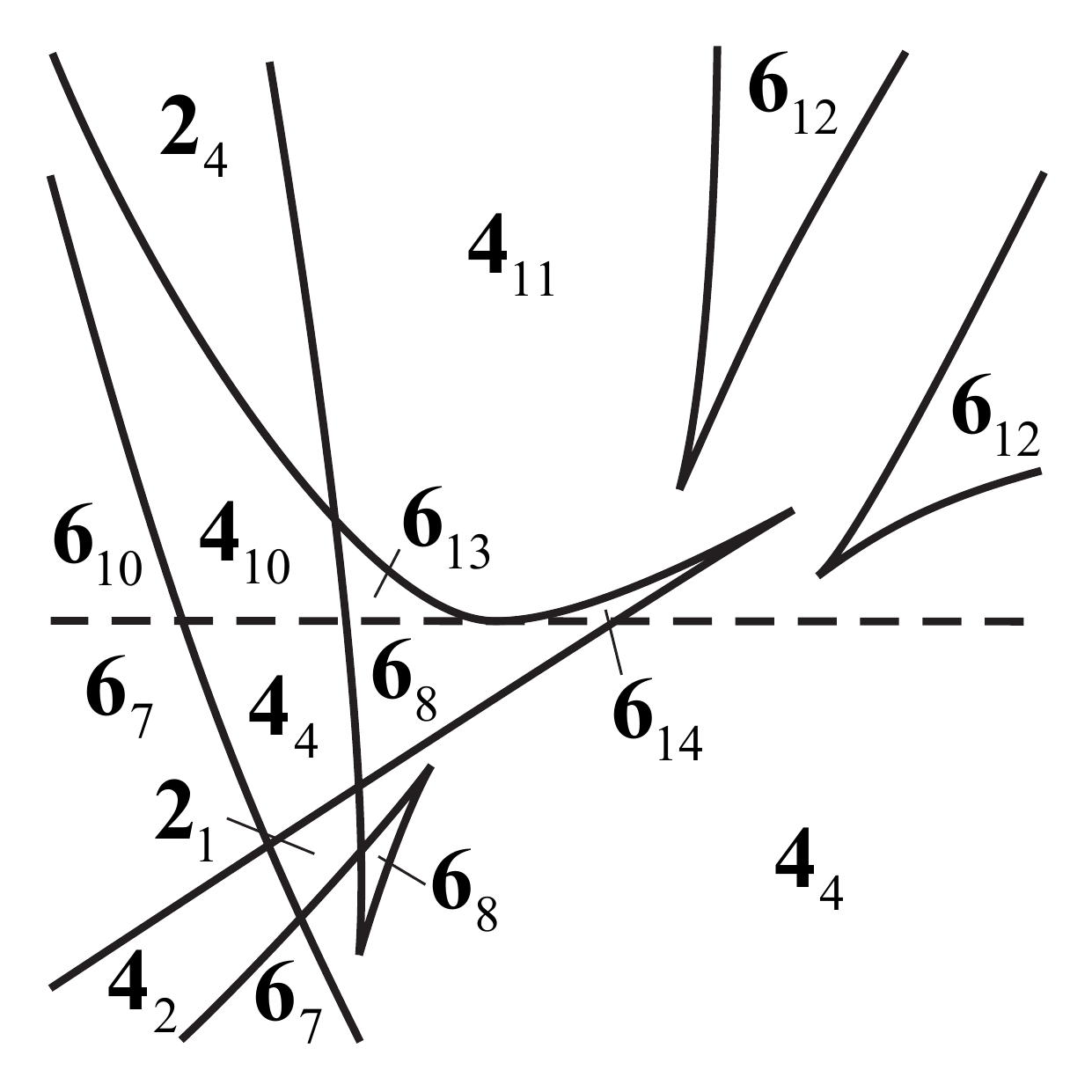}&
89)\includegraphics[width=4cm]{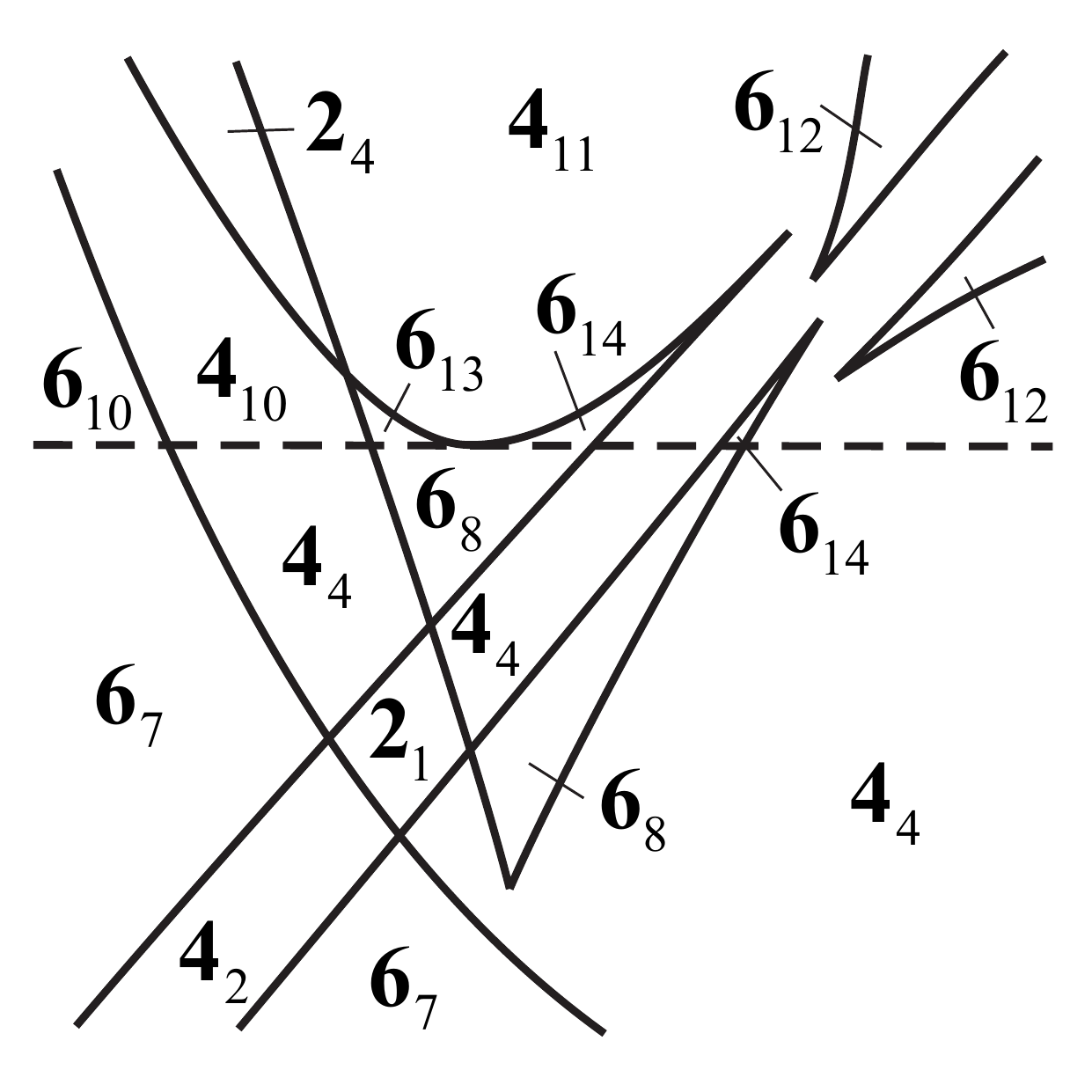}&
90)\includegraphics[width=4cm]{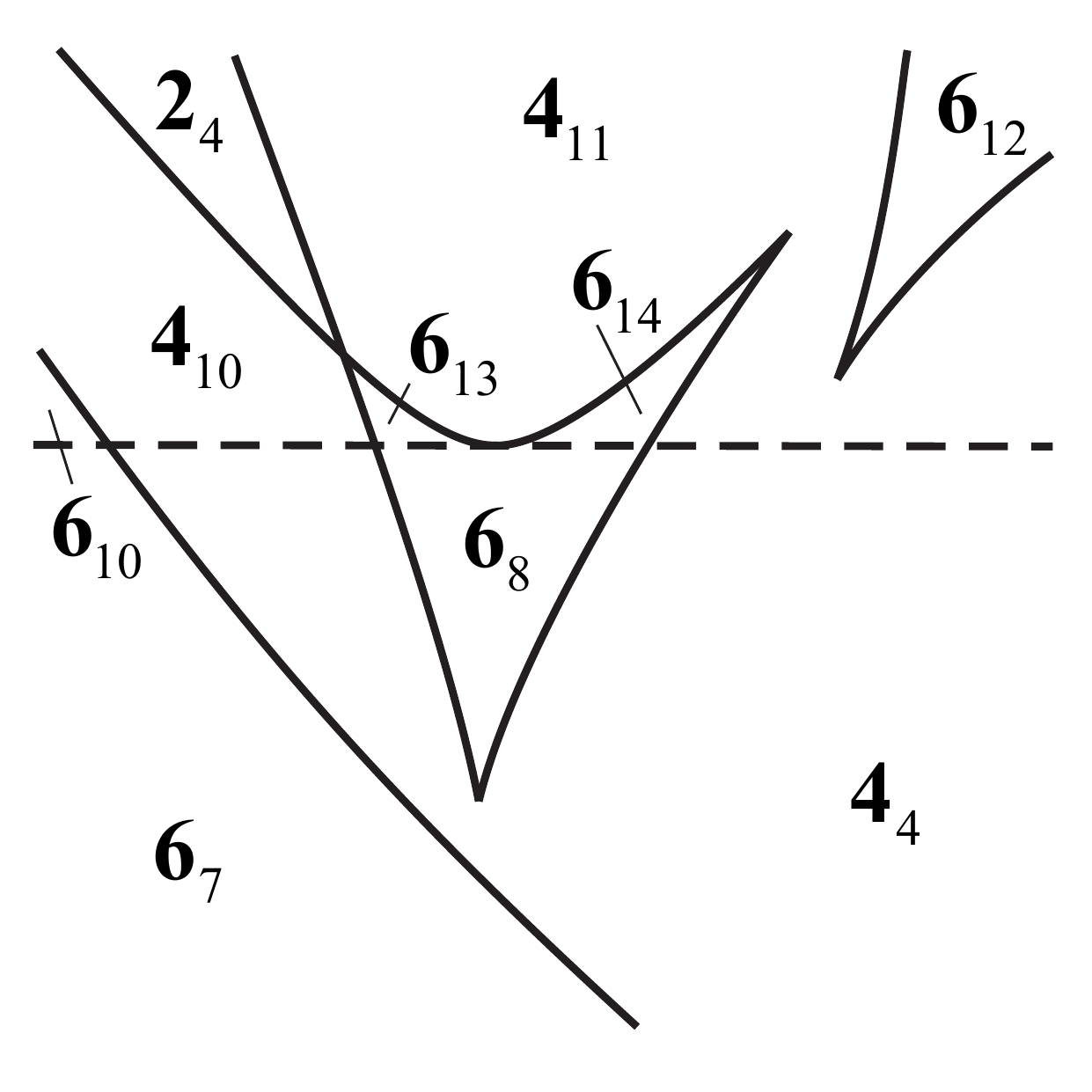}
\end{tabular}
\caption{The skeletons $\Sigma_{S_4,t_1,q_5},S_4\neq0$ for domains $76 - 90$.}
\label{zona76-90}
\end{center}
\end{figure}

\begin{figure}
\begin{center}
\begin{tabular}{ccc}
91)\includegraphics[width=4cm]{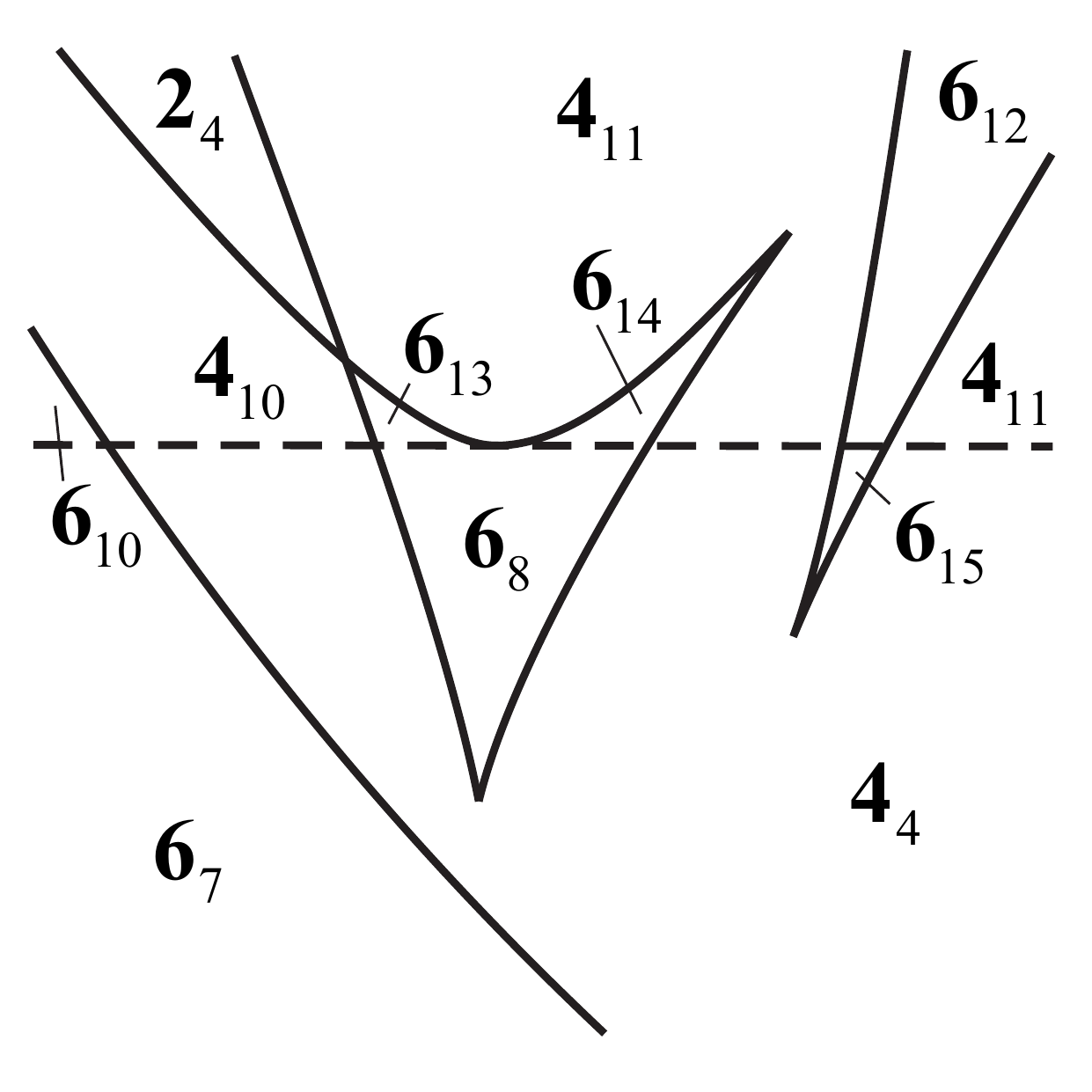}&
92)\includegraphics[width=4cm]{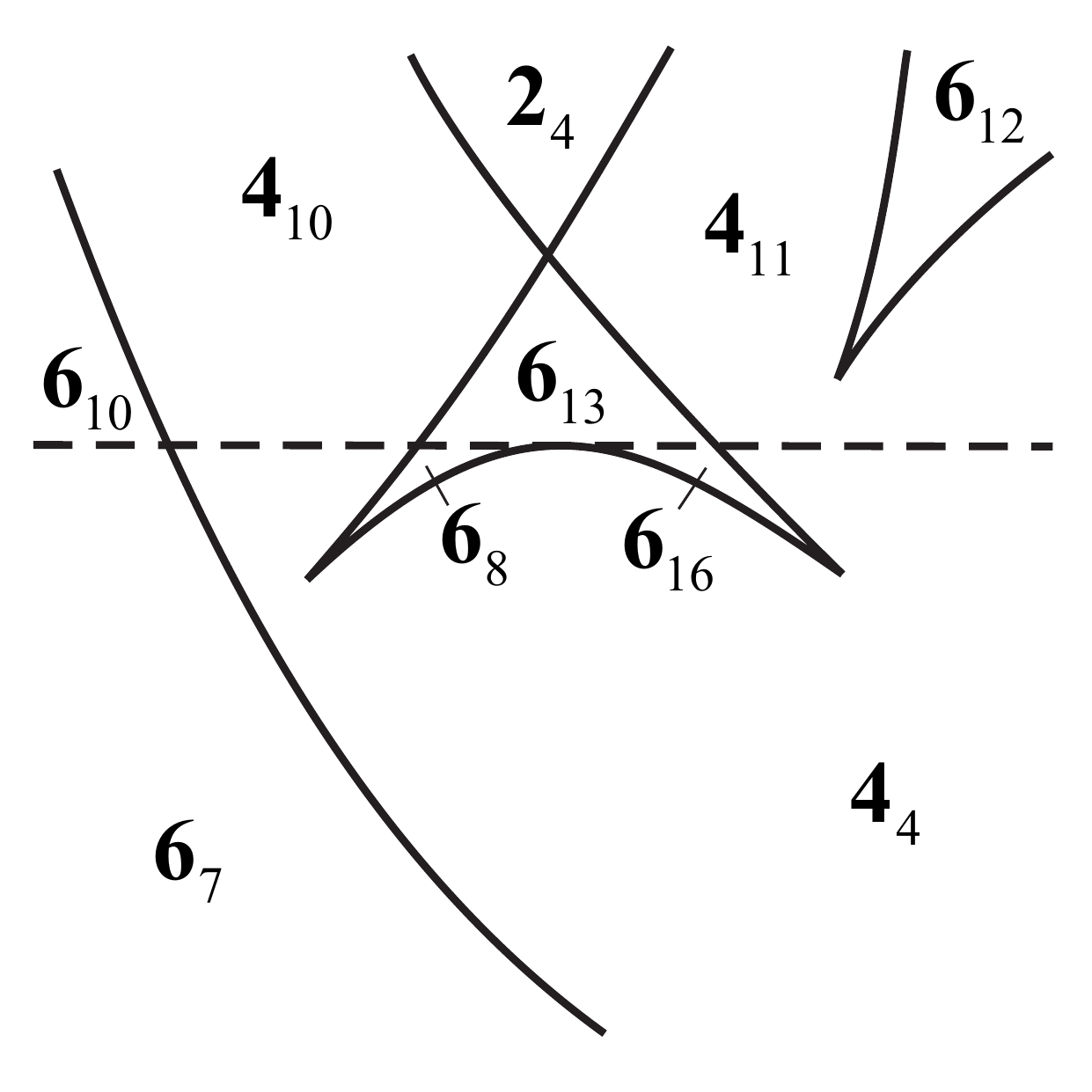}&
93)\includegraphics[width=4cm]{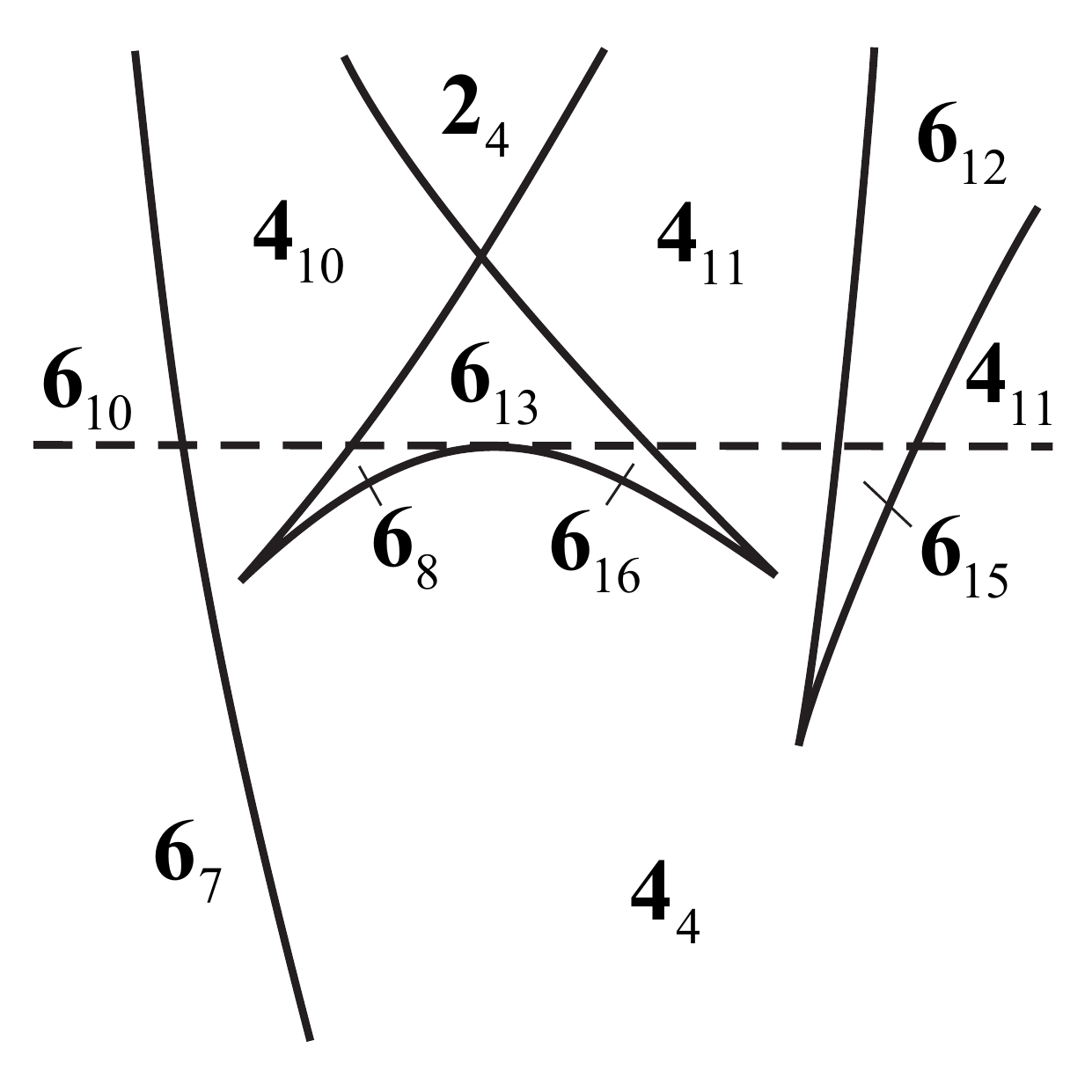}
\end{tabular}
\caption{The skeletons $\Sigma_{S_4,t_1,q_5},S_4\neq0$ for domains $91 - 93$.}
\label{zona91-93}
\end{center}
\end{figure}

\begin{figure}
\begin{center}
\begin{tabular}{ccc}
\includegraphics[width=4cm]{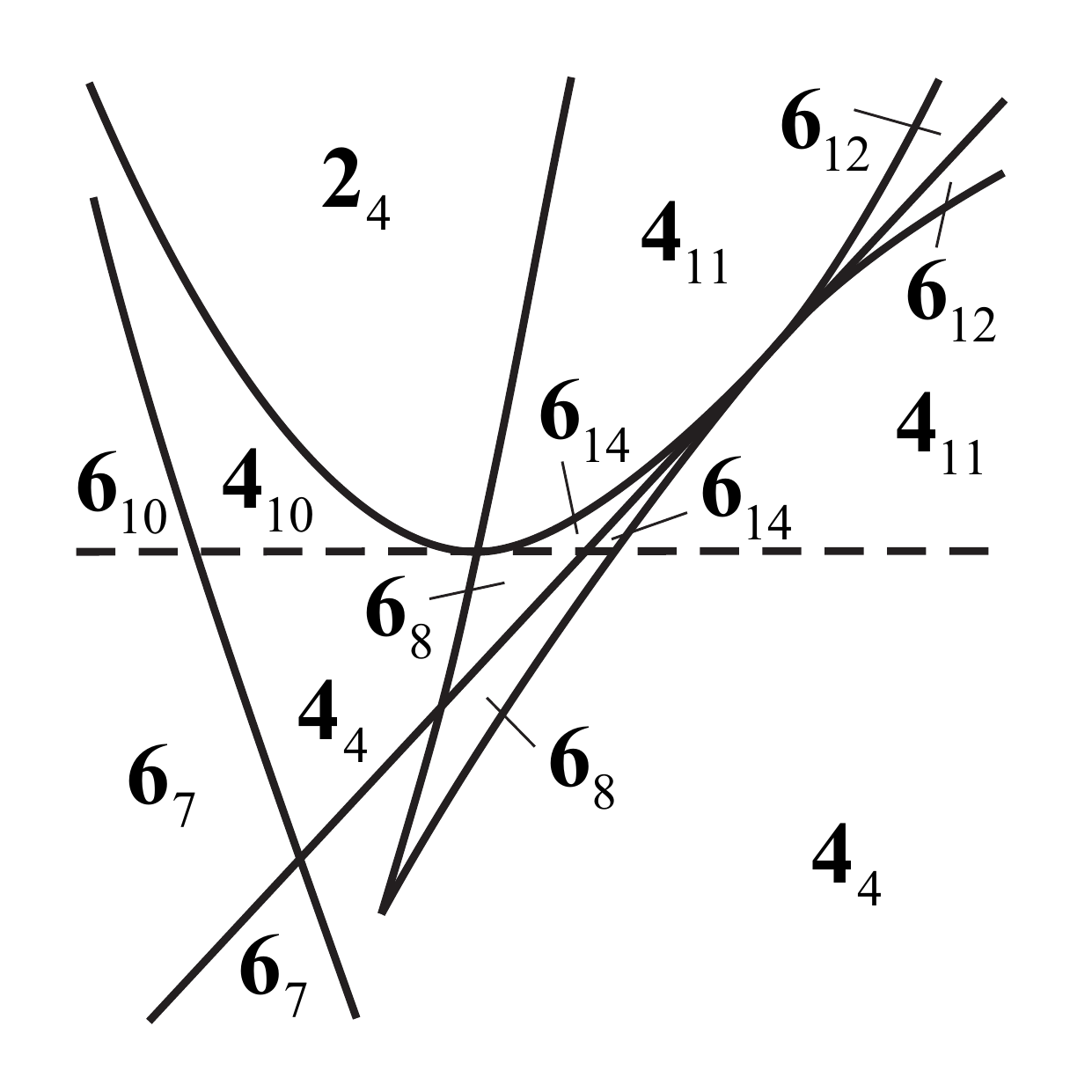}&\includegraphics[width=4cm]{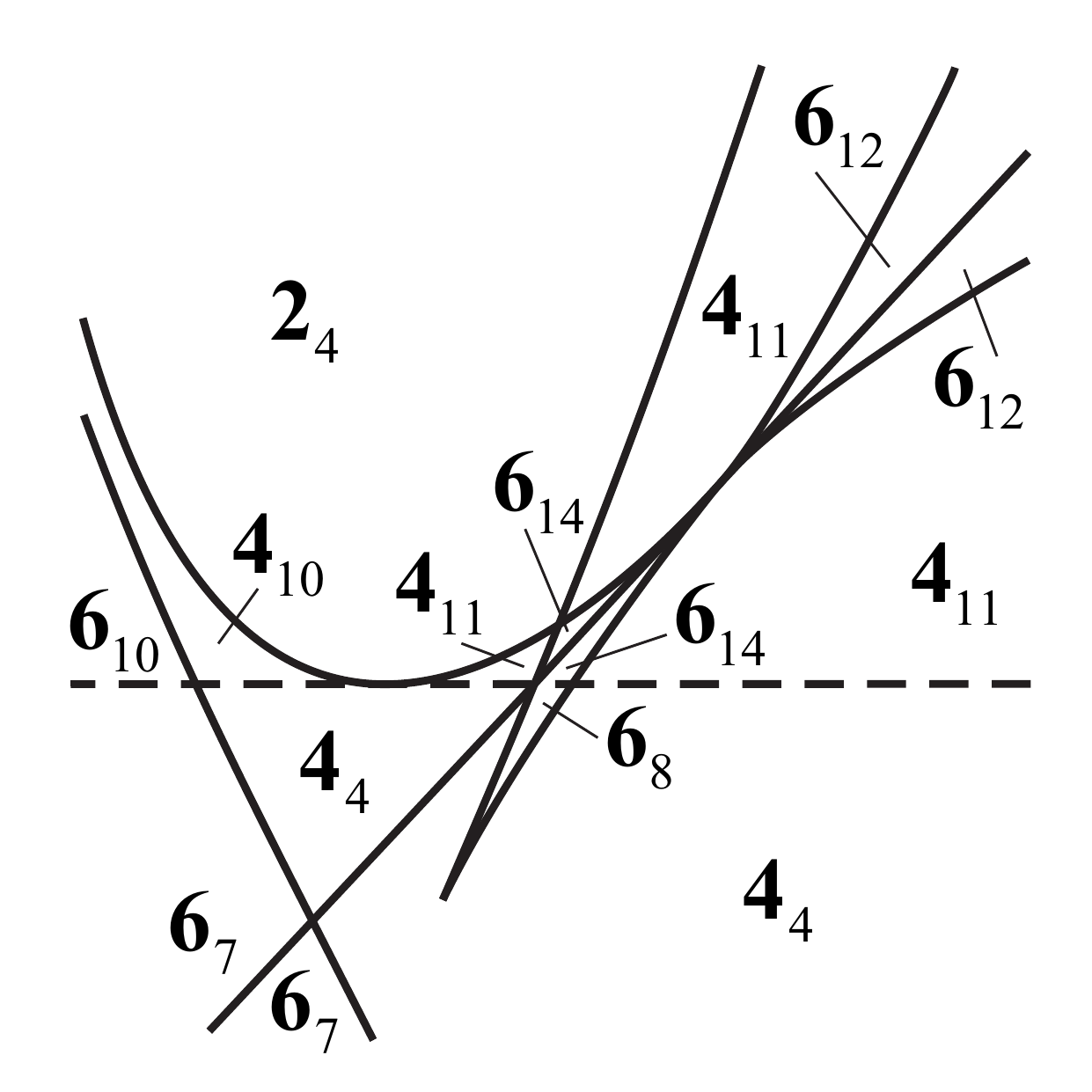}&
\includegraphics[width=4cm]{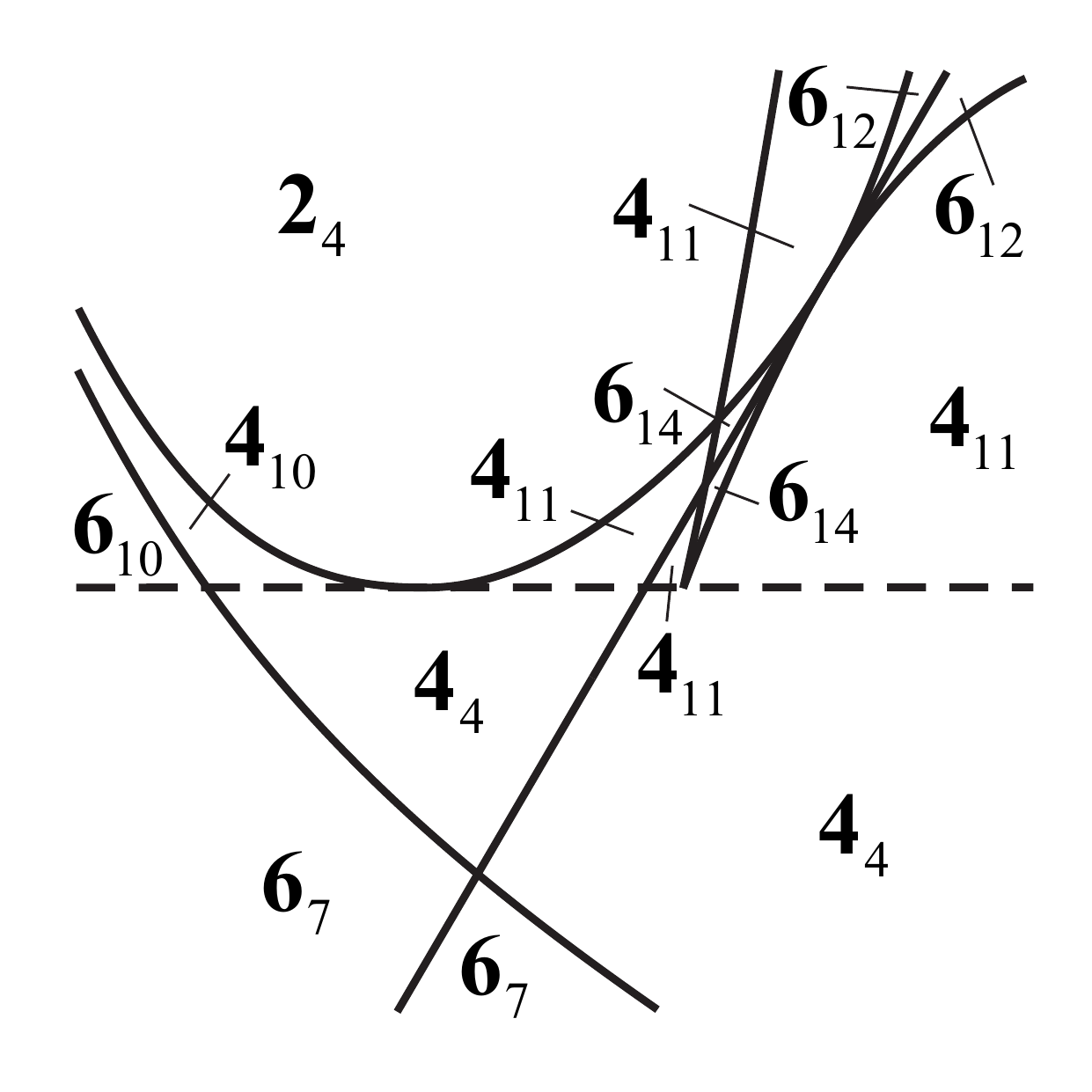}\\
$(t_1,q_5)\equiv(\ref{spoint8})$&$(t_1,q_5)\equiv(\ref{spoint5})$&
$(t_1,q_5)\equiv(\ref{trper-xi5-g1})$\\
\includegraphics[width=4cm]{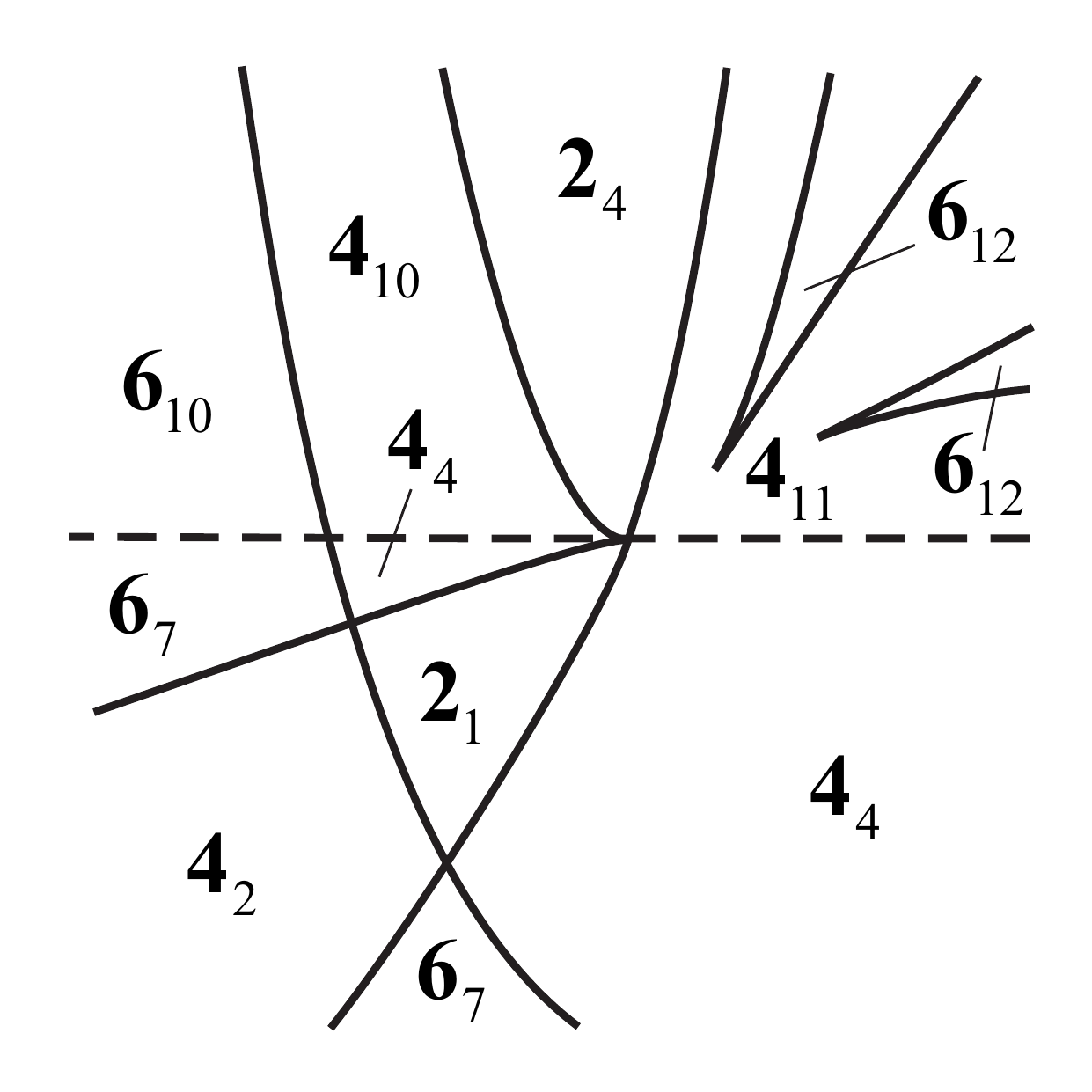}&\includegraphics[width=4cm]{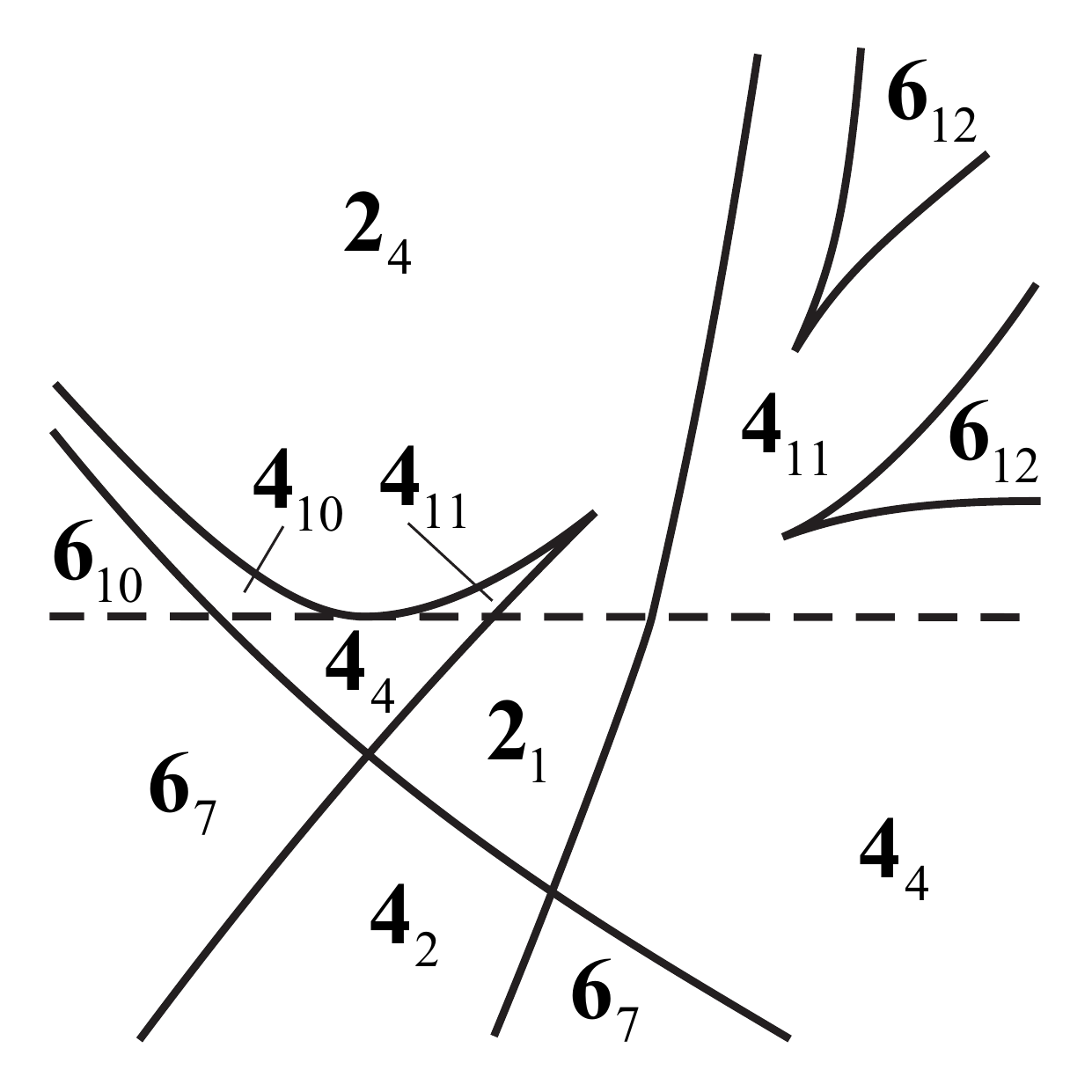}&\includegraphics[width=4cm]{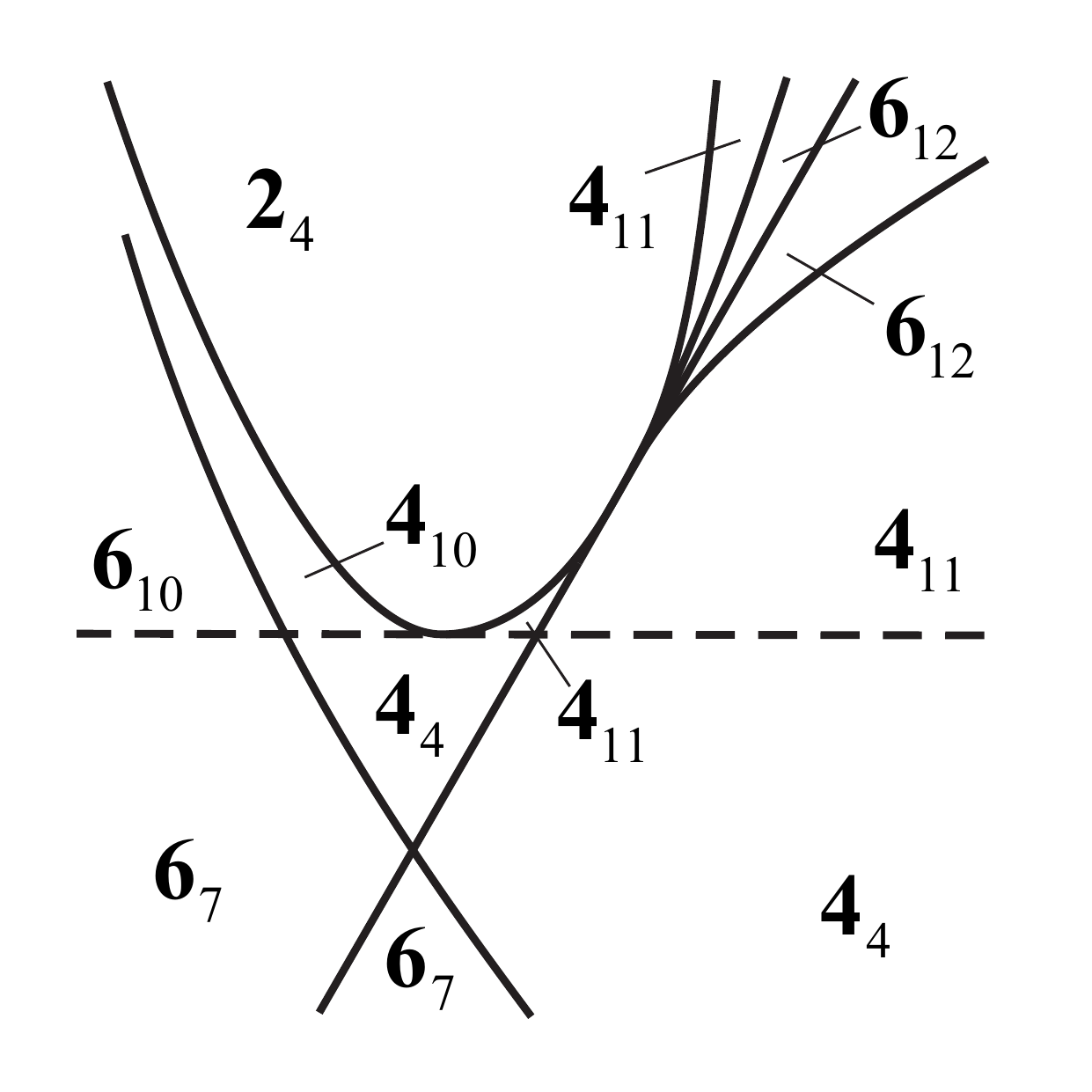}\\
$(t_1,q_5)\equiv(\ref{spoint6})$&$(t_1,q_5)\equiv(\ref{spoint4})$&$(t_1,q_5)\equiv(\ref{vozvrat})$\\
\includegraphics[width=4cm]{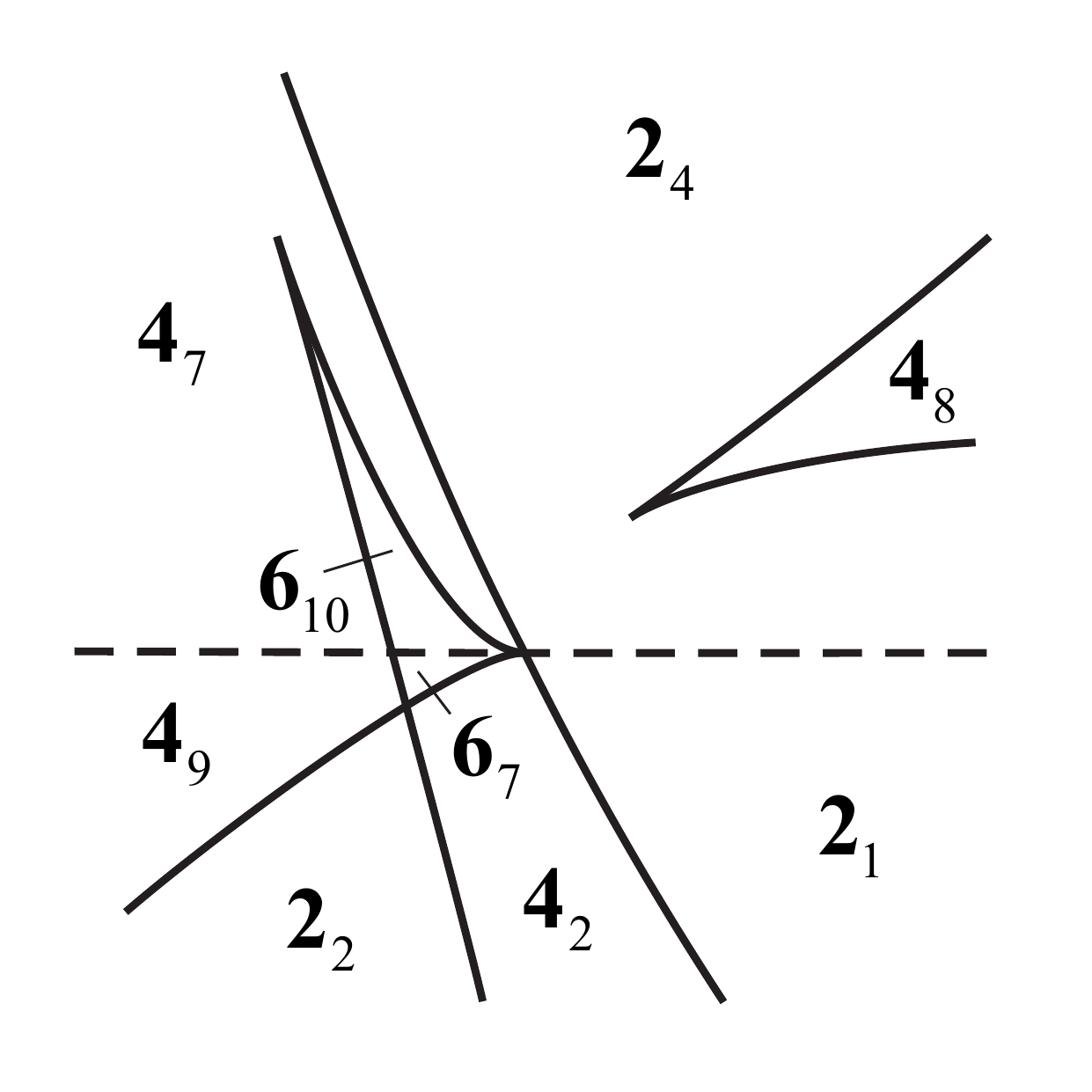}&\includegraphics[width=4cm]{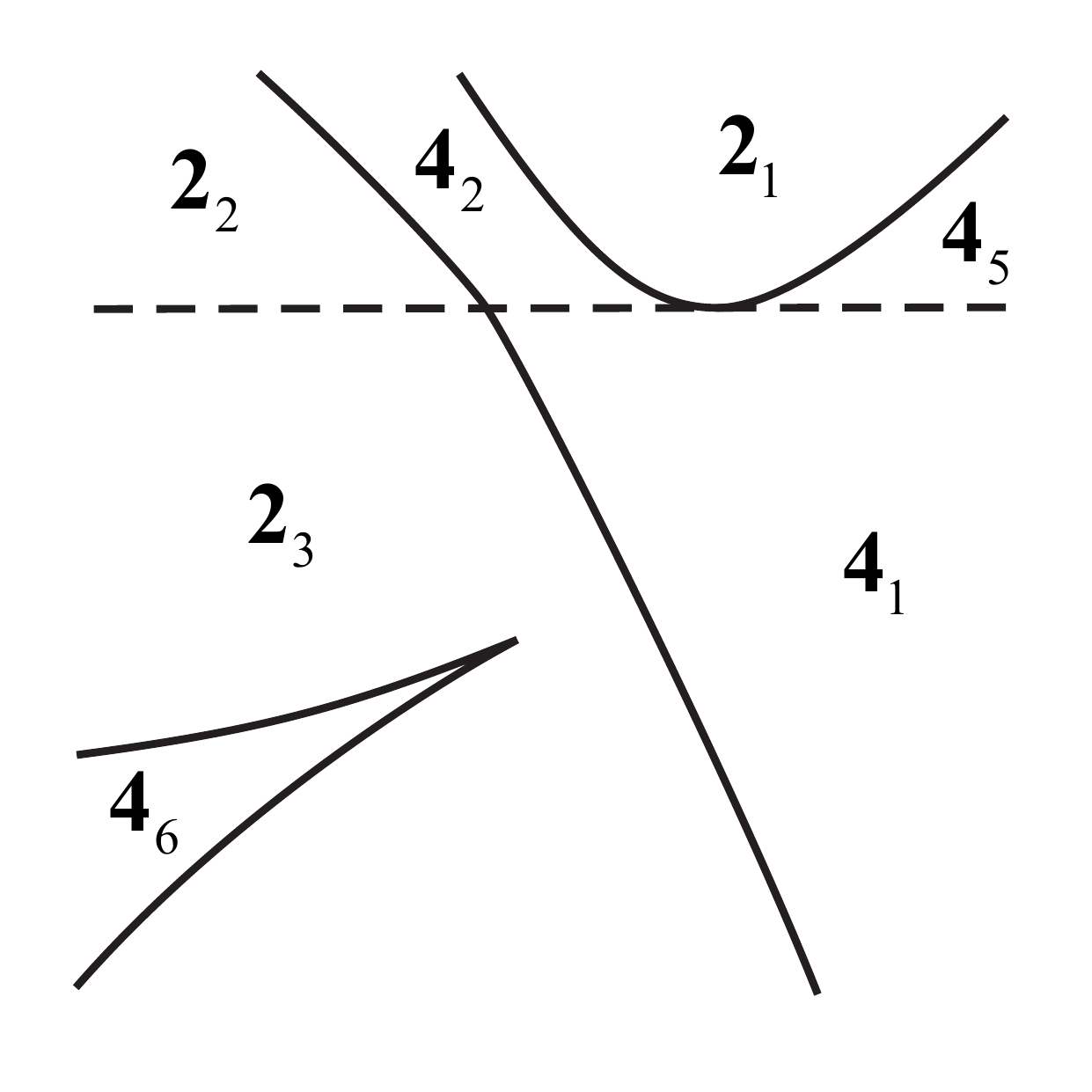}&
\includegraphics[width=4cm]{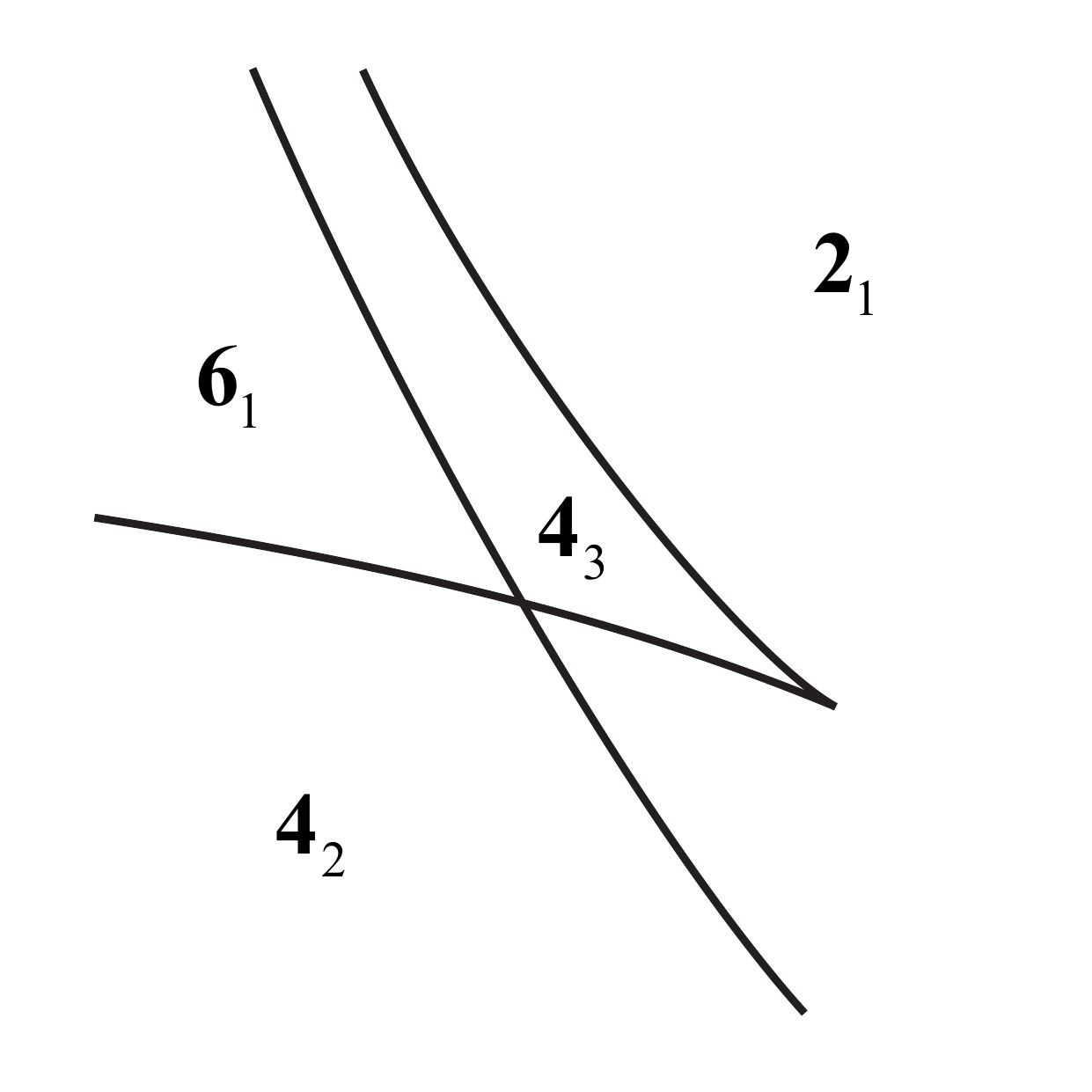}\\
$(t_1,q_5)\equiv(\ref{spoint2})$&$(t_1,q_5)\equiv(\ref{spoint1})$&$t_1=q_5=0$\\
\\
\end{tabular}
\caption{The skeletons $\Sigma_{S_4,t_1,q_5},S_4\neq0$ for some singular points of $B_{S_4}$.}
\label{sing-sigma}
\end{center}
\end{figure}

\begin{figure}
\begin{center}
\begin{tabular}{ccc}
15)\includegraphics[width=4cm]{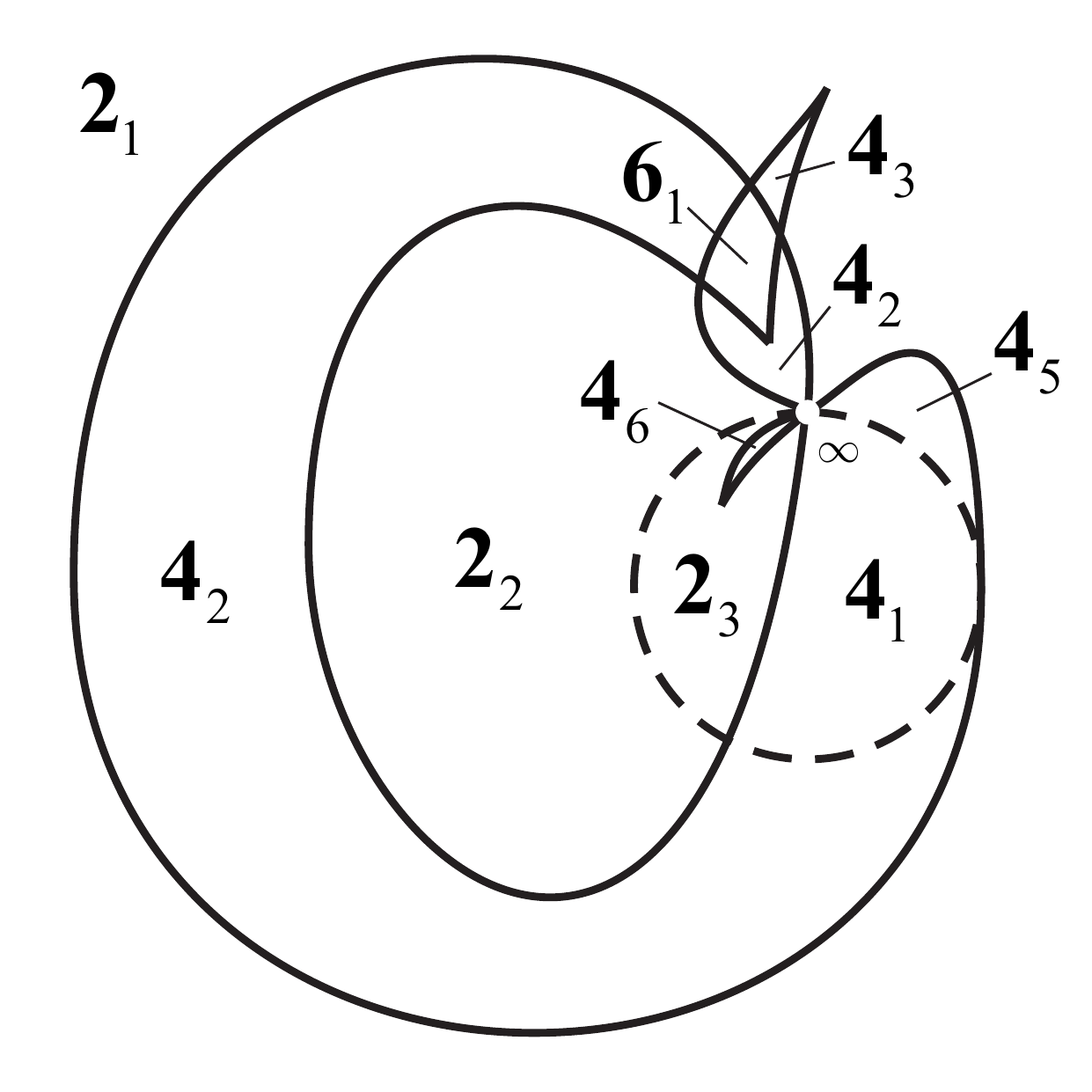}&&
34)\includegraphics[width=4cm]{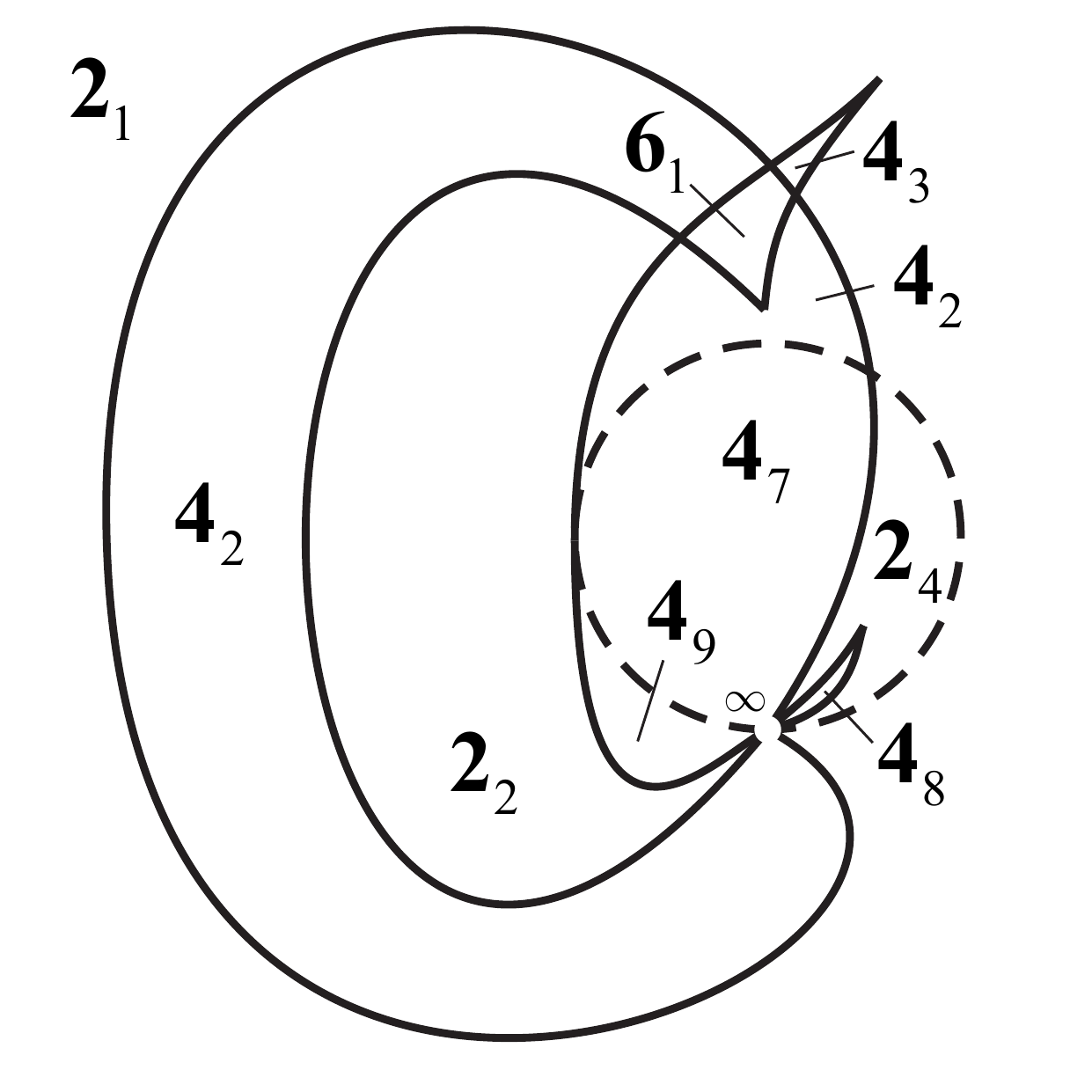}\\
\\
&$\stackrel{{\cal P}_{10}}{\longleftrightarrow}$&\\
\\
12)\includegraphics[width=4cm]{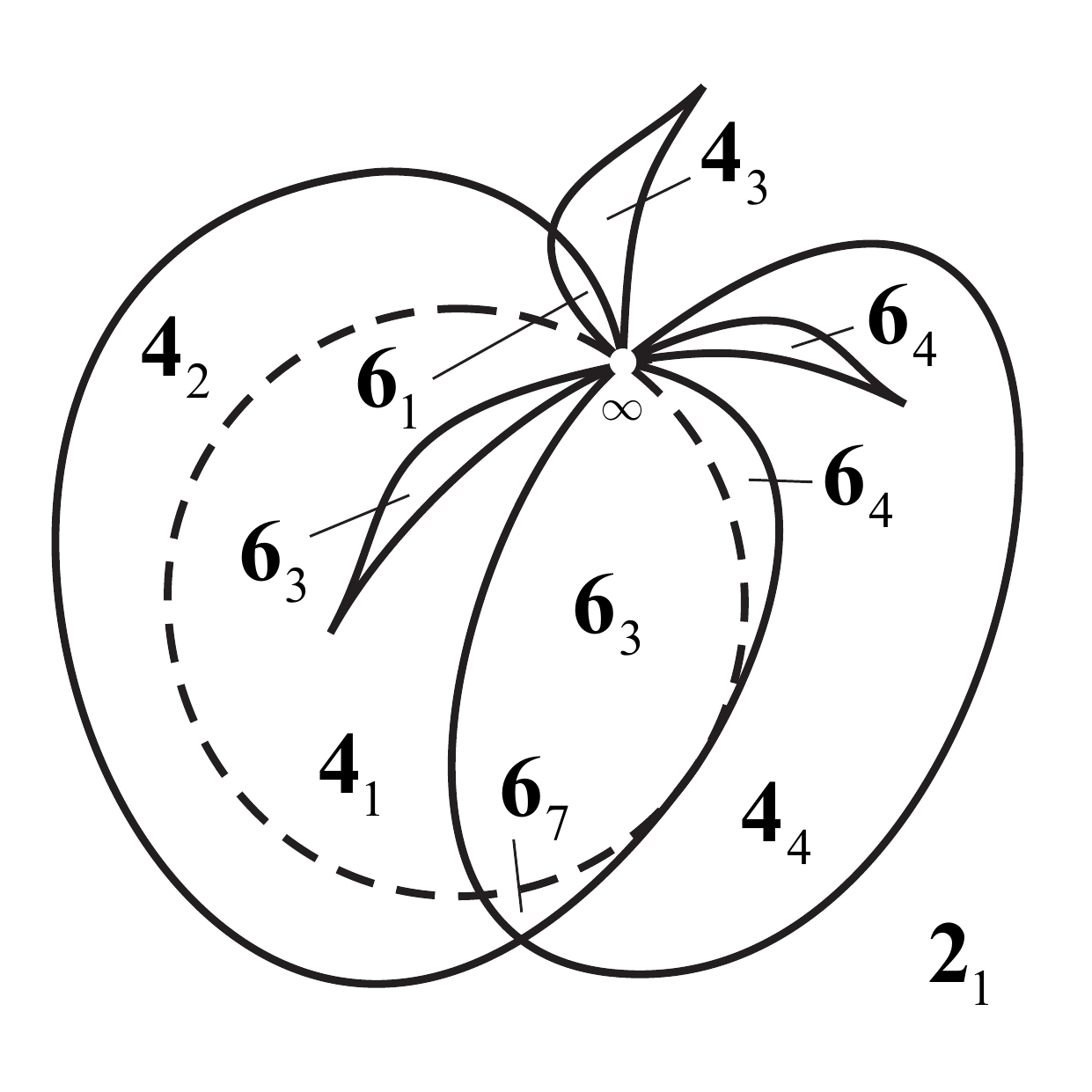}&&
51)\includegraphics[width=4cm]{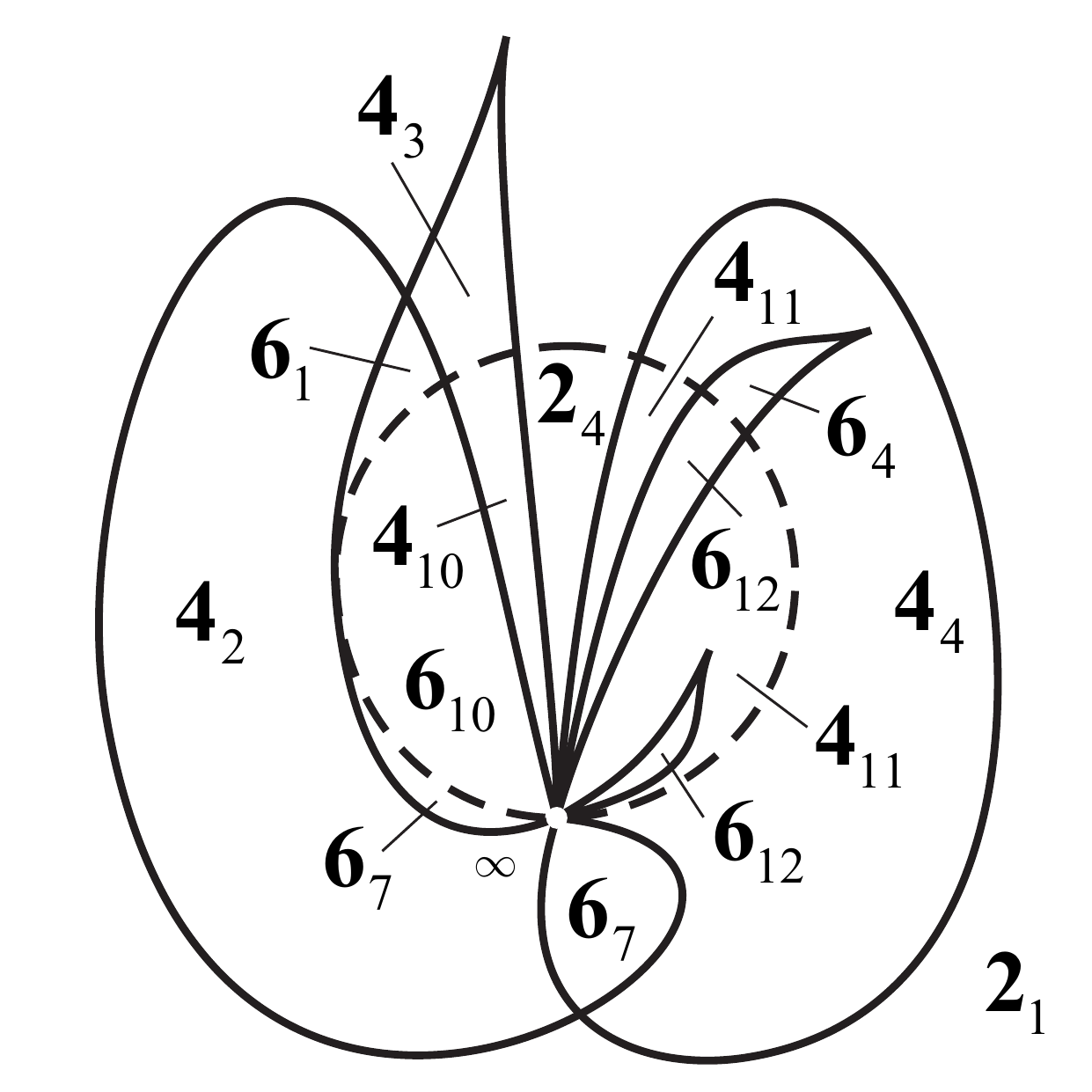}\\
\\
\hline
\\
10)\includegraphics[width=4cm]{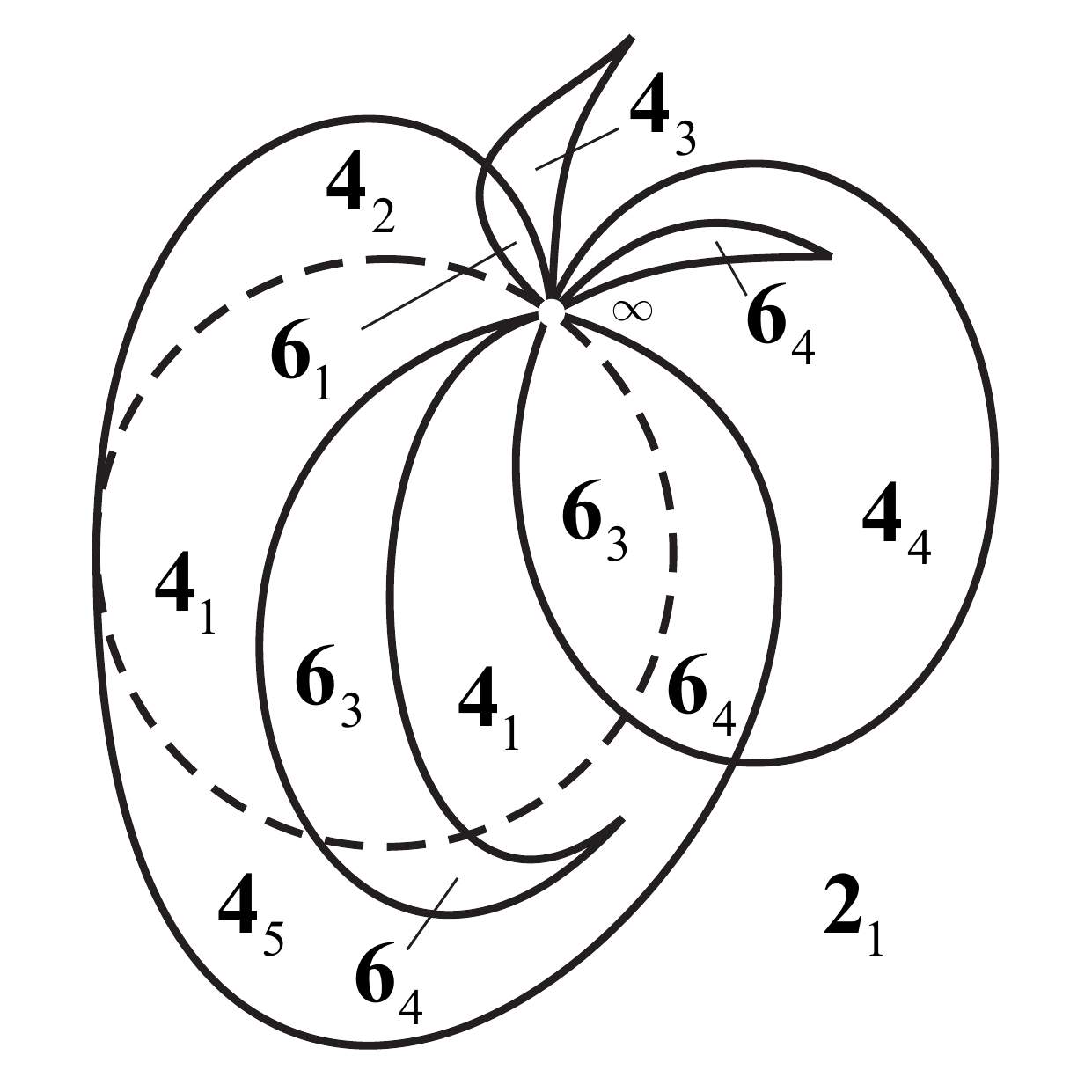}&&13)\includegraphics[width=4cm]{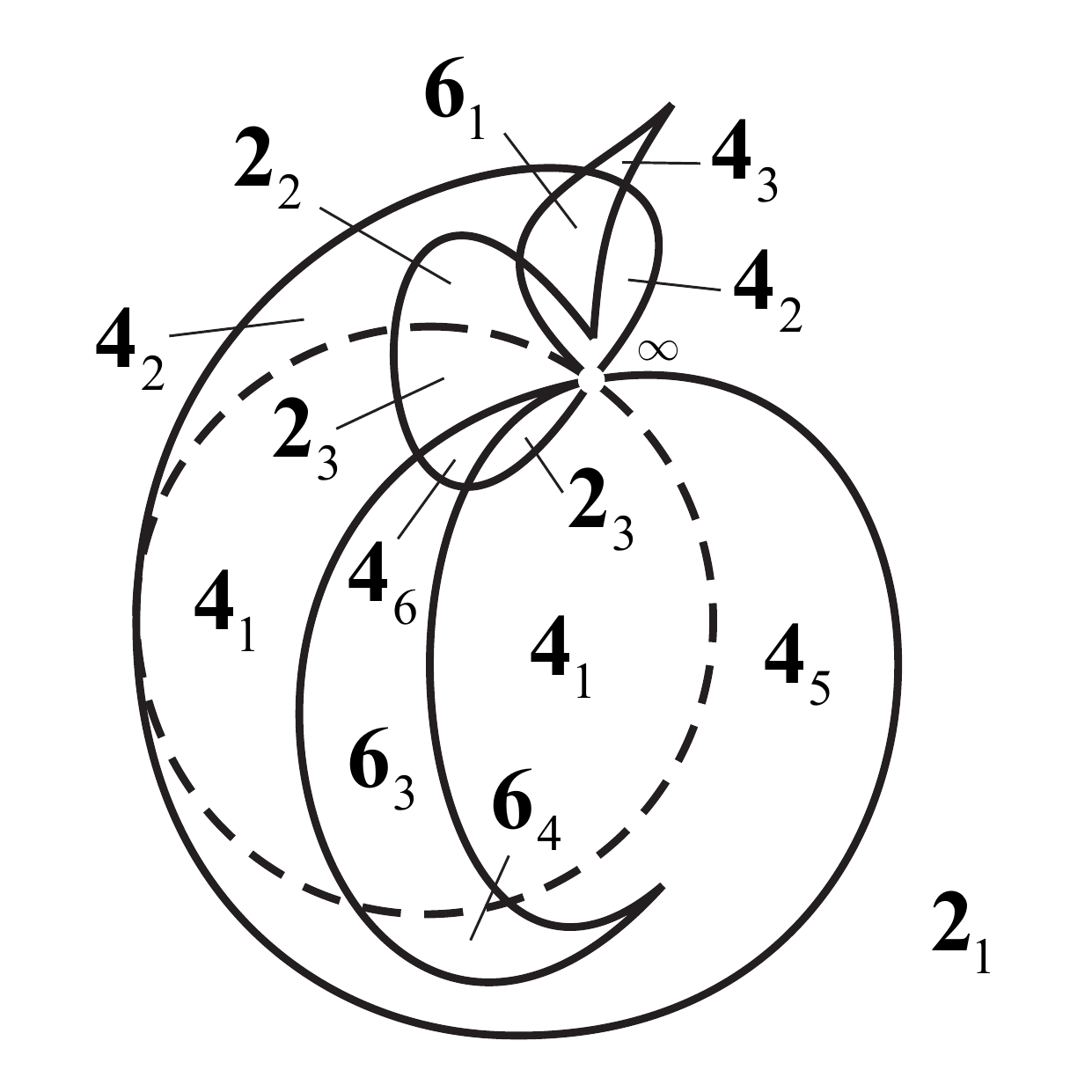}
\\
\\
&$\stackrel{{\cal P}_{11}}{\longleftrightarrow}$&\\
\\
48)\includegraphics[width=4cm]{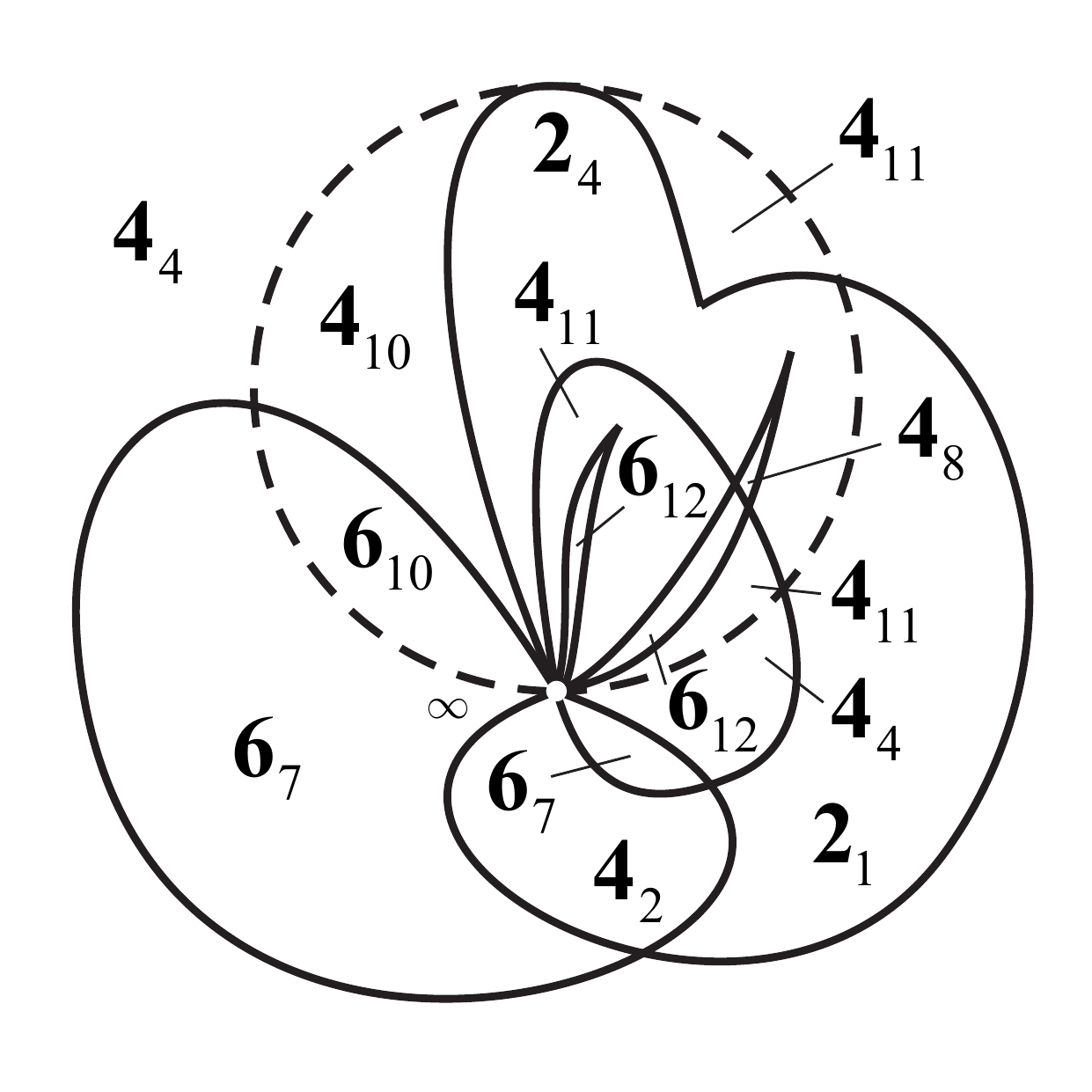}&&42)\includegraphics[width=4cm]{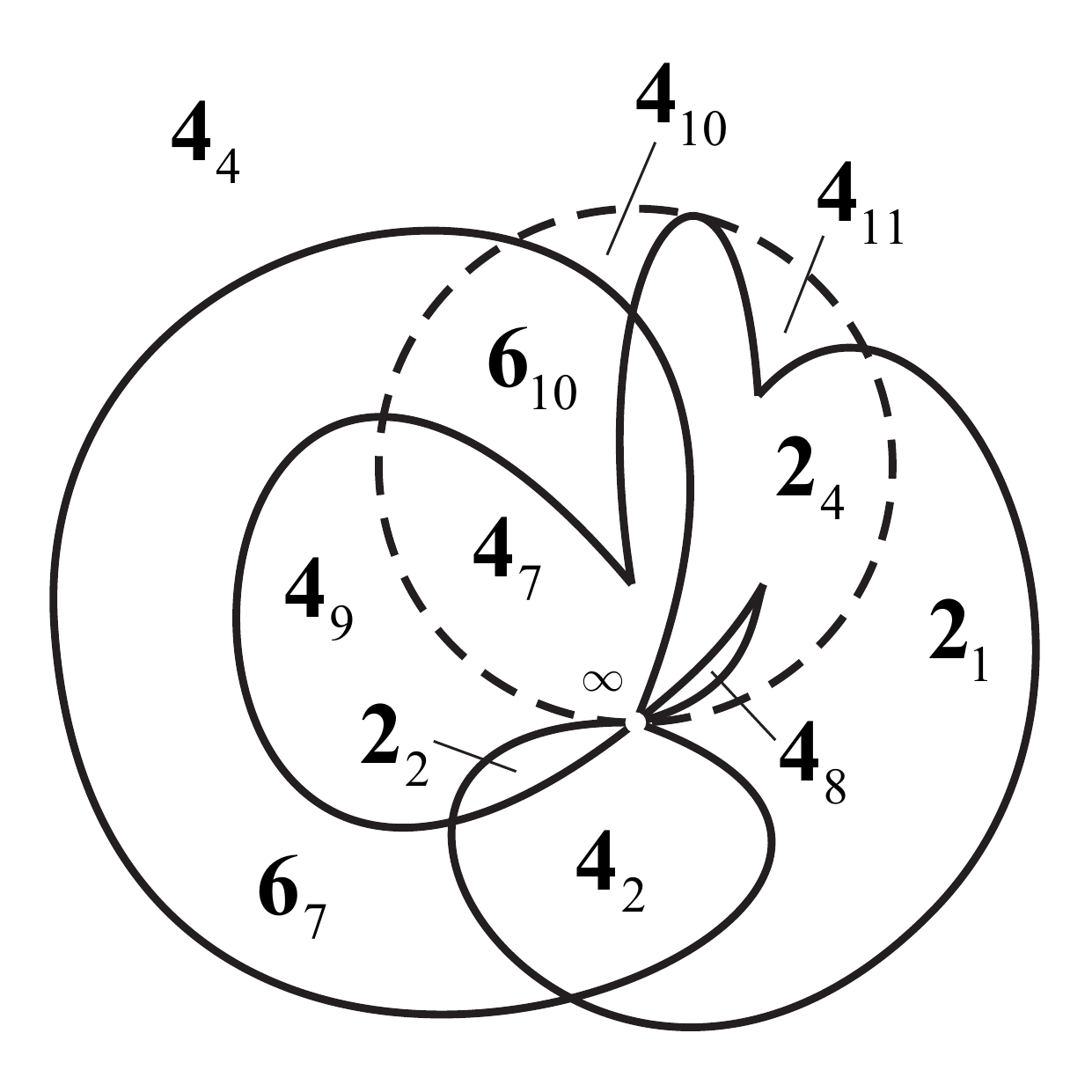}
\\
\\
\end{tabular}
\caption{Transformations of strata under perestroikas ${\cal P}_{10}$ and ${\cal P}_{11}$.}
\label{zony-infty}
\end{center}
\end{figure}

\medskip
Two-dimensional strata of stratifications $\mathfrak{S}_{S_4,t_1,q_5}$ are enumerated
(see Fig. \ref{zona1-15} -- \ref{sing-sigma}). This enumeration is not one-to-one but it determines the one-to-one enumeration for the four-dimensional strata of the stratification $\mathfrak{S}_{S_4}$.

The enumeration starts in domain $1$ and extends to the other domains $2-93$ using properties of perestroikas ${\cal P}_{1} - {\cal P}_{11}$ of $\Sigma_{S_4,t_1,q_5}$ when the point $(t_1,q_5)$ passes into the neighboring domains. A stratum of type $A_1^k$ with number $i$ is denoted by ${\bf k}_i$. The type of a stratum is not changed when $(t_1,q_5)$ transversally intersects the dashed line (at a point that does not lie on a solid line). As a result of a sequence of perestroikas, it may turn out that strata ${\bf k}_i$ and ${\bf k}_j$, where $i<j$, are different connected components of a section of the same stratum of type $A_1^k$ of the stratification $\mathfrak{S}_{S_4}$. Then the notation ${\bf k}_j$ changes everywhere to ${\bf k}_i$, and the number $j$ is not used further for enumeration of type $A_1^k$ strata. Similarly, we set ${\bf k}_j={\bf k}_i$ if the passages to a given point $(t_1,q_5)$ along two different smooth curves lead to different notations ${\bf k}_i$ and ${\bf k}_j$ for the same type $A_1^k$ stratum of $\mathfrak{S}_{S_4,t_1,q_5}$. Finally, we shift the enumeration to the left using the missing numbers. This enumeration naturally extends to two-dimensional strata of $\mathfrak{S}_{S_4,t_1,q_5}$ for $(t_1,q_5)\in B_{S_4}$. The stratum ${\bf k}_i$ of $\mathfrak{S}_{S_4}$ consists of points $(t_1,t_2,S_3,S_4,q_5)$ such that $(t_2,S_3)$ belongs to a stratum ${\bf k}_i$ of $\mathfrak{S}_{S_4,t_1,q_5}$.

\medskip
\predlo\label{predl-straty-S4} {\it The stratification $\mathfrak{S}_{S_4}, S_4\neq0$ has $31$ four-dimensional strata. Namely it has four strata ${\bf 2}_1 - {\bf 2}_4$ of type $A_1^2$, eleven strata ${\bf 4}_1 - {\bf 4}_{11}$ of type $A_1^4$, and sixteen strata ${\bf 6}_{1} - {\bf 6}_{16}$ of type $A_1^6$.}

\medskip
\remark To control the enumeration process in the third type perestroikas ${\cal P}_{10},{\cal P}_{11}$, we drew complete images of some skeletons $\Sigma_{S_4,t_1,q_5}$ under inversions with centers at various points of $\mathbb{R}^2=\{(t_2,S_3)\}$ (see Fig. \ref{zony-infty}).

\medskip
Two-dimensional strata of any stratification $\mathfrak{S}_{S_4,t_1,q_5},S_4\neq0$ are contractible. Next, we study the question on the contractibility of four-dimensional strata of $\mathfrak{S}_{S_4}$.

\section{Local properties of the complement to $\Sigma_{S_4},S_4\neq0$}\label{lok-property-1}

Let the map $f$ has a multisingularity of type $\mathcal{B}$ at a point $y$ of the target space. Since $f$ is Lagrangian stable, we see that for any $\mathcal{A}\in\mathbb{S}^+$ the intersection of the manifold $\mathcal{A}_f$ with the open ball of sufficiently small radius centred at $y$ is a smooth manifold, and the equivalence class of this manifold under diffeomorphisms depends only on $\mathcal{A}$ and $\mathcal{B}$. We denote this manifold by $\Xi_{\mathcal{A}}(\mathcal{B})$. Its Euler characteristic is called the adjacency index of a multisingularity of type $\mathcal{B}$ to a multisingularity of type $\mathcal{A}$ and is denoted by $J_{\mathcal{A}}(\mathcal{B})$ (see \cite{Sed2015}). 

Let us consider an arbitrary point $x=(t_1,t_2,S_3,S_4,q_5)$ in the source space of $f$. By $\Sigma^c(x;r)$ denote the intersection of the complement to $\Sigma$ with an open ball of radius $r>0$ centered at $x$. As before, there exists $r=r(x)$ such that $\Sigma^c(x;r)$ is diffeomorphic to $\Sigma^c(x;r(x))$ for any $r\leq r(x)$. The set $\Sigma^c(x;r(x))$ is denoted by $\Sigma^c(x)$. A connected component of the set $\Sigma^c(x)$ is called a component of type $A_1^k$ if it belongs to a type $A_1^k$ stratum of the stratification $\mathfrak{S}$.

Now from \cite[Theorem 5.1]{Sed2015} we obtain the following statement.

\medskip
\lemma\label{Ox} {\it Let $f$ have a singularity of type $X$ at $x$ and a multisingularity of type $X\mathcal{A}$ at $f(x)$. Then the Euler characteristic of the union of type $A_1^k$ components of the set $\Sigma^c(x)$ is equal to:
$$
kJ_{A_1^k}(X)J_{\mathbf{1}}(\mathcal{A})+(k-2)J_{A_1^{k-2}}(X)J_{A_1^2}(\mathcal{A})+(k-4)J_{A_1^{k-4}}(X)J_{A_1^4}(\mathcal{A})+\dots
$$
if $\deg X\equiv k\,(\mathrm{mod}\,2)$, and
$$
(k-1)J_{A_1^{k-1}}(X)J_{A_1}(\mathcal{A})+(k-3)J_{A_1^{k-3}}(X)J_{A_1^3}(\mathcal{A})+(k-5)J_{A_1^{k-5}}(X)J_{A_1^5}(\mathcal{A})+\dots
$$
if $\deg X\equiv k-1\,(\mathrm{mod}\,2)$.}

\medskip
Let us fix an arbitrary $S_4\neq0$. By $\Sigma^c_{S_4}(x)$ denote the intersection of the set $\Sigma^c(x)$ with the hyperplane $S_4=\mathrm{const}$. A connected component of the set $\Sigma^c_{S_4}(x)$ is called a component of type $A_1^k$ if it belongs to a type $A_1^k$ stratum of the stratification $\mathfrak{S}_{S_4}$.

By Lemma \ref{sech}, it follows that the set $\Sigma^c(x)$ is diffeomorphic to the direct product $\Sigma^c_{S_4}(x)\times (0,+\infty)$. Lemmas $\ref{raspad-E6},\ref{Ox}$ and the description \cite{Sed2015} of adjacencies of Lagrangian monosingularities of types $A_{\mu}^\pm,D_{\mu}^\pm$ imply the following properties of the set $\Sigma^c_{S_4}(x)$ (see Fig. \ref{zona1-15} -- \ref{sing-sigma}).

\medskip
\lemma\label{lok-dzeta} {\it Let $S_4\neq0$ and $\alpha$ be any of the curves $\zeta_2,\zeta_3^+$,
\begin{equation}
\zeta_1: \delta t_1<-\sqrt[5]{\frac{(13\sqrt{5}+29)S_4^4}{108}}\quad (\mbox{bounded by $(\ref{spoint1})$}),
\label{zeta1-l}
\end{equation}
\begin{equation}
\zeta_1: \delta t_1>-\sqrt[5]{\frac{(13\sqrt{5}+29)S_4^4}{108}},
\label{zeta1-r}
\end{equation}
\begin{equation}
\zeta_3^-: q_5<-\sqrt[5]{6(\sqrt{5}-1)S_4^2}\quad (\mbox{bounded by $(\ref{spoint4})$}),
\label{dzeta3-l}
\end{equation}
\begin{equation}
\zeta_3^-: -\sqrt[5]{6(\sqrt{5}-1)S_4^2}<q_5<-\sqrt[5]{\frac{3S_4^2}{4}}\quad (\mbox{between $(\ref{spoint4})$ and $(\ref{vozvrat})$}).
\label{dzeta3-r}
\end{equation}
Suppose $(t_1,q_5)\in\alpha$ and the curve $(\ref{psi-2})$ has a singularity of type $4/3$ at $(t_2,S_3)$. Then $f$ is non-singular at $x$ and has a multisingularity of type $A_4A_1^2$ at $f(x)$. The set $\Sigma^c_{S_4}(x)$ consists of three contractible components; there is one component of each type $A_1^2,A_1^4,A_1^6$.}

\medskip
\lemma\label{lok-eta2} {\it Let $S_4\neq0$ and $\alpha$ be any of the curves $\eta_3^+$,
\begin{equation}
\eta_1: q_5<\sqrt[5]{12 \left(\sqrt{5}-2\right)S_4^2}\quad (\mbox{bounded by $(\ref{spoint2})$}),
\label{eta1-l}
\end{equation}
\begin{equation}
\eta_1 : q_5>\sqrt[5]{12 \left(\sqrt{5}-2\right)S_4^2},
\label{eta1-r}
\end{equation}
\begin{equation}
\eta_2: \delta t_1<\frac{\sqrt[5]{144\left(\sqrt{5}+2\right)^2S_4^4}}{6\left(\sqrt{5}+3\right)}\quad (\mbox{bounded by $(\ref{spoint6})$}),
\label{eta2-l}
\end{equation}
\begin{equation}
\eta_2: \frac{\sqrt[5]{144\left(\sqrt{5}+2\right)^2S_4^4}}{6\left(\sqrt{5}+3\right)}<\delta t_1<\sqrt[5]{\frac{S_4^4}{432}}\quad (\mbox{between $(\ref{spoint6})$ and $(\ref{vozvrat})$}),
\label{eta2-m}
\end{equation}
\begin{equation}
\eta_2:
\delta t_1>\sqrt[5]{\frac{S_4^4}{432}}\quad (\mbox{bounded by $(\ref{vozvrat})$}),
\label{eta2-r}
\end{equation}
\begin{equation}
\eta_3^-: q_5<0,
\label{eta3-r}
\end{equation}
\begin{equation}
\eta_3^-: q_5>0.
\label{eta3-l}
\end{equation}
Suppose $(t_1,q_5)\in\alpha$ and
\begin{equation}
\begin{tabular}{c}
$(t_2,S_3)$ is \mbox{the intersection point of two branches of the curve $(\ref{psi-2})$; one of them }\\ \mbox{is smooth and the other has semicubical cusp at this point.}
\end{tabular}
\label{usl-1}
\end{equation}
Then $f$ is non-singular at $x$ and has a multisingularity of type $A_3^{\pm}A_2A_1$ at $f(x)$. The set $\Sigma^c_{S_4}(x)$ consists of $1,2,1$ contractible components of types $A_1^2,A_1^4,A_1^6$, respectively. The components of type $A_1^4$ belong to the same stratum ${\bf 4}_{11}$ of $\mathfrak{S}_{S_4}$ if $\alpha\equiv(\ref{eta2-m})$ and to different strata in other cases.}

\medskip
\lemma\label{lok-ksi1} {\it Let $S_4\neq0$ and $\alpha$ be any of the curves $\xi_1,\xi_2,\xi_4^\pm,\xi_5^\pm$,
\begin{equation}
\xi_3 : q_5<0,
\label{xi3-b}
\end{equation}
\begin{equation}
\xi_3 : q_5>0.
\label{xi3-t}
\end{equation}
Suppose $(t_1,q_5)\in\alpha$ and
\begin{equation}
\begin{tabular}{c}
$(t_2,S_3)$ is \mbox{a semicubical cusp of the curve $(\ref{psi-2})$ lying on $(\ref{crit-tochki})$.}\\
\end{tabular}
\label{usl-2}
\end{equation}
Then $f$ has a singularity of type $A_2$ at $x$ and a multisingularity of type $A_3^{\pm}A_2A_1$ at $f(x)$. The set $\Sigma^c_{S_4}(x)$ consists of four contractible components; there are two components of each type $A_1^4,A_1^6$. They belong to different strata of $\mathfrak{S}_{S_4}$.}

\medskip
\lemma\label{lok-tau} {\it Let $S_4\neq0$ and $\alpha$ be any of the curves $\tau_1$,
\begin{equation}
\tau_2 : \delta t_1<\frac{\sqrt[5]{144\left(\sqrt{5}+2\right)^2S_4^4}}{6\left(\sqrt{5}+3\right)}\quad (\mbox{bounded by $(\ref{spoint6})$}),
\label{tau2-b}
\end{equation}
\begin{equation}
\tau_2: \delta t_1>\frac{\sqrt[5]{144\left(\sqrt{5}+2\right)^2S_4^4}}{6\left(\sqrt{5}+3\right)},
\label{tau2-r}
\end{equation}
\begin{equation}
\tau_3 : \delta t_1<\frac{\sqrt[5]{144\left(\sqrt{5}-2\right)^2S_4^4}}{6\left(3-\sqrt{5}\right)}\quad (\mbox{bounded by $(\ref{spoint2})$}),
\label{tau3-l}
\end{equation}
\begin{equation}
\tau_3 : \delta t_1>\frac{\sqrt[5]{144\left(\sqrt{5}-2\right)^2S_4^4}}{6\left(3-\sqrt{5}\right)},
\label{tau3-r}
\end{equation}
\begin{equation}
\tau_4 : q_5<0,
\label{tau4-}
\end{equation}
\begin{equation}
\tau_4 : q_5>0.
\label{tau4+}
\end{equation}
Suppose $(t_1,q_5)\in\alpha$ and $(t_2,S_3)$ is the intersection point of two smooth branches of the curve $(\ref{psi-2})$, one of which has a simple tangency with $(\ref{crit-tochki})$ at this point. Then $f$ has a singularity of type $A_3^{\pm}$ at $x$ and a multisingularity of type $A_3^{\pm}A_2A_1$ at $f(x)$. The set $\Sigma^c_{S_4}(x)$ consists of $1,4,3$ contractible components of types $A_1^2,A_1^4,A_1^6$, respectively. The components of type $A_1^6$ belong to different strata of $\mathfrak{S}_{S_4}$. The components of type $A_1^4$ belong to different strata if $\alpha$ is one of the curves $\tau_1,(\ref{tau3-l}),(\ref{tau3-r}),(\ref{tau4-})$. In the cases $\alpha\equiv(\ref{tau2-b}),(\ref{tau2-r}),(\ref{tau4+})$, they belong to three strata, one of which contains two components and each of the other two contains exactly one; strata with two components are ${\bf 4}_4,{\bf 4}_{11},{\bf 4}_2$, respectively.}

\medskip
\lemma\label{lok-gamma2} {\it Let $S_4\neq0$ and $\alpha$ be any of the curves $\gamma_2^-$,
\begin{equation}
\gamma_2^+: q_5<-\sqrt[5]{12\left(\sqrt{5}+2\right)S_4^2}\quad (\mbox{bounded by $(\ref{spoint6})$}),
\label{gamma2+b}
\end{equation}
\begin{equation}
\gamma_2^+: -\sqrt[5]{12\left(\sqrt{5}+2\right)S_4^2}<q_5<\sqrt[5]{12 \left(\sqrt{5}-2\right)S_4^2}\quad (\mbox{between $(\ref{spoint6})$ and $(\ref{spoint2})$}),
\label{gamma2-m}
\end{equation}
\begin{equation}
\gamma_2^+: q_5>\sqrt[5]{12 \left(\sqrt{5}-2\right)S_4^2}\quad (\mbox{bounded by $(\ref{spoint2})$}).
\label{gam2+t}
\end{equation}
Suppose $(t_1,q_5)\in\alpha$ and $(t_2,S_3)\equiv(\ref{spec1})$. Then:

{\rm1)} if $\alpha\equiv(\ref{gamma2-m})$, then $f$ has a monosingularity o type $A_4$ at $f(x)$; the set $\Sigma^c_{S_4}(x)$ consists of six contractible components;  there are two components of type $A_1^2$ and four components of type $A_1^4$; they belong to different strata of $\mathfrak{S}_{S_4}$;

{\rm2)} if $\alpha\not\equiv(\ref{gamma2-m})$, then $f$ has a singularity of type $A_4$ at $x$ and a multisingularity of type $A_4A_1^2$ at $f(x)$; the set $\Sigma^c_{S_4}(x)$ consists of six contractible components; there are two components of type $A_1^4$ and four components of type $A_1^6$; they belong to different strata of $\mathfrak{S}_{S_4}$.}

\medskip
\lemma\label{lok-34} {\it Let $S_4\neq0$, $(t_1,q_5)$ be any of points $(\ref{spoint2}),(\ref{spoint6})$, and $(t_2,S_3)\equiv(\ref{spec1})$. Then $f$ has a singularity of type $A_4$ at $x$ and a multisingularity of type $A_4A_2$ at $f(x)$. The set $\Sigma^c_{S_4}(x)$ consists of $2,6,4$ contractible components of types $A_1^2,A_1^4,A_1^6$, respectively. They belong to different strata of $\mathfrak{S}_{S_4}$ if $(t_1,q_5)\equiv(\ref{spoint2})$. In the case $(t_1,q_5)\equiv(\ref{spoint6})$, the components of types $A_1^2,A_1^6$ belong to different strata; the components of type $A_1^4$ belong to four strata; two of these strata are ${\bf 4}_4,{\bf 4}_{11}$; each of them contains exactly two components; each of the other two strata contains exactly one component.}

\medskip
\lemma\label{gamm3} {\it Let $S_4\neq0$ and $\alpha$ be any of the curves $\gamma_3^-$,
\begin{equation}
\gamma_3^+: q_5<-\sqrt[5]{12\left(\sqrt{5}+2\right)S_4^2}\quad (\mbox{bounded by $(\ref{spoint6})$}),
\label{gam3+bb}
\end{equation}
\begin{equation}
\gamma_3^+: -\sqrt[5]{12\left(\sqrt{5}+2\right)S_4^2}<q_5<-\sqrt[5]{12S_4^2}\quad (\mbox{between $(\ref{spoint6})$ and $(\ref{spoint5})$}),
\label{gamma3+-m}
\end{equation}
\begin{equation}
\gamma_3^+: -\sqrt[5]{12S_4^2}<q_5<-\sqrt[5]{6(\sqrt{5}-1)S_4^2},\quad (\mbox{between $(\ref{spoint5})$ and $(\ref{spoint4})$}),
\label{gam3+bt}
\end{equation}
\begin{equation}
\gamma_3^+: 0<q_5<\sqrt[5]{12 \left(\sqrt{5}-2\right)S_4^2},\quad (\mbox{between the origin and $(\ref{spoint2})$}),
\label{gam3+m}
\end{equation}
\begin{equation}
\gamma_3^+: q_5>\sqrt[5]{12 \left(\sqrt{5}-2\right)S_4^2}\quad (\mbox{bounded by $(\ref{spoint2})$}).
\label{gam3+t}
\end{equation}
Suppose $(t_1,q_5)\in\alpha$; $(t_2,S_3)$ is a pairwise transversal intersection point of the line $(\ref{crit-tochki})$ and two smooth branches of the curve $(\ref{psi-2})$. Then $f$ has a singularity of type $A_2$ at $x$ and a multisingularity of type $A_2^3$ at $f(x)$. The set $\Sigma^c_{S_4}(x)$ consists of $2,4,2$ contractible components of types $A_1^2,A_1^4,A_1^6$, respectively. The components of types $A_1^2,A_1^6$ belong to different strata of $\mathfrak{S}_{S_4}$. The components of type $A_1^4$ belong to different strata if $\alpha$ is one of the curves $(\ref{gam3+bb}),(\ref{gam3+m}),(\ref{gam3+t})$. They belong to strata ${\bf 4}_1,{\bf 4}_{2}$ if $\alpha\equiv\gamma_3^-$ and to ${\bf 4}_4,{\bf 4}_{11}$ if $\alpha$ is one of the curves $(\ref{gamma3+-m}),(\ref{gam3+bt})$; each of the strata ${\bf 4}_1,{\bf 4}_{2},{\bf 4}_4,{\bf 4}_{11}$ in the mentioned cases contains exactly two components.}

\medskip
\lemma\label{lok-36} {\it Let $S_4\neq0$, $(t_1,q_5)$ be any of points $(\ref{spoint1}),(\ref{spoint4})$, and $(t_2,S_3)$ be determined by formulas $(\ref{psi-2})$ for $u$ given by $(\ref{u-36}),(\ref{u-37})$, respectively. Then $f$ has a singularity of type $A_2$ at $x$ and a multisingularity of type $A_4A_2$ at $f(x)$. The set $\Sigma^c_{S_4}(x)$ consists of six contractible components; there are two components of each type $A_1^2,A_1^4,A_1^6$; they belong to different strata of $\mathfrak{S}_{S_4}$.}

\medskip
\lemma\label{lok-gamma-} {\it Let $S_4\neq0$, $(t_1,q_5)\in\gamma_1^-$, and $(t_2,S_3)\in(\ref{psi-4})$. Then:

{\rm1)} if $(t_2,S_3)$ does not coincide with $(\ref{spec2})$ and $(\ref{spec-kas})$, then $f$ is non-singular at $x$ and has a multisingularity of type $A_2A_1^4$ at $f(x)$; the set $\Sigma^c_{S_4}(x)$ consists of two contractible components; there is one component of each type $A_1^4,A_1^6$;

{\rm2)} if $(t_2,S_3)\equiv(\ref{spec2})$, then $f$ is non-singular at $x$ and has a multisingularity of type $D_4^+A_1^2$ at $f(x)$; the set $\Sigma^c_{S_4}(x)$ consists of $1,2,1$ contractible components of types $A_1^2,A_1^4,A_1^6$, respectively; they belong to different strata of $\mathfrak{S}_{S_4}$;

{\rm3)} if $(t_2,S_3)\equiv(\ref{spec-kas})$, then $f$ has a singularity of type $A_2$ at $x$ and a multisingularity of type $A_2^2A_1^2$ at $f(x)$; the set $\Sigma^c_{S_4}(x)$ consists of four contractible components; there are two components of each type $A_1^4,A_1^6$; they belong to different strata of $\mathfrak{S}_{S_4}$.}

\medskip
\lemma\label{lok-gamma+} {\it Let $S_4\neq0$, $(t_1,q_5)\in\gamma_1^+$, and $(t_2,S_3)\in(\ref{psi-4})$. Then:

{\rm1)} if $(t_2,S_3)$ lies outside the segment with ends $(\ref{peresech1}),(\ref{peresech2})$ and does not coincide with $(\ref{spec2})$, $(\ref{spec-kas})$, then $f$ is non-singular at $x$ and has a multisingularity of type $A_2A_1^4$ at $f(x)$; the set $\Sigma^c_{S_4}(x)$ consists of two contractible components; there is one component of each type $A_1^4,A_1^6$;

{\rm2)} if $(t_2,S_3)$ lies between $(\ref{peresech1}),(\ref{peresech2})$ and $(t_2,S_3)\not\equiv(\ref{spec-kas})$, then $f$ is non-singular at $x$ and has a multisingularity of type $A_2A_1^2$ at $f(x)$; the set $\Sigma^c_{S_4}(x)$ consists of two contractible components; there is one component of each type $A_1^2,A_1^4$;

{\rm3)} if $(t_1,q_5)\in(\ref{gam1+b})\cup(\ref{gamma1+part})\cup(\ref{gam1+t})$ and $(t_2,S_3)$ is any of the points $(\ref{peresech1}),(\ref{peresech2})$, then $f$ is non-singular at $x$ and has a multisingularity of type $A_2^2A_1^2$ at $f(x)$; the set $\Sigma^c_{S_4}(x)$ consists of $1,2,1$ contractible components of types $A_1^2,A_1^4,A_1^6$, respectively; the components of type $A_1^4$ belong to different strata of $\mathfrak{S}_{S_4}$ if $(t_2,S_3)\equiv(\ref{peresech1})$; in the case $(t_2,S_3)\equiv(\ref{peresech2})$, these components belong to different strata if $(t_1,q_5)\in(\ref{gam1+t})$, to the stratum ${\bf 4}_{11}$ if $(t_1,q_5)\in(\ref{gamma1+part})$, and to ${\bf 4}_{4}$ if $(t_1,q_5)\in(\ref{gam1+b})$;

{\rm4)} if $(t_1,q_5)$ is any of the points $(\ref{spoint5}),(\ref{vozvrat})$ and $(t_2,S_3)\equiv(\ref{peresech1})$, then $f$ is non-singular at $x$ and has a multisingularity of type $A_2^2A_1^2$ at $f(x)$; the set $\Sigma^c_{S_4}(x)$ consists of $1,2,1$ contractible components of types $A_1^2,A_1^4,A_1^6$, respectively; the components of type $A_1^4$ belong to different strata of $\mathfrak{S}_{S_4}$;

{\rm5)} if $(t_1,q_5)\equiv(\ref{spoint5})$ and $(t_2,S_3)\equiv(\ref{peresech2})$, then $f$ has a singularity of type $A_2$ at $x$ and a multisingularity of type $A_2^3$ at $f(x)$; the set $\Sigma^c_{S_4}(x)$ consists of $2,4,2$ contractible components of types $A_1^2,A_1^4,A_1^6$, respectively; the components of types $A_1^2,A_1^6$ belong to different strata of $\mathfrak{S}_{S_4}$; the components of type $A_1^4$ belong to strata ${\bf 4}_4,{\bf 4}_{11}$, each of which contains exactly two components of $\Sigma^c_{S_4}(x)$;

{\rm6)} if $(t_1,q_5)\in(\ref{gam1+b})$ and $(t_2,S_3)\equiv(\ref{spec-kas})$, then $f$ has a singularity of type $A_2$ at $x$ and a multisingularity of type $A_2^2A_1^2$ at $f(x)$; the set $\Sigma^c_{S_4}(x)$ consists of four contractible components; there are two components of each type $A_1^4,A_1^6$; they belong to different strata of $\mathfrak{S}_{S_4}$;

{\rm7)} if $(t_1,q_5)\in(\ref{gamma1+part})\cup(\ref{gam1+t})\cup(\ref{vozvrat})$ and $(t_2,S_3)\equiv(\ref{spec-kas})$, then $f$ has a singularity of type $A_2$ at $x$ and a multisingularity of type $A_2^2$ at $f(x)$; the set $\Sigma^c_{S_4}(x)$ consists of four contractible components; there are two components of each type $A_1^2,A_1^4$; they belong to different strata of $\mathfrak{S}_{S_4}$;

{\rm8)} if $(t_1,q_5)\in(\ref{gam1+t})$ and $(t_2,S_3)\equiv(\ref{spec2})$, then $f$ is non-singular at $x$ and has a multisingularity of type $D_4^+A_1^2$ at $f(x)$; the set $\Sigma^c_{S_4}(x)$ consists of $1,2,1$ contractible components of types $A_1^2,A_1^4$, $A_1^6$, respectively; the components of type $A_1^4$ belong to different strata of $\mathfrak{S}_{S_4}$;

{\rm9)} if $(t_1,q_5)\in(\ref{gamma1+bv})$ and $(t_2,S_3)\equiv(\ref{spec2})$, then $f$ is non-singular at $x$ and has a multisingularity of type $D_4^-A_1^2$ at $f(x)$; the set $\Sigma^c_{S_4}(x)$ consists of three components; there are one component of type $A_1^4$ and two components of type $A_1^6$; the component of type $A_1^4$ is homotopy equivalent to $S^1$ and belongs to the stratum ${\bf 4}_{11}$; the components of type $A_1^6$ are contractible and belong to different strata of $\mathfrak{S}_{S_4}$;

{\rm10)} if $(t_1,q_5)\equiv(\ref{vozvrat})$ and $(t_2,S_3)\equiv(\ref{spec2})$, then $f$ is non-singular at $x$ and has a multisingularity of type $D_5^{\delta}A_1$ at $f(x)$; the set $\Sigma^c_{S_4}(x)$ consists of $1,2,2$ components of type $A_1^2,A_1^4,A_1^6$, respectively; they belong to different strata of $\mathfrak{S}_{S_4}$; the components of types $A_1^2,A_1^6$ are contractible; one component of type $A_1^4$ is contractible and the other one is homotopy equivalent to $S^1$; the non-contractible component belongs to the stratum  ${\bf 4}_{11}$.}

\medskip
\remark The non-contractible component of the set $\Sigma^c_{S_4}(x)$ in item 9 of Lemma \ref{lok-gamma+} is homotopy equivalent to a non-contractible connected component of the complement in $\mathbb{R}^{3}$ to the caustic of a Lagrangian germ of type $D_4^-$. As a curve defining the generator of its fundamental group for $x=(\delta,\frac{7}{6}\delta,10\delta,-6\delta,-3)$, we can take, for example, a circle
\begin{equation}
\delta t_1=\frac1{18}(18+\cos\varphi),\quad \delta t_2=\frac1{18}(21+\sin\varphi),\quad
\delta S_3=10,\quad \delta S_4=-6,\quad q_5=-3.
\label{fgr-D4}
\end{equation}

\medskip
\lemma\label{lok-int} {\it Let $S_4\neq0$, $(t_1,q_5)\equiv(\ref{trper-xi5-g1})$ and $(t_2,S_3)$ be an interior point of the segment with ends at $(\ref{spec-kas})$ and the point that satisfies the condition $(\ref{usl-2})$. Then $f$ has a singularity of type $A_2$ at $x$ and a multisingularity of type $A_2A_1^2$ at $f(x)$. The set $\Sigma^c_{S_4}(x)$ consists of two contractible components of $A_1^4$; they belong to different strata of $\mathfrak{S}_{S_4}$.}

\medskip
\lemma\label{lok-00} {\it The skeleton $\Sigma_{S_4,t_1,q_5}$ for $S_4\neq0$ and $t_1=q_5=0$ does not have non-ordinary singular points. Connected components of the intersection of any stratum of $\mathfrak{S}_{S_4}$ with a neighbourhood of infinity are contractible.}

\section{The topology of type $A_1^k$ strata of $\mathfrak{S}_{S_4},S_4\neq0$}

Let $x=(t_1,t_2,S_3,S_4,q_5)$. We fix $S_4\neq0$. By $R_{S_4}({\bf k}_i)$ denote the set of points $(t_1,q_5)\in\mathbb{R}^2$ such that $\mathfrak{S}_{S_4,t_1,q_5}$ has  at least one stratum ${\bf k}_i$.

\medskip
{\sc 1) Strata of type $A_1^2$.}
\medskip

\lemma\label{A2-1-styag} {\it The stratum ${\bf 2}_1\in\mathfrak{S}_{S_4},S_4\neq0$ is contractible.}

{\sc Proof.} Let $R_1$ be the domain in $\mathbb{R}^2$ that is bounded by the curve $\gamma_1^+$ and contains the origin; the subset $R_2\subset R_1$ be the union of the curve $(\ref{eta2-l})$ and the domain bounded by the line $(\ref{eta2-l})\cup(\ref{spoint6})\cup(\ref{gam3+bb})$. Then
$R_{S_4}({\bf 2}_1)=R_1$. The stratification $\mathfrak{S}_{S_4,t_1,q_5}$ has two strata ${\bf 2}_1$ if $(t_1,q_5)\in R_2$ (including domains $1-7,73-76,83$) and one stratum ${\bf 2}_1$ if $(t_1,q_5)\in R_1\setminus R_2$ (including domains $8-49,51-58,61,64-70,77-82,84-89$).

The sets $R_1,R_2$ are contractible. Strata ${\bf 2}_{1}\in\mathfrak{S}_{S_4,t_1,q_5}$ have two forms ${\mathcal R}_1,{\mathcal R}_2$. A stratum has the form ${\mathcal R}_i$ if: 1) $(t_1,q_5)\in R_i$; 2) the boundary of a stratum of the form ${\mathcal R}_2$ has one or three singular points; one of them is a transversal intersection point of two smooth branches of the curve (\ref{psi-2}) and each of the other two (if any) is a transversal intersection point of the line (\ref{crit-tochki}) with a smooth branch of the curve (\ref{psi-2}).

The stratification $\mathfrak{S}_{S_4,t_1,q_5}$ has one stratum of the form ${\mathcal R}_{i}$ for any $(t_1,q_5)\in R_i$. Each stratum of the form ${\mathcal R}_{i}$ is contracted to a point $(t_2,S_3)$ that continuously depends on $(t_1,q_5)$ (see Fig. \ref{zona1-15} -- \ref{sing-sigma}). Hence, the set of points $x$ such that $(t_1,q_5)\in R_i$ and $(t_2,S_3)$ belongs to a stratum of the form ${\mathcal R}_{i}$ is the total space of a locally trivial bundle. This space is contractible. It is called the part ${\mathcal R}_{i}$ of the stratum ${\bf 2}_{1}\in\mathfrak{S}_{S_4}$.

Consider the set $\beta$ of points $x$ such that $(t_1,q_5)$ belongs to the closure of the curve (\ref{eta2-l}) and $(t_2,S_3)$ satisfies the condition (\ref{usl-1}). It is closed and contractible. By Lemmas \ref{lok-eta2}, \ref{lok-34}, it follows that the intersection of the set $\Sigma^c_{S_4}(x)$ with the stratum ${\bf 2}_{1}$ is contractible for any $x\in\beta$. But Whitney stratification is locally trivial (see \cite{Mather}). Hence, there is a neighbourhood $U$ of the set $\beta$ such that the intersection $U\cap{\bf 2}_{1}$ is contractible. It remains to note that parts ${\mathcal R}_1$ and ${\mathcal R}_2$ of the stratum ${\bf 2}_1$ are contracted into $U$. $\Box$

\medskip
\lemma\label{A2-2-styag} {\it The stratum ${\bf 2}_2\in\mathfrak{S}_{S_4},S_4\neq0$ is contractible.}

{\sc Proof.} The set $R=R_{S_4}({\bf 2}_2)$ is the domain $q_5>0$. The stratification $\mathfrak{S}_{S_4,t_1,q_5}$ has one stratum ${\bf 2}_2$ for any point $(t_1,q_5)\in R$ (including domains $13-45$). This stratum is contracted to a point that continuously depends on $(t_1,q_5)$. Hence, the stratum ${\bf 2}_{2}\in\mathfrak{S}_{S_4}$ is homotopy equivalent to $R$. $\Box$

\medskip
\lemma\label{A2-3-styag} {\it The strata ${\bf 2}_3,{\bf 2}_4\in\mathfrak{S}_{S_4},S_4\neq0$ are contractible.}

{\sc Proof.} ${\bf 2}_3$) Let $R_1$ be the domain $\delta t_1<0,q_5>0$; the subset $R_2\subset R_1$ be the union of the curve $(\ref{eta3-l})$ and the domain bounded by the line $(\ref{lucht1-})\cup(0,0)\cup(\ref{eta3-l})$. Then $R_{S_4}({\bf 2}_3)=R_1$. The stratification $\mathfrak{S}_{S_4,t_1,q_5}$ has two strata ${\bf 2}_3$ if $(t_1,q_5)\in R_2$ (including domains $13,14,16,17,19$, $20,22,23,25,27$) and one stratum ${\bf 2}_3$ if $(t_1,q_5)\in R_1\setminus R_2$ (including domains $15,18,21,24,26,28$).

The sets $R_1,R_2$ are contractible. Strata ${\bf 2}_3\in\mathfrak{S}_{S_4,t_1,q_5}$ have two forms ${\mathcal R}_1,{\mathcal R}_2$. A stratum has the form ${\mathcal R}_i$ if: 1) $(t_1,q_5)\in R_i$; 2) the boundary of a stratum of the form ${\mathcal R}_2$ has one singular point; it is a transversal intersection point of two smooth branches of the curve (\ref{psi-2}). Next we repeat the arguments from the proof of Lemma \ref{A2-1-styag}. The set $\beta$ consists of $x$ such that $(t_1,q_5)\in(\ref{eta3-l})$ and $(t_2,S_3)$ satisfies the condition (\ref{usl-1}).

${\bf 2}_4$) Let $R_1$ be the half-plane $\delta t_1>0$; the subset $R_2\subset R_1$ be the union of the curve $(\ref{eta2-m})$ and the domain bounded by the line $(\ref{vozvrat})\cup(\ref{eta2-m})\cup(\ref{spoint6})\cup(\ref{gamma3+-m})\cup(\ref{spoint5})\cup(\ref{gamma1+part})$; the subset $R_3\subset R_1$ be the union of the curve $(\ref{eta2-r})$ and the domain bounded by the line $(\ref{gam1+t})\cup(\ref{vozvrat})\cup(\ref{eta2-r})$. Then $R_{S_4}({\bf 2}_4)=R_1$. The stratification $\mathfrak{S}_{S_4,t_1,q_5}$ has two strata ${\bf 2}_4$ if $(t_1,q_5)$ belongs to $R_2$ (including domains $61,64-67$) or $R_3$ (including domains $38,44,45,49$); one stratum ${\bf 2}_4$ if $(t_1,q_5)\in R_1\setminus (R_2\cup R_3)$ (including domains $29-37,39-43,46-48,50-60,62,63,68-93$).

The sets $R_1 - R_3$ are contractible. Strata ${\bf 2}_4\in\mathfrak{S}_{S_4,t_1,q_5}$ have three forms ${\mathcal R}_1,{\mathcal R}_2,{\mathcal R}_3$. A stratum has the form ${\mathcal R}_i$ if: 1) $(t_1,q_5)\in R_i$; 2) strata of the form ${\mathcal R}_2$ are bounded; 3) if a stratum of the form ${\mathcal R}_3$ is unbounded, then its boundary is connected and contains unbounded interval of the line (\ref{crit-tochki}). The set $\beta$ in this case consists of $x$ such that $(t_1,q_5)$ belongs to the closure of the union $(\ref{eta2-m})\cup(\ref{eta2-r})$ and $(t_2,S_3)$ satisfies the condition (\ref{usl-1}). $\Box$

\medskip
Lemmas \ref{A2-1-styag} -- \ref{A2-3-styag} imply:

\medskip
\predlo\label{styg-A_1^2} {\it The strata ${\bf 2}_1 - {\bf 2}_4$ of $\mathfrak{S}_{S_4},S_4\neq0$ are contractible.}

\medskip
{\sc 2) Strata of type $A_1^4$.}
\medskip

\lemma\label{A4-1-styag} {\it The strata ${\bf 4}_1,{\bf 4}_2,{\bf 4}_4,{\bf 4}_7\in\mathfrak{S}_{S_4},S_4\neq0$ are contractible.}

{\sc Proof.} ${\bf 4}_1$) Let $R_1$ be the domain in $\mathbb{R}^2$ that is bounded by the curves $\gamma_1^-$ and $t_1=0$; the set $R_2$ be the union of the curve $\xi_1$ and the domain that is bounded by $\xi_1$ and does not contain the origin; the subset $R_3\subset R_1$ be the union of the curve $\xi_4^+$ and the domain bounded by the line $\xi_4^+\cup(\ref{spoint1})\cup\gamma_3^-$. Then $R_{S_4}({\bf 4}_{1})=R_1\cup R_2$ that is the half-plane $\delta t_1<0$. The stratification $\mathfrak{S}_{S_4,t_1,q_5}$ has three strata ${\bf 4}_1$ if $(t_1,q_5)\in R_2\cap R_3$ (including domain $22$); two strata ${\bf 4}_1$ if $(t_1,q_5)$ belongs to $(R_1\cap R_2)\setminus R_3$ (including domains $4,10,13,16,19$) or $R_3\setminus R_2$ (including domains $23,24$); one stratum ${\bf 4}_1$ if $(t_1,q_5)$ belongs to $R_2\setminus R_1$ (including domains $1-3,8,9$) or $R_1\setminus (R_2\cup R_3)$ (including domains $5-7,11,12,14,15,17,18,20,21,25-28$).

The sets $R_1$ -- $R_3$ are contractible. Strata ${\bf 4}_1\in\mathfrak{S}_{S_4,t_1,q_5}$ have three forms ${\mathcal R}_1,{\mathcal R}_2,{\mathcal R}_3$. A stratum has the form ${\mathcal R}_i$ if: 1) $(t_1,q_5)\in R_i$; 2) a stratum of the form ${\mathcal R}_3$ is bounded and its boundary has three singular points; 3) other bounded strata ${\bf 4}_1$ have the form ${\mathcal R}_2$; 4) the boundary of any unbounded stratum of the form ${\mathcal R}_2$ has one singular point. Next we repeat the arguments from the proof of Lemma \ref{A2-1-styag}. 

The set $\beta$ has two connected components $\beta_1,\beta_2$. They consists of $x$ such that $(t_1,q_5)$ belongs to the closure of the curve $\xi_4^+$ in the case $\beta_1$ and to $\xi_1$ in the case $\beta_2$; $(t_2,S_3)$ satisfies the condition (\ref{usl-2}) if $(t_1,q_5)\not\equiv(\ref{spoint1})$ and is determined by the formula $(\ref{psi-2})$ for $u$ given by $(\ref{u-36})$ if $(t_1,q_5)\equiv(\ref{spoint1})$. Each of the components is closed and contractible. By Lemmas \ref{lok-ksi1}, \ref{lok-36}, it follows that the intersection of the set $\Sigma^c_{S_4}(x)$ with the stratum ${\bf 4}_{1}\in\mathfrak{S}_{S_4}$ is contractible for any $x\in\beta$. Hence, there are disjoint neighbourhoods $U_1,U_2$ of the sets $\beta_1,\beta_2$ such that the intersections $U_1\cap{\bf 4}_{1},U_2\cap{\bf 4}_{1}$ are contractible and are contracted onto each other in the considered stratum. Parts ${\mathcal R}_3,{\mathcal R}_2$ of this stratum are contracted into $U_1,U_2$, respectively. Part ${\mathcal R}_1$ is contracted into any of neighbourhoods $U_1,U_2$.

${\bf 4}_2$) Let $R_1$ be the domain that is bounded by the curve $\gamma_1^+$ and contains the origin; the subset $R_2\subset R_1$ is the union of the curve (\ref{eta1-l}) and the domain bounded by the line $(\ref{lucht1-})\cup(0,0)\cup(\ref{gam3+m})\cup(\ref{spoint2})\cup(\ref{eta1-l})$; the subset $R_3\subset R_1$ is the union of the curve $\xi_4^-$ and the domain bounded by the line $\xi_4^-\cup(\ref{spoint1})\cup\gamma_3^-$; the subset $R_4\subset R_1$ is the union of curve $(\ref{xi3-t})$ and the domain bounded by the line $(\ref{xi3-t})\cup(0,0)\cup(\ref{tau4+})$. Then $R_{S_4}({\bf 4}_{2})=R_1$. The stratification $\mathfrak{S}_{S_4,t_1,q_5}$ has three strata ${\bf 4}_2$ if $(t_1,q_5)\in R_2\cap R_4$ (including domain $35$); two strata ${\bf 4}_2$ if $(t_1,q_5)$ belongs to $R_2\setminus R_4$ (including domains $13-15,34,36$) or $R_3$ (including domains $19-21$) or $R_4\setminus R_2$ (including domains $31,39$); one stratum ${\bf 4}_2$ if $(t_1,q_5)\in R_1\setminus (R_2\cap R_3\cup R_4)$ (including domains $1-12,16-18,22-30,32,33,37,38,40-49,51-58,61,64-70,73-89$).

The sets $R_1$ -- $R_4$ are contractible. Strata ${\bf 4}_2\in\mathfrak{S}_{S_4,t_1,q_5}$ have four forms ${\mathcal R}_i,i=1,\dots,4$. A stratum has the form ${\mathcal R}_i$ if: 1) $(t_1,q_5)\in R_i$; 2) the boundary of a stratum of the form ${\mathcal R}_4$ has three singular points; one of them is the tangency point of the line (\ref{crit-tochki}) with the curve (\ref{psi-2}); 3) the boundary of a stratum of the form ${\mathcal R}_2$ has an odd number of singular points but this stratum is not a stratum of the form ${\mathcal R}_4$; 4) a stratum of the form ${\mathcal R}_3$ is bounded. The set $\beta$ in this case has three connected components $\beta_1,\beta_2,\beta_3$. They consist of $x$ such that:

$\beta_1$: $(t_1,q_5)$ belongs to the closure of the curve $\xi_4^-$, $(t_2,S_3)$ satisfies the condition (\ref{usl-2}) if $(t_1,q_5)\not\equiv(\ref{spoint1})$ and is determined by the formula $(\ref{psi-2})$ for $u$ given by $(\ref{u-36})$ if $(t_1,q_5)\equiv(\ref{spoint1})$;

$\beta_2$: $(t_1,q_5)$ belongs to the closure of the curve (\ref{eta1-l}), $(t_2,S_3)$ satisfies the condition (\ref{usl-1});

$\beta_3$: $(t_1,q_5)\in(\ref{xi3-t})$ and $(t_2,S_3)$ satisfies the condition (\ref{usl-2}).

Each of the components is closed and contractible. There are pairwise disjoint neighbourhoods $U_1,U_2,U_3$ of $\beta_1,\beta_2,\beta_3$ such that their intersections with the stratum ${\bf 4}_{2}\in\mathfrak{S} _{S_4}$ are contractible. Parts ${\mathcal R}_2,{\mathcal R}_3,{\mathcal R}_4$ of this stratum are contracted into $U_1,U_2,U_3$, respectively. Part ${\mathcal R}_1$ is contracted into any of them.

${\bf 4}_4$) Let $R_1$ be the half-plane $q_5<0$; the set $R_2$ be the union of the curve $\gamma_1^+$ and the domain bounded by $\gamma_1^+$ and the line $(\ref{tau3-l})\cup(\ref{spoint2})\cup(\ref{gam3+t})$; the subset $R_3\subset R_1\cap R_2$ be the union of the curve $\xi_5^-$ and the domain bounded by the line $\xi_5^-\cup(\ref{spoint4})\cup(\ref{gam3+bt})\cup(\ref{spoint5})\cup(\ref{gam1+b})$. Then $R_{S_4}({\bf 4}_{4})=R_1\cup R_2$. The stratification $\mathfrak{S}_{S_4,t_1,q_5}$ has three strata ${\bf 4}_4$ if $(t_1,q_5)\in R_3$ (including domains $56,66,70,89$); two strata ${\bf 4}_4$ if $(t_1,q_5)\in(R_1\cap R_2)\setminus R_3$ (including domains $47-49,53-55,57,58,61,64,65,67-69,74,76,78-82,85-88$); one stratum ${\bf 4}_4$ if $(t_1,q_5)$ belongs to $R_1\setminus R_2$ (including domains $1-12,46,50-52,59,60,62,63,71-73,75,77,83,84,90-93$) or $R_2\setminus R_1$ (including domains $41-45$).

The sets $R_1 - R_3$ are contractible. Strata ${\bf 4}_4\in\mathfrak{S}_{S_4,t_1,q_5}$ have three forms ${\mathcal R}_1,{\mathcal R}_2,{\mathcal R}_3$. A stratum has the form ${\mathcal R}_i$ if: 1) $(t_1,q_5)\in R_i$; 2) a stratum of the form ${\mathcal R}_1$ is unbounded; 3) a stratum of the form ${\mathcal R}_2$ is bounded; it borders an unbounded two-dimensional stratum along a smooth arc of the curve (\ref{psi-2}). The set $\beta$ has two connected components $\beta_1,\beta_2$. They consist of $x$ such that:

$\beta_1$: $(t_1,q_5)\in\gamma_1^+$, $(t_2,S_3)$ belongs to the intersection of two segments with ends (\ref{peresech1}),(\ref{spec-kas}) and (\ref{peresech1}),(\ref{peresech2});

$\beta_2$: $(t_1,q_5)$ belongs to the closure of the curve $\xi_5^-$, $(t_2,S_3)$ satisfies the condition (\ref{usl-2}) if $(t_1,q_5)\not\equiv(\ref{spoint4})$ and is determined by the formula $(\ref{psi-2})$ for $u$ given by $(\ref{u-37})$ if $(t_1,q_5)\equiv(\ref{spoint4})$.

Each of the components is closed and contractible. There are disjoint neighbourhoods $U_1,U_2$ of $\beta_1,\beta_2$ such that their intersections with the stratum ${\bf 4}_{4}\in\mathfrak{S} _{S_4}$ are contractible. Parts ${\mathcal R}_2,{\mathcal R}_3$ of this stratum are contracted into $U_1,U_2$, respectively. Part ${\mathcal R}_1$ is contracted into $U_1$.

${\bf 4}_7$) Let $R_1$ be the domain $\delta t_1>0,q_5>0$; the subset $R_2\subset R_1$ be the union of the curve (\ref{eta1-r}) and the domain bounded by the line $(\ref{eta1-r})\cup(\ref{spoint2})\cup(\ref{gam3+t})$. Then $R_{S_4}({\bf 4}_7)=R_1$. The stratification $\mathfrak{S}_{S_4,t_1,q_5}$ has two strata ${\bf 4}_7$ if $(t_1,q_5)\in R_2$ (including domains $37,38$); one stratum ${\bf 4}_7$ if $(t_1,q_5)\in R_1\setminus R_2$ (including domains $29-36,39-45$). The sets $R_1,R_2$ are contractible. Strata ${\bf 4}_7\in\mathfrak{S}_{S_4,t_1,q_5}$ have two forms ${\mathcal R}_1,{\mathcal R}_2$. A stratum has the form ${\mathcal R}_i$ if: 1) $(t_1,q_5)\in R_i$; 2) a stratum of the form ${\mathcal R}_2$ is bounded. The set $\beta$ consists of $x$ such that $(t_1,q_5)$ belongs to the closure of the curve (\ref{eta1-r}) and $(t_2,S_3)$ satisfies the condition (\ref{usl-1}). $\Box$

\medskip
\lemma\label{A4-3-styag} {\it The strata ${\bf 4}_3,{\bf 4}_5,{\bf 4}_6,{\bf 4}_8,{\bf 4}_9,{\bf 4}_{10}\in\mathfrak{S}_{S_4},S_4\neq0$ are contractible.}

{\sc Proof.} The sets $R_{S_4}({\bf 4}_3),R_{S_4}({\bf 4}_5),R_{S_4}({\bf 4}_6),R_{S_4}({\bf 4}_8),R_{S_4}({\bf 4}_9),R_{S_4}({\bf 4}_{10})$ are contractible do\-mains in $\mathbb{R}^2$. Namely:

$R_{S_4}({\bf 4}_3)$ is bounded by the line $(\ref{eta1-l})\cup(\ref{spoint2})\cup(\ref{gamma2-m})\cup(\ref{spoint6})\cup(\ref{tau2-b})$ and contains the origin; the stratification $\mathfrak{S}_{S_4,t_1,q_5}$ has one stratum ${\bf 4}_3$ if $(t_1,q_5)\in R_{S_4}({\bf 4}_3)$ (including domains $1-15,34-36,39-41,46,47,51-53,73-79,83-85$);

$R_{S_4}({\bf 4}_5)$ is bounded by the curve $\tau_1$ and the line $(\ref{luchq5+})\cup(0,0)\cup(\ref{tau4-})$; the stratification $\mathfrak{S}_{S_4,t_1,q_5}$ has one stratum ${\bf 4}_5$ if $(t_1,q_5)\in R_{S_4}({\bf 4}_5)$ (including domains $3-5,9-11,13-28$);

$R_{S_4}({\bf 4}_6)$ is the domain $\delta t_1<0,q_5>0$; the stratification $\mathfrak{S}_{S_4,t_1,q_5}$ has one stratum ${\bf 4}_6$ if $(t_1,q_5)\in R_{S_4}({\bf 4}_6)$ (including domains $13-28$);

$R_{S_4}({\bf 4}_8)$ lies in the half-plane $\delta t_1>0$ and is bounded by the line $(\ref{luchq5+})\cup(0,0)\cup\eta_3^+\cup(\ref{vozvrat})\cup(\ref{eta3-r})$; the stratification $\mathfrak{S}_{S_4,t_1,q_5}$ has one stratum ${\bf 4}_8$ if $(t_1,q_5)\in R_{S_4}({\bf 4}_8)$ (including domains $29-50$);

$R_{S_4}({\bf 4}_9)$ is the domain $\delta t_1>0,q_5>0$; the stratification $\mathfrak{S}_{S_4,t_1,q_5}$ has one stratum ${\bf 4}_9$ if $(t_1,q_5)\in R_{S_4}({\bf 4}_9)$ (including domains $29-45$);

$R_{S_4}({\bf 4}_{10})$ lies in the half-plane $\delta t_1>0$ and is bounded by the line $(\ref{luchq5-})\cup(0,0)\cup(\ref{gam3+m})\cup(\ref{spoint2})\cup(\ref{tau3-r})$; the stratification $\mathfrak{S}_{S_4,t_1,q_5}$ has one stratum ${\bf 4}_{10}$ if $(t_1,q_5)\in R_{S_4}({\bf 4}_{10})$ (including domains $39-42,45-93$).

Next we repeat the arguments from the proof of Lemma \ref{A2-2-styag}. $\Box$

\medskip
\lemma\label{A4-11-styag} {\it The stratum ${\bf 4}_{11}\in\mathfrak{S}_{S_4},S_4\neq0$ is homotopy equivalent to $S^1$.}

{\sc Proof.} Consider the curves
\begin{equation}
\xi_5^+: q_5>-\frac{3}{2}\sqrt[5]{\frac{9S_4^2}{8}}\quad (\mbox{between (\ref{spoint4}) and (\ref{trper-xi5-g1})}),
\label{xi5+l}
\end{equation}
\begin{equation}
\xi_5^+: q_5<-\frac{3}{2}\sqrt[5]{\frac{9S_4^2}{8}}\quad (\mbox{bounded by (\ref{trper-xi5-g1})}),
\label{xi5+r}
\end{equation}
\begin{equation}
\gamma_1^+: -\sqrt[5]{12S_4^2}<q_5<-\frac{3}{2}\sqrt[5]{\frac{9S_4^2}{8}}\quad (\mbox{between (\ref{spoint5}) and (\ref{trper-xi5-g1})}),
\label{gam1+bm}
\end{equation}
\begin{equation}
\gamma_1^+: q_5>-\frac{3}{2}\sqrt[5]{\frac{9S_4^2}{8}}\quad (\mbox{bounded by (\ref{trper-xi5-g1})}).
\label{gam1+tm}
\end{equation}
Let $R_1$ be the domain $\delta t_1>0,q_5<0$; the set $R_2$ be the union of the curve $(\ref{xi5+r})\cup(\ref{trper-xi5-g1})\cup(\ref{gam1+tm})$ and the domain bounded by this curve and the line $(\ref{eta1-r})\cup(\ref{spoint2})\cup(\ref{gamma2-m})\cup(\ref{spoint6})\cup(\ref{tau2-r})$; the subset $R_3\subset R_1$ be the union of the curve $\xi_2$ and the domain that is bounded by $\xi_2$ and does not contain the origin; the subset $R_4\subset R_1$ be the union of the curve $(\ref{xi3-b})$ and the domain bounded by the line $(\ref{luchq5-})\cup(0,0)\cup(\ref{xi3-b})$; the subset $R_5\subset R_1$ be the union of the curve $(\ref{eta3-r})\cup(\ref{vozvrat})\cup\eta_3^+$ and the domain bounded by this curve and the line $(\ref{lucht1+})\cup(0,0)\cup(\ref{eta3-r})\cup(\ref{vozvrat})\cup\eta_3^+$; the subset $R_6\subset R_1\cap R_2$ be the union of the curve $(\ref{xi5+l})\cup(\ref{trper-xi5-g1})\cup(\ref{gam1+bm})$ and the domain bounded by the line $(\ref{spoint4})\cup(\ref{xi5+l})\cup(\ref{trper-xi5-g1})\cup(\ref{gam1+bm})\cup(\ref{spoint5})\cup(\ref{gam3+bt})$. Then $R_{S_4}({\bf 4}_{11})=R_1\cup R_2$. The stratification $\mathfrak{S}_{S_4,t_1,q_5}$ has four strata ${\bf 4}_{11}$ if $(t_1,q_5)\equiv(\ref{trper-xi5-g1})$; three strata ${\bf 4}_{11}$ if $(t_1,q_5)$ belongs to $R_2\cap R_3$ (including domain $72$) or $R_2\cap R_5$ (including domains $48,49$) or $R_6$ (including domains $57,67$); two strata ${\bf 4}_{11}$ if $(t_1,q_5)$ belongs to $R_5\setminus R_2$ (including domains $46,47,50$) or $R_4$ (including domains $51,73,74,77,78,80,83-86$) or $R_2\cap R_1$ (including domains $54-56,58,61,64-66,68-71$) or $R_3\setminus R_2$ (including domains $60,63,91,93$); one stratum ${\bf 4}_{11}$ if $(t_1,q_5)$ belongs to $R_2\setminus R_1$ (including domains $37,38,42-45$) or $R_1\setminus (R_2\cup R_3\cup R_4\cup R_5)$ (including domains $52,53,59,62,75,76,79,81,82,87-90,92$).

The sets $R_1 - R_6$ are contractible. Strata ${\bf 4}_{11}\in\mathfrak{S}_{S_4,t_1,q_5}$ have six forms ${\mathcal R}_i,i=1,\dots,6$. A stratum has the form ${\mathcal R}_i$ if: 1) $(t_1,q_5)\in R_i$; 2) a stratum of the form ${\mathcal R}_6$ is bounded; its boundary has three singular points; each of them is a transversal intersection point of two smooth branches of $\Sigma_{S_4,t_1,q_5}$; 3) a stratum of the form ${\mathcal R}_2$ is bounded but is not a stratum of the form ${\mathcal R}_6$; 4) a stratum of the form ${\mathcal R}_1$ is unbounded; its boundary contains an unbounded interval of the line (\ref{crit-tochki}); 5) a stratum of any of the forms ${\mathcal R}_3,{\mathcal R}_4,{\mathcal R}_5$ is unbounded but it is not a stratum of the form ${\mathcal R}_1$.

The stratification $\mathfrak{S}_{S_4,t_1,q_5}$ has one stratum ${\bf 4}_{11}$ of the form ${\mathcal R}_{i}$ for any $(t_1,q_5)\in R_i$. Let $\widetilde{R_i}$ be $R_i$ if $i\neq1$ and be the complement to the closure of the curve (\ref{gamma1+bv}) in $R_1$ if $i=1$. Then each stratum of the form ${\mathcal R}_{i}$ is contracted to a point $(t_2,S_3)$ that continuously depends on $(t_1,q_5)\in\widetilde{R_i}$. The set of points $x$ such that $(t_1,q_5)\in\widetilde{R_i}$ and $(t_2,S_3)$ is a point of a stratum ${\bf 4}_{11}$ of the form ${\mathcal R}_{i}$ is the total space of a locally trivial bundle. This space is contractible. It is called the part ${\mathcal R}_{i}$ of the stratum ${\bf 4}_{11}\in\mathfrak{S}_{S_4}$.

The set $\beta$ has four connected components $\beta_1=\beta_1^1\cup\beta_1^2,\beta_2=\beta_2^1\cup\beta_2^2\cup\beta_2^3,\beta_3,\beta_4$. They consist of $x$ such that:

$\beta_1^1$: $(t_1,q_5)$ belongs to the closure of the curve (\ref{gamma1+bv}), $(t_2,S_3)\equiv(\ref{spec2})$;

$\beta_1^2$: $(t_1,q_5)$ belongs to the union $(\ref{eta3-r})\cup\eta_3^+$, $(t_2,S_3)$ satisfies the condition (\ref{usl-1});

$\beta_2^1$: $(t_1,q_5)$ belongs to the closure of the union $(\ref{gam1+bm})\cup(\ref{gam1+tm})$, $(t_2,S_3)\equiv(\ref{spec-kas})$;

$\beta_2^2$: $(t_1,q_5)$ belongs to the closure of the curve $\xi_5^+$, $(t_2,S_3)$ satisfies the condition (\ref{usl-2}) if $(t_1,q_5)\not\equiv(\ref{spoint4})$ and is determined by the formula $(\ref{psi-2})$ for $u$ given by $(\ref{u-37})$ if $(t_1,q_5)\equiv(\ref{spoint4})$;

$\beta_2^3$: $(t_1,q_5)\equiv(\ref{trper-xi5-g1})$, $(t_2,S_3)$ belongs to the segment with ends at (\ref{spec-kas}) and the point that satisfies the condition (\ref{usl-2}); 

$\beta_3$: $(t_1,q_5)\in\xi_2$, $(t_2,S_3)$ satisfies the condition (\ref{usl-2});

$\beta_4$: $(t_1,q_5)\in(\ref{xi3-b})$, $(t_2,S_3)$ satisfies the condition (\ref{usl-2}).

Each of the components is closed and contractible. By Lemmas \ref{lok-eta2}, \ref{lok-ksi1}, \ref{lok-36}, \ref{lok-gamma+}, \ref{lok-int}, \ref{lok-00} it follows that the intersection of the set $\Sigma^c_{S_4}(x)$ with the stratum ${\bf 4}_{11}\in\mathfrak{S}_{S_4}$ is contractible for any $x\in\beta\setminus \beta_1^1$. If $x\in\beta_1^1$, then $\Sigma^c_{S_4}(x)$ does not intersect parts ${\mathcal R}_{3},{\mathcal R}_{4},{\mathcal R}_{6}$ of this stratum. The intersections of $\Sigma^c_{S_4}(x)$ with parts ${\mathcal R}_{2},{\mathcal R}_{5}$ are not empty if and only if $(t_1,q_5)\equiv(\ref{vozvrat})$. These intersections are contractible. Finally, the intersection of $\Sigma^c_{S_4}(x)$ with part ${\mathcal R}_{1}$ is homotopy equivalent to $S^1$. Hence, there are pairwise disjoint neighbourhoods $U_1,U_2,U_3,U_4$ of components $\beta_1,\beta_2,\beta_3,\beta_4$ such that the intersection $U_1\cap{\bf 4}_{11}$ is homotopy equivalent to $S^1$ and all other intersections $U_i\cap{\bf 4}_{11},i\neq1$ are contractible. Parts ${\mathcal R}_3,{\mathcal R}_4$ of the stratum ${\bf 4}_{11}$ are contracted into $U_3,U_4$, respectively. The union ${\mathcal R}_2\cup{\mathcal R}_6$ is contracted into $U_2$, and ${\mathcal R}_1\cup{\mathcal R}_5$ is contracted into $U_1$. $\Box$

\medskip
Lemmas \ref{A4-1-styag} -- \ref{A4-11-styag} imply:

\medskip
\predlo\label{styg-A_1^4} {\it The strata ${\bf 4}_1 - {\bf 4}_{10}$ of the stratification $\mathfrak{S}_{S_4},S_4\neq0$ are contractible. The stratum ${\bf 4}_{11}$ is homotopy equivalent to $S^1$.}

\medskip
{\sc 3) Strata of type $A_1^6$.}
\medskip

\lemma\label{A6-1-styag} {\it The strata ${\bf 6}_1,{\bf 6}_2,{\bf 6}_5,{\bf 6}_6,{\bf 6}_9,{\bf 6}_{10},{\bf 6}_{11},{\bf 6}_{13},{\bf 6}_{15},{\bf 6}_{16}\in\mathfrak{S}_{S_4},S_4\neq0$ are contrac\-tible.}

{\sc Proof.} The sets $
R_{S_4}({\bf 6}_1), R_{S_4}({\bf 6}_2), R_{S_4}({\bf 6}_5), R_{S_4}({\bf 6}_6), R_{S_4}({\bf 6}_9), R_{S_4}({\bf 6}_{10}),R_{S_4}({\bf 6}_{11}),R_{S_4}({\bf 6}_{13})$, $R_{S_4}({\bf 6}_{15}), R_{S_4}({\bf 6}_{16})
$
are contractible domains in $\mathbb{R}^2$. Namely:

$R_{S_4}({\bf 6}_1)$ is bounded by the lines $\xi_4^+\cup(\ref{spoint1})\cup(\ref{zeta1-r})$ and $(\ref{gam2+t})\cup(\ref{spoint2})\cup(\ref{tau3-l})$; the stratification $\mathfrak{S}_{S_4,t_1,q_5}$ has one stratum ${\bf 6}_1$ if $(t_1,q_5)\in R_{S_4}({\bf 6}_1)$ (including domains $1-24,30-32,34-36,39,40,46,51,52$,$73,75,77,83,84$);

$R_{S_4}({\bf 6}_2)$ lies in the half-plane $\delta t_1<0$ and is bounded by the curve $\tau_1$; the stratification $\mathfrak{S}_{S_4,t_1,q_5}$ has one stratum ${\bf 6}_2$ if $(t_1,q_5)\in R_{S_4}({\bf 6}_2)$ (including domains $1,2,8$);

$R_{S_4}({\bf 6}_5)$ lies in the half-plane $\delta t_1<0$ and is bounded by the curve $\zeta_2$; it coincides with domain $1$; the stratification $\mathfrak{S}_{S_4,t_1,q_5}$ has one stratum ${\bf 6}_5$ if $(t_1,q_5)\in R_{S_4}({\bf 6}_5)$;

$R_{S_4}({\bf 6}_6)$ lies in the half-plane $q_5<0$ and is bounded by the line $(\ref{eta2-l})\cup(\ref{spoint6})\cup(\ref{gamma2+b})$; the stratification $\mathfrak{S}_{S_4,t_1,q_5}$ has one stratum ${\bf 6}_6$ if $(t_1,q_5)\in R_{S_4}({\bf 6}_6)$ (including domains $1-7,73-81,83-87$);

$R_{S_4}({\bf 6}_9)$ lies in the half-plane $\delta t_1<0$ and is bounded by the line $\xi_4^-\cup(\ref{spoint1})\cup(\ref{zeta1-l})$; the stratification $\mathfrak{S}_{S_4,t_1,q_5}$ has one stratum ${\bf 6}_9$ if $(t_1,q_5)\in R_{S_4}({\bf 6}_9)$ (including domains $19-26$);

$R_{S_4}({\bf 6}_{10})$ lies in the half-plane $\delta t_1>0$ and is bounded by the line $(\ref{luchq5-})\cup(0,0)\cup(\ref{xi3-t})$; the stratification $\mathfrak{S}_{S_4,t_1,q_5}$ has one stratum ${\bf 6}_{10}$ if $(t_1,q_5)\in R_{S_4}({\bf 6}_{10})$ (including domains $31-33,35-93$);

$R_{S_4}({\bf 6}_{11})$ lies in the half-plane $\delta t_1>0$ and is bounded by the line $(\ref{tau3-r})\cup(\ref{spoint2})\cup(\ref{gam2+t})$; the stratification $\mathfrak{S}_{S_4,t_1,q_5}$ has one stratum ${\bf 6}_{11}$ if $(t_1,q_5)\in R_{S_4}({\bf 6}_{11})$ (including domains $33,37,38$, $43,44$);

$R_{S_4}({\bf 6}_{13})$ lies in the half-plane $\delta t_1>0$ and is bounded by the line $(\ref{gam3+bb})\cup(\ref{spoint6})\cup(\ref{tau2-r})$; the stratification $\mathfrak{S}_{S_4,t_1,q_5}$ has one stratum ${\bf 6}_{13}$ if $(t_1,q_5)\in R_{S_4}({\bf 6}_{13})$ (including domains $77-82,84-93$);

$R_{S_4}({\bf 6}_{15})$ lies in the half-plane $\delta t_1>0$ and is bounded by the curve $\xi_2$; the stratification $\mathfrak{S}_{S_4,t_1,q_5}$ has one stratum ${\bf 6}_{15}$ if $(t_1,q_5)\in R_{S_4}({\bf 6}_{15})$ (including domains $60,63,72,91,93$);

$R_{S_4}({\bf 6}_{16})$ lies in the half-plane $\delta t_1>0$ and is bounded by the curve $\gamma_2^-$; the stratification $\mathfrak{S}_{S_4,t_1,q_5}$ has one stratum ${\bf 6}_{16}$ if $(t_1,q_5)\in R_{S_4}({\bf 6}_{16})$ (including domains $92,93$).

Next we repeat the arguments from the proof of Lemma \ref{A2-2-styag}. $\Box$

\medskip
\lemma\label{A6-3-styag} {\it The strata ${\bf 6}_3,{\bf 6}_4,{\bf 6}_7,{\bf 6}_8,{\bf 6}_{12},{\bf 6}_{14}\in\mathfrak{S}_{S_4},S_4\neq0$ are contractible.}

{\sc Proof.} ${\bf 6}_{3}$) Let $R_1$ be the domain $\delta t_1<0,q_5<0$; the set $R_2$ be the union of the curve $\gamma_1^-$ and the domain bounded by $\gamma_1^-$ and the line $(\ref{luchq5-})\cup(0,0)\cup(\ref{eta3-l})$. Then $R_{S_4}({\bf 6}_{3})=R_1\cup R_2$. The stratification $\mathfrak{S}_{S_4,t_1,q_5}$ has two strata ${\bf 6}_{3}$ if $(t_1,q_5)\in R_1\cap R_2$ (including domains $4-7,10-12$); one stratum ${\bf 6}_{3}$ if $(t_1,q_5)$ belongs to $R_1\setminus R_2$ (including domains $1-3,8,9$) or $R_2\setminus R_1$ (including domains $13,14,16,17,19,20,22,23,25,27$).

The sets $R_1,R_2$ are contractible. Strata ${\bf 6}_3\in\mathfrak{S}_{S_4,t_1,q_5}$ have two forms ${\mathcal R}_1,{\mathcal R}_2$. A stratum has the form ${\mathcal R}_i$ if: 1) $(t_1,q_5)\in R_i$; 2) the boundary of a stratum of the form ${\mathcal R}_1$ contains an unbounded interval of the line (\ref{crit-tochki}). The set $\beta$ consists of $x$ such that $(t_1,q_5)\in\gamma_1^-$ and $(t_2,S_3)$ belongs to the closed infinite interval of the line $(\ref{psi-4})$ that is bounded by the point $(\ref{spec-kas})$ and does not contain $(\ref{spec2})$. It is closed and contractible. Next we repeat the arguments from the proof of Lemma \ref{A2-1-styag}.

${\bf 6}_{4}$) Let $R_1$ be the domain $\delta t_1<0,q_5<0$; the set $R_2$ be the union of the curve $\gamma_1^-$ and the domain bounded by the curves $\gamma_1^-$ and $\xi_1$; the set $R_3$ be the union of $\gamma_1^-$ and the domain bounded by $\gamma_1^-$ and the line $(\ref{lucht1-})\cup(0,0)\cup(\ref{xi3-b})$. Then $R_{S_4}({\bf 6}_{4})=R_1\cup R_2\cup R_3$. The stratification $\mathfrak{S}_{S_4,t_1,q_5}$ has three strata ${\bf 6}_{4}$ if $(t_1,q_5)\in R_1\cap R_2\cap R_3$ (including domains $4,10$); two strata ${\bf 6}_{4}$ if $(t_1,q_5)\in R_1\cap R_3$ (including domains $5-7,11,12$); one stratum ${\bf 6}_{4}$ if $(t_1,q_5)$ belongs to $R_1\setminus R_3$ (including domains $1-3,8,9$) or $R_2\setminus R_1$ (including domains $13,16,19,22$) or $R_3\setminus R_1$ (including domains $51,73,74,77,78,80,83-86$).

The sets $R_1 - R_3$ are contractible. Strata ${\bf 6}_4\in\mathfrak{S}_{S_4,t_1,q_5}$ have three forms ${\mathcal R}_1,{\mathcal R}_2,{\mathcal R}_3$. A stratum has the form ${\mathcal R}_i$ if: 1) $(t_1,q_5)\in R_i$; 2) the boundary of a stratum of the form ${\mathcal R}_1$ contains an unbounded interval of the line (\ref{crit-tochki}); 3) a stratum of the form ${\mathcal R}_2$ is bounded. The set $\beta$ consists of $x$ such that $(t_1,q_5)\in\gamma_1^-$ and $(t_2,S_3)$ belongs to the closed infinite interval of the line $(\ref{psi-4})$ that is bounded by the point $(\ref{spec-kas})$ and contains $(\ref{spec2})$. It is closed and contractible.

${\bf 6}_{7}$) Let $R_1$ be the domain that lies in the half-plane $\delta t_1>0$ and is bounded by the line $(\ref{luchq5-})\cup(0,0)\cup(\ref{tau4+})$; the set $R_2$ be the union of the curve $\gamma_1^+$ and the domain bounded by $\gamma_1^+$ and the line $(\ref{lucht1+})\cup(0,0)\cup(\ref{tau4-})$. Then $R_{S_4}({\bf 6}_{7})=R_1\cup R_2$. The stratification $\mathfrak{S}_{S_4,t_1,q_5}$ has two strata ${\bf 6}_{7}$ if $(t_1,q_5)\in R_1\cap R_2$ (including domains $46-49,51-58,61,64-70,73-89$); one stratum ${\bf 6}_{7}$ if $(t_1,q_5)$ belongs to $R_1\setminus R_2$ (including domains $32,33,36-38,40-45,50,59,60,62,63,71,72,90-93$) or $R_2\setminus R_1$ (including domains $6,7,12$).

The sets $R_1,R_2$ are contractible. Strata ${\bf 6}_7\in\mathfrak{S}_{S_4,t_1,q_5}$ have two forms ${\mathcal R}_1,{\mathcal R}_2$. A stratum has the form ${\mathcal R}_i$ if: 1) $(t_1,q_5)\in R_i$; 2) a stratum of the form ${\mathcal R}_1$ is either bounded or its boundary contains an unbounded interval of the line (\ref{crit-tochki}). The set $\beta$ consists of $x$ such that $(t_1,q_5)\in\gamma_1^+$ and $(t_2,S_3)$ belongs to the closed infinite interval of the line $(\ref{psi-4})$ that is bounded by the point $(\ref{peresech1})$ and does not contain $(\ref{spec2})$. It is closed and contractible.

${\bf 6}_{8}$) Let $R_1$ be the domain that lies in the half-plane $q_5<0$ and is bounded by the line $\xi_5^+\cup(\ref{spoint4})\cup(\ref{dzeta3-l})$; the set $R_2$ be the union of the curve (\ref{gam1+b}) and the domain that lies in the half-plane $\delta t_1>0$ and is bounded by the line $(\ref{tau2-b})\cup(\ref{spoint6})\cup(\ref{gamma3+-m})\cup(\ref{spoint5})\cup(\ref{gam1+b})$. Then $R_{S_4}({\bf 6}_{8})=R_1\cup R_2$. The stratification $\mathfrak{S}_{S_4,t_1,q_5}$ has two strata ${\bf 6}_{8}$ if $(t_1,q_5)\in R_1\cap R_2$ (including domains $69,70,86-89$); one stratum ${\bf 6}_{8}$ if $(t_1,q_5)$ belongs to $R_1\setminus R_2$ (including domains $7,55-57,65-67,71,72,83-85,90-93$) or $R_2\setminus R_1$ (including domains $68,80-82$).

The sets $R_1,R_2$ are contractible. Strata ${\bf 6}_8\in\mathfrak{S}_{S_4,t_1,q_5}$ have two forms ${\mathcal R}_1,{\mathcal R}_2$. A stratum has the form ${\mathcal R}_i$ if: 1) $(t_1,q_5)\in R_i$; 2) at least one of singular points of the boundary of a stratum of the form ${\mathcal R}_1$ is a semicubical cusp. The set $\beta$ consists of $x$ such that $(t_1,q_5)$ belongs to the closure of the curve (\ref{gam1+b}) and $(t_2,S_3)$ belongs to the segment with ends $(\ref{peresech2}),(\ref{spec-kas})$. It is closed and contractible.

${\bf 6}_{12}$) Let $R_1$ be the domain $\delta t_1>0,q_5<0$; the set $R_2$ be the union of the curve (\ref{gam1+t}) and the domain that lies in the half-plane $\delta t_1>0$ and is bounded by the line $(\ref{eta2-r})\cup(\ref{vozvrat})\cup(\ref{gam1+t})$; the subset $R_3\subset R_1$ be the union of the curve $\gamma_1^+$ and the domain bounded by $\gamma_1^+$ and the line $(\ref{luchq5-})\cup(0,0)\cup(\ref{lucht1+})$. Then $R_{S_4}({\bf 6}_{12})=R_1\cup R_2$. The stratification $\mathfrak{S}_{S_4,t_1,q_5}$ has three strata ${\bf 6}_{12}$ if $(t_1,q_5)\in R_2\cap R_3$ (including domain $49$); two strata ${\bf 6}_{12}$ if $(t_1,q_5)\in R_3\setminus R_2$ (including domains $46-48,51-58,61,64-70,73-89$); one stratum ${\bf 6}_{12}$ if $(t_1,q_5)$ belongs to $R_2\setminus R_3$ (including domains $38,44,45$) or $R_1\setminus R_3$ (including domains $50,59,60,62,63,71,72,90-93$).

The sets $R_1 - R_3$ are contractible. Strata ${\bf 6}_{12}\in\mathfrak{S}_{S_4,t_1,q_5}$ have three forms ${\mathcal R}_1,{\mathcal R}_2,{\mathcal R}_3$. A stratum has the form ${\mathcal R}_i$ if: 1) $(t_1,q_5)\in R_i$; 2) a stratum of the form ${\mathcal R}_2$ is bounded; 3) a stratum of the form ${\mathcal R}_1$ is unbounded; one of smooth curves bounding it and an infinite interval of the line (\ref{crit-tochki}) bound a domain in the intersection of a stratum ${\bf 4}_{11}$ with a neighbourhood of infinity. The set $\beta$ is the union of two sets $\beta_1,\beta_2$. The set $\beta_1$ consists of $x$ such that $(t_1,q_5)\in\gamma_1^+$ and $(t_2,S_3)$ belongs to the closed infinite interval of the line $(\ref{psi-4})$ that is bounded by the point (\ref{spec2}) and does not contain (\ref{spec-kas}). The set $\beta_2$ consists of $x$ such that $(t_1,q_5)\in(\ref{gam1+t})$ and $(t_2,S_3)$ belongs to the segment with ends (\ref{peresech2}),(\ref{spec2}). The set $\beta$ is closed and contractible.

${\bf 6}_{14}$) Let $R_1$ be the domain bounded by the curve $\gamma_2^-$ and the line $\zeta_3^+\cup(\ref{vozvrat})\cup(\ref{spoint4})\cup(\ref{dzeta3-r})\cup\xi_5^-$; the set $R_2$ be the union of the curve (\ref{gamma1+bv}) and the domain that lies in the half-plane $\delta t_1>0$ and is bounded by the line $(\ref{gamma2+b})\cup(\ref{spoint6})\cup(\ref{eta2-m})\cup(\ref{vozvrat})\cup(\ref{gamma1+bv})$. Then $R_{S_4}({\bf 6}_{14})=R_1\cup R_2$. The stratification $\mathfrak{S}_{S_4,t_1,q_5}$ has two strata ${\bf 6}_{14}$ if $(t_1,q_5)\in R_1\cap R_2$ (including domains $61,66,67,70,89$); one stratum ${\bf 6}_{14}$ if $(t_1,q_5)$ belongs to $R_1\setminus R_2$ (including domains $56-58,62,63,71,74,90,91$) or $R_2\setminus R_1$ (including domains $64,65,68,69,82,88$).

The sets $R_1,R_2$ are contractible. Strata ${\bf 6}_{14}\in\mathfrak{S}_{S_4,t_1,q_5}$ have two forms ${\mathcal R}_1,{\mathcal R}_2$. A stratum has the form ${\mathcal R}_i$ if: 1) $(t_1,q_5)\in R_i$; 2) a stratum of the form ${\mathcal R}_1$ belongs to the connected component of the complement to the curve (\ref{psi-2}) that is diffeomorphic to the bounded connected component of the complement to the swallowtail section by a plane transversally intersecting the self-intersection line of this swallowtail near its vertex. The set $\beta$ consists of $x$ such that $(t_1,q_5)$ belongs to the closure of the curve (\ref{gamma1+bv}) and $(t_2,S_3)$ belongs to the intersection of two segments with ends (\ref{spec-kas}),(\ref{spec2}) and (\ref{peresech2}),(\ref{spec2}). It is closed and contractible. $\Box$

\medskip
Lemmas \ref{A6-1-styag} and \ref{A6-3-styag} imply:

\medskip
\predlo\label{styg-A_1^6} {\it The strata ${\bf 6}_{1} - {\bf 6}_{16}$ of the stratification $\mathfrak{S}_{S_4},S_4\neq0$ are contractible.}

\section{The bifurcation diagram $B_{0}$}

The equations of curves that compose the bifurcation diagram $B_{S_4}$ depend continuously on the parameter $S_4$. The limit position of these curves and domains $1-93$  in Fig. \ref{razbienie} as $S_4\rightarrow0$ can be found using Lemmas \ref{samoperesechS4-2} -- \ref{samoperesechS4-5}. Simple calculations show that limit curves are smooth everywhere except the origin. Each of them is part of the union of the curves given by equations
\begin{equation}
q_5^2=12|t_1|,\quad q_5=\sqrt{\frac32|t_1|},\quad q_5=0,\quad t_1=0.
\label{predel}
\end{equation}

\medskip
\lemma\label{povedenie} {\it Ten domains $1,8,22,25,28,29,33,38,60,93$ are turned into connected components of the complement to the union of curves {\rm(\ref{predel})}. All other domains from the list $1-93$ disappear as $S_4\rightarrow0$.}

\medskip
The bifurcation diagram $B_{0}$ is shown in Fig. \ref{razbienie0}. The origin divides each of curves (\ref{predel}) into two arcs. Each arc belongs to the limit of several curves from Lemmas \ref{samoperesechS4-2} -- \ref{samoperesechS4-5}. The notations of two such curves are indicated in the figure near the corresponding arcs. The limit positions of other curves are determined by the comparing of Fig. \ref{razbienie} and \ref{razbienie0} using Squeeze Theorem. The connected components of the complement to $B_0$ are denoted by the same numbers as the corresponding domains in Fig. \ref{razbienie}.

\begin{figure}[h]
\begin{center}
\includegraphics[width=9cm]{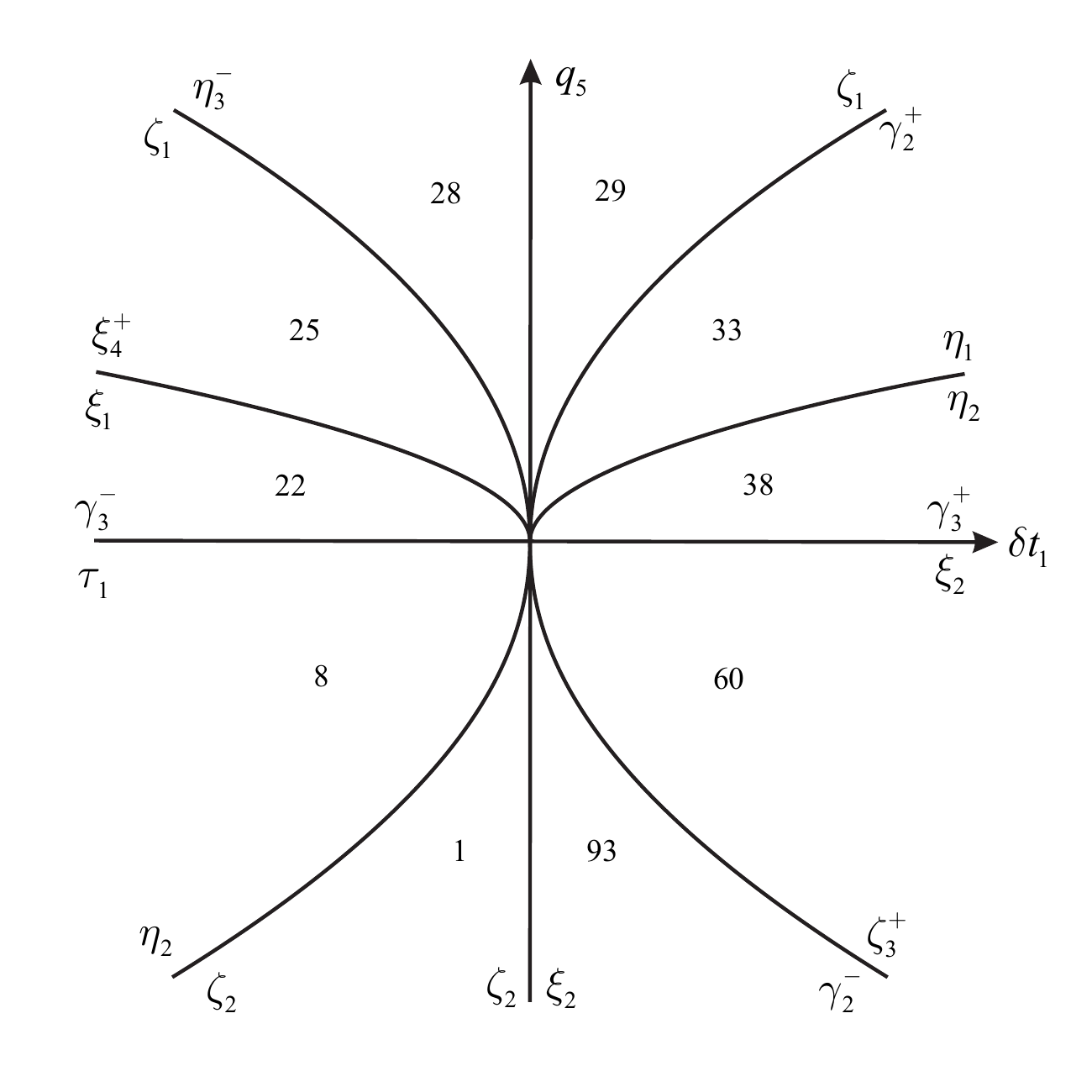}
\caption{The bifurcation diagram $B_{0}$.}
\label{razbienie0}
\end{center}
\end{figure}

\section{Type $A_1^k$ strata of the stratification $\mathfrak{S}$}

Three-dimensional space $S_4=t_1=0$ belongs to the skeleton $\Sigma_0$. By analogy with Lemma \ref{sech}, all sections of $\Sigma_0$ by hyperplanes $t_1=\mathrm{const}>0$ are diffeomorphic.

\medskip
\lemma\label{simS4} {\it The skeleton $\Sigma_0$ is invariant under diffeomorphism}
\begin{equation}
(t_1,t_2,S_3,S_4,q_5)\mapsto (-t_1,t_2,S_3-4q_5t_1,S_4,q_5).
\label{diffeo-S4-0}
\end{equation}

{\sc Proof.} Consider the generating family
$
\Phi(v,t,q)+\lambda_1v_1+\lambda_2v_2
$
of the Lagrangian map $F$ from Lemma \ref{proizv-psi}. The diffeomorphism
\begin{equation}
(v_1,v_2,t_1,t_2,S_3,S_4,q_5,\lambda_1,\lambda_2)\mapsto
(v_1+2t_1,v_2,-t_1,t_2,S_3-4q_5t_1,S_4,q_5,\lambda_1,\lambda_2)
\label{diffeo-S4-0-1}
\end{equation}
takes this family to
$\Phi(v,t,q)+\lambda_1v_1+(\lambda_2+2t_1S_4)v_2+2t_1(\lambda_1-2t_1)$.
The restrictions of both families to the hyperplane $S_4=0$ in the space of parameters $(t_1,t_2,S_3,S_4,q_5,\lambda_1,\lambda_2)$ differ by a function that depends only on the parameter. The diffeomorphism (\ref{diffeo-S4-0-1}) preserves the subspace $\lambda_1=\lambda_2=0$. Hence, the section of the caustic of the map $F$ by the subspace $S_4=\lambda_1=\lambda_2=0$ is invariant under (\ref{diffeo-S4-0}). $\Box$

\medskip
Thus all sections of the skeleton $\Sigma_0$ by hyperplanes $t_1=\mathrm{const}\neq0$ are diffeomorphic.

\medskip
\lemma\label{samoperesechS4} {\it Let $S_4=0,t_1\neq0$. Then the following statements hold.

{\rm 1)} The skeleton $\Sigma_{0,t_1,q_5}$ is symmetric with respect to the line $t_2=0$. If $q_5\neq0$, then $\Sigma_{0,t_1,q_5}$ consists of the curve {\rm(\ref{psi-2})} and two parallel lines $S_3=0$, $S_3=4q_5t_1$. If $q_5=0$, then $\Sigma_{0,t_1,q_5}$ consists of parabola $S_3=9\delta t_2^2$ and its tangent $S_3=0$.

{\rm 2)} If $q_5\neq0$, then the curve {\rm(\ref{psi-2})} is tangent to the line $S_3=0$ at $t_2=0$. The tangency is simple if $q_5^2\neq12\delta t_1$; in the case $q_5^2=12\delta t_1$, the curve {\rm(\ref{psi-2})} has a singularity of type $4/3$ at the point of tangency.

{\rm 3)} If
$$
q_5\in\left(-2\sqrt{3|t_1|},0
\right)\cup\left(2\sqrt{3|t_1|},+\infty\right),
$$
then the curve {\rm(\ref{psi-2})} is smooth. In the case
$$
q_5\in\left(-\infty,-2\sqrt{3|t_1|}
\right)\cup\left(0,2\sqrt{3|t_1|}\right),
$$
it has one self-intersection point and two semicubical cusps. The self-intersection point has coordinates
$$
t_2=0,\quad S_3=-\frac{(q_5^2 - 12\delta t_1)^2}{12\delta q_5}.
$$
It is the transversal intersection point of two smooth branches of the curve {\rm(\ref{psi-2})} corresponding to the values
$$
u=\pm\sqrt{\frac{q_5\left(12\delta t_1-q_5^2\right)}{3\left(12\delta t_1+q_5^2\right)}}.
$$
The cusps correspond to the values
\begin{equation}
u=\pm\sqrt{\frac{q_5\left(\sqrt[3]{12\delta t_1}-\sqrt[3]{q_5^2}\right)}{3\left(\sqrt[3]{12\delta t_1}+\sqrt[3]{q_5^2}\right)}}.
\label{osobye-S4-0}
\end{equation}
If $q_5=\sqrt{-\frac32\delta t_1}$, then the cusps lie on the line $S_3=0$. If $q_5=\sqrt{\frac32\delta t_1}$, then they lie on the line $S_3=4q_5t_1$.

{\rm 4)} The closure of the curve {\rm(\ref{psi-2})} {\rm(}as a subset of $\mathbb{R}^2${\rm)} is obtained by adding the limit point
$$
t_2=0,\quad S_3=4q_5t_1
$$
as $u\rightarrow\infty$. The closure is tangent to the line $S_3=4q_5t_1$ at $t_2=0$. The tangency is simple if $q_5^2\neq-12\delta t_1$; in the case $q_5^2=-12\delta t_1$, the closure of the curve {\rm(\ref{psi-2})} has a singularity of type $4/3$ at the point of tangency.}

\medskip
The proof of Lemma \ref{samoperesechS4} is a direct calculation. The formula (\ref{osobye-S4-0}) follows from Lemma \ref{vyr-krit-psi}. The skeletons $\Sigma_{0,t_1,q_5}$ for points $(t_1,q_5)$ of domains $1,8,22,25,28,29,33,38,60,93$ in Fig. \ref{razbienie0} are shown in the middle pictures of the triptychs presented in Fig. \ref{predel3-2chetv} and \ref{predel1-4chetv} up to diffeomorphism. The $\delta t_2$ axis is horizontal and is directed to the right. The $\delta S_3$ axis is vertical and is directed upward. The line (\ref{crit-tochki}) is dashed.

The pictures on the left and on the right of each triptych shows the limit positions of the skeleton $\Sigma_{S_4,t_1,q_5}$ as $\delta S_4\rightarrow-0$ and $\delta S_4\rightarrow+0$, respectively. The symbols ${\bf k}_i^{\pm}$ denote the limits of the strata ${\bf k}_i\in\mathfrak{S}_{S_4,t_1,q_5}$ as $\delta S_4\rightarrow\pm0$. The dotted line denotes (\ref{peresech-3}) (in the middle pictures of the triptychs) and the limit positions of some branches of the curve (\ref{psi-2}) as $\delta S_4\rightarrow\mp0$ (in the pictures on the left and on the right).

Type $A_1^k$ strata of the stratification $\mathfrak{S}_{0,t_1,q_5}$ are denoted by ${\bf k}_j^0$. The type of a stratum is not changed when $(t_1,q_5)$ transversally intersects the dashed or dotted line (at a point that does not lie on a solid line). The coincidence of the subscripts of the strata ${\bf k}_j^0$ means that they are different connected components of the section of the same type $A_1^k$ stratum of the stratification $\mathfrak{S}$ by the hyperplane $S_4=0$. This stratum is denoted by ${\bf k}_j^{\mathfrak{S}}$. A description of the strata ${\bf k}_j^{\mathfrak{S}}$  for all possible ${\bf k}$ and $j$ is given in lemmas below.

\begin{figure}
\begin{center}
\begin{tabular}{ccc}
$-0$\includegraphics[width=4cm]{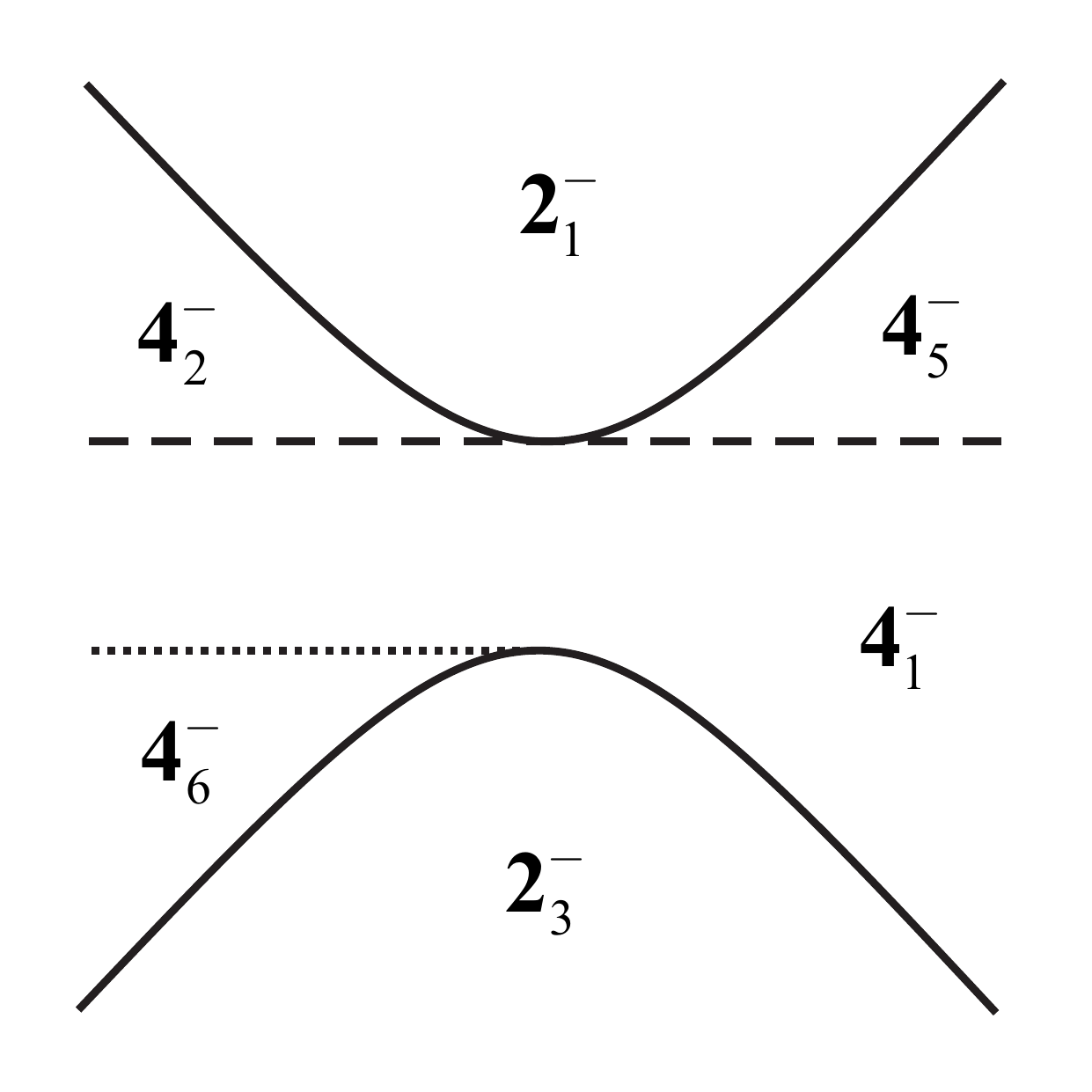}&
28)\includegraphics[width=4cm]{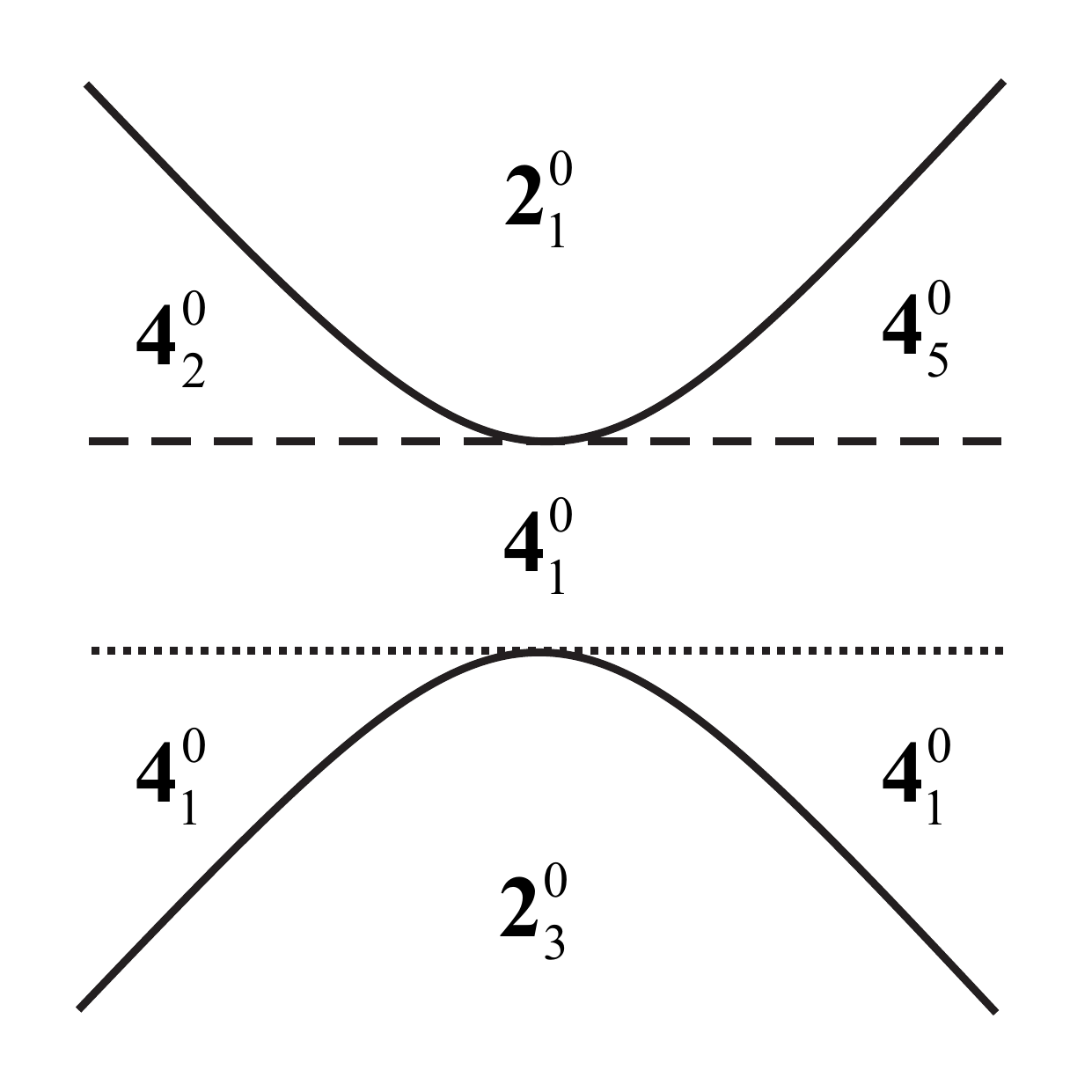}&
$+0$\includegraphics[width=4cm]{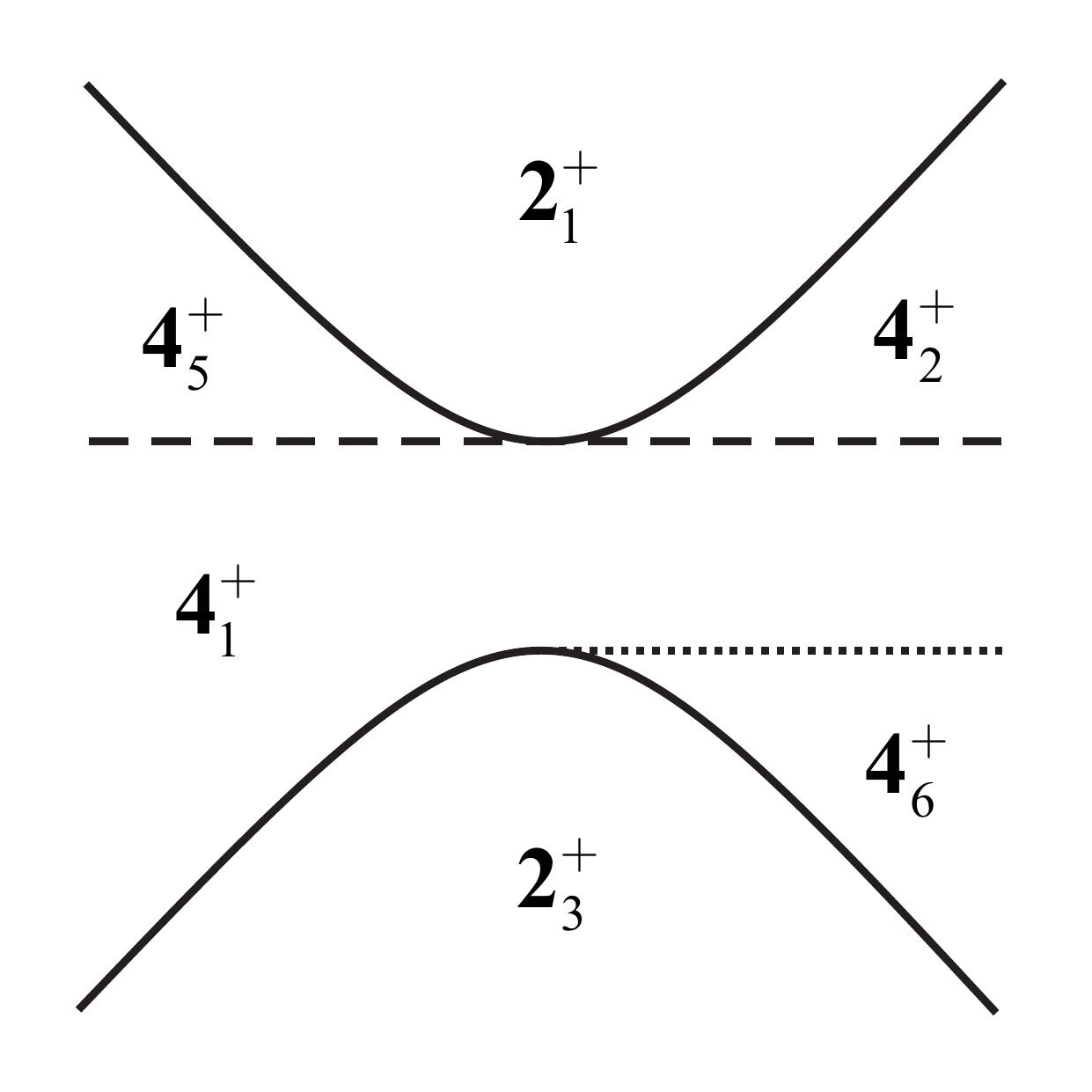}\\
$-0$\includegraphics[width=4cm]{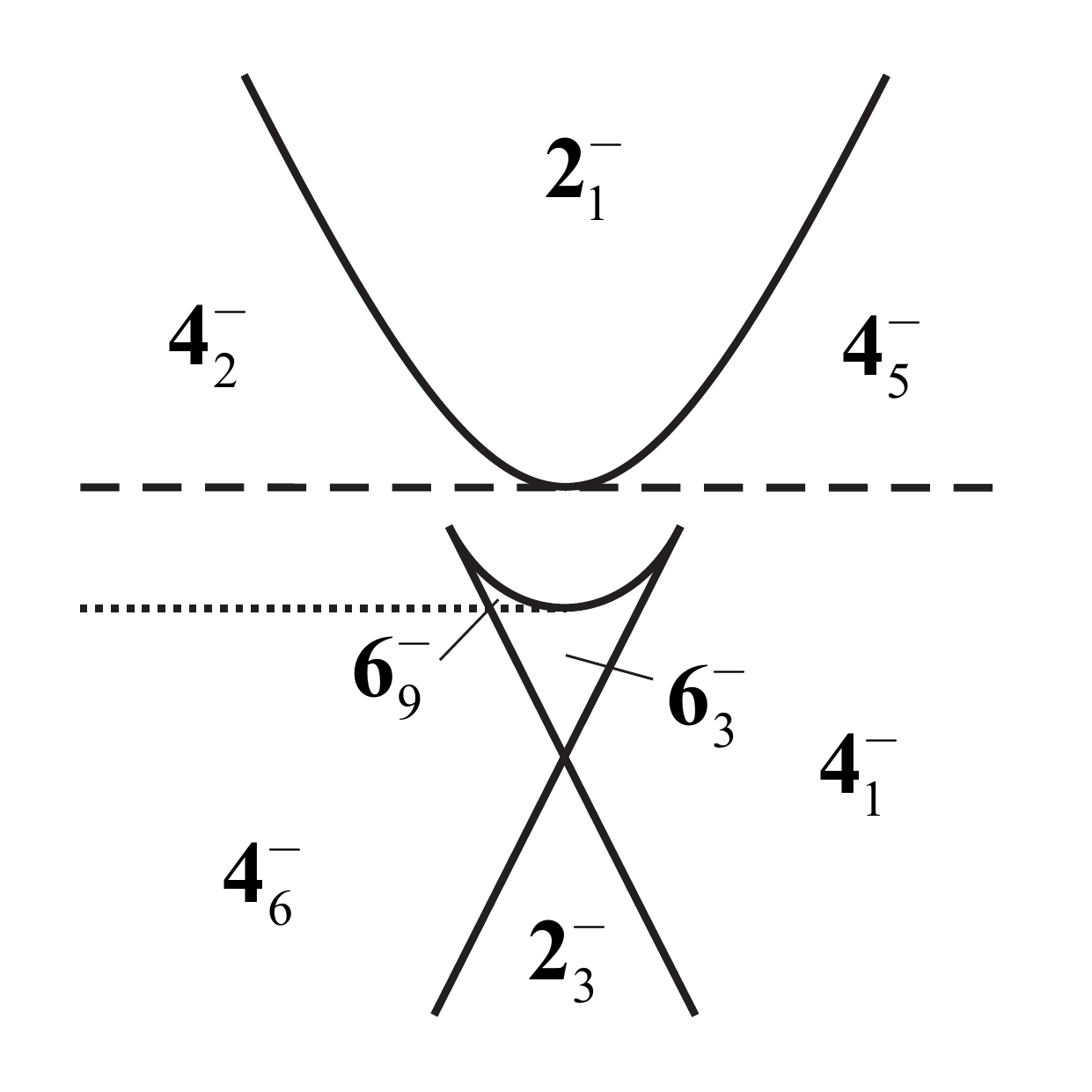}&
25)\includegraphics[width=4cm]{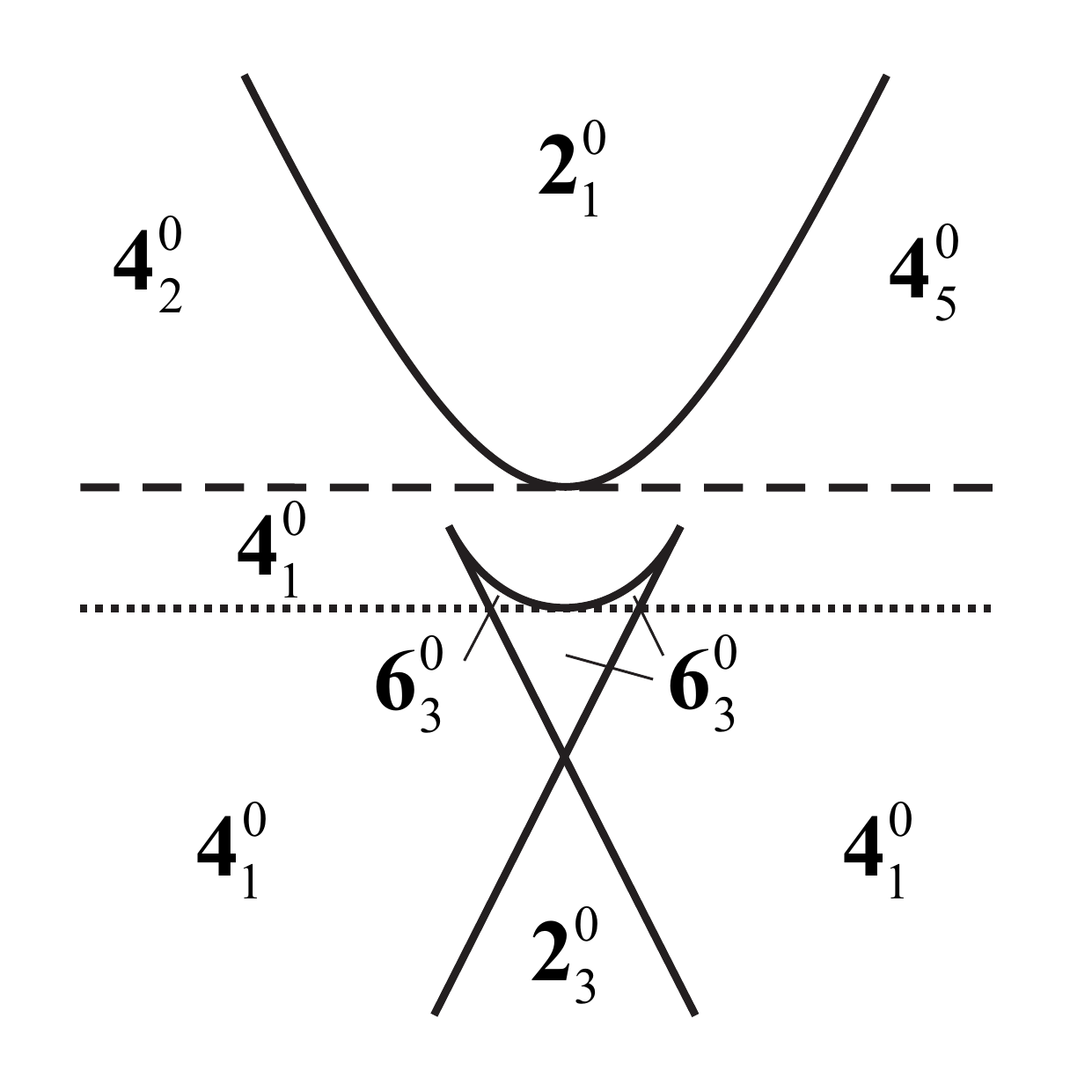}&
$+0$\includegraphics[width=4cm]{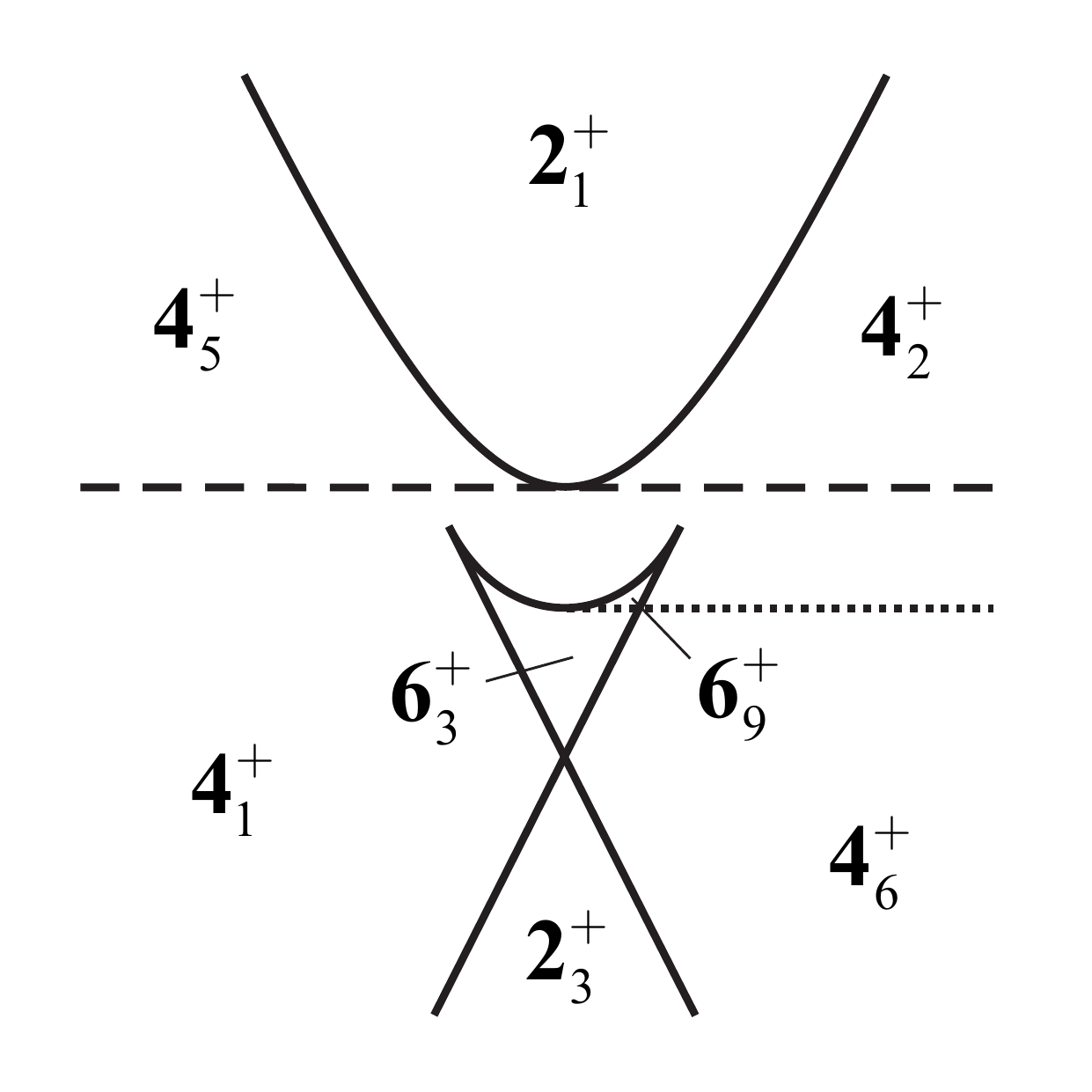}\\
$-0$\includegraphics[width=4cm]{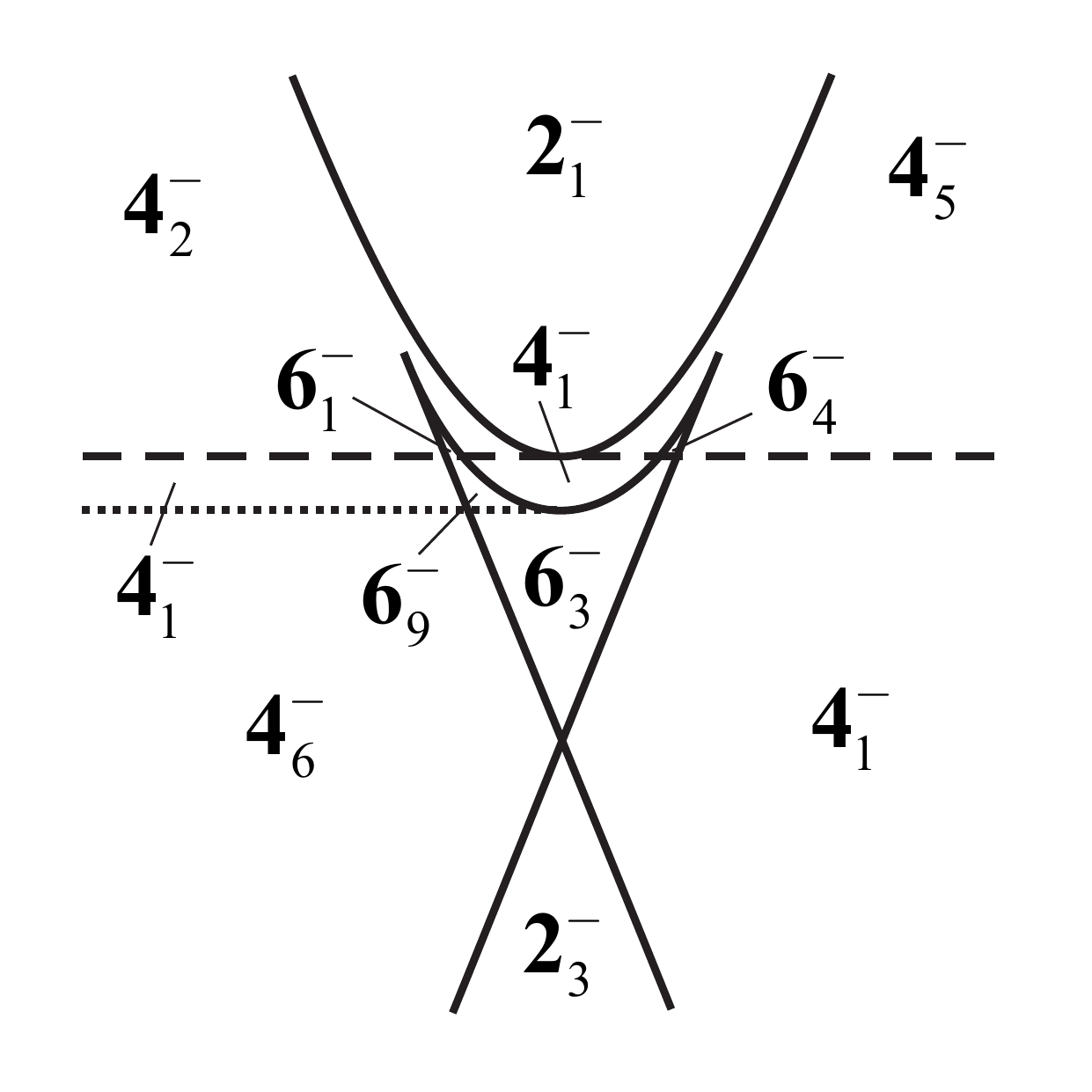}&
22)\includegraphics[width=4cm]{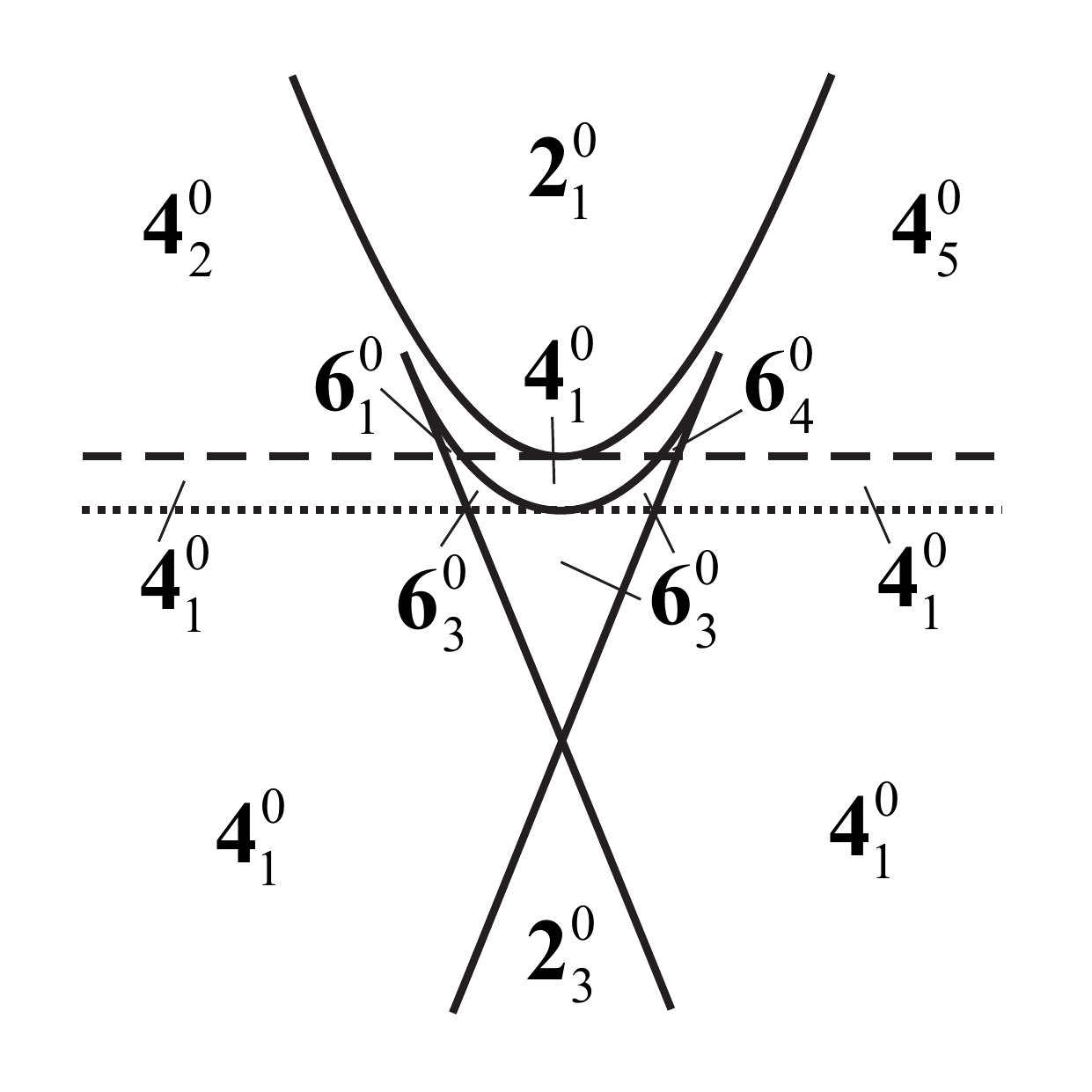}&
$+0$\includegraphics[width=4cm]{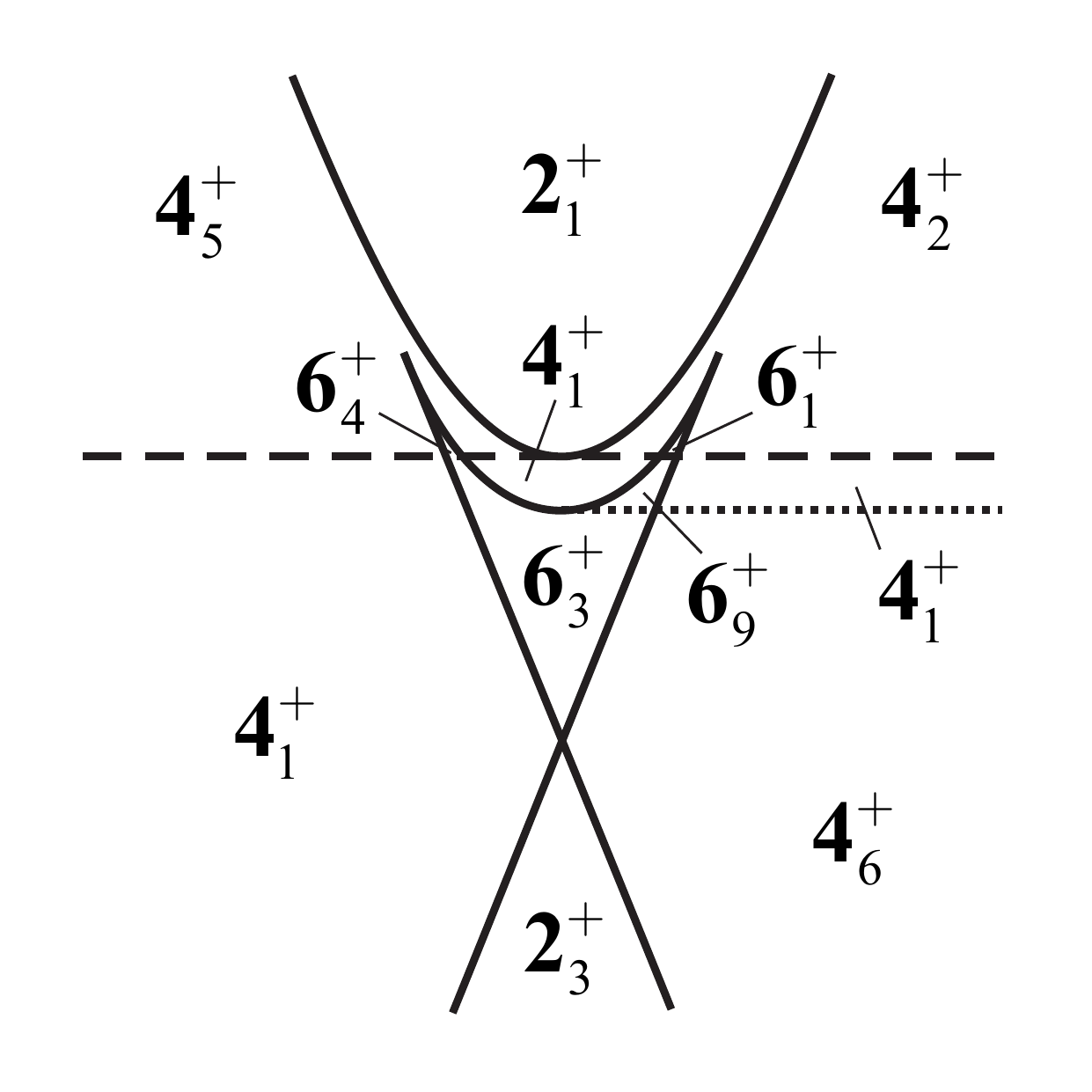}\\
$-0$\includegraphics[width=4cm]{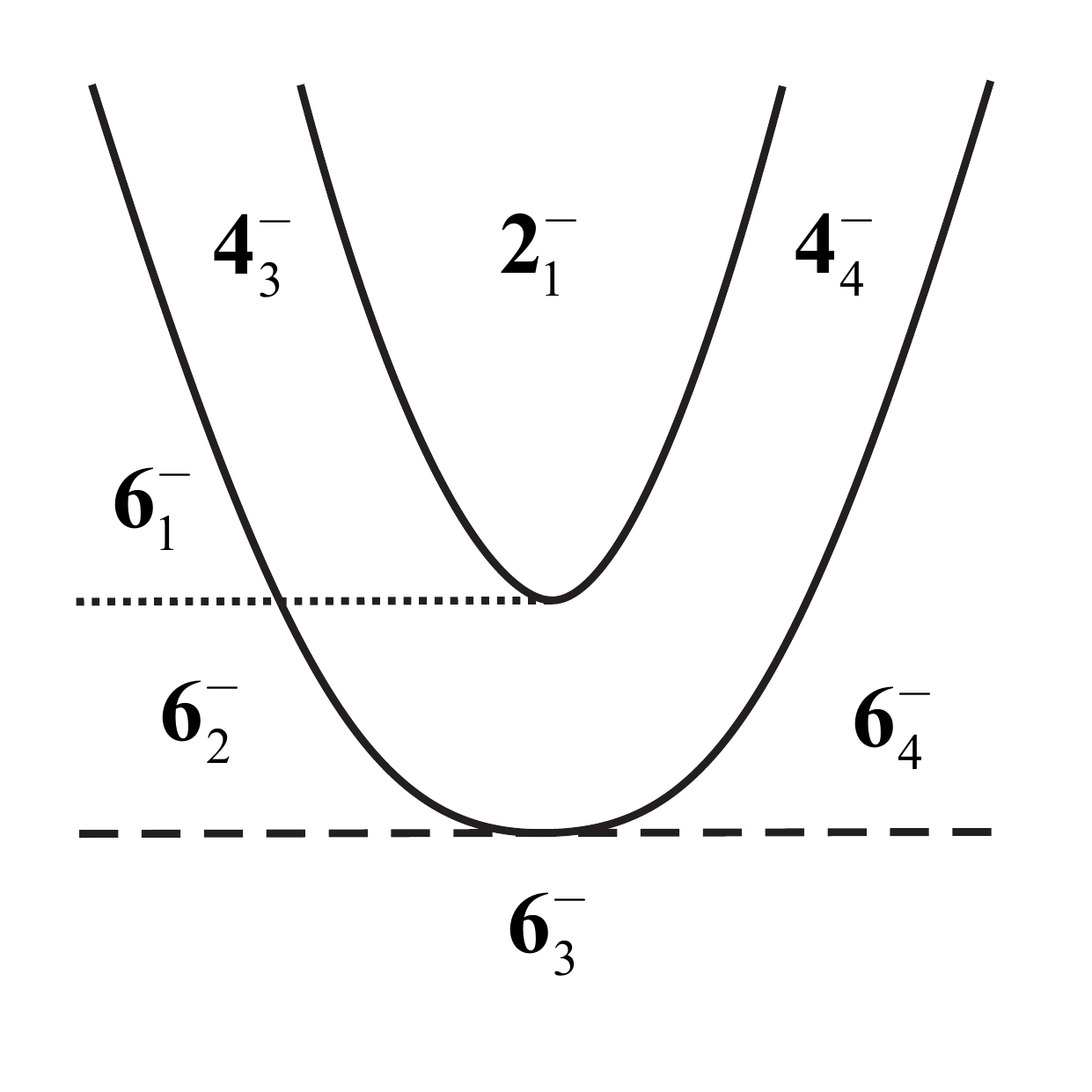}&
8)\includegraphics[width=4cm]{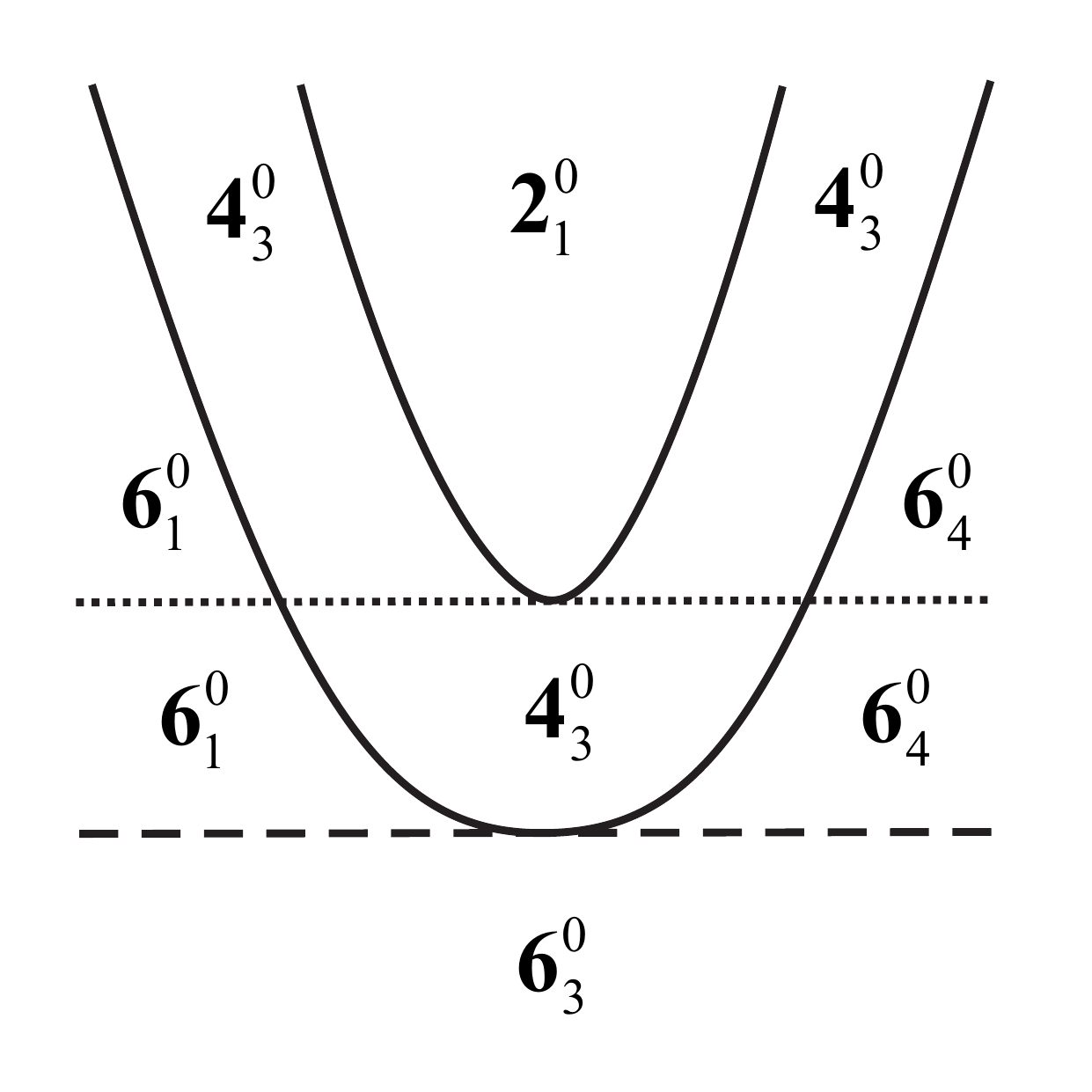}&
$+0$\includegraphics[width=4cm]{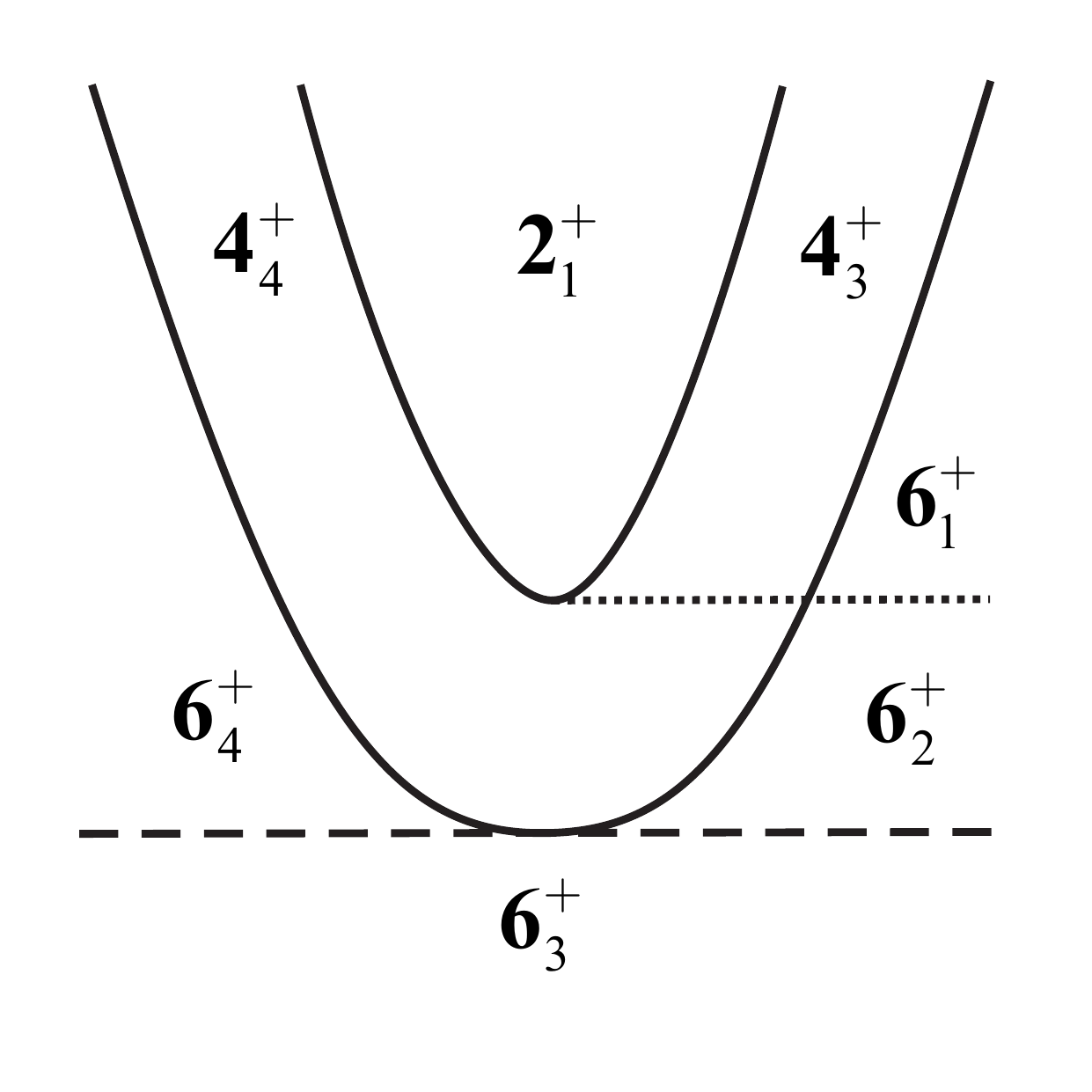}\\
$-0$\includegraphics[width=4cm]{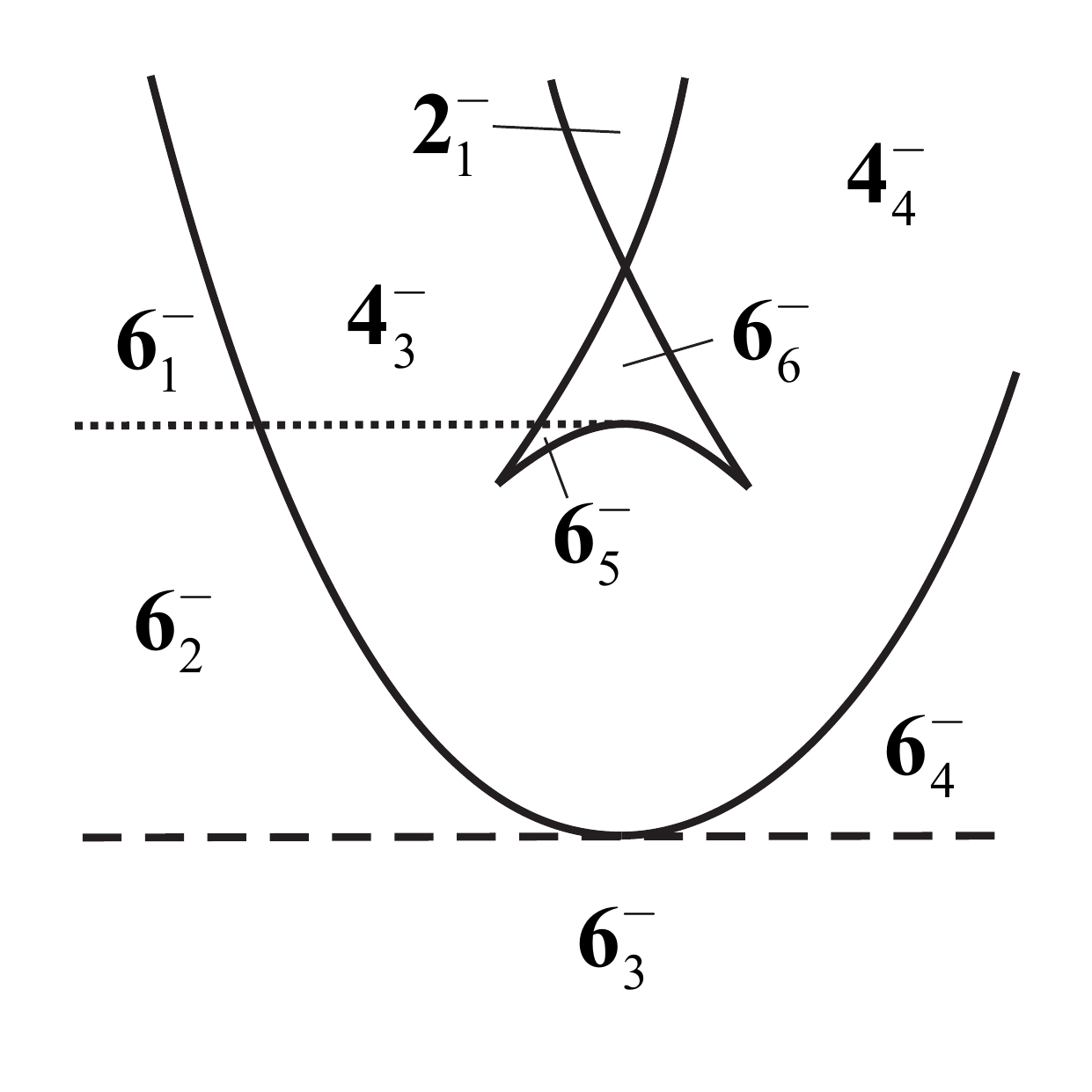}&
1)\includegraphics[width=4cm]{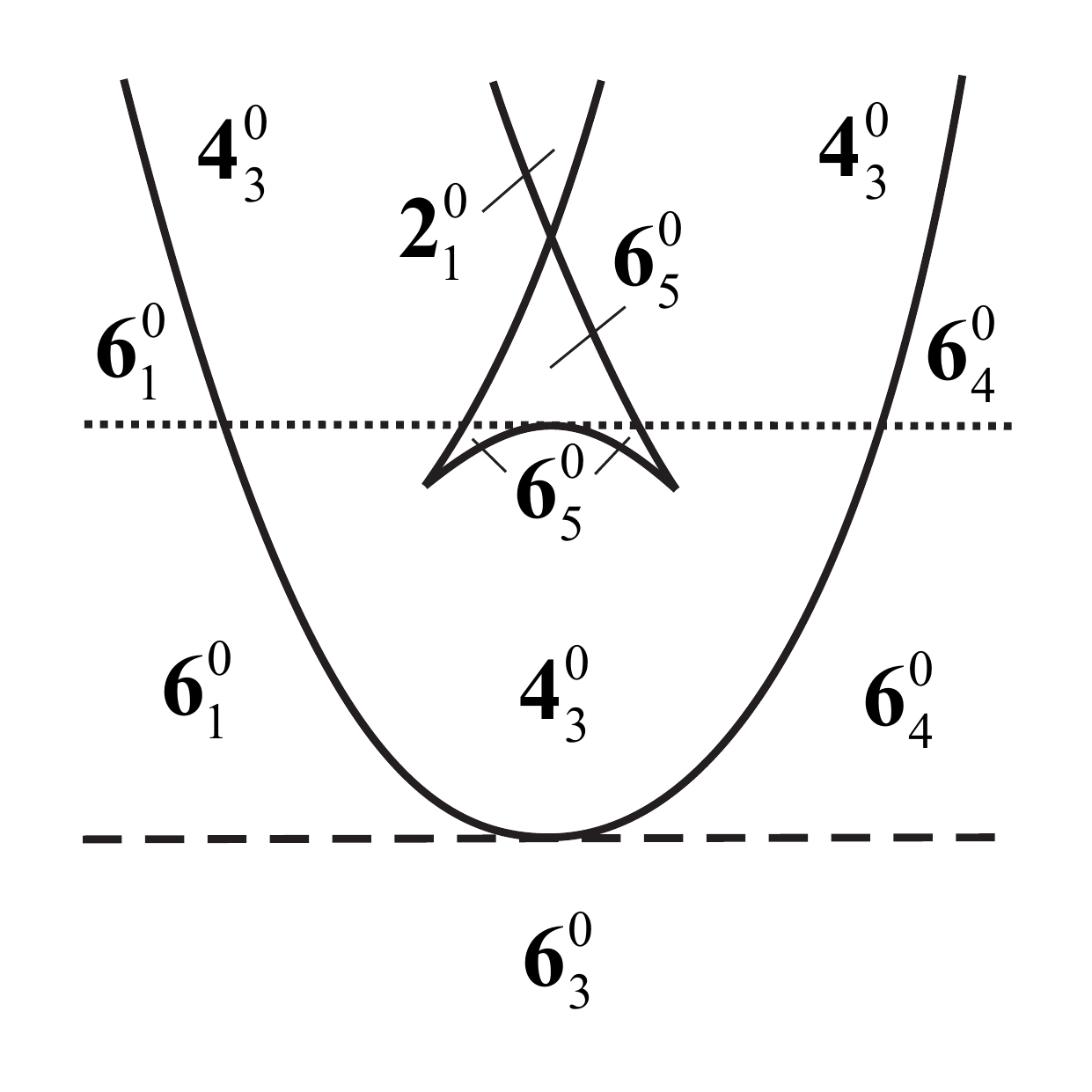}&
$+0$\includegraphics[width=4cm]{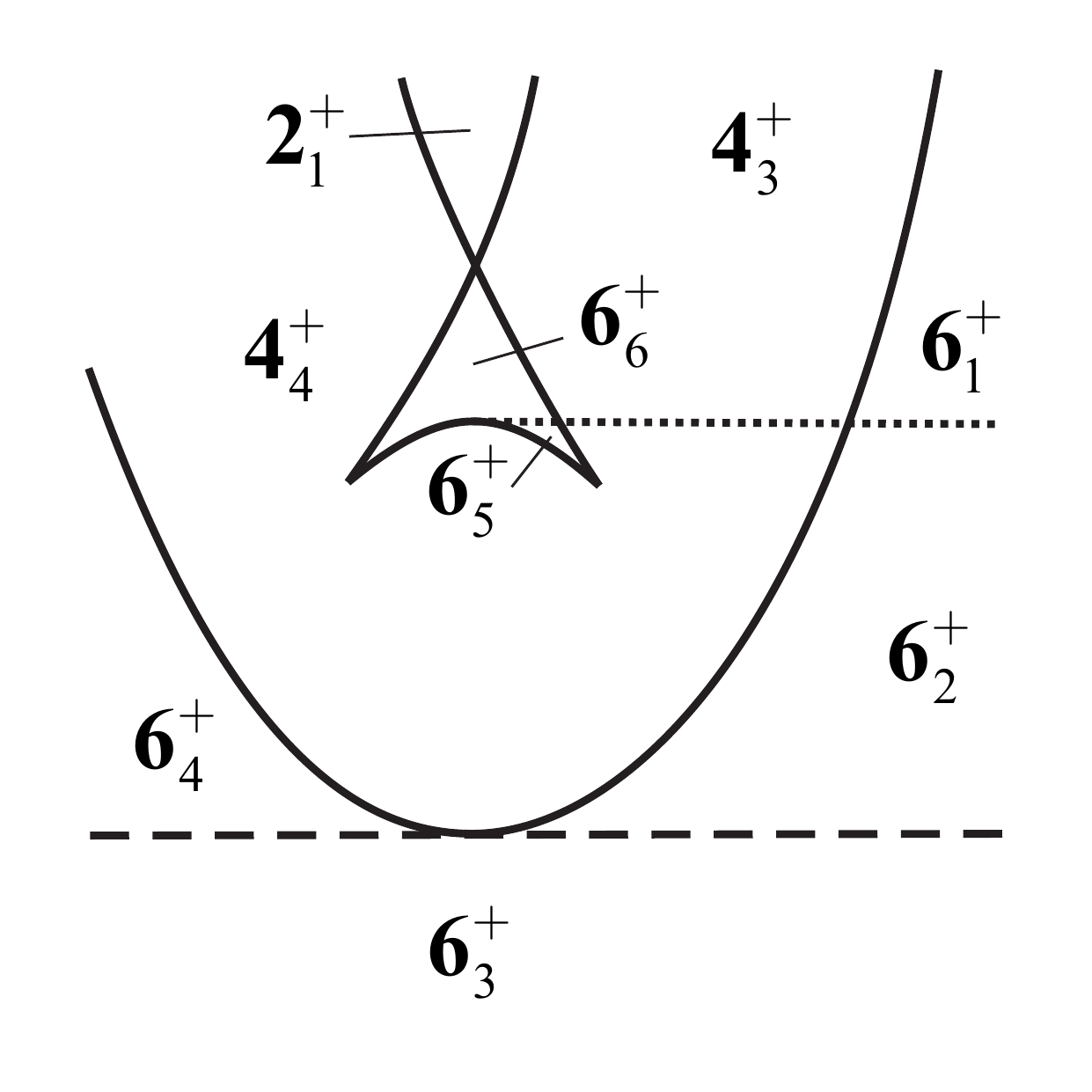}
\end{tabular}
\caption{The skeletons $\Sigma_{0,t_1,q_5}$ for domains $1,8,22,25,28$.}
\label{predel3-2chetv}
\end{center}
\end{figure}

\begin{figure}
\begin{center}
\begin{tabular}{ccc}
$-0$\includegraphics[width=4cm]{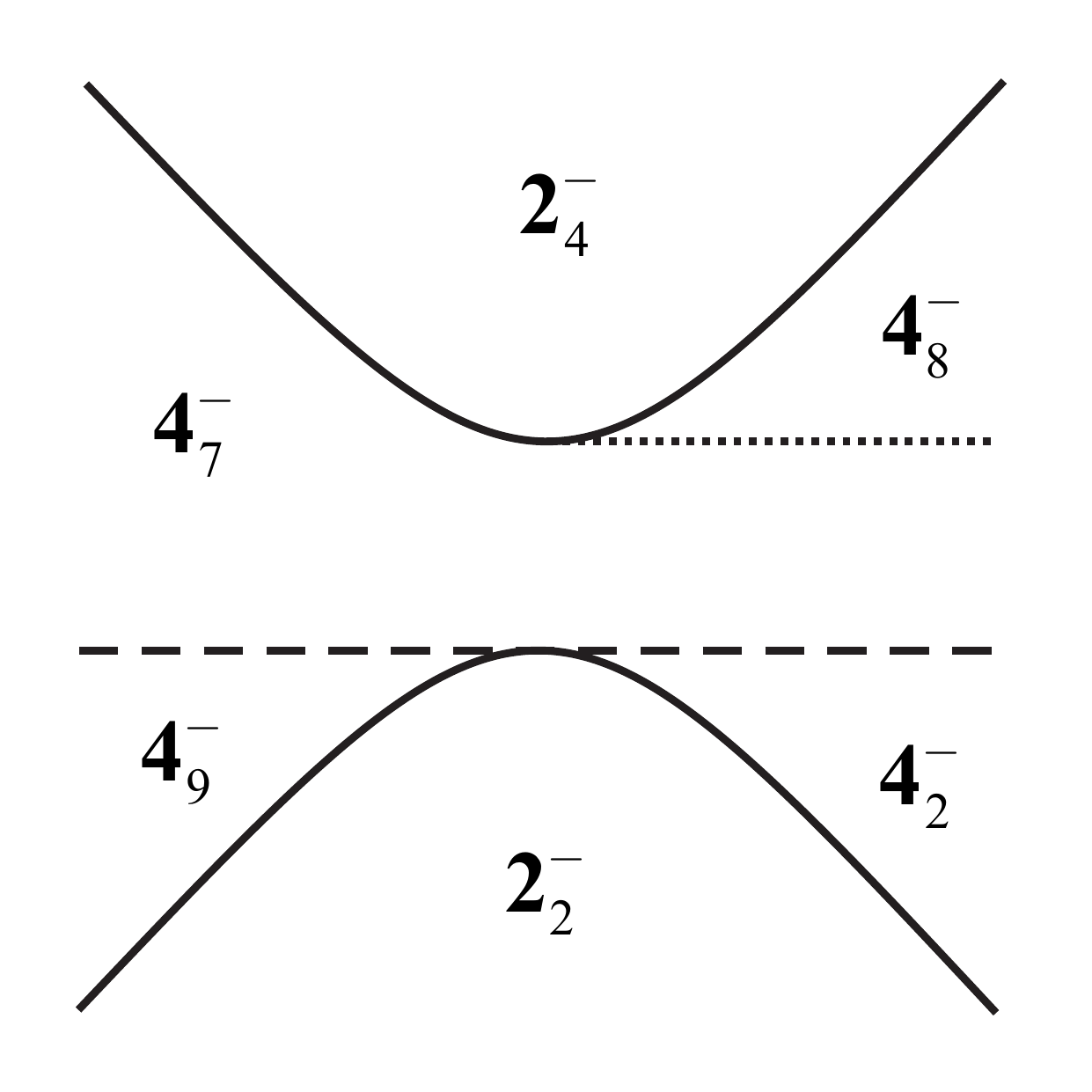}&
29)\includegraphics[width=4cm]{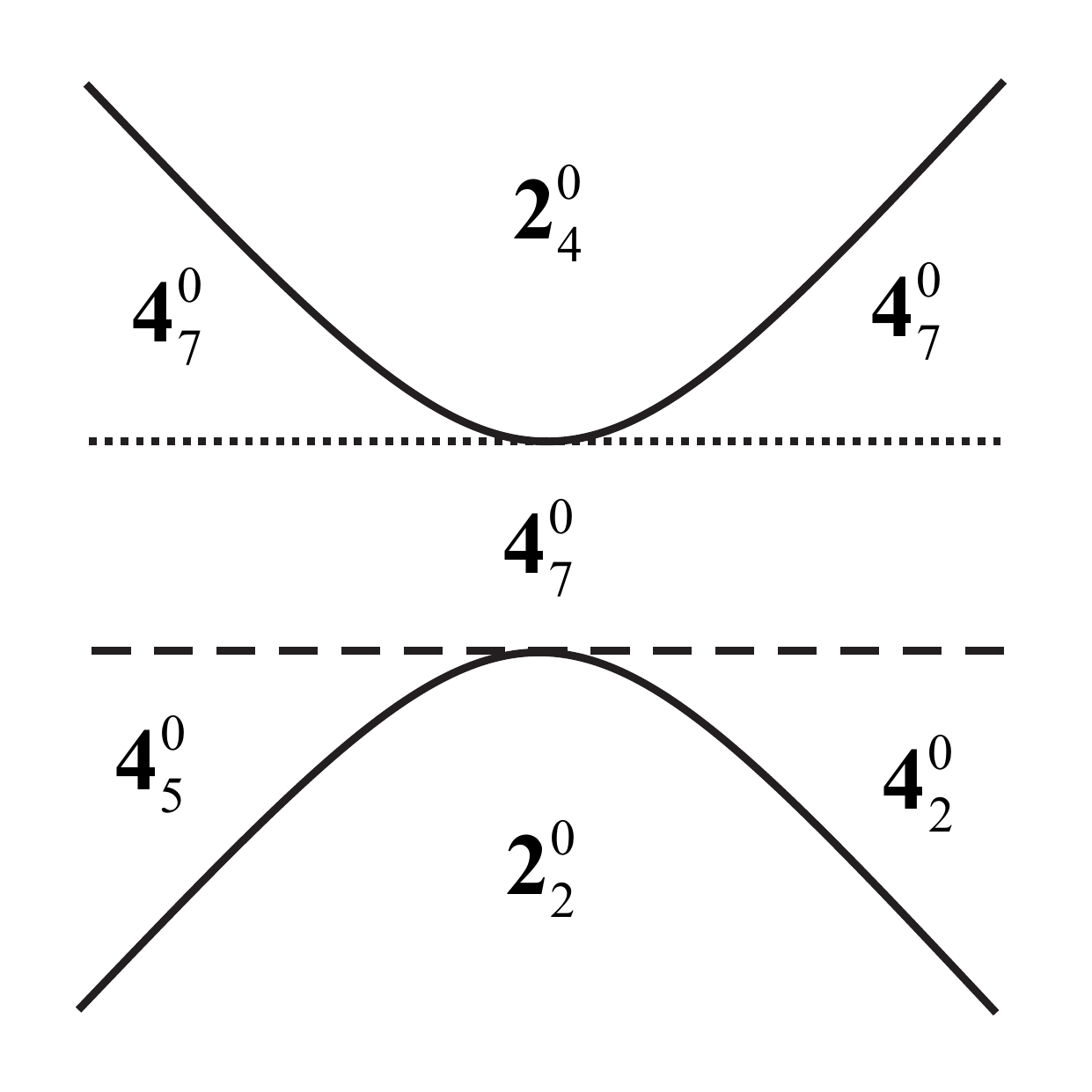}&
$+0$\includegraphics[width=4cm]{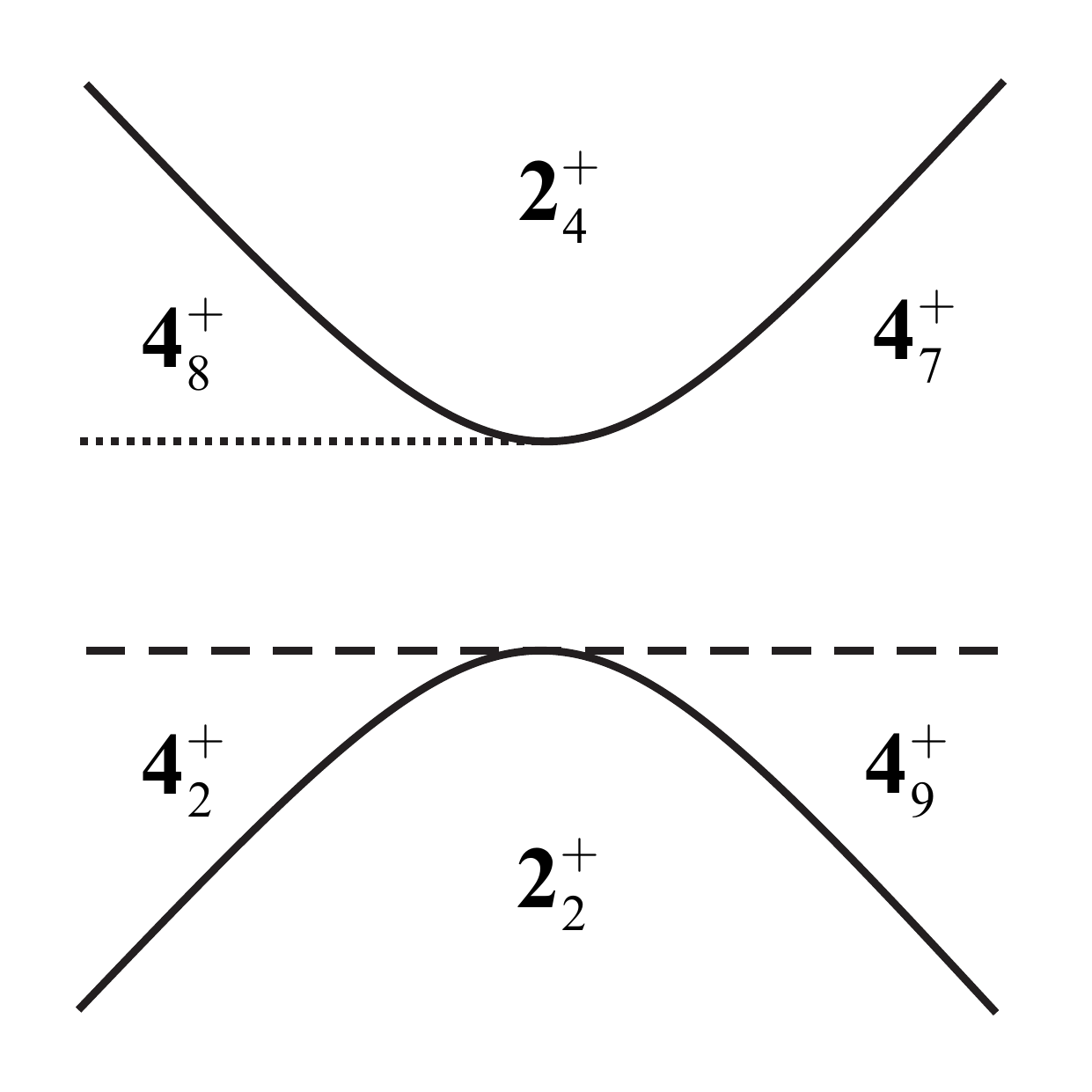}\\
$-0$\includegraphics[width=4cm]{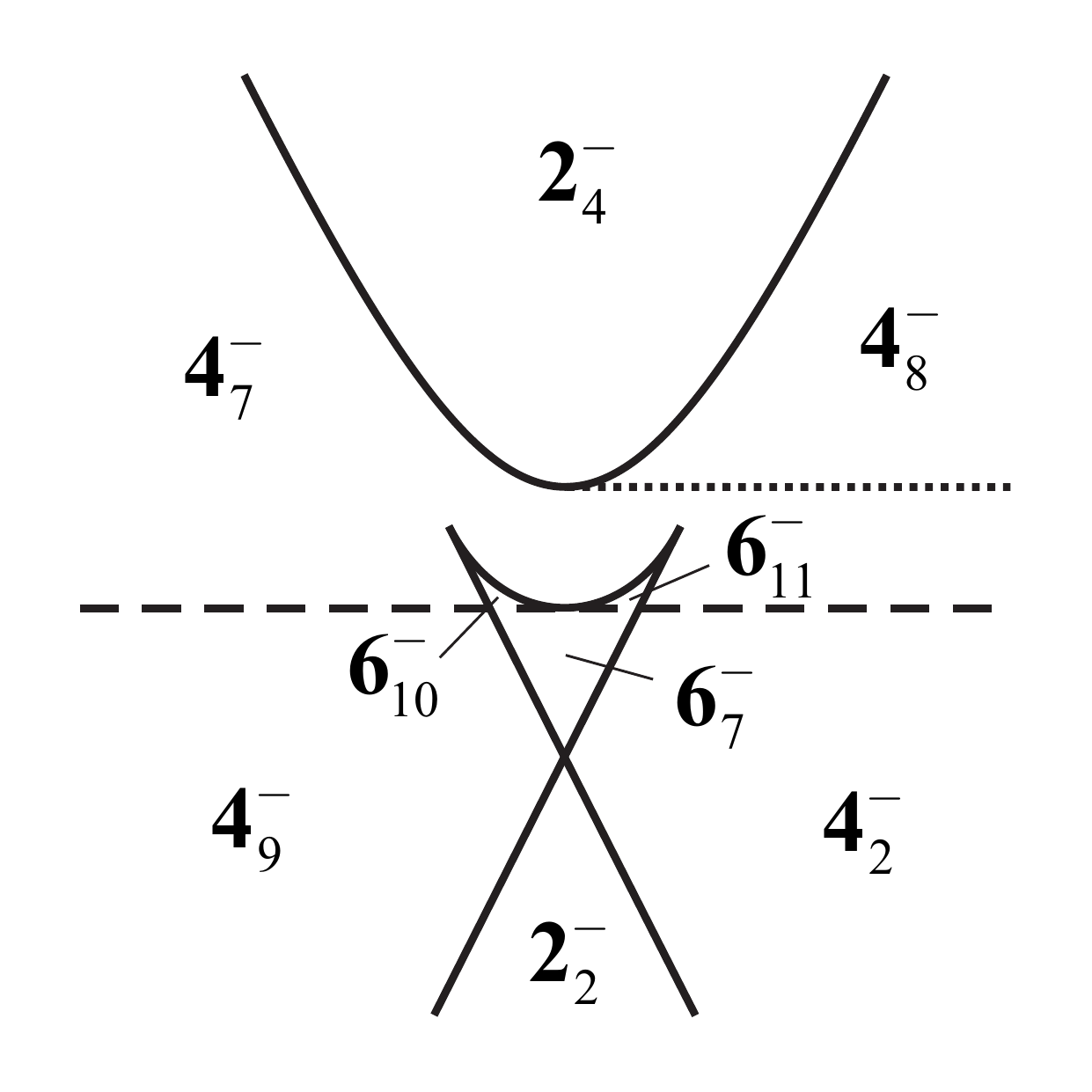}&
33)\includegraphics[width=4cm]{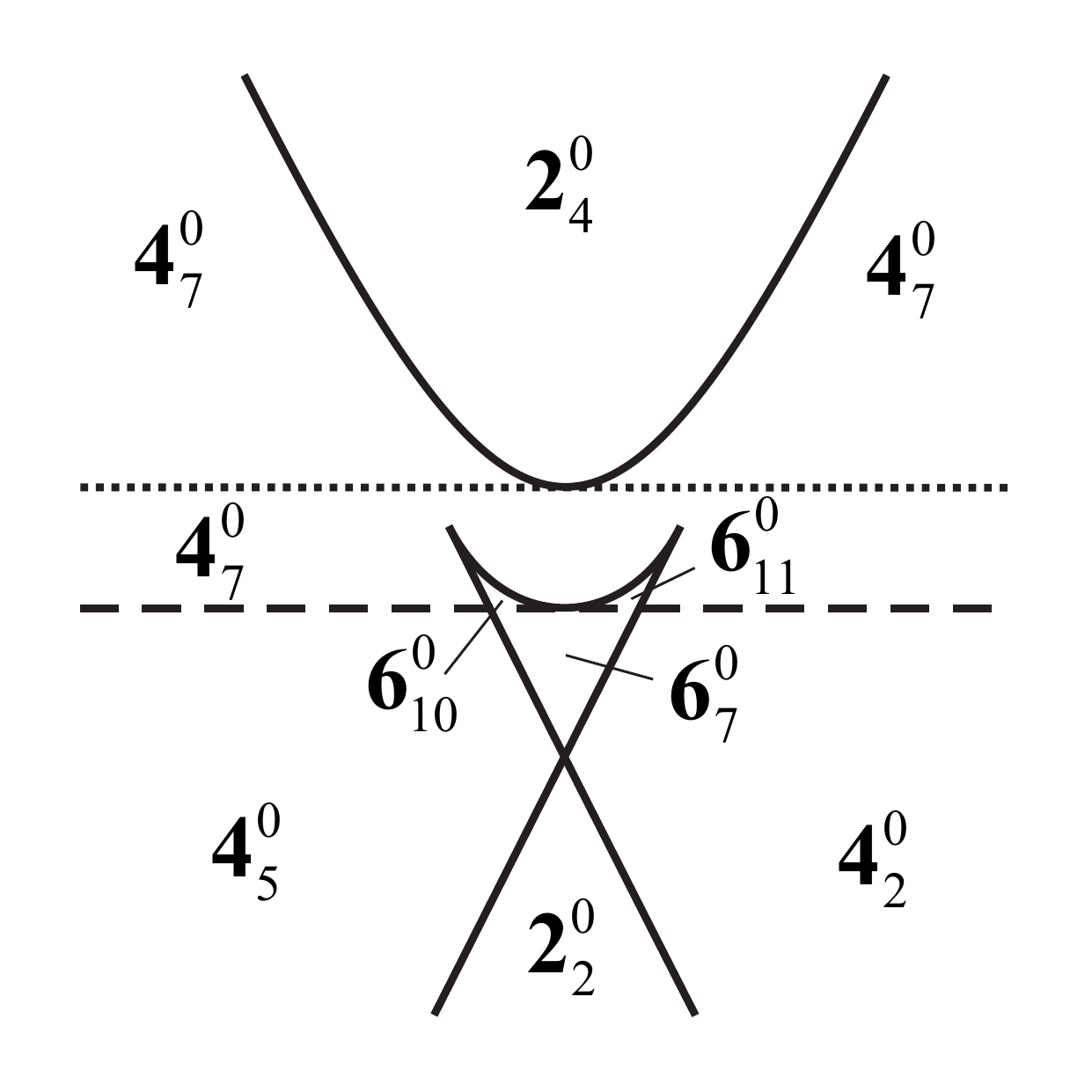}&
$+0$\includegraphics[width=4cm]{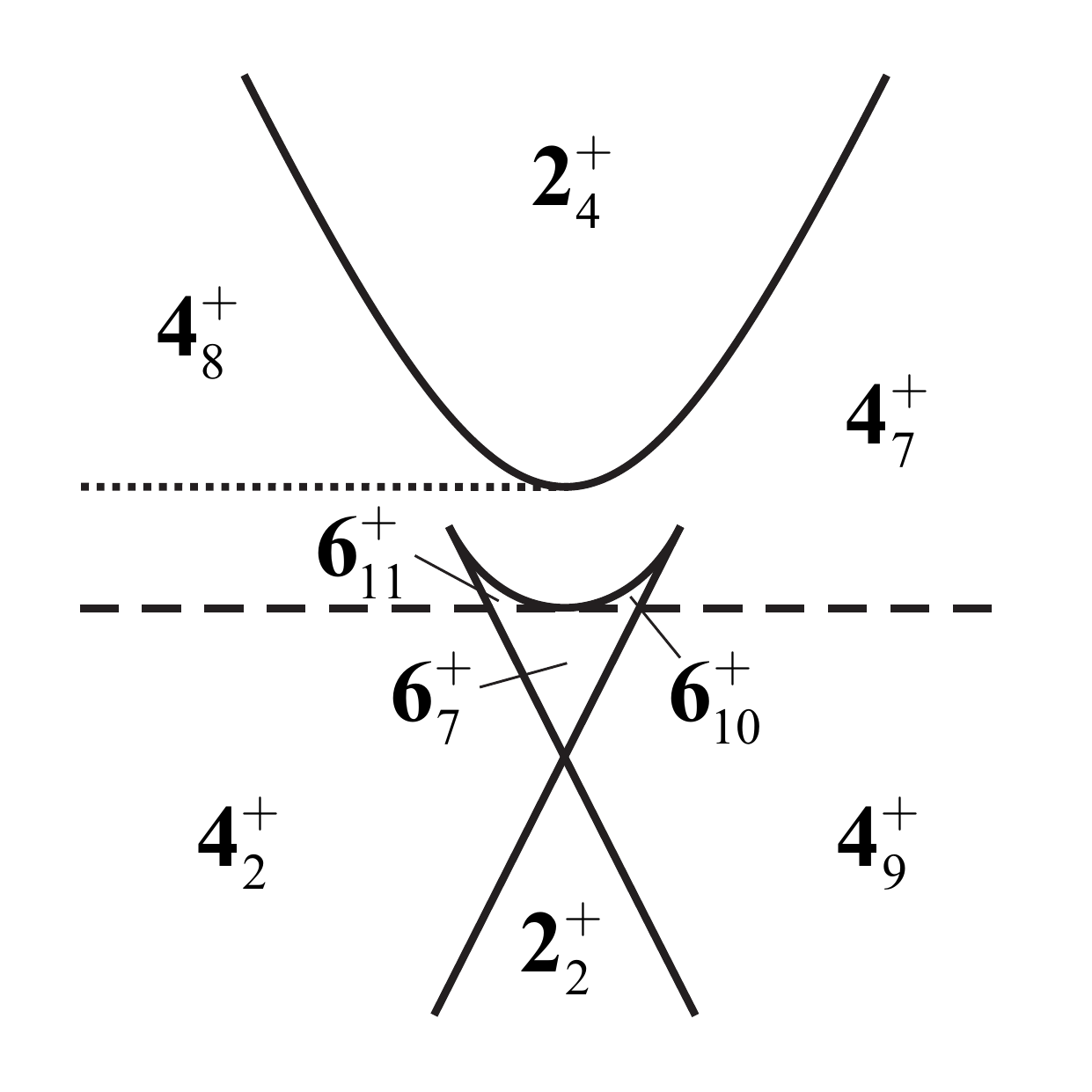}\\
$-0$\includegraphics[width=4cm]{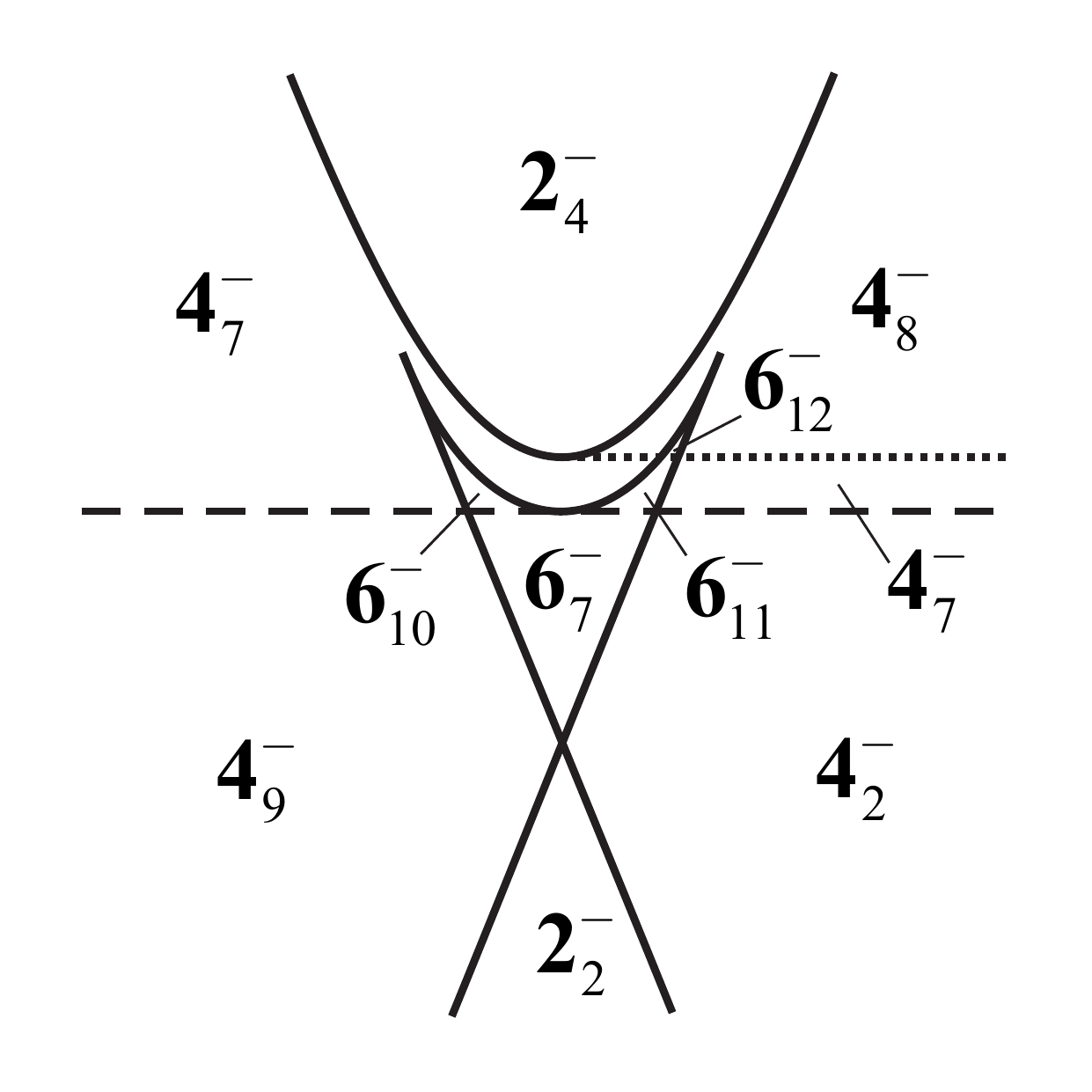}&
38)\includegraphics[width=4cm]{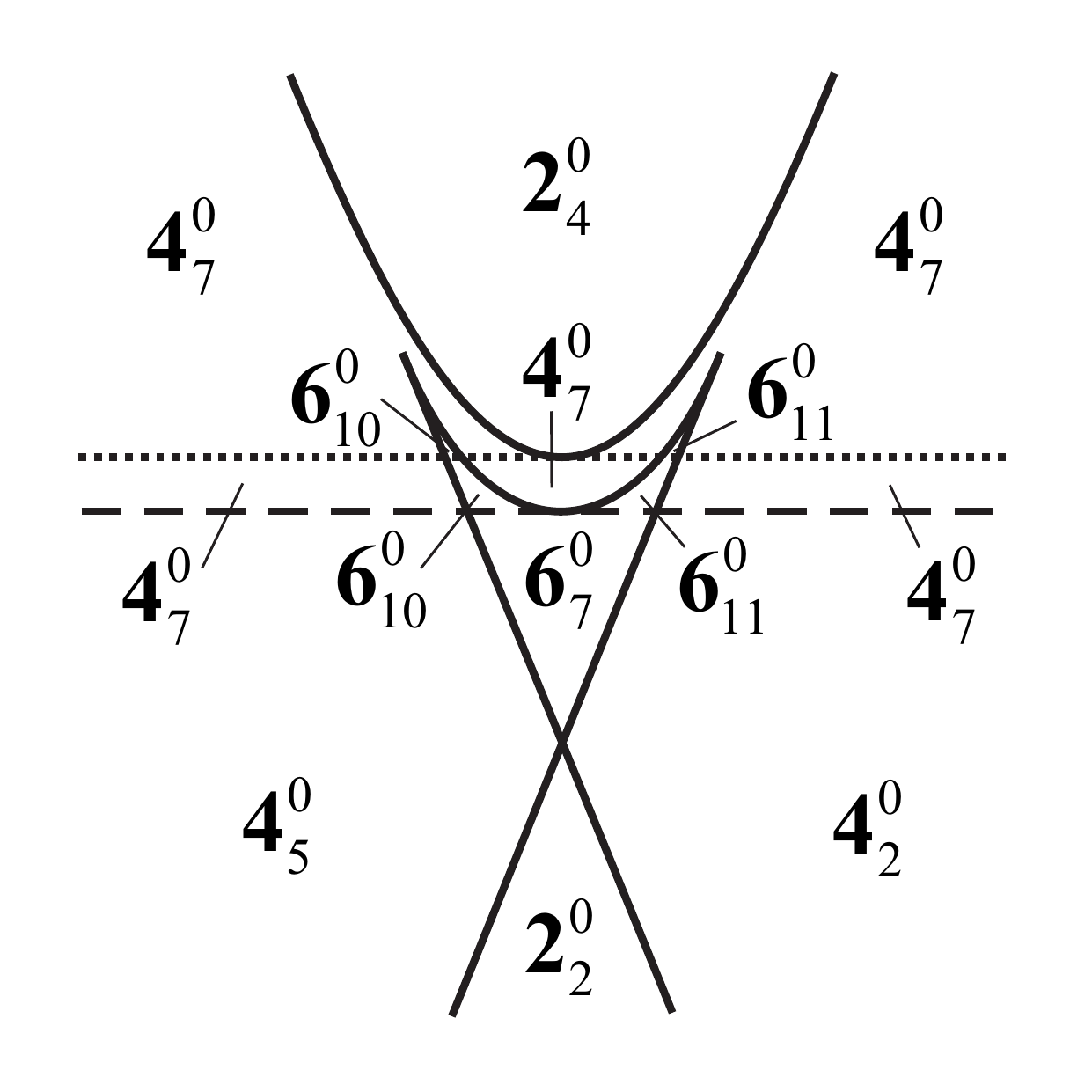}&
$+0$\includegraphics[width=4cm]{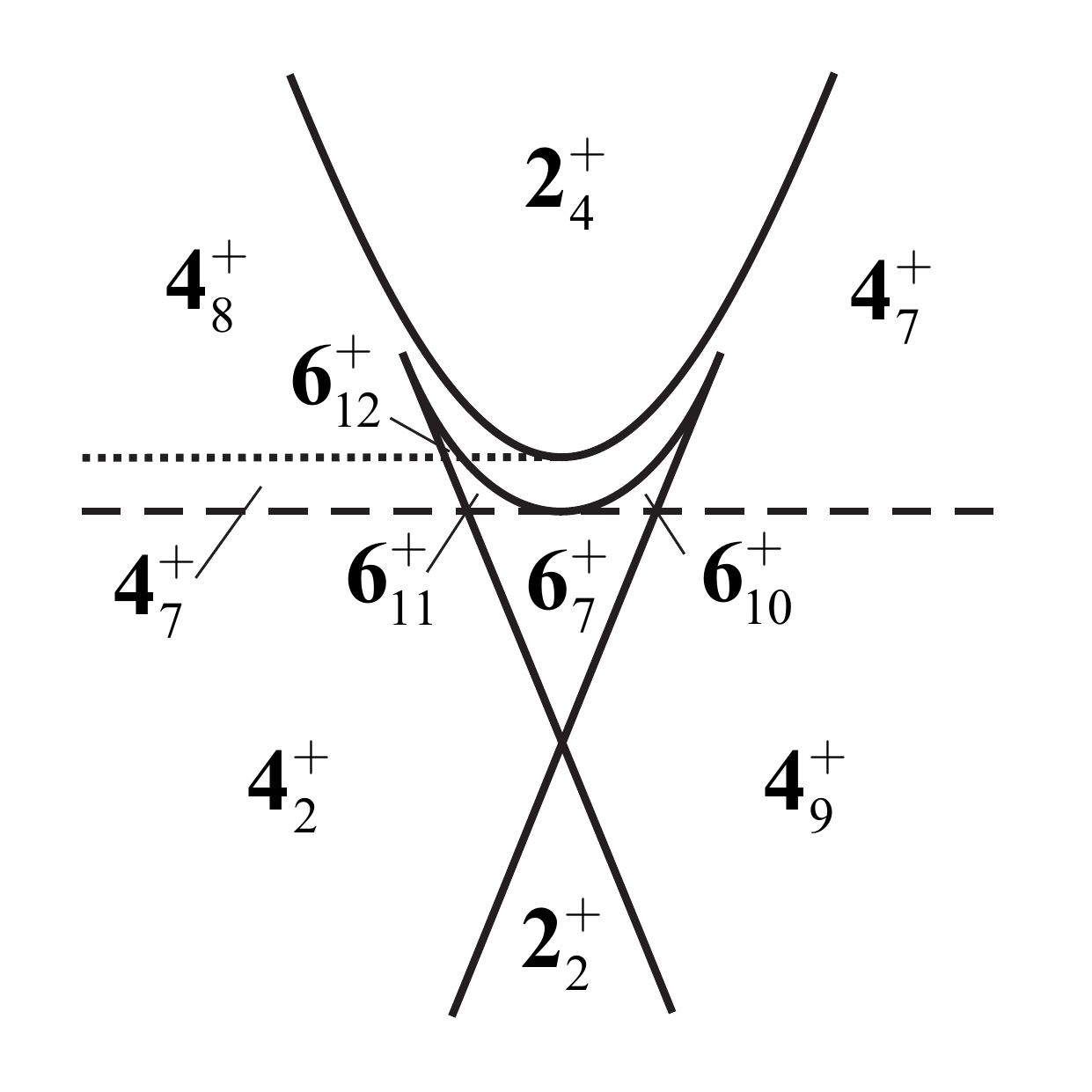}\\
$-0$\includegraphics[width=4cm]{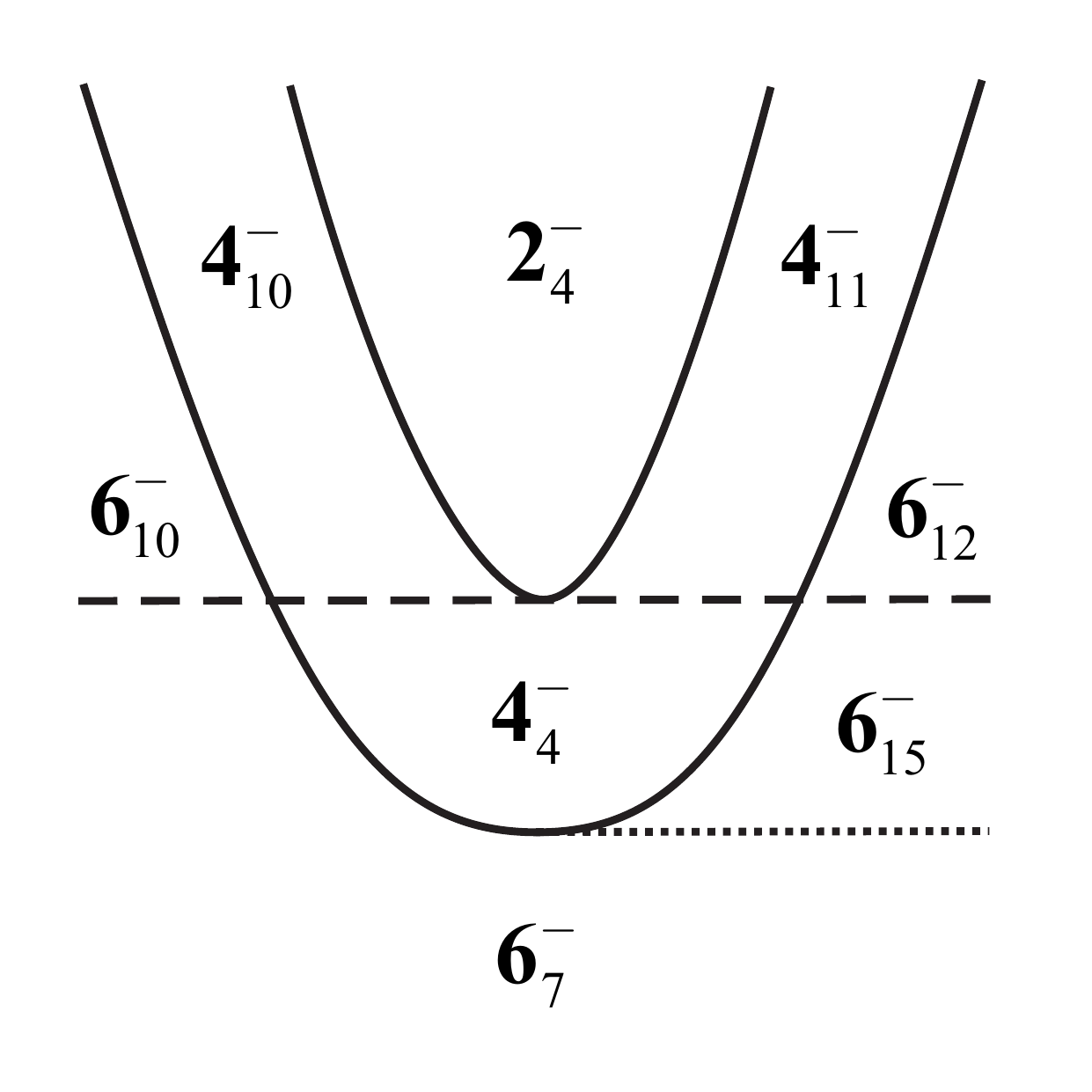}&
60)\includegraphics[width=4cm]{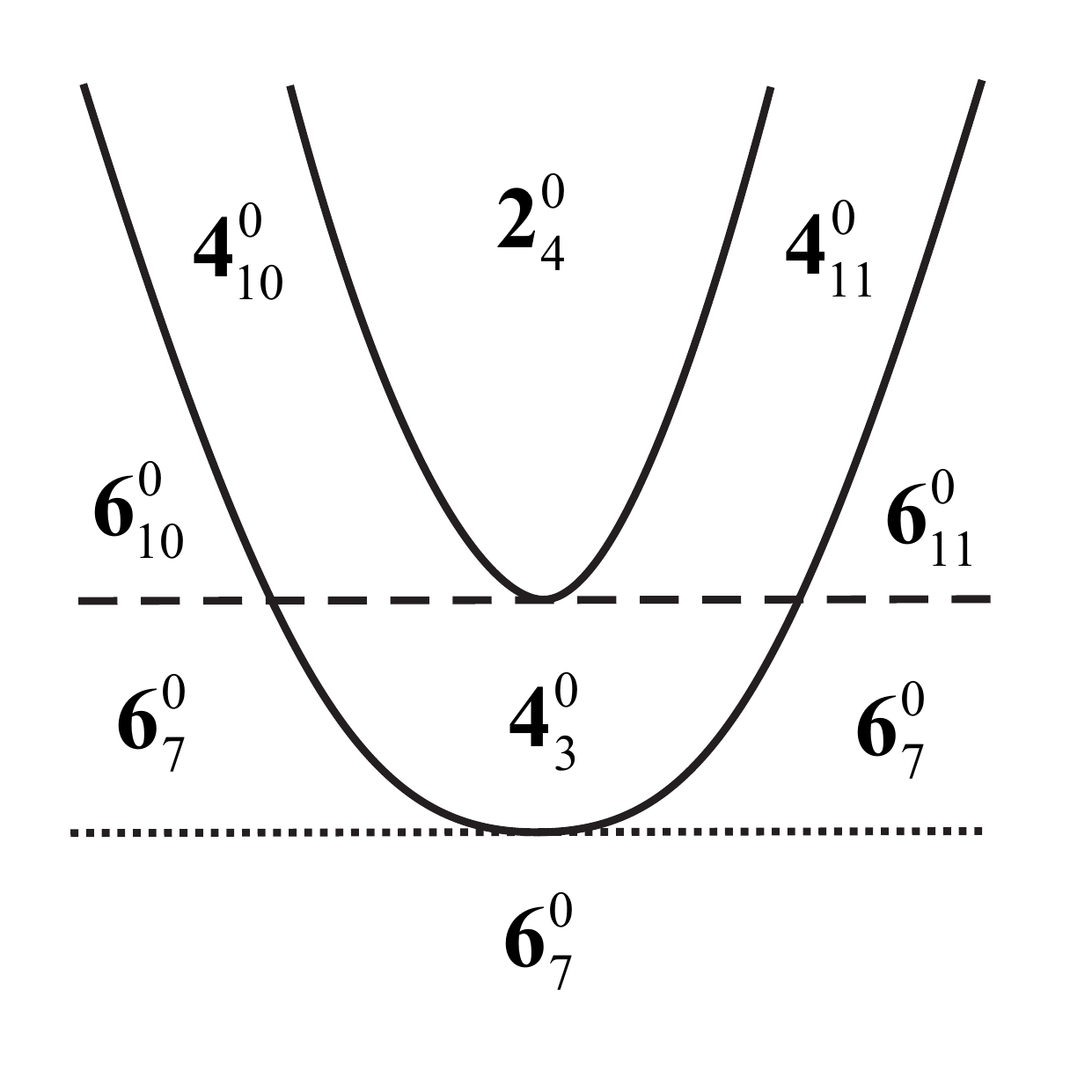}&
$+0$\includegraphics[width=4cm]{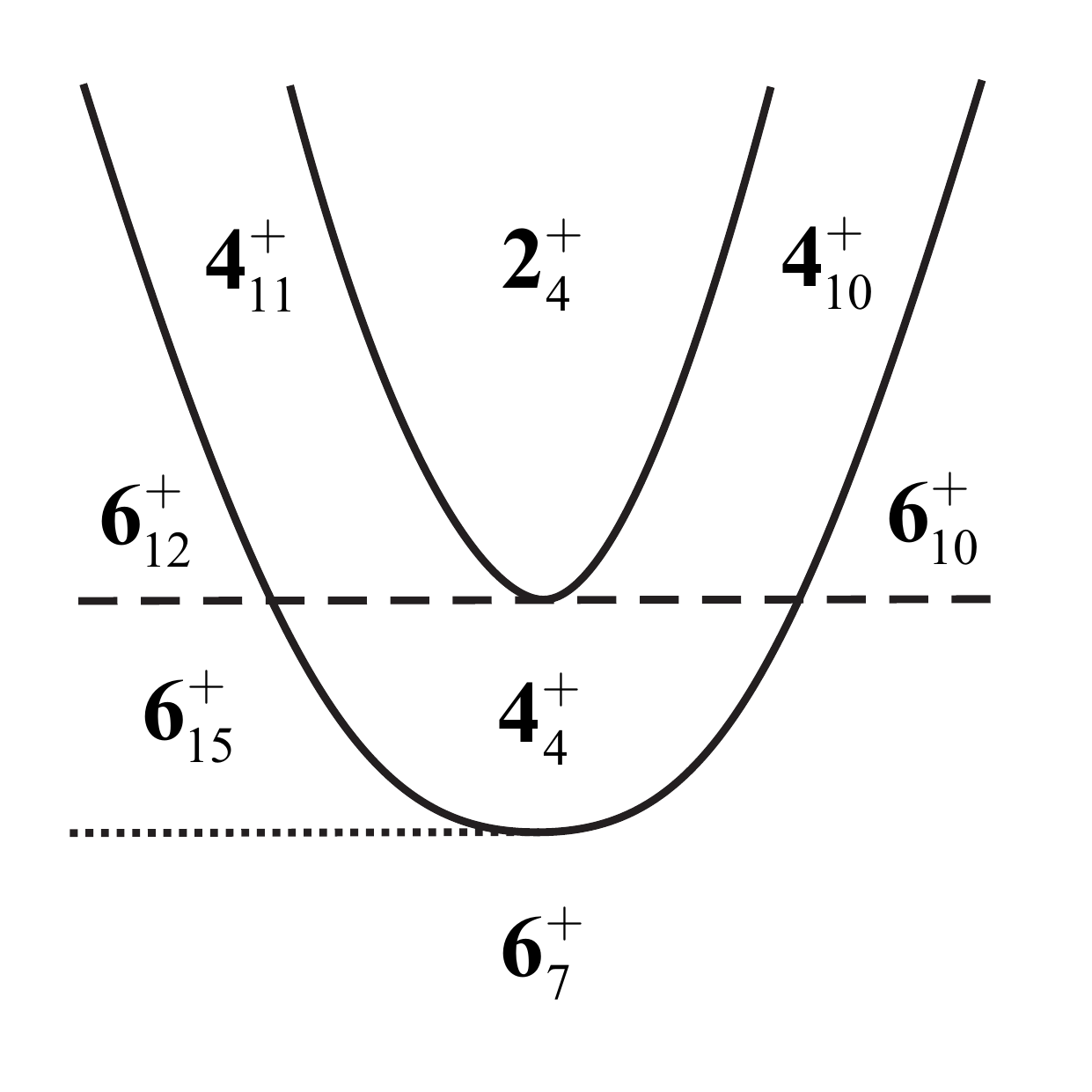}\\
$-0$\includegraphics[width=4cm]{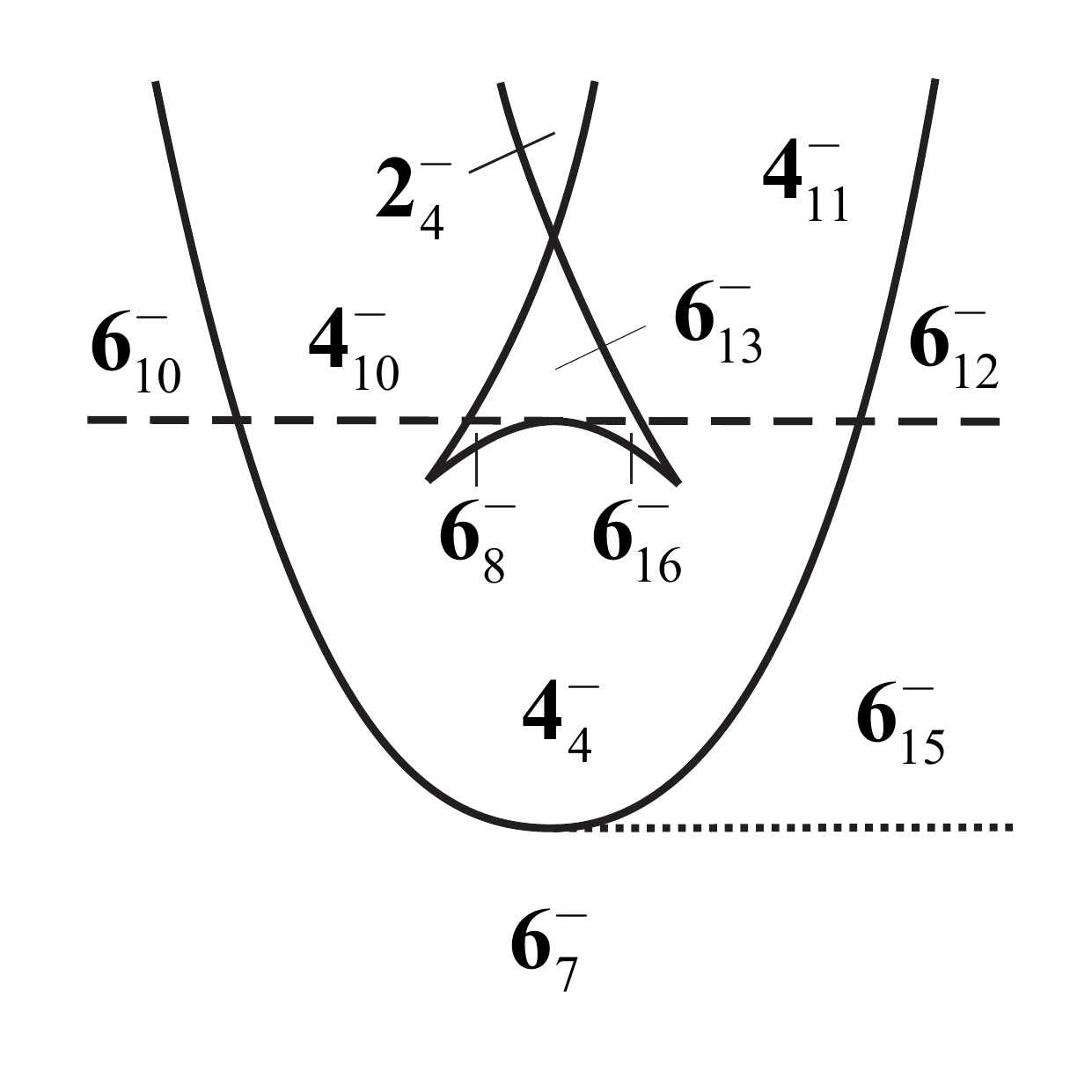}&
93)\includegraphics[width=4cm]{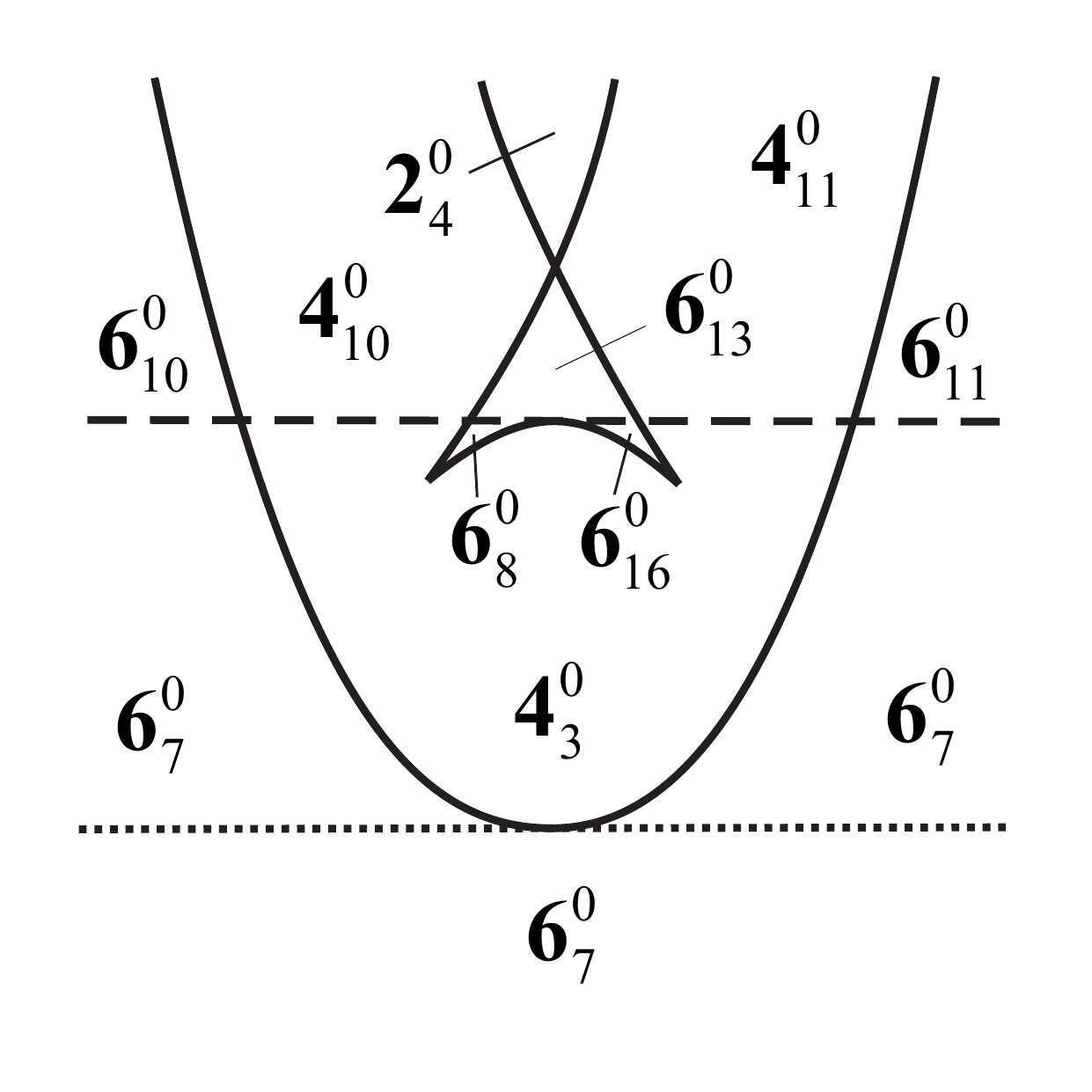}&
$+0$\includegraphics[width=4cm]{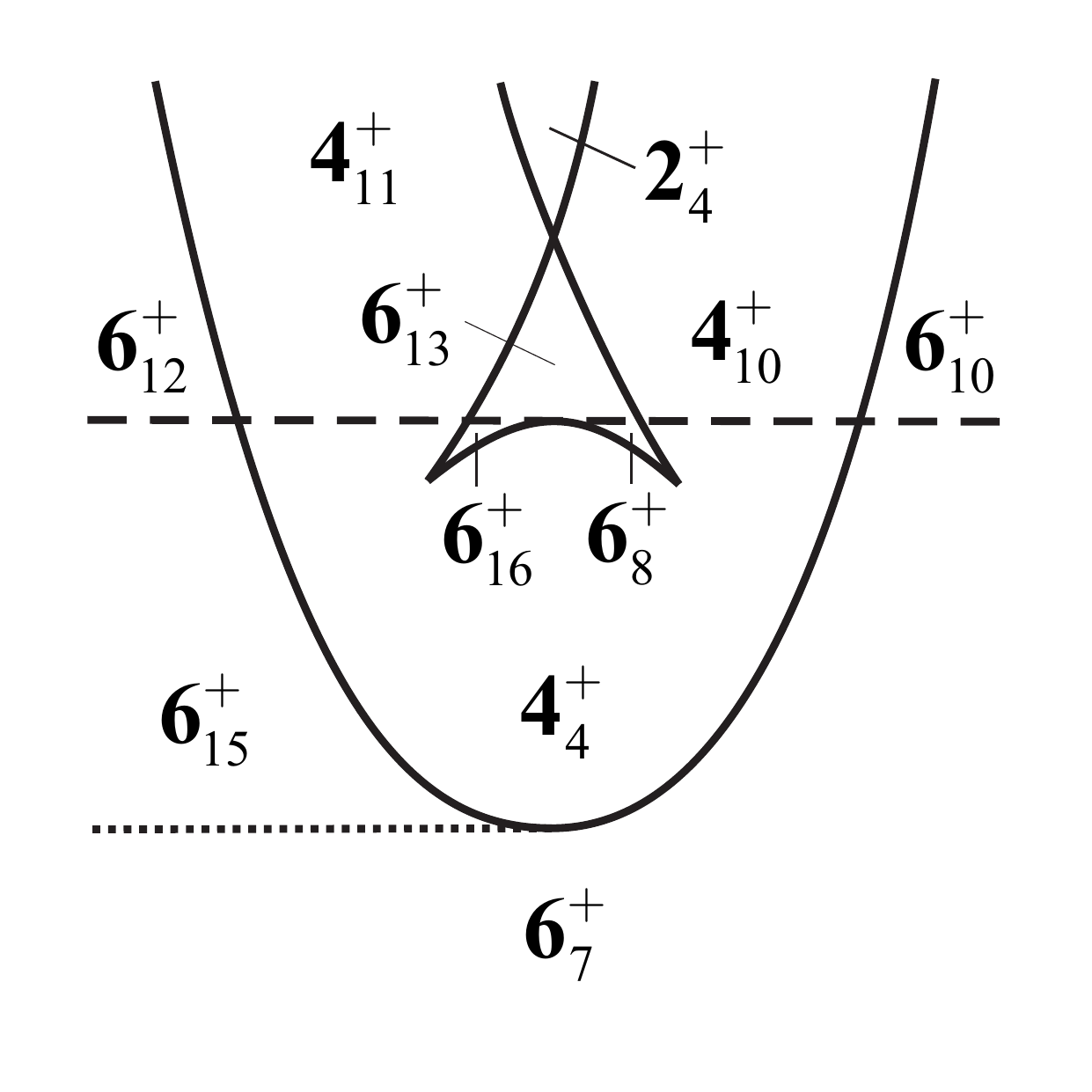}
\end{tabular}
\caption{The skeletons $\Sigma_{0,t_1,q_5}$ for domains $29,33,38,60,93$.}
\label{predel1-4chetv}
\end{center}
\end{figure}

\medskip
\lemma\label{straty-2} {\it For any $i=1,\dots,4$, the set ${\bf 2}_i^{\mathfrak{S}}$ of points $x=(t_1,t_2,S_3,S_4,q_5)$ such that $(t_2,S_3)$ belongs to the stratum ${\bf 2}_{i}\in\mathfrak{S}_{S_4,t_1,q_5}$ if $S_4\neq0$ and to the domain ${\bf 2}_{i}^-$ if $S_4=0$ is a type $A_1^2$ stratum of $\mathfrak{S}$.}

{\sc Proof.} For any point $(t_1,q_5)$ of the domain $\delta t_1<0$, we have the equality ${\bf 2}_{1}^-={\bf 2}_{1}^+$ (see Fig. \ref{predel3-2chetv}). It follows that the set ${\bf 2}_1^{\mathfrak{S}}$ is open, connected, and its boundary belongs to $\Sigma$. The other cases are considered similarly:
${\bf 2}_{2}^-={\bf 2}_{2}^+$ for $\delta t_1>0,q_5>0$; ${\bf 2}_{3}^-={\bf 2}_{3}^+$ for $\delta t_1<0,q_5>0$; ${\bf 2}_{4}^-={\bf 2}_{4}^+$ for $\delta t_1>0$. $\Box$

\medskip
\lemma\label{straty-4-1} {\it The following sets are type $A_1^4$ strata of $\mathfrak{S}$:

${\bf 4}_{1}^{\mathfrak{S}}=\{\,x\,|$ $(t_2,S_3)$ belongs to the union of strata ${\bf 4}_{1},{\bf 4}_{6}\in\mathfrak{S}_{S_4,t_1,q_5}$ if $S_4\neq0$ and to the union of domains ${\bf 4}_{6}^-,{\bf 4}_{6}^+,{\bf 4}_{1}^-\setminus \overline{{\bf 4}_{6}^+}$ if $S_4=0\}$\footnote{Here and below, by $\overline{U}$ we denote the closure of a subset $U\subset\mathbb{R}^{n}$.};

${\bf 4}_{7}^{\mathfrak{S}}=\{\,x\,|$ $(t_2,S_3)$ belongs to the union of strata ${\bf 4}_{7},{\bf 4}_{8}\in\mathfrak{S}_{S_4,t_1,q_5}$ if $S_4\neq0$ and to the union of domains ${\bf 4}_{8}^-,{\bf 4}_{8}^+,{\bf 4}_{7}^-\setminus \overline{{\bf 4}_{8}^+}$ if $S_4=0\}$;

${\bf 4}_{3}^{\mathfrak{S}}=\{\,x\,|$ $(t_2,S_3)$ belongs to the union of strata ${\bf 4}_{3},{\bf 4}_{4}\in\mathfrak{S}_{S_4,t_1,q_5}$ if $S_4\neq0$ and to the union of domains ${\bf 4}_{3}^-,{\bf 4}_{3}^+,{\bf 4}_{4}^-\setminus \overline{{\bf 4}_{3}^+}$ if $S_4=0\}$;

${\bf 4}_{2}^{\mathfrak{S}}=\{\,x\,|$ $(t_2,S_3)$ belongs to the stratum ${\bf 4}_{2}\in\mathfrak{S}_{S_4,t_1,q_5}$ if $\delta S_4<0$, to the union of strata ${\bf 4}_{5},{\bf 4}_{9}$ if $\delta S_4>0$, and to the domain ${\bf 4}_{2}^-$ if $S_4=0\}$;

${\bf 4}_{5}^{\mathfrak{S}}=\{\,x\,|$ $(t_2,S_3)$ belongs to the union of strata ${\bf 4}_{5},{\bf 4}_{9}\in\mathfrak{S}_{S_4,t_1,q_5}$ if $\delta S_4<0$, to the stratum ${\bf 4}_{2}$ if $\delta S_4>0$, and to the domain ${\bf 4}_{2}^+$ if $S_4=0\}$;

${\bf 4}_{10}^{\mathfrak{S}}=\{\,x\,|$ $(t_2,S_3)$ belongs to the stratum ${\bf 4}_{10}\in\mathfrak{S}_{S_4,t_1,q_5}$ if $\delta S_4<0$, to the stratum ${\bf 4}_{11}$ if $\delta S_4>0$, and to the domain ${\bf 4}_{10}^-$ if $S_4=0\}$;

${\bf 4}_{11}^{\mathfrak{S}}=\{\,x\,|$ $(t_2,S_3)$ belongs to the stratum ${\bf 4}_{11}\in\mathfrak{S}_{S_4,t_1,q_5}$ if $\delta S_4<0$, to the stratum ${\bf 4}_{10}$ if $\delta S_4>0$, and to the domain ${\bf 4}_{11}^-$ if $S_4=0\}$.}

{\sc Proof.} We have the inclusions ${\bf 4}_{6}^-\subset{\bf 4}_{1}^+,{\bf 4}_{6}^+\subset{\bf 4}_{1}^-$ and the equality ${\bf 4}_{1}^-\setminus \overline{{\bf 4}_{6}^+}={\bf 4}_{1}^+\setminus\overline{{\bf 4}_{6}^-}$ for any point $(t_1,q_5)$ of the domain $\delta t_1<0,q_5>0$. It now follows that the set ${\bf 4}_{1}^{\mathfrak{S}}$ is open, connected, and its boundary belongs to $\Sigma$. The other cases are considered similarly: ${\bf 4}_{8}^-\subset{\bf 4}_{7}^+,{\bf 4}_{8}^+\subset{\bf 4}_{7}^-$ and ${\bf 4}_{7}^-\setminus \overline{{\bf 4}_{8}^+}={\bf 4}_{7}^+\setminus \overline{{\bf 4}_{8}^-}$ for $\delta t_1>0,q_5>0$; ${\bf 4}_{3}^-\subset{\bf 4}_{4}^+$, ${\bf 4}_{3}^+\subset{\bf 4}_{4}^-$, ${\bf 4}_{4}^-\setminus\overline{{\bf 4}_{3}^+}={\bf 4}_{4}^+\setminus\overline{{\bf 4}_{3}^-}$ for $\delta t_1<0,q_5<0$ and ${\bf 4}_{4}^-={\bf 4}_{4}^+$ for $\delta t_1>0,q_5<0$; ${\bf 4}_{2}^-={\bf 4}_{5}^+,{\bf 4}_{5}^-={\bf 4}_{2}^+$ for $\delta t_1<0,q_5>0$ and ${\bf 4}_{2}^-={\bf 4}_{9}^+,{\bf 4}_{9}^-={\bf 4}_{2}^+$ for $\delta t_1>0,q_5>0$; ${\bf 4}_{10}^-={\bf 4}_{11}^+$ and ${\bf 4}_{11}^-={\bf 4}_{10}^+$ for $\delta t_1>0,q_5<0$. $\Box$

\medskip
\lemma\label{straty-6-1} {\it The following sets are type $A_1^6$ strata of $\mathfrak{S}$:

${\bf 6}_1^{\mathfrak{S}}=\{\,x\,|$ $(t_2,S_3)$ belongs to the union of strata ${\bf 6}_{1},{\bf 6}_{2}\in\mathfrak{S}_{S_4,t_1,q_5}$ if $\delta S_4<0$, to the stratum ${\bf 6}_{4}$ if $\delta S_4>0$, and to the union of domains ${\bf 6}_{1}^-,{\bf 6}_{2}^-$ if $S_4=0\}$;

${\bf 6}_4^{\mathfrak{S}}=\{\,x\,|$ $(t_2,S_3)$ belongs to the stratum ${\bf 6}_{4}\in\mathfrak{S}_{S_4,t_1,q_5}$ if $\delta S_4<0$, to the union of strata ${\bf 6}_{1}$,${\bf 6}_{2}$ if $\delta S_4>0$, and to the union of domains ${\bf 6}_{1}^+,{\bf 6}_{2}^+$ if $S_4=0\}$;

${\bf 6}_{10}^{\mathfrak{S}}=\{\,x\,|$ $(t_2,S_3)$ belongs to the stratum ${\bf 6}_{10}\in\mathfrak{S}_{S_4,t_1,q_5}$ if $\delta S_4<0$, to the union of strata ${\bf 6}_{11},{\bf 6}_{12}$ if $\delta S_4>0$, and to the union of domains ${\bf 6}_{11}^+,{\bf 6}_{12}^+$ if $S_4=0\}$;

${\bf 6}_{11}^{\mathfrak{S}}=\{\,x\,|$ $(t_2,S_3)$ belongs to the union of strata ${\bf 6}_{11}$,${\bf 6}_{12}\in\mathfrak{S}_{S_4,t_1,q_5}$ if $\delta S_4<0$, to the stratum ${\bf 6}_{10}$ if $\delta S_4>0$, and to the union of domains ${\bf 6}_{11}^-,{\bf 6}_{12}^-$ if $S_4=0\}$;

${\bf 6}_{3}^{\mathfrak{S}}=\{\,x\,|$ $(t_2,S_3)$ belongs to the union of strata ${\bf 6}_{3},{\bf 6}_{9}\in\mathfrak{S}_{S_4,t_1,q_5}$ if $S_4\neq0$ and to the union of domains ${\bf 6}_{9}^-$,${\bf 6}_{9}^+$,${\bf 6}_{3}^-\setminus \overline{{\bf 6}_{9}^+}$ if $S_4=0\}$;

${\bf 6}_{5}^{\mathfrak{S}}=\{\,x\,|$ $(t_2,S_3)$ belongs to the union of strata ${\bf 6}_{5},{\bf 6}_{6}\in\mathfrak{S}_{S_4,t_1,q_5}$ if $S_4\neq0$ and to the union of domains ${\bf 6}_{5}^-,{\bf 6}_{5}^+,{\bf 6}_{6}^-\setminus\overline{{\bf 6}_{5}^+}$ if $S_4=0\}$;

${\bf 6}_{7}^{\mathfrak{S}}=\{\,x\,|$ $(t_2,S_3)$ belongs to the union of strata ${\bf 6}_{7},{\bf 6}_{15}\in\mathfrak{S}_{S_4,t_1,q_5}$ if $S_4\neq0$ and to the union of domains ${\bf 6}_{15}^-,{\bf 6}_{15}^+,{\bf 6}_{7}^-\setminus \overline{{\bf 6}_{15}^+}$ if $S_4=0\}$;

${\bf 6}_8^{\mathfrak{S}}=\{\,x\,|$ $(t_2,S_3)$ belongs to the stratum ${\bf 6}_{8}\in\mathfrak{S}_{S_4,t_1,q_5}$ if $\delta S_4<0$, to the stratum ${\bf 6}_{16}$ if $\delta S_4>0$, and to the domain ${\bf 6}_{8}^-$ if $S_4=0\}$;

${\bf 6}_{16}^{\mathfrak{S}}=\{\,x\,|$ $(t_2,S_3)$ belongs to the stratum ${\bf 6}_{16}\in\mathfrak{S}_{S_4,t_1,q_5}$ if $\delta S_4<0$, to the stratum ${\bf 6}_{8}$ if $\delta S_4>0$, and to the domain ${\bf 6}_{16}^-$ if $S_4=0\}$;

${\bf 6}_{13}^{\mathfrak{S}}=\{\,x\,|$ $(t_2,S_3)$ belongs to the stratum ${\bf 6}_{13}\in\mathfrak{S}_{S_4,t_1,q_5}$ if $S_4\neq0$ and to the domain ${\bf 6}_{13}^-$ if $S_4=0\}$.}

{\sc Proof.} We have the inclusions ${\bf 6}_1^-\cup{\bf 6}_2^-\subset{\bf 6}_4^+,{\bf 6}_4^-\supset{\bf 6}_1^+\cup{\bf 6}_2^+$ for any point $(t_1,q_5)$ of the domain $\delta t_1<0,q_5<\sqrt{\frac32|t_1|}$. It now follows that the sets ${\bf 6}_{1}^{\mathfrak{S}}$ and ${\bf 6}_{4}^{\mathfrak{S}}$ are open, connected and their boundaries belong to $\Sigma$. The other cases are considered similarly: ${\bf 6}_{10}^-\supset{\bf 6}_{11}^+\cup{\bf 6}_{12}^+,{\bf 6}_{11}^-\cup{\bf 6}_{12}^-\subset{\bf 6}_{10}^+$ for $\delta t_1>0,q_5<2\sqrt{3|t_1|}$; ${\bf 6}_{9}^-\subset{\bf 6}_{3}^+,{\bf 6}_{3}^-\supset{\bf 6}_{9}^+$ and ${\bf 6}_{3}^-\setminus \overline{{\bf 6}_{9}^+}={\bf 6}_{3}^+\setminus \overline{{\bf 6}_{9}^-}$ for $\delta t_1<0,q_5<2\sqrt{3|t_1|}$; ${\bf 6}_{5}^-\subset{\bf 6}_{6}^+,{\bf 6}_{6}^-\supset{\bf 6}_{5}^+$ and ${\bf 6}_{6}^-\setminus \overline{{\bf 6}_{5}^+}={\bf 6}_{6}^+\setminus \overline{{\bf 6}_{5}^-}$ for $\delta t_1<0,q_5<-2\sqrt{3|t_1|}$; ${\bf 6}_{7}^-\supset{\bf 6}_{15}^+,{\bf 6}_{15}^-\subset{\bf 6}_{7}^+$ and ${\bf 6}_{7}^-\setminus \overline{{\bf 6}_{15}^+}={\bf 6}_{7}^+\setminus \overline{{\bf 6}_{15}^-}$ for $\delta t_1>0,q_5<2\sqrt{3|t_1|}$; ${\bf 6}_{8}^-={\bf 6}_{16}^+,{\bf 6}_{16}^-={\bf 6}_{8}^+,{\bf 6}_{13}^-={\bf 6}_{13}^+$ for $\delta t_1>0,q_5<-2\sqrt{3|t_1|}$. $\Box$

\medskip
\lemma\label{straty-14} {\it The sets ${\bf 6}_{14}^{\mathfrak{S},\pm}$ of points $x$ such that $(t_2,S_3)$ belongs to the stratum ${\bf 6}_{14}\in\mathfrak{S}_{S_4,t_1,q_5}$ and $\pm \delta S_4>0$, respectively, are type $A_1^6$ strata of $\mathfrak{S}$.}

{\sc Proof.} By Lemma \ref{sech}, it follows that the sets ${\bf 6}_{14}^{\mathfrak{S},\pm}$ are open, connected, and their boundary points that do not lie in the hyperplane $S_4=0$ belong to $\Sigma$. It remains to note that if a limit point $x$ of the set ${\bf 6}_{14}^{\mathfrak{S},\pm}$ lies in the hyperplane $S_4=0$, then $(t_1,q_5)$ belongs to the boundary of domain $93$ in Fig. \ref{razbienie0}, i.e. $x\in\Sigma$. $\Box$

\medskip
Lemmas \ref{straty-2} -- \ref{straty-14} imply:

\medskip
\predlo\label{straty-vse} {\it The stratification $\mathfrak{S}$ has $23$ five-dimensional strata. Namely it has four strata ${\bf 2}_1^{\mathfrak{S}} - {\bf 2}_4^{\mathfrak{S}}$ of type $A_1^2$, seven strata ${\bf 4}_1^{\mathfrak{S}} - {\bf 4}_{3}^{\mathfrak{S}},{\bf 4}_{5}^{\mathfrak{S}},{\bf 4}_{7}^{\mathfrak{S}},{\bf 4}_{10}^{\mathfrak{S}},{\bf 4}_{11}^{\mathfrak{S}}$ of type $A_1^4$ and twelve strata ${\bf 6}_{1}^{\mathfrak{S}},{\bf 6}_{3}^{\mathfrak{S}} - {\bf 6}_{5}^{\mathfrak{S}},{\bf 6}_{7}^{\mathfrak{S}},{\bf 6}_{8}^{\mathfrak{S}},{\bf 6}_{10}^{\mathfrak{S}},{\bf 6}_{11}^{\mathfrak{S}},{\bf 6}_{13}^{\mathfrak{S}}$, ${\bf 6}_{14}^{\mathfrak{S},\pm},{\bf 6}_{16}^{\mathfrak{S}}$ of type $A_1^6$.}

\section{Local properties of the complement to $\Sigma$}

This paragraph is a continuation of the section \ref{lok-property-1}. Now we describe local properties of the complement to the skeleton $\Sigma$ in a neighbourhood of a point $x=(t_1,t_2,S_3,S_4,q_5)$ such that $S_4=0$. The number $\ast$ in the formulas $(t_2,S_3)\in\overline{{\bf k}_{\ast}^0}$ below can be the subscript of any stratum ${\bf k}_{\ast}^{\mathfrak{S}}$ from the list in Proposition \ref{straty-vse}.

\medskip
\lemma\label{lok-0-1} {\it Let $S_4=0,q_5t_1\neq0$, and $(t_2,S_3)$ be a nonsingular point of $\Sigma_{0,t_1,q_5}$ that does not belong to the lines $(\ref{crit-tochki}),(\ref{peresech-3})$. Then:

{\rm1)} if $(t_2,S_3)\in\overline{{\bf 2}_{\ast}^0}$, then $f$ is non-singular at $x$ and has a multisingularity of type $A_2A_1^2$ at $f(x)$; the set $\Sigma^c(x)$ consists of two contractible components; there is one component of each type $A_1^2,A_1^4$;

{\rm2)} if $(t_2,S_3)\in\overline{{\bf 6}_{\ast}^0}$, then $f$ is non-singular at $x$ and has a multisingularity of type $A_2A_1^4$ at $f(x)$; the set $\Sigma^c(x)$ consists of two contractible components; there is one component of each type $A_1^4,A_1^6$.}

\medskip
\lemma\label{lok-0-2} {\it Let $S_4=0,q_5t_1\neq0$, and $(t_2,S_3)\in(\ref{peresech-3})$ be a nonsingular point of $\Sigma_{0,t_1,q_5}$. Then:

{\rm1)} if $(t_2,S_3)\in\overline{{\bf 4}_{\ast}^0}$, then $f$ is non-singular at $x$ and has a multisingularity of type $A_2A_1^2$ at $f(x)$; the set $\Sigma^c(x)$ consists of two contractible components; there is one component of each type $A_1^2,A_1^4$;

{\rm2)} if $(t_2,S_3)\in\overline{{\bf 6}_{\ast}^0}$, then $f$ is non-singular at $x$ and has a multisingularity of type $A_2A_1^4$ at $f(x)$; the set $\Sigma^c(x)$ consists of two contractible components; there is one component of each type $A_1^4,A_1^6$.}

\medskip
\lemma\label{lok-0-3} {\it Let $S_4=0,q_5t_1\neq0$, and $(t_2,S_3)\in(\ref{crit-tochki})$ be a nonsingular point of $\Sigma_{0,t_1,q_5}$. Then:

{\rm1)} if $(t_2,S_3)\in\overline{{\bf 4}_{\ast}^0}$, then $f$ has a singularity of type $A_2$ at $x$ and a multisingularity of type $A_2A_1^2$ at $f(x)$; the set $\Sigma^c(x)$ consists of two contractible components of type $A_1^4$; they belong to different strata of $\mathfrak{S}$;

{\rm2)} if $(t_2,S_3)\in\overline{{\bf 6}_{\ast}^0}$, then $f$ has a singularity of type $A_2$ at $x$ and a multisingularity of type $A_2A_1^4$ at $f(x)$; the set $\Sigma^c(x)$ consists of two contractible components of type $A_1^6$; they belong to different strata of $\mathfrak{S}$.}

\medskip
\lemma\label{lok-0-4} {\it Let $S_4=0,q_5t_1\neq0$, and $(t_2,S_3)$ be a semicubical cusp of the curve {\rm(\ref{psi-2})}. Then:

{\rm1)} if $(t_2,S_3)$ does not belong to the lines $(\ref{crit-tochki}),(\ref{peresech-3})$, then $f$ is non-singular at $x$ and has a multisingularity of type $A_3^\pm A_1^3$ at $f(x)$; the set $\Sigma^c(x)$ consists of two contractible components; there is one component of each type $A_1^4,A_1^6$;

{\rm2)} if $(t_2,S_3)$ belongs to the line $(\ref{peresech-3})$, then $f$ is non-singular at $x$ and has a multisingularity of type $A_3^\pm A_2A_1$ at $f(x)$; the set $\Sigma^c(x)$ consists of $1,2,1$ contractible components of types $A_1^2,A_1^4,A_1^6$, respectively; the components of type $A_1^4$ belong to the different strata ${\bf 4}_7^{\mathfrak{S}},{\bf 4}_{10}^{\mathfrak{S}}$ for $\delta t_2<0$ and to ${\bf 4}_7^{\mathfrak{S}},{\bf 4}_{11}^{\mathfrak{S}}$ for $\delta t_2>0$;

{\rm3)} if $(t_2,S_3)$ belongs to the line $(\ref{crit-tochki})$, then $f$ has a singularity of type $A_2$ at $x$ and a multisingularity of type $A_3^\pm A_2A_1$ at $f(x)$; the set $\Sigma^c(x)$ consists of four contractible components; there are two components of each type $A_1^4,A_1^6$; they belong to different strata of $\mathfrak{S}$.}

\medskip
\lemma\label{lok-0-5} {\it Let $S_4=0,q_5t_1\neq0$, and $(t_2,S_3)$ be a self-intersection point of the curve {\rm(\ref{psi-2})}. Then $f$ is non-singular at $x$ and has a multisingularity of type $A_2^2A_1^2$ at $f(x)$. The set $\Sigma^c(x)$ consists of $1,2,1$ contractible components of types $A_1^2,A_1^4,A_1^6$, respectively. If $\delta t_1>0$, then the components of type $A_1^4$ belong to different strata of $\mathfrak{S}$. If $\delta t_1<0$, then they belong to the same stratum ${\bf 4}_{1}^{\mathfrak{S}}$ for $q_5>0$ and to ${\bf 4}_{3}^{\mathfrak{S}}$ for $q_5<0$.}

\medskip
\lemma\label{lok-0-6} {\it Let $S_4=0,q_5t_1\neq0$, and $(t_2,S_3)\in(\ref{peresech-3})$ be a singular point of $\Sigma_{0,t_1,q_5}$. Then:

{\rm1)} if $(t_2,S_3)$ is a transversal intersection point of the line $(\ref{peresech-3})$ with a smooth branch of the curve {\rm(\ref{psi-2})}, then $f$ is non-singular at $x$ and has a multisingularity of type $A_2^2A_1^2$ at $f(x)$; the set $\Sigma^c(x)$ consists of $1,2,1$ contractible components of types $A_1^2,A_1^4,A_1^6$, respectively; the components of type $A_1^4$ belong to the different strata ${\bf 4}_{2}^{\mathfrak{S}},{\bf 4}_3^{\mathfrak{S}}$ for $\delta t_1<0,\delta t_2<0,q_5<0$; to ${\bf 4}_{3}^{\mathfrak{S}},{\bf 4}_5^{\mathfrak{S}}$ for $\delta t_1<0,\delta t_2>0,q_5<0$; to ${\bf 4}_{7}^{\mathfrak{S}},{\bf 4}_{10}^{\mathfrak{S}}$ for $\delta t_1>0,\delta t_2<0,q_5>0$; to ${\bf 4}_{7}^{\mathfrak{S}},{\bf 4}_{11}^{\mathfrak{S}}$ for $\delta t_1>0,\delta t_2>0,q_5>0$; and to the same stratum ${\bf 4}_{1}^{\mathfrak{S}}$ for $\delta t_1<0,q_5>0$;

{\rm2)} if $(t_2,S_3)\in\overline{{\bf 2}_{\ast}^0}$ is a simple tangency point of the line $(\ref{peresech-3})$ with the closure of the curve {\rm(\ref{psi-2})}, then $f$ is non-singular at $x$ and has a multisingularity of type $A_3^\pm A_1$ at $f(x)$; the set $\Sigma^c(x)$ consists of two contractible components; there is one component of each type $A_1^2,A_1^4$;

{\rm3)} if $(t_2,S_3)\in\overline{{\bf 6}_{\ast}^0}$ is a simple tangency point of the line $(\ref{peresech-3})$ with the closure of the curve {\rm(\ref{psi-2})}, then $f$ is non-singular at $x$ and has a multisingularity of type $A_3^\pm A_1^3$ at $f(x)$; the set $\Sigma^c(x)$ consists of two contractible components; there is one component of each type $A_1^4,A_1^6$;

{\rm4)} if $t_2=0,q_5^2=-12\delta t_1$, then $f$ is non-singular at $x$ and has a multisingularity of type $A_5^\pm A_1$ at $f(x)$; the set $\Sigma^c(x)$ consists of three contractible components; there is one component of each type $A_1^2,A_1^4,A_1^6$.}

\medskip
\lemma\label{lok-0-7} {\it Let $S_4=0,q_5t_1\neq0$, $(t_2,S_3)\in(\ref{crit-tochki})$ be a singular point of $\Sigma_{0,t_1,q_5}$. Then:

{\rm1)} if $(t_2,S_3)$ is a transversal intersection point of the line $(\ref{crit-tochki})$ with a smooth branch of the curve
{\rm(\ref{psi-2})}, then $f$ has a singularity of type $A_2$ at $x$ and multisingularity of type $A_2^2A_1^2$ at $f(x)$; the set $\Sigma^c(x)$ consists of four contractible components; there are two components of each type $A_1^4,A_1^6$; they belong to different strata of $\mathfrak{S}$;

{\rm2)} if $(t_2,S_3)\in\overline{{\bf 2}_{\ast}^0}$ is a simple tangency point of the line $(\ref{crit-tochki})$ with a smooth branch of the curve {\rm(\ref{psi-2})}, then $f$ has a singularity of type $A_3^\pm$ at $x$ and a multisingularity of type $A_3^\pm A_1$ at $f(x)$; the set $\Sigma^c(x)$ consists of four contractible components; there are one component of type $A_1^2$ and three components of type $A_1^4$; they belong to different strata of $\mathfrak{S}$;

{\rm3)} if $(t_2,S_3)\in\overline{{\bf 6}_{\ast}^0}$ is a simple tangency point of the line $(\ref{crit-tochki})$ with a smooth branch of the curve {\rm(\ref{psi-2})}, then $f$ has a singularity of type $A_3^\pm$ at $x$ and a multisingularity of type $A_3^\pm A_1^3$ at $f(x)$; the set $\Sigma^c(x)$ consists of four contractible components; there are one component of type $A_1^4$ and three components of type $A_1^6$; they belong to different strata of $\mathfrak{S}$;

{\rm4)} if $t_2=0,q_5^2=12\delta t_1$, then $f$ has a singularity of type $A_5^\pm$ at $x$ and a multisingularity of type $A_5^\pm A_1$ at $f(x)$; the set $\Sigma^c(x)$ consists of $1,3,5$ contractible components of types $A_1^2,A_1^4,A_1^6$, respectively; they belong to the different strata ${\bf 2}_2^{\mathfrak{S}},{\bf 4}_2^{\mathfrak{S}},{\bf 4}_5^{\mathfrak{S}},{\bf 4}_7^{\mathfrak{S}},{\bf 6}_1^{\mathfrak{S}},{\bf 6}_4^{\mathfrak{S}},{\bf 6}_7^{\mathfrak{S}},{\bf 6}_{10}^{\mathfrak{S}},{\bf 6}_{11}^{\mathfrak{S}}$ for $q_5>0$ and to ${\bf 2}_4^{\mathfrak{S}},{\bf 4}_3^{\mathfrak{S}},{\bf 4}_{10}^{\mathfrak{S}}$, ${\bf 4}_{11}^{\mathfrak{S}},{\bf 6}_8^{\mathfrak{S}},{\bf 6}_{13}^{\mathfrak{S}},{\bf 6}_{14}^{\mathfrak{S},\pm},{\bf 6}_{16}^{\mathfrak{S}}$ for $q_5<0$.}

\medskip
\lemma\label{lok-0-8} {\it Let $S_4=q_5=0,t_1\neq0$. Then:

{\rm1)} if $\delta S_3=9t_2^2,t_2\neq0$, then
$f$ is non-singular at $x$ and has a multisingularity of type $A_2^2A_1^2$ at $f(x)$; the set $\Sigma^c(x)$ consists of $1,2,1$ contractible components of types $A_1^2,A_1^4,A_1^6$, respectively; the components of type $A_1^4$ belong to the different strata ${\bf 4}_{2}^{\mathfrak{S}},{\bf 4}_3^{\mathfrak{S}}$ for $\delta t_1<0,\delta t_2<0$; to ${\bf 4}_{3}^{\mathfrak{S}},{\bf 4}_5^{\mathfrak{S}}$ for $\delta t_1<0,\delta t_2>0$; to ${\bf 4}_{7}^{\mathfrak{S}},{\bf 4}_{10}^{\mathfrak{S}}$ for $\delta t_1>0,\delta t_2<0$; and to ${\bf 4}_{7}^{\mathfrak{S}},{\bf 4}_{11}^{\mathfrak{S}}$ for $\delta t_1>0,\delta t_2>0$;

{\rm2)} if $S_3=0,t_2\neq0$, then $f$ has a singularity of type $A_2$ at $x$ and a multisingularity of type $A_2^2A_1^2$ at $f(x)$; the set $\Sigma^c(x)$ consists of four contractible components; there are two components of each type $A_1^4,A_1^6$; they belong to different strata of $\mathfrak{S}$; the components of type $A_1^4$ belong to strata ${\bf 4}_1^{\mathfrak{S}},{\bf 4}_{2}^{\mathfrak{S}}$ for $\delta t_1<0,\delta t_2<0$; to ${\bf 4}_1^{\mathfrak{S}},{\bf 4}_{5}^{\mathfrak{S}}$ for $\delta t_1<0,\delta t_2>0$; to ${\bf 4}_3^{\mathfrak{S}},{\bf 4}_{10}^{\mathfrak{S}}$ for $\delta t_1>0,\delta t_2<0$; and to ${\bf 4}_3^{\mathfrak{S}},{\bf 4}_{11}^{\mathfrak{S}}$ for $\delta t_1>0,\delta t_2>0$;

{\rm3)} if $S_3=t_2=0$, then $f$ has a singularity of type $A_3^\delta$ at $x$ and a multisingularity of type $(A_3^\delta)^2$ at $f(x)$; the set $\Sigma^c(x)$ consists of $1,4,3$ contractible components of types $A_1^2,A_1^4,A_1^6$, respectively; they belong to different strata of $\mathfrak{S}$; the components of type $A_1^4$ belong to strata ${\bf 4}_1^{\mathfrak{S}},{\bf 4}_{2}^{\mathfrak{S}},{\bf 4}_3^{\mathfrak{S}},{\bf 4}_{5}^{\mathfrak{S}}$ for $\delta t_1<0$ and to ${\bf 4}_3^{\mathfrak{S}},{\bf 4}_7^{\mathfrak{S}},{\bf 4}_{10}^{\mathfrak{S}},{\bf 4}_{11}^{\mathfrak{S}}$ for $\delta t_1>0$.}

\medskip
\remark The limit positions of the skeleton $\Sigma_{0,t_1,q_5}$ for fixed $q_5$ as $t_1\rightarrow\pm0$ are shown by solid and dotted lines in Fig. \ref{S4=t1=0}. The dotted line represents the line $S_3=0$. The solid line shows parabolas
\begin{equation}
S_3=9\delta \left(t_2\pm\frac16\sqrt{-\frac{q_5^3}{3}}\right)^2.
\label{paraboly}
\end{equation}
Their vertices are denoted by bold points. In the case $q_5=0$, parabolas (\ref{paraboly}) coincide. There is only the line $S_3=0$ for $q_5>0$. The limit positions of two-dimensional strata of $\mathfrak{S}_{0,t_1,q_5}$ are denoted by the same symbols as the strata themselves.

\begin{figure}
\begin{center}
\begin{tabular}{ccc}
&\includegraphics[width=4cm]{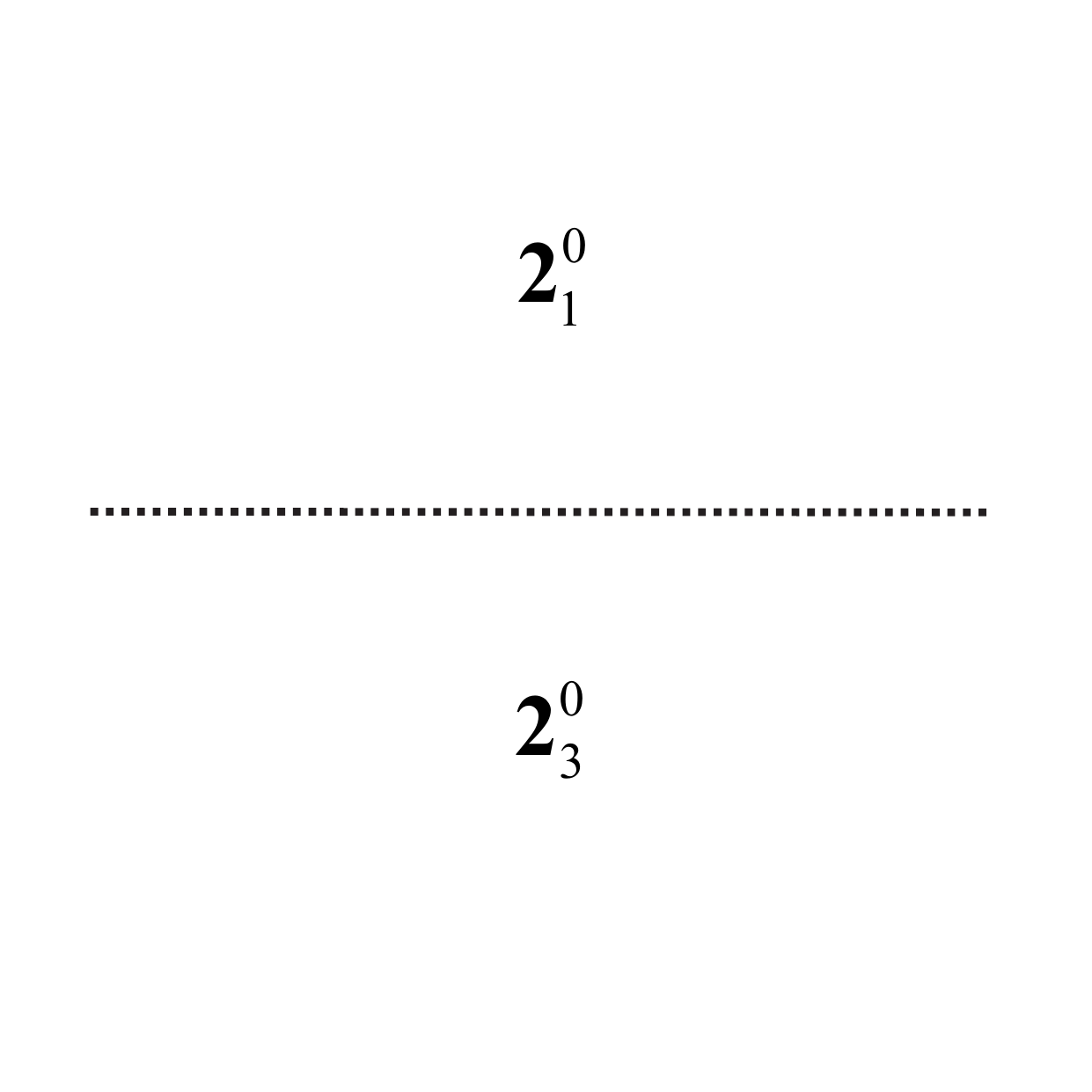}&
\includegraphics[width=4cm]{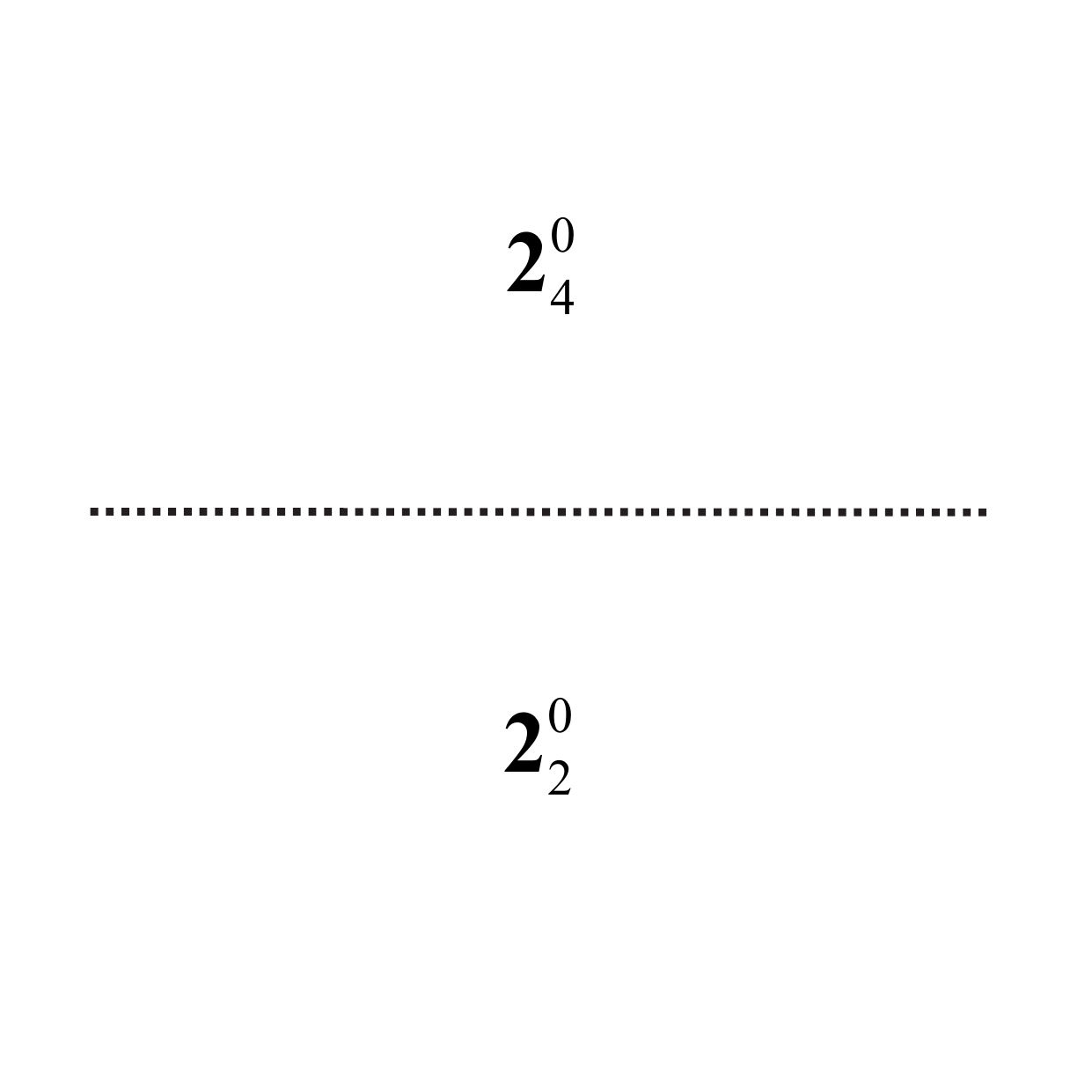}\\
$q_5>0:$&$-0$&$+0$
\\
\\
&\includegraphics[width=4cm]{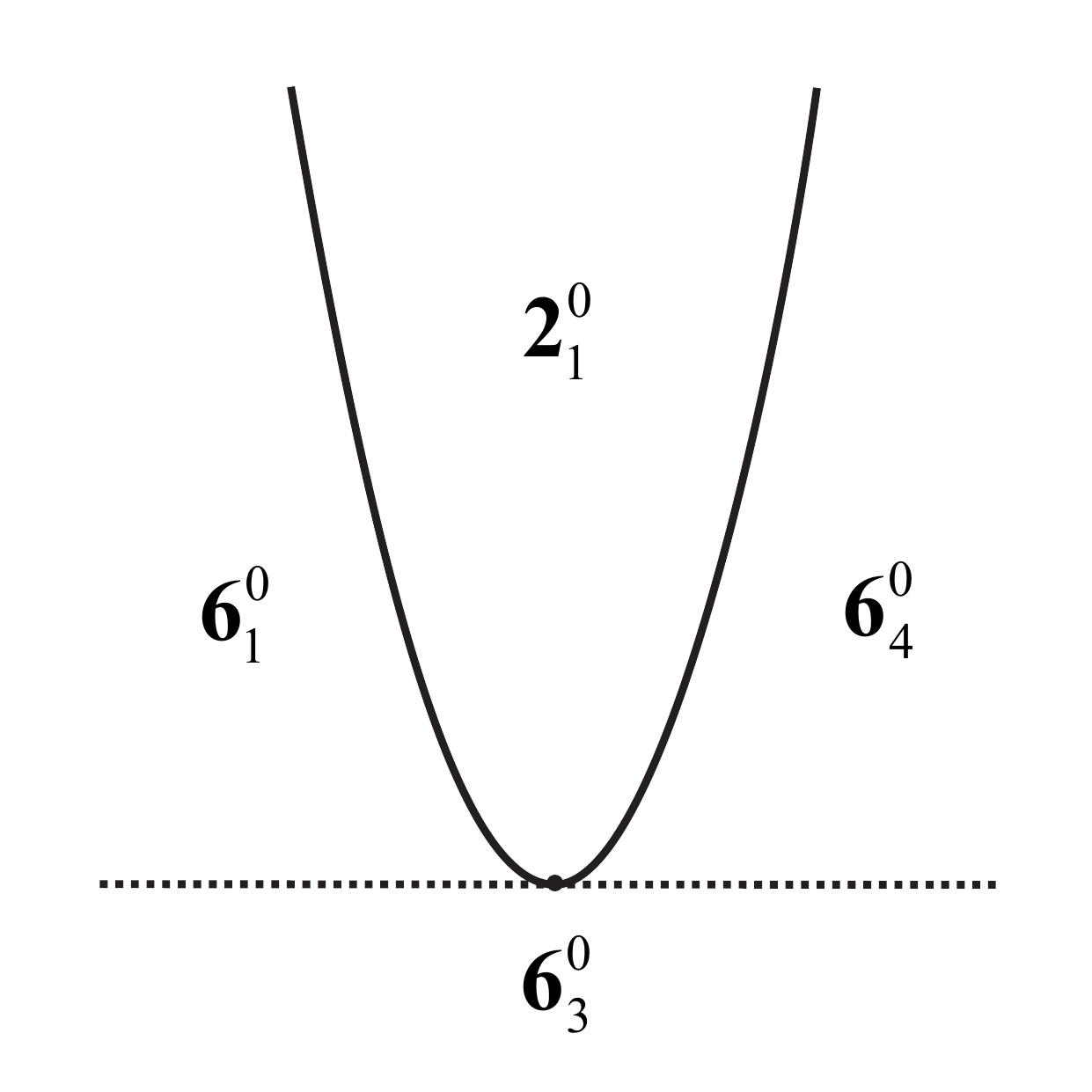}&
\includegraphics[width=4cm]{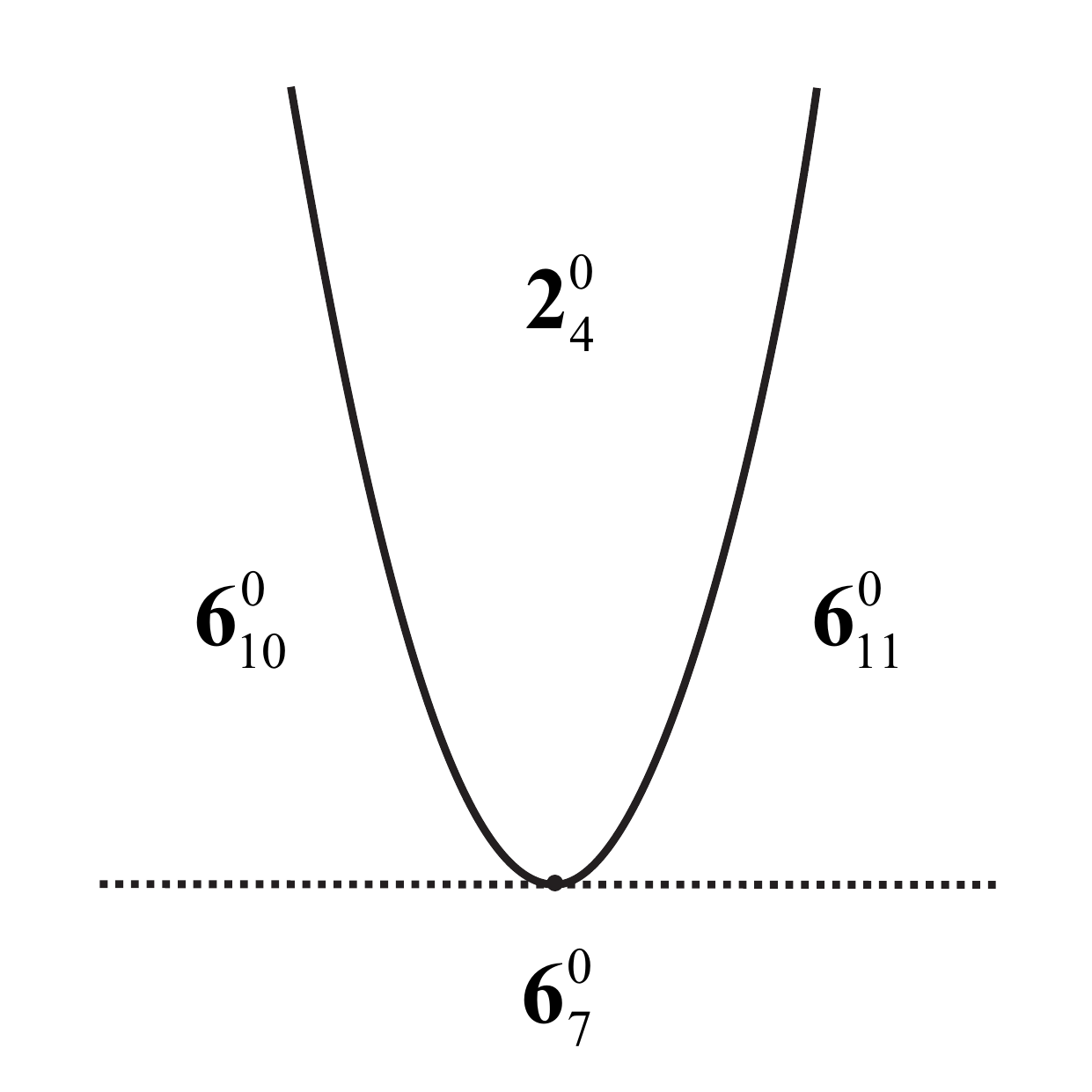}\\
\\
$q_5=0:$&$-0$&$+0$
\\
\\
&\includegraphics[width=4cm]{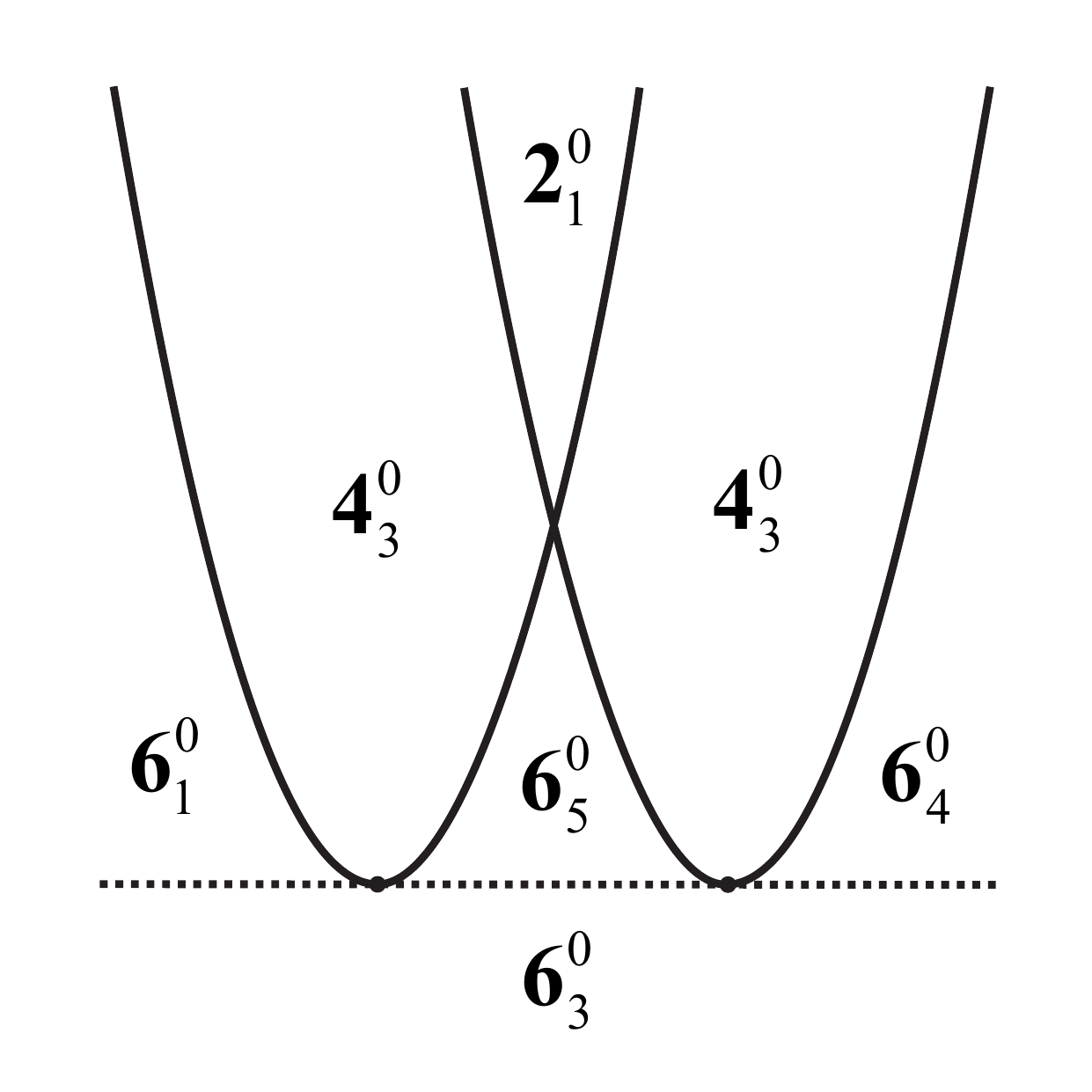}&
\includegraphics[width=4cm]{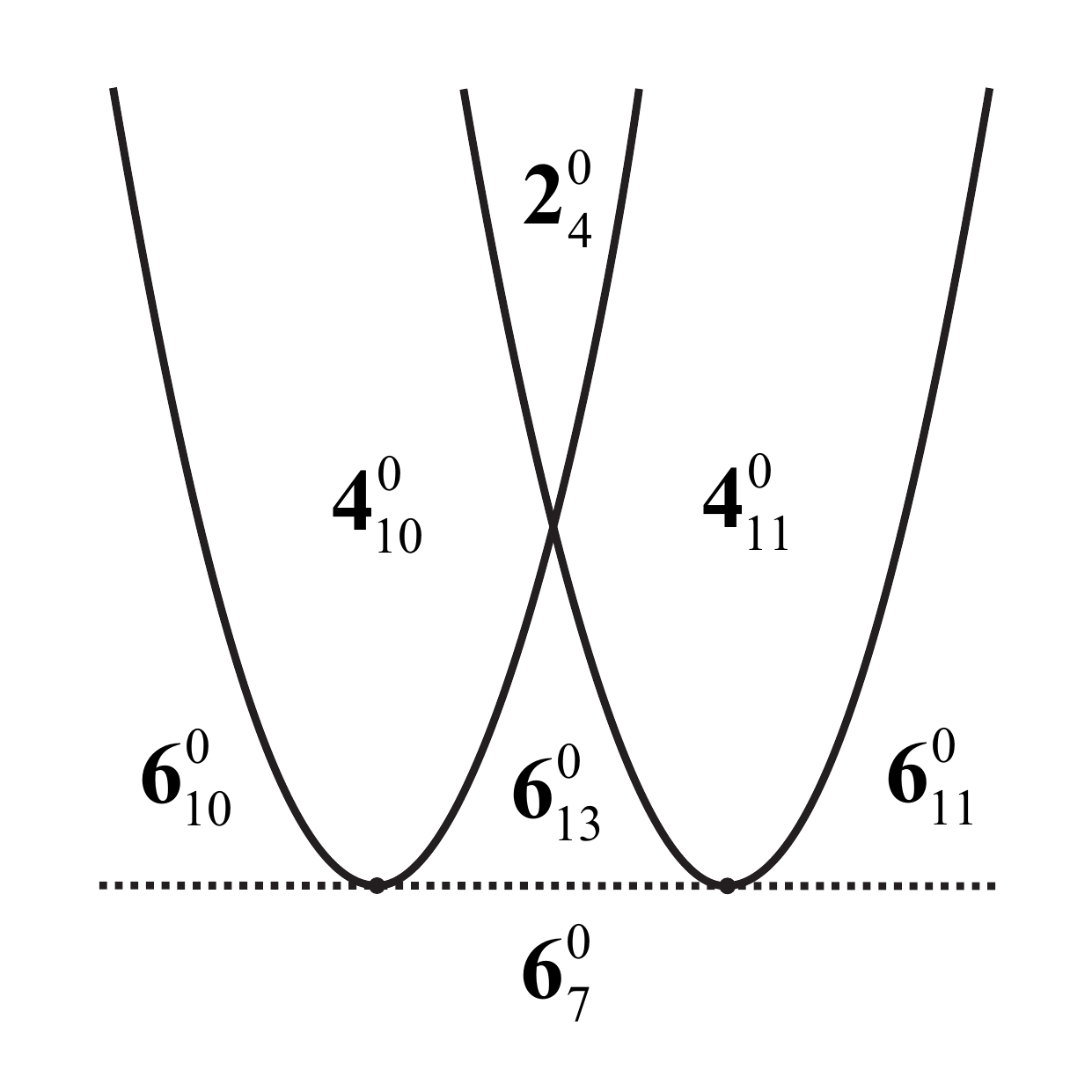}\\
\\
$q_5<0:$&$-0$&$+0$
\\
\\
\end{tabular}
\caption{The limit positions of $\Sigma_{0,t_1,q_5},q_5=\mathrm{const}$ as $\delta t_1\rightarrow\pm0$.}
\label{S4=t1=0}
\end{center}
\end{figure}

\medskip
\lemma\label{lok-0-9} {\it Let $S_4=t_1=0,q_5>0$. Then:

{\rm1)} if $S_3\neq0$, then $f$ has a singularity of type $A_2$ at $x$ and a multisingularity of type $A_2$ at $f(x)$; the set $\Sigma^c(x)$ consists of two contractible components of type $A_1^2$; they belong to different strata ${\bf 2}_1^{\mathfrak{S}},{\bf 2}_4^{\mathfrak{S}}$ for $\delta S_3>0$ and to ${\bf 2}_2^{\mathfrak{S}},{\bf 2}_3^{\mathfrak{S}}$ for $\delta S_3<0$;

{\rm2)} if $S_3=0$, then $f$ has a singularity of type $D_4^{+}$ at $x$ and a monosingularity of type $D_4^{+}$ at $f(x)$; the set $\Sigma^c(x)$ consists of eight contractible components; there are four components of each type $A_1^2,A_1^4$; they belong to different strata ${\bf 2}_1^{\mathfrak{S}} - {\bf 2}_4^{\mathfrak{S}},{\bf 4}_1^{\mathfrak{S}},{\bf 4}_{2}^{\mathfrak{S}},{\bf 4}_5^{\mathfrak{S}},{\bf 4}_{7}^{\mathfrak{S}}$.}

\medskip
\lemma\label{lok-0-10} {\it Let $S_4=t_1=q_5=0,(t_2,S_3)\neq0$. Then:

{\rm1)} if $\delta S_3>9t_2^2$, then $f$ has a singularity of type $A_2$ at $x$ and a monosingularity of type $A_2$ at $f(x)$; the set $\Sigma^c(x)$ consists of two contractible components of type $A_1^2$; they belong to different strata ${\bf 2}_1^{\mathfrak{S}},{\bf 2}_4^{\mathfrak{S}}$;

{\rm2)} if $\delta S_3<9t_2^2,S_3\neq0$, then $f$ has a singularity of type $A_2$ at $x$ and a multisingularity of type $A_2^3$ at $f(x)$; the set $\Sigma^c(x)$ consists of $2,4,2$ contractible components of types $A_1^2,A_1^4,A_1^6$, respectively; in the case $\delta S_3<0$, they belong to strata ${\bf 2}_{2}^{\mathfrak{S}},{\bf 2}_{3}^{\mathfrak{S}},{\bf 4}_1^{\mathfrak{S}},{\bf 4}_2^{\mathfrak{S}},{\bf 4}_5^{\mathfrak{S}}, {\bf 6}_{3}^{\mathfrak{S}},{\bf 6}_{7}^{\mathfrak{S}}$, where ${\bf 4}_1^{\mathfrak{S}}$ contains two components and each of the others contains exactly one; in the case $\delta S_3>0$, the components belong to different strata ${\bf 2}_{1}^{\mathfrak{S}},{\bf 2}_{4}^{\mathfrak{S}},{\bf 4}_2^{\mathfrak{S}},{\bf 4}_3^{\mathfrak{S}},{\bf 4}_7^{\mathfrak{S}},{\bf 4}_{10}^{\mathfrak{S}},{\bf 6}_{1}^{\mathfrak{S}},{\bf 6}_{10}^{\mathfrak{S}}$ for $\delta t_2<0$ and to ${\bf 2}_{1}^{\mathfrak{S}},{\bf 2}_{4}^{\mathfrak{S}},{\bf 4}_3^{\mathfrak{S}},{\bf 4}_5^{\mathfrak{S}},{\bf 4}_7^{\mathfrak{S}},{\bf 4}_{11}^{\mathfrak{S}},{\bf 6}_{4}^{\mathfrak{S}},{\bf 6}_{11}^{\mathfrak{S}}$ for $\delta t_2>0$;

{\rm3)} if $S_3=9\delta t_2^2,t_2\neq0$, then $f$ has a singularity of type $A_2$ at $x$ and a multisingularity of type $D_4^{+}A_2$ at $f(x)$; the set $\Sigma^c(x)$ consists of $2,4,2$ contractible components of types $A_1^2,A_1^4,A_1^6$, respectively; they belong to different strata ${\bf 2}_{1}^{\mathfrak{S}},{\bf 2}_{4}^{\mathfrak{S}},{\bf 4}_{2}^{\mathfrak{S}},{\bf 4}_{3}^{\mathfrak{S}},{\bf 4}_{7}^{\mathfrak{S}},{\bf 4}_{10}^{\mathfrak{S}},{\bf 6}_{1}^{\mathfrak{S}},{\bf 6}_{10}^{\mathfrak{S}}$ for $\delta t_2<0$ and to ${\bf 2}_{1}^{\mathfrak{S}},{\bf 2}_{4}^{\mathfrak{S}},{\bf 4}_{3}^{\mathfrak{S}},{\bf 4}_{5}^{\mathfrak{S}},{\bf 4}_{7}^{\mathfrak{S}},{\bf 4}_{11}^{\mathfrak{S}},{\bf 6}_{4}^{\mathfrak{S}},{\bf 6}_{11}^{\mathfrak{S}}$ for $\delta t_2>0$;

{\rm4)} if $S_3=0,t_2\neq0$, then $f$ has a singularity of type $D_4^{+}$ at $x$ and a multisingularity of type $D_4^{+}A_2$ at $f(x)$; the set $\Sigma^c(x)$ consists of $4,8,4$ contractible components of types $A_1^2,A_1^4,A_1^6$, respectively; they belong to strata ${\bf 2}_1^{\mathfrak{S}} - {\bf 2}_4^{\mathfrak{S}},{\bf 4}_1^{\mathfrak{S}}-{\bf 4}_3^{\mathfrak{S}},{\bf 4}_{5}^{\mathfrak{S}},{\bf 4}_{7}^{\mathfrak{S}},{\bf 6}_{3}^{\mathfrak{S}},{\bf 6}_{7}^{\mathfrak{S}}$ as well as to ${\bf 4}_{10}^{\mathfrak{S}},{\bf 6}_{1}^{\mathfrak{S}},{\bf 6}_{10}^{\mathfrak{S}}$ for $\delta t_2<0$ and to ${\bf 4}_{11}^{\mathfrak{S}},{\bf 6}_{4}^{\mathfrak{S}},{\bf 6}_{11}^{\mathfrak{S}}$ for $\delta t_2>0$; here each of strata ${\bf 4}_1^{\mathfrak{S}},{\bf 4}_2^{\mathfrak{S}}$ for $\delta t_2<0$ and ${\bf 4}_1^{\mathfrak{S}},{\bf 4}_5^{\mathfrak{S}}$ for $\delta t_2>0$ contains two components; each of the others contains exactly one.}

\medskip
\remark 1) The item 2 of Lemma \ref{lok-0-10} follow from the description of connected components of the complement to the set $\Sigma_{S_4,t_1,q_5},S_4\neq0$ at a neighbourhood of a transversal intersection point of two smooth branches of the curve (\ref{psi-2}) (see Fig. \ref{zona1-15} -- \ref{sing-sigma}). For example, type $A_1^4$ components of the set $\Sigma^c(x)$ for $\delta=1$ and $x=(0,0,-1,0,0)$ are numbered by the points
$$
(-0.05;0;-1;\pm0.1;0)\in{\bf 4}_1^{\mathfrak{S}},\quad (0.05;0;-1;-0.1;0)\in{\bf 4}_2^{\mathfrak{S}},\quad (0.05;0;-1;0.1;0)\in{\bf 4}_5^{\mathfrak{S}}.
$$

2) Type $A_1^4$ components of the set $\Sigma^c(x)$ in the item 4 of Lemma \ref{lok-0-10} for $\delta=1$ and $x=(0,-1,0,0,0)$ are numbered by the points
$$
(-1;-1;-1;\pm1;0.1)\in{\bf 4}_1^{\mathfrak{S}},\quad (-0.5;-1;1;\pm1;0.1)\in{\bf 4}_2^{\mathfrak{S}},\quad (1;-1;-1;0.2;0.1)\in{\bf 4}_3^{\mathfrak{S}},
$$
$$
(1;-1;-1;-1;0.1)\in{\bf 4}_5^{\mathfrak{S}},\quad (0.1;-1;1;-0.1;0.1)\in{\bf 4}_7^{\mathfrak{S}},\quad (1;-1;1;0.2;0.1)\in{\bf 4}_{10}^{\mathfrak{S}}.
$$

\medskip
\lemma\label{lok-0-11} {\it Let $S_4=t_1=0,q_5<0$. Then:

{\rm1)} if $(t_2,S_3)\in\overline{{\bf 2}_{\ast}^0}$ is a point of one of parabolas {\rm(\ref{paraboly})} and $t_2S_3\neq0$, then $f$ has a singularity of type $A_2$ at $x$ and a multisingularity of type $A_2^2$ at $f(x)$; the set $\Sigma^c(x)$ consists of four contractible components; there are two components of each type $A_1^2,A_1^4$; they belong to different strata of $\mathfrak{S}$;

{\rm2)} if $(t_2,S_3)\in\overline{{\bf 6}_{\ast}^0}$ is a point of one of parabolas {\rm(\ref{paraboly})} and $t_2S_3\neq0$, then $f$ has a singularity of type $A_2$ at $x$ and a multisingularity of type $A_2^2A_1^2$ at $f(x)$; the set $\Sigma^c(x)$ consists of four contractible components; there are two components of each type $A_1^4,A_1^6$; they belong to different strata of $\mathfrak{S}$;

{\rm3)} if $(t_2,S_3)$ is the intersection point of parabolas {\rm(\ref{paraboly})}, then $f$ has a singularity of type $A_2$ at $x$ and a multisingularity of type $A_2^3$ at $f(x)$; the set $\Sigma^c(x)$ consists of $2,4,2$ contractible components of types $A_1^2,A_1^4,A_1^6$, respectively; the components of types $A_1^2,A_1^6$ belong to different strata of $\mathfrak{S}$; two components of type $A_1^4$ belong to different strata ${\bf 4}_{10}^{\mathfrak{S}},{\bf 4}_{11}^{\mathfrak{S}}$ and the other two belong to the same stratum ${\bf 4}_{3}^{\mathfrak{S}}$;

{\rm4)} if $(t_2,S_3)$ is a point of the line $S_3=0$ lying outside the segment with ends at vertices of the parabolas {\rm(\ref{paraboly})}, then $f$ has a singularity of type $D_4^{+}$ at $x$ and a multisingularity of type $D_4^{+}A_1^2$ at $f(x)$; the set $\Sigma^c(x)$ consists of eight contractible components; there are four components of each type $A_1^4,A_1^6$; they belong to different strata ${\bf 4}_1^{\mathfrak{S}} - {\bf 4}_3^{\mathfrak{S}},{\bf 4}_{10}^{\mathfrak{S}},{\bf 6}_1^{\mathfrak{S}},{\bf 6}_{3}^{\mathfrak{S}},{\bf 6}_{7}^{\mathfrak{S}},{\bf 6}_{10}^{\mathfrak{S}}$ for $\delta t_2<0$ and to ${\bf 4}_1^{\mathfrak{S}},{\bf 4}_3^{\mathfrak{S}},{\bf 4}_{5}^{\mathfrak{S}},{\bf 4}_{11}^{\mathfrak{S}},{\bf 6}_3^{\mathfrak{S}},{\bf 6}_{4}^{\mathfrak{S}},{\bf 6}_{7}^{\mathfrak{S}},{\bf 6}_{11}^{\mathfrak{S}}$ for $\delta t_2>0$;

{\rm5)} if $(t_2,S_3)$ is a point of the line $S_3=0$ lying inside the segment with ends at vertices of the parabolas {\rm(\ref{paraboly})}, then $f$ has a singularity of type $D_4^{-}$ at $x$ and a multisingularity of type $D_4^{-}A_1^2$ at $f(x)$; the set $\Sigma^c(x)$ consists of nine connected components; there are one component of type $A_1^4$ and eight components of type $A_1^6$; the component of type $A_1^4$ is homotopy equivalent to $S^1$ and belongs to the stratum ${\bf 4}_{3}^{\mathfrak{S}}$; the components of type $A_1^6$ are contractible and belong to different strata ${\bf 6}_1^{\mathfrak{S}},{\bf 6}_3^{\mathfrak{S}} - {\bf 6}_{5}^{\mathfrak{S}},{\bf 6}_{7}^{\mathfrak{S}},{\bf 6}_{8}^{\mathfrak{S}},{\bf 6}_{13}^{\mathfrak{S}},{\bf 6}_{16}^{\mathfrak{S}}$;

{\rm6)} if $(t_2,S_3)$ is the vertex of one of parabolas {\rm(\ref{paraboly})}, then $f$ has a singularity of type $D_5^{\delta}$ at $x$ and a multisingularity of type $D_5^{\delta}A_1$ at $f(x)$; the set $\Sigma^c(x)$ consists of $1,5,10$ connected components of types $A_1^2,A_1^4,A_1^6$, respectively; the component of type $A_1^2$ is contractible and belongs to the stratum ${\bf 2}_{1}^{\mathfrak{S}}$; three components of type $A_1^4$ are contractible and the other two are homotopy equivalent to $S^1$; the contractible components of type $A_1^4$ belong to different strata ${\bf 4}_1^{\mathfrak{S}},{\bf 4}_2^{\mathfrak{S}},{\bf 4}_7^{\mathfrak{S}}$ for $\delta t_2<0$ and to ${\bf 4}_1^{\mathfrak{S}},{\bf 4}_5^{\mathfrak{S}},{\bf 4}_7^{\mathfrak{S}}$ for $\delta t_2>0$; the non-contractible components of type $A_1^4$ belong to different strata ${\bf 4}_3^{\mathfrak{S}},{\bf 4}_{10}^{\mathfrak{S}}$ for $\delta t_2<0$ and to ${\bf 4}_3^{\mathfrak{S}},{\bf 4}_{11}^{\mathfrak{S}}$ for $\delta t_2>0$; the components of type $A_1^6$ are contractible and belong to different strata ${\bf 6}_1^{\mathfrak{S}},{\bf 6}_{3}^{\mathfrak{S}},{\bf 6}_{4}^{\mathfrak{S}},{\bf 6}_{5}^{\mathfrak{S}},{\bf 6}_{7}^{\mathfrak{S}},{\bf 6}_{8}^{\mathfrak{S}},{\bf 6}_{13}^{\mathfrak{S}},{\bf 6}_{16}^{\mathfrak{S}}$ as well as to ${\bf 6}_{10}^{\mathfrak{S}},{\bf 6}_{14}^{\mathfrak{S},+}$ for $\delta t_2<0$ and to ${\bf 6}_{11}^{\mathfrak{S}},{\bf 6}_{14}^{\mathfrak{S},-}$ for $\delta t_2>0$.}

\begin{figure}[h]
\begin{center}
\includegraphics[width=11cm]{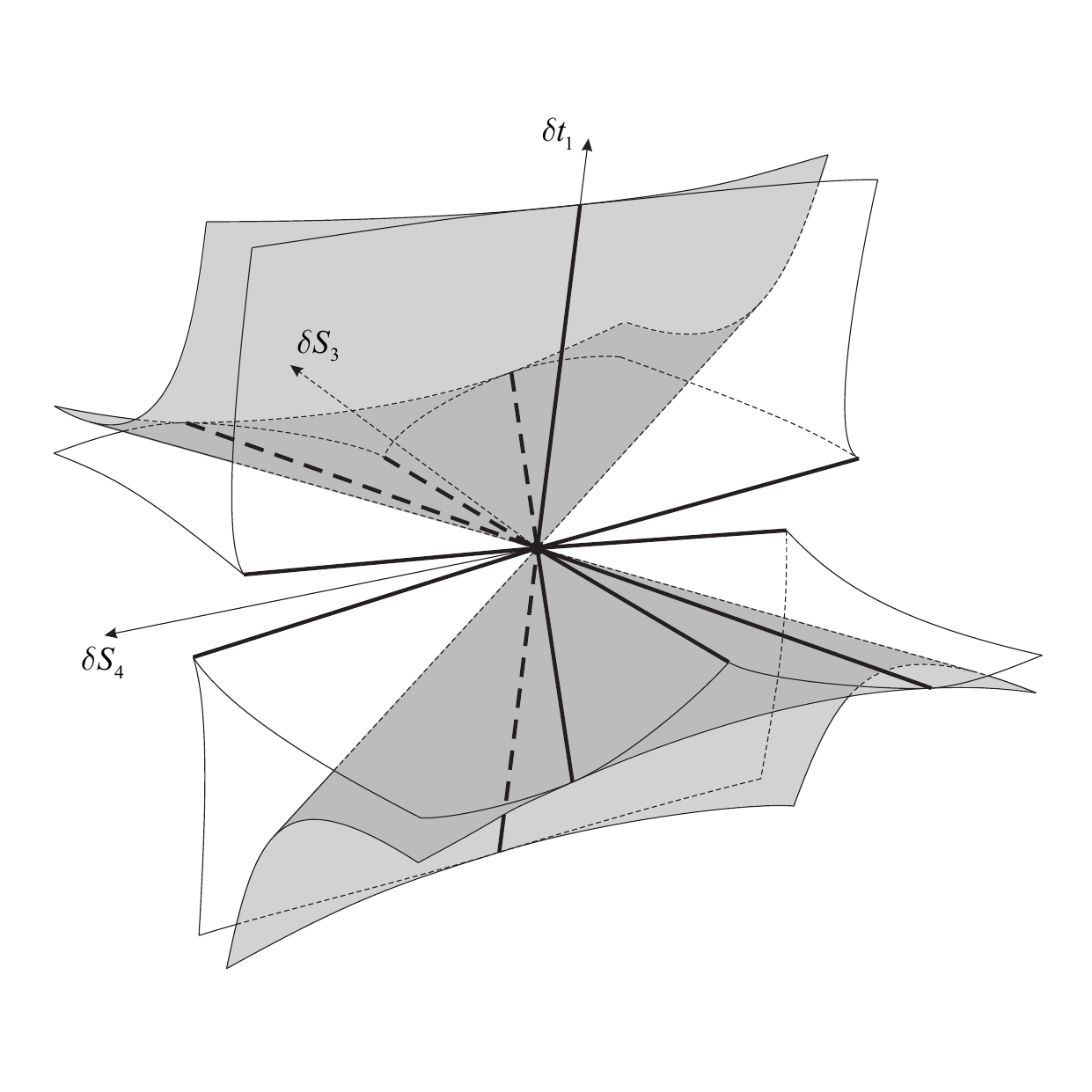}
\caption{A section germ of the skeleton $\Sigma$ by the space $t_2=0,q_5=-3$.}
\label{S1-komponenta}
\end{center}
\end{figure}

\medskip
\remark 1) The type $A_1^4$ non-contractible component of the set $\Sigma^c(x)$ in the item 5 of Lemma \ref{lok-0-11} for $x=(0,0,0,0,-3)$ is homotopy equivalent to the non-contractible connected component of the complement to a germ of the surface in Fig. \ref{S1-komponenta} at the origin. This germ is a section germ of the skeleton $\Sigma$ by the space $t_2=0,q_5=-3$. The semi-transparent surface represents a germ of the cone $6t_1S_3=S_4^2$. The dotted lines will become invisible if $\Sigma$ is made opaque.

2) The description of type $A_1^4$ components in the item 6 of Lemma \ref{lok-0-11} is based on a corollary of the results from \cite{Sed2015}. Namely the manifold $\Xi_{A_1^4}(D_5^{\delta}A_1)$ is homotopy equivalent to $\Xi_{A_1^3}(D_5^{\delta})$. Therefore it has two connected components, one of which is homotopy equivalent to $S^1$ and the other is contractible. The inverse image of the contractible component under the map $f$ has three contractible connected components in a sufficiently small neighbourhood of the point $x$ (at which $f$ has a singularity of type $D_5^{\delta}$). The inverse image of the non-contractible component consists of connected components that are homotopy equivalent to $S^1$. The restriction of $f$ to one of them is a two-sheeted covering since a Lagrangian monosingularity of type $D_5^{\delta}$ is adjacent to a multisingularity of type $D_4^-A_1$. Hence, there is one more non-contractible connected component of the inverse image in a neighbourhood of the point $x$; the restriction of $f$ to this component is a diffeomorphism.

\section{The topology of type $A_1^k$ strata of $\mathfrak{S}$}

The hyperplane $S_4=0$ intersects all strata of the stratification $\mathfrak{S}$ enumerated in Proposition \ref{straty-vse} except ${\bf 6}_{14}^{\mathfrak{S},\pm}$. By the quasi-homogeneity of the map $f$, the strata ${\bf 6}_{14}^{\mathfrak{S},\pm}$ are homotopy equivalent to their intersections with the hyperplanes $\pm\delta S_4=1$, respectively. Thus using Lemma \ref{A6-3-styag}, we get

\medskip
\predlo\label{6-14} {\it The strata ${\bf 6}_{14}^{\mathfrak{S},\pm}$ are contractible.}

\medskip
Let $M$ be a sphere of radius $1$ centered at the origin in the source space $\mathbb{R}^5$ of the map $f$. Then each stratum of $\mathfrak{S}$ containing at least one point  $x=(t_1,t_2,S_3,S_4,q_5)\neq0$ is homotopy equivalent to its intersection with $M$. The intersection of such a stratum with the half-sphere $M\cap\{\pm \delta S_4>0\}$ is diffeomorphic to its intersection with the hyperplane $\pm \delta S_4=1$, respectively. 

\medskip
\predlo\label{2-1} {\it The strata ${\bf 2}_1^{\mathfrak{S}} - {\bf 2}_4^{\mathfrak{S}}, {\bf 4}_1^{\mathfrak{S}}, {\bf 4}_{2}^{\mathfrak{S}}, {\bf 4}_{5}^{\mathfrak{S}}, {\bf 4}_{7}^{\mathfrak{S}}, {\bf 6}_{1}^{\mathfrak{S}}, {\bf 6}_{3}^{\mathfrak{S}} - {\bf 6}_{5}^{\mathfrak{S}}, {\bf 6}_{7}^{\mathfrak{S}}, {\bf 6}_{8}^{\mathfrak{S}},{\bf 6}_{10}^{\mathfrak{S}}, {\bf 6}_{11}^{\mathfrak{S}}, {\bf 6}_{13}^{\mathfrak{S}}, {\bf 6}_{16}^{\mathfrak{S}}$
 are contractible.}

{\sc Proof.} Let $\sigma$ be the intersection of the hyperplane $S_4=0$ with the stratum ${\bf 2}_1^{\mathfrak{S}}$. Then the set $\overline{\sigma}\cap M$ is contractible. Indeed, $\overline{\sigma}$ is homotopy equivalent to the union of the following subsets of $\mathbb{R}^4=\{(t_1,t_2,S_3,q_5)\}$: 
$$
\begin{aligned}
&\{\delta t_1\leq0,t_2=0,S_3=0,q_5\geq0\},\\
&\{\delta t_1\leq0,t_2=0,S_3=4q_5t_1,-2\sqrt{3|t_1|}\leq q_5\leq0\},\\
&\{\delta t_1\leq0,t_2=0,12\delta q_5S_3=-(q_5^2-12\delta t_1)^2,q_5\leq-2\sqrt{3|t_1|}\}.
\end{aligned}
$$
This follows from Fig. \ref{predel3-2chetv}, \ref{S4=t1=0} and Lemma \ref{samoperesechS4}.
The mentioned union is a continuous surface homotopically equivalent to 
$$
\mathcal{D}=\{\delta t_1\leq0,t_2=0,S_3=0\}.
$$ 

Now by Lemmas \ref{lok-0-1}, \ref{lok-0-2}, \ref{lok-0-5} -- \ref{lok-0-8}, \ref{lok-0-9}, \ref{lok-0-10}, \ref{lok-0-11}, it follows that each point of $\overline{\sigma}\cap M$ has a neighbourhood in $M$ such that its intersection with ${\bf 2}_1^{\mathfrak{S}}$ is contractible. But Whitney stratification is locally trivial (see \cite{Mather}). Hence, there is a neighbourhood $U\subset M$ of the set $\overline{\sigma}\cap M$ such that the intersection ${\bf 2}_1^{\mathfrak{S}}\cap U$ is contractible. It remains to note that the intersection ${\bf 2}_1^{\mathfrak{S}}\cap M$ is contracted into $U$ by Lemma \ref{A2-1-styag}.

The contractibility of the strata ${\bf 2}_2^{\mathfrak{S}} - {\bf 2}_4^{\mathfrak{S}},{\bf 4}_1^{\mathfrak{S}},{\bf 4}_{2}^{\mathfrak{S}},{\bf 4}_{5}^{\mathfrak{S}},{\bf 4}_{7}^{\mathfrak{S}},{\bf 6}_1^{\mathfrak{S}},{\bf 6}_{3}^{\mathfrak{S}}-{\bf 6}_{5}^{\mathfrak{S}},{\bf 6}_{7}^{\mathfrak{S}},{\bf 6}_{8}^{\mathfrak{S}},{\bf 6}_{10}^{\mathfrak{S}},{\bf 6}_{11}^{\mathfrak{S}}, {\bf 6}_{13}^{\mathfrak{S}}$, ${\bf 6}_{16}^{\mathfrak{S}}$ is proved similarly. In these cases, the surface $\mathcal{D}$ is determined by the following conditions: 

${\bf 2}_4^{\mathfrak{S}}$) $\delta t_1\geq0,t_2=0,S_3=0$; 

${\bf 2}_2^{\mathfrak{S}},{\bf 2}_3^{\mathfrak{S}},{\bf 4}_1^{\mathfrak{S}}$) $t_1=0,t_2=0,\delta S_3\leq0,q_5\geq0$; 

${\bf 4}_{2}^{\mathfrak{S}},{\bf 4}_{5}^{\mathfrak{S}}$) $t_2=0,S_3=0,q_5\geq0$; 

${\bf 4}_7^{\mathfrak{S}}$) $\delta t_1\geq0,t_2=0,S_3=0,q_5\geq0$; 

${\bf 6}_1^{\mathfrak{S}},{\bf 6}_4^{\mathfrak{S}}$) $\delta t_1\leq0,t_2=0,S_3=0,q_5\leq\sqrt{\frac32|t_1|}$; 

${\bf 6}_3^{\mathfrak{S}}$) $\delta t_1\leq0,t_2=0,S_3=0,q_5\leq2\sqrt{3|t_1|}$; 

${\bf 6}_7^{\mathfrak{S}},{\bf 6}_{10}^{\mathfrak{S}},{\bf 6}_{11}^{\mathfrak{S}}$) $\delta t_1\geq0,t_2=0,S_3=0,q_5\leq2\sqrt{3|t_1|}$; 

${\bf 6}_5^{\mathfrak{S}}$) $\delta t_1\leq0,t_2=0,S_3=0,q_5\leq-2\sqrt{3|t_1|}$; 

${\bf 6}_{8}^{\mathfrak{S}},{\bf 6}_{13}^{\mathfrak{S}},{\bf 6}_{16}^{\mathfrak{S}}$) $\delta t_1\geq0,t_2=0,S_3=0,q_5\leq-2\sqrt{3|t_1|}$. $\Box$

\medskip
\predlo\label{4-3} {\it The strata ${\bf 4}_{3}^{\mathfrak{S}},{\bf 4}_{10}^{\mathfrak{S}},{\bf 4}_{11}^{\mathfrak{S}}$ are homotopy equivalent to $S^1$.}

{\sc Proof.} Let $\sigma$ be the intersection of the hyperplane $S_4=0$ with the stratum ${\bf 4}_{3}^{\mathfrak{S}}$. Then the set $\overline{\sigma}\cap M$ is contractible since $\overline{\sigma}$ is homotopy equivalent to the surface 
$$
\mathcal{D}=\{t_2=0,S_3=0,q_5\leq0\}
$$ 
in $\mathbb{R}^4=\{(t_1,t_2,S_3,q_5)\}$.

By $\mathcal{I}$ denote the subset  of the hyperplane $S_4=0$ given by the following conditions: 
$$
t_1=0,S_3=0,q_5^3+108t_2^2\leq0.
$$ 
The intersection $\mathcal{I}\cap M$ is diffeomorphic to a segment. By Lemmas \ref{lok-0-1} -- \ref{lok-0-8}, \ref{lok-0-10}, \ref{lok-0-11}, it follows that $\mathcal{I}$ lies in $\overline{\sigma}$ and each point $x\in\overline{\sigma}\cap M$ has a neighbourhood in $\mathbb{R}^5$ such that its intersection with the stratum ${\bf 4}_3^{\mathfrak{S}}$ is homotopy equivalent to $S^1$ if $x\in\mathcal{I}$ and consists of contractible connected components if $x\notin\mathcal{I}$. Hence, there is a neighbourhood $U\subset M$ of the set $\overline{\sigma}\cap M$ such that the intersection ${\bf 4}_3^{\mathfrak{S}}\cap U$ is homotopy equivalent to $S^1$. It remains to note that the intersection ${\bf 4}_3^{\mathfrak{S}}\cap M$ is contracted into $U$ by Lemma \ref{A4-3-styag}.

The strata ${\bf 4}_{10}^{\mathfrak{S}}$ and ${\bf 4}_{11}^{\mathfrak{S}}$ are considered similarly. In both cases, the surface $\mathcal{D}$ is determined by the conditions $\delta t_1\geq0,t_2=0,S_3=0,q_5\leq0$. The set $\mathcal{I}$ is given by the following formulas:

${\bf 4}_{10}^{\mathfrak{S}}$) $t_1=0,q_5^3+108t_2^2=0,\delta t_2\leq0,S_3=0,q_5\leq0$;

${\bf 4}_{11}^{\mathfrak{S}}$) $t_1=0,q_5^3+108t_2^2=0,\delta t_2\geq0,S_3=0,q_5\leq0$. $\Box$

\section{Proof of Theorem \ref{th}}

{\sc 1) Connected components of the manifold $(A_1^2)^f$.}

\medskip
The point $x=(t_1,t_2,S_3,S_4,q_5)=(-\delta,0,\delta,0,4)$ belongs to the stratum ${\bf 2}_1^{\mathfrak{S}}$. Its image under the map $f$ has exactly two preimages: $x$ and $\widehat{x}=(\delta,0,17\delta,0,4)$. The point $\widehat{x}$ belongs to the stratum ${\bf 2}_4^{\mathfrak{S}}$. Hence by Lemma \ref{kratnost}, Propositions \ref{straty-vse} and \ref{2-1}, it follows that the manifold of type $A_1^2$ multisingularities of the map $f$ has two connected components
$$
f({\bf 2}_1^{\mathfrak{S}})=f({\bf 2}_4^{\mathfrak{S}}) \quad \mbox{ and } \quad f({\bf 2}_2^{\mathfrak{S}})=f({\bf 2}_3^{\mathfrak{S}}).
$$
Both of them are contractible. The restriction of $f$ to each of the strata ${\bf 2}_1^{\mathfrak{S}} - {\bf 2}_4^{\mathfrak{S}}$ is an embedding.

\medskip
{\sc 2) Connected components of the manifold $(A_1^4)^f$.}
\medskip

By Lemma \ref{kratnost}, Propositions \ref{straty-vse} and \ref{4-3}, it follows that the manifold of type $A_1^4$ multisingularities of the map $f$ has two connected components
$$
f({\bf 4}_{1}^{\mathfrak{S}})=f({\bf 4}_{2}^{\mathfrak{S}})=f({\bf 4}_{5}^{\mathfrak{S}})=f({\bf 4}_{7}^{\mathfrak{S}}) \quad \mbox{ and } \quad f({\bf 4}_{3}^{\mathfrak{S}})=f({\bf 4}_{10}^{\mathfrak{S}})=f({\bf 4}_{11}^{\mathfrak{S}}).
$$
The first one is contractible and the second is homotopy equivalent to $S^1$. The restriction of $f$ to each of the strata ${\bf 4}_1^{\mathfrak{S}},{\bf 4}_2^{\mathfrak{S}},{\bf 4}_5^{\mathfrak{S}},{\bf 4}_7^{\mathfrak{S}}$ is an embedding.

Consider the point $x=(-\delta,0,6\delta,0,-3)$. It belongs to the stratum ${\bf 4}_3^{\mathfrak{S}}$. Its image under the map $f$ has one more preimage in ${\bf 4}_3^{\mathfrak{S}}$, namely $\widehat{x}=(\delta,0,-6\delta,0,-3)$. Hence the restriction of $f$ to the stratum ${\bf 4}_3^{\mathfrak{S}}$ is a two-sheeted covering. The restriction of $f$ to each of the strata ${\bf 4}_{10}^{\mathfrak{S}},{\bf 4}_{11}^{\mathfrak{S}}$ is an embedding.

\medskip
{\sc 3) Connected components of the manifold $(A_1^6)^f$.}
\medskip

The restriction of $f$ to a type $A_1^6$ stratum of $\mathfrak{S}$ preserves the sign of the coordinate $q_5$. Now it follows from Lemma \ref{kratnost}, Propositions \ref{straty-vse}, \ref{6-14} and \ref{2-1} that the manifold of type $A_1^6$ multisingularities of the map $f$ has two connected components
$$
f({\bf 6}_{5}^{\mathfrak{S}})=f({\bf 6}_{8}^{\mathfrak{S}})=f({\bf 6}_{13}^{\mathfrak{S}})=f({\bf 6}_{14}^{\mathfrak{S},-})=f({\bf 6}_{14}^{\mathfrak{S},+})=f({\bf 6}_{16}^{\mathfrak{S}})
$$
(it lies in the half-space $q_5<0$) and
$$
f({\bf 6}_{1}^{\mathfrak{S}})=f({\bf 6}_{3}^{\mathfrak{S}})=f({\bf 6}_{4}^{\mathfrak{S}})=f({\bf 6}_{7}^{\mathfrak{S}})=f({\bf 6}_{10}^{\mathfrak{S}})=f({\bf 6}_{11}^{\mathfrak{S}}).
$$
They are contractible. The restriction of $f$ to each of type $A_1^6$ strata of $\mathfrak{S}$ is an embedding.

\medskip
\remark The image of the circle (\ref{fgr-D4}) under the map $f$ is the curve
\begin{equation*}
\begin{aligned}
q_1=&-\frac1{54}(5+18\cos\varphi+3\sin\varphi),\\
q_2=&\frac1{2916}(234-162\cos\varphi+2703\sin\varphi-126\cos2\varphi-27\sin2\varphi-\sin3\varphi),\\
q_3=&\frac{\delta}{108}(-19+18\cos\varphi-84\sin\varphi+\cos2\varphi),\quad
q_4=\frac{\delta}3(3+\sin\varphi),\quad q_5=-3.
\end{aligned}
\end{equation*}
It determines the generator for the fundamental group of the noncontractible component of the complement to the caustic of the map $f$.

\medskip

\begin{flushleft}

Gubkin University, Moscow, Russia

\bigskip

vdsedykh@gmail.com

sedykh@mccme.ru

\end{flushleft}

\end{document}